\pdfoutput=1 % ensure arXiv direct PDF output

\documentclass[11pt]{amsart}

\usepackage{amssymb,amsthm,enumitem,colonequals,tikz-cd,microtype}
\usepackage[normalem]{ulem} 
\usepackage{mathtools} % For the slashed arrow
\usepackage{graphicx}

\usepackage[osf]{Baskervaldx}
\usepackage[baskervaldx]{newtxmath}
\usepackage[cal=boondoxo]{mathalfa} % mathcal from STIX, unslanted a bit

\usepackage[top=3.75cm, bottom=3cm, left=3.5cm, right=3.5cm]{geometry}

\usepackage{xcolor}
\colorlet{darkblue}{blue!55!black}
\colorlet{darkcyan}{cyan!50!black}
\colorlet{darkgreen}{green!60!black}

\PassOptionsToPackage{hyphens}{url}
\usepackage{hyperref}
\hypersetup{
    colorlinks=true, 
    linkcolor=darkblue,
    urlcolor=darkcyan,
    citecolor=darkgreen,
}

\def\eqref#1{\textcolor{darkblue}{(\ref{#1})}}

\usepackage[nameinlink,capitalize]{cleveref} 
\Crefformat{section}{#2\S#1#3}
\Crefmultiformat{section}{#2\S\S#1#3}{ and~#2#1#3}{, #2#1#3}{, and~#2#1#3}

\usepackage[pagewise]{lineno}
\RequirePackage{tikz-cd}
\RequirePackage{amssymb}
\usetikzlibrary{calc}
\usetikzlibrary{decorations.pathmorphing}
\usetikzlibrary{spath3}
\usetikzlibrary{nfold}

\tikzset{curve/.style={settings={#1},to path={(\tikztostart)
    .. controls ($(\tikztostart)!\pv{pos}!(\tikztotarget)!\pv{height}!270:(\tikztotarget)$)
    and ($(\tikztostart)!1-\pv{pos}!(\tikztotarget)!\pv{height}!270:(\tikztotarget)$)
    .. (\tikztotarget)\tikztonodes}},
    settings/.code={\tikzset{quiver/.cd,#1}
        \def\pv##1{\pgfkeysvalueof{/tikz/quiver/##1}}},
    quiver/.cd,pos/.initial=0.35,height/.initial=0}

\tikzset{between/.style n args={2}{/tikz/execute at end to={
    \tikzset{spath/split at keep middle={current}{#1}{#2}}
}}}

\tikzset{
    tail reversed/.code={\pgfsetarrowsstart{tikzcd to}},
    2tail/.code={\pgfsetarrowsstart{Implies[reversed]}},
    2tail reversed/.code={\pgfsetarrowsstart{Implies}},
    coil/.code={\pgfsetarrowsstart{Glyph[glyph math command=looparrowleft]}},
    coil'/.code={\pgfsetarrowsstart{Glyph[glyph math command=looparrowleft, swap]}},
    multimap/.code={\pgfsetarrowsend{Glyph[glyph math command=multimap]}},
    no body/.style={/tikz/dash pattern=on 0 off 1mm},
}

\let\oldequation\equation
\let\oldendequation\endequation
\renewenvironment{equation}{\linenomathNonumbers\oldequation}{\oldendequation\endlinenomath}
\expandafter\let\expandafter\oldequationstar\csname equation*\endcsname
\expandafter\let\expandafter\oldendequationstar\csname endequation*\endcsname
\renewenvironment{equation*}{\linenomathNonumbers\oldequationstar}{\oldendequationstar\endlinenomath}
\let\oldalign\align
\let\oldendalign\endalign
\renewenvironment{align}{\linenomathNonumbers\oldalign}{\oldendalign\endlinenomath}
\expandafter\let\expandafter\oldalignstar\csname align*\endcsname
\expandafter\let\expandafter\oldendalignstar\csname endalign*\endcsname
\renewenvironment{align*}{\linenomathNonumbers\oldalignstar}{\oldendalignstar\endlinenomath}

\theoremstyle{plain}
\newtheorem{theorem}{Theorem}[section]
\newtheorem{lemma}[theorem]{Lemma}
\newtheorem{corollary}[theorem]{Corollary}
\newtheorem{proposition}[theorem]{Proposition}

\theoremstyle{definition}
\newtheorem{definition}[theorem]{Definition}
\newtheorem{example}[theorem]{Example}
\newtheorem{remark}[theorem]{Remark}
\newtheorem{setup}[theorem]{Setup}

\newtheorem{reminder}[theorem]{Reminder}
\newtheorem{question}[theorem]{Question}
\newtheorem{lemma-definition}[theorem]{Lemma-Definition}
\newtheorem{example-claim}[theorem]{Example-Claim}

\newtheorem*{ack}{Acknowledgments}

\AddToHook{env/conjecture/begin}{\crefalias{theorem}{conjecture}}
\AddToHook{env/lemma/begin}{\crefalias{theorem}{lemma}}
\AddToHook{env/corollary/begin}{\crefalias{theorem}{corollary}}
\AddToHook{env/proposition/begin}{\crefalias{theorem}{proposition}}
\AddToHook{env/definition/begin}{\crefalias{theorem}{definition}}
\AddToHook{env/remark/begin}{\crefalias{theorem}{remark}}
\AddToHook{env/setup/begin}{\crefalias{theorem}{setup}}
\AddToHook{env/example/begin}{\crefalias{theorem}{example}}
\AddToHook{env/conjecture/begin}{\crefalias{theorem}{conjecture}}
\AddToHook{env/notation/begin}{\crefalias{theorem}{notation}}
\AddToHook{env/reminder/begin}{\crefalias{theorem}{reminder}}
\AddToHook{env/question/begin}{\crefalias{theorem}{question}}
\AddToHook{env/lemma-definition/begin}{\crefalias{theorem}{lemma-definition}}
\AddToHook{env/example-claim/begin}{\crefalias{theorem}{example-claim}}

\makeatletter
\def\slashedarrowfill@#1#2#3#4#5{%
  $\m@th\thickmuskip0mu\medmuskip\thickmuskip\thinmuskip\thickmuskip
   \relax#5#1\mkern-7mu%
   \cleaders\hbox{$#5\mkern-2mu#2\mkern-2mu$}\hfill
   \mathclap{#3}\mathclap{#2}%
   \cleaders\hbox{$#5\mkern-2mu#2\mkern-2mu$}\hfill
   \mkern-7mu#4$%
}
\def\rightslashedarrowfill@{%
  \slashedarrowfill@\relbar\relbar\mapstochar\rightarrow}
\newcommand\xslashedrightarrow[2][]{%
  \ext@arrow 0055{\rightslashedarrowfill@}{#1}{#2}}
\makeatother

\def\slashedrightarrow{\xslashedrightarrow{}}

\numberwithin{equation}{section}
\numberwithin{theorem}{section}

\title[Relative Fourier--Mukai theory for algebraic spaces]{Relative Fourier--Mukai theory for algebraic spaces}

\author[E.~Guisado Villalgordo]{El\'{i}as Guisado Villalgordo}
\address{E.~Guisado Villalgordo,
BCAM -- Basque Center for Applied Mathematics, Mazarredo 14, 48009
Bilbao, Basque Country -- Spain}
\email{eguisado@bcamath.org}

\author[P.~Lank]{Pat Lank}
\address{P.~Lank,
Dipartimento di Matematica `F. Enriques', Universit\`{a} degli Studi di Milano, Via Cesare
Saldini 50, 20133 Milano, Italy}
\email{lankp.math@gmail.com}

\date{\today}


\begin{document}

\begin{abstract}
	We study the variation of Fourier--Mukai equivalences in proper flat families of Noetherian algebraic spaces. 
	For bounded pseudocoherent kernels, we prove that the equivalence locus is open and commutes with arbitrary Noetherian base change.
	To do so, we develop a theory of relative perfectness and relative quasi-perfectness for algebraic spaces.
	This includes a fiberwise criteria for integral transforms with nonperfect kernels. 
	As an application, we construct an explicit autoequivalence for elliptic fibrations.
\end{abstract}

\maketitle

\tableofcontents

%%%%%%%%%%%%%%%%%%%%%%%%%%%%%%%%%%%
\setcounter{part}{-1}
\part{Introduction and preliminaries}
\label{sec:part0}
%%%%%%%%%%%%%%%%%%%%%%%%%%%%%%%%%%%

%%%%%%%%%%%%%%%%%%%%%%%%%%%%%%%%%%%
\section{Introduction}
\label{sec:intro}

%%%%%%%%%%%%%%%%%%%%%%%%%%%%%%%%%%%
\subsection{Motivation}
\label{sec:intro_motivation}
%%%%%%%%%%%%%%%%%%%%%%%%%%%%%%%%%%%

Algebraic spaces sit at the intersection for those interested in varieties, quotients, and moduli theory.
Derived categories homologically describe these geometric objects. 
Fourier--Mukai theory, which originates in \cite{Mukai:1981}, compares shapes through their derived categories, and its relative version compares families \cite{Ruiperez/Hernandez/Martin/SanchodeSalas:2009}.

This article develops a relative Fourier--Mukai theory for proper flat families of Noetherian algebraic spaces.
It allows for singular fibers and nonperfect kernels. 
The main question is whether the fibers on which an integral transform is an equivalence form an open locus, and how this locus behaves under base change.
Our arguments use the $2$-category of triangulated categories.
\Cref{app:higher_categories} compares our results with higher-categorical approaches, highlighting their overlap and differences.
See \Cref{cmp:ex:unit,cmp:rem:enhanced}.

We fix notation.
Let $f_i\colon Y_i \to S$ be proper flat morphisms of Noetherian algebraic spaces for $i=1,2$.
Consider the fibered square
\begin{displaymath}
	% https://q.uiver.app/#q=WzAsNCxbMSwxLCJcXG1hdGhjYWx7U30uIl0sWzAsMSwiXFxtYXRoY2Fse1l9XzEiXSxbMSwwLCJcXG1hdGhjYWx7WX1fMiJdLFswLDAsIlxcbWF0aGNhbHtZfV8xXFx0aW1lc197XFxtYXRoY2Fse1N9fSBcXG1hdGhjYWx7WX1fMiJdLFsxLDAsImZfMSIsMl0sWzIsMCwiZl8yIl0sWzMsMiwiZl5cXHByaW1lXzEiLDAseyJsYWJlbF9wb3NpdGlvbiI6MzB9XSxbMywxLCJmXlxccHJpbWVfMiIsMl1d
	\begin{tikzcd}
		{Y_1\times_{S}Y_2} & {Y_2} \\
		{Y_1} & {S.}
		\arrow["{f_1^\prime}", from=1-1, to=1-2]
		\arrow["{f_2^\prime}"', from=1-1, to=2-1]
		\arrow["{f_2}", from=1-2, to=2-2]
		\arrow["{f_1}"', from=2-1, to=2-2]
	\end{tikzcd}
\end{displaymath}
For any $K\in D_{\operatorname{qc}}(Y_1\times_S Y_2)$, the \textbf{integral transform} with kernel $K$ is the functor
\begin{displaymath}
	\Phi_K\colon D_{\operatorname{qc}}(Y_1)\to D_{\operatorname{qc}}(Y_2),
	\qquad A\mapsto \mathbf{R}(f_1^\prime)_\ast(\mathbf{L}(f_2^\prime)^\ast A\otimes^{\mathbf{L}}K).
\end{displaymath}
We say that $Y_1$ and $Y_2$ are \textbf{Fourier--Mukai $S$-partners} if there exists $K\in D^b_{\operatorname{coh}}(Y_1\times_S Y_2)$ such that $\Phi_K$ restricts to an equivalence $D^b_{\operatorname{coh}}(Y_1)\to D^b_{\operatorname{coh}}(Y_2)$.
Let $p\in |S|$.
Choose a representative $t\colon \operatorname{Spec}(k)\to S$ of $p$.
Denote by $t^\prime\colon Y_1\times_S Y_2\times_S \operatorname{Spec}(k)\to Y_1\times_S Y_2$ the canonical morphism.
The \textbf{fiber} of $\Phi_K$ at $t$ is the integral transform
\begin{displaymath}
	\Phi_{\mathbf{L}(t^\prime)^\ast K}\colon D_{\operatorname{qc}}(Y_1\times_S \operatorname{Spec}(k))\to D_{\operatorname{qc}}(Y_2\times_S \operatorname{Spec}(k)).
\end{displaymath}

\begin{question}
	\label{question:main}
	Let $K\in D^b_{\operatorname{coh}}(Y_1\times_S Y_2)$. 
	How do full faithfulness and equivalences of the fibers of $\Phi_K$ vary in $|S|$?
\end{question}

For $K\in D^b_{\operatorname{coh}}(Y_1\times_S Y_2)$, write $\operatorname{FM}(K)$ for the set of points $p\in |S|$ such that $\Phi_{\mathbf{L}(t^\prime)^\ast K}$ is an equivalence on $D_{\operatorname{qc}}$ for some representative $t$ of $p$. 
Similarly, $\operatorname{fm}(K)$ denotes the corresponding full faithfulness locus.
See \Cref{def:fm_locus}. 
If $K$ is relatively perfect over each $Y_i$, then the choice of representative does not matter.
See \Cref{lem:independence_for_rep_with_ff_or_eq}.

Our main result is the following.

\begin{theorem}
	\label{thm:fm_locus}
	Consider proper flat morphisms $f_1\colon Y_1 \to S$ and $f_2\colon Y_2 \to S$ of Noetherian algebraic spaces.
	Let $K\in D_{\operatorname{coh}}^b(Y_1 \times_{S}  Y_2 )$.
	Then:
	\begin{enumerate}
		\item \label{thm:fm_locus1} $\operatorname{FM}(K)$ is an open subset of $|S|$
		\item \label{thm:fm_locus2} $\operatorname{FM}(K)$ commutes with arbitrary Noetherian base change.
	\end{enumerate}
	If $K$ is relatively perfect over both $Y_i$, then a similar statement holds for $\operatorname{fm}(K)$.
\end{theorem}

For the equivalence locus, no perfectness or relative perfectness hypothesis on $K$ is required. 
The base change assertion means that for every morphism $t\colon T\to S$ of Noetherian algebraic spaces, with induced map $t^\prime$ on the product, $t^{-1}(\operatorname{FM}(K))=\operatorname{FM}(\mathbf{L}(t^\prime)^\ast K)$. 
See \Cref{thm:fm_locus_weak_approx}. 
For relatively perfect kernels, analogs with $\operatorname{fm}(K)$ follow from \Cref{thm:openness_for_equivalence_fullfaithful_dbcoh,cor:deform}.

These base change results show that if $Y_1$ and $Y_2$ are Fourier--Mukai $S$-partners, then $Y_1\times_S T$ and $Y_2\times_S T$ are Fourier--Mukai $T$-partners for every Noetherian algebraic space $T$ over $S$. 
No flatness assumption on $T\to S$ is needed.
See \Cref{introthm:base_change,introthm:descent_ascent}.
By \Cref{ex:fm_locus_counterexample}, boundedness cannot be dropped.

%%%%%%%%%%%%%%%%%%%%%%%%%%%%%%%%%%%
\subsection{Background and relation to earlier work}
\label{sec:intro_background}
\label{sec:intro_novelty}
%%%%%%%%%%%%%%%%%%%%%%%%%%%%%%%%%%%

Orlov showed that an equivalence given by an integral transform between smooth projective varieties persists under field extension \cite[Lemma 2.12]{Orlov:2002}. 
The kernels in Orlov's work are perfect complexes.

More generally, in the case of perfect kernels, Hall--Priver prove openness of full faithfulness under additional hypotheses on algebraic stacks \cite[Corollary 5.9]{Hall/Priver:2024} (i.e.\ proper tame morphisms with Gorenstein fibers). 
Loc.\ cit.\ bootstraps \cite{Anno/Logvinenko:2012}. 
However, the work of Anno--Logvinenko do not extend directly to nonperfect kernels.
Indeed, \cite[Lemma 2.3]{Anno/Logvinenko:2012} is used in the proof of \cite[Theorem 3.1]{Anno/Logvinenko:2012}, and the proof of \cite[Appendix A, Eq.\ (A.18)]{Anno/Logvinenko:2012} explicitly uses perfectness of the kernel.

Nonperfect kernels are inevitable. 
In particular, for proper flat families with singular fibers, every kernel of a Fourier--Mukai partnership cannot be perfect \cite{GuisadoVillalgordo/Lank:2026b}.
In fact, smoothness of a morphism is completely characterized by perfectness of the kernels. 

Criteria for integral transforms on singular varieties were developed by Hern\'{a}ndez--L\'{o}pez--Sancho \cite{Ruiperez/Martin/deSalas:2007,Ruiperez/Hernandez/Martin/SanchodeSalas:2009}. 
The proof of \cite[Proposition 2.15]{Ruiperez/Hernandez/Martin/SanchodeSalas:2009} is incomplete.
It was corrected in our scheme-theoretic predecessor \cite{GuisadoVillalgordo/Lank/ManaliRahul/Pavic:2025}.
Loc.\ cit.\ proves fiberwise criteria for proper flat families over separated Noetherian schemes.
The calculus of mates was leveraged to control adjunctions under base change. 

We further build upon \cite{GuisadoVillalgordo/Lank/ManaliRahul/Pavic:2025} and work with algebraic spaces.
Our new results concern further base change and openness properties.
The manuscript provides a full theory of relative Fourier--Mukai theory that works for arbitrary Noetherian algebraic spaces.
It does not require separatedness of the base, and establishes new observations.

%%%%%%%%%%%%%%%%%%%%%%%%%%%%%%%%%%%
\subsection{Relative perfectness}
\label{sec:intro_relative_perfectness}
%%%%%%%%%%%%%%%%%%%%%%%%%%%%%%%%%%%

The proof of \Cref{thm:fm_locus} starts with relative perfectness.
This notion was first developed for schemes by Illusie, and later in \cite{Lieblich:2006,AlonsoTarrio/JeremiasLopez/SanchodeSalas:2023,StacksProject}.
We develop these for algebraic spaces in two ways.
The first notion is motivated by \cite{AlonsoTarrio/JeremiasLopez/SanchodeSalas:2023}, and the second by \cite[\S 3]{Ballard:2009}.

For a morphism $f\colon Y\to X$ of quasi-compact quasi-separated algebraic spaces, an object $E\in D_{\operatorname{qc}}(Y)$ is \textbf{$f$-perfect} if
\begin{displaymath}
	E\otimes^{\mathbf{L}}\mathbf{L}f^\ast A\in D^b_{\operatorname{qc}}(Y)
	\quad\textrm{for every }A\in D^b_{\operatorname{qc}}(X),
\end{displaymath}
and \textbf{$f$-quasi-perfect} if
\begin{displaymath}
	\mathbf{R}f_\ast(E\otimes^{\mathbf{L}}P)\in\operatorname{Perf}(X)
	\quad\textrm{for every }P\in\operatorname{Perf}(Y).
\end{displaymath}
See \Cref{def:relative_perf_stacks,def:quasi_perf}. 
These differ in general (see \Cref{ex:relative_perfect_and_quasi_do_not_coincide}), but coincide for pseudocoherent complexes under a proper morphism.

\begin{theorem}
	\label{thm:stacky_all_coincide}
	Let $f\colon Y\to X$ be a proper morphism of Noetherian algebraic spaces.
	For any $E\in D_{\operatorname{qc}}(Y)$ pseudocoherent, the following are equivalent:
	\begin{enumerate}
		\item \label{cor:perfectness_relative_via_coherence1} $E$ is $f$-perfect
		\item \label{cor:perfectness_relative_via_coherence2} $E\otimes^{\mathbf{L}}\mathbf{L}f^\ast A\in D^b_{\operatorname{coh}}(Y)$ for all $A\in D^b_{\operatorname{coh}}(X)$
		\item \label{cor:perfectness_relative_via_coherence3} $E$ is $f$-quasi-perfect.
	\end{enumerate}
\end{theorem}

This extends \cite[Theorem 4.3]{AlonsoTarrio/JeremiasLopez/SanchodeSalas:2023} to algebraic spaces.
See \Cref{sec:properness}. 
Relative perfectness are \'{e}tale local on both source and target (\Cref{thm:relative_perf_iff_smooth_locally}). 
For $Y_1\times_S Y_2$, we say a kernel is relatively perfect over $Y_i$ if it is perfect relative to the projection onto $Y_i$.

Integral transforms preserving bounded pseudocoherence or perfectness are characterized as follows.

\begin{corollary}
	\label{cor:preservation}
	Let $S$ be a Noetherian algebraic space.
	Consider proper morphisms of algebraic spaces $f_1\colon Y_1 \to S$ and $f_2\colon Y_2 \to S$.
	Denote by $p_i\colon Y_1 \times_{S} Y_2 \to Y_i$ the natural projections. For any $E\in D^-_{\operatorname{coh}}(Y_1 \times_{S} Y_2)$, the following are equivalent:
	\begin{enumerate}
		\item $E$ is $p_1$-perfect (resp.\ $p_2$-perfect)
		\item $\Phi_E (D^b_{\operatorname{coh}}(Y_1))\subseteq D^b_{\operatorname{coh}}(Y_2)$ (resp.\ $\Phi_E (\operatorname{Perf}(Y_1))\subseteq \operatorname{Perf}(Y_2)$).
	\end{enumerate}
\end{corollary}

See \Cref{sec:properness}. 
In the bounded pseudocoherent case, such kernels enjoy many properties.

\begin{theorem}
	\label{introthm:base_change}
	Let $f_i\colon Y_i \to S$ be proper flat morphisms of Noetherian algebraic spaces.
	For any $K\in D^b_{\operatorname{coh}}(Y_1\times_S Y_2)$, the following are equivalent:
	\begin{enumerate}
		\item \label{introthm:base_change1} $K$ is relatively perfect over $Y_1$
		\item \label{introthm:base_change2} for every closed point $p\in |S|$ there exists a representative $t\colon \operatorname{Spec}(k)\to S$ of $p$ such that $\mathbf{L}(t^\prime)^\ast K$ is relatively perfect over $Y_1\times_S \operatorname{Spec}(k)$
		\item \label{introthm:base_change3} there exists a faithfully flat morphism $T\to S$ of Noetherian algebraic spaces such that the derived pullback of $K$ to $Y_1\times_S Y_2\times_S T$ is relatively perfect over $Y_1\times_S T$.
	\end{enumerate}
	The analogous statement holds for $Y_2$ in place of $Y_1$.
	Moreover, if $K$ is relatively perfect over each $Y_i$, then for every morphism $T\to S$ of Noetherian algebraic spaces, the derived pullback of $K$ to $Y_1\times_S Y_2\times_S T$ is relatively perfect over each $Y_i\times_S T$.
\end{theorem}

See \Cref{prop:reflecting_bounded_pseudocoherence_perfectness,thm:bounded_pseudocoherence_perfectness_faithfully_flat_affine,cor:noetherian_base_change_for_relative_perfectness}. 
Boundedness is essential (\Cref{ex:relative_perf_counterexample}).
See \Cref{rem:comparison_stacks_project}.
In particular, base change along arbitrary morphisms is new (\Cref{cor:noetherian_base_change_for_relative_perfectness}), and does not appear in \cite{GuisadoVillalgordo/Lank/ManaliRahul/Pavic:2025}.

%%%%%%%%%%%%%%%%%%%%%%%%%%%%%%%%%%%
\subsection{Fiberwise criteria}
\label{sec:intro_fiberwise}
%%%%%%%%%%%%%%%%%%%%%%%%%%%%%%%%%%%

The previous results allow us to compare integral transforms with their fibers.
This applies if the kernels are not perfect. 
Full faithfulness and equivalence can be tested at closed points.

\begin{theorem}
	\label{introthm:descent_ascent}
	Let $S$ be a Noetherian algebraic space.
	Consider proper flat morphisms of algebraic spaces $f_1\colon Y_1 \to S$ and $f_2\colon Y_2 \to S$.
	Suppose $K\in D^b_{\operatorname{coh}}(Y_1 \times_{S}  Y_2 )$ is relatively perfect over each $Y_i$.
	Then the following are equivalent:
	\begin{enumerate}
		\item $\Phi_{K}$ is fully faithful (resp.\ an equivalence) on $D^b_{\operatorname{coh}}$
		\item $\Phi_{\mathbf{L}(t^\prime)^\ast K}$ is fully faithful (resp.\ an equivalence) on $D^b_{\operatorname{coh}}$ for every $t\colon \operatorname{Spec}(k)\to S$ of finite type that represents any closed point $p\in |S|$
		\item for any closed point $p\in |S|$ there exists a representative $t\colon \operatorname{Spec}(k)\to S$ of $p$ such that $\Phi_{\mathbf{L}(t^\prime)^\ast K}$ is fully faithful (resp.\ an equivalence) on $D^b_{\operatorname{coh}}$.
		%%NOTE: $t$ need not be of finite type for this condition
	\end{enumerate}
\end{theorem}

See \Cref{thm:descent_ascent}. 
This extends the scheme-theoretic criterion of \cite[Theorem 6.11]{GuisadoVillalgordo/Lank/ManaliRahul/Pavic:2025}. 
The compatibility of units and counits with base change is established in \Cref{prop:pullbacks_is_a_morphism_of_adjoint_integral_transforms}.
It is not a purely formal consequence of base change for the functors. 
\Cref{introthm:descent_ascent,introthm:base_change} and \Cref{cor:noetherian_base_change_for_relative_perfectness} are key inputs for \Cref{thm:fm_locus}.

%%%%%%%%%%%%%%%%%%%%%%%%%%%%%%%%%%%
\subsection{An autoequivalence}
\label{sec:intro_applications}
\label{sec:intro_applications_autoequivalence}
%%%%%%%%%%%%%%%%%%%%%%%%%%%%%%%%%%%

As an application, we construct an autoequivalence of an elliptic fibration. 
Here we allow singular fibers. 
Recall an elliptic $X$-fibration is a proper Gorenstein morphism $f\colon Y\to X$ of Noetherian algebraic spaces with geometrically integral genus-one fibers. See \Cref{rem:elliptic_fibration}.

\begin{theorem}
	\label{thm:elliptic_fibration_autoequivalence}
	Let $f\colon Y \to X$ be an elliptic $X$-fibration.
	Then $\Phi_{\mathcal{I}_{\Delta_f}}$ induces an autoequivalence of $D^b_{\operatorname{coh}}(Y)$ where $\mathcal{I}_{\Delta_f}$ is the ideal sheaf of the diagonal $\Delta_f\colon Y \to Y \times_X Y$.
\end{theorem}

This extends \cite[Proposition 2.16]{Ruiperez/Hernandez/Martin/SanchodeSalas:2009} from varieties over an algebraically closed field to proper algebraic spaces.
See \Cref{sec:autoequivalence}.
Combined with \Cref{thm:fm_locus}, this has arithmetic consequences.

\begin{example}
	\label{ex:generic_elliptic_fibration}
	Let $S$ be a Noetherian Dedekind scheme with generic point $\eta$.
	For any geometrically connected elliptic curve $C$ over $\kappa(\eta)$, there exists a proper flat morphism $f\colon Y \to S$ from a scheme of Krull dimension two satisfying $Y\times_S \operatorname{Spec}(\kappa(\eta))\cong C$ \cite[\S 10.1.1]{Liu:2002}; see e.g.\ Proposition 1.8 and Remark 1.11 of loc.\ cit.
	Denote by $\mathcal{I}_{\Delta_f}$ the ideal sheaf of the diagonal $\Delta_f\colon Y \to Y\times_S Y$.
	By \Cref{thm:elliptic_fibration_autoequivalence,prop:derived_pullback_ideal_sheaf},
	the fiber of $\Phi_{\mathcal{I}_{\Delta_f}}$ at $\eta$ induces an autoequivalence.
	Hence, by \Cref{thm:fm_locus}, $\operatorname{FM}(\mathcal{I}_{\Delta_f})$ is an open subset of $|S|$ containing $\eta$, and so it contains all but finitely many closed points of $S$.
	For instance, if $S=\operatorname{Spec}(\mathbb{Z})$, then the fiber of $\Phi_{\mathcal{I}_{\Delta_f}}$ is an autoequivalence for all but finitely many reductions modulo primes.
\end{example}

%%%%%%%%%%%%%%%%%%%%%%%%%%%%%%%%%%%
\subsection{Methodology and related approaches}
\label{sec:intro_methodology}
%%%%%%%%%%%%%%%%%%%%%%%%%%%%%%%%%%%

The proof of relative perfectness uses \'{e}tale localization and compact generators. 
For the fiberwise criteria, we describe the adjoints of integral transforms through Grothendieck duality, extending \cite{Rizzardo:2017} and using \cite{Neeman:2023}. 

The calculus of mates \cite{Kelly/Street:1974} identifies the base changes of the units and counits with those of the fiberwise adjunctions (\Cref{prop:pullbacks_is_a_morphism_of_adjoint_integral_transforms,rem:ff_descent}). 
Mates are also used by Bergh--Schn\"{u}rer to study semiorthogonal decompositions under base change \cite{Bergh/Schnurer:2020}.
However, in loc.\ cit.\ perfectness of the kernels are used, see e.g.\ the proof of \cite[Proposition 4.9]{Bergh/Schnurer:2020}.

Openness follows by expressing the failure of full faithfulness or equivalence through supports of cones of these transformations (\Cref{thm:openness_for_equivalence_fullfaithful_dbcoh}).
This further requires restricting to the open locus of relative perfectness (\Cref{prop:relatively_perfect_locus,thm:fm_locus_weak_approx}).
Related questions for algebraic stacks are studied in \cite{GuisadoVillalgordo/Lank:2026a}.

%%%%%%%%%%%%%%%%%%%%%%%%%%%%%%%%%%%
\subsection{Organization}
\label{sec:intro_organization}
%%%%%%%%%%%%%%%%%%%%%%%%%%%%%%%%%%%

\Cref{sec:part0} contains the introduction and preliminaries, including duality for algebraic spaces. 
\Cref{sec:part1} develops relative perfectness and its behavior under base change. 
\Cref{sec:part2} proves the fiberwise criteria using mates, derived reflexivity, and base change for adjoints. 
\Cref{sec:part3} treats elliptic fibrations and openness of Fourier--Mukai loci. 
The first two appendices review preferred equivalence classes and the $2$-category of triangulated categories.
The third compares our results with enhanced approaches.

\begin{ack}
	The authors thank Rom\'{a}n Ahumada Gialanella, Matthew R.\ Ballard, Arend Bayer, Ana Cristina L{\'{o}}pez Mart{\'{\i}}n, Jack Hall, Andres Fernandez Herrero, Emilio Franco, Ian Grojnowski, Kabeer Manali Rahul,  Nebojsa Pavic, Fei Peng, Fernando Sancho de Salas, Dylan Spence, and Paolo Stellari for discussions.
	Lank was supported by the ERC Advanced Grant 101095900-TriCatApp.
	Guisado Villalgordo was supported under the grant PRE2022-104796 by the Spanish Ministry of Science and Innovation.
\end{ack}

%%%%%%%%%%%%%%%%%%%%%%%%%%%%%%%%%%%%%
\paragraph*{\textbf{AI declaration}}
%%%%%%%%%%%%%%%%%%%%%%%%%%%%%%%%%%%%%

The authors used large language models (Claude Opus 5.5 and GPT Astra 6.1) to assist with proofreading (e.g.\ grammar, notation, and exposition) and checking arguments (e.g.\ identifying possible gaps or errors). 
The main results and methods are due to the authors.
\Cref{app:higher_categories} was developed through discussions with these models, following the authors' proposal to adapt the construction of \cite[\S 2.1]{Canonaco/Stellari:2012}. 
These discussions contributed to the example and the comparison between enhanced and triangulated integral transforms.
The authors take full responsibility for the contents of the manuscript.

%%%%%%%%%%%%%%%%%%%%%%%%%%%%%%%%%%%

\section{Preliminaries}
\label{sec:preliminaries}
%%%%%%%%%%%%%%%%%%%%%%%%%%%%%%%%%%%%%

By a `natural' morphism, we mean something which is functorial.
By a `canonical' morphism, we mean something which is constructed and might not necessarily be functorial.

%%%%%%%%%%%%%%%%%%%%%%%%%%%%%%%%%%%%%
\subsection{Generation}
\label{sec:preliminaries_generation}
%%%%%%%%%%%%%%%%%%%%%%%%%%%%%%%%%%%%%

We discuss a form of generation for triangulated categories.
See \cite{Bondal/VandenBergh:2003} for details.
Let $\mathcal{T}$ be a triangulated category with shift functor $[1]\colon \mathcal{T} \to \mathcal{T}$.
Consider a subcategory $S \subseteq \mathcal{T}$.
A triangulated subcategory of $\mathcal{T}$ is called \textbf{thick} if it is closed under direct summands.

Denote by $\langle S \rangle$ the smallest thick subcategory of $\mathcal{T}$ containing $S$.
If $S$ consists of a single object $G$, we write $\langle S \rangle := \langle G \rangle$.
Set $\operatorname{add}(S)$ to be the smallest strictly full subcategory of $\mathcal{T}$ containing $S$ that is closed under shifts, finite coproducts, and direct summands.
Inductively, let $\langle S \rangle_0$ consist of all zero objects in $\mathcal{T}$, $\langle S \rangle_1 := \operatorname{add}(S)$, and
\begin{displaymath}
	\langle S \rangle_n := \operatorname{add} \{ \operatorname{cone}(\phi) : \phi \in \operatorname{Hom}_{\mathcal{T}} (\langle S \rangle_{n-1}, \langle S \rangle_1) \}.
\end{displaymath}
It can be checked that $\langle S \rangle = \bigcup_{n=0}^\infty \langle S \rangle_n$.
We say $E$ is \textbf{finitely built by $S$} if $E\in \langle S\rangle$.

\begin{example}
	\label{ex:regular_ring}
	Let $X=\operatorname{Spec}(R)$ where $R$ is a regular ring.
	Then $\langle \mathcal{O}_X \rangle = D^b_{\operatorname{coh}}(X)$.
	More generally, let $X$ be an affine Noetherian scheme and $Z\subseteq X$ be closed.
	If $P\in \operatorname{Perf}(X)$ satisfies $\operatorname{Supp}(P)=Z$, then $\langle P \rangle = \operatorname{Perf}(X)\cap D^b_{\operatorname{coh},Z}(X)$. See \cite[Lemma 1.2]{Neeman:1992}.
\end{example}

If $\mathcal{T}$ admits small coproducts, then the collection of compact objects in $\mathcal{T}$ is denoted by $\mathcal{T}^c$.
These form a triangulated subcategory of $\mathcal{T}$.
We say that $\mathcal{T}$ is \textbf{compactly generated} if it coincides with the smallest triangulated subcategory of $\mathcal{T}$ containing $\mathcal{T}^c$ and closed under small coproducts.
Equivalently, $\mathcal{T}$ is compactly generated if, for any $E \in \mathcal{T}$ satisfying $\operatorname{Hom}(P, E) = 0$ for all $P \in \mathcal{T}^c$, one has $E \cong 0$ \cite[Lemma 2.2.1]{Schwede/Shipley:2003}.
Note that classical generators for $\mathcal{T}^c$ coincide with compact generators for $\mathcal{T}$ \cite[\href{https://stacks.math.columbia.edu/tag/09SR}{Tag 09SR}]{StacksProject}.

\begin{example}
	Let $X$ be a quasi-compact quasi-separated scheme.
	Then $D_{\operatorname{qc}}(X)^c=\operatorname{Perf}(X)$, and moreover, $\operatorname{Perf}(X)$ admits a classical generator.
	See \cite[Theorem 3.1.1]{Bondal/VandenBergh:2003}.
	More generally, let $Z\subseteq X$ be closed with quasi-compact complement.
	By \cite[Theorem 6.8]{Rouquier:2008}, $D_{\operatorname{qc},Z}(X)$ is compactly generated by a single object.
	Hence, $\operatorname{Perf}(X)\cap D_{\operatorname{qc},Z}(X)$ admits a classical generator.
	See also \Cref{sec:prelim_stacks_perfects}.
\end{example}

We still keep the hypothesis that $\mathcal{T}$ has small coproducts.
Then a triangulated subcategory of $\mathcal{T}$ is \textbf{localizing} when it's closed under small coproducts.
Given a set $\mathcal{S}\subseteq \mathcal{T}$, set $\overline{\langle \mathcal{S} \rangle}$ to be the smallest localizing subcategory of $\mathcal{T}$ containing $\mathcal{S}$.
A localizing subcategory is thick,
see \cite[Remark 1.4 and Proposition 3.2]{Bokstedt/Neeman:1993} or \cite[Example 3.5]{Hall/Rydh:2017}.

%%%%%%%%%%%%%%%%%%%%%%%%%%%%%%%%%%%%%%%%%%%
\subsubsection{Generation criteria for functors}
%%%%%%%%%%%%%%%%%%%%%%%%%%%%%%%%%%%
We say that a triangulated functor $F\colon \mathcal{A}\to\mathcal{B}$ between triangulated categories is \textbf{essentially dense} if for every object $B\in\mathcal{B}$ there is an object $A\in\mathcal{A}$ such that
$\langle F(A)\rangle_1\ni B$
(this notion also makes sense for an additive functor between additive categories).
The following is known but we record it for later use;
see \Cref{app:triang_is_2-cat} for the definitions.

\begin{lemma}\hfill
	\begin{enumerate}[
		label=(\arabic*),
		ref=\thelemma(\arabic*)
		]
		\item \label[lemma]{lem:subcat_wher_nat_is_iso_is_closed_under_coprod_and_cofacts}
		Let $\eta\colon F\Rightarrow G$ be a natural transformation of functors $\mathcal{A}\overset{F}{\underset{G}{\rightrightarrows}}\mathcal{B}$ of categories.
		Let $I$ be a set.
		Assume $\mathcal{A}$ has $I$-coproducts and that $F,G$ preserve them.
		Then the full subcategory
		\begin{displaymath}
			\mathcal{A}_\eta :=
			\{
			A\in \mathcal{A} : \eta_A\textrm{ is an isomorphism}
			\}
		\end{displaymath}
		is strictly full and is closed under $I$-coproducts and $I$-cofactors.\footnote{This last thing means if $\coprod_{i\in I}A_i$ lies in $\mathcal{A}_\eta$, then so does each $A_i$.}
		\item \label[lemma]{lem:subcat_where_trinat_is_iso_is_localizing}
		Let $\eta\colon F\Rightarrow G$ be a trinatural transformation of triangulated functors $\mathcal{A}\overset{F}{\underset{G}{\rightrightarrows}}\mathcal{B}$.
		The full subcategory
		\begin{displaymath}
			\mathcal{A}_\eta :=
			\{
			A\in \mathcal{A} : \eta_A\textrm{ is an isomorphism}
			\}
		\end{displaymath}
		is a thick triangulated subcategory of $\mathcal{A}$.
		If $\mathcal{A}$ has small coproducts and $F,G$ preserve small coproducts, then $\mathcal{A}_\eta$ is localizing.
		\item \label[lemma]{lem:prewhiskering_by_ess_dense_reflects_isos}
		Let
		\begin{displaymath}
			% https://q.uiver.app/#q=WzAsMyxbMCwwLCJcXG1hdGhjYWx7QX0iXSxbMSwwLCJcXG1hdGhjYWx7Qn0iXSxbMiwwLCJcXG1hdGhjYWx7Q30iXSxbMSwyLCJHXzEiLDAseyJjdXJ2ZSI6LTF9XSxbMSwyLCJHXzIiLDIseyJjdXJ2ZSI6MX1dLFswLDEsIkYiXSxbMyw0LCJcXGV0YSIsMCx7InNob3J0ZW4iOnsic291cmNlIjoyMCwidGFyZ2V0IjoyMH19XV0=
			\begin{tikzcd}[cramped]
				{\mathcal{A}} & {\mathcal{B}} & {\mathcal{C}}
				\arrow["F", from=1-1, to=1-2]
				\arrow[""{name=0, anchor=center, inner sep=0}, "{G_1}", bend left=30, from=1-2, to=1-3]
				\arrow[""{name=1, anchor=center, inner sep=0}, "{G_2}"', bend right=30, from=1-2, to=1-3]
				\arrow["\eta", shorten <=2pt, shorten >=2pt, Rightarrow, from=0, to=1]
			\end{tikzcd}
		\end{displaymath}
		be a diagram in the $2$-category of triangulated categories.
		That, is $F,G_1,G_2$ is are triangulated functors of triangulated categories and $\eta$ is a trinatural transformation.
		Suppose $F$ is essentially dense.
		Then $\eta$ is an isomorphism if, and only if, $\eta F$ is an isomorphism.
	\end{enumerate}
\end{lemma}

\begin{proof}
	We prove the first point.
	Strict fullness comes from naturality of $\eta$.
	Let $A_i\in\mathcal{A}$, $i\in I$, be objects.
	Consider the following diagram:
	\begin{displaymath}
		% https://q.uiver.app/#q=WzAsNCxbMCwwLCJcXGNvcHJvZF97aVxcaW4gSX1GQV9pIl0sWzAsMSwiXFxjb3Byb2Rfe2lcXGluIEl9R0FfaSJdLFsxLDAsIkZcXGxlZnQoXFxjb3Byb2Rfe2lcXGluIEl9QV9pXFxyaWdodCkiXSxbMSwxLCJHXFxsZWZ0KFxcY29wcm9kX3tpXFxpbiBJfUFfaVxccmlnaHQpIl0sWzAsMSwiXFxjb3Byb2Rfe2lcXGluIEl9XFxldGFfe0FfaX0iLDJdLFsyLDMsIlxcZXRhX3tcXGNvcHJvZF97aVxcaW4gSX1BX2l9Il0sWzAsMl0sWzEsM11d
		\begin{tikzcd}[cramped]
			{\coprod_{i\in I}FA_i} & {F\left(\coprod_{i\in I}A_i\right)} \\
			{\coprod_{i\in I}GA_i} & {G\left(\coprod_{i\in I}A_i\right)}
			\arrow[from=1-1, to=1-2]
			\arrow["{\coprod_{i\in I}\eta_{A_i}}"', from=1-1, to=2-1]
			\arrow["{\eta_{\coprod_{i\in I}A_i}}", from=1-2, to=2-2]
			\arrow[from=2-1, to=2-2]
		\end{tikzcd}
	\end{displaymath}
	where the composite $HA_j\hookrightarrow\coprod_{i\in I}HA_i\to H\left(\coprod_{i\in I}A_i\right)$ equals $H(A_j\hookrightarrow\coprod_{i\in I}A_i)$, for $H\in\{F,G\}$.
	The square commutes and the hypothesis means the top and bottom map are isomorphisms.
	It follows that $\mathcal{A}_\eta$ is closed under $I$-coproducts.
	On the other hand, since $\operatorname{Hom}_\mathcal{A}(\coprod_{i\in I}\eta_{A_i},-)\cong\prod_{i\in I}\operatorname{Hom}_\mathcal{A}(\eta_{A_i},-)$,
	Yoneda lemma implies $\mathcal{A}_\eta$ is closed under $I$-cofactors.

	We show the second point.
	Closedness under shifts follow from the transnaturality condition \eqref{diag:transnat},
	see \Cref{lem:transnat_is_trinat}.
	Closedness under cones follow from the trinaturality condition \eqref{diag:trinat}.
	Point 1 gives strict fullness.
	Hence $\mathcal{A}_\eta\subseteq\mathcal{A}$ is a triangulated subcategory.
	Also, Point 1 gives the claim on thickness (resp.\ localizingness) when applied to $I$ a finite set (resp.\ $I$ an arbitrary set).

	Lastly, the third point is an immediate application of the second one.
\end{proof}

\begin{lemma}
	\label{lem:localizing_iff_cpt_gen}
	Assume $\mathcal{T}$ is compactly generated.
	Then $G\in \mathcal{T}^c$ compactly generates $\mathcal{T}$ if, and only if, $\overline{\langle G \rangle} = \mathcal{T}$.
\end{lemma}

\begin{proof}
	By \cite[Remark 1.16, Theorem 1.17, \& Proposition 9.1.19]{Neeman:2001}, the Verdier localization $\pi\colon \mathcal{T}\to \mathcal{T}/\overline{\langle G \rangle}$ admits a right adjoint $Q\colon \mathcal{T}/\overline{\langle G \rangle}\to \mathcal{T}$.
	Hence, for every $E\in \mathcal{T}$, there exists a distinguished triangle
	\begin{displaymath}
		S_E \to E \to (Q\circ \pi)(E) \to S_E [1]
	\end{displaymath}
	with $S_E \in \overline{\langle G \rangle}$ and $(Q\circ \pi)(E)\in \overline{\langle G \rangle}^\perp$ (e.g.\ $A\in \mathcal{T}$ such that $\operatorname{Hom}(G,A)\cong 0$ for all $G\in \overline{\langle G \rangle}$).
	Therefore, the claim follows.
\end{proof}

\begin{lemma}
	\label{lem:equivalence_or_fully_faithful_via_compacts}
	Suppose $F\colon \mathcal{T} \to \mathcal{S}$ is an exact functor between triangulated categories which are compactly generated by a single object.
	Assume $F$ has a right adjoint $G$ that preserves small coproducts.
	Then $F$ is an equivalence if, and only if, $F$ restricts to an equivalence $\mathcal{T}^c \to \mathcal{S}^c$.
\end{lemma}

\begin{proof}
	We may assume $F\dashv G$ in the $2$-category of triangulated categories, \Cref{lem:triadjoint_iff_adjoint}.
	By \cite[Theorem 5.1]{Neeman:1996}, $F$ preserves compacts.
	Recall that an equivalence $\mathcal{T}\to\mathcal{S}$ induces an equivalence on compact objects $\mathcal{T}^c \to \mathcal{S}^c$.
	Conversely, if $F$ restricts to an equivalence $F\colon \mathcal{T}^c \to \mathcal{S}^c$, then $F\colon \mathcal{T}^c \to \mathcal{S}^c$ has a right adjoint.
	Hence, \cite[Lemma 2.6]{Balmer/Dell'Ambrogio/Sanders:2016} implies $G$ preserves compact objects, and the adjunction of $F$ and $G$ on ambient categories restricts to an adjunction on categories of compact objects.
	To show $F$ is an equivalence, we show that $F$ and $G$ are both fully faithful.
	Denote the counit and unit of the adjunction respectively by $\epsilon$ and $\eta$.
	These are trinatural transformations, see \Cref{lem:equiv_transadj,lem:transnat_is_trinat}.
	Set $\mathcal{T}^\prime\subseteq\mathcal{T}$ and $\mathcal{S}^\prime\subseteq\mathcal{S}$ to be the strictly full subcategories  consisting respectively of objects where $\eta$ and $\epsilon$ are isomorphisms contain $\mathcal{T}^c$ and $\mathcal{S}^c$.
	By \Cref{lem:subcat_where_trinat_is_iso_is_localizing}, $\mathcal{T}^\prime$ and $\mathcal{S}^\prime$ are closed under shifts, extensions, and small coproducts (hence, homotopy colimits).
	By \cite[\href{https://stacks.math.columbia.edu/tag/09SN}{Tag 09SN}]{StacksProject}, each object in $\mathcal{T}$ and $\mathcal{S}$ can be obtained respectively by $\mathcal{T}^c$ and $\mathcal{S}^c$ using these operations.
	Hence, $\mathcal{T} \subseteq \mathcal{T}^\prime$ and $\mathcal{S} \subseteq \mathcal{S}^\prime$, which completes the proof.
\end{proof}

\begin{remark}
	\label{rmk:equivalences_induced}
	For $i=1,2$, let $Y_i$ be Noetherian algebraic spaces.
	By \cite[Corollary 1.2]{Lank:2026d}, each $D_{\operatorname{qc}}(Y_i)$ is a weakly approximable triangulated category which is coherent.
	See \cite[Definition 1.25]{Neeman:2026} and \cite[Definition 10.1.1]{Canonaco/Haesemeyer/Neeman/Stellari:2024} for details.
	Let $\Phi\colon D_{\operatorname{qc}}(Y_1)\to D_{\operatorname{qc}}(Y_2)$ be an exact functor that preserves coproducts.
	The following are equivalent:
	\begin{enumerate}
		\item \label{rmk:equivalences_induced1} $\Phi$ is an equivalence
		\item \label{rmk:equivalences_induced2} $\Phi$ restricts to an equivalence on $D^-_{\operatorname{coh}}$
		\item \label{rmk:equivalences_induced3} $\Phi$ restricts to an equivalence on $D^b_{\operatorname{coh}}$
		\item \label{rmk:equivalences_induced4} $\Phi$ restricts to an equivalence on $\operatorname{Perf}$.
	\end{enumerate}
	Indeed, by \cite[Theorem B]{Canonaco/Haesemeyer/Neeman/Stellari:2024}, $\eqref{rmk:equivalences_induced1} \implies \eqref{rmk:equivalences_induced2} \implies \eqref{rmk:equivalences_induced3} \implies \eqref{rmk:equivalences_induced4}$.
	Additionally, if $\eqref{rmk:equivalences_induced4}$ holds, then $\Phi$ admits a right adjoint that preserves small coproducts \cite[Theorem 5.1]{Neeman:1996}.
	Hence, \Cref{lem:equivalence_or_fully_faithful_via_compacts} tells us $\eqref{rmk:equivalences_induced4} \implies \eqref{rmk:equivalences_induced1}$.
	In the particular case $\Phi=\Phi_K$ is a relative Fourier--Mukai transform with a relative perfect kernel $K$ over each $Y_i$, the fact that \eqref{rmk:equivalences_induced1}, \eqref{rmk:equivalences_induced3}, \eqref{rmk:equivalences_induced4} are equivalent can be proven without resorting to weak approximability, see \Cref{lem:equivalences_induced}.
\end{remark}

\begin{remark}
	\label{rem:rel_perfect_is_free}
	Let $f_i\colon Y_i\to S$ be proper morphisms of Noetherian algebraic spaces.
	Let $K\in D^-_{\operatorname{coh}}(Y_1\times_SY_2)$.
	If $\Phi_K$ induces an equivalence on any of $D_{\operatorname{qc}}$, $D_{\operatorname{coh}}^-$, $D_{\operatorname{coh}}^b$ or $\operatorname{Perf}$, then \Cref{cor:preservation} and \Cref{rmk:equivalences_induced} imply that $K$ is relatively perfect over each $Y_i$.
\end{remark}

%%%%%%%%%%%%%%%%%%%%%%%%%%%%%%%%%%%%%
\subsection{Algebraic spaces}
\label{sec:preliminaries_stacks}
%%%%%%%%%%%%%%%%%%%%%%%%%%%%%%%%%%%%%

We follow \cite{StacksProject} for conventions on algebraic spaces.
Unless stated otherwise, we work over a fixed base scheme.
In this subsection, let $X$ be a quasi-compact quasi-separated algebraic space.

%%%%%%%%%%%%%%%%%%%%%%%%%%%%%%%%%%%%%%%%%%%%
\subsubsection{Notions}
\label{sec:prelim_stacks_notions}
%%%%%%%%%%%%%%%%%%%%%%%%%%%%%%%%%%%%%%%%%%%%

An \textbf{\'{e}tale presentation} of $X$ is an \'{e}tale, finitely presented, surjective morphism to $X$ from a scheme.
The underlying topological space of $X$ is given by equivalence classes of morphisms from fields to the algebraic space \cite[\href{https://stacks.math.columbia.edu/tag/03BT}{Tag 03BT}]{StacksProject}.
We denote it by $|X|$.
Recall that a morphism from a field to a quasi-separated algebraic space is an affine morphism \cite[\href{https://stacks.math.columbia.edu/tag/09TF}{Tag 09TF}]{StacksProject}.

%%%%%%%%%%%%%%%%%%%%%%%%%%%%%%%%%%%%%%%%%%%%
\subsubsection{Categories}
\label{sec:prelim_stacks_categories}
%%%%%%%%%%%%%%%%%%%%%%%%%%%%%%%%%%%%%%%%%%%%

$\operatorname{Mod}(X)$ is the Grothendieck abelian category of sheaves of $\mathcal{O}_X$-modules on the small \'{e}tale site $X_{\textrm{\'{e}tale}}$ of $X$ \cite[\href{https://stacks.math.columbia.edu/tag/03LT}{Tag 03LT}]{StacksProject}.
$\operatorname{Qcoh}(X)$ is the full subcategory of $\operatorname{Mod}(X)$ consisting of quasi-coherent sheaves.
$D(X) := D(\operatorname{Mod}(X))$ is the derived category of $\operatorname{Mod}(X)$.
$D_{\operatorname{qc}}(X)$ is the full subcategory of $D(X)$ consisting of complexes with quasi-coherent cohomology sheaves.
$\operatorname{Perf}(X)$ is the full subcategory of perfect complexes in $D_{\operatorname{qc}}(X)$.
If $X$ is Noetherian, then $\operatorname{coh}(X)$ is the full subcategory of $\operatorname{Mod}(X)$ consisting of coherent sheaves and $D^b_{\operatorname{coh}}(X)$ denotes the full subcategory of $D(X)$ consisting of bounded pseudocoherent complexes.

%%%%%%%%%%%%%%%%%%%%%%%%%%%%%%%%%%%
\subsubsection{Support}
\label{sec:preliminaries_stacks_support}
%%%%%%%%%%%%%%%%%%%%%%%%%%%%%%%%%%%

For any $M \in \operatorname{Qcoh}(X)$, define $\operatorname{Supp}(M) := p(\operatorname{Supp}(p^\ast M)) \subseteq |X|$ where $p\colon U \to X$ is any \'{e}tale presentation from a scheme; this is independent of $p$ \cite[\href{https://stacks.math.columbia.edu/tag/07TY}{Tag 07TY}]{StacksProject}.
This coincides with the notion of support for abelian sheaves on the small \'{e}tale site \cite[\href{https://stacks.math.columbia.edu/tag/04K7}{Tag 04K7}]{StacksProject}.
For $E \in D_{\operatorname{qc}}(X)$, set
\begin{displaymath}
	\operatorname{Supp}(E) := \bigcup_{j \in \mathbb{Z}} \operatorname{Supp}(\mathcal{H}^j(E)) \subseteq |X|.
\end{displaymath}
Given a closed subset $Z \subseteq |X|$, we say $E$ is \textbf{supported on $Z$} if $\operatorname{Supp}(E) \subseteq Z$. If $E\in \operatorname{Qcoh}(X)$ and is of finite type, then $\operatorname{Supp}(E)\subseteq |X|$ is closed \cite[\href{https://stacks.math.columbia.edu/tag/07TZ}{Tag 07TZ}]{StacksProject}.
On a locally Noetherian algebraic space, quasi-coherent sheaves of finite type coincide with coherent sheaves \cite[\href{https://stacks.math.columbia.edu/tag/07UB}{Tag 07UB}]{StacksProject}.
Denote by $D_{\operatorname{qc},Z}(X)$ the strictly full subcategory of $D_{\operatorname{qc}}(X)$ consisting of objects supported on $Z$.

\begin{lemma}
	\label{lem:factor_residual_gerbes}
	Let $X$ be a decent algebraic space (e.g.\ quasi-separated \cite[\href{https://stacks.math.columbia.edu/tag/03I7}{Tag 03I7}]{StacksProject}).
	If $t\colon \operatorname{Spec}(k)\to X$ is a morphism from a field such that $t(\operatorname{Spec}(k))=p$, then $t$ represents $p$ (see \cite[\href{https://stacks.math.columbia.edu/tag/03BT}{Tag 03BT}]{StacksProject}).
\end{lemma}

\begin{proof}
	Let $i\colon Z_p\to X$ be the residual space of $X$ at $p$ (see \cite[\href{https://stacks.math.columbia.edu/tag/06QZ}{Tag 06QZ} and \href{https://stacks.math.columbia.edu/tag/06R0}{06R0}]{StacksProject}).
	Choose any $h\colon \operatorname{Spec}(\ell)\to X$ that represents $p$.
	By \cite[\href{https://stacks.math.columbia.edu/tag/06QZ}{Tag 06QZ}]{StacksProject} (e.g.\ use $X$ is decent), there exists a commutative diagram
	\begin{displaymath}
		% https://q.uiver.app/#q=WzAsMyxbMSwwLCJcXG1hdGhjYWx7Wn1fcCJdLFsxLDEsIlxcbWF0aGNhbHtYfS4iXSxbMCwwLCJcXG9wZXJhdG9ybmFtZXtTcGVjfShrKSJdLFswLDEsImkiXSxbMiwxLCJ0IiwyXSxbMiwwLCJ0XlxccHJpbWUiXV0=
		\begin{tikzcd}
			{\operatorname{Spec}(k)} & {Z_p} \\
			& {X.}
			\arrow["{t^\prime}", from=1-1, to=1-2]
			\arrow["t"', from=1-1, to=2-2]
			\arrow["i", from=1-2, to=2-2]
		\end{tikzcd}
	\end{displaymath}
	Moreover, from \cite[\href{https://stacks.math.columbia.edu/tag/0H1T}{Tag 0H1T}]{StacksProject}, there exists a commutative diagram
	\begin{displaymath}
		% https://q.uiver.app/#q=WzAsMyxbMCwwLCJcXG1hdGhjYWx7Wn1fcCJdLFswLDEsIlxcbWF0aGNhbHtYfS4iXSxbMSwwLCJcXG9wZXJhdG9ybmFtZXtTcGVjfShcXGVsbCkiXSxbMCwxLCJpIiwyXSxbMiwxLCJoIl0sWzIsMCwiaF5cXHByaW1lIiwyXV0=
		\begin{tikzcd}
			{Z_p} & {\operatorname{Spec}(\ell)} \\
			{X.}
			\arrow["i"', from=1-1, to=2-1]
			\arrow["{h^\prime}"', from=1-2, to=1-1]
			\arrow["h", from=1-2, to=2-1]
		\end{tikzcd}
	\end{displaymath}
	Consider the fiber product
	\begin{displaymath}
		% https://q.uiver.app/#q=WzAsNCxbMSwxLCJcXG1hdGhjYWx7Wn1fcC4iXSxbMCwxLCJcXG9wZXJhdG9ybmFtZXtTcGVjfShrKSJdLFsxLDAsIlxcb3BlcmF0b3JuYW1le1NwZWN9KFxcZWxsKSJdLFswLDAsIlxcb3BlcmF0b3JuYW1le1NwZWN9KGspIFxcdGltZXNfe1xcbWF0aGNhbHtafV9wfSBcXG9wZXJhdG9ybmFtZXtTcGVjfShcXGVsbCkgIl0sWzEsMCwidF5cXHByaW1lIiwyXSxbMiwwLCJoXlxccHJpbWUiXSxbMywxLCJwXzEiLDJdLFszLDIsInBfMiJdXQ==
		\begin{tikzcd}
			{\operatorname{Spec}(k) \times_{Z_p} \operatorname{Spec}(\ell) } & {\operatorname{Spec}(\ell)} \\
			{\operatorname{Spec}(k)} & {Z_p.}
			\arrow["{p_2}", from=1-1, to=1-2]
			\arrow["{p_1}"', from=1-1, to=2-1]
			\arrow["{h^\prime}", from=1-2, to=2-2]
			\arrow["{t^\prime}"', from=2-1, to=2-2]
		\end{tikzcd}
	\end{displaymath}
	Since $|Z_p|$ is a singleton, $t^\prime$ is surjective.
	Hence, by base change, $p_2$ is surjective.
	Thus, $\operatorname{Spec}(k) \times_{Z_p} \operatorname{Spec}(\ell)$ is nonempty.
	Here $\operatorname{Spec}(k) \times_{Z_p} \operatorname{Spec}(\ell)$ is an algebraic space.
	Choose an \'{e}tale presentation $s\colon U\to \operatorname{Spec}(k) \times_{Z_p} \operatorname{Spec}(\ell)$.
	Fix some $p^\prime \in |\operatorname{Spec}(k) \times_{Z_p} \operatorname{Spec}(\ell)|$.
	Consider any $g\colon \operatorname{Spec}(K)\to \operatorname{Spec}(k) \times_{Z_p} \operatorname{Spec}(\ell)$ that represents $p^\prime$.
	Then we have a commutative diagram
	\begin{displaymath}
		% https://q.uiver.app/#q=WzAsNCxbMSwxLCJcXG1hdGhjYWx7WH0uIl0sWzAsMSwiXFxvcGVyYXRvcm5hbWV7U3BlY30oaykiXSxbMSwwLCJcXG9wZXJhdG9ybmFtZXtTcGVjfShcXGVsbCkiXSxbMCwwLCJcXG9wZXJhdG9ybmFtZXtTcGVjfShLKSJdLFsxLDAsInQiLDJdLFsyLDAsImgiXSxbMywxLCJwXzEgXFxjaXJjIGciLDJdLFszLDIsInBfMiBcXGNpcmMgZyJdXQ==
		\begin{tikzcd}
			{\operatorname{Spec}(K)} & {\operatorname{Spec}(\ell)} \\
			{\operatorname{Spec}(k)} & {X.}
			\arrow["{p_2 \circ g}", from=1-1, to=1-2]
			\arrow["{p_1 \circ g}"', from=1-1, to=2-1]
			\arrow["h", from=1-2, to=2-2]
			\arrow["t"', from=2-1, to=2-2]
		\end{tikzcd}
	\end{displaymath}
	This completes the proof.
\end{proof}

\begin{lemma}
	\label{lem:support_is_cohomological_for_finite_type_cohomology}
	Let $X$ be a Noetherian algebraic space and $E\in D_{\operatorname{qc}}(X)$ have coherent cohomology.
	Then $\operatorname{Supp}(E)$ coincides with the $p\in |X|$ such that $\mathbf{L}i^\ast E \not\cong 0$ where $i$ is a representative of $p$; in fact, this is independent of the representative.
\end{lemma}

\begin{proof}
	Let $s\colon U \to X$ be an \'{e}tale presentation.
	The first claim follows from \cite[\href{https://stacks.math.columbia.edu/tag/07TZ}{Tag 07TZ}]{StacksProject} and the string of equalities:
	\begin{displaymath}
		\begin{aligned}
			\operatorname{Supp}(E)
			&= \bigcup_n \operatorname{Supp}(\mathcal{H}^n (E))
			\\&= \bigcup_n s(\operatorname{Supp}(s^\ast \mathcal{H}^n (E)))
			\\&= \bigcup_n s(\operatorname{Supp}(\mathcal{H}^n (s^\ast E)))
			\\&= \bigcup_n s(\operatorname{Supp}(\mathcal{H}^n (\mathbf{L}s^\ast E)))
			\\&= s(\bigcup_n \operatorname{Supp}(\mathcal{H}^n (\mathbf{L}s^\ast E)))
			\\&= s(\operatorname{Supp}(\mathbf{L}s^\ast E)).
		\end{aligned}
	\end{displaymath}
	By \cite[\href{https://stacks.math.columbia.edu/tag/056J}{Tag 056J}]{StacksProject}, $\operatorname{Supp}(s^\ast \mathcal{H}^n (E))$ is the set of $ q\in U$ such that $i^\ast s^\ast \mathcal{H}^n (E)$ is nonzero where $i\colon \operatorname{Spec}(\kappa(q))\to U$ is the canonical morphism.
	Choose $p\in X$ such that $\mathbf{L}i^\ast E \not\cong 0$ where $i$ is a representative of $p$.
	Since $s$ is surjective, \Cref{lem:factor_residual_gerbes} allows us to choose a representative $i^\prime$ of some $q\in s^{-1}(p)$.
	We claim that $\mathbf{L}(s\circ i^\prime)^\ast E\not\cong0$.
	Indeed, there exists a field $\ell$ and commutative diagram
	\begin{displaymath}
		% https://q.uiver.app/#q=WzAsNCxbMSwwLCJcXG9wZXJhdG9ybmFtZXtTcGVjfShrKSJdLFsxLDEsIlxcbWF0aGNhbHtTfS4iXSxbMCwxLCJcXG9wZXJhdG9ybmFtZXtTcGVjfShcXGthcHBhKHEpKSJdLFswLDAsIlxcb3BlcmF0b3JuYW1le1NwZWN9KFxcZWxsKSJdLFswLDEsImkiXSxbMiwxLCJzXFxjaXJjIGleXFxwcmltZSIsMl0sWzMsMCwiZyJdLFszLDIsImgiLDJdXQ==
		\begin{tikzcd}
			{\operatorname{Spec}(\ell)} & {\operatorname{Spec}(k)} \\
			{\operatorname{Spec}(\kappa(q))} & {X.}
			\arrow["g", from=1-1, to=1-2]
			\arrow["h"', from=1-1, to=2-1]
			\arrow["i", from=1-2, to=2-2]
			\arrow["{s\circ i^\prime}"', from=2-1, to=2-2]
		\end{tikzcd}
	\end{displaymath}
	If $\mathbf{L}i^\ast E \not\cong 0$, then $\mathbf{L}(i\circ g)^\ast E \not\cong 0$.
	Hence, by \Cref{lem:zero_object_via_covers}, $\mathbf{L}(s\circ i^\prime)^\ast E$ is nonzero.
	Notice that this shows independence of the choice of representative for $p$ regarding whether the derived pullback of $E$ is nonzero.
	Now, from \cite[Lemma A.1]{Iyengar/Lipman/Neeman:2015}, $q\in \operatorname{Supp}(\mathbf{L}s^\ast E)$, and so, $p\in \operatorname{Supp}(E)$.
	We check the reverse inclusion.
	Choose $p\in \operatorname{Supp}(E)$.
	As $s$ is surjective, we can find a representative $\operatorname{Spec}(\kappa(q))\to X$ given by the composition of $s$ with the canonical morphism $i\colon \operatorname{Spec}(\kappa(q))\to U$ for some $q\in \operatorname{Supp}(\mathbf{L}s^\ast E)$.
	%%NOTE: Use union above to realize supp = s (supp).
	By \cite[Lemma A.1]{Iyengar/Lipman/Neeman:2015}, we have that $\mathbf{L}(s\circ i)^\ast E\not\cong 0$, and so we are done.
	Observe that \cite{Iyengar/Lipman/Neeman:2015} proves results on the small Zariski site of $U$.
	However, it is valid for the small \'{e}tale site of $U$.
	Indeed, we can apply \cite[\href{https://stacks.math.columbia.edu/tag/071Q}{Tag 071Q}]{StacksProject} to realize that the morphism of ringed topoi $\epsilon\colon (U_{\textrm{\'{e}tale}} , \mathcal{O}_{U_{\textrm{\'{e}tale}}}) \to (U_{\textrm{Zar}}, \mathcal{O}_{\textrm{Zar}})$ yields a $t$-exact equivalence
	\begin{displaymath}
		\mathbf{L}\epsilon^\ast\colon D_{\operatorname{qc}} (\mathcal{O}_{\textrm{Zar}}) \leftrightarrows D_{\operatorname{qc}} (\mathcal{O}_{\textrm{\'{e}tale}})\colon \mathbf{R}\epsilon_\ast.
	\end{displaymath}
	By \cite[\href{https://stacks.math.columbia.edu/tag/07TY}{Tag 07TY}]{StacksProject}, the support of an abelian sheaf on $U_{\textrm{\'{e}tale}}$ agrees with the notion as defined for quasi-coherent modules on schemes.
	This completes the proof.
\end{proof}

%%%%%%%%%%%%%%%%%%%%%%%%%%%%%%%%%%%%%%%%%%%%
\subsubsection{Perfect complexes}
\label{sec:prelim_stacks_perfects}
%%%%%%%%%%%%%%%%%%%%%%%%%%%%%%%%%%%%%%%%%%%%

Perfect complexes are defined on any ringed site \cite[\href{https://stacks.math.columbia.edu/tag/08G4}{Tag 08G4}]{StacksProject}, e.g.\ on the small \'{e}tale site of $X$.
A complex is \textbf{strictly perfect} if it is a bounded complex whose terms are direct summands of finite free modules.
A complex is \textbf{perfect} if it is locally strictly perfect.
For any $Z\subseteq |X|$ closed on a quasi-compact quasi-separated algebraic space such that $|X|\setminus Z$ is quasi-compact, perfect complexes of $D_{\operatorname{qc},Z}(X)$ coincide with compacts of $D_{\operatorname{qc},Z}(X)$, and in particular they are perfect complexes in $D_{\operatorname{qc}}(X)$.
Moreover, $D_{\operatorname{qc},Z}(X)$ is singly compactly generated for such $Z$.
See \cite[\href{https://stacks.math.columbia.edu/tag/0AED}{Tag 0AED}, \href{https://stacks.math.columbia.edu/tag/09M8}{09M8}, and \href{https://stacks.math.columbia.edu/tag/0AEC}{0AEC}]{StacksProject}.

%%%%%%%%%%%%%%%%%%%%%%%%%%%%%%%%%%%%%
\subsubsection{Internal homs}
\label{sec:prelim_stacks_internal_homs}
%%%%%%%%%%%%%%%%%%%%%%%%%%%%%%%%%%%%%

Let $E,G\in D(X)$.
There exists the derived tensor product $E\otimes^{\mathbf{L}} G \in D(X)$ and derived sheaf Hom $\operatorname{\mathbb{R}\mathcal{H}\! \mathit{om}}(E,G)\in D(X)$.
Here, $(-)\otimes^{\mathbf{L}} E$ is left adjoint to $\operatorname{\mathbb{R}\mathcal{H}\! \mathit{om}}(E,-)$ on $D(X)$. Applying \cite[\href{https://stacks.math.columbia.edu/tag/08F5}{Tag 08F5} and \href{https://stacks.math.columbia.edu/tag/0FPQ}{0FPQ}]{StacksProject}, $D_{\operatorname{qc}}$ is symmetric monoidal, i.e.\ $E\otimes^{\mathbf{L}} G \in D_{\operatorname{qc}}(X)$ for all $E,G\in D_{\operatorname{qc}}(X)$.
However, this need not be the case for $\operatorname{\mathbb{R}\mathcal{H}\! \mathit{om}}(E,G)$ despite the formation of $\operatorname{\mathbb{R}\mathcal{H}\! \mathit{om}}(E,-)$ being \'{e}tale local \cite[\href{https://stacks.math.columbia.edu/tag/04LX}{Tag 04LX} and \href{https://stacks.math.columbia.edu/tag/08JB}{08JB}]{StacksProject}. By \cite[\href{https://stacks.math.columbia.edu/tag/09IY}{Tag 09IY}]{StacksProject}, the category $D_{\operatorname{qc}}(X)$ is singly compactly generated.
Also, the endofunctor $(-)\otimes^{\mathbf{L}} E$ on $D_{\operatorname{qc}}(X)$ preserves small coproducts.
Hence, \cite[Theorem 8.4.4]{Neeman:2001} (or \cite[Theorem 4.1]{Neeman:1996}) implies that the endofunctor admits a right adjoint $\operatorname{\mathbf{R}\mathcal{H}\! \mathit{om}}(E,-)$ on $D_{\operatorname{qc}}(X)$.
This equips $D_{\operatorname{qc}}(X)$ with the structure of a closed symmetric monoidal triangulated category.
Denote by $i\colon D_{\operatorname{qc}}(X)\to D(X)$ for the natural inclusion.
It admits a right adjoint $Q\colon D(X) \to D_{\operatorname{qc}}(X)$ by reasoning above because $i$ preserves small coproducts.

\begin{lemma}
	\label{lem:rhom_tensor_morphism}
	Let $X$ be an algebraic space.
	If $K$ is perfect and $L,M\in D(X)$, then there exists a functorial isomorphism
	\begin{displaymath}
		K\otimes^{\mathbf{L}} \operatorname{\mathbb{R}\mathcal{H}\! \mathit{om}}(M,L) \to \operatorname{\mathbb{R}\mathcal{H}\! \mathit{om}} ( \operatorname{\mathbb{R}\mathcal{H}\! \mathit{om}}(K,M), L).
	\end{displaymath}
\end{lemma}

\begin{proof}
	We compute this directly:
	\begin{displaymath}
		\begin{aligned}
			K\otimes^{\mathbf{L}} \operatorname{\mathbb{R}\mathcal{H}\! \mathit{om}}(M,L)
			& \cong \operatorname{\mathbb{R}\mathcal{H}\! \mathit{om}} ( \operatorname{\mathbb{R}\mathcal{H}\! \mathit{om}}(K,\mathcal{O}_X), \operatorname{\mathbb{R}\mathcal{H}\! \mathit{om}}(M,L)) && \textrm{(\cite[\href{https://stacks.math.columbia.edu/tag/08JJ}{Tag 08JJ}]{StacksProject})}
			\\& \cong \operatorname{\mathbb{R}\mathcal{H}\! \mathit{om}} ( \operatorname{\mathbb{R}\mathcal{H}\! \mathit{om}}(K,\mathcal{O}_X) \otimes^{\mathbf{L}} M , L) && \textrm{(\cite[\href{https://stacks.math.columbia.edu/tag/08J9}{Tag 08J9}]{StacksProject})}
			\\& \cong \operatorname{\mathbb{R}\mathcal{H}\! \mathit{om}} ( \operatorname{\mathbb{R}\mathcal{H}\! \mathit{om}}(K,M), L) && \textrm{(\cite[\href{https://stacks.math.columbia.edu/tag/08JJ}{Tag 08JJ}]{StacksProject})}.
		\end{aligned}
	\end{displaymath}
\end{proof}

\begin{lemma}
	\label{lem:internal_hom}
	Let $X$ be an algebraic space.
	For any $E,L\in D_{\operatorname{qc}}(X)$, there is a canonical isomorphism
	\begin{displaymath}
		\operatorname{\mathbf{R}\mathcal{H}\! \mathit{om}} (E,L) \to Q(\operatorname{\mathbb{R}\mathcal{H}\! \mathit{om}} (i(E),i(L))).
	\end{displaymath}
	In particular, if $\operatorname{\mathbb{R}\mathcal{H}\! \mathit{om}} (i(E),i(L))$ has quasi-coherent cohomology, then
	\begin{displaymath}
		\operatorname{\mathbf{R}\mathcal{H}\! \mathit{om}} (i(E),i(L))
		\cong \operatorname{\mathbb{R}\mathcal{H}\! \mathit{om}} (i(E),i(L)).
	\end{displaymath}
	This happens if $E$ is pseudocoherent and $L$ is locally bounded below, see \cite[\href{https://stacks.math.columbia.edu/tag/0A8A}{Tag 0A8A}]{StacksProject}.
\end{lemma}

\begin{proof}
	Note that there exist adjunctions
	\begin{displaymath}
		% https://q.uiver.app/#q=WzAsMyxbMCwwLCJEX3tcXG9wZXJhdG9ybmFtZXtxY319KFxcbWF0aGNhbHtYfSkiXSxbMSwwLCJEKFxcbWF0aGNhbHtYfSkiXSxbMiwwLCJEKFxcbWF0aGNhbHtYfSkuIl0sWzEsMiwiKCAoLSkgXFxvdGltZXNee1xcbWF0aGJme0x9fSBpKEUpKSIsMCx7ImN1cnZlIjotMn1dLFsyLDEsIlxcb3BlcmF0b3JuYW1le1xcbWF0aGJie1J9XFxtYXRoY2Fse0h9XFwhIFxcbWF0aGl0e29tfX0gKGkoRSksLSkpIiwwLHsiY3VydmUiOi0yfV0sWzAsMSwiaSIsMCx7ImN1cnZlIjotMn1dLFsxLDAsIlEiLDAseyJjdXJ2ZSI6LTJ9XV0=
		\begin{tikzcd}
			{D_{\operatorname{qc}}(X)} & {D(X)} & {D(X).}
			\arrow["i", bend left=30, from=1-1, to=1-2]
			\arrow["Q", bend left=30, from=1-2, to=1-1]
			\arrow["{( (-) \otimes^{\mathbf{L}} i(E))}", bend left=30, from=1-2, to=1-3]
			\arrow["{\operatorname{\mathbb{R}\mathcal{H}\! \mathit{om}} (i(E),-)}", bend left=30, from=1-3, to=1-2]
		\end{tikzcd}
	\end{displaymath}
	By composing, we obtain
	\begin{displaymath}
		% https://q.uiver.app/#q=WzAsMixbMCwwLCJEX3tcXG9wZXJhdG9ybmFtZXtxY319KFxcbWF0aGNhbHtYfSkiXSxbMiwwLCJEKFxcbWF0aGNhbHtYfSkuIl0sWzAsMSwiKCAoLSkgXFxvdGltZXNee1xcbWF0aGJme0x9fSBpKEUpKSBcXGNpcmMgaSIsMCx7ImN1cnZlIjotMn1dLFsxLDAsIlEgKFxcb3BlcmF0b3JuYW1le1xcbWF0aGJie1J9XFxtYXRoY2Fse0h9XFwhIFxcbWF0aGl0e29tfX0gKGkoRSksLSkpIiwwLHsiY3VydmUiOi0yfV1d
		\begin{tikzcd}
			{D_{\operatorname{qc}}(X)} && {D(X).}
			\arrow["{( (-) \otimes^{\mathbf{L}} i(E)) \circ i}", bend left=30, from=1-1, to=1-3]
			\arrow["{Q (\operatorname{\mathbb{R}\mathcal{H}\! \mathit{om}} (i(E),-))}", bend left=30, from=1-3, to=1-1]
		\end{tikzcd}
	\end{displaymath}
	Using that $i$ is monoidal, it follows that $( (-) \otimes^{\mathbf{L}} i(E)) \circ i$ induces an endofunctor on $D_{\operatorname{qc}}(X)$, whose right adjoint is $\operatorname{\mathbf{R}\mathcal{H}\! \mathit{om}} (E,-)$.
	Therefore, restriction gives an adjunction
	\begin{displaymath}
		% https://q.uiver.app/#q=WzAsMixbMCwwLCJEX3tcXG9wZXJhdG9ybmFtZXtxY319KFxcbWF0aGNhbHtYfSkiXSxbMiwwLCJEX3tcXG9wZXJhdG9ybmFtZXtxY319KFxcbWF0aGNhbHtYfSkuIl0sWzAsMSwiKCAoLSkgXFxvdGltZXNee1xcbWF0aGJme0x9fSBpKEUpKSBcXGNpcmMgaSIsMCx7ImN1cnZlIjotMn1dLFsxLDAsIlEgKFxcb3BlcmF0b3JuYW1le1xcbWF0aGJie1J9XFxtYXRoY2Fse0h9XFwhIFxcbWF0aGl0e29tfX0gKGkoRSksLSkpIiwwLHsiY3VydmUiOi0yfV1d
		\begin{tikzcd}
			{D_{\operatorname{qc}}(X)} && {D_{\operatorname{qc}}(X).}
			\arrow["{( (-) \otimes^{\mathbf{L}} i(E)) \circ i}", bend left=30, from=1-1, to=1-3]
			\arrow["{Q (\operatorname{\mathbb{R}\mathcal{H}\! \mathit{om}} (i(E),-))}", bend left=30, from=1-3, to=1-1]
		\end{tikzcd}
	\end{displaymath}
	Thus, the desired claim follows from uniqueness of adjoints.
\end{proof}

\begin{remark}
	It is no surprise that the two internal homs, $\operatorname{\mathbb{R}\mathcal{H}\! \mathit{om}}$ and $\operatorname{\mathbf{R}\mathcal{H}\! \mathit{om}}$, will have very different behaviour.
	Here's a practical example:
	Given a morphism of quasi-compact quasi-separated algebraic spaces, for $L\in D_{\operatorname{qc}}(X)$ and $K\in D_{\operatorname{qc}}(Y)$, we have canonical maps
	\begin{displaymath}
		\begin{aligned}
			\mathbf{R}f_\ast \operatorname{\mathbb{R}\mathcal{H}\! \mathit{om}}(L,f^\times K)
			&\to\operatorname{\mathbb{R}\mathcal{H}\! \mathit{om}}(\mathbf{R}f_\ast L,K),
			\\ \mathbf{R}f_\ast \operatorname{\mathbf{R}\mathcal{H}\! \mathit{om}}(L,f^\times K)
			&\to \operatorname{\mathbf{R}\mathcal{H}\! \mathit{om}}(\mathbf{R}f_\ast L,K).
		\end{aligned}
	\end{displaymath}
	The latter is obtained by mapping the former under the coherator, see \cite[\href{https://stacks.math.columbia.edu/tag/0GG4}{Tag 0GG4}]{StacksProject} and use \Cref{lem:internal_hom}.
	The latter is an isomorphism \cite[\href{https://stacks.math.columbia.edu/tag/0E58}{Tag 0E58}]{StacksProject},
	but the former is not an isomorphism in general \cite[\href{https://stacks.math.columbia.edu/tag/0GEU}{Tag 0GEU}]{StacksProject}.
\end{remark}

%%%%%%%%%%%%%%%%%%%%%%%%%%%%%%%%%%%%%
\subsubsection{Functors}
\label{sec:prelim_stacks_functors}
%%%%%%%%%%%%%%%%%%%%%%%%%%%%%%%%%%%%%

Let $f\colon Y \to X$ be a morphism of quasi-compact quasi-separated algebraic spaces. Recall that there exists an adjoint pair
\begin{displaymath}
	\mathbf{L}f^\ast\colon D(X) \leftrightarrows D(Y)\colon \mathbf{R}f_\ast.
\end{displaymath}
See \cite[\href{https://stacks.math.columbia.edu/tag/07A6}{Tag 07A6}]{StacksProject}. By \cite[\href{https://stacks.math.columbia.edu/tag/08FA}{Tag 08FA} and \href{https://stacks.math.columbia.edu/tag/08F4}{08F4}]{StacksProject}, $\mathbf{L}f^\ast$ and $\mathbf{R}f_\ast$ preserve quasi-coherent cohomology.
Also, from \cite[\href{https://stacks.math.columbia.edu/tag/07A6}{Tag 07A6}]{StacksProject}, the restrictions of $\mathbf{L}f^\ast$ and $\mathbf{R}f_\ast$ form an adjoint pair on $D_{\operatorname{qc}}$. Using \cite[\href{https://stacks.math.columbia.edu/tag/07A4}{Tag 07A4}]{StacksProject}, it follows that $\mathbf{L}f^\ast$ is monoidal. Lastly, $\mathbf{R}f_\ast$ admits a right adjoint $f^\times$ \cite[\href{https://stacks.math.columbia.edu/tag/0E55}{Tag 0E55}]{StacksProject}.

%%%%%%%%%%%%%%%%%%%%%%%%%%%%%%%%%%%
\subsubsection{Adjoint pseudofunctors}
\label{sec:preliminaries_stacks_lipman_hashimoto}
%%%%%%%%%%%%%%%%%%%%%%%%%%%%%%%%%%%

We need the following formalism. See \cite[\S 3.6]{Lipman:2009} for details.

\begin{lemma}
	\label{lem:spacey_Lipman_Hashimoto}
	Let $\operatorname{Sp}$ be the category of quasi-compact quasi-separated algebraic spaces over a base scheme (e.g.\ $\mathbb{Z}$).
	For each $X\in \operatorname{Sp}$, set $X^\ast = X_\ast = D_{\operatorname{qc}}(X)$, which is a closed $\Delta$-category with product $\otimes$ and internal Hom $\operatorname{\mathbf{R}\mathcal{H}\! \mathit{om}}(E,-)$.
	For $g\colon W \to Y$ and $f\colon Y \to X$ in $\operatorname{Sp}$, write $f^\ast := \mathbf{L}f^\ast\colon D_{\operatorname{qc}}(X) \to D_{\operatorname{qc}}(Y)$ and $f_\ast := \mathbf{R}f_\ast\colon D_{\operatorname{qc}}(Y) \to D_{\operatorname{qc}}(X)$, and $d_{f,g}$ and $c_{f,g}$ be respectively as in \cite[\href{https://stacks.math.columbia.edu/tag/0D6D}{Tag 0D6D} and \href{https://stacks.math.columbia.edu/tag/0D6E}{0D6E}]{StacksProject}.
	In the sense of \cite[(3.6.10)]{Lipman:2009}, this defines an adjoint pair $({}^\ast,{}_\ast)$ of monoidal $\Delta$-pseudofunctors on $\operatorname{Sp}$.
\end{lemma}

\begin{proof}
	We verify the hypotheses of \cite[(3.6.10)]{Lipman:2009}, following the proof of the \cite[Scholium, (3.6.10)]{Lipman:2009} where appropriate.
	Denote by $\operatorname{TriCat}$ the $2$-category of triangulated categories.
	\begin{enumerate}
		\item The functor $(-)_\ast\colon \operatorname{Sp} \to \operatorname{TriCat}$ is covariant in the sense of \cite[(3.6.5)]{Lipman:2009}.
		It is easy to see that $c_{1,g} = c_{f,1}$ are identities.
		To see that \cite[(3.6.5.1)]{Lipman:2009} commutes, observe that \cite[$(3.6.3)_\ast$]{Lipman:2009} holds for underived pushforward of complexes of modules on the small \'{e}tale sites.
		Indeed, the proof of \cite[$(3.6.3)_\ast$]{Lipman:2009} depends only on functoriality of composition of direct images of morphisms of ringed topoi.
		Now, \cite[\href{https://stacks.math.columbia.edu/tag/0D6E}{Tag 0D6E} and \href{https://stacks.math.columbia.edu/tag/07A5}{07A5}]{StacksProject} can be used to show that the analog of \cite[(3.6.4.1)]{Lipman:2009} holds in our setting.
		Hence, \cite[(3.6.5.1)]{Lipman:2009} commutes.
		\item The functor $(-)^\ast\colon \operatorname{Sp} \to \operatorname{TriCat}$ is contravariant in the sense of \cite[(3.6.5)]{Lipman:2009}.
		It is easy to see that $d_{1,g} = d_{f,1}$ are identities.
		We check that the analog of \cite[(3.6.5.2)]{Lipman:2009} commutes.
		Note that \cite[\href{https://stacks.math.columbia.edu/tag/0D6D}{Tag 0D6D}]{StacksProject} gives the analog of \cite[$(3.6.4)^\ast$]{Lipman:2009}.
		Since derived pullback preserves $K$-flat resolutions \cite[\href{https://stacks.math.columbia.edu/tag/0G7E}{Tag 0G7E}]{StacksProject}, the analog of \cite[$(3.6.3)^\ast$]{Lipman:2009} with derived pullbacks holds, yielding \cite[(3.6.5.2)]{Lipman:2009}.
		\item We have explained \cite[$(3.6.7)(a,b,c)$ \& $(3.6.7)(d(i,ii))$]{Lipman:2009}.
		However, \cite[(3.6.7.2)]{Lipman:2009} is \cite[\href{https://stacks.math.columbia.edu/tag/0FPN}{Tag 0FPN}]{StacksProject}, \cite[(3.6.7.1)]{Lipman:2009} is formal because the proof depends only on properties of closed symmetric monoidal categories, and \cite[Definition 3.4.2]{Lipman:2009} are \cite[\href{https://stacks.math.columbia.edu/tag/0FPK}{Tag 0FPK}, \href{https://stacks.math.columbia.edu/tag/0FPL}{0FPL}, \href{https://stacks.math.columbia.edu/tag/0FPM}{0FPM}, and \href{https://stacks.math.columbia.edu/tag/0FPN}{0FPN}]{StacksProject}.
		\item \cite[$(3.6.7)(d(iii))$]{Lipman:2009} follows by applying \cite[Lemma-Definition 3.3.5]{Lipman:2009} to see that the derived functor versions of \cite[$(3.6.2)$ \& $(3.6.2)^{op}$]{Lipman:2009} hold.
		\item \cite[$3.6.7(d(iv))$]{Lipman:2009} follows from \cite[\href{https://stacks.math.columbia.edu/tag/0B6C}{Tag 0B6C}]{StacksProject}, whereas \cite[$3.6.7(d(v))$]{Lipman:2009} is \cite[\href{https://stacks.math.columbia.edu/tag/07A6}{Tag 07A6}]{StacksProject}.\qedhere
	\end{enumerate}
\end{proof}

\begin{lemma}
	\label{lem:faithfully_flat_is_conservative}
	Consider a faithfully flat morphism $f\colon Y \to X$ of algebraic spaces.
	Then $\mathbf{L}f^\ast\colon D_{\operatorname{qc}}(X)\to D_{\operatorname{qc}}(Y)$ is conservative.
\end{lemma}

\begin{proof}
	Let $E\in D_{\operatorname{qc}}(X)$ satisfy $\mathbf{L}f^\ast E\cong 0$.
	Since $f$ is flat, we have $f^\ast \mathcal{H}^j (E)\cong \mathcal{H}^j (\mathbf{L}f^\ast E)$ for all $j\in \mathbb{Z}$.
	Thus, it suffices to check that $f^\ast \mathcal{H}^j (E)\cong 0$ implies $\mathcal{H}^j (E)\cong 0$ for all $j\in \mathbb{Z}$.
	Fix an \'{e}tale presentation $s\colon U \to X$ from a scheme.
	Consider the projection morphisms $f^\prime\colon Y \times_{X} U \to U$ and $s^\prime\colon Y\times_{X} U \to Y$.
	By base change, $f^\prime$ is faithfully flat and $s^\prime$ is faithfully flat.
	Let $t\colon V \to Y\times_{X} U$ be an \'{e}tale presentation from a scheme.
	Then $t^\ast (s^\prime)^\ast f^\ast \mathcal{H}^j (E)\cong 0$ for all $j\in \mathbb{Z}$, which is equivalent to $t^\ast (f^\prime)^\ast s^\ast \mathcal{H}^j (E)\cong 0$ for all $j\in \mathbb{Z}$.
	Since $f^\prime \circ t$ is a faithfully flat morphism of schemes, it follows that $s^\ast \mathcal{H}^j (E)\cong 0$ for all $j\in \mathbb{Z}$ (see e.g.\ \cite[Proposition 14.11]{Gortz/Wedhorn:2020}).
	Finally, since $s$ is an \'{e}tale presentation, $\mathcal{H}^j (E)\cong 0$ for all $j \in \mathbb{Z}$ \cite[\href{https://stacks.math.columbia.edu/tag/07TY}{Tag 07TY}]{StacksProject}.
\end{proof}

\begin{lemma}
	\label{lem:zero_object_via_covers}
	Let $X$ be a quasi-compact quasi-separated algebraic space.
	For any $E\in D_{\operatorname{qc}}(X)$, the following are equivalent:
	\begin{enumerate}
		\item \label{lem:zero_object_via_covers1} $E\cong 0$
		\item \label{lem:zero_object_via_covers2} $\mathbf{L}f^\ast E\cong 0$ for all flat morphisms $f\colon Y\to X$ of algebraic spaces
		\item \label{lem:zero_object_via_covers3} $\mathbf{L}f^\ast E\cong 0$ for some smooth surjective morphism $f\colon Y\to X$ of algebraic spaces
		\item \label{lem:zero_object_via_covers4} there exists a flat surjective morphism $f\colon Y\to X$ of algebraic spaces such that $\mathbf{L}f^\ast E\cong 0$.
	\end{enumerate}
	Moreover, if $\mathcal{H}^j (\mathbf{L} f^\ast E)\cong 0$ for some $j$, then $\mathcal{H}^j (E)\cong 0$ whenever $f\colon Y \to X$ is faithfully flat.
\end{lemma}

\begin{proof}
	It is straightforward to check that $\eqref{lem:zero_object_via_covers1} \implies \eqref{lem:zero_object_via_covers2} \implies \eqref{lem:zero_object_via_covers3} \implies \eqref{lem:zero_object_via_covers4}$.
	%%NOTE: $\eqref{lem:zero_object_via_covers2} \implies \eqref{lem:zero_object_via_covers3}$ follows by taking $1_{X}$; $\eqref{lem:zero_object_via_covers3} \implies \eqref{lem:zero_object_via_covers4}$ follows because smooth implies flat.
	We show that $\eqref{lem:zero_object_via_covers4} \implies \eqref{lem:zero_object_via_covers1}$.
	Assume there exists a flat surjective morphism $f\colon Y\to X$ of algebraic spaces such that $\mathbf{L}f^\ast E\cong 0$.
	Let $s\colon U \to X$ be a flat surjective morphism from a scheme.
	Consider the fibered square
	\begin{displaymath}
		% https://q.uiver.app/#q=WzAsNCxbMCwwLCJcXG1hdGhjYWx7WX1cXHRpbWVzX3tcXG1hdGhjYWx7WH19IFxcbWF0aGNhbHtVfSJdLFswLDEsIlxcbWF0aGNhbHtZfSJdLFsxLDEsIlxcbWF0aGNhbHtYfS4iXSxbMSwwLCJVIl0sWzAsMSwic15cXHByaW1lIiwyXSxbMSwyLCJmIiwyXSxbMywyLCJzIl0sWzAsMywiZl5cXHByaW1lIl1d
		\begin{tikzcd}
			{Y\times_{X} U} & U \\
			{Y} & {X.}
			\arrow["{f^\prime}", from=1-1, to=1-2]
			\arrow["{s^\prime}"', from=1-1, to=2-1]
			\arrow["s", from=1-2, to=2-2]
			\arrow["f"', from=2-1, to=2-2]
		\end{tikzcd}
	\end{displaymath}
	Choose an \'{e}tale presentation $t\colon V \to Y\times_{X} U$.
	By base change, $f^\prime$ is flat and surjective.
	Hence, $f^\prime\circ t$ is a flat surjective morphism of schemes.
	Since $\mathbf{L}f^\ast E\cong 0$, it follows that $\mathbf{L}(f\circ s^\prime \circ t)^\ast E\cong 0$.
	However, $f\circ s^\prime \circ t = s\circ f^\prime \circ t$, and so $\mathbf{L}(s\circ f^\prime \circ t)^\ast E\cong 0$. By \Cref{lem:faithfully_flat_is_conservative}, it follows that $E\cong 0$ because $s\circ f^\prime \circ t$ is faithfully flat. The last claim follows from the fact that $\mathcal{H}^j (\mathbf{L} f^\ast E) \cong f^\ast \mathcal{H}^j (E)$ and \Cref{lem:faithfully_flat_is_conservative}.
\end{proof}

\begin{lemma}
	\label{lem:gortz_wedhorn_internal_hom}
	Let $f\colon Y \to X$ be a morphism of quasi-compact quasi-separated algebraic spaces.
	The canonical morphism
	\begin{displaymath}
		\mathbf{L}f^\ast \operatorname{\mathbf{R}\mathcal{H}\! \mathit{om}}(E,G) \to \operatorname{\mathbf{R}\mathcal{H}\! \mathit{om}}(\mathbf{L}f^\ast E, \mathbf{L}f^\ast G)
	\end{displaymath}
	is an isomorphism if any of the following hold:
	\begin{enumerate}
		\item \label{lem:gortz_wedhorn_internal_hom1} $E \in \operatorname{Perf}(X)$ and $G\in D_{\operatorname{qc}}(X)$
		\item \label{lem:gortz_wedhorn_internal_hom2} $E$ is pseudocoherent on $X$, $G \in D_{\operatorname{qc}}^+ (X)$, and $f$ is of finite tor-dimension.
	\end{enumerate}
\end{lemma}

\begin{proof}
	By \cite[\href{https://stacks.math.columbia.edu/tag/08JF}{Tag 08JF}]{StacksProject}, there exists a canonical morphism
	\begin{displaymath}
		\rho_f\colon \mathbf{L}f^\ast \operatorname{\mathbb{R}\mathcal{H}\! \mathit{om}}(E,G) \to \operatorname{\mathbb{R}\mathcal{H}\! \mathit{om}}(\mathbf{L}f^\ast E, \mathbf{L}f^\ast G).
	\end{displaymath}
	By \Cref{lem:internal_hom}, there exists a canonical isomorphism
	\begin{displaymath}
		\operatorname{\mathbf{R}\mathcal{H}\! \mathit{om}} (E,G) \to Q(\operatorname{\mathbb{R}\mathcal{H}\! \mathit{om}} (i(E),i(G))).
	\end{displaymath}
	In both cases, \cite[\href{https://stacks.math.columbia.edu/tag/0A8A}{Tag 0A8A}]{StacksProject} implies $\operatorname{\mathbb{R}\mathcal{H}\! \mathit{om}}(E,G)\in D_{\operatorname{qc}}$, and so we have an isomorphism
	\begin{displaymath}
		Q(\operatorname{\mathbb{R}\mathcal{H}\! \mathit{om}} (i(E),i(G))) \to \operatorname{\mathbb{R}\mathcal{H}\! \mathit{om}}(E,G).
	\end{displaymath}
	By \cite[\href{https://stacks.math.columbia.edu/tag/08H4}{Tag 08H4} and \href{https://stacks.math.columbia.edu/tag/08H6}{08H6}]{StacksProject}, $\mathbf{L}f^\ast$ preserves perfectness and pseudocoherence.
	If $f$ is of finite tor-dimension, then $\mathbf{L}f^\ast$ preserves bounded below complexes.
	Then a similar argument for both cases give a chain of isomorphisms
	\begin{displaymath}
		\operatorname{\mathbf{R}\mathcal{H}\! \mathit{om}} (\mathbf{L}f^\ast E,\mathbf{L}f^\ast G)
		\to Q(\operatorname{\mathbb{R}\mathcal{H}\! \mathit{om}} (i(\mathbf{L}f^\ast E),i(\mathbf{L}f^\ast G))) \to\operatorname{\mathbb{R}\mathcal{H}\! \mathit{om}}(\mathbf{L}f^\ast E,\mathbf{L}f^\ast G).
	\end{displaymath}
	Hence, we have a morphism
	\begin{displaymath}
		% https://q.uiver.app/#q=WzAsNCxbMCwxLCJcXG1hdGhiZntMfWZeXFxhc3QgXFxvcGVyYXRvcm5hbWV7XFxtYXRoYmZ7Un1cXG1hdGhjYWx7SH1cXCEgXFxtYXRoaXR7b219fShFLEcpIl0sWzAsMCwiXFxtYXRoYmZ7TH1mXlxcYXN0IFxcb3BlcmF0b3JuYW1le1xcbWF0aGJie1J9XFxtYXRoY2Fse0h9XFwhIFxcbWF0aGl0e29tfX0oRSxHKSJdLFsxLDAsIlxcb3BlcmF0b3JuYW1le1xcbWF0aGJie1J9XFxtYXRoY2Fse0h9XFwhIFxcbWF0aGl0e29tfX0gKFxcbWF0aGJme0x9Zl5cXGFzdCBFLFxcbWF0aGJme0x9Zl5cXGFzdCBHKSAiXSxbMSwxLCJcXG9wZXJhdG9ybmFtZXtcXG1hdGhiZntSfVxcbWF0aGNhbHtIfVxcISBcXG1hdGhpdHtvbX19IChcXG1hdGhiZntMfWZeXFxhc3QgRSxcXG1hdGhiZntMfWZeXFxhc3QgRykgLiJdLFswLDEsIlxcY29uZyJdLFsxLDIsIlxccmhvX2YiXSxbMiwzLCJcXGNvbmciXV0=
		\begin{tikzcd}
			{\mathbf{L}f^\ast \operatorname{\mathbb{R}\mathcal{H}\! \mathit{om}}(E,G)} & {\operatorname{\mathbb{R}\mathcal{H}\! \mathit{om}} (\mathbf{L}f^\ast E,\mathbf{L}f^\ast G) } \\
			{\mathbf{L}f^\ast \operatorname{\mathbf{R}\mathcal{H}\! \mathit{om}}(E,G)} & {\operatorname{\mathbf{R}\mathcal{H}\! \mathit{om}} (\mathbf{L}f^\ast E,\mathbf{L}f^\ast G) .}
			\arrow["{\rho_f}", from=1-1, to=1-2]
			\arrow["\cong", from=1-2, to=2-2]
			\arrow["\cong", from=2-1, to=1-1]
		\end{tikzcd}
	\end{displaymath}
	It suffices to prove that $\rho_f$ is an isomorphism.

	Let $p\colon U \to X$ be an \'{e}tale presentation.
	Form the fibered square
	\begin{displaymath}
		% https://q.uiver.app/#q=WzAsNCxbMCwwLCJZXFx0aW1lc19YIFUiXSxbMCwxLCJZIl0sWzEsMCwiVSJdLFsxLDEsIlguIl0sWzAsMSwicF5cXHByaW1lIiwyXSxbMiwzLCJwIl0sWzEsMywiZiIsMl0sWzAsMiwiZl5cXHByaW1lIl1d
		\begin{tikzcd}
			{Y\times_X U} & U \\
			Y & {X.}
			\arrow["{f^\prime}", from=1-1, to=1-2]
			\arrow["{p^\prime}"', from=1-1, to=2-1]
			\arrow["p", from=1-2, to=2-2]
			\arrow["f"', from=2-1, to=2-2]
		\end{tikzcd}
	\end{displaymath}
	Choose an \'{e}tale presentation $p^{\prime\prime}\colon V\to Y\times_X U$.
	Set $q := p^\prime \circ p^{\prime\prime}$ and $g := f^\prime \circ p^{\prime\prime}$.
	Note that $q$ is an \'{e}tale presentation.
	By \Cref{lem:spacey_Lipman_Hashimoto} and \cite[Exercise 3.7.1.1]{Lipman:2009}, there exists a commutative diagram,
	\begin{displaymath}
		% https://q.uiver.app/#q=WzAsNixbMCwwLCJcXG1hdGhiZntMfXFeXFxhc3QgXFxtYXRoYmZ7TH0gZl5cXGFzdCBcXG9wZXJhdG9ybmFtZXtcXG1hdGhiYntSfVxcbWF0aGNhbHtIfVxcISBcXG1hdGhpdHtvbX19KEUsRykiXSxbMSwwLCJcXG1hdGhiZntMfXFeXFxhc3QgXFxvcGVyYXRvcm5hbWV7XFxtYXRoYmJ7Un1cXG1hdGhjYWx7SH1cXCEgXFxtYXRoaXR7b219fShcXG1hdGhiZntMfWZeXFxhc3QgRSwgXFxtYXRoYmZ7TH1mXlxcYXN0IEcpIl0sWzEsMSwiXFxvcGVyYXRvcm5hbWV7XFxtYXRoYmJ7Un1cXG1hdGhjYWx7SH1cXCEgXFxtYXRoaXR7b219fShcXG1hdGhiZntMfXFeXFxhc3QgXFxtYXRoYmZ7TH0gZl5cXGFzdCBFLFxcbWF0aGJme0x9cV5cXGFzdCBcXG1hdGhiZntMfSBmXlxcYXN0IEcpIl0sWzEsMiwiXFxvcGVyYXRvcm5hbWV7XFxtYXRoYmJ7Un1cXG1hdGhjYWx7SH1cXCEgXFxtYXRoaXR7b219fShcXG1hdGhiZntMfWdeXFxhc3QgXFxtYXRoYmZ7TH0gcF5cXGFzdCBFLFxcbWF0aGJme0x9Z15cXGFzdCBcXG1hdGhiZntMfSBwXlxcYXN0IEcpLiJdLFswLDEsIlxcbWF0aGJme0x9Z15cXGFzdCBcXG1hdGhiZntMfSBwXlxcYXN0IFxcb3BlcmF0b3JuYW1le1xcbWF0aGJie1J9XFxtYXRoY2Fse0h9XFwhIFxcbWF0aGl0e29tfX0oRSxHKSJdLFswLDIsIlxcbWF0aGJme0x9Z15cXGFzdCBcXG9wZXJhdG9ybmFtZXtcXG1hdGhiYntSfVxcbWF0aGNhbHtIfVxcISBcXG1hdGhpdHtvbX19KFxcbWF0aGJme0x9cF5cXGFzdCBFLCBcXG1hdGhiZntMfXBeXFxhc3QgRykiXSxbMCwxLCJcXG1hdGhiZntMfXFeXFxhc3QgXFxyaG9fZiJdLFsxLDIsIlxccmhvX3EiXSxbMiwzLCJcXGNvbmciXSxbNSwzLCJcXHJob19nIiwyXSxbMCw0LCJcXGNvbmciLDJdLFs0LDUsIlxcbWF0aGJme0x9Z15cXGFzdCBcXHJob19wIiwyXV0=
		\begin{tikzcd}
			{\mathbf{L}q^\ast \mathbf{L} f^\ast \operatorname{\mathbb{R}\mathcal{H}\! \mathit{om}}(E,G)} & {\mathbf{L}q^\ast \operatorname{\mathbb{R}\mathcal{H}\! \mathit{om}}(\mathbf{L}f^\ast E, \mathbf{L}f^\ast G)} \\
			{\mathbf{L}g^\ast \mathbf{L} p^\ast \operatorname{\mathbb{R}\mathcal{H}\! \mathit{om}}(E,G)} & {\operatorname{\mathbb{R}\mathcal{H}\! \mathit{om}}(\mathbf{L}q^\ast \mathbf{L} f^\ast E,\mathbf{L}q^\ast \mathbf{L} f^\ast G)} \\
			{\mathbf{L}g^\ast \operatorname{\mathbb{R}\mathcal{H}\! \mathit{om}}(\mathbf{L}p^\ast E, \mathbf{L}p^\ast G)} & {\operatorname{\mathbb{R}\mathcal{H}\! \mathit{om}}(\mathbf{L}g^\ast \mathbf{L} p^\ast E,\mathbf{L}g^\ast \mathbf{L} p^\ast G).}
			\arrow["{\mathbf{L}q^\ast \rho_f}", from=1-1, to=1-2]
			\arrow["\cong"', from=1-1, to=2-1]
			\arrow["{\rho_q}", from=1-2, to=2-2]
			\arrow["{\mathbf{L}g^\ast \rho_p}"', from=2-1, to=3-1]
			\arrow["\cong", from=2-2, to=3-2]
			\arrow["{\rho_g}"', from=3-1, to=3-2]
		\end{tikzcd}
	\end{displaymath}
	By \cite[\href{https://stacks.math.columbia.edu/tag/04LX}{Tag 04LX} and \href{https://stacks.math.columbia.edu/tag/08JB}{08JB}]{StacksProject}, $\rho_q$ and $\rho_p$ are isomorphisms.
	Since $q$ is faithfully flat, \Cref{lem:zero_object_via_covers} says it suffices to prove that $\rho_g$ is an isomorphism.
	Indeed, this implies $\mathbf{L}q^\ast \rho_f$ is an isomorphism, and hence $\mathbf{L}q^\ast \operatorname{cone}(\rho_f)\cong 0$.
	In such a case, we would obtain that $\operatorname{cone}(\rho_f)\cong 0$, and hence, $\rho_f$ is an isomorphism.
	% If $f$ is \'{e}tale, then by change, so is $f^\prime$.
	% Hence, $g$ is \'{e}tale with scheme source-and-target.
	% By \cite[\href{https://stacks.math.columbia.edu/tag/04LX}{Tags 04LX} \& \href{https://stacks.math.columbia.edu/tag/08JB}{08JB}]{StacksProject}, $\rho_g$ is an isomorphism whenever $f$ is \'{e}tale.
	% This proves \eqref{lem:gortz_wedhorn_internal_hom1}.
	Applying \cite[\href{https://stacks.math.columbia.edu/tag/08GH}{Tag 08GH}]{StacksProject}, we can reduce to checking the claim on the small Zariski sites of $U$ and $V$.
	Therefore, the desired claim follows from \cite[Proposition 22.70]{Gortz/Wedhorn:2023}.
\end{proof}

The following shows where \Cref{lem:gortz_wedhorn_internal_hom} can fail beyond the mentioned cases.

\begin{example}
	Let $k$ be a field.
	Set $A:=k[\varepsilon]/(\varepsilon^2)$, $f:=1_{\operatorname{Spec}(A)}\colon \operatorname{Spec}(A)\to\operatorname{Spec}(A)$, $E:=k=A/(\varepsilon)$, and $F:=A$.
	Then $f$ is flat and finitely presented, $E$ is pseudocoherent,
	and $F$ is perfect relative to $\operatorname{Spec}(A)$.
	Consider the closed immersion $i\colon \operatorname{Spec}(k)\to\operatorname{Spec}(A)$.
	There exists a periodic free resolution
	\begin{displaymath}
		\cdots
		\xrightarrow{\varepsilon}
		A
		\xrightarrow{\varepsilon}
		A
		\xrightarrow{\varepsilon}
		A
		\to k
		\to 0.
	\end{displaymath}
	Applying $\operatorname{Hom}_A(-,A)$ gives the complex
	\begin{displaymath}
		0
		\to A
		\xrightarrow{\varepsilon} A
		\xrightarrow{\varepsilon} A
		\xrightarrow{\varepsilon}
		\cdots.
	\end{displaymath}
	Its degree zero cohomology is $k$, whereas all its positive cohomology vanishes.
	Consequently, $\mathbf{R}\operatorname{Hom}_A(k,A) \cong k$.
	It follows that $\mathbf{L}i^\ast \mathbf{R}\operatorname{Hom}_A(k,A) \cong k\otimes_A^{\mathbf{L}}k$.
	Tensoring the periodic resolution of $k$ with $k$ shows that $k\otimes_A^{\mathbf{L}}k$ is represented by
	\begin{displaymath}
		\cdots
		\xrightarrow{0}
		k
		\xrightarrow{0}
		k
		\xrightarrow{0}
		k
		\to 0
	\end{displaymath}
	with copies of $k$ in nonpositive degrees.
	Hence, for every $j\geq 0$,
	\begin{displaymath}
		\mathcal{H}^{-j} ( \mathbf{L}i^\ast  \mathbf{R}\operatorname{Hom}_A(k,A)) \cong k.
	\end{displaymath}
	On the other hand, $\mathbf{L}i^\ast k \cong k\otimes_A^{\mathbf{L}}k$
	and $\mathbf{L}i^\ast A\cong k$.
	Since $k$ is a field, $\mathbf{R}\operatorname{Hom}_k ( k\otimes_A^{\mathbf{L}}k,  k)$ is represented by
	\begin{displaymath}
		0
		\to k
		\xrightarrow{0} k
		\xrightarrow{0} k
		\xrightarrow{0} \cdots
	\end{displaymath}
	with copies of $k$ in nonnegative degrees.
	In particular,
	\begin{displaymath}
		\mathcal{H}^{-1}(\mathbf{R}\operatorname{Hom}_k ( \mathbf{L}i^\ast k, \mathbf{L}i^\ast A))
		\cong 0,
	\end{displaymath}
	whereas
	\begin{displaymath}
		\mathcal{H}^{-1} (\mathbf{L}i^\ast \mathbf{R}\operatorname{Hom}_A(k,A)
		) \cong k.
	\end{displaymath}
	Thus, the canonical base change morphism
	\begin{displaymath}
		\mathbf{L}i^\ast
		\mathbf{R}\operatorname{Hom}_A(k,A)
		\to
		\mathbf{R}\operatorname{Hom}_k
		(\mathbf{L}i^\ast k, \mathbf{L}i^\ast A)
	\end{displaymath}
	is not an isomorphism.
\end{example}

%%%%%%%%%%%%%%%%%%%%%%%%%%%%%%%%%%%
\subsubsection{Quasi-affine diagonal}
\label{sec:preliminaries_quasi-affine_diagonal}
%%%%%%%%%%%%%%%%%%%%%%%%%%%%%%%%%%%

\begin{lemma}
	\label{lem:quasi-affine_via_diagonal}
	Consider a commutative diagram of algebraic spaces
	\begin{displaymath}
		% https://q.uiver.app/#q=WzAsMyxbMCwwLCJcXG1hdGhjYWx7Wn0iXSxbMSwwLCJcXG1hdGhjYWx7WX0iXSxbMSwxLCJcXG1hdGhjYWx7WH0uIl0sWzAsMiwiZyIsMl0sWzEsMiwiZiJdLFswLDEsImgiXV0=
		\begin{tikzcd}
			{Z} & {Y} \\
			& {X.}
			\arrow["g", from=1-1, to=1-2]
			\arrow["h"', from=1-1, to=2-2]
			\arrow["f", from=1-2, to=2-2]
		\end{tikzcd}
	\end{displaymath}
	If both $h$ and $\Delta_f$ are quasi-affine, then $g$ is quasi-affine.
\end{lemma}

\begin{proof}
	There exists a fibered square
	\begin{displaymath}
		% https://q.uiver.app/#q=WzAsNCxbMCwxLCJcXG1hdGhjYWx7Wn0iXSxbMSwwLCJcXG1hdGhjYWx7WX0iXSxbMSwxLCJcXG1hdGhjYWx7WH0uIl0sWzAsMCwiXFxtYXRoY2Fse1p9XFx0aW1lc197XFxtYXRoY2Fse1h9fSBcXG1hdGhjYWx7WX0iXSxbMCwyLCJnIiwyXSxbMSwyLCJmIl0sWzMsMSwicF8yIl0sWzMsMCwicF8xIiwyXV0=
		\begin{tikzcd}
			{Z\times_{X} Y} & {Y} \\
			{Z} & {X.}
			\arrow["{p_2}", from=1-1, to=1-2]
			\arrow["{p_1}"', from=1-1, to=2-1]
			\arrow["f", from=1-2, to=2-2]
			\arrow["h"', from=2-1, to=2-2]
		\end{tikzcd}
	\end{displaymath}
	By base change, $p_2$ is quasi-affine \cite[\href{https://stacks.math.columbia.edu/tag/03WO}{Tag 03WO}]{StacksProject}.
	Moreover, the morphism $(1,g)\colon Z \to Z\times_{X} Y$ is the base change of the diagonal $\Delta_f\colon Y\to Y\times_{X}Y$ by the morphism $Z\times_{X}Y \to Y\times_{X}Y$.
	%%NOTE: See Lemma 1 of 'Magic Squares' in notes folder or Gortz--Wedhorn, Volume I, Prop 9.3(2).
	Again, by base change, $(1,g)\colon Z \to Z\times_{X} Y$ is quasi-affine.
	Since $g= p_2 \circ (1,g)$, it must be itself quasi-affine \cite[\href{https://stacks.math.columbia.edu/tag/03WN}{Tag 03WN}]{StacksProject}.
\end{proof}

\begin{lemma}
	\label{lem:quasi_affine_diagonal}
	Let $X$ be a quasi-separated algebraic space.
	Then every morphism to $X$ from a quasi-affine scheme is quasi-affine.
\end{lemma}

\begin{proof}
	Since quasi-separated algebraic spaces have quasi-affine diagonal \cite[\href{https://stacks.math.columbia.edu/tag/03HK}{Tag 03HK}, \href{https://stacks.math.columbia.edu/tag/05W7}{05W7}, and \href{https://stacks.math.columbia.edu/tag/082J}{082J}]{StacksProject}, the claim follows from \Cref{lem:quasi-affine_via_diagonal}.
\end{proof}

%%%%%%%%%%%%%%%%%%%%%%%%%%%%%%%%%%%%%
\subsection{Integral transforms and base change notation}
\label{sec:prelim_fm_setup}
%%%%%%%%%%%%%%%%%%%%%%%%%%%%%%%%%%%%%

For a point $p\in |S|$, a representative is a morphism $t\colon\operatorname{Spec}(k)\to S$ in the equivalence class defining $p$.
Denote by $\kappa(p)$ its residue field.
Set $\mathbf{L}(t^\prime)^\ast K$ to be the base changed kernel.
We use the following throughout. 

\begin{setup}
	\label{setup:fm_cube_pullback}
	Let $t\colon T \to S$ be a morphism of Noetherian algebraic spaces.
	Suppose $f_1\colon Y_1 \to S$ and $f_2\colon Y_2 \to S$ are proper flat morphisms of algebraic spaces.
	Consider the commutative diagram
	\begin{equation}
		\label{diag:fm_cube_pullback}
		% https://q.uiver.app/#q=WzAsOCxbNCwxLCJcXG1hdGhjYWx7VH0iXSxbNCwzLCJcXG1hdGhjYWx7U30iXSxbMSwzLCJcXG1hdGhjYWx7WX1fMSJdLFsyLDIsIlxcbWF0aGNhbHtZfV8yIl0sWzEsMSwiXFxtYXRoY2Fse1l9XzFcXHRpbWVzX3tcXG1hdGhjYWx7U319IFxcbWF0aGNhbHtUfSJdLFsyLDAsIlxcbWF0aGNhbHtZfV8yXFx0aW1lc19cXG1hdGhjYWx7U30gXFxtYXRoY2Fse1R9Il0sWzAsMCwiXFxtYXRoY2Fse1l9XzFcXHRpbWVzX1xcbWF0aGNhbHtTfSBcXG1hdGhjYWx7WX1fMlxcdGltZXNfXFxtYXRoY2Fse1N9IFxcbWF0aGNhbHtUfSJdLFswLDIsIlxcbWF0aGNhbHtZfV8xXFx0aW1lc197XFxtYXRoY2Fse1N9fSBcXG1hdGhjYWx7WX1fMiJdLFswLDEsInQiXSxbMiwxLCJmXzEiLDJdLFszLDEsImZfMiJdLFs0LDAsImdfMSJdLFs1LDAsImdfMiJdLFs0LDIsInRfMSIsMCx7ImxhYmVsX3Bvc2l0aW9uIjo3MH1dLFs1LDMsInRfMiIsMix7ImxhYmVsX3Bvc2l0aW9uIjo4MCwic3R5bGUiOnsiYm9keSI6eyJuYW1lIjoiZGFzaGVkIn19fV0sWzYsNCwiZ18yXlxccHJpbWUiLDJdLFs2LDUsImdfMV5cXHByaW1lIiwyXSxbNywzLCJmXlxccHJpbWVfMSIsMCx7ImxhYmVsX3Bvc2l0aW9uIjozMCwic3R5bGUiOnsiYm9keSI6eyJuYW1lIjoiZGFzaGVkIn19fV0sWzcsMiwiZl5cXHByaW1lXzIiLDJdLFs2LDcsInReXFxwcmltZSIsMl1d
		\begin{tikzcd}
			{Y_1\times_S Y_2\times_S T} && {Y_2\times_S T} && \\
			& {Y_1\times_{S} T} &&& {T} \\
			{Y_1\times_{S} Y_2} && {Y_2} \\
			& {Y_1} &&& {S}
			\arrow["{g_1^\prime}"', from=1-1, to=1-3]
			\arrow["{g_2^\prime}"', from=1-1, to=2-2]
			\arrow["{t^\prime}"', from=1-1, to=3-1]
			\arrow["{g_2}", from=1-3, to=2-5]
			\arrow["{t_2}"'{pos=0.8}, dashed, from=1-3, to=3-3]
			\arrow["{g_1}", from=2-2, to=2-5]
			\arrow["{t_1}"{pos=0.7}, from=2-2, to=4-2]
			\arrow["t", from=2-5, to=4-5]
			\arrow["{f^\prime_1}"{pos=0.3}, dashed, from=3-1, to=3-3]
			\arrow["{f^\prime_2}"', from=3-1, to=4-2]
			\arrow["{f_2}", from=3-3, to=4-5]
			\arrow["{f_1}"', from=4-2, to=4-5]
		\end{tikzcd}
	\end{equation}
	obtained by base change along $t$.
\end{setup}

%%%%%%%%%%%%%%%%%%%%%%%%%%%%%%%%%%%
\subsection{Duality for spaces}
\label{sec:duality}
%%%%%%%%%%%%%%%%%%%%%%%%%%%%%%%%%%%

We use Neeman's construction of Grothendieck duality \cite{Neeman:2023} directly for $D_{\operatorname{qc}}$ defined on the small \'{e}tale sites.
See also \cite[\href{https://stacks.math.columbia.edu/tag/0E4V}{Tag 0E4V}]{StacksProject}.
The pullback-pushforward formalism and the compact generation results used in the construction of loc.\ cit.\ are taken on these sites.
The resulting functors and comparison morphisms have the properties recorded.

%%%%%%%%%%%%%%%%%%%%%%%%%%%%%%%%%%%
\subsubsection{Categories}
\label{sec:duality_2cats}
%%%%%%%%%%%%%%%%%%%%%%%%%%%%%%%%%%%

Fix a base scheme $S$, regarded as an object of the big fppf site \cite[\href{https://stacks.math.columbia.edu/tag/021R}{Tag 021R}]{StacksProject}.
Let $(\operatorname{Sch}/S)_{fppf}$ be the big fppf site of $S$ \cite[\href{https://stacks.math.columbia.edu/tag/021S}{Tag 021S}]{StacksProject}.
We define the $(2,1)$-category $\operatorname{Sp}$ of algebraic spaces over $S$ as follows.

Every algebraic space $X$ over $S$ determines a category $\mathcal{S}_X$ fibered in sets over $(\operatorname{Sch}/S)_{fppf}$ \cite[\href{https://stacks.math.columbia.edu/tag/04TM}{Tag 04TM} and \href{https://stacks.math.columbia.edu/tag/02Y2}{02Y2}]{StacksProject}, and hence  an object in the $2$-category of categories fibered over $(\operatorname{Sch}/S)_{fppf}$ \cite[\href{https://stacks.math.columbia.edu/tag/04S8}{Tag 04S8}]{StacksProject}.

A morphism of algebraic spaces $f\colon Y \to X$ induces a $1$-morphism $\mathcal{S}_Y \to \mathcal{S}_X$ in the $2$-category of categories fibered over $(\operatorname{Sch}/S)_{fppf}$.
Thus algebraic spaces over $S$ form a $(2,1)$-category, denoted $\operatorname{Sp}/S$, whose $2$-morphisms are invertible \cite[\href{https://stacks.math.columbia.edu/tag/003S}{Tag 003S}]{StacksProject}.
In particular, $\operatorname{Sp}/S$ is a $2$-subcategory of the $(2,1)$-category of algebraic stacks over $S$ \cite[\href{https://stacks.math.columbia.edu/tag/03YP}{Tag 03YP}]{StacksProject}.

Denote by $(\operatorname{Sp}/S)_{Noeth,sep,fp}$ the $2$-subcategory of $\operatorname{Sp}/S$ consisting of Noetherian algebraic spaces, separated finitely presented morphisms, and all $2$-morphisms.
Every $1$-morphism $f\colon Y \to X$ of $(\operatorname{Sp}/S)_{Noeth,sep,fp}$ admits a Nagata compactification $f=p\circ j$ where $p\colon Y^\prime \to X$ is proper and $j\colon Y\to Y^\prime$ is an open immersion \cite{Conrad/Lieblich/Olsson:2012,Raoult:1971, Raoult:1974}.
Let $\operatorname{TriCat}$ be the $2$-category whose objects are triangulated categories, $1$-morphisms are exact functors, and $2$-morphisms are natural transformations.

%%%%%%%%%%%%%%%%%%%%%%%%%%%%%%%%%%%
\subsubsection{Functors}
\label{sec:duality_functors}
%%%%%%%%%%%%%%%%%%%%%%%%%%%%%%%%%%%

There exist two contravariant pseudofunctors
\begin{displaymath}
	\mathbf{L}(-)^\ast, (-)^\times\colon (\operatorname{Sp}/S)_{Noeth,sep,fp} \to \operatorname{TriCat},
\end{displaymath}
and a functor
\begin{displaymath}
	(-)^!\colon (\operatorname{Sp}/S)_{Noeth,sep,fp} \to \operatorname{TriCat}
\end{displaymath}
which is contravariant and oplax. We define them as follows:
\begin{itemize}
	\item On objects: $\mathbf{L}(X)^\ast = (X)^\times = (X)^!=D_{\operatorname{qc}}(X)$ for all objects $X$ in $(\operatorname{Sp}/S)_{Noeth,sep,fp}$.
	\item On $1$-morphisms: For a $1$-morphism $f\colon Y\to X$ in $(\operatorname{Sp}/S)_{Noeth,sep,fp}$, $\mathbf{L}f^\ast$ is the derived pullback functor, $f^\times$ is the right adjoint to pushforward $\mathbf{R}f_\ast$, and $f^!$ is defined as in \cite[Construction 9.4]{Neeman:2023} (if $f=p\circ j$ is a Nagata compactification, then $f^! := \mathbf{L} j^\ast \circ p^\times$).
	\item On $2$-morphisms: Since $(\operatorname{Sp}/S)_{Noeth,sep,fp}$ is a $(2,1)$-category, each $2$-morphism is invertible, and so we obtain natural isomorphisms between functors.
\end{itemize}
The functors $\mathbf{L}(-)^\ast$ and $(-)^\times$ are pseudofunctorial; that is, if $f\colon Z \to Y$ and $g\colon Y \to X$ are composable $1$-morphisms, then there exist natural isomorphisms $\tau(f,g)\colon \mathbf{L}(g\circ f)^\ast \implies \mathbf{L}f^\ast \circ \mathbf{L}g^\ast$ and $\delta(f,g)\colon  (g\circ f)^\times \implies f^\times \circ g^\times$.
In contrast, $(-)^!$ is not in general pseudofunctorial.
Rather, $(-)^!$ is oplax, equipped with natural morphisms $\rho(f,g)\colon (g\circ f)^! \implies f^! \circ g^!$. Moreover, $(-)^\times$ and $(-)^!$ are oplax modules over $\mathbf{L}(-)^\ast$. This means there exist oplax natural transformations $\xi\colon \mathbf{L}(-)^\ast \times (-)^\times\to (-)^\times$ and $\sigma\colon \mathbf{L}(-)^\ast \times (-)^! \to (-)^!$ satisfying an associativity property. See \cite[Theorem 1.8(i \textrm{to} v)]{Neeman:2023}.

%%%%%%%%%%%%%%%%%%%%%%%%%%%%%%%%%%%
\subsubsection{Relations}
\label{sec:duality_relations}
%%%%%%%%%%%%%%%%%%%%%%%%%%%%%%%%%%%

We now collect base change formulas and properties of these functors. The following results are \cite[Theorem 1.8.6 to 1.8.9 \& Lemma 0.1]{Neeman:2023}.

\begin{lemma}
	\label{lem:neeman186}
	Let $f\colon Z \to Y$ and $g\colon Y \to X$ be a pair of composable $1$-morphisms in $(\operatorname{Sp}/S)_{Noeth,sep,fp}$. Then the morphism $\rho(f,g)\colon (g \circ f)^! \to f^! \circ g^!$ is an isomorphism if one of the following holds:
	\begin{enumerate}
		\item $f$ is of finite tor-dimension (e.g.\ flat)
		\item $g$ is proper
		\item $g\circ f$ is proper
		\item We restrict to the subcategory $D^+_{\operatorname{qc}}(X)$.
	\end{enumerate}
\end{lemma}

\begin{lemma}
	\label{lem:neeman187}
	Let $f\colon Y \to X$ be a $1$-morphism in $(\operatorname{Sp}/S)_{Noeth,sep,fp}$. If $f$ is proper, then there exists a natural isomorphism $f^\times \to f^!$ on $D_{\operatorname{qc}}$.
\end{lemma}

\begin{lemma}
	\label{lem:neeman188}
	Consider a $2$-fibered square in $(\operatorname{Sp}/S)_{Noeth,sep,fp}$
	\begin{displaymath}
		% https://q.uiver.app/#q=WzAsNCxbMCwwLCJXIl0sWzEsMCwiWCJdLFsxLDEsIloiXSxbMCwxLCJZIl0sWzAsMSwidSJdLFsxLDIsImciXSxbMCwzLCJmIiwyXSxbMywyLCJ2IiwyXV0=
		\begin{tikzcd}
			W & X \\
			Y & Z
			\arrow["u", from=1-1, to=1-2]
			\arrow["f"', from=1-1, to=2-1]
			\arrow["g", from=1-2, to=2-2]
			\arrow["v"', from=2-1, to=2-2]
		\end{tikzcd}
	\end{displaymath}
	where $u,v$ are flat. Then the base change morphism $\mathbf{L}u^\ast g^! \to f^! \mathbf{L} v^\ast$ is an isomorphism if one of the following conditions holds:
	\begin{enumerate}
		\item $f$ is of finite tor-dimension
		\item We restrict to the subcategory $D^+_{\operatorname{qc}}(Z)$.
	\end{enumerate}
	More generally, consider a fibered square of Noetherian algebraic spaces
	\begin{displaymath}
		% https://q.uiver.app/#q=WzAsNCxbMCwwLCJXIl0sWzEsMCwiWCJdLFsxLDEsIloiXSxbMCwxLCJZIl0sWzAsMSwidSJdLFsxLDIsImciXSxbMCwzLCJmIiwyXSxbMywyLCJ2IiwyXV0=
		\begin{tikzcd}
			W & X \\
			Y & Z
			\arrow["u", from=1-1, to=1-2]
			\arrow["f"', from=1-1, to=2-1]
			\arrow["g", from=1-2, to=2-2]
			\arrow["v"', from=2-1, to=2-2]
		\end{tikzcd}
	\end{displaymath}
	where $v$ is flat and $g$ is proper. Then the base change morphism $\mathbf{L}u^\ast g^\times E \to f^\times \mathbf{L} v^\ast E$ is an isomorphism if one of the following conditions holds:
	\begin{enumerate}
		\item $f$ is of finite tor-dimension and $E\in D_{\operatorname{qc}}(Z)$
		\item We restrict to the subcategory $D^+_{\operatorname{qc}}(Z)$.
	\end{enumerate}
\end{lemma}

\begin{lemma}
	\label{lem:neeman189}
	Let $f\colon Y \to X$ be a $1$-morphism in $(\operatorname{Sp}/S)_{Noeth,sep,fp}$.
	Choose $E,G\in D_{\operatorname{qc}}(X)$.
	There exists bifunctorial morphisms
	\begin{displaymath}
		\chi(f,E,G)\colon \mathbf{L}f^\ast E \otimes^{\mathbf{L}} f^\times G \to f^\times (E\otimes^{\mathbf{L}} G)
	\end{displaymath}
	and
	\begin{displaymath}
		\sigma(f,E,G)\colon \mathbf{L}f^\ast E \otimes^{\mathbf{L}} f^! G \to f^! (E\otimes^{\mathbf{L}} G).
	\end{displaymath}
	If $E$ is perfect, then both $\chi(f,E,G)$ and $\sigma(f,E,G)$ are isomorphisms.
	If $f$ is of finite tor-dimension (e.g.\ flat), then $\sigma(f,E,G)$ is an isomorphism.
\end{lemma}

%%%%%%%%%%%%%%%%%%%%%%%%%%%%%%%%%%%
\part{Relatively perfect complexes on algebraic spaces}
\label{sec:part1}
%%%%%%%%%%%%%%%%%%%%%%%%%%%%%%%%%%%

%%%%%%%%%%%%%%%%%%%%%%%%%%%%%%%%%%%
\section{Characterization of perfectness}
\label{sec:perfectness}
%%%%%%%%%%%%%%%%%%%%%%%%%%%%%%%%%%%

\begin{lemma}
	\label{lem:23}
	Let $X$ be a quasi-compact quasi-separated algebraic space.
	Then $\mathbf{R}i_\ast E\in D^b_{\operatorname{qc}}(X)$ for all $E\in D^b_{\operatorname{qc}}(\operatorname{Spec}(k))$ where $i\colon \operatorname{Spec}(k) \to X$ represents some $p\in |X|$.
\end{lemma}

\begin{proof}
	Since $X$ is quasi-separated it follows that $\mathbf{R} i_\ast \mathcal{O}_{\operatorname{Spec}(k)}\in D^b_{\operatorname{qc}}(X)$.
	Indeed, $i$ must be an affine morphism \cite[\href{https://stacks.math.columbia.edu/tag/09TF}{Tag 09TF}]{StacksProject}, so $\mathbf{R}^j i_\ast \mathcal{O}_{\operatorname{Spec}(k)}\cong 0$ for $j\not=0$.
	%%NOTE: This can be checked on schemes. Indeed, $i$ is representable by schemes, so apply Tag 0G9R.
	Now, any $E\in D^b_{\operatorname{qc}}(\operatorname{Spec}(k))$ is a small coproduct of shifts of $\mathcal{O}_{\operatorname{Spec}(k)}$.
	In particular, $E\cong \oplus_{\alpha\in I} \mathcal{O}_{\operatorname{Spec}(k)}^{\oplus \alpha}[n_\alpha]$ where $I$ is some finite set of cardinals.
	As $\mathbf{R} i_\ast \mathcal{O}_{\operatorname{Spec}(k)}^{\oplus \alpha} \in D^b_{\operatorname{qc}}(X)$ for each $\alpha\in I$, it follows that $\mathbf{R} i_\ast E \in D^b_{\operatorname{qc}}(X)$.
\end{proof}

\begin{lemma}
	\label{lem:quasi-affine_is_conservative}
	Let $f\colon Y \to X$ be a quasi-affine morphism of quasi-compact quasi-separated algebraic spaces.
	Then $\mathbf{R}f_\ast\colon D_{\operatorname{qc}}(Y) \to D_{\operatorname{qc}}(X)$ is conservative (i.e.\ $\mathbf{R}f_\ast E \cong 0$ implies $E\cong 0$).
\end{lemma}

\begin{proof}
	There exists the canonical factorization
	\begin{displaymath}
		Y \xrightarrow{j} \underline{\operatorname{Spec}}_X (f_\ast \mathcal{O}_Y) \xrightarrow{f^\prime} X.
	\end{displaymath}
	See \cite[\href{https://stacks.math.columbia.edu/tag/081X}{Tag 081X}]{StacksProject}.
	Since $f^\prime$ is affine, $\mathbf{R}f^\prime_\ast$ is conservative \cite[\href{https://stacks.math.columbia.edu/tag/0E4R}{Tag 0E4R}]{StacksProject}.
	Moreover, $f$ is quasi-affine, $Y \to \underline{\operatorname{Spec}}_X (f_\ast \mathcal{O}_Y)$ is an open immersion \cite[\href{https://stacks.math.columbia.edu/tag/086S}{Tag 086S}]{StacksProject}.
	Hence, $\mathbf{R}j_\ast$ is conservative since $\mathbf{L}j^\ast \mathbf{R}j_\ast E \cong E$ for all $E\in D_{\operatorname{qc}}(Y)$.
	%%NOTE: One way to see this is look at flat base change of $j$ along $j$, fiber product is $Y$ and projections are identities.
	Thus, the claim follows.
\end{proof}

\begin{lemma}
	\label{lem:34}
	Let $X$ be a quasi-compact quasi-separated algebraic space.
	Suppose $i\colon \operatorname{Spec}(k) \to X$ represents some $p\in |X|$.
	If $E\otimes^{\mathbf{L}} \mathbf{R} i_\ast \mathcal{O}_{\operatorname{Spec}(k)} \in D^b_{\operatorname{qc}}(X)$, then $\mathbf{L}i^\ast E\in D^b_{\operatorname{qc}}(\operatorname{Spec}(k))$.
\end{lemma}

\begin{proof}
	The projection formula yields
	\begin{displaymath}
		\mathbf{R} i_\ast \mathbf{L}i^\ast E  \cong E\otimes^{\mathbf{L}} \mathbf{R} i_\ast \mathcal{O}_{\operatorname{Spec}(k)}.
	\end{displaymath}
	Hence, by hypothesis, $\mathbf{R} i_\ast \mathbf{L}i^\ast E \in D^b_{\operatorname{qc}}(X)$.
	Moreover, affineness of $i$ implies $\mathcal{H}^j (\mathbf{R} i_\ast A) \cong i_\ast \mathcal{H}^j (A)$ for all $j\in \mathbb{Z}$.
	Thus, \Cref{lem:quasi-affine_is_conservative} gives the boundedness of $\mathbf{L}i^\ast E$.
\end{proof}

\begin{lemma}
	\label{lem:45}
	Let $X$ be a quasi-compact quasi-separated algebraic space.
	Choose $E\in D^b_{\operatorname{qc}}(X)$ pseudocoherent.
	If $\mathbf{L}i^\ast E\in D^b_{\operatorname{qc}}(\operatorname{Spec}(k))$ for all $i\colon \operatorname{Spec}(k) \to X$ that represents any $p\in |X|$, then $\operatorname{\mathbf{R}\mathcal{H}\! \mathit{om}}(E,A)\in D^b_{\operatorname{qc}}(X)$ for all $A\in D^b_{\operatorname{qc}}(X)$.
\end{lemma}

\begin{proof}
	Consider an \'{e}tale presentation $s\colon U \to X$ from an affine scheme. Observe that for any $q\in U$, with the canonical morphism $i\colon \operatorname{Spec}(\kappa(q))\to U$, the hypothesis says $\mathbf{L}(s\circ i)^\ast E \in D^b_{\operatorname{qc}}(\operatorname{Spec}(\kappa(q)))$.
	Hence, by \cite[Theorem 2.3]{AlonsoTarrio/JeremiasLopez/SanchodeSalas:2023}, $\operatorname{\mathbb{R}\mathcal{H}\! \mathit{om}}(\mathbf{L}s^\ast E, \mathbf{L}s^\ast A) \in D^b_{\operatorname{qc}}(U)$ for any $A\in D^b_{\operatorname{qc}}(X)$ (e.g.\ apply \cite[\href{https://stacks.math.columbia.edu/tag/071Q}{Tag 071Q}, \href{https://stacks.math.columbia.edu/tag/08GH}{08GH}, and \href{https://stacks.math.columbia.edu/tag/08JP}{08JP}]{StacksProject}).
	%%NOTE: See \cite[\S 1.2]{AlonsoTarrio/JeremiasLopez/SanchodeSalas:2023}, which shows internal hom is taken to be in $D(X)$ for their notation.
	It follows that $\operatorname{\mathbb{R}\mathcal{H}\! \mathit{om}} (E,A) \in D^b_{\operatorname{qc}}(X)$.
	Thus, \Cref{lem:internal_hom} shows that $\operatorname{\mathbf{R}\mathcal{H}\! \mathit{om}} (E,A)\cong \operatorname{\mathbb{R}\mathcal{H}\! \mathit{om}} (E,A)$, which completes the proof.
\end{proof}

\begin{lemma}
	\label{lem:perf_pullback_characterize}
	Let $X$ be a quasi-compact quasi-separated algebraic space.
	Then $E\in D(X)$ is perfect if, and only if, there exists an \'{e}tale presentation $s\colon U \to X$ such that $\mathbf{L}s^\ast E$ is perfect in $D(U)$.
\end{lemma}

\begin{proof}
	By \cite[\href{https://stacks.math.columbia.edu/tag/08HG}{Tag 08HG}]{StacksProject}, perfectness of complexes can be detected on the small Zariski site or small \'{e}tale site of a scheme.
	Let $E\in D(X)$ be perfect.
	Choose any \'{e}tale presentation $s\colon U \to X$. By \cite[\href{https://stacks.math.columbia.edu/tag/08H6}{Tag 08H6}]{StacksProject}, $\mathbf{L}s^\ast E$ is perfect.
	We check the converse.

	Assume there exists an \'{e}tale presentation $s\colon U \to X$ such that $\mathbf{L}s^\ast E$ is perfect in $D(U)$.
	Choose an \'{e}tale morphism $t\colon V \to X$.
	Consider the projections $p\colon V\times_X U \to V$ and $q\colon V\times_X U \to U$ from the fiber product.
	Since $\mathbf{L}s^\ast E$ is perfect, it follows that $\mathbf{L}(s\circ q)^\ast E$ is perfect. Hence, $\mathbf{L}(t\circ p)^\ast E$ is perfect.
	Let $r\colon W \to V\times_X U$ be an \'{e}tale presentation.
	Then $\mathbf{L}(t\circ p \circ r)^\ast E$ is perfect.
	However, $p\circ r$ is an \'{e}tale covering of $V$ in $X_{\textrm{\'{e}tale}}$, and so $\mathbf{L}t^\ast E$ is perfect \cite[\href{https://stacks.math.columbia.edu/tag/04LX}{Tag 04LX} and \href{https://stacks.math.columbia.edu/tag/08G5}{08G5}]{StacksProject}.
\end{proof}

\begin{lemma}
	\label{lem:61}
	Let $X$ be a quasi-compact quasi-separated algebraic space.
	Consider $E\in D^b_{\operatorname{qc}}(X)$ pseudocoherent.
	If $\operatorname{\mathbf{R}\mathcal{H}\! \mathit{om}}(E, \mathbf{R} i_\ast \mathcal{O}_{\operatorname{Spec}(k)})\in D^b_{\operatorname{qc}}(X)$ for any $p\in |X|$ and for $i\colon \operatorname{Spec}(k) \to X$ that represents $p$, then $E$ is perfect.
\end{lemma}

\begin{proof}
	Consider an \'{e}tale presentation $s\colon U \to X$ from an affine scheme.
	Fix $q\in U$.
	Let $i\colon \operatorname{Spec}(k)\to U$ be the canonical morphism where $k:=\kappa(q)$.
	There exists a commutative diagram obtained by pullbacks and canonical morphisms,
	\begin{displaymath}
		% https://q.uiver.app/#q=WzAsNyxbMywxLCJVIl0sWzMsMiwiXFxtYXRoY2Fse1h9LiJdLFsyLDIsIlUiXSxbMiwxLCJVXFx0aW1lc197XFxtYXRoY2Fse1h9fSBVIl0sWzEsMiwiXFxvcGVyYXRvcm5hbWV7U3BlY30oaykiXSxbMSwxLCJcXG9wZXJhdG9ybmFtZXtTcGVjfShrKVxcdGltZXNfe1xcbWF0aGNhbHtYfX0gVSJdLFswLDAsIlxcb3BlcmF0b3JuYW1le1NwZWN9KGspIl0sWzAsMSwicyJdLFsyLDEsInMiLDJdLFszLDIsInFfMSIsMl0sWzMsMCwicV8yIiwyXSxbNCwyLCJpIiwyXSxbNSwzLCJpXlxccHJpbWUiXSxbNSw0LCJzXlxccHJpbWUiLDJdLFs2LDQsIjFfe1xcb3BlcmF0b3JuYW1le1NwZWN9KGspfSIsMix7ImN1cnZlIjozfV0sWzYsMCwiaSIsMCx7ImN1cnZlIjotM31dLFs2LDUsImgiXV0=
		\begin{tikzcd}
			{\operatorname{Spec}(k)} &&& \\
			& {\operatorname{Spec}(k)\times_{X} U} & {U\times_{X} U} & U \\
			& {\operatorname{Spec}(k)} & U & {X.}
			\arrow["h", from=1-1, to=2-2]
			\arrow["i", bend left=30, from=1-1, to=2-4]
			\arrow["{1_{\operatorname{Spec}(k)}}"', bend right=30, from=1-1, to=3-2]
			\arrow["{i^\prime}", from=2-2, to=2-3]
			\arrow["{s^\prime}"', from=2-2, to=3-2]
			\arrow["{q_2}"', from=2-3, to=2-4]
			\arrow["{q_1}"', from=2-3, to=3-3]
			\arrow["s", from=2-4, to=3-4]
			\arrow["i"', from=3-2, to=3-3]
			\arrow["s"', from=3-3, to=3-4]
		\end{tikzcd}
	\end{displaymath}
	Since $X$ is quasi-separated, $s\circ i$ is affine \cite[\href{https://stacks.math.columbia.edu/tag/09TF}{Tag 09TF}]{StacksProject}.
	Hence, $q_2 \circ i^\prime$ is affine by base change, and thus, $\operatorname{Spec}(k)\times_{X} U$ is an affine scheme.
	%%NOTE: $U$ is affine
	This implies that $s^\prime$ is an affine,
	%%NOTE: It is a morphism between affine schemes
	and hence, separated morphism.
	As $1_{\operatorname{Spec}(k)}$ is proper and $s^\prime$ is separated, it follows that $h$ is proper.
	%%NOTE: Use Tag 04NX
	Consequently, $\mathbf{R} h_\ast \mathcal{O}_{\operatorname{Spec}(k)}\in D^b_{\operatorname{coh}}(\operatorname{Spec}(k)\times_{X} U)$.
	Moreover, $s^\prime$ is smooth because it is the base change of $s$ along $s\circ i$. Since $\operatorname{Spec}(k)\times_{X} U$ is affine and regular, we have that $\mathbf{R} h_\ast \mathcal{O}_{\operatorname{Spec}(k)} \in \langle \mathcal{O}_{\operatorname{Spec}(k)\times_{X} U} \rangle$.
	Then affineness of $q_2 \circ i^\prime$ yields
	\begin{displaymath}
		\mathbf{R}(q_2 \circ i^\prime)_\ast \mathcal{O}_{\operatorname{Spec}(k)\times_{X} U} \cong (q_2 \circ i^\prime)_\ast \mathcal{O}_{\operatorname{Spec}(k)\times_{X} U}\in D^b_{\operatorname{qc}}(U).
	\end{displaymath}
	By flat base change,
	\begin{displaymath}
		\begin{aligned}
			\mathbf{R}(q_2 \circ i^\prime)_\ast \mathcal{O}_{\operatorname{Spec}(k)\times_{X} U}
			&\cong \mathbf{R}(q_2 \circ i^\prime)_\ast \mathbf{L}(s^\prime)^\ast \mathcal{O}_{\operatorname{Spec}(k)}
			\\&\cong \mathbf{L}s^\ast \mathbf{R}(s\circ i)_\ast \mathcal{O}_{\operatorname{Spec}(k)}.
		\end{aligned}
	\end{displaymath}

	Since $E$ is pseudocoherent and $\mathbf{R}(s\circ i)_\ast \mathcal{O}_{\operatorname{Spec}(k)}\in D^b_{\operatorname{qc}}(X)$ (hence, is bounded below), \Cref{lem:gortz_wedhorn_internal_hom} yields an isomorphism (e.g.\ use that $s$ is \'{e}tale),
	\begin{displaymath}
		\mathbf{L}s^\ast \operatorname{\mathbf{R}\mathcal{H}\! \mathit{om}}(E,  \mathbf{R}(s\circ i)_\ast \mathcal{O}_{\operatorname{Spec}(k)})
		\to \operatorname{\mathbf{R}\mathcal{H}\! \mathit{om}}( \mathbf{L}s^\ast E , \mathbf{L}s^\ast \mathbf{R}(s\circ i)_\ast \mathcal{O}_{\operatorname{Spec}(k)}).
	\end{displaymath}
	It follows that
	\begin{displaymath}
		\begin{aligned}
			\operatorname{\mathbf{R}\mathcal{H}\! \mathit{om}}( \mathbf{L}s^\ast E , \mathbf{R}i_\ast \mathcal{O}_{\operatorname{Spec}(k)} )
			&\cong \operatorname{\mathbf{R}\mathcal{H}\! \mathit{om}}( \mathbf{L}s^\ast E , \mathbf{R}(q_2 \circ i^\prime)_\ast \mathbf{R} h_\ast \mathcal{O}_{\operatorname{Spec}(k)})
			\\&\in \langle \operatorname{\mathbf{R}\mathcal{H}\! \mathit{om}}( \mathbf{L}s^\ast E , \mathbf{R}(q_2 \circ i^\prime)_\ast \mathcal{O}_{\operatorname{Spec}(k)\times_{X} U} ) \rangle
			\\&\in \langle \operatorname{\mathbf{R}\mathcal{H}\! \mathit{om}}( \mathbf{L}s^\ast E , \mathbf{L}s^\ast \mathbf{R}(s\circ i)_\ast \mathcal{O}_{\operatorname{Spec}(k)} ) \rangle
			\\&\in \langle \mathbf{L}s^\ast \operatorname{\mathbf{R}\mathcal{H}\! \mathit{om}}(E,  \mathbf{R}(s\circ i)_\ast \mathcal{O}_{\operatorname{Spec}(k)}) \rangle.
		\end{aligned}
	\end{displaymath}
	By hypothesis,
	\begin{displaymath}
		\operatorname{\mathbf{R}\mathcal{H}\! \mathit{om}}(E,  \mathbf{R}(s\circ i)_\ast \mathcal{O}_{\operatorname{Spec}(k)})\in D^b_{\operatorname{qc}}(X),
	\end{displaymath}
	and hence,
	\begin{displaymath}
		\mathbf{L}s^\ast \operatorname{\mathbf{R}\mathcal{H}\! \mathit{om}}(E,  \mathbf{R}(s\circ i)_\ast \mathcal{O}_{\operatorname{Spec}(k)})  \in D^b_{\operatorname{qc}}(U).
	\end{displaymath}
	However, $D^b_{\operatorname{qc}}(U)$ is thick in $D_{\operatorname{qc}}(U)$, and so,
	\begin{displaymath}
		\operatorname{\mathbf{R}\mathcal{H}\! \mathit{om}}( \mathbf{L}s^\ast E , \mathbf{R}i_\ast \mathcal{O}_{\operatorname{Spec}(k)} )\in D^b_{\operatorname{qc}}(U).
	\end{displaymath}
	As $q\in U$ was arbitrary, \cite[Theorem 2.3]{AlonsoTarrio/JeremiasLopez/SanchodeSalas:2023} shows that $\mathbf{L}s^\ast E$ is perfect (e.g.\ apply \cite[\href{https://stacks.math.columbia.edu/tag/071Q}{Tag 071Q}, \href{https://stacks.math.columbia.edu/tag/08GH}{08GH}, \href{https://stacks.math.columbia.edu/tag/08HG}{08HG}, and \href{https://stacks.math.columbia.edu/tag/08JP}{08JP}]{StacksProject}).
	Thus, via \Cref{lem:perf_pullback_characterize}, $E$ is perfect.
\end{proof}

\begin{lemma}
	\label{lem:pseudocoherence_affine_faithfully_flat_cover}
	Let $f\colon Y \to X$ be a morphism of Noetherian algebraic spaces and $E\in D_{\operatorname{qc}}(X)$.
	If $E$ is pseudocoherent (resp.\ $\operatorname{Perf}(X)$), then $\mathbf{L}f^\ast E$ is pseudocoherent (resp.\ $\mathbf{L}f^\ast E\in \operatorname{Perf}(Y)$).
	Additionally, if $f$ is faithfully flat, then the converse holds.
\end{lemma}

\begin{proof}
	The first claim follows from \cite[\href{https://stacks.math.columbia.edu/tag/08H4}{Tag 08H4} and \href{https://stacks.math.columbia.edu/tag/08H6}{08H6}]{StacksProject}.
	Choose any \'{e}tale morphism $s\colon U \to X$ from a scheme.
	Consider the fibered square
	\begin{displaymath}
		% https://q.uiver.app/#q=WzAsNCxbMSwwLCJVIl0sWzEsMSwiWC4iXSxbMCwxLCJZIl0sWzAsMCwiWVxcdGltZXNfWCBVIl0sWzAsMSwicyJdLFsyLDEsImYiLDJdLFszLDIsInNeXFxwcmltZSIsMl0sWzMsMCwiZl5cXHByaW1lIl1d
		\begin{tikzcd}
			{Y\times_X U} & U \\
			Y & {X.}
			\arrow["{f^\prime}", from=1-1, to=1-2]
			\arrow["{s^\prime}"', from=1-1, to=2-1]
			\arrow["s", from=1-2, to=2-2]
			\arrow["f"', from=2-1, to=2-2]
		\end{tikzcd}
	\end{displaymath}
	Let $t\colon V \to Y\times_X U$ be an \'{e}tale presentation.
	By the first claim, $\mathbf{L}(f\circ s^\prime \circ t)^\ast E$ is pseudocoherent (resp.\ perfect).
	Hence, $\mathbf{L}(f^\prime \circ t)^\ast \mathbf{L}s^\ast E$ is pseudocoherent (resp.\ perfect).
	By \cite[Proposition 22.52]{Gortz/Wedhorn:2023}, $\mathbf{L}s^\ast E$ is pseudocoherent (resp.\ perfect).
	Then \cite[\href{https://stacks.math.columbia.edu/tag/08FU}{Tag 08FU}]{StacksProject} implies $E$ is pseudocoherent (resp.\ perfect).
\end{proof}

\begin{lemma}
	\label{lem:pullback_is_perfect_iff_for_some_or_all_in_equivalence_class}
	Let $X$ be a quasi-compact quasi-separated algebraic space.
	Fix $p\in |X|$.
	For any pseudocoherent complex $E$ on $X$, the following are equivalent:
	\begin{enumerate}
		\item \label{lem:pullback_is_perfect_iff_for_some_or_all_in_equivalence_class1} $\mathbf{L}i^\ast E\in D^b_{\operatorname{coh}}(k)$ for some $i\colon \operatorname{Spec}(k) \to X$ that represents $p$
		\item \label{lem:pullback_is_perfect_iff_for_some_or_all_in_equivalence_class2} $\mathbf{L}t^\ast E\in D^b_{\operatorname{coh}}(L)$ for any $t\colon \operatorname{Spec}(L) \to X$ that represents $p$.
	\end{enumerate}
\end{lemma}

\begin{proof}
	It is straightforward to check that $\eqref{lem:pullback_is_perfect_iff_for_some_or_all_in_equivalence_class2} \implies \eqref{lem:pullback_is_perfect_iff_for_some_or_all_in_equivalence_class1}$.
	We prove the converse.
	Assume $\mathbf{L}i^\ast E\in D^b_{\operatorname{coh}}(k)$ for some $i\colon \operatorname{Spec}(k) \to X$ that represents some $p\in |X|$.
	Let $t\colon \operatorname{Spec}(L) \to X$ be any representative of $p\in |X|$.
	Then there exist a field $\ell$ and a commutative square
	\begin{displaymath}
		% https://q.uiver.app/#q=WzAsNCxbMSwwLCJcXG9wZXJhdG9ybmFtZXtTcGVjfShMKSJdLFsxLDEsIlxcbWF0aGNhbHtYfSJdLFswLDEsIlxcb3BlcmF0b3JuYW1le1NwZWN9KGspIl0sWzAsMCwiXFxvcGVyYXRvcm5hbWV7U3BlY30oXFxlbGwpIl0sWzAsMSwidCJdLFsyLDEsImkiLDJdLFszLDAsImleXFxwcmltZSJdLFszLDIsInReXFxwcmltZSIsMl1d
		\begin{tikzcd}
			{\operatorname{Spec}(\ell)} & {\operatorname{Spec}(L)} \\
			{\operatorname{Spec}(k)} & {X}
			\arrow["{i^\prime}", from=1-1, to=1-2]
			\arrow["{t^\prime}"', from=1-1, to=2-1]
			\arrow["t", from=1-2, to=2-2]
			\arrow["i"', from=2-1, to=2-2]
		\end{tikzcd}
	\end{displaymath}
	because $i$ and $t$ are in the same equivalence class.
	If $\mathbf{L}i^\ast E\in D^b_{\operatorname{coh}}(k)$, then
	\begin{displaymath}
		\mathbf{L}(t \circ i^\prime)^\ast E = \mathbf{L}(i\circ t^\prime)^\ast E\in D^b_{\operatorname{coh}} (\ell).
	\end{displaymath}
	Since $i^\prime$ is faithfully flat, \Cref{lem:zero_object_via_covers,lem:pseudocoherence_affine_faithfully_flat_cover} imply $\mathbf{L}t^\ast E\in D^b_{\operatorname{coh}} (L)$.
\end{proof}

\begin{lemma}
	\label{lem:pullback_is_perfect_iff_closed_pt}
	Let $X$ be a quasi-compact quasi-separated algebraic space.
	Suppose $p\in |X|$ is a closed point.
	For any $E\in D^b_{\operatorname{qc}}(X)$ pseudocoherent, the following are equivalent:
	\begin{enumerate}
		\item \label{lem:pullback_is_perfect_iff_closed_pt1} $\mathbf{L}i^\ast E\in D^b_{\operatorname{coh}}(k)$ for some $i\colon \operatorname{Spec}(k) \to X$ that represents $p$
		\item \label{lem:pullback_is_perfect_iff_closed_pt2} $\mathbf{L}t^\ast E\in D^b_{\operatorname{coh}}(L)$ for any $t\colon \operatorname{Spec}(L) \to X$ that represents a $q\in |X|$ which is a generalization of $p$ (i.e.\ $p\in \overline{\{q\}}$).
	\end{enumerate}
\end{lemma}

\begin{proof}
	It is clear that $\eqref{lem:pullback_is_perfect_iff_closed_pt2}\implies \eqref{lem:pullback_is_perfect_iff_closed_pt1}$.
	We prove the converse.
	Let $s\colon U \to X$ be an \'{e}tale presentation.
	Choose any $q\in |X|$ which is a generalization of $p$.
	By \cite[\href{https://stacks.math.columbia.edu/tag/03JX}{Tag 03JX} and \href{https://stacks.math.columbia.edu/tag/03K2}{03K2}]{StacksProject}, there exist $u,v\in U$ such that $s(u)=p$, $s(v)=q$, and $v$ is a generalization of $u$.
	Consider the canonical morphisms
	\begin{displaymath}
		\operatorname{Spec}(\kappa(u)) \xrightarrow{i_u} \operatorname{Spec}(\mathcal{O}_{U,u}) \xrightarrow{s_u} U.
	\end{displaymath}
	If $\mathbf{L}i^\ast E\in D^b_{\operatorname{coh}}(k)$, then \Cref{lem:pullback_is_perfect_iff_for_some_or_all_in_equivalence_class} implies that $\mathbf{L}(s \circ s_u \circ i_u)^\ast E\in D^b_{\operatorname{coh}}(\kappa(u))$.
	It follows that $\mathbf{L}(s \circ s_u)^\ast E \in \operatorname{Perf}(\mathcal{O}_{U,u})$ (see e.g.\ \cite[\href{https://stacks.math.columbia.edu/tag/068W}{Tag 068W} and \href{https://stacks.math.columbia.edu/tag/08HG}{08HG}]{StacksProject}).
	%%NOTE: Use that $u$ is only closed point of spec of the local ring
	Denote by $\sigma\colon \operatorname{Spec}(\mathcal{O}_{U,v}) \to \operatorname{Spec}(\mathcal{O}_{U,u})$ and $i_v\colon \operatorname{Spec}(\kappa(v))\to \operatorname{Spec}(\mathcal{O}_{U,v})$ the canonical morphisms.
	Then $\mathbf{L}(s \circ s_u \circ \sigma)^\ast E \in \operatorname{Perf}(\mathcal{O}_{U,v})$ and $\mathbf{L}(s \circ s_u \circ \sigma \circ i_v)^\ast E \in D^b_{\operatorname{coh}}(\kappa(v))$.
	Note that $s_u \circ \sigma \circ i_v$ is the canonical morphism $\operatorname{Spec}(\kappa(v))\to U$.
	Since $(s \circ s_u \circ \sigma \circ i_v)(v)=q$, \Cref{lem:pullback_is_perfect_iff_for_some_or_all_in_equivalence_class} tells us $\mathbf{L}t^\ast E\in D^b_{\operatorname{coh}}(L)$ for any $t\colon \operatorname{Spec}(L) \to X$ that represents $q\in |X|$.
\end{proof}

\begin{definition}
	Let $X$ be an algebraic space.
	An object $E\in D_{\operatorname{qc}}(X)$ is said to have \textbf{finite flat dimension} if $E\otimes^{\mathbf{L}} B\in D^b_{\operatorname{qc}}(X)$ for all $B\in D^b_{\operatorname{qc}}(X)$.
\end{definition}

\begin{remark}
	\label{rmk:topological_fact_for_stacks}
	By \cite[\href{https://stacks.math.columbia.edu/tag/03I7}{Tag 03I7}]{StacksProject}, a quasi-separated algebraic space is decent.
	Moreover, \cite[\href{https://stacks.math.columbia.edu/tag/03K3}{Tag 03K3}]{StacksProject} shows that a decent algebraic space has a Kolmogorov underlying topological space.
	Also, from \cite[\href{https://stacks.math.columbia.edu/tag/005E}{Tag 005E}]{StacksProject}, any nonempty quasi-compact Kolmogorov topological space contains a closed point.
	Hence, any $p\in |X|$ is the generalization of a closed point $p^\prime \in |X|$.
\end{remark}

\begin{proposition}
	\label{prop:perfectness}
	Let $X$ be a quasi-compact quasi-separated algebraic space.
	For any $E\in D^b_{\operatorname{qc}}(X)$ pseudocoherent, the following are equivalent:
	\begin{enumerate}
		\item \label{prop:perfectness1} $E$ is perfect
		\item \label{prop:perfectness2} $E$ has finite flat dimension
		\item \label{prop:perfectness3} $E\otimes^{\mathbf{L}} \mathbf{R} i_\ast \mathcal{O}_{\operatorname{Spec}(k)} \in D^b_{\operatorname{qc}}(X)$ for all $i\colon \operatorname{Spec}(k) \to X$ that represents any $p\in |X|$
		\item \label{prop:perfectness4} $\mathbf{L}i^\ast E\in D^b_{\operatorname{qc}}(\operatorname{Spec}(k))$ for all $i\colon \operatorname{Spec}(k) \to X$ that represents any $p\in |X|$
		\item \label{prop:perfectness4_some} for any $p\in |X|$ there exists a representative  $i\colon \operatorname{Spec}(k) \to X$ such that $\mathbf{L}i^\ast E\in D^b_{\operatorname{qc}}(\operatorname{Spec}(k))$
		\item \label{prop:perfectness4_any_closed} $\mathbf{L}i^\ast E\in D^b_{\operatorname{qc}}(\operatorname{Spec}(k))$ for any $i\colon \operatorname{Spec}(k) \to X$ that represents any closed $p\in |X|$
		\item \label{prop:perfectness4_some_closed} for any closed $p\in |X|$  there exists a representative $i\colon \operatorname{Spec}(k) \to X$ such that $\mathbf{L}i^\ast E\in D^b_{\operatorname{qc}}(\operatorname{Spec}(k))$
		\item \label{prop:perfectness5} $\operatorname{\mathbf{R}\mathcal{H}\! \mathit{om}}(E,A)\in D^b_{\operatorname{qc}}(X)$ for all $A\in D^b_{\operatorname{qc}}(X)$
		\item \label{prop:perfectness6} $\operatorname{\mathbf{R}\mathcal{H}\! \mathit{om}}(E, \mathbf{R} i_\ast \mathcal{O}_{\operatorname{Spec}(k)})\in D^b_{\operatorname{qc}}(X)$ for all $i\colon \operatorname{Spec}(k) \to X$ that represents any $p\in |X|$.
	\end{enumerate}
\end{proposition}

\begin{proof}
	If $E$ is perfect, then it has finite flat dimension.
	Indeed, let $s\colon U \to X$ be an \'{e}tale presentation from an affine scheme.
	For any $A\in D^b_{\operatorname{qc}}(X)$, we know that $\mathbf{L}s^\ast E \otimes^{\mathbf{L}} \mathbf{L}s^\ast A\in D^b_{\operatorname{qc}}(U)$ because $\mathbf{L}s^\ast E$ is perfect (e.g.\ use \Cref{lem:perf_pullback_characterize} and \cite[Theorem 2.3]{AlonsoTarrio/JeremiasLopez/SanchodeSalas:2023}).
	Hence, $\mathbf{L}s^\ast E \otimes^{\mathbf{L}} \mathbf{L}s^\ast A \cong \mathbf{L}s^\ast (E \otimes^{\mathbf{L}}A)$ implies that $E \otimes^{\mathbf{L}}A$ has bounded cohomology, and so $E \otimes^{\mathbf{L}}A \in D^b_{\operatorname{qc}}(X)$.
	Thus, $\eqref{prop:perfectness1} \implies \eqref{prop:perfectness2}$.

	By \Cref{lem:23}, we know that $\mathbf{R} i_\ast \mathcal{O}_{\operatorname{Spec}(k)}\in D^b_{\operatorname{qc}}(X)$ because $X$ is quasi-separated.
	Hence, $\eqref{prop:perfectness2} \implies \eqref{prop:perfectness3}$.
	Moreover, \Cref{lem:34} tells us $\eqref{prop:perfectness3} \implies \eqref{prop:perfectness4}$ since $X$ is quasi-separated.
	As for $\eqref{prop:perfectness4} \implies \eqref{prop:perfectness5}$, this is \Cref{lem:45}.
	Clearly, $\eqref{prop:perfectness5} \implies \eqref{prop:perfectness6}$, whereas $\eqref{prop:perfectness6} \implies \eqref{prop:perfectness1}$ is \Cref{lem:61}.

	Lastly, \Cref{lem:pullback_is_perfect_iff_closed_pt,lem:pullback_is_perfect_iff_for_some_or_all_in_equivalence_class} imply that $\eqref{prop:perfectness4} \iff \eqref{prop:perfectness4_some} \iff \eqref{prop:perfectness4_any_closed} \iff \eqref{prop:perfectness4_some_closed}$; see \cite[\href{https://stacks.math.columbia.edu/tag/08H4}{Tag 08H4}]{StacksProject} and \Cref{rmk:topological_fact_for_stacks}.
\end{proof}

%%%%%%%%%%%%%%%%%%%%%%%%%%%%%%%%%%%
\section{Relatively perfect complexes}
\label{sec:relative_perfect}
%%%%%%%%%%%%%%%%%%%%%%%%%%%%%%%%%%%

\begin{definition}
	\label{def:relative_perf_stacks}
	Let $f\colon Y\to X$ be a morphism of quasi-compact quasi-separated algebraic spaces.
	We say that $E\in D_{\operatorname{qc}}(Y)$ is \textbf{relatively perfect over $X$}, or \textbf{$f$-perfect}, if $E\otimes^{\mathbf{L}} \mathbf{L}f^\ast A \in  D^b_{\operatorname{qc}}(Y)$ for all $A\in D^b_{\operatorname{qc}}(X)$.
\end{definition}

\begin{remark}
	\label{rmk:relative_perfects}
	\hfill
	\begin{enumerate}
		\item \label{rmk:relative_perfects1} Any object $E\in D_{\operatorname{qc}}(Y)$ which is $f$-perfect belongs to $D^b_{\operatorname{qc}}(Y)$.
		\item \label{rmk:relative_perfects2} Whereas \cite[3.1]{AlonsoTarrio/JeremiasLopez/SanchodeSalas:2023} requires $E$ to be pseudocoherent, we don't do this in \Cref{def:relative_perf_stacks}.
		Indeed, the results of this section hold without this extra hypothesis.
		Later on we will only consider relatively perfect complexes that are also pseudocoherent.
		\item \label{rmk:relative_perfects3} Let $f\colon Y \to X$ be a morphism of quasi-compact quasi-separated algebraic spaces.
		We claim that for any $E$ is $f$-perfect if, and only if, there exists $[a,b]\subseteq \mathbb{Z}$ such that $\mathcal{H}^j (E\otimes^{\mathbf{L}} \mathbf{L} f^\ast M)\cong 0$ for all $j\not\in [a,b]$ and $M\in \operatorname{Qcoh}(X)$.
		This is akin to \cite[Lemma 3.1]{GuisadoVillalgordo/Lank/ManaliRahul/Pavic:2025}. 
		First, suppose that $E\otimes^{\mathbf{L}} \mathbf{L}f^\ast  D^b_{\operatorname{qc}}(X) \subseteq D^b_{\operatorname{qc}}(Y)$.
		Assume the contrary.
		Then for each $n \geq 1$ there exists an $M_n\in \operatorname{Qcoh}(X)$ such that $\mathcal{H}^j (E\otimes^{\mathbf{L}} \mathbf{L} f^\ast M_n)\not\cong 0$ for some $j\not\in [-n,n]$.
		Set $M:= \oplus_{n\geq 1} M_n$.
		By hypothesis, $E\otimes^{\mathbf{L}} \mathbf{L}f^\ast  M \in D^b_{\operatorname{qc}}(Y)$, and hence there exists $[a,b]\subseteq \mathbb{Z}$ such that $\mathcal{H}^j (E\otimes^{\mathbf{L}} \mathbf{L} f^\ast M)\cong 0$ for all $j\not\in [a,b]$.
		Since $\mathcal{H}^j (E\otimes^{\mathbf{L}} \mathbf{L} f^\ast M_n)$ is a direct summand of $\mathcal{H}^j (E\otimes^{\mathbf{L}} \mathbf{L} f^\ast M)$ for each $n\geq 1$, $\mathcal{H}^j (E\otimes^{\mathbf{L}} \mathbf{L} f^\ast M_n)\cong 0$ for all $j\not\in [a,b]$ and $n\geq 1$.
		This leads to a contradiction.
		Indeed, for $n> \max\{|a|,|b|\}$ there exists $j\not\in [-n,n]$ with
		$\mathcal{H}^j(E\otimes^{\mathbf{L}} \mathbf{L}f^\ast M_n)\not\cong 0$.
		Next, suppose there exists $[a,b]\subseteq \mathbb{Z}$ such that $\mathcal{H}^j (E\otimes^{\mathbf{L}} \mathbf{L} f^\ast M)\cong 0$ for all $j\not\in [a,b]$ and $M\in \operatorname{Qcoh}(X)$.
		It is straightforward to check that $B$ is finitely built by its cohomology sheaves $\mathcal{H}^i (B)$ (see e.g.\ \cite[Lemma 2.5]{GuisadoVillalgordo/Lank/ManaliRahul/Pavic:2025}).
		Moreover, the hypothesis implies $E\otimes^{\mathbf{L}} \mathbf{L}f^\ast \mathcal{H}^i (B) \in D^b_{\operatorname{qc}}(Y)$ for all $i\in \mathbb{Z}$.
		%%NOTE: $D^b_{\operatorname{qc}}(Y)$ is thick in $D_{qc}$
		Since $B\in \langle \{ \mathcal{H}^i (B) \}_{i\in\mathbb{Z}} \rangle$ and $D^b_{\operatorname{qc}}(Y)$ is a thick subcategory of $D_{\operatorname{qc}}(Y)$, it follows that
		\begin{displaymath}
			E\otimes^{\mathbf{L}} \mathbf{L} f^\ast B\in \langle \{ E\otimes^{\mathbf{L}} \mathbf{L} f^\ast \mathcal{H}^i (B) \}_{i\in \mathbb{Z}} \rangle\subseteq D^b_{\operatorname{qc}}(Y).
		\end{displaymath}
		This completes the proof.
		\item \label{rmk:relative_perfects4} Recall that a morphism $f\colon Y\to X$ of quasi-compact quasi-separated algebraic spaces has \textbf{finite tor-dimension} if there exists $[a,b]\subseteq \mathbb{Z}$ such that $\mathcal{H}^j (E\otimes^{\mathbf{L}} \mathbf{L} f^\ast M)\cong 0$ for all $j\not\in [a,b]$ and $M\in \operatorname{Qcoh}(X)$.
		Hence, by \eqref{rmk:relative_perfects3}, $\mathcal{O}_Y$ is $f$-perfect if, and only if, $f$ has finite tor-dimension. 
		In such a case, it follows from \Cref{prop:perfectness} that every perfect complex on $Y$ is $f$-perfect.
		Coversely, if every perfect complex on $Y$ is $f$-perfect, then so is $\mathcal{O}_Y$, and hence $f$ is of finite tor-dimension.
	\end{enumerate}
\end{remark}

\begin{lemma}
	\label{lem:finite_coh_dim_implies_restriction_to_bounded}
	Let $f\colon Y\to X$ be a morphism of quasi-compact quasi-separated algebraic spaces.
	Then $\mathbf{R}f_\ast D^b_{\operatorname{qc}}(Y)\subseteq D^b_{\operatorname{qc}}(X)$.
\end{lemma}

\begin{proof}
	By the truncation triangles for the bounded standard $t$-structure, every object of $D^b_{\operatorname{qc}}(Y)$ is finitely built by its cohomology sheaves.
	Hence, we can reduce to checking the claim for every quasi-coherent $\mathcal{O}_Y$-module.
	Applying \cite[\href{https://stacks.math.columbia.edu/tag/08FA}{Tag 08FA}]{StacksProject}, we can find some $N\geq 0$ such that $\mathbf{R}^i f_\ast M \cong 0$ for all $M\in \operatorname{Qcoh}(Y)$ and $i>N$.
	This completes the proof.
\end{proof}

\begin{lemma}
	\label{lem:quasi_affine_reflects_boundedness}
	Let $f\colon Y \to X$ be a quasi-affine morphism of quasi-compact quasi-separated algebraic spaces.
	If $E\in D_{\operatorname{qc}}(Y)$ satisfies $\mathbf{R}f_\ast E \in D^b_{\operatorname{qc}}(X)$, then $E\in D^b_{\operatorname{qc}}(Y)$.
\end{lemma}

\begin{proof}
	By \cite[\href{https://stacks.math.columbia.edu/tag/081X}{Tag 081X}]{StacksProject}, $f$ has a canonical factorization
	\begin{displaymath}
		Y \xrightarrow{j} \underline{\operatorname{Spec}}_X (f_\ast \mathcal{O}_Y) \xrightarrow{f^\prime} X.
	\end{displaymath}
	Here $j$ is an open immersion and $f^\prime$ is affine.
	Since $\mathbf{R}f_\ast E \in D^b_{\operatorname{qc}}(X)$ and $f^\prime$ is affine, it follows that $\mathbf{R}j_\ast E\in D^b_{\operatorname{qc}}(\underline{\operatorname{Spec}}_X (f_\ast \mathcal{O}_Y))$.
	Indeed, by \Cref{lem:quasi-affine_is_conservative} and \cite[\href{https://stacks.math.columbia.edu/tag/073H}{Tag 073H}]{StacksProject}, $\mathbf{R}f^\prime_\ast$ is conservative and $t$-exact with respect to the standard $t$-structures.
	However, $j$ is an open immersion, and so the counit $\mathbf{L}j^\ast \mathbf{R}j_\ast E \to E$ is an isomorphism.
	See \Cref{lem:spacey_recollement}.
	As $j$ is flat and $\mathbf{R}j_\ast E\in D^b_{\operatorname{qc}}(\underline{\operatorname{Spec}}_X (f_\ast \mathcal{O}_Y))$, it follows that $E\in D^b_{\operatorname{qc}}(Y)$.
\end{proof}

\begin{proposition}
	\label{prop:relative_perfect_smooth_locality_ascent}
	Let $f\colon Y\to X$ be a morphism of quasi-compact quasi-separated algebraic spaces.
	Suppose that $s\colon U\to X$ is a flat quasi-affine morphism from a quasi-compact quasi-separated algebraic space.
	Consider the tor-independent fibered square
	\begin{displaymath}
		% https://q.uiver.app/#q=WzAsNCxbMSwwLCJcXG1hdGhjYWx7VX0iXSxbMSwxLCJcXG1hdGhjYWx7WH0uIl0sWzAsMCwiXFxtYXRoY2Fse1l9XFx0aW1lc197XFxtYXRoY2Fse1h9fSBcXG1hdGhjYWx7VX0iXSxbMCwxLCJcXG1hdGhjYWx7WX0iXSxbMCwxLCJzIl0sWzIsMCwiZl5cXHByaW1lIl0sWzIsMywic15cXHByaW1lIiwyXSxbMywxLCJmIiwyXV0=
		\begin{tikzcd}
			{Y\times_{X} U} & {U} \\
			{Y} & {X.}
			\arrow["{f^\prime}", from=1-1, to=1-2]
			\arrow["{s^\prime}"', from=1-1, to=2-1]
			\arrow["s", from=1-2, to=2-2]
			\arrow["f"', from=2-1, to=2-2]
		\end{tikzcd}
	\end{displaymath}
	If $E\in D_{\operatorname{qc}}(Y)$ is $f$-perfect, then $\mathbf{L}(s^\prime)^\ast E$ is $f^\prime$-perfect.
\end{proposition}

\begin{proof}
	Note that $E\in D^b_{\operatorname{qc}}(Y)$.
	By \Cref{lem:finite_coh_dim_implies_restriction_to_bounded}, we have $\mathbf{R}s_\ast D^b_{\operatorname{qc}}(U)\subseteq D^b_{\operatorname{qc}}(X)$.
	Since $E$ is $f$-perfect, it follows that
	\begin{displaymath}
		E\otimes^{\mathbf{L}} \mathbf{L}f^\ast \mathbf{R}s_\ast D^b_{\operatorname{qc}}(U)\subseteq D^b_{\operatorname{qc}}(Y).
	\end{displaymath}
	From tor-independent base change \cite[\href{https://stacks.math.columbia.edu/tag/08IR}{Tag 08IR}]{StacksProject}, we obtain that
	\begin{displaymath}
		E\otimes^{\mathbf{L}} \mathbf{R}s^\prime_\ast \mathbf{L}(f^\prime)^\ast D^b_{\operatorname{qc}}(U)\subseteq D^b_{\operatorname{qc}}(Y).
	\end{displaymath}
	By the projection formula, for any $A\in D^b_{\operatorname{qc}}(U)$, there exists an isomorphism
	\begin{displaymath}
		E\otimes^{\mathbf{L}} \mathbf{R}s^\prime_\ast \mathbf{L}(f^\prime)^\ast A \cong \mathbf{R}s^\prime_\ast ( \mathbf{L}(f^\prime)^\ast A \otimes^{\mathbf{L}} \mathbf{L}(s^\prime)^\ast E).
	\end{displaymath}
	By base change, $s^\prime$ is quasi-affine.
	By \Cref{lem:quasi_affine_reflects_boundedness},
	\begin{displaymath}
		\mathbf{R}s^\prime_\ast ( \mathbf{L}(f^\prime)^\ast A \otimes^{\mathbf{L}} \mathbf{L}(s^\prime)^\ast E)\in D^b_{\operatorname{qc}}(Y)
	\end{displaymath}
	implies
	\begin{displaymath}
		\mathbf{L}(f^\prime)^\ast A \otimes^{\mathbf{L}} \mathbf{L}(s^\prime)^\ast E\in D^b_{\operatorname{qc}}(Y\times_{X} U).
	\end{displaymath}
	Therefore, $\mathbf{L}(s^\prime)^\ast E$ is $f^\prime$-perfect.
\end{proof}

\begin{proposition}
	\label{prop:relative_perfect_smooth_locality_descent}
	Let $f\colon Y\to X$ be a morphism of quasi-compact quasi-separated algebraic spaces.
	Suppose that $s\colon U\to X$ is a flat surjective morphism of algebraic spaces.
	Consider the fibered square
	\begin{displaymath}
		% https://q.uiver.app/#q=WzAsNCxbMSwwLCJcXG1hdGhjYWx7VX0iXSxbMSwxLCJcXG1hdGhjYWx7WH0uIl0sWzAsMCwiXFxtYXRoY2Fse1l9XFx0aW1lc197XFxtYXRoY2Fse1h9fSBcXG1hdGhjYWx7VX0iXSxbMCwxLCJcXG1hdGhjYWx7WX0iXSxbMCwxLCJzIl0sWzIsMCwiZl5cXHByaW1lIl0sWzIsMywic15cXHByaW1lIiwyXSxbMywxLCJmIiwyXV0=
		\begin{tikzcd}
			{Y\times_{X} U} & {U} \\
			{Y} & {X.}
			\arrow["{f^\prime}", from=1-1, to=1-2]
			\arrow["{s^\prime}"', from=1-1, to=2-1]
			\arrow["s", from=1-2, to=2-2]
			\arrow["f"', from=2-1, to=2-2]
		\end{tikzcd}
	\end{displaymath}
	If $\mathbf{L}(s^\prime)^\ast E$ is $f^\prime$-perfect where $E\in D_{\operatorname{qc}}(Y)$, then $E$ is $f$-perfect.
\end{proposition}

\begin{proof}
	By base change, $s^\prime$ is flat and surjective.
	As $\mathbf{L}(s^\prime)^\ast E$ is $f^\prime$-perfect, $E\in D^b_{\operatorname{qc}}(Y)$.
	Let $A\in D^b_{\operatorname{qc}}(X)$.
	Since $\mathbf{L}(s^\prime)^\ast E$ is $f^\prime$-perfect, it follows that
	\begin{displaymath}
		\mathbf{L}(f^\prime)^\ast \mathbf{L}s^\ast A \otimes^{\mathbf{L}} \mathbf{L}(s^\prime)^\ast E\in D^b_{\operatorname{qc}}(Y\times_{X} U).
	\end{displaymath}
	This implies
	\begin{displaymath}
		\mathbf{L}(f\circ s^\prime)^\ast A \otimes^{\mathbf{L}} \mathbf{L}(s^\prime)^\ast E\in D^b_{\operatorname{qc}}(Y\times_{X} U).
	\end{displaymath}
	Yet, by \Cref{lem:zero_object_via_covers}, $E\otimes^{\mathbf{L}} \mathbf{L}f^\ast A \in D^b_{\operatorname{qc}}(Y)$. Hence, $E$ is $f$-perfect.
\end{proof}

\begin{corollary}
	\label{cor:smooth_locally_f_perfect_implies_f_perfect}
	Let $f\colon Y\to X$ be a morphism of quasi-compact quasi-separated algebraic spaces.
	Choose $E\in D_{\operatorname{qc}}(Y)$.
	Assume $\mathbf{L}q^\ast E$ is $g$-perfect for every commutative diagram
	\begin{displaymath}
		% https://q.uiver.app/#q=WzAsNCxbMCwwLCJWIl0sWzAsMSwiXFxtYXRoY2Fse1l9Il0sWzEsMCwiVSJdLFsxLDEsIlxcbWF0aGNhbHtYfSJdLFswLDEsInEiLDJdLFsyLDMsInAiXSxbMSwzLCJmIiwyXSxbMCwyLCJnIl1d
		\begin{tikzcd}
			V & U \\
			{Y} & {X}
			\arrow["g", from=1-1, to=1-2]
			\arrow["q"', from=1-1, to=2-1]
			\arrow["p", from=1-2, to=2-2]
			\arrow["f"', from=2-1, to=2-2]
		\end{tikzcd}
	\end{displaymath}
	where $p$ and $q$ are flat surjective morphisms from quasi-compact quasi-separated schemes.
	Then $E$ is $f$-perfect.
\end{corollary}

\begin{proof}
	Note that $E\in D^b_{\operatorname{qc}} (Y)$.
	Let $p\colon U \to X$ be an \'{e}tale presentation.
	Consider the fibered square
	\begin{displaymath}
		% https://q.uiver.app/#q=WzAsNCxbMCwwLCJcXG1hdGhjYWx7WX1cXHRpbWVzX3tcXG1hdGhjYWx7WH19IFUiXSxbMCwxLCJcXG1hdGhjYWx7WX0iXSxbMSwxLCJcXG1hdGhjYWx7WH0uIl0sWzEsMCwiVSJdLFswLDEsInBeXFxwcmltZSIsMl0sWzEsMiwiZiIsMl0sWzMsMiwicCJdLFswLDMsImZeXFxwcmltZSJdXQ==
		\begin{tikzcd}
			{Y\times_{X} U} & U \\
			{Y} & {X.}
			\arrow["{f^\prime}", from=1-1, to=1-2]
			\arrow["{p^\prime}"', from=1-1, to=2-1]
			\arrow["p", from=1-2, to=2-2]
			\arrow["f"', from=2-1, to=2-2]
		\end{tikzcd}
	\end{displaymath}
	Choose an \'{e}tale presentation $t\colon V \to Y\times_{X} U$.
	By assumption, we know that $\mathbf{L}(p^\prime \circ t)^\ast E$ is $(f^\prime \circ t)$-perfect.
	Hence, for every $B\in D^b_{\operatorname{qc}}(X)$, it follows that
	\begin{displaymath}
		\mathbf{L}(p^\prime \circ t)^\ast E \otimes^{\mathbf{L}} \mathbf{L}(f\circ p^\prime \circ t)^\ast B \in D^b_{\operatorname{qc}}(V).
	\end{displaymath}
	Since $t$ is flat and surjective, \Cref{lem:zero_object_via_covers} gives
	\begin{displaymath}
		\mathbf{L}(p^\prime)^\ast E \otimes^{\mathbf{L}} \mathbf{L}(f\circ p^\prime)^\ast B \in D^b_{\operatorname{qc}}(Y\times_{X} U).
	\end{displaymath}
	However, by similar reasoning, we obtain that $E\otimes^{\mathbf{L}} \mathbf{L}f^\ast B\in D^b_{\operatorname{qc}}(Y)$.
\end{proof}

\begin{theorem}
	\label{thm:relative_perf_iff_smooth_locally}
	Let $X$ be a quasi-compact quasi-separated algebraic space.
	Consider a quasi-compact quasi-separated morphism $f\colon Y\to X$ of algebraic spaces.
	Then $E\in D_{\operatorname{qc}}(Y)$ is $f$-perfect if, and only if, $\mathbf{L}q^\ast E$ is $g$-perfect for every commutative diagram
	\begin{displaymath}
		% https://q.uiver.app/#q=WzAsNCxbMCwwLCJWIl0sWzAsMSwiXFxtYXRoY2Fse1l9Il0sWzEsMCwiVSJdLFsxLDEsIlxcbWF0aGNhbHtYfSJdLFswLDEsInEiLDJdLFsyLDMsInAiXSxbMSwzLCJmIiwyXSxbMCwyLCJnIl1d
		\begin{tikzcd}
			V & U \\
			{Y} & {X}
			\arrow["g", from=1-1, to=1-2]
			\arrow["q"', from=1-1, to=2-1]
			\arrow["p", from=1-2, to=2-2]
			\arrow["f"', from=2-1, to=2-2]
		\end{tikzcd}
	\end{displaymath}
	where $p$ and $q$ are \'{e}tale surjective morphisms from quasi-compact quasi-separated schemes.
\end{theorem}

\begin{proof}
	In both cases, $E$ has bounded cohomology, and so we can impose this without loss of generality in the proof.
	First, let $E\in D^b_{\operatorname{qc}}(Y)$ be $f$-perfect.
	Form a commutative diagram
	\begin{displaymath}
		% https://q.uiver.app/#q=WzAsNCxbMCwwLCJWIl0sWzAsMSwiXFxtYXRoY2Fse1l9Il0sWzEsMCwiVSJdLFsxLDEsIlxcbWF0aGNhbHtYfSJdLFswLDEsInEiLDJdLFsyLDMsInAiXSxbMSwzLCJmIiwyXSxbMCwyLCJnIl1d
		\begin{tikzcd}
			V & U \\
			{Y} & {X}
			\arrow["g", from=1-1, to=1-2]
			\arrow["q"', from=1-1, to=2-1]
			\arrow["p", from=1-2, to=2-2]
			\arrow["f"', from=2-1, to=2-2]
		\end{tikzcd}
	\end{displaymath}
	where $p$ and $q$ are \'{e}tale surjective morphisms from quasi-compact quasi-separated schemes.
	If needed, choose $s\colon U^\prime \to U$ an \'{e}tale presentation from an affine scheme.
	We show that $\mathbf{L}q^\ast E$ is $g$-perfect.
	Consider the commutative diagram
	\begin{displaymath}
		% https://q.uiver.app/#q=WzAsNixbMCwxLCJcXG1hdGhjYWx7WX1cXHRpbWVzX3tcXG1hdGhjYWx7WH19IFUiXSxbMCwyLCJcXG1hdGhjYWx7WX0iXSxbMSwyLCJcXG1hdGhjYWx7WH0uIl0sWzEsMSwiVSJdLFsxLDAsIlVeXFxwcmltZSJdLFswLDAsIlxcbWF0aGNhbHtZfVxcdGltZXNfe1xcbWF0aGNhbHtYfX0gVV5cXHByaW1lIl0sWzAsMSwicF5cXHByaW1lIiwyXSxbMSwyLCJmIiwyXSxbMywyLCJwIl0sWzAsMywiZl5cXHByaW1lIl0sWzQsMywicyJdLFs1LDAsInNeXFxwcmltZSIsMl0sWzUsNCwiZl57XFxwcmltZSBcXHByaW1lfSJdXQ==
		\begin{tikzcd}
			{Y\times_{X} U^\prime} & {U^\prime} \\
			{Y\times_{X} U} & U \\
			{Y} & {X.}
			\arrow["{f^{\prime\prime}}", from=1-1, to=1-2]
			\arrow["{s^\prime}"', from=1-1, to=2-1]
			\arrow["s", from=1-2, to=2-2]
			\arrow["{f^\prime}", from=2-1, to=2-2]
			\arrow["{p^\prime}"', from=2-1, to=3-1]
			\arrow["p", from=2-2, to=3-2]
			\arrow["f"', from=3-1, to=3-2]
		\end{tikzcd}
	\end{displaymath}
	By \Cref{lem:quasi_affine_diagonal}, $p\circ s$ is quasi-affine.
	Moreover, as $p\circ s$ is flat,
	%%NOTE: It is a composition of smooth morphisms. Also, it is qcqs with scheme source.
	\Cref{prop:relative_perfect_smooth_locality_ascent} implies that $\mathbf{L}(p^\prime \circ s^\prime)^\ast E$ is $f^{\prime\prime}$-perfect.
	Next, \Cref{prop:relative_perfect_smooth_locality_descent} applied to $f^\prime$ and $s$, shows $\mathbf{L}(p^\prime)^\ast E$ is $f^\prime$-perfect.
	Since $g$ factors through $f^\prime$ by a morphism $h\colon V\to Y\times_{X} U$, it follows that $\mathbf{L}(p^\prime \circ h)^\ast E$ is $f^\prime \circ h$-perfect because $h$ is \'{e}tale \cite[\href{https://stacks.math.columbia.edu/tag/05W3}{Tag 05W3}]{StacksProject}, which shows the desired claim.
	%%NOTE: Use $g = f^\prime \circ h$.
	The converse direction follows from \Cref{cor:smooth_locally_f_perfect_implies_f_perfect}.
\end{proof}

%%%%%%%%%%%%%%%%%%%%%%%%%%%%%%%%%%%
\section{Relatively quasi-perfect complexes}
\label{sec:quasi-perfect}
%%%%%%%%%%%%%%%%%%%%%%%%%%%%%%%%%%%

\begin{definition}
	\label{def:quasi_perf}
	Let $f\colon Y\to X$ be a morphism of quasi-compact quasi-separated algebraic spaces.
	We say that $E\in D_{\operatorname{qc}}(Y)$ is \textbf{relatively quasi-perfect over $X$}, or \textbf{$f$-quasi-perfect}, if $\mathbf{R} f_\ast (E\otimes^{\mathbf{L}} P)\in \operatorname{Perf}(X)$ for all $P\in \operatorname{Perf} (Y)$.
\end{definition}

\begin{remark}
	\label{rmk:quasi-perfect_morphisms}
	Recall that $f$ is \textbf{quasi-perfect} if $\mathbf{R}f_\ast$ sends perfect complexes to perfect complexes \cite{Lipman/Neeman:2007}.
	Hence, $f$ is quasi-perfect if, and only if, $\mathcal{O}_Y$ is $f$-quasi-perfect.
	In other words, \Cref{def:quasi_perf} is a relative version of this notion for objects of $D_{\operatorname{qc}}(Y)$.
	For proper morphisms of Noetherian algebraic spaces, quasi-perfectness can be checked on the \'{e}tale local rings of the target, their completions, and their (strict) Henselizations, and the locus over which a proper morphism is quasi-perfect is open \cite{DeDeyn/Lank/ManaliRahul:2025b}.
\end{remark}

\begin{lemma}
	\label{lem:quasi-perfect_bounded}
	Let $f\colon Y\to X$ be a morphism of quasi-compact quasi-separated algebraic spaces.
	If $E$ is $f$-quasi-perfect, then $E\otimes^{\mathbf{L}} Q$ is $f$-quasi-perfect for all $Q\in \operatorname{Perf}(Y)$.
	Moreover, if $E$ is pseudocoherent, then $\operatorname{\mathbb{R}\mathcal{H}\! \mathit{om}}(E, f^\times \mathcal{O}_X)\in D^b_{\operatorname{qc}}(Y)$.
	Lastly, if $f^\times \mathcal{O}_X$ has coherent cohomology and $Y$ is Noetherian, then $\operatorname{\mathbb{R}\mathcal{H}\! \mathit{om}}(E, f^\times \mathcal{O}_X)$ is pseudocoherent (and hence, is in $D^b_{\operatorname{coh}}(Y)$).
\end{lemma}

\begin{proof}
	If $E$ is $f$-quasi-perfect, then $E\otimes^{\mathbf{L}} Q$ is $f$-quasi-perfect for all $Q\in \operatorname{Perf}(Y)$, e.g.\ $f$-quasi-perfectness is equivalent to
	\begin{displaymath}
		\mathbf{R} f_\ast (E\otimes^{\mathbf{L}} \operatorname{Perf}(Y))\subseteq \operatorname{Perf}(X).
	\end{displaymath}
	We check the second claim.

	By \cite[\href{https://stacks.math.columbia.edu/tag/0E56}{Tag 0E56}]{StacksProject}, we have $f^\times \mathcal{O}_X \in D^+_{\operatorname{qc}}(Y)$.
	Moreover, if $E$ is pseudocoherent, it follows that $\operatorname{\mathbb{R}\mathcal{H}\! \mathit{om}}(E, f^\times \mathcal{O}_X)\in D_{\operatorname{qc}}(Y)$ \cite[\href{https://stacks.math.columbia.edu/tag/0A8A}{Tag 0A8A}]{StacksProject}.
	Furthermore, $\operatorname{\mathbb{R}\mathcal{H}\! \mathit{om}}$ is \'{e}tale local, meaning it commutes with derived pullback along \'{e}tale presentations \cite[\href{https://stacks.math.columbia.edu/tag/04LX}{Tag 04LX} and \href{https://stacks.math.columbia.edu/tag/08JB}{08JB}]{StacksProject}.
	Hence, in the case $E$ is pseudocoherent, it follows that $\operatorname{\mathbb{R}\mathcal{H}\! \mathit{om}}(E, f^\times \mathcal{O}_X)\in D^+_{\operatorname{qc}}(Y)$ \cite[\href{https://stacks.math.columbia.edu/tag/0A6H}{Tag 0A6H}]{StacksProject}.

	Next, we prove that $\operatorname{\mathbb{R}\mathcal{H}\! \mathit{om}}(E, f^\times \mathcal{O}_X)\in D^-_{\operatorname{qc}}(Y)$ if $E$ is pseudocoherent.
	In other words, $\operatorname{\mathbb{R}\mathcal{H}\! \mathit{om}}(E, f^\times \mathcal{O}_X)$ is bounded for any $E$ pseudocoherent and $f$-quasi-perfect.
	To check this claim, let $G$ be a compact generator for $D_{\operatorname{qc}}(Y)$ and $n\in \mathbb{Z}$.
	By adjunction,
	\begin{displaymath}
		\begin{aligned}
			\operatorname{Hom}(G[n],\operatorname{\mathbb{R}\mathcal{H}\! \mathit{om}}(E, f^\times \mathcal{O}_X))
			&\cong \operatorname{Hom}(G[n] \otimes^{\mathbf{L}} E , f^\times \mathcal{O}_X) && (\textrm{tensor/hom})
			\\&\cong \operatorname{Hom}(\mathbf{R}f_\ast (G \otimes^{\mathbf{L}} E)[n], \mathcal{O}_X) && (\textrm{right adj. of push}).
		\end{aligned}
	\end{displaymath}
	Since $E$ is $f$-quasi-perfect, $\mathbf{R}f_\ast (G \otimes^{\mathbf{L}} E)\in \operatorname{Perf}(X)$.
	Then \cite[\href{https://stacks.math.columbia.edu/tag/0GFH}{Tag 0GFH}]{StacksProject} asserts that $\operatorname{Hom}(\mathbf{R}f_\ast (G \otimes^{\mathbf{L}} E)[-N], \mathcal{O}_X)\cong 0$ for $N \gg 0$.
	A further application of \cite[\href{https://stacks.math.columbia.edu/tag/0GFH}{Tag 0GFH}]{StacksProject} shows $\operatorname{\mathbb{R}\mathcal{H}\! \mathit{om}}(E, f^\times \mathcal{O}_X)\in D^b_{\operatorname{qc}}(Y)$.

	Lastly, we show if $f^\times \mathcal{O}_X$ has coherent cohomology and $Y$ is Noetherian, it follows that $\operatorname{\mathbb{R}\mathcal{H}\! \mathit{om}}(E, f^\times \mathcal{O}_X)$ has coherent cohomology.
	This follows from \cite[\href{https://stacks.math.columbia.edu/tag/0D0S}{Tag 0D0S}]{StacksProject}.
\end{proof}

\begin{lemma}
	\label{lem:neeman2023cor5_10}
	Let $f_i\colon Y_i \to S$ be morphisms of quasi-compact quasi-separated algebraic spaces.
	If $P_i$ compactly generates $D_{\operatorname{qc}}(Y_i)$ for each $i$, then $D_{\operatorname{qc}}(Y_1\times_{S} Y_2)$ is compactly generated by $\mathbf{L}g_1^\ast P_1 \otimes^{\mathbf{L}} \mathbf{L}g_2^\ast P_2$ where $g_i\colon Y_1\times_{S} Y_2 \to Y_i$ are the natural projections.
\end{lemma}

\begin{remark}
	\Cref{lem:neeman2023cor5_10} is a version of \cite[Corollary 5.10]{Neeman:2023} for the small \'{e}tale site.
\end{remark}

\begin{proof}
	[Proof of \Cref{lem:neeman2023cor5_10}]
	Since $\mathbf{L}g_1^\ast P_1 \otimes^{\mathbf{L}} \mathbf{L}g_2^\ast P_2$ is perfect, it is compact.
	Let $E\in D_{\operatorname{qc}}(Y_1\times_S Y_2)$ satisfy for all $n\in \mathbb{Z}$,
	\begin{displaymath}
		\operatorname{Hom}(\mathbf{L}g_1^\ast P_1 \otimes^{\mathbf{L}} \mathbf{L}g_2^\ast P_2 [n] , E)\cong 0.
	\end{displaymath}
	By tensor--Hom adjunction, and using that $P_1$ compactly generates $D_{\operatorname{qc}}(Y_1)$, 
	\begin{displaymath}
		\mathbf{R}(g_1)_\ast
		\operatorname{\mathbf{R}\mathcal{H}\! \mathit{om}}
		(\mathbf{L}g_2^\ast P_2,E)\cong 0.
	\end{displaymath}
	Choose an \'{e}tale surjection $s\colon U\to Y_1$ with $U$ affine.
	Set $t\colon U\times_S Y_2\to Y_1\times_S Y_2$ to be its base change. Denote by $h_1$ and $h_2$ respectively the projections from $U\times_S Y_2$ to $U$ and $Y_2$.
	Flat base change and the perfectness of $P_2$ yield
	\begin{displaymath}
		\mathbf{R}(h_1)_\ast
		\operatorname{\mathbf{R}\mathcal{H}\! \mathit{om}}
		(\mathbf{L}h_2^\ast P_2,\mathbf{L}t^\ast E)
		\cong 0.
	\end{displaymath}
	%%NOTE: Flat base change left most internal square for U and Y_2 over S with U on bottom. Use pullback commutes with Rhom if left most input is perfect. 
	As $U$ is affine, adjunction yields $\operatorname{Hom} (P_2,\mathbf{R}(h_2)_\ast\mathbf{L}t^\ast E[n])=0$ for every $n\in\mathbb{Z}$.
	%%NOTE: Hom (O_U, (h_2)_\ast blah) \cong 0 so apply h_2^\ast O_U = O, use tensor hom, then push/pull
	Since $P_2$ compactly generates $D_{\operatorname{qc}}(Y_2)$, it follows that $\mathbf{R}(h_2)_\ast\mathbf{L}t^\ast E\cong 0$.
	By \Cref{lem:quasi_affine_diagonal}, $U\to S$ is quasi-affine, and hence so is $h_2$.
	By \Cref{lem:quasi-affine_is_conservative}, $\mathbf{L}t^\ast E\cong 0$.
	As $t$ is faithfully flat, $E\cong 0$.
\end{proof}

\begin{lemma}
	\label{lem:preservation_of_perfect}
	Let $S$ be a quasi-compact quasi-separated algebraic space.
	Suppose $f_1\colon Y_1 \to S$ and $f_2\colon Y_2 \to S$ are morphisms of algebraic spaces.
	Denote by $p_i\colon Y_1 \times_{S} Y_2 \to Y_i$ the natural projection. For any $E\in D_{\operatorname{qc}}(Y_1 \times_{S} Y_2)$, the following are equivalent:
	\begin{enumerate}
		\item $E$ is $p_2$-quasi-perfect
		\item $\Phi_E (\operatorname{Perf}(Y_1))\subseteq \operatorname{Perf}(Y_2)$.
	\end{enumerate}
\end{lemma}

\begin{proof}
	First, assume $E$ is $p_2$-quasi-perfect.
	Since $\mathbf{L}p_1^\ast \operatorname{Perf}(Y_1)\subseteq \operatorname{Perf}(Y_1 \times_{S} Y_2)$, we have $\Phi_E (\operatorname{Perf}(Y_1))\subseteq \operatorname{Perf}(Y_2)$. We check the converse.

	By \cite[\href{https://stacks.math.columbia.edu/tag/09IY}{Tag 09IY}]{StacksProject}, each $D_{\operatorname{qc}}(Y_i)$ is compactly generated by an object $G_i$.
	By \Cref{lem:neeman2023cor5_10}, $\mathbf{L}p_1^\ast G_1 \otimes^{\mathbf{L}} \mathbf{L}p_2^\ast G_2$ is a compact generator for $D_{\operatorname{qc}}(Y_1 \times_{S} Y_2)$.
	Hence, $\mathbf{L}p_1^\ast G_1 \otimes^{\mathbf{L}} \mathbf{L}p_2^\ast G_2$ is a classical generator for $\operatorname{Perf}(Y_1 \times_{S} Y_2)$ (see \cite[\href{https://stacks.math.columbia.edu/tag/09SR}{Tag 09SR} and \href{https://stacks.math.columbia.edu/tag/09M8}{09M8}]{StacksProject}). Therefore,
	\begin{displaymath}
		\Phi_E (G_1)\otimes^{\mathbf{L}} G_2 \cong \mathbf{R}(p_2)_\ast (E \otimes^{\mathbf{L}} \mathbf{L}p_1^\ast G_1 \otimes^{\mathbf{L}} \mathbf{L}p_2^\ast G_2)\in \operatorname{Perf}(Y_2).
	\end{displaymath}
	It follows that $E$ is $p_2$-quasi-perfect.
\end{proof}

\begin{lemma}
	\label{lem:f_perf_implies_perf_preserved}
	Let $f\colon Y\to X$ be a proper morphism of Noetherian algebraic spaces. Let $E\in D_{\operatorname{qc}}(Y)$ be pseudocoherent.
	If $E$ is $f$-perfect, then $E$ is $f$-quasi-perfect.
\end{lemma}

\begin{proof}
	Throughout, we freely appeal to \Cref{lem:finite_coh_dim_implies_restriction_to_bounded}.
	Since $E$ is $f$-perfect and pseudocoherent, we know that $E \in D^b_{\operatorname{coh}}(Y)$.
	Let $P\in \operatorname{Perf}(Y)$. Then $E \otimes^{\mathbf{L}} P \in D^b_{\operatorname{coh}}(Y)$ because $(-)\otimes^{\mathbf{L}} P$ is an endofunctor on $D^b_{\operatorname{coh}}(Y)$.
	Since $f$ is proper, it follows that $\mathbf{R}f_\ast (E \otimes^{\mathbf{L}} P) \in D^b_{\operatorname{coh}}(X)$.
	Now let $i\colon \operatorname{Spec}(k) \to X$ be a representative of some $p\in |X|$.
	Then $\mathbf{R} i_\ast \mathcal{O}_{\operatorname{Spec}(k)} \in D^b_{\operatorname{qc}}(X)$.
	By projection formula, we have
	\begin{displaymath}
		\mathbf{R}f_\ast (E \otimes^{\mathbf{L}} P \otimes^{\mathbf{L}} \mathbf{L}f^\ast \mathbf{R}i_\ast \mathcal{O}_{\operatorname{Spec}(k)})\cong \mathbf{R}f_\ast (E\otimes^{\mathbf{L}} P) \otimes^{\mathbf{L}} \mathbf{R}i_\ast \mathcal{O}_{\operatorname{Spec}(k)}.
	\end{displaymath}
	Since $E$ is $f$-perfect, we obtain that $E \otimes^{\mathbf{L}} \mathbf{L}f^\ast \mathbf{R}i_\ast \mathcal{O}_{\operatorname{Spec}(k)}\in D^b_{\operatorname{qc}}(Y)$, and so
	\begin{displaymath}
		E \otimes^{\mathbf{L}} P \otimes^{\mathbf{L}} \mathbf{L}f^\ast \mathbf{R}i_\ast \mathcal{O}_{\operatorname{Spec}(k)}\in D^b_{\operatorname{qc}}(Y).
	\end{displaymath}
	By \Cref{lem:finite_coh_dim_implies_restriction_to_bounded}, it follows that
	\begin{displaymath}
		\mathbf{R}f_\ast (E \otimes^{\mathbf{L}} P \otimes^{\mathbf{L}} \mathbf{L}f^\ast \mathbf{R}i_\ast \mathcal{O}_{\operatorname{Spec}(k)})\in D^b_{\operatorname{qc}}(X).
	\end{displaymath}
	As $p$ was arbitrary, \eqref{prop:perfectness3} of \Cref{prop:perfectness} implies $\mathbf{R}f_\ast (E\otimes^{\mathbf{L}} P)\in \operatorname{Perf}(X)$.
\end{proof}

\begin{example}
	\label{ex:relative_perfect_and_quasi_do_not_coincide}
	The two notions of $f$-perfect and $f$-quasi-perfectness do not coincide in general.
	For example, any morphism to a field has structure sheaves being $f$-perfect.
	However, if the morphism is not proper, structure sheaves need not be $f$-quasi-perfect.
\end{example}

\begin{lemma}
	\label{lem:neeman23lem5_3}
	Let $f\colon Y \to X$ be a morphism of quasi-compact quasi-separated algebraic spaces.
	For any $E\in D_{\operatorname{qc}}(Y)$ and $P\in D_{\operatorname{qc}}(X)$, there exists an isomorphism
	\begin{displaymath}
		\mathbf{R}f_\ast \operatorname{\mathbf{R}\mathcal{H}\! \mathit{om}} (E ,f^\times P) \cong \operatorname{\mathbf{R}\mathcal{H}\! \mathit{om}}( \mathbf{R}f_\ast E, P).
	\end{displaymath}
\end{lemma}

\begin{proof}
	We mimic \cite[1.6.1]{Iyengar/Lipman/Neeman:2015}.
	Let $A\in D_{\operatorname{qc}}(X)$.
	By adjunctions,
	\begin{displaymath}
		\begin{aligned}
			\operatorname{Hom}(A, \mathbf{R}f_\ast \operatorname{\mathbf{R}\mathcal{H}\! \mathit{om}} (E ,f^\times P))
			&\cong \operatorname{Hom}(\mathbf{L}f^\ast  A, \operatorname{\mathbf{R}\mathcal{H}\! \mathit{om}} (E ,f^\times P)) && \textrm{(pull/push)}
			\\&\cong \operatorname{Hom}(\mathbf{L}f^\ast  A \otimes^{\mathbf{L}} E , f^\times P) && \textrm{(tensor/hom)}
			\\&\cong \operatorname{Hom}(\mathbf{R}f_\ast (\mathbf{L}f^\ast  A \otimes^{\mathbf{L}} E) , P) && \textrm{(right adjoint of push)}
			\\&\cong \operatorname{Hom}(\mathbf{R}f_\ast E \otimes^{\mathbf{L}} A, P) && \textrm{(projection formula)}
			\\&\cong \operatorname{Hom}(A,  \operatorname{\mathbf{R}\mathcal{H}\! \mathit{om}}( \mathbf{R}f_\ast E, P)) && \textrm{(tensor/hom)}.
		\end{aligned}
	\end{displaymath}
	This finishes the proof.
\end{proof}

\begin{lemma}
	\label{lem:Ballard_relative_perf_implies_dual_is_such}
	Let $f\colon Y\to X$ be a morphism of quasi-compact quasi-separated algebraic spaces.
	Consider an $f$-quasi-perfect pseudocoherent complex $E$ on $Y$.
	Then $\operatorname{\mathbf{R}\mathcal{H}\! \mathit{om}}(E, f^\times P)$ is $f$-quasi-perfect for all $P\in \operatorname{Perf}(X)$.
\end{lemma}

\begin{proof}
	By \cite[\href{https://stacks.math.columbia.edu/tag/0E56}{Tag 0E56}]{StacksProject}, we know that $f^\times P \in D^+_{\operatorname{qc}}(Y)$.
	Since $E$ is pseudocoherent, it follows that $\operatorname{\mathbb{R}\mathcal{H}\! \mathit{om}}(E, f^\times P)\in D_{\operatorname{qc}}(Y)$ \cite[\href{https://stacks.math.columbia.edu/tag/0A8A}{Tag 0A8A}]{StacksProject}.
	In this situation, \Cref{lem:internal_hom} yields
	\begin{displaymath}
		\operatorname{\mathbf{R}\mathcal{H}\! \mathit{om}}(E, f^\times P) \cong \operatorname{\mathbb{R}\mathcal{H}\! \mathit{om}}(E, f^\times P).
	\end{displaymath}
	Choose $Q\in \operatorname{Perf}(Y)$.
	From \cite[\href{https://stacks.math.columbia.edu/tag/08JJ}{Tag 08JJ} and \href{https://stacks.math.columbia.edu/tag/0A8A}{0A8A}]{StacksProject}, the double duality morphism $Q\to (Q^\vee)^\vee$ is an isomorphism where $Q^\vee := \operatorname{\mathbf{R}\mathcal{H}\! \mathit{om}}(Q, \mathcal{O}_{Y})$, and there exists a natural isomorphism of functors $\operatorname{\mathbf{R}\mathcal{H}\! \mathit{om}}(Q, -)\cong \operatorname{\mathbb{R}\mathcal{H}\! \mathit{om}}(Q, -)$.
	Then
	\begin{displaymath}
		\begin{aligned}
			\operatorname{\mathbb{R}\mathcal{H}\! \mathit{om}}(E, f^\times P) \otimes^{\mathbf{L}} Q
			&\cong \operatorname{\mathbb{R}\mathcal{H}\! \mathit{om}}(E, f^\times P) \otimes^{\mathbf{L}} (Q^\vee)^\vee
			\\&\cong \operatorname{\mathbb{R}\mathcal{H}\! \mathit{om}} (Q^\vee, \operatorname{\mathbb{R}\mathcal{H}\! \mathit{om}}(E, f^\times P) ) && (\textrm{\cite[\href{https://stacks.math.columbia.edu/tag/08JJ}{Tag 08JJ}]{StacksProject}})
			\\&\cong \operatorname{\mathbb{R}\mathcal{H}\! \mathit{om}} (Q^\vee \otimes^{\mathbf{L}} E ,f^\times P) && (\textrm{\cite[\href{https://stacks.math.columbia.edu/tag/08J9}{Tag 08J9}]{StacksProject}}).
		\end{aligned}
	\end{displaymath}
	Applying \Cref{lem:internal_hom},
	\begin{displaymath}
		\operatorname{\mathbb{R}\mathcal{H}\! \mathit{om}} (Q^\vee \otimes^{\mathbf{L}} E ,f^\times P) \cong \operatorname{\mathbf{R}\mathcal{H}\! \mathit{om}} (Q^\vee \otimes^{\mathbf{L}} E ,f^\times P).
	\end{displaymath}
	Moreover, \Cref{lem:neeman23lem5_3} yields an isomorphism
	\begin{displaymath}
		\mathbf{R}f_\ast \operatorname{\mathbf{R}\mathcal{H}\! \mathit{om}} (Q^\vee \otimes^{\mathbf{L}} E ,f^\times P) \cong \operatorname{\mathbf{R}\mathcal{H}\! \mathit{om}}( \mathbf{R}f_\ast (Q^\vee \otimes^{\mathbf{L}} E), P).
	\end{displaymath}
	Since $Q^\vee$ is perfect, the hypothesis implies that $\mathbf{R}f_\ast (Q^\vee \otimes^{\mathbf{L}} E)$ is perfect.
	As $P$ is perfect, we conclude that $\operatorname{\mathbf{R}\mathcal{H}\! \mathit{om}}( \mathbf{R}f_\ast (Q^\vee \otimes^{\mathbf{L}} E), P)$ is perfect, which completes the proof.
\end{proof}

\begin{lemma}
	\label{lem:f-quasi-perfect_via_small_coproducts}
	Let $f\colon Y\to X$ be a morphism of quasi-compact quasi-separated algebraic spaces.
	For any $E\in D_{\operatorname{qc}}(Y)$, $E$ is $f$-quasi-perfect if, and only if, $\operatorname{\mathbf{R}\mathcal{H}\! \mathit{om}} (E,f^\times (-))\colon D_{\operatorname{qc}}(X) \to D_{\operatorname{qc}}(Y)$ preserves small coproducts.
\end{lemma}

\begin{proof}
	There exists an adjunction
	\begin{displaymath}
		\mathbf{R}f_\ast (E\otimes^{\mathbf{L}} (-))\colon D_{\operatorname{qc}}(Y) \leftrightarrows D_{\operatorname{qc}}(X)\colon \operatorname{\mathbf{R}\mathcal{H}\! \mathit{om}} (E,f^\times (-)) .
	\end{displaymath}
	Hence, the desired claim follows from \cite[Theorem 5.1]{Neeman:1996}.
\end{proof}

\begin{proposition}
	\label{prop:stacky_f_perf_gives_iso_between_uppershriek_and_pullback_up_to_tensor}
	Let $f\colon Y\to X$ be a morphism of quasi-compact quasi-separated algebraic spaces.
	If $E\in D_{\operatorname{qc}}(Y)$ is $f$-quasi-perfect, then there exists an isomorphism for all $G,B\in D_{\operatorname{qc}}(X)$,
	\begin{displaymath}
		\operatorname{\mathbf{R}\mathcal{H}\! \mathit{om}} (E,f^\times G) \otimes^{\mathbf{L}} \mathbf{L}f^\ast B \to \operatorname{\mathbf{R}\mathcal{H}\! \mathit{om}} (E, f^\times (G\otimes^{\mathbf{L}} B)).
	\end{displaymath}
\end{proposition}

\begin{proof}
	We mimic the proof of \cite[Lemma 3.5]{Ballard:2009} and add details.
	Recall that there exists the adjunction
	\begin{displaymath}
		\mathbf{R}f_\ast (E\otimes^{\mathbf{L}}(-))
		\colon D_{\operatorname{qc}}(Y)
		\leftrightarrows D_{\operatorname{qc}}(X)
		\colon \operatorname{\mathbf{R}\mathcal{H}\! \mathit{om}} (E,f^\times(-)).
	\end{displaymath}
	Consider the counit
	\begin{displaymath}
		\mathbf{R}f_\ast(E\otimes^{\mathbf{L}} \operatorname{\mathbf{R}\mathcal{H}\! \mathit{om}} (E,f^\times G)) \to G.
	\end{displaymath}
	Tensoring with $B$ gives a morphism
	\begin{equation}
		\label{eq:stacky_f_perf_gives_iso_between_uppershriek_and_pullback_up_to_tensoring1}
		\mathbf{R}f_\ast(E\otimes^{\mathbf{L}} \operatorname{\mathbf{R}\mathcal{H}\! \mathit{om}}(E,f^\times G)) \otimes^{\mathbf{L}}B
		\to G\otimes^{\mathbf{L}}B.
	\end{equation}
	By projection formula, there exists a canonical isomorphism
	\begin{equation}
		\label{eq:stacky_f_perf_gives_iso_between_uppershriek_and_pullback_up_to_tensoring2}
		\begin{aligned}
			&\mathbf{R}f_\ast ( E\otimes^{\mathbf{L}} \operatorname{\mathbf{R}\mathcal{H}\! \mathit{om}} (E,f^\times G) ) \otimes^{\mathbf{L}}B
			\\&\cong\mathbf{R}f_\ast ( E\otimes^{\mathbf{L}} \operatorname{\mathbf{R}\mathcal{H}\! \mathit{om}} (E,f^\times G) \otimes^{\mathbf{L}}\mathbf{L}f^\ast B ).
		\end{aligned}
	\end{equation}
	Composing the inverse of \eqref{eq:stacky_f_perf_gives_iso_between_uppershriek_and_pullback_up_to_tensoring2} with \eqref{eq:stacky_f_perf_gives_iso_between_uppershriek_and_pullback_up_to_tensoring1} yields
	\begin{displaymath}
		\mathbf{R}f_\ast ( E\otimes^{\mathbf{L}} \operatorname{\mathbf{R}\mathcal{H}\! \mathit{om}} (E,f^\times G) \otimes^{\mathbf{L}}\mathbf{L}f^\ast B )
		\to G\otimes^{\mathbf{L}}B.
	\end{displaymath}
	By adjunctions, this determines the desired morphism
	\begin{displaymath}
		\operatorname{\mathbf{R}\mathcal{H}\! \mathit{om}} (E,f^\times G) \otimes^{\mathbf{L}}\mathbf{L}f^\ast B
		\to \operatorname{\mathbf{R}\mathcal{H}\! \mathit{om}} (E,f^\times(G\otimes^{\mathbf{L}}B)).
	\end{displaymath}

	We first assume $B\in\operatorname{Perf}(X)$.
	For every $A\in D_{\operatorname{qc}}(Y)$, there is a sequence
	of natural isomorphisms
	\begin{displaymath}
		\begin{aligned}
			\operatorname{Hom} & ( A , \operatorname{\mathbf{R}\mathcal{H}\! \mathit{om}} (E,f^\times G)\otimes \mathbf{L}f^\ast B )
			\\&\cong \operatorname{Hom} ( A , \operatorname{\mathbf{R}\mathcal{H}\! \mathit{om}} ( \operatorname{\mathbf{R}\mathcal{H}\! \mathit{om}} ( \mathbf{L} f^\ast B, \mathcal{O}_{Y} ) , \operatorname{\mathbf{R}\mathcal{H}\! \mathit{om}} ( E , f^\times G) ) ) && (\textrm{\cite[\href{https://stacks.math.columbia.edu/tag/08JJ}{Tag 08JJ}]{StacksProject}})
			\\&\cong \operatorname{Hom} ( A \otimes^{\mathbf{L}} \operatorname{\mathbf{R}\mathcal{H}\! \mathit{om}} ( \mathbf{L} f^\ast B, \mathcal{O}_{Y} ) , \operatorname{\mathbf{R}\mathcal{H}\! \mathit{om}} ( E , f^\times G)  ) && \textrm{(tensor/hom)}
			\\&\cong \operatorname{Hom} ( A \otimes^{\mathbf{L}} \operatorname{\mathbf{R}\mathcal{H}\! \mathit{om}} ( \mathbf{L} f^\ast B, \mathcal{O}_{Y} ) \otimes^{\mathbf{L}} E , f^\times G  ) && \textrm{(tensor/hom)}
			\\&\cong \operatorname{Hom} ( \mathbf{R}f_\ast (A \otimes^{\mathbf{L}} \operatorname{\mathbf{R}\mathcal{H}\! \mathit{om}} ( \mathbf{L} f^\ast B, \mathcal{O}_{Y} ) \otimes^{\mathbf{L}} E ), G  ) && (\textrm{right adj. of push})
			\\&\cong \operatorname{Hom} ( \mathbf{R}f_\ast (A \otimes^{\mathbf{L}} \mathbf{L} f^\ast (\operatorname{\mathbf{R}\mathcal{H}\! \mathit{om}} ( B , \mathcal{O}_{X} )) \otimes^{\mathbf{L}} E ), G  ) && \textrm{(\Cref{lem:gortz_wedhorn_internal_hom})}
			\\&\cong \operatorname{Hom} ( \mathbf{R}f_\ast (A \otimes^{\mathbf{L}}  E ) \otimes^{\mathbf{L}}  \operatorname{\mathbf{R}\mathcal{H}\! \mathit{om}} ( B , \mathcal{O}_{X} ), G  ) && \textrm{(projection formula)}
			\\&\cong \operatorname{Hom} ( \mathbf{R}f_\ast (A \otimes^{\mathbf{L}}  E ) , \operatorname{\mathbf{R}\mathcal{H}\! \mathit{om}} (\operatorname{\mathbf{R}\mathcal{H}\! \mathit{om}} ( B , \mathcal{O}_{X} ) , G)  ) && \textrm{(tensor/hom)}
			\\&\cong \operatorname{Hom} ( \mathbf{R}f_\ast (A \otimes^{\mathbf{L}}  E ) , B\otimes^{\mathbf{L}} G )  && (\textrm{\cite[\href{https://stacks.math.columbia.edu/tag/08JJ}{Tag 08JJ}]{StacksProject}})
			\\&\cong \operatorname{Hom} ( A \otimes^{\mathbf{L}}  E  , f^\times (B\otimes^{\mathbf{L}} G) ) && (\textrm{right adj. of push})
			\\&\cong \operatorname{Hom} ( A , \operatorname{\mathbf{R}\mathcal{H}\! \mathit{om}}  (E, f^\times (B\otimes^{\mathbf{L}} G)) ) && \textrm{(tensor/hom)}.
		\end{aligned}
	\end{displaymath}
	Under this sequence of adjunction isomorphisms, postcomposition
	with the canonical morphism constructed above is the resulting
	bijection between the first and last Hom-sets.
	Indeed, the canonical morphism was defined as the adjunct of the
	counit
	\begin{displaymath}
		\mathbf{R}f_\ast(E\otimes^{\mathbf{L}} \operatorname{\mathbf{R}\mathcal{H}\! \mathit{om}} (E,f^\times G))
		\to G
	\end{displaymath}
	after tensoring with $B$ and applying the projection formula.
	Unwinding the adjunctions in the preceding sequence shows that postcomposition with this morphism inserts precisely that counit.
	The triangle identity for
	\begin{displaymath}
		\mathbf{R}f_\ast(E\otimes^{\mathbf{L}}(-))
		\dashv
		\operatorname{\mathbf{R}\mathcal{H}\! \mathit{om}}
		(E,f^\times(-))
	\end{displaymath}
	then gives the identity under the common Hom-set description.
	Consequently, postcomposition with the canonical morphism is a bijection for every $A\in D_{\operatorname{qc}}(Y)$.
	Hence, the canonical morphism is an isomorphism whenever $B$ is perfect.

	Fix now $G\in D_{\operatorname{qc}}(X)$, and let $\mathcal{T}\subseteq D_{\operatorname{qc}}(X)$ be the strictly full subcategory consisting of those $B$ for which the canonical morphism is an isomorphism.
	By \Cref{lem:f-quasi-perfect_via_small_coproducts}, $\operatorname{\mathbf{R}\mathcal{H}\! \mathit{om}} (E,f^\times(-))$ preserves small coproducts.
	Since $\mathbf{L}f^\ast$ and $\otimes^{\mathbf{L}}$ preserve
	small coproducts as well, both sides of the canonical morphism,
	viewed as functors of $B$, preserve small coproducts.
	Therefore, $\mathcal{T}$ is a localizing subcategory.

	The argument above shows $\operatorname{Perf}(X)\subseteq\mathcal{T}$.
	Since $\operatorname{Perf}(X)$ compactly generates $D_{\operatorname{qc}}(X)$, it follows that $\mathcal{T}=D_{\operatorname{qc}}(X)$.
	This completes the proof.
\end{proof}

\begin{lemma}
	\label{lem:zero_iff_pushforward_zero}
	Let $f\colon Y\to X$ be a morphism of quasi-compact quasi-separated algebraic spaces.
	Let $P\in D_{\operatorname{qc}}(Y)$ be a compact generator.
	Then $E\in D_{\operatorname{qc}}(Y)$ is the zero object if, and only if, $\mathbf{R}f_\ast (E\otimes^{\mathbf{L}} P)\cong 0$.
\end{lemma}

\begin{proof}
	If $E\in D_{\operatorname{qc}}(Y)$ is the zero object, then $\mathbf{R}f_\ast (E\otimes^{\mathbf{L}} P)\cong 0$.
	We prove the converse.
	By adjunction,
	\begin{displaymath}
		\begin{aligned}
			0
			&\cong \operatorname{Hom}(\mathcal{O}_{X} , \mathbf{R}f_\ast (E\otimes^{\mathbf{L}} P)[n])\\
			&\cong \operatorname{Hom}(\mathcal{O}_{X} , \mathbf{R}f_\ast (E[n]\otimes^{\mathbf{L}} P))
			\\&\cong \operatorname{Hom}(\mathbf{L}f^\ast \mathcal{O}_{X} , E[n]\otimes^{\mathbf{L}} P) && \textrm{(pull/push)}
			\\&\cong \operatorname{Hom}(\mathcal{O}_{Y} , \operatorname{\mathbf{R}\mathcal{H}\! \mathit{om}}(\operatorname{\mathbf{R}\mathcal{H}\! \mathit{om}}(P,\mathcal{O}_{Y}) , E[n]) ) && \textrm{(\cite[\href{https://stacks.math.columbia.edu/tag/08JJ}{Tag 08JJ}]{StacksProject})}
			%%NOTE: P = P double dual
			\\&\cong \operatorname{Hom}(\operatorname{\mathbf{R}\mathcal{H}\! \mathit{om}}(P,\mathcal{O}_{Y}), E[n]) && \textrm{(\cite[\href{https://stacks.math.columbia.edu/tag/08J7}{Tag 08J7}]{StacksProject}).}
		\end{aligned}
	\end{displaymath}
	%% OLD PROOF: It follows from the assumption that $E\cong 0$.
	% Indeed, use that perfect complexes are dualizable in $D_{\operatorname{qc}}(Y)$ \cite[\href{https://stacks.math.columbia.edu/tag/0FPU}{Tags 0FPU}, \href{https://stacks.math.columbia.edu/tag/0FPV}{0FPV}, \href{https://stacks.math.columbia.edu/tag/09M8}{09M8}]{StacksProject}.
	% In particular, we have shown that $\operatorname{Hom}(Q,E[n])\cong 0$ for any compact generator $Q$ of $D_{\operatorname{qc}}(Y)$ and $n\in \mathbb{Z}$.
	% To see, the computation above shows
	% \begin{displaymath}
	%     \operatorname{Hom}(\operatorname{\mathbf{R}\mathcal{H}\! \mathit{om}}(Q,\mathcal{O}_{Y})[n], E) \cong 0
	% \end{displaymath}
	% for all $\in \mathbb{Z}$ and fixed compact generator $Q$.
	However, $\operatorname{\mathbf{R}\mathcal{H}\! \mathit{om}}(P,\mathcal{O}_{Y})$ is a compact generator, and so the claim follows.
\end{proof}

\begin{lemma}
	\label{lem:upper_shriek_preserves_Dplusqc}
	Let $f\colon Y\to X$ be a separated finitely presented morphism of Noetherian algebraic spaces.
	Then $f^! (D^+_{\operatorname{qc}}(X))\subseteq D^+_{\operatorname{qc}}(Y)$.
\end{lemma}

\begin{proof}
	Choose a Nagata compactification $f= p \circ j$ with $p\colon Y^\prime \to X$ proper and $j\colon Y \to Y^\prime$ an open immersion.
	By construction, $f^!$ is naturally isomorphic to $\mathbf{L}j^\ast \circ p^\times$.
	As $j$ is an open immersion, $j^\ast$ is exact, and so $\mathbf{L}j^\ast D^+_{\operatorname{qc}}(Y^\prime) \subseteq D^+_{\operatorname{qc}}(Y)$.
	By \cite[\href{https://stacks.math.columbia.edu/tag/0E56}{Tag 0E56}]{StacksProject}, $p^\times D^+_{\operatorname{qc}}(X) \subseteq D^+_{\operatorname{qc}}(Y^\prime)$.
	Thus, the claim follows.
\end{proof}

\begin{lemma}
	\label{lem:upper_shriek_coherent_cohomology}
	Let $f\colon Y\to X$ be a separated finitely presented morphism of Noetherian algebraic spaces. Then $f^! D^+_{\operatorname{coh}}(X)\subseteq D^+_{\operatorname{coh}}(Y)$.
\end{lemma}

\begin{proof}
	Let $E\in D^+_{\operatorname{coh}}(X)$.
	By \Cref{lem:upper_shriek_preserves_Dplusqc}, $f^! E \in D^+_{\operatorname{qc}}(Y)$
	It suffices to check that $f^! E$ has coherent cohomology.

	Choose an \'{e}tale presentation $s\colon U \to X$ from an affine scheme. Consider the fibered square
	\begin{displaymath}
		% https://q.uiver.app/#q=WzAsNCxbMSwwLCJVIl0sWzEsMSwiXFxtYXRoY2Fse1h9LiJdLFswLDEsIlxcbWF0aGNhbHtZfSJdLFswLDAsIlxcbWF0aGNhbHtZfVxcdGltZXNfe1xcbWF0aGNhbHtYfX0gVSJdLFswLDEsInMiXSxbMiwxLCJmIiwyXSxbMywwLCJmXlxccHJpbWUiXSxbMywyLCJzXlxccHJpbWUiLDJdXQ==
		\begin{tikzcd}
			{Y\times_{X} U} & U \\
			{Y} & {X.}
			\arrow["{f^\prime}", from=1-1, to=1-2]
			\arrow["{s^\prime}"', from=1-1, to=2-1]
			\arrow["s", from=1-2, to=2-2]
			\arrow["f"', from=2-1, to=2-2]
		\end{tikzcd}
	\end{displaymath}
	By \Cref{lem:quasi_affine_diagonal,lem:quasi-affine_via_diagonal}, $s$ is separated.
	Base change implies each morphism in the diagram is finitely presented, and each algebraic space is Noetherian.
	By \Cref{lem:neeman188},
	%%NOTE: Use D^+_{qc}
	$\mathbf{L}(s^\prime)^\ast f^! E \cong (f^\prime)^! \mathbf{L}s^\ast E$.
	Note that $\mathbf{L}s^\ast E\in D^+_{\operatorname{qc}}(U)$.

	Choose an \'{e}tale surjective morphism $t\colon V \to Y\times_{X} U$ from an affine scheme.
	Consider the fibered square
	\begin{displaymath}
		% https://q.uiver.app/#q=WzAsNCxbMSwwLCJWIl0sWzEsMSwiXFxtYXRoY2Fse1l9XFx0aW1lc197XFxtYXRoY2Fse1h9fSBVLiJdLFswLDEsIlYiXSxbMCwwLCJWXlxccHJpbWUiXSxbMCwxLCJ0Il0sWzIsMSwidCIsMl0sWzMsMiwidF8xIiwyXSxbMywwLCJ0XzIiXV0=
		\begin{tikzcd}
			{V^\prime} & V \\
			V & {Y\times_{X} U.}
			\arrow["{t_2}", from=1-1, to=1-2]
			\arrow["{t_1}"', from=1-1, to=2-1]
			\arrow["t", from=1-2, to=2-2]
			\arrow["t"', from=2-1, to=2-2]
		\end{tikzcd}
	\end{displaymath}
	By \Cref{lem:quasi_affine_diagonal}, $t$ is quasi-affine.
	In particular, $t$ is representable by schemes.
	As $V$ is affine, each $t_i$ is separated and \'{e}tale with scheme source.

	Applying \Cref{lem:neeman188},
	%%NOTE: Use that t and t_i are of finite tor-dimension
	there exists an isomorphism
	\begin{displaymath}
		t^!_2 \mathbf{L}t^\ast (f^\prime)^! \mathbf{L}s^\ast E \cong \mathbf{L}t_1^\ast t^! (f^\prime)^! \mathbf{L}s^\ast E.
	\end{displaymath}
	Since $t_i$ is \'{e}tale, \cite[\href{https://stacks.math.columbia.edu/tag/0FWI}{Tag 0FWI}]{StacksProject} shows that $t^!_i \cong \mathbf{L}t_i^\ast$.
	Combining gives the identifications,
	\begin{displaymath}
		\mathbf{L} t^\ast_2 \mathbf{L}t^\ast (f^\prime)^! \mathbf{L}s^\ast E
		\cong
		t^!_2 \mathbf{L}t^\ast (f^\prime)^! \mathbf{L}s^\ast E
		\cong \mathbf{L}t_1^\ast t^! (f^\prime)^! \mathbf{L}s^\ast E.
	\end{displaymath}
	Note that \Cref{lem:neeman186}
	%%NOTE: Use D^+_{qc}
	says $(f^\prime\circ t)^! \mathbf{L} s^\ast E \cong t^! (f^\prime)^! \mathbf{L} s^\ast E$.
	By \cite[\href{https://stacks.math.columbia.edu/tag/0AU1}{Tag 0AU1}]{StacksProject}, $(f^\prime\circ t)^! \mathbf{L} s^\ast E\in D^+_{\operatorname{coh}}(V)$.
	Hence,
	\begin{displaymath}
		\mathbf{L} t^\ast_2 \mathbf{L}t^\ast (f^\prime)^! \mathbf{L}s^\ast E  \cong \mathbf{L}t_1^\ast t^! (f^\prime)^! \mathbf{L}s^\ast E\in D^+_{\operatorname{coh}}(V^\prime).
	\end{displaymath}
	Since $t \circ t_2$ is an \'{e}tale presentation, the derived pullback is $t$-exact with respect to the standard $t$-structures, and so \cite[\href{https://stacks.math.columbia.edu/tag/07UB}{Tag 07UB}]{StacksProject} shows that $(f^\prime)^! \mathbf{L}s^\ast E \in D^+_{\operatorname{coh}}(Y\times_{X}U)$.
	Similar reasoning shows $f^! E\in D^+_{\operatorname{coh}}(Y)$ as desired.
\end{proof}

\begin{proposition}
	\label{prop:ballard_quasi_perfect_involution}
	Let $f\colon Y\to X$ be a proper morphism of Noetherian algebraic spaces.
	If $E\in D_{\operatorname{qc}}(Y)$ is $f$-quasi-perfect and pseudocoherent, then the canonical morphism
	\begin{displaymath}
		\nu\colon E \to \operatorname{\mathbf{R}\mathcal{H}\! \mathit{om}} (\operatorname{\mathbf{R}\mathcal{H}\! \mathit{om}} (E,f^\times \mathcal{O}_{X}) , f^\times \mathcal{O}_{X})
	\end{displaymath}
	is an isomorphism.
\end{proposition}

\begin{proof}
	We mimic the proof \cite[Lemma 3.9]{Ballard:2009} and add
	details.
	Recall that the canonical morphism is the morphism corresponding to the identity under the adjunctions
	\begin{displaymath}
		\begin{aligned}
			&\operatorname{Hom}(E , \operatorname{\mathbf{R}\mathcal{H}\! \mathit{om}} (\operatorname{\mathbf{R}\mathcal{H}\! \mathit{om}} (E,f^\times \mathcal{O}_{X}) , f^\times \mathcal{O}_{X}))
			\\&\cong \operatorname{Hom}(E \otimes^{\mathbf{L}} \operatorname{\mathbf{R}\mathcal{H}\! \mathit{om}} (E,f^\times \mathcal{O}_{X}),  f^\times \mathcal{O}_{X}) && (\textrm{tensor/hom})
			\\&\cong \operatorname{Hom}(\operatorname{\mathbf{R}\mathcal{H}\! \mathit{om}} (E, f^\times \mathcal{O}_{X}) \otimes^{\mathbf{L}} E ,  f^\times \mathcal{O}_{X}) && (\textrm{switch})
			\\&\cong \operatorname{Hom}(\operatorname{\mathbf{R}\mathcal{H}\! \mathit{om}} (E, f^\times \mathcal{O}_{X}) , \operatorname{\mathbf{R}\mathcal{H}\! \mathit{om}} (E , f^\times \mathcal{O}_{X})) && (\textrm{tensor/hom}).
		\end{aligned}
	\end{displaymath}
	For each $P\in\operatorname{Perf}(Y)$, naturality of evaluation gives a commutative diagram
	\begin{displaymath}
		% https://q.uiver.app/#q=WzAsNCxbMCwwLCJFIFxcb3RpbWVzXntcXG1hdGhiZntMfX0gUCJdLFsxLDAsIlxcb3BlcmF0b3JuYW1le1xcbWF0aGJme1J9XFxtYXRoY2Fse0h9XFwhIFxcbWF0aGl0e29tfX0gKFxcb3BlcmF0b3JuYW1le1xcbWF0aGJme1J9XFxtYXRoY2Fse0h9XFwhIFxcbWF0aGl0e29tfX0gKEUsZl5cXHRpbWVzIFxcbWF0aGNhbHtPfV97XFxtYXRoY2Fse1h9fSkgLCBmXlxcdGltZXMgXFxtYXRoY2Fse099X3tcXG1hdGhjYWx7WH19KSBcXG90aW1lc157XFxtYXRoYmZ7TH19IFAiXSxbMCwxLCJFIFxcb3RpbWVzXntcXG1hdGhiZntMfX0gUCJdLFsxLDEsIlxcb3BlcmF0b3JuYW1le1xcbWF0aGJme1J9XFxtYXRoY2Fse0h9XFwhIFxcbWF0aGl0e29tfX0gKFxcb3BlcmF0b3JuYW1le1xcbWF0aGJme1J9XFxtYXRoY2Fse0h9XFwhIFxcbWF0aGl0e29tfX0gKEUgXFxvdGltZXNee1xcbWF0aGJme0x9fSBQICxmXlxcdGltZXMgXFxtYXRoY2Fse099X3tcXG1hdGhjYWx7WH19KSAsIGZeXFx0aW1lcyBcXG1hdGhjYWx7T31fe1xcbWF0aGNhbHtYfX0pLiJdLFswLDEsIlxcbnVfRSBcXG90aW1lc157XFxtYXRoYmZ7TH19IFAiXSxbMCwyLCIxX3tFXFxvdGltZXNee1xcbWF0aGJme0x9fSBQfSIsMl0sWzEsM10sWzIsMywiXFxudV97RVxcb3RpbWVzXntcXG1hdGhiZntMfX0gUH0iLDJdXQ==
		\begin{tikzcd}
			{E \otimes^{\mathbf{L}} P} & {\operatorname{\mathbf{R}\mathcal{H}\! \mathit{om}} (\operatorname{\mathbf{R}\mathcal{H}\! \mathit{om}} (E,f^\times \mathcal{O}_{X}) , f^\times \mathcal{O}_{X}) \otimes^{\mathbf{L}} P} \\
			{E \otimes^{\mathbf{L}} P} & {\operatorname{\mathbf{R}\mathcal{H}\! \mathit{om}} (\operatorname{\mathbf{R}\mathcal{H}\! \mathit{om}} (E \otimes^{\mathbf{L}} P ,f^\times \mathcal{O}_{X}) , f^\times \mathcal{O}_{X}).}
			\arrow["{\nu_E \otimes^{\mathbf{L}} P}", from=1-1, to=1-2]
			\arrow["{1_{E\otimes^{\mathbf{L}} P}}"', from=1-1, to=2-1]
			\arrow[from=1-2, to=2-2]
			\arrow["{\nu_{E\otimes^{\mathbf{L}} P}}"', from=2-1, to=2-2]
		\end{tikzcd}
	\end{displaymath}
	The right vertical morphism is induced by \Cref{lem:rhom_tensor_morphism} and \cite[\href{https://stacks.math.columbia.edu/tag/08J9}{Tag 08J9}]{StacksProject}, and is an isomorphism because $P$ is perfect.

	We justify that all the internal Homs above have the required bounded and quasi-coherent cohomology.
	By \Cref{lem:neeman187}, there exists a natural isomorphism $f^\times\cong f^!$.
	Hence, \Cref{lem:upper_shriek_coherent_cohomology} gives $f^\times\mathcal{O}_X\in D^+_{\operatorname{coh}}(Y)$.
	By \Cref{lem:Ballard_relative_perf_implies_dual_is_such}, $\operatorname{\mathbf{R}\mathcal{H}\! \mathit{om}}(E,f^\times\mathcal{O}_X)$ is $f$-quasi-perfect.
	Then \Cref{lem:quasi-perfect_bounded} gives the required boundedness, and \Cref{lem:internal_hom} identifies the internal Homs with the corresponding derived sheaf Homs.
	The same reasoning applies to $E\otimes^{\mathbf{L}}P$.

	By the preceding commutative diagram, for every $Q\in\operatorname{Perf}(Y)$, $\operatorname{cone}(\nu_E) \otimes^{\mathbf{L}}Q \cong \operatorname{cone} (\nu_{E\otimes^{\mathbf{L}}Q})$.
	Hence, by \Cref{lem:zero_iff_pushforward_zero}, it suffices to prove that $\mathbf{R}f_\ast (\nu_{E\otimes^{\mathbf{L}}Q})$ is an isomorphism for every $Q\in\operatorname{Perf}(Y)$.

	By \Cref{lem:neeman23lem5_3}, there is a natural isomorphism
	\begin{displaymath}
		\eta_{E\otimes^{\mathbf{L}}Q}\colon
		\mathbf{R}f_\ast
		\operatorname{\mathbf{R}\mathcal{H}\! \mathit{om}}
		(E\otimes^{\mathbf{L}}Q,f^\times\mathcal{O}_X)
		\to \operatorname{\mathbf{R}\mathcal{H}\! \mathit{om}} (\mathbf{R}f_\ast(E\otimes^{\mathbf{L}}Q), \mathcal{O}_X).
	\end{displaymath}
	Applying the same isomorphism to $\operatorname{\mathbf{R}\mathcal{H}\! \mathit{om}}(E\otimes^{\mathbf{L}}Q,f^\times\mathcal{O}_X)$.
	Using the contravariance of $\operatorname{\mathbf{R}\mathcal{H}\! \mathit{om}} (-,\mathcal{O}_X)$ gives a commutative diagram
	\begin{displaymath}
		% https://q.uiver.app/#q=WzAsNCxbMCwwLCJcXG1hdGhiZntSfWZfXFxhc3QgKEUgXFxvdGltZXNee1xcbWF0aGJme0x9fSBRKSJdLFswLDEsIlxcbWF0aGJme1J9Zl9cXGFzdCBcXGJpZ2coIFxcb3BlcmF0b3JuYW1le1xcbWF0aGJme1J9XFxtYXRoY2Fse0h9XFwhIFxcbWF0aGl0e29tfX0gKFxcb3BlcmF0b3JuYW1le1xcbWF0aGJme1J9XFxtYXRoY2Fse0h9XFwhIFxcbWF0aGl0e29tfX0gKEUgIFxcb3RpbWVzXntcXG1hdGhiZntMfX0gUSAsZl5cXHRpbWVzIFxcbWF0aGNhbHtPfV97XFxtYXRoY2Fse1h9fSkgLCBmXlxcdGltZXMgXFxtYXRoY2Fse099X3tcXG1hdGhjYWx7WH19KSBcXGJpZ2cpIl0sWzAsMiwiXFxvcGVyYXRvcm5hbWV7XFxtYXRoYmZ7Un1cXG1hdGhjYWx7SH1cXCEgXFxtYXRoaXR7b219fSAoXFxtYXRoYmZ7Un1mX1xcYXN0IFxcb3BlcmF0b3JuYW1le1xcbWF0aGJme1J9XFxtYXRoY2Fse0h9XFwhIFxcbWF0aGl0e29tfX0gKEUgXFxvdGltZXNee1xcbWF0aGJme0x9fSBRICxmXlxcdGltZXMgXFxtYXRoY2Fse099X3tcXG1hdGhjYWx7WH19KSAsIFxcbWF0aGNhbHtPfV97XFxtYXRoY2Fse1h9fSkiXSxbMCwzLCJcXG9wZXJhdG9ybmFtZXtcXG1hdGhiZntSfVxcbWF0aGNhbHtIfVxcISBcXG1hdGhpdHtvbX19IChcXG9wZXJhdG9ybmFtZXtcXG1hdGhiZntSfVxcbWF0aGNhbHtIfVxcISBcXG1hdGhpdHtvbX19IChcXG1hdGhiZntSfWZfXFxhc3QgIChFIFxcb3RpbWVzXntcXG1hdGhiZntMfX0gUSkgICwgXFxtYXRoY2Fse099X3tcXG1hdGhjYWx7WH19KSAsIFxcbWF0aGNhbHtPfV97XFxtYXRoY2Fse1h9fSkuIl0sWzAsMSwiXFxtYXRoYmZ7Un1mX1xcYXN0IChcXG51X3tFIFxcb3RpbWVzXntcXG1hdGhiZntMfX0gUX0pIl0sWzEsMiwiXFxldGFfe1xcb3BlcmF0b3JuYW1le1xcbWF0aGJme1J9XFxtYXRoY2Fse0h9XFwhIFxcbWF0aGl0e29tfX0gKEUgIFxcb3RpbWVzXntcXG1hdGhiZntMfX0gUSAsZl5cXHRpbWVzIFxcbWF0aGNhbHtPfV97XFxtYXRoY2Fse1h9fSl9Il0sWzIsMywiXFxvcGVyYXRvcm5hbWV7XFxtYXRoYmZ7Un1cXG1hdGhjYWx7SH1cXCEgXFxtYXRoaXR7b219fSAoXFxldGFfe0VcXG90aW1lc157XFxtYXRoYmZ7TH19IFF9XnstMX0sIFxcbWF0aGNhbHtPfV97XFxtYXRoY2Fse1h9fSkiXV0=
		\begin{tikzcd}
			{\mathbf{R}f_\ast (E \otimes^{\mathbf{L}} Q)} \\
			{\mathbf{R}f_\ast ( \operatorname{\mathbf{R}\mathcal{H}\! \mathit{om}} (\operatorname{\mathbf{R}\mathcal{H}\! \mathit{om}} (E  \otimes^{\mathbf{L}} Q ,f^\times \mathcal{O}_{X}) , f^\times \mathcal{O}_{X}) )} \\
			{\operatorname{\mathbf{R}\mathcal{H}\! \mathit{om}} (\mathbf{R}f_\ast \operatorname{\mathbf{R}\mathcal{H}\! \mathit{om}} (E \otimes^{\mathbf{L}} Q ,f^\times \mathcal{O}_{X}) , \mathcal{O}_{X})} \\
			{\operatorname{\mathbf{R}\mathcal{H}\! \mathit{om}} (\operatorname{\mathbf{R}\mathcal{H}\! \mathit{om}} (\mathbf{R}f_\ast  (E \otimes^{\mathbf{L}} Q)  , \mathcal{O}_{X}) , \mathcal{O}_{X}).}
			\arrow["{\mathbf{R}f_\ast (\nu_{E \otimes^{\mathbf{L}} Q})}", from=1-1, to=2-1]
			\arrow["{\eta_{\operatorname{\mathbf{R}\mathcal{H}\! \mathit{om}} (E  \otimes^{\mathbf{L}} Q ,f^\times \mathcal{O}_{X})}}", from=2-1, to=3-1]
			\arrow["{\operatorname{\mathbf{R}\mathcal{H}\! \mathit{om}} (\eta_{E\otimes^{\mathbf{L}} Q}^{-1}, \mathcal{O}_{X})}", from=3-1, to=4-1]
		\end{tikzcd}
	\end{displaymath}
	We claim that the composite in this diagram is the canonical
	biduality morphism \cite[\href{https://stacks.math.columbia.edu/tag/0A97}{Tag 0A97} and \href{https://stacks.math.columbia.edu/tag/08JJ}{08JJ}]{StacksProject},
	\begin{displaymath}
		\mathbf{R}f_\ast(E\otimes^{\mathbf{L}}Q)
		\to \operatorname{\mathbf{R}\mathcal{H}\! \mathit{om}} ( \operatorname{\mathbf{R}\mathcal{H}\! \mathit{om}}( \mathbf{R}f_\ast(E\otimes^{\mathbf{L}}Q), \mathcal{O}_X ), \mathcal{O}_X ).
	\end{displaymath}
	To see this, it suffices to identify its adjunct under tensor-Hom adjunction \cite[\href{https://stacks.math.columbia.edu/tag/08J7}{Tag 08J7}]{StacksProject}.
	Unwinding the natural isomorphism of \Cref{lem:neeman23lem5_3}, the adjunct of the outer composite is
	\begin{displaymath}
		\begin{aligned}
			&\operatorname{\mathbf{R}\mathcal{H}\! \mathit{om}}(\mathbf{R}f_\ast(E\otimes^{\mathbf{L}}Q), \mathcal{O}_X)\otimes^{\mathbf{L}} \mathbf{R}f_\ast(E\otimes^{\mathbf{L}}Q)
			\\&\xrightarrow{\eta_{E\otimes^{\mathbf{L}}Q}^{-1}\otimes1}
			\mathbf{R}f_\ast\operatorname{\mathbf{R}\mathcal{H}\! \mathit{om}}
			(E\otimes^{\mathbf{L}}Q,f^\times\mathcal{O}_X)\otimes^{\mathbf{L}} \mathbf{R}f_\ast(E\otimes^{\mathbf{L}}Q)
			\\&\to\mathbf{R}f_\ast( \operatorname{\mathbf{R}\mathcal{H}\! \mathit{om}} (E\otimes^{\mathbf{L}}Q,f^\times\mathcal{O}_X)  \otimes^{\mathbf{L}} (E\otimes^{\mathbf{L}}Q))
			\\&\to \mathbf{R}f_\ast f^\times\mathcal{O}_X
			\to \mathcal{O}_X.
		\end{aligned}
	\end{displaymath}
	Here the second morphism is the relative cup product, the third is evaluation, and the last is the counit of $\mathbf{R}f_\ast\dashv f^\times$.
	By the construction of the natural isomorphism in \Cref{lem:neeman23lem5_3}, this composite is precisely the evaluation morphism
	\begin{displaymath}
		\operatorname{\mathbf{R}\mathcal{H}\! \mathit{om}}
		(\mathbf{R}f_\ast(E\otimes^{\mathbf{L}}Q), \mathcal{O}_X ) \otimes^{\mathbf{L}} \mathbf{R}f_\ast(E\otimes^{\mathbf{L}}Q)
		\to \mathcal{O}_X.
	\end{displaymath}
	Its adjunct is the canonical biduality morphism, as
	claimed.

	Since $E$ is $f$-quasi-perfect and $Q\in\operatorname{Perf}(Y)$, $\mathbf{R}f_\ast(E\otimes^{\mathbf{L}}Q) \in\operatorname{Perf}(X)$.
	Hence, its canonical biduality morphism is an isomorphism
	\cite[\href{https://stacks.math.columbia.edu/tag/08JJ}{Tag 08JJ} and \href{https://stacks.math.columbia.edu/tag/0A8A}{0A8A}]{StacksProject}.
	The two morphisms denoted by $\eta$ above are isomorphisms, so $\mathbf{R}f_\ast(\nu_{E\otimes^{\mathbf{L}}Q})$ is an isomorphism.
	This holds for every $Q\in\operatorname{Perf}(Y)$.
	Therefore, \Cref{lem:zero_iff_pushforward_zero} implies $\operatorname{cone}(\nu_E)\cong0$, and hence, $\nu_E$ is an isomorphism.
\end{proof}

\begin{remark}
	A variation of \Cref{prop:ballard_quasi_perfect_involution} for $f^!$ is proved in \Cref{lem:involution_for_f_perfect}.
\end{remark}

\begin{theorem}
	\label{thm:adjoint_dqc_for_f-quasi-perfect}
	Let $f\colon Y\to X$ be a proper morphism of Noetherian algebraic spaces.
	If $E\in D_{\operatorname{qc}}(Y)$ is $f$-quasi-perfect and pseudocoherent, then there exists an adjunction
	\begin{displaymath}
		\mathbf{R} f_\ast ( (-) \otimes^{\mathbf{L}} \operatorname{\mathbf{R}\mathcal{H}\! \mathit{om}}(E, f^\times \mathcal{O}_{X}) )\colon D_{\operatorname{qc}}(Y) \leftrightarrows D_{\operatorname{qc}}(X)\colon E \otimes^{\mathbf{L}} \mathbf{L}f^\ast (-).
	\end{displaymath}
\end{theorem}

\begin{proof}
	This follows from natural isomorphisms for all $A\in D_{\operatorname{qc}}(Y)$ and $G\in D_{\operatorname{qc}}(X)$:
	\begin{displaymath}
		\begin{aligned}
			\operatorname{Hom} & (A , E \otimes^{\mathbf{L}} \mathbf{L}f^\ast G)
			\\&\cong \operatorname{Hom}(A , \operatorname{\mathbf{R}\mathcal{H}\! \mathit{om}} ( \operatorname{\mathbf{R}\mathcal{H}\! \mathit{om}}(E, f^\times \mathcal{O}_{X}) , f^\times \mathcal{O}_{X}) \otimes^{\mathbf{L}} \mathbf{L}f^\ast G) && (\textrm{\Cref{prop:ballard_quasi_perfect_involution}})
			\\&\cong \operatorname{Hom}(A , \operatorname{\mathbf{R}\mathcal{H}\! \mathit{om}} ( \operatorname{\mathbf{R}\mathcal{H}\! \mathit{om}}(E, f^\times \mathcal{O}_{X}) , f^\times G)) && (\textrm{\Cref{prop:stacky_f_perf_gives_iso_between_uppershriek_and_pullback_up_to_tensor}})
			\\&\cong \operatorname{Hom}(A \otimes^{\mathbf{L}} \operatorname{\mathbf{R}\mathcal{H}\! \mathit{om}}(E, f^\times \mathcal{O}_{X}) ,  f^\times G) && (\textrm{tensor/hom})
			\\&\cong \operatorname{Hom}(\mathbf{R} f_\ast ( A \otimes^{\mathbf{L}} \operatorname{\mathbf{R}\mathcal{H}\! \mathit{om}}(E, f^\times \mathcal{O}_{X}) ) ,  G) && (\textrm{right adj. of push}).
		\end{aligned}
	\end{displaymath}
\end{proof}

%%%%%%%%%%%%%%%%%%%%%%%%%%%%%%%%%%%
\section{Properness and base change}
\label{sec:properness}
%%%%%%%%%%%%%%%%%%%%%%%%%%%%%%%%%%%

\begin{lemma}
	[cf.\ {\cite[\href{https://stacks.math.columbia.edu/tag/0GFH}{Tag 0GFH} and \href{https://stacks.math.columbia.edu/tag/0GFI}{0GFI}]{StacksProject}}]
	\label{lem:boundedness_via_coaisle}
	Let $X$ be a quasi-compact quasi-separated algebraic space.
	Consider a closed subset $Z\subseteq |X|$ with quasi-compact complement.
	Fix a compact generator $G$ for $D_{\operatorname{qc},Z}(X)$.
	For any $E\in D_{\operatorname{qc},Z}(X)$, the following are equivalent:
	\begin{enumerate}
		\item \label{lem:boundedness_via_coaisle1} $E\in D^+_{\operatorname{qc},Z}(X)$
		\item \label{lem:boundedness_via_coaisle2} there exists $i\geq 0$ such that $\operatorname{Hom}(G[t],E)\cong 0$ if $t\geq i$.
	\end{enumerate}
	A similar argument holds for objects in $D^-_{\operatorname{qc}}(X)$ where we take $[-t]$.
\end{lemma}

\begin{proof}
	We prove the case for $Z=|X|$.
	The general case can be argued analogously.
	It is possible to argue directly from \cite[\href{https://stacks.math.columbia.edu/tag/0GFH}{Tag 0GFH} and \href{https://stacks.math.columbia.edu/tag/0GFI}{0GFI}]{StacksProject} but we provide an alternate proof.

	Suppose $E\in D^+_{\operatorname{qc}}(X)$.
	Then there exists an $n\geq 0$ such that $E\in D^{\geq -n}_{\operatorname{qc}}(X)$.
	Hence, $E[-n]\in D^{\geq 0}_{\operatorname{qc}}(X)= D^{\geq -n}_{\operatorname{qc}}(X)[-n]$.
	%%NOTE: D^{\geq -n}_{\operatorname{qc}}(X) = D^{\geq 0}_{\operatorname{qc}}(X) [n]
	It follows that $\operatorname{Hom}(A,E[-n-1])=0$ for all $A\in D^{\leq 0}_{\operatorname{qc}}(X)$ because $D^{\geq 0}_{\operatorname{qc}}(X)[-1] = (D^{\leq 0}_{\operatorname{qc}}(X))^\perp$.
	Since $D^{\leq 0}_{\operatorname{qc}}(X)$ lies in the preferred equivalence class, there exists an $N\geq 0$ such that $\overline{\langle G \rangle}^{(-\infty,0]}[N]\subseteq D^{\leq 0}_{\operatorname{qc}}(X)$.
	%%NOTE: \overline{\langle G \rangle}^{(-\infty,0]}[N] is $\leq - N$ of cocomplete $t$-structure generated by $G$.
	Thus, $\operatorname{Hom}(G[N],E[-n-1])=0$, and hence, $\operatorname{Hom}(G[N+a],E[-n-1])=0$ for all $a\geq 0$ because $D^{\leq 0}_{\operatorname{qc}}(X)$ is closed under positive shifts.
	As $n\geq 0$, it follows $\operatorname{Hom}(G[N+a+n+1],E)=0$ for all $a\geq 0$.
	In other words, $\operatorname{Hom}(G[i],E) \cong 0$ for all $i \gg 0$, which shows that $\eqref{lem:boundedness_via_coaisle1} \implies \eqref{lem:boundedness_via_coaisle2}$.

	Conversely, suppose that $\operatorname{Hom}(G[i],E)\cong 0$ for all $i \gg 0$.
	Then there exists a smallest $n\geq 0$ such that this vanishing occurs. Set $n$ to be this integer.
	By \Cref{lem:preferred_eq_class}, there exists $N\geq 0$ satisfying $D^{\leq -N}_{\operatorname{qc}}(X) = D^{\leq 0}_{\operatorname{qc}}(X)[N] \subseteq \overline{\langle G[n] \rangle}^{(-\infty,0]}$.
	Hence, for all $A\in D^{\leq 0}_{\operatorname{qc}}(X)[N]$, we have $\operatorname{Hom}(A,E)=0$.
	It follows that $E\in (D^{\leq 0}_{\operatorname{qc}}(X)[N])^\perp = D^{\geq -N}_{\operatorname{qc}}(X)[-1]$.
	Since $N\geq 0$, we obtain that $D^{\leq -N - a}_{\operatorname{qc}}(X)\subseteq D^{\leq -N}_{\operatorname{qc}}(X)$ for all $a\geq 0$.
	%%NOTE: Use aisles closed under shifting by positive shifts
	If necessary, choose $N > 1$.
	This implies that $E \in D^{\geq - N}_{\operatorname{qc}}(X)$ for some $N \gg 1$.
	Consequently, $E\in D^+_{\operatorname{qc}}(X)$.

	The last claim is argued similarly the argument above.
\end{proof}

\begin{proposition}
	\label{prop:stacky_all_coincide}
	Let $f\colon Y\to X$ be a proper morphism of Noetherian algebraic spaces.
	Choose $E\in D^-_{\operatorname{coh}}(Y)$.
	Then $E$ is $f$-perfect if, and only if, it is $f$-quasi-perfect.
	In such a case, $E$ has bounded cohomology.
\end{proposition}

\begin{proof}
	The last claim follows from the fact $E$ is pseudocoherent and $f$-perfect. Moreover, \Cref{lem:f_perf_implies_perf_preserved} shows $f$-perfectness implies $f$-quasi-perfectness.
	We check the converse.
	Choose an \'{e}tale presentation $s\colon U \to X$ from an affine scheme.
	Consider the fibered square
	\begin{displaymath}
		% https://q.uiver.app/#q=WzAsNCxbMSwwLCJVIl0sWzEsMSwiXFxtYXRoY2Fse1h9LiJdLFswLDEsIlxcbWF0aGNhbHtZfSJdLFswLDAsIlxcbWF0aGNhbHtZfVxcdGltZXNfe1xcbWF0aGNhbHtYfX0gVSJdLFswLDEsInMiXSxbMiwxLCJmIiwyXSxbMywwLCJmXlxccHJpbWUiXSxbMywyLCJzXlxccHJpbWUiLDJdXQ==
		\begin{tikzcd}
			{Y\times_{X} U} & U \\
			{Y} & {X.}
			\arrow["{f^\prime}", from=1-1, to=1-2]
			\arrow["{s^\prime}"', from=1-1, to=2-1]
			\arrow["s", from=1-2, to=2-2]
			\arrow["f"', from=2-1, to=2-2]
		\end{tikzcd}
	\end{displaymath}
	Fix $A\in D^b_{\operatorname{qc}}(X)$.
	To show $E\otimes^{\mathbf{L}} \mathbf{L}f^\ast A \in D^b_{\operatorname{qc}}(Y)$, it suffices to check that
	\begin{displaymath}
		\mathbf{L}(s^\prime)^\ast (E\otimes^{\mathbf{L}} \mathbf{L} f^\ast A)
		\cong \mathbf{L}(s^\prime)^\ast E\otimes^{\mathbf{L}} \mathbf{L}(s^\prime)^\ast \mathbf{L} f^\ast A \in D^b_{\operatorname{qc}}(Y\times_{X} U)
	\end{displaymath}
	By \Cref{lem:quasi-affine_via_diagonal}, $s$ is quasi-affine.
	Hence, $s^\prime$ must be quasi-affine.
	Then \Cref{lem:pullback_compact_generator} shows $\mathbf{L}(s^\prime)^\ast \operatorname{Perf}(Y)$ compactly generates $D_{\operatorname{qc}} (Y\times_{X} U)$.
	Applying flat base change, it follows that
	\begin{displaymath}
		\begin{aligned}
			\mathbf{R}f^\prime_\ast (\mathbf{L}(s^\prime)^\ast E \otimes^{\mathbf{L}} \operatorname{Perf}(Y\times_{X} U))
			&\subseteq \mathbf{R}f^\prime_\ast (\mathbf{L}(s^\prime)^\ast E \otimes^{\mathbf{L}} \langle \mathbf{L}(s^\prime)^\ast \operatorname{Perf}(Y) \rangle )
			\\&\subseteq \langle \mathbf{R}f^\prime_\ast (\mathbf{L}(s^\prime)^\ast E \otimes^{\mathbf{L}}  \mathbf{L}(s^\prime)^\ast \operatorname{Perf}(Y) ) \rangle
			\\&\subseteq \langle \mathbf{R}f^\prime_\ast \mathbf{L}(s^\prime)^\ast ( E \otimes^{\mathbf{L}}  \operatorname{Perf}(Y) ) \rangle
			\\&\subseteq \langle \mathbf{L}s^\ast \mathbf{R}f_\ast ( E \otimes^{\mathbf{L}}  \operatorname{Perf}(Y) ) \rangle.
		\end{aligned}
	\end{displaymath}
	As $E$ is $f$-quasi-perfect, $\mathbf{R}f_\ast (E\otimes^{\mathbf{L}} \operatorname{Perf}(Y))\subseteq \operatorname{Perf}(X)$.
	Hence, we obtain
	\begin{displaymath}
		\langle \mathbf{L}s^\ast \mathbf{R}f_\ast ( E \otimes^{\mathbf{L}}  \operatorname{Perf}(Y) ) \rangle\subseteq \operatorname{Perf}(U),
	\end{displaymath}
	and so $\mathbf{L}(s^\prime)^\ast E$ is $f^\prime$-quasi-perfect (use \cite[\href{https://stacks.math.columbia.edu/tag/08H4}{Tag 08H4}]{StacksProject} for pseudocoherence).

	Let $G$ be a compact generator for $D_{\operatorname{qc}}(Y\times_{X} U)$.
	By \Cref{thm:adjoint_dqc_for_f-quasi-perfect},
	we know that
	\begin{displaymath}
		\begin{aligned}
			\operatorname{Hom} & (G , \mathbf{L}(s^\prime)^\ast E \otimes^{\mathbf{L}} \mathbf{L}(f^\prime)^\ast \mathbf{L}s^\ast A)
			\\&\cong \operatorname{Hom}(\mathbf{R}f^\prime_\ast ( G \otimes^{\mathbf{L}} \operatorname{\mathbf{R}\mathcal{H}\! \mathit{om}}(\mathbf{L}(s^\prime)^\ast E, (f^\prime)^\times \mathcal{O}_U) ) ,  \mathbf{L}(s^\ast A)).
		\end{aligned}
	\end{displaymath}
	However, \Cref{lem:Ballard_relative_perf_implies_dual_is_such} says
	\begin{displaymath}
		\mathbf{R}f^\prime_\ast ( G \otimes^{\mathbf{L}} \operatorname{\mathbf{R}\mathcal{H}\! \mathit{om}}(\mathbf{L}(s^\prime)^\ast E, (f^\prime)^\times \mathcal{O}_U) ) \subseteq \operatorname{Perf}(U).
	\end{displaymath}
	Since $U$ is affine, $D_{\operatorname{qc}}(U)$ is compactly generated by
	\begin{displaymath}
		\mathbf{R}f^\prime_\ast ( G \otimes^{\mathbf{L}} \operatorname{\mathbf{R}\mathcal{H}\! \mathit{om}}(\mathbf{L}(s^\prime)^\ast E, (f^\prime)^\times \mathcal{O}_U) ) \oplus \mathcal{O}_U.
	\end{displaymath}
	By \Cref{lem:neeman187}, $b^! \cong b^\times$ for all proper morphisms $b$ of Noetherian algebraic spaces.
	We apply this below.
	Now, $\mathbf{L}s^\ast A \in D^b_{\operatorname{qc}}(U)$, and so \Cref{lem:boundedness_via_coaisle} implies
	\begin{displaymath}
		\operatorname{Hom}((\mathbf{R}f^\prime_\ast ( G \otimes^{\mathbf{L}} \operatorname{\mathbf{R}\mathcal{H}\! \mathit{om}}(\mathbf{L}(s^\prime)^\ast E , (f^\prime)^! \mathcal{O}_U) ) \oplus \mathcal{O}_U) [i],  \mathbf{L}(s^\ast A)) \cong 0
	\end{displaymath}
	for all $i\gg 0$. Hence,
	\begin{displaymath}
		\operatorname{Hom}(\mathbf{R}f^\prime_\ast ( G \otimes^{\mathbf{L}} \operatorname{\mathbf{R}\mathcal{H}\! \mathit{om}}(\mathbf{L}(s^\prime)^\ast E, (f^\prime)^! \mathcal{O}_U) ) [i],  \mathbf{L}(s^\ast A)) \cong 0.
	\end{displaymath}
	By adjunction, we have
	\begin{displaymath}
		\operatorname{Hom} (G [i] , \mathbf{L}(s^\prime)^\ast E \otimes^{\mathbf{L}} \mathbf{L}(f^\prime)^\ast \mathbf{L}s^\ast A) \cong 0
	\end{displaymath}
	for all $i\gg 0$. Then \Cref{lem:boundedness_via_coaisle} shows that
	%%NOTE: Use that this is in the preferred equivalence class for the affine scheme $U$.
	$\mathbf{L}(s^\prime)^\ast E \otimes^{\mathbf{L}} \mathbf{L}(f^\prime)^\ast \mathbf{L}s^\ast A\in D^+_{\operatorname{qc}}(Y\times_{X} U)$. A similar argument shows $\mathbf{L}(s^\prime)^\ast E \otimes^{\mathbf{L}} \mathbf{L}(f^\prime)^\ast \mathbf{L}s^\ast A\in D^-_{\operatorname{qc}}(Y\times_{X} U)$. Therefore, $\mathbf{L}(s^\prime)^\ast E \otimes^{\mathbf{L}} \mathbf{L}(f^\prime)^\ast \mathbf{L}s^\ast A\in D^b_{\operatorname{qc}}(Y\times_{X} U)$, which completes the proof.
\end{proof}

\begin{remark}
	Let $f_i\colon Y_i\to S$ be proper morphisms of Noetherian algebraic spaces for $i=1,2$.
	Let $K\in D_{\operatorname{coh}}^-(Y_1\times_SY_2)$.
	If the Fourier--Mukai transform $\Phi_K\colon D_{\operatorname{qc}}(Y_1)\to D_{\operatorname{qc}}(Y_2)$ is an equivalence,
	then $K\in D_{\operatorname{coh}}^b(Y_1\times_SY_2)$.
	Indeed, an equivalence preserves compact objects, which on $Y_i$ coincide with the perfect complexes \cite[\href{https://stacks.math.columbia.edu/tag/09M8}{Tag 09M8}]{StacksProject}.
	Now combine \Cref{lem:preservation_of_perfect} and \Cref{prop:stacky_all_coincide}.
\end{remark}

\begin{question}
	Let $f_i\colon Y_i\to S$ be proper morphisms of Noetherian algebraic spaces for $i=1,2$.
	If $K\in D_{\operatorname{qc}}(Y_1\times_SY_2)$ induces an equivalence $\Phi_K$, does it follow $K\in D_{\operatorname{coh}}^b(Y_1\times_SY_2)$?
\end{question}

\begin{lemma}
	\label{lem:perfectness_via_coherence}
	Let $X$ be a Noetherian scheme.
	Then $E\in D^b_{\operatorname{coh}}(X)$ is perfect if, and only if, $E\otimes^{\mathbf{L}}\mathbf{R}i_\ast \mathcal{O}_{\operatorname{Spec}(\kappa(p))}\in D^b_{\operatorname{coh}}(X)$ for all morphisms $i\colon \operatorname{Spec}(\kappa(p)) \to X$ associated to a closed point $p\in X$.
\end{lemma}

\begin{proof}
	It suffices to prove the claim on the small Zariski site \cite[\href{https://stacks.math.columbia.edu/tag/071Q}{Tag 071Q}]{StacksProject}.
	By \Cref{prop:perfectness}, if $E$ is perfect, then $E\otimes^{\mathbf{L}}\mathbf{R}i_\ast \mathcal{O}_{\operatorname{Spec}(\kappa(p))}\in D^b_{\operatorname{coh}}(X)$ for all morphisms $i\colon \operatorname{Spec}(\kappa(p)) \to X$ associated to a closed point $p\in X$.
	We check the converse.
	Choose a closed point $p\in X$.
	Denote by $i\colon \operatorname{Spec}(\kappa(p)) \to X$ the associated closed immersion.
	The morphism $i$ factors via canonical morphisms,
	\begin{displaymath}
		% https://q.uiver.app/#q=WzAsMyxbMCwwLCJcXG9wZXJhdG9ybmFtZXtTcGVjfShcXGthcHBhKHApKSJdLFsxLDAsIlxcb3BlcmF0b3JuYW1le1NwZWN9KFxcbWF0aGNhbHtPfV97WCxwfSkiXSxbMiwwLCJYLiJdLFswLDEsImleXFxwcmltZSJdLFsxLDIsInMiXV0=
		\begin{tikzcd}
			{\operatorname{Spec}(\kappa(p))} & {\operatorname{Spec}(\mathcal{O}_{X,p})} & {X.}
			\arrow["{i^\prime}", from=1-1, to=1-2]
			\arrow["s", from=1-2, to=1-3]
		\end{tikzcd}
	\end{displaymath}
	It follows that
	\begin{displaymath}
		\mathbf{L}s^\ast (E\otimes^{\mathbf{L}}\mathbf{R}i_\ast \mathcal{O}_{\operatorname{Spec}(\kappa(p))}) \cong
		E_p \otimes^{\mathbf{L}}\mathbf{R}i^\prime_\ast \mathcal{O}_{\operatorname{Spec}(\kappa(p))}\in D^b_{\operatorname{coh}}(\mathcal{O}_{X,p}).
	\end{displaymath}
	Hence, $E_p$ is perfect.
	Since $p\in X$ was an arbitrary closed point, this implies $E$ is perfect (see e.g.\ \cite[\href{https://stacks.math.columbia.edu/tag/068W}{Tag 068W}]{StacksProject}).
	%%NOTE: That is, stalks at closed points is perfect.
\end{proof}

\begin{remark}
	\Cref{lem:perfectness_via_coherence} implicitly appears in the proof of \cite[Lemma 3.4]{Dutta/Lank/ManaliRahul:2025} when loc.\ cit.\ references \cite[Theorem 2.3(3)]{AlonsoTarrio/JeremiasLopez/SanchodeSalas:2023}.
\end{remark}

\begin{lemma}
	\label{lem:perfect_locus_open}
	Let $X$ be a Noetherian algebraic space.
	For any $E\in D^b_{\operatorname{coh}}(X)$, the collection of $p\in |X|$ such that $\mathbf{L}s^\ast E$ is perfect, where $s\colon \operatorname{Spec}(k)\to X$ is a representative of $p$, is open in $|X|$.
\end{lemma}

\begin{proof}
	Choose an \'{e}tale presentation $\pi\colon U \to X$ from an affine scheme.
	Let $S$ be the collection of $p\in |X|$ such that $\mathbf{L}j_p^\ast E$ is perfect where $j_p\colon \operatorname{Spec}(k)\to X$ represents $p$.
	Set $S^\prime$ to be the collection of $u\in U$ such that $(\mathbf{L}\pi^\ast E)_u\in \langle \mathcal{O}_{U,u} \rangle$.
	Applying \Cref{lem:pullback_is_perfect_iff_for_some_or_all_in_equivalence_class}, we see that this condition is independent of the choice of representative.

	Let $q\in U$ satisfy $\pi(q)\in S$.
	Consider the canonical morphisms
	\begin{displaymath}
		\operatorname{Spec}(\kappa(q))\xrightarrow{i_q} \operatorname{Spec}(\mathcal{O}_{U,q}) \xrightarrow{t_q} U.
	\end{displaymath}
	By hypothesis, $\mathbf{L}(\pi \circ t_q \circ i_q)^\ast E$ is perfect.
	This implies $\mathbf{L}(\pi \circ t_q)^\ast E\in \operatorname{Perf}(\mathcal{O}_{U,q})$ (see e.g.\ \cite[\href{https://stacks.math.columbia.edu/tag/068W}{Tag 068W}]{StacksProject}).
	Hence, $\pi^{-1}(S)\subseteq S^\prime$.

	Conversely, let $u\in S^\prime$.
	Then $\mathbf{L}(\pi \circ t_u)^\ast E\in \operatorname{Perf}(\mathcal{O}_{U,u})$, and so $\mathbf{L}(\pi \circ t_u \circ i_u)^\ast E$ is perfect.
	Thus, $\pi^{-1}(S) = S^\prime$.

	Since $\pi$ is surjective, it follows that $\pi(S^\prime)=S$.
	By \cite[Proposition 2.5]{Letz:2021}, $S^\prime$ is Zariski open in $U$ (which can be checked on the small Zariski site \cite[\href{https://stacks.math.columbia.edu/tag/071Q}{Tag 071Q}]{StacksProject}).
	As $\pi$ is a smooth morphism, it is an open morphism, and so $S$ must be open.
\end{proof}

\begin{proof}
	[Proof of \Cref{thm:stacky_all_coincide}]
	By \Cref{prop:stacky_all_coincide}, we have $\eqref{cor:perfectness_relative_via_coherence1} \iff \eqref{cor:perfectness_relative_via_coherence3}$.
	It is immediate that $\eqref{cor:perfectness_relative_via_coherence1} \implies \eqref{cor:perfectness_relative_via_coherence2}$.
	We prove $\eqref{cor:perfectness_relative_via_coherence2} \implies \eqref{cor:perfectness_relative_via_coherence3}$.
	Fix a closed point $p\in |X|$.
	Let $i\colon Z_p \to X$ be the residual space at $p$ (see \cite[\href{https://stacks.math.columbia.edu/tag/06QZ}{Tag 06QZ} and \href{https://stacks.math.columbia.edu/tag/06R0}{06R0}]{StacksProject}).
	By \cite[\href{https://stacks.math.columbia.edu/tag/0H1R}{Tag 0H1R}]{StacksProject}, $Z_p$ can be represented by $\operatorname{Spec}(k)$ for some field $k$, i.e.\ we can choose $Z_p = \operatorname{Spec}(k)$.
	Since $p$ is closed, $i$ is a closed immersion \cite[\href{https://stacks.math.columbia.edu/tag/0H1U}{Tag 0H1U}]{StacksProject}.
	In particular, $\mathbf{R}i_\ast\mathcal{O}_{\operatorname{Spec}(k)} \in D^b_{\operatorname{coh}}(X)$.

	Fix any $Q\in \operatorname{Perf}(Y)$.
	Since $E$ is pseudocoherent, $\eqref{cor:perfectness_relative_via_coherence2}$ implies $E$ has bounded cohomology.
	Hence, $\mathbf{R}f_\ast (E\otimes^{\mathbf{L}} Q)\in D^b_{\operatorname{coh}}(X)$.
	By $\eqref{cor:perfectness_relative_via_coherence2}$, it follows that $E\otimes^{\mathbf{L}} \mathbf{L}f^\ast \mathbf{R}i_\ast \mathcal{O}_{\operatorname{Spec}(k)}\in D^b_{\operatorname{coh}}(Y)$.
	Since $Q$ is perfect, we have that
	\begin{displaymath}
		Q\otimes^{\mathbf{L}} E\otimes^{\mathbf{L}}\mathbf{L}f^\ast \mathbf{R}i_\ast \mathcal{O}_{\operatorname{Spec}(k)}\in D^b_{\operatorname{coh}}(Y).
	\end{displaymath}
	Moreover, $f$ being proper implies
	\begin{displaymath}
		\mathbf{R}f_\ast (Q\otimes^{\mathbf{L}} E\otimes^{\mathbf{L}} \mathbf{L}f^\ast \mathbf{R}i_\ast \mathcal{O}_{\operatorname{Spec}(k)})\in D^b_{\operatorname{coh}}(X).
	\end{displaymath}
	By projection formula, we obtain
	\begin{displaymath}
		\mathbf{R}f_\ast (Q\otimes^{\mathbf{L}} E) \otimes^{\mathbf{L}} \mathbf{R}i_\ast \mathcal{O}_{\operatorname{Spec}(k)} \in D^b_{\operatorname{coh}}(X).
	\end{displaymath}
	A further application of projection tells us that
	\begin{displaymath}
		\mathbf{R}f_\ast (Q\otimes^{\mathbf{L}} E) \otimes^{\mathbf{L}} \mathbf{R}i_\ast \mathcal{O}_{\operatorname{Spec}(k)}
		\cong \mathbf{R}i_\ast \mathbf{L}i^\ast \mathbf{R}f_\ast (Q\otimes^{\mathbf{L}} E).
	\end{displaymath}
	Then \Cref{lem:34} implies $\mathbf{L}i^\ast \mathbf{R}f_\ast (Q\otimes^{\mathbf{L}} E) \in D^b_{\operatorname{coh}}(\operatorname{Spec}(k))$.
	By \Cref{lem:perfect_locus_open}, the collection $B$ of $p\in |X|$ such that
	\begin{displaymath}
		\mathbf{L}s^\ast (\mathbf{R}f_\ast (Q\otimes^{\mathbf{L}} E))
	\end{displaymath}
	is perfect, where $s\colon \operatorname{Spec}(k)\to X$ is a representative of $p$, is open in $|X|$.
	However, the work above shows every closed point belongs to this collection.
	By \Cref{rmk:topological_fact_for_stacks}, every point in $|X|$ is a generalization of a closed point, and so $B = |X|$.
	Therefore, \Cref{prop:perfectness} implies $\mathbf{R}f_\ast (Q\otimes^{\mathbf{L}} E) \in \operatorname{Perf}(X)$.
\end{proof}

\begin{proposition}
	\label{prop:preservation_of_dbcoh}
	Let $S$ be a Noetherian algebraic space.
	Consider proper morphisms of algebraic spaces $f_1\colon Y_1 \to S$ and $f_2\colon Y_2 \to S$.
	Denote by $p_i\colon Y_1 \times_{S} Y_2 \to Y_i$ the natural projections.
	For any $E\in D_{\operatorname{qc}}(Y_1 \times_{S} Y_2)$ pseudocoherent, the following are equivalent:
	\begin{enumerate}
		\item $E$ is $p_1$-perfect
		\item $\Phi_E (D^b_{\operatorname{coh}}(Y_1))\subseteq D^b_{\operatorname{coh}}(Y_2)$.
	\end{enumerate}
\end{proposition}

\begin{proof}
	First, assume $E$ is $p_1$-perfect.
	Base change says each $p_i$ is proper.
	It follows that
	\begin{displaymath}
		E\otimes^{\mathbf{L}} \mathbf{L}p_1^\ast D^b_{\operatorname{coh}}(Y_1)\subseteq D^b_{\operatorname{coh}}(Y_1 \times_{S} Y_2).
	\end{displaymath}
	Since $p_2$ is proper, we have
	\begin{displaymath}
		\mathbf{R}(p_2)_\ast (E\otimes^{\mathbf{L}} \mathbf{L}p_1^\ast D^b_{\operatorname{coh}}(Y_1))\subseteq D^b_{\operatorname{coh}}(Y_2).
	\end{displaymath}
	See  \cite[\href{https://stacks.math.columbia.edu/tag/08GK}{Tag 08GK}]{StacksProject}. We prove the converse.

	Recall that each $D_{\operatorname{qc}}(Y_i)$ is compactly generated by an object $G_i$.
	Moreover, from \Cref{lem:neeman2023cor5_10}, $\mathbf{L}p_1^\ast G_1 \otimes^{\mathbf{L}} \mathbf{L}p_2^\ast G_2$ is a compact generator for $D_{\operatorname{qc}}(Y_1 \times_{S} Y_2)$.
	Choose $B\in D^b_{\operatorname{qc}}(Y_1)$.
	We want to show that $E\otimes^{\mathbf{L}} \mathbf{L}p_1^\ast B\in D^b_{\operatorname{qc}}(Y_1 \times_{S} Y_2)$.
	Applying \Cref{thm:stacky_all_coincide}, we can test for $p_1$-perfectness by restricting to $B\in D^b_{\operatorname{coh}}(Y_1)$.
	In this case, as $E$ and $\mathbf{L}p_1^\ast B$ are pseudocoherent, it follows that $E\otimes^{\mathbf{L}} \mathbf{L}p_1^\ast B$ is pseudocoherent, and so $E\otimes^{\mathbf{L}} \mathbf{L}p_1^\ast B\in D^-_{\operatorname{coh}}(Y_1 \times_{S} Y_2)$.
	To show that $E\otimes^{\mathbf{L}} \mathbf{L}p_1^\ast B\in D^b_{\operatorname{coh}}(Y_1 \times_{S} Y_2)$, it suffices to check that $E\otimes^{\mathbf{L}} \mathbf{L}p_1^\ast B\in D^+_{\operatorname{qc}}(Y_1 \times_{S} Y_2)$.
	By adjunctions,
	\begin{displaymath}
		\begin{aligned}
			&\operatorname{Ext}^n  (\mathbf{L}p_1^\ast G_1 \otimes^{\mathbf{L}} \mathbf{L}p_2^\ast G_2, E\otimes^\mathbf{L} \mathbf{L}p_1^\ast B)
			\\&\cong \operatorname{Ext}^n (\mathbf{L}p_2^\ast G_2, \operatorname{\mathbf{R}\mathcal{H}\! \mathit{om}}(\mathbf{L}p_1^\ast G_1, E\otimes^\mathbf{L} \mathbf{L}p_1^\ast B))  && (\textrm{tensor/hom})
			\\&\cong \operatorname{Ext}^n (\mathbf{L}p_2^\ast G_2, \operatorname{\mathbf{R}\mathcal{H}\! \mathit{om}}(\mathbf{L}p_1^\ast G_1,\mathcal{O}_{Y_1\times_{S} Y_2}) \otimes^{\mathbf{L}} E\otimes^\mathbf{L} \mathbf{L}p_1^\ast B) && \textrm{(\cite[\href{https://stacks.math.columbia.edu/tag/08JJ}{Tag 08JJ}]{StacksProject})}
			\\&\cong \operatorname{Ext}^n (\mathbf{L}p_2^\ast G_2, \operatorname{\mathbf{R}\mathcal{H}\! \mathit{om}}(\mathbf{L}p_1^\ast G_1,\mathbf{L}p_1^\ast \mathcal{O}_{Y_1}) \otimes^{\mathbf{L}} E\otimes^\mathbf{L} \mathbf{L}p_1^\ast B ) && (\mathbf{L}p_1^\ast \mathcal{O}_{Y_1} =\mathcal{O}_{Y_1\times_{S} Y_2})
			\\&\cong \operatorname{Ext}^n (\mathbf{L}p_2^\ast G_2, \mathbf{L}p_1^\ast (\operatorname{\mathbf{R}\mathcal{H}\! \mathit{om}}( G_1,\mathcal{O}_{Y_1})) \otimes^{\mathbf{L}} E\otimes^\mathbf{L} \mathbf{L}p_1^\ast B ) && \textrm{(\Cref{lem:gortz_wedhorn_internal_hom})}
			\\&\cong \operatorname{Ext}^n (\mathbf{L}p_2^\ast G_2, \mathbf{L}p_1^\ast (\operatorname{\mathbf{R}\mathcal{H}\! \mathit{om}}( G_1,\mathcal{O}_{Y_1})\otimes^{\mathbf{L}} B) \otimes^{\mathbf{L}} E ) && \textrm{(\cite[\href{https://stacks.math.columbia.edu/tag/07A4}{Tag 07A4}]{StacksProject})}
			\\&\cong \operatorname{Ext}^n (G_2, \mathbf{R}(p_2)_\ast (\mathbf{L}p_1^\ast (\operatorname{\mathbf{R}\mathcal{H}\! \mathit{om}}( G_1,\mathcal{O}_{Y_1})\otimes^{\mathbf{L}} B) \otimes^{\mathbf{L}} E) ) && \textrm{(pull/push)}
			\\&\cong \operatorname{Ext}^n (G_2, \Phi_E (\operatorname{\mathbf{R}\mathcal{H}\! \mathit{om}}( G_1,\mathcal{O}_{Y_1})\otimes^{\mathbf{L}} B)).
		\end{aligned}
	\end{displaymath}
	Clearly, $\operatorname{\mathbf{R}\mathcal{H}\! \mathit{om}}( G_1,\mathcal{O}_{Y_1})$ is perfect, and so
	\begin{displaymath}
		\operatorname{\mathbf{R}\mathcal{H}\! \mathit{om}}( G_1,\mathcal{O}_{Y_1})\otimes^{\mathbf{L}} B\in D^b_{\operatorname{coh}}(Y_1).
	\end{displaymath}
	Since
	\begin{displaymath}
		\Phi_E (\operatorname{\mathbf{R}\mathcal{H}\! \mathit{om}}( G_1,\mathcal{O}_{Y_1})\otimes^{\mathbf{L}} B)\in D^b_{\operatorname{coh}}(Y_2),
	\end{displaymath}
	we know that
	\begin{displaymath}
		\Phi_E (\operatorname{\mathbf{R}\mathcal{H}\! \mathit{om}}( G_1,\mathcal{O}_{Y_1})\otimes^{\mathbf{L}} B)\in D^+_{\operatorname{qc}}(Y_2).
	\end{displaymath}
	Hence, by \Cref{lem:boundedness_via_coaisle}, there exists an $i\gg 0$ such that for all $t\geq i$, one has
	\begin{displaymath}
		\operatorname{Hom}(G_2 [t], \Phi_E (\operatorname{\mathbf{R}\mathcal{H}\! \mathit{om}}( G_1,\mathcal{O}_{Y_1})\otimes^{\mathbf{L}} B)) \cong 0.
	\end{displaymath}
	Then the chain of isomorphisms above imply the claim via \Cref{lem:boundedness_via_coaisle}.
\end{proof}

\begin{proof}
	[Proof of \Cref{cor:preservation}]
	This follows from \Cref{lem:preservation_of_perfect}, \Cref{prop:preservation_of_dbcoh}, and \Cref{thm:stacky_all_coincide}.
\end{proof}

\begin{remark}
	\label{rmk:fm_cube_pullback}
	Consider \Cref{setup:fm_cube_pullback} with $t$ a quasi-affine morphism. By base change, $t^\prime$ and each $t_i$ are quasi-affine.
	Moreover, base change implies each algebraic space in the cube is Noetherian.
	Note that $D_{\operatorname{qc}}(Y_1)$ and $D_{\operatorname{qc}}(Y_2)$ are respectively compactly generated by some $G_1$ and $G_2$.
	Hence, \Cref{lem:neeman2023cor5_10} shows that $D_{\operatorname{qc}}(Y_1\times_{S} Y_2)$ is compactly generated by $\mathbf{L}(f^\prime_2)^\ast G_1 \otimes^{\mathbf{L}} \mathbf{L}(f^\prime_1)^\ast G_2$.
	Furthermore, from \Cref{lem:pullback_compact_generator}, we have that $D_{\operatorname{qc}}(Y_1\times_{S} T)$, $D_{\operatorname{qc}}(Y_2\times_{S} T)$, and $D_{\operatorname{qc}}(Y_1\times_{S} T\times_{S} Y_2)$ are respectively compactly generated by $\mathbf{L}t_1^\ast G_1$, $\mathbf{L}t_2^\ast G_2$, and $\mathbf{L}(t^\prime)^\ast (\mathbf{L}(f^\prime_2)^\ast G_1 \otimes^{\mathbf{L}} \mathbf{L}(f^\prime_1)^\ast G_2)$.
\end{remark}

\begin{lemma}
	\label{lem:fm_cube_locally_finite_tor_dimension}
	Consider \Cref{setup:fm_cube_pullback}.
	Then each face of \eqref{diag:fm_cube_pullback} is a tor-independent square.
\end{lemma}

\begin{proof}
	This follows from the fact that $f_i$, $f^\prime_i$, $g_i$, and $g^\prime_i$ are flat.
	This can be checked after taking \'{e}tale presentations, applying \cite[\href{https://stacks.math.columbia.edu/tag/08IQ}{Tag 08IQ}]{StacksProject} and reducing to the scheme case (see e.g.\ \cite[Lemma 22.93 \& Remark 22.94]{Gortz/Wedhorn:2020}).
	%%NOTE: In Tag 08IQ, apply (3) and use that base change along projections are flat at times.
\end{proof}

\begin{lemma}
	\label{lem:fm_cube_flat_base_change_natural_isomorphism}
	Consider \Cref{setup:fm_cube_pullback}.
	Let $K\in D_{\operatorname{qc}}(Y_1 \times_{S} Y_2)$.
	Then on $D_{\operatorname{qc}}$ there is a natural isomorphism
	\begin{displaymath}
		\alpha^K\colon \mathbf{L} t_2^\ast \circ \Phi_{K} \to \Phi_{\mathbf{L}(t^\prime)^\ast K} \circ \mathbf{L} t_1^\ast.
	\end{displaymath}
\end{lemma}

\begin{proof}
	There exists a string of natural isomorphisms for each $E\in D_{\operatorname{qc}}(Y_1)$:
	\begin{displaymath}
		\begin{aligned}
			\mathbf{L} t_2^\ast \Phi_{K} (E)
			&= \mathbf{L} t_2^\ast \mathbf{R}(f^\prime_1)_\ast (\mathbf{L} (f^\prime_2)^\ast E \otimes^{\mathbf{L}} K)
			\\&\cong \mathbf{R}(g^\prime_1)_\ast \mathbf{L}(t^\prime)^\ast (\mathbf{L} (f^\prime_2)^\ast E \otimes^{\mathbf{L}} K)
			&&\textrm{(flat base change)}
			\\&\cong \mathbf{R}(g^\prime_1)_\ast ( \mathbf{L}(t^\prime)^\ast  \mathbf{L} (f^\prime_2)^\ast E \otimes^{\mathbf{L}}  \mathbf{L}(t^\prime)^\ast K) && \textrm{(monoidality)}
			\\&\cong \mathbf{R}(g^\prime_1)_\ast ( \mathbf{L} (g^\prime_2)^\ast  \mathbf{L} (t_1)^\ast E \otimes^{\mathbf{L}}  \mathbf{L}(t^\prime)^\ast K) && \textrm{(pseudofunctoriality)}
			\\&= \Phi_{\mathbf{L}(t^\prime)^\ast K}( \mathbf{L} t_1^\ast E).
		\end{aligned}
	\end{displaymath}
	This completes the proof.
\end{proof}

\begin{lemma}
	\label{lem:pushforward_coherent_support_is_image}
	Let $f\colon Y \to X$ be a proper flat morphism of Noetherian algebraic spaces. Choose a closed subset $Z\subseteq |Y|$. Then $\mathbf{R}f_\ast D^b_{\operatorname{coh},Z}(Y) \subseteq D^b_{\operatorname{coh},f(Z)}(X)$.
\end{lemma}

\begin{proof}
	For any $E\in D^b_{\operatorname{coh},Z}(Y)$, define $W_E := \operatorname{Supp}(\mathbf{R}f_\ast E)$.
	Choose $p\in W_E$, and let $t\colon \operatorname{Spec}(k)\to X$ be a representative of $p$.
	Consider the fibered square
	\begin{displaymath}
		% https://q.uiver.app/#q=WzAsNCxbMSwwLCJcXG9wZXJhdG9ybmFtZXtTcGVjfShrKSJdLFsxLDEsIlguIl0sWzAsMSwiWSJdLFswLDAsIllcXHRpbWVzX1hcXG9wZXJhdG9ybmFtZXtTcGVjfShrKSJdLFswLDEsInQiXSxbMiwxLCJmIiwyXSxbMywyLCJ0XlxccHJpbWUiLDJdLFszLDAsImZeXFxwcmltZSJdXQ==
		\begin{tikzcd}
			{Y\times_X\operatorname{Spec}(k)} & {\operatorname{Spec}(k)} \\
			Y & {X.}
			\arrow["{f^\prime}", from=1-1, to=1-2]
			\arrow["{t^\prime}"', from=1-1, to=2-1]
			\arrow["t", from=1-2, to=2-2]
			\arrow["f"', from=2-1, to=2-2]
		\end{tikzcd}
	\end{displaymath}
	By \Cref{lem:support_is_cohomological_for_finite_type_cohomology}, $\mathbf{L}t^\ast \mathbf{R}f_\ast E \not\cong 0$.
	Applying flat base change, it follows that $\mathbf{R}f^\prime_\ast \mathbf{L}(t^\prime)^\ast E \not\cong 0$.
	This ensures that $Y\times_X\operatorname{Spec}(k)$ is nonempty (otherwise, $\mathbf{R}f^\prime_\ast \mathbf{L}(t^\prime)^\ast E \cong 0$), and so $f^\prime$ is surjective.
	Moreover, it follows that $\mathbf{L}(t^\prime)^\ast E \not\cong 0$.
	Choose any point $q\in \operatorname{Supp}(\mathbf{L}(t^\prime)^\ast E)$, and representative $h\colon \operatorname{Spec}(\ell) \to Y\times_X\operatorname{Spec}(k)$ of $q$.
	By \Cref{lem:support_is_cohomological_for_finite_type_cohomology}, $\mathbf{L}(t^\prime \circ h)^\ast E \not\cong 0$.
	Using \Cref{lem:factor_residual_gerbes}, $(t^\prime\circ h)$ represents $t^\prime (q)$.
	By \Cref{lem:support_is_cohomological_for_finite_type_cohomology}, $t^\prime (q) \in \operatorname{Supp}(E)\subseteq Z$.
	Therefore, $p = (f \circ t^\prime)(q) \in f(Z)$.
\end{proof}

\begin{lemma}
	\label{lem:reflecting_bounded_pseudocoherence_perfectness_trivial_case}
	Consider \Cref{setup:fm_cube_pullback}.
	Let $K\in D^b_{\operatorname{coh}}(Y_1 \times_{S} Y_2)$.
	For any $p\not \in (f_1 \circ f^\prime_2)(\operatorname{Supp}(K))$, and for any representative $t\colon \operatorname{Spec}(k)\to S$ of $p$, one has that $\mathbf{L}(t^\prime)^\ast K$ is relatively perfect over $Y_1 \times_{S} \operatorname{Spec}(k)$ (resp.\ $Y_2 \times_{S} \operatorname{Spec}(k)$).
	Moreover, $\Phi_{\mathbf{L}(t^\prime)^\ast K}(E)\cong 0$ for all $E\in D_{\operatorname{qc}}(Y_1\times_S \operatorname{Spec}(k))$.
\end{lemma}

\begin{proof}
	As $K$ is a bounded pseudocoherent complex, $\operatorname{Supp}(K)$ is closed.
	Then properness of $f_1 \circ f^\prime_2$ implies $(f_1 \circ f^\prime_2)(\operatorname{Supp}(K)) \subseteq |S|$ is closed.
	We claim that $\mathbf{L}(t^\prime)^\ast K\cong 0$ for all such $p$.
	Assume the contrary.
	Choose a representative $t\colon \operatorname{Spec}(k)\to S$ of $p\in |S|$ such that $\mathbf{L}(t^\prime)^\ast K \not= 0$.
	By \Cref{lem:support_is_cohomological_for_finite_type_cohomology}, there exists $q\in |Y_1 \times_S Y_2 \times_S \operatorname{Spec}(k)|$ and a representative $a\colon \operatorname{Spec}(k^\prime)\to Y_1 \times_S Y_2 \times_S \operatorname{Spec}(k)$ of $q$ such that $\mathbf{L}(t^\prime \circ a)^\ast K \not\cong 0$.
	Then \Cref{lem:factor_residual_gerbes,lem:support_is_cohomological_for_finite_type_cohomology} imply $t^\prime (q)\in \operatorname{Supp}(K)$.
	However,
	\begin{displaymath}
		(f_1 \circ f^\prime_2 \circ t^\prime \circ a)(q) = p \in (f_1 \circ f^\prime_2)(\operatorname{Supp}(K)),
	\end{displaymath}
	which is absurd.
	Therefore, $\mathbf{L}(t^\prime)^\ast K\cong 0$ for all $p\not\in (f_1 \circ f^\prime_2)(\operatorname{Supp}(K))$.
	In fact, the last claim becomes clear.
\end{proof}

\begin{lemma}
	\label{lem:induced_dbcoh_perf_preservation_upon_pullback_with_t_affine}
	Consider \Cref{setup:fm_cube_pullback} with $t$ a quasi-affine morphism.
	Let $K\in D^-_{\operatorname{coh}}(Y_1 \times_{S} Y_2)$.
	If $K$ is relatively perfect over $Y_1$ (resp.\ over $Y_2$), then $\Phi_{\mathbf{L} (t^\prime)^\ast K}$ restricts to an exact functor on $D^b_{\operatorname{coh}}$ (resp.\ on $\operatorname{Perf}$).
\end{lemma}

\begin{proof}
	Suppose $K$ is relatively perfect over $Y_2$. As $\Phi_K (\operatorname{Perf}(Y_1 ))\subseteq \operatorname{Perf}(Y_2 )$, \Cref{cor:preservation} implies
	\begin{displaymath}
		\mathbf{R}(f_1^\prime)_\ast (K \otimes^{\mathbf{L} }\operatorname{Perf}(Y_1 \times_{S} Y_2 )) \subseteq \operatorname{Perf}(Y_2 ).
	\end{displaymath}
	By base change, we know that $t_1,t_2,t^\prime$ are quasi-affine morphisms. Choose $G\in \operatorname{Perf}(Y_1 \times_{S} Y_2 )$ such that $\operatorname{Perf}(Y_1 \times_{S} Y_2 ) = \langle G \rangle$.
	Now, $\mathbf{L}(t^\prime)^\ast G$ satisfies
	\begin{displaymath}
		\operatorname{Perf}(Y_1 \times_{S} Y_2 \times_{S} T) = \langle \mathbf{L} (t^\prime)^\ast G \rangle,
	\end{displaymath}
	see \cite[\href{https://stacks.math.columbia.edu/tag/09SR}{Tag 09SR}]{StacksProject} and \Cref{rmk:fm_cube_pullback}.
	By base change, $g_1,g_2, f^\prime_2,f^\prime_1,g_2^\prime,g_1^\prime$ are proper flat morphisms.
	Using flat base change, we have that
	\begin{displaymath}
		\begin{aligned}
			\mathbf{R} (g^\prime_1)_\ast
			& (\mathbf{L}(t^\prime)^\ast K \otimes^\mathbf{L} \operatorname{Perf}(Y_1 \times_{S} Y_2  \times_{S} T))
			\\&\subseteq \mathbf{R} (g^\prime_1)_\ast \langle \mathbf{L}(t^\prime)^\ast (K \otimes^\mathbf{L} G) \rangle
			\\&\subseteq \langle \mathbf{R} (g^\prime_1)_\ast \mathbf{L}(t^\prime)^\ast (K \otimes^\mathbf{L} G) \rangle
			\\&\subseteq \langle \mathbf{L}t_2^\ast \mathbf{R}(f_1^\prime)_\ast (K \otimes^\mathbf{L} G) \rangle
			\\&\subseteq \operatorname{Perf}(Y_2 \times_{S} T).
		\end{aligned}
	\end{displaymath}
	By \Cref{lem:pseudocoherence_affine_faithfully_flat_cover}, $\mathbf{L}(t^\prime)^\ast K$ is pseudocoherent.
	Thus, \Cref{cor:preservation} says $\Phi_{\mathbf{L}(t^\prime)^\ast K}$ restricts to an exact functor on $\operatorname{Perf}$.
	A similar argument proves the case for $K$ being relatively perfect over $Y_1$.
\end{proof}

\begin{proposition}
	\label[proposition]{prop:reflecting_bounded_pseudocoherence_perfectness}
	Consider \Cref{setup:fm_cube_pullback}.
	Let $K\in D^b_{\operatorname{coh}}(Y_1 \times_{S} Y_2)$.
	\begin{enumerate}[label=(\arabic*), ref=\theremark(\arabic*)]
		\item \label[proposition]{prop:reflecting_bounded_pseudocoherence_perfectness1} If $\mathbf{L}(t^\prime)^\ast K$ is relatively perfect over $Y_1 \times_{S} T$ (resp.\ $Y_2 \times_{S} T$) where $t$ is faithfully flat, then $K$ is relatively perfect over $Y_1$ (resp.\ $Y_2$).
		\item \label[proposition]{prop:reflecting_bounded_pseudocoherence_perfectness2} If $\mathbf{L}(t^\prime)^\ast K$ is relatively perfect over $Y_1 \times_{S} \operatorname{Spec}(k)$ (resp.\ $Y_2 \times_{S} \operatorname{Spec}(k)$) for all morphisms $t\colon \operatorname{Spec}(k)\to S$ from a field, then $K$ is relatively perfect over $Y_1$ (resp.\ $Y_2$).
		\item \label[proposition]{prop:reflecting_bounded_pseudocoherence_perfectness3} If for every closed point $p\in |S|$ there exists a representative $t\colon \operatorname{Spec}(k)\to S$ of $p$ such that $\mathbf{L}(t^\prime)^\ast K$ is relatively perfect over $Y_1\times_{S} \operatorname{Spec}(k)$ (resp.\ $Y_2 \times_{S} \operatorname{Spec}(k)$), then $K$ is relatively perfect over $Y_1$ (resp.\ $Y_2$).
	\end{enumerate}
\end{proposition}

\begin{proof}
	We prove the case $\mathbf{L}(t^\prime)^\ast K$ is relatively perfect over $Y_1\times_{S} T$.
	The other case is argued similarly using \Cref{lem:pseudocoherence_affine_faithfully_flat_cover}.

	We check the first claim.
	Let $E\in D^b_{\operatorname{coh}}(Y_1)$.
	As $t^\prime$ is flat, it follows that $\mathbf{L} t_1^\ast E\in D^b_{\operatorname{coh}}(Y_1\times_{S} T)$.
	By \Cref{lem:fm_cube_flat_base_change_natural_isomorphism}, we know that $\mathbf{L} t_2^\ast \Phi_K (E)\cong \Phi_{\mathbf{L}(t^\prime)^\ast K} (\mathbf{L} t_1^\ast E)$.
	If $\mathbf{L}(t^\prime)^\ast K$ is relatively perfect over $Y_1 \times_{S} T$, then \Cref{cor:preservation} shows $\Phi_{\mathbf{L}(t^\prime)^\ast K} (\mathbf{L} t_1^\ast E) \in D^b_{\operatorname{coh}}(Y_2\times_{S} T)$.
	Since $t_2$ is faithfully flat, we know for each $j\in \mathbb{Z}$ that $\mathcal{H}^j (\mathbf{L} t_2^\ast \Phi_K (E))\cong 0$ implies $\mathcal{H}^j (\Phi_K (E))\cong 0$ (see \Cref{lem:zero_object_via_covers}).
	Consequently, we have $\Phi_K (D^b_{\operatorname{coh}}(Y_1)) \subseteq D^b_{\operatorname{coh}}(Y_2)$, which implies $K$ is relatively perfect over $Y_1$ by \Cref{cor:preservation}.

	Next, we check the second claim.
	Fix a compact generator $G\in D_{\operatorname{qc}}(Y_1\times_S Y_2)$.
	By \Cref{lem:pushforward_coherent_support_is_image},
	\begin{displaymath}
		\operatorname{Supp}(\mathbf{R}(f^\prime_2)_\ast (K \otimes^{\mathbf{L}} G)) \subseteq f^\prime_2 (\operatorname{Supp} (K \otimes^{\mathbf{L}} G)).
	\end{displaymath}
	To show $\mathbf{R}(f^\prime_2)_\ast (K \otimes^{\mathbf{L}} G)$ is perfect, \Cref{prop:perfectness} says we can test for boundedness of the derived pullback of $\mathbf{R}(f^\prime_2)_\ast (K \otimes^{\mathbf{L}} G)$ along representatives of points in $|Y_1|$.
	By \Cref{lem:support_is_cohomological_for_finite_type_cohomology}, it suffices to check this condition at all points in $\operatorname{Supp}(\mathbf{R}(f^\prime_2)_\ast (K \otimes^{\mathbf{L}} G))$ since the derived pullback is the zero object otherwise.
	Choose any $p\in f^\prime_2(\operatorname{Supp}(K\otimes^{\mathbf{L}} G))$, and let $h_1\colon \operatorname{Spec}(k)\to Y_1$ be a representative of $p$.
	Denote by $t\colon \operatorname{Spec}(k)\to S$ the composition $f_1\circ h_1$.
	There exists a commutative diagram
	\begin{displaymath}
		\begin{tikzcd}
			{\operatorname{Spec}(k)} && \\
			& {Y_1\times_{S} \operatorname{Spec}(k)} & {\operatorname{Spec}(k)} \\
			& {Y_1} & {S.}
			\arrow["h"{description}, from=1-1, to=2-2]
			\arrow["{1_{\operatorname{Spec}(k)}}", bend left=30, from=1-1, to=2-3]
			\arrow["{h_1}"', bend right=30, from=1-1, to=3-2]
			\arrow["{g_1}", from=2-2, to=2-3]
			\arrow["{t_1}"', from=2-2, to=3-2]
			\arrow["t", from=2-3, to=3-3]
			\arrow["{f_1}", from=3-2, to=3-3]
		\end{tikzcd}
	\end{displaymath}
	As $\mathbf{L}(t^\prime)^\ast K$ is relatively perfect over $Y_1 \times_S \operatorname{Spec}(k)$, it follows that
	\begin{displaymath}
		\begin{aligned}
			\mathbf{R}(g^\prime_2)_\ast
			& \mathbf{L}(t^\prime)^\ast ( K \otimes^{\mathbf{L}} G )
			\\&\cong \mathbf{R}(g^\prime_2)_\ast (\mathbf{L}(t^\prime)^\ast K \otimes^{\mathbf{L}} \mathbf{L}(t^\prime)^\ast G)
			\\&\in \mathbf{R}(g^\prime_2)_\ast (\mathbf{L}(t^\prime)^\ast K \otimes^{\mathbf{L}} \operatorname{Perf}(Y_1\times_{S} Y_2 \times_{S} \operatorname{Spec}(k))
			\\&\subseteq \operatorname{Perf}(Y_1\times_{S} \operatorname{Spec}(k)).
		\end{aligned}
	\end{displaymath}
	There is a fibered square
	\begin{displaymath}
		\begin{tikzcd}
			{Y_1\times_{S} Y_2\times_{S} \operatorname{Spec}(k)} & {Y_1\times_{S} \operatorname{Spec}(k)} \\
			{Y_1\times_{S} Y_2} & {Y_1}
			\arrow["{g_2^\prime}"', from=1-1, to=1-2]
			\arrow["{t^\prime}"', from=1-1, to=2-1]
			\arrow["{t_1}", from=1-2, to=2-2]
			\arrow["{f^\prime_2}"', from=2-1, to=2-2]
		\end{tikzcd}
	\end{displaymath}
	By flat base change,
	\begin{displaymath}
		\begin{aligned}
			\mathbf{L}t_1^\ast \mathbf{R}(f^\prime_2)_\ast ( K \otimes^{\mathbf{L}} G)
			&\cong \mathbf{R}(g^\prime_2)_\ast \mathbf{L}(t^\prime)^\ast ( K \otimes^{\mathbf{L}} G)
			\\&\in \operatorname{Perf}(Y_1\times_{S} \operatorname{Spec}(k)).
		\end{aligned}
	\end{displaymath}
	This implies that
	\begin{displaymath}
		\mathbf{L}h_1^\ast \mathbf{R}(f^\prime_2)_\ast ( K \otimes^{\mathbf{L}} G)
		\cong \mathbf{L}h^\ast \mathbf{L}t_1^\ast \mathbf{R}(f^\prime_2)_\ast ( K \otimes^{\mathbf{L}} G)
		\in D^b_{\operatorname{coh}}(k),
	\end{displaymath}
	because the derived pullback of perfect complexes remains perfect.
	% Using \Cref{prop:perfectness}, we deduce that $\mathbf{R}(f^\prime_2)_\ast ( K \otimes^{\mathbf{L}} P) \in \operatorname{Perf}(Y_1)$.
	% To see, note that $\mathbf{R}(f^\prime_2)_\ast ( K \otimes^{\mathbf{L}} P)$ is bounded and pseudocoherent because $K \otimes^{\mathbf{L}} P$ is a bounded pseudocoherent complex and $f^\prime_2$ is proper.
	% Indeed, we have shown that $\mathbf{L}b^\ast \mathbf{R}(f^\prime_2)_\ast ( K \otimes^{\mathbf{L}} P)$ is bounded for all $q\in Y_1$ with representative morphism $b\colon \operatorname{Spec}(\ell) \to Y_1$ that represents $q$.
	Thus, $K$ is relatively perfect over $Y_1$.
	Indeed, $\langle G \rangle = \operatorname{Perf}(Y_1\times_S Y_2)$, and so $\mathbf{R}(f^\prime_2)_\ast (K\otimes^{\mathbf{L}} \operatorname{Perf}(Y_1\times_S Y_2))\subseteq \langle \mathbf{R}(f^\prime_2)_\ast (K\otimes^{\mathbf{L}}G) \rangle$, whereas the work above shows the latter complex is perfect.

	Lastly, we check the third claim.
	Assume for every closed point $p\in |S|$ there exists a representative $t\colon \operatorname{Spec}(k)\to S$ of $p$ such that $\mathbf{L}(t^\prime)^\ast K$ is relatively perfect over $Y_1\times_{S} \operatorname{Spec}(k)$.
	Fix a compact generator $G\in D_{\operatorname{qc}}(Y_1\times_S Y_2)$.
	It suffices to check that $\mathbf{R}(f^\prime_2)_\ast (K\otimes^{\mathbf{L}} G)$ is perfect.
	By \Cref{lem:perfect_locus_open}, we only need to show that every closed point of $|Y_1|$ belongs to the collection of $p\in |Y_1|$ such that $\mathbf{L}s^\ast \mathbf{R}(f^\prime_2)_\ast (K\otimes^{\mathbf{L}} G)$ is perfect, where $s\colon \operatorname{Spec}(k)\to Y_1$ is a representative of $p$.
	Choose a closed point $p\in |Y_1|$.
	Let $h_1\colon \operatorname{Spec}(k)\to Y_1$ be a representative of $p$.
	Denote by $h\colon \operatorname{Spec}(k)\to S$ the composition $f_1\circ h_1$.
	Since $f_1$ is closed, $f_1 (p) \in |S|$ is closed.
	By the hypotheses, there exists a representative $\tilde{t}\colon \operatorname{Spec}(k^\prime)\to S$ of $f_1 (p)$ such that $\mathbf{L}(t^\prime)^\ast K$ is relatively perfect over $Y_1\times_{S} \operatorname{Spec}(k^\prime)$.
	Applying \Cref{lem:factor_residual_gerbes}, $h$ and $\tilde{t}$ represent $f_1 (p)$.
	Hence, there exist a field $\ell$ and a commutative diagram,
	\begin{displaymath}
		% https://q.uiver.app/#q=WzAsNCxbMSwwLCJcXG9wZXJhdG9ybmFtZXtTcGVjfShrXlxccHJpbWUpIl0sWzEsMSwiUy4iXSxbMCwxLCJcXG9wZXJhdG9ybmFtZXtTcGVjfShrKSJdLFswLDAsIlxcb3BlcmF0b3JuYW1le1NwZWN9KFxcZWxsKSJdLFswLDEsInQiXSxbMiwxLCJoIiwyXSxbMywyLCJiIiwyXSxbMywwLCJhIl1d
		\begin{tikzcd}
			{\operatorname{Spec}(\ell)} & {\operatorname{Spec}(k^\prime)} \\
			{\operatorname{Spec}(k)} & {S.}
			\arrow["a", from=1-1, to=1-2]
			\arrow["b"', from=1-1, to=2-1]
			\arrow["\tilde{t}", from=1-2, to=2-2]
			\arrow["h"', from=2-1, to=2-2]
		\end{tikzcd}
	\end{displaymath}
	The hypothesis says the derived pullback of $K$ along $\tilde{t}$ is relatively perfect over $Y_1 \times_S \operatorname{Spec}(k^\prime)$.
	As $a$ is affine, \Cref{lem:induced_dbcoh_perf_preservation_upon_pullback_with_t_affine} with \Cref{cor:preservation} implies the derived pullback of $K$ along $t\circ a$ is relatively perfect over $Y_1 \times_S \operatorname{Spec}(\ell)$.
	However, $b$ is faithfully flat and affine, so the first claim $\eqref{prop:reflecting_bounded_pseudocoherence_perfectness1}$ implies the derived pullback of $K$ along $h$ is relatively perfect over $Y_1 \times_S \operatorname{Spec}(k)$ (e.g.\ use that $t \circ a = h\circ b$).
	Consequently, one can now apply the argument from the second claim using $h$.
\end{proof}

\begin{theorem}
	\label[theorem]{thm:bounded_pseudocoherence_perfectness_faithfully_flat_affine}
	Consider \Cref{setup:fm_cube_pullback}.
	Let $K\in D^b_{\operatorname{coh}}(Y_1 \times_{S} Y_2)$.
	Then the following are equivalent:
	\begin{enumerate}
		\item \label[theorem]{thm:bounded_pseudocoherence_perfectness_faithfully_flat_affine1} $K$ is relatively perfect over $Y_1$ (resp.\ $Y_2$)
		\item \label[theorem]{thm:bounded_pseudocoherence_perfectness_faithfully_flat_affine2} $\mathbf{L}(t^\prime)^\ast K$ is relatively perfect over $Y_1\times_{S} T$ (resp.\ $Y_2 \times_{S} T$) for every quasi-affine morphism $t$ from a Noetherian algebraic space
		\item \label[theorem]{thm:bounded_pseudocoherence_perfectness_faithfully_flat_affine3} $\mathbf{L}(t^\prime)^\ast K$ is relatively perfect over $Y_1\times_{S} \operatorname{Spec}(k)$ (resp.\ $Y_2\times_{S} \operatorname{Spec}(k)$) for all morphisms $t\colon \operatorname{Spec}(k)\to S$ from a field
		\item \label[theorem]{thm:bounded_pseudocoherence_perfectness_faithfully_flat_affine3_closed} $\mathbf{L}(t^\prime)^\ast K$ is relatively perfect over $Y_1\times_{S} \operatorname{Spec}(k)$ (resp.\ $Y_2\times_{S} \operatorname{Spec}(k)$) for all morphisms $t\colon \operatorname{Spec}(k)\to S$ that represents a closed point $p\in |S|$
		\item \label[theorem]{thm:bounded_pseudocoherence_perfectness_faithfully_flat_affine4} $\mathbf{L}(t^\prime)^\ast K$ is relatively perfect over $Y_1\times_{S} T$ (resp.\ $Y_2\times_{S} T$) for some faithfully flat morphism $t\colon T\to S$ of Noetherian algebraic spaces.
	\end{enumerate}
\end{theorem}

\begin{proof}
	By \Cref{lem:pseudocoherence_affine_faithfully_flat_cover}, the derived pullback of $K$ along any morphism above is pseudocoherent.
	We only prove the case of being relatively perfect over $Y_1$, since the other follows analogously.
	Applying \Cref{lem:induced_dbcoh_perf_preservation_upon_pullback_with_t_affine} with \Cref{cor:preservation}, \eqref{thm:bounded_pseudocoherence_perfectness_faithfully_flat_affine1} implies \eqref{thm:bounded_pseudocoherence_perfectness_faithfully_flat_affine2}. Moreover, \eqref{thm:bounded_pseudocoherence_perfectness_faithfully_flat_affine2} implies \eqref{thm:bounded_pseudocoherence_perfectness_faithfully_flat_affine3} because $S$ is quasi-separated \cite[\href{https://stacks.math.columbia.edu/tag/09TF}{Tag 09TF}]{StacksProject}.
	Clearly, \eqref{thm:bounded_pseudocoherence_perfectness_faithfully_flat_affine3} implies \eqref{thm:bounded_pseudocoherence_perfectness_faithfully_flat_affine3_closed}. Furthermore, to see that \eqref{thm:bounded_pseudocoherence_perfectness_faithfully_flat_affine3_closed} implies \eqref{thm:bounded_pseudocoherence_perfectness_faithfully_flat_affine4}, note that \Cref{prop:reflecting_bounded_pseudocoherence_perfectness3} shows that $K$ is relatively perfect over $Y_1$, i.e.\ take $t = 1_{S}$.
	Lastly, by \Cref{prop:reflecting_bounded_pseudocoherence_perfectness1}, \eqref{thm:bounded_pseudocoherence_perfectness_faithfully_flat_affine4} implies \eqref{thm:bounded_pseudocoherence_perfectness_faithfully_flat_affine1}.
\end{proof}

\begin{example}
	\label{ex:relative_perf_counterexample}
	Set $S=\operatorname{Spec}(\mathbb{Z})$.
    Let $Y_1=Y_2=S$ with $f_1$ and $f_2$ the identity.
	Choose an infinite sequence of prime numbers $p_1 < p_2 < \cdots$.
	For each $n\geq 1$, denote by $t_n\colon \operatorname{Spec}(\mathbb{F}_{p_n}) \to S$ the associated closed immersion.
	Consider
	\begin{displaymath}
		K := \mathcal{O}_S \oplus \bigoplus_{n\geq 1} (t_n)_\ast \mathcal{O}_{\operatorname{Spec}(\mathbb{F}_{p_n})} [n] \in D^-_{\operatorname{coh}}(S).
	\end{displaymath}
	Let $t\colon \operatorname{Spec}(k)\to S$ be a representative of a point $p\in |S|$.
	If $p\not\in \{p_1,p_2,\dots\}$, then $\mathbf{L}t^\ast K\cong \mathcal{O}_{\operatorname{Spec}(k)}$.
	If $p=p_n$, then $\mathbf{L}t^\ast K\cong \mathcal{O}_{\operatorname{Spec}(k)} \oplus \mathcal{O}_{\operatorname{Spec}(k)}[n] \oplus \mathcal{O}_{\operatorname{Spec}(k)}[n+1]$.
	In either case, $\mathbf{L}t^\ast K \in D^b_{\operatorname{coh}}(\operatorname{Spec}(k))$, and so $\mathbf{L}t^\ast K$ is relatively perfect over each factor $\operatorname{Spec}(k)$.
	However, $K$ is not bounded below, and so it is not relatively perfect over either factor $S$ (take $A = \mathcal{O}_S$ in \Cref{def:relative_perf_stacks}).
	Thus, the boundedness hypothesis in \Cref{prop:reflecting_bounded_pseudocoherence_perfectness} and \Cref{thm:bounded_pseudocoherence_perfectness_faithfully_flat_affine} cannot be dropped.
\end{example}

\begin{corollary}
	[Noetherian base change for relative perfectness]
	\label{cor:noetherian_base_change_for_relative_perfectness}
	Consider \Cref{setup:fm_cube_pullback}.
	Let $K\in D^b_{\operatorname{coh}}(Y_1 \times_{S} Y_2)$ which is relatively perfect over each $Y_i$.
	Then $\mathbf{L}(t^\prime)^\ast K$ is relatively perfect over $Y_1 \times_S T$ (resp.\ $Y_2 \times_S T$).
\end{corollary}

\begin{proof}
	Choose an \'{e}tale presentation $s\colon U \to T$ from an affine scheme.
	By \Cref{lem:quasi_affine_diagonal}, $t \circ s$ is quasi-affine.
	Denote by $s^\prime\colon Y_1 \times_S Y_2 \times_S U \to Y_1 \times_S Y_2 \times_S T$ the natural morphism.
	Then \Cref{cor:preservation} and \Cref{lem:induced_dbcoh_perf_preservation_upon_pullback_with_t_affine} implies $\mathbf{L}(t^\prime \circ s^\prime)^\ast K$ is relatively perfect over $Y_1 \times_S U$ (resp.\ $Y_2 \times_S U$)
	By \cite[\href{https://stacks.math.columbia.edu/tag/08H4}{Tag 08H4}]{StacksProject}, $\mathbf{L}(t^\prime \circ s^\prime)^\ast K$ and $\mathbf{L}(t^\prime)^\ast K$ are pseudocoherent.
	By \Cref{cor:preservation}, $\mathbf{L}(t^\prime \circ s^\prime)^\ast K$ being relatively perfect over $Y_1 \times_S U$ (resp.\ $Y_2 \times_S U$) implies $\mathbf{L}(t^\prime \circ s^\prime)^\ast K$ has bounded cohomology.
	Hence, $\mathbf{L}(t^\prime \circ s^\prime)^\ast K \in D^b_{\operatorname{coh}}(Y_1 \times_S Y_2 \times_S U)$.
	Since $s$ is an \'{e}tale presentation, \Cref{lem:zero_object_via_covers} asserts that $\mathbf{L}(t^\prime)^\ast K$ has bounded cohomology.
	Applying \Cref{thm:bounded_pseudocoherence_perfectness_faithfully_flat_affine}, it follows that $\mathbf{L}(t^\prime)^\ast K$ is relatively perfect over $Y_1 \times_S T$ (resp.\ $Y_2 \times_S T$).
\end{proof}

\begin{remark}
	\label{rem:comparison_stacks_project}
	Let $f\colon Y\to X$ be a morphism of quasi-compact quasi-separated algebraic spaces which is flat and locally of finite presentation.
	In \cite[\href{https://stacks.math.columbia.edu/tag/0DKN}{Tag 0DKN}]{StacksProject} (see also \cite{Lieblich:2006}), a pseudocoherent $E\in D(\mathcal{O}_Y)$ is called \textbf{$X$-perfect} if it locally has finite tor dimension as an object of $D(f^{-1}\mathcal{O}_X)$.
	This agrees with $f$-perfectness.
	Indeed, both notions are \'{e}tale local (\Cref{thm:relative_perf_iff_smooth_locally} and \cite[\href{https://stacks.math.columbia.edu/tag/0DKQ}{Tag 0DKQ}]{StacksProject}), and so one can apply \cite[Proposition 3.2]{GuisadoVillalgordo/Lank/ManaliRahul/Pavic:2025}.
	Unlike the Stacks Project, \Cref{def:relative_perf_stacks} requires neither flatness nor pseudocoherence.
	In the flat case, \cite[\href{https://stacks.math.columbia.edu/tag/0DKT}{Tags 0DKT}, \href{https://stacks.math.columbia.edu/tag/0DKX}{0DKX}, \href{https://stacks.math.columbia.edu/tag/0GFK}{0GFK}, and \href{https://stacks.math.columbia.edu/tag/0GFL}{0GFL}]{StacksProject} yield alternative proofs of parts of \Cref{thm:stacky_all_coincide,introthm:base_change}.
	Our proofs are independent and work with $t$-structures.
\end{remark}

%%%%%%%%%%%%%%%%%%%%%%%%%%%%%%%%%%%
\part{Fiberwise criteria}
\label{sec:part2}
%%%%%%%%%%%%%%%%%%%%%%%%%%%%%%%%%%%

%%%%%%%%%%%%%%%%%%%%%%%%%%%%%%%%%%%
\section{Mates in geometry}
\label{sec:mates_in_geometry}
%%%%%%%%%%%%%%%%%%%%%%%%%%%%%%%%%%%

In this section, we state the general results about mates and the applications to geometry required in the rest of the article.

\subsection{Mates}
\label{sec:mates}
We refer \cite[\S 2]{Kelly/Street:1974}.
Recall that given a $2$-category $\mathcal{T}$, an adjunction in $\mathcal{T}$ is a tuple $(F,G,\eta,\epsilon)$, where $F\colon \mathcal{C}\to\mathcal{C}^\prime$ and $G\colon \mathcal{C}^\prime\to\mathcal{C}$ are $1$-morphisms in $\mathcal{T}$ and $\eta\colon \operatorname{id}_{\mathcal{C}}\Rightarrow GF$ and $\epsilon\colon FG\Rightarrow\operatorname{id}_{\mathcal{C}^\prime}$ are $2$-morphisms, respectively referred to as the unit and the counit of the adjunction, that satisfy the triangle identities:
\begin{displaymath}
	\begin{aligned}
		(\epsilon F)\circ (F\eta)
		&=\operatorname{id}_F,\\
		(G\epsilon)\circ(\eta G)&=\operatorname{id}_G.
	\end{aligned}
\end{displaymath}
We also say `$F\dashv G$ in $\mathcal{T}$ with unit and counit $\eta$ and $\epsilon$.'

\begin{setup}
	\label{setup:mates}
	Let $\mathcal{T}$ be a $2$-category.
	Consider a diagram in $\mathcal{T}$:
	\begin{equation}
		\label{diag:morph_adj}
		% https://q.uiver.app/#q=WzAsNCxbMCwwLCJcXG1hdGhjYWx7Q30iXSxbMSwwLCJcXG1hdGhjYWx7Q31eXFxwcmltZSJdLFswLDEsIlxcbWF0aGNhbHtEfSJdLFsxLDEsIlxcbWF0aGNhbHtEfV5cXHByaW1lIl0sWzAsMSwiSCJdLFsyLDMsIksiLDJdLFsyLDAsIkciLDIseyJjdXJ2ZSI6MX1dLFswLDIsIkYiLDIseyJjdXJ2ZSI6MX1dLFszLDEsIkdeXFxwcmltZSIsMix7ImN1cnZlIjoxfV0sWzEsMywiRl5cXHByaW1lIiwyLHsiY3VydmUiOjF9XSxbOSw4LCJcXGRhc2h2IiwzLHsic2hvcnRlbiI6eyJzb3VyY2UiOjMwLCJ0YXJnZXQiOjMwfSwibGV2ZWwiOjEsInN0eWxlIjp7ImJvZHkiOnsibmFtZSI6Im5vbmUifSwiaGVhZCI6eyJuYW1lIjoibm9uZSJ9fX1dLFs3LDYsIlxcZGFzaHYiLDMseyJzaG9ydGVuIjp7InNvdXJjZSI6MzAsInRhcmdldCI6MzB9LCJsZXZlbCI6MSwic3R5bGUiOnsiYm9keSI6eyJuYW1lIjoibm9uZSJ9LCJoZWFkIjp7Im5hbWUiOiJub25lIn19fV1d
		\begin{tikzcd}
			{\mathcal{C}} & {\mathcal{C}^\prime} \\
			{\mathcal{D}} & {\mathcal{D}^\prime,x}
			\arrow["H", from=1-1, to=1-2]
			\arrow[""{name=0, anchor=center, inner sep=0}, "F"', bend right=30, from=1-1, to=2-1]
			\arrow[""{name=1, anchor=center, inner sep=0}, "{F^\prime}"', bend right=30, from=1-2, to=2-2]
			\arrow[""{name=2, anchor=center, inner sep=0}, "G"', bend right=30, from=2-1, to=1-1]
			\arrow["K"', from=2-1, to=2-2]
			\arrow[""{name=3, anchor=center, inner sep=0}, "{G^\prime}"', bend right=30, from=2-2, to=1-2]
			\arrow["\dashv"{description}, draw=none, from=2, to=0]
			\arrow["\dashv"{description}, draw=none, from=3, to=1]
		\end{tikzcd}
	\end{equation}
	where we have adjunctions $(F,G,\eta,\epsilon)$ and $(F^\prime,G^\prime,\eta^\prime,\epsilon^\prime)$.
\end{setup}

\begin{lemma}
	\label[lemma]{lem:mates_bijection}
	Assume \Cref{setup:mates}.
	\begin{enumerate}[label=(\arabic*), ref=\thelemma(\arabic*)]
		\item\label[lemma]{lem:mates_bijection1}
		The following diagram commutes:
		\begin{equation}
			\label{diag:mates_bijection}
			% https://q.uiver.app/#q=WzAsNCxbMCwwLCJcXG9wZXJhdG9ybmFtZXtIb219KEYnSCxLRikiXSxbMSwwLCJcXG9wZXJhdG9ybmFtZXtIb219KEgsRydLRikiXSxbMSwxLCJcXG9wZXJhdG9ybmFtZXtIb219KEhHLEcnSykiXSxbMCwxLCJcXG9wZXJhdG9ybmFtZXtIb219KEYnSEcsSykiXSxbMCwxXSxbMSwyXSxbMCwzXSxbMywyXV0=
			\begin{tikzcd}
				{\operatorname{Hom}(F^\prime H,KF)} & {\operatorname{Hom}(H,G^\prime KF)} \\
				{\operatorname{Hom}(F^\prime HG,K)} & {\operatorname{Hom}(HG,G^\prime K)}
				\arrow["\cong", from=1-1, to=1-2]
				\arrow["\cong"', from=1-1, to=2-1]
				\arrow["\cong", from=1-2, to=2-2]
				\arrow["\cong", from=2-1, to=2-2]
			\end{tikzcd}
		\end{equation}
		where the bijections are given by \cite[Proposition 4.3.7, Remark 4.3.9]{Riehl:26}.
		\item\label[lemma]{lem:mates_bijection2}
		There is a canonical bijection
		\begin{align*}
			\{ \textrm{$2$-morphisms } \alpha\colon F^\prime H\to KF \}
			&\cong \{ \textrm{$2$-morphisms } \beta\colon HG\to G^\prime K \}
			\intertext{given by the mutually inverse maps}
			\alpha
			&\mapsto \overline{\alpha}
			=(G^\prime K\epsilon) \circ (G^\prime\alpha G) \circ (\eta^\prime HG)\\
			(\epsilon^\prime KF)\circ(F^\prime\beta F)\circ(F^\prime H\eta)
			=\overline{\beta}
			&\mapsfrom\beta
		\end{align*}
	\end{enumerate}
\end{lemma}

\begin{proof}
	The following proof follows \cite[Proposition 4.3.8 and footnote 8]{Riehl:26}.
	See also \cite[Proposition 2.1]{Kelly/Street:1974}, \cite[Lemma 4.2]{Bergh/Schnurer:2020}, \cite[Sect. 1]{Cheng/Gurski/Riehl:2014}.

	To prove \ref{lem:mates_bijection1}, observe that given $\alpha\colon F^\prime H\to KF$, we have
	\begin{equation}
		\label{eq:alpha_bar}
		G^\prime K\epsilon\circ(G^\prime\alpha\circ\eta^\prime H)G
		=(G^\prime K\epsilon) \circ (G^\prime\alpha G) \circ (\eta^\prime HG)
		=G^\prime(K\epsilon\circ\alpha G)\circ\eta^\prime HG,
	\end{equation}
	by \cite[Appendix, \S 1, (I), (II), (III)]{Godement:1958}.\footnote{Even though Godement states his formulas for the category of categories, these are valid in any $2$-category.}
	The LHS and RHS are the image of $\alpha$ under the clockwise and counterclockwise composites of the square \eqref{diag:mates_bijection}, respectively.

	The proof of \ref{lem:mates_bijection2} follows by noting \eqref{eq:alpha_bar} is $\overline{\alpha}$ and that the image of $\beta\colon HG\to G^\prime K$ under the inverse commutative square to \eqref{diag:mates_bijection} gives the stated formula for $\overline{\beta}$.
\end{proof}

\begin{definition}
	Consider \Cref{setup:mates}.
	\begin{enumerate}
		\item In \Cref{lem:mates_bijection}, if $\beta=\overline{\alpha}$ (equivalently, $\alpha=\overline{\beta}$), then we say that $(\alpha,\beta)$ are \textbf{mates} under the diagram \eqref{diag:morph_adj} (with respect to the adjunctions $(F^\prime,G^\prime,\eta^\prime,\epsilon^\prime)$ and $(F,G,\eta,\epsilon)$).
		We also say $\alpha$ is the \textbf{left mate} (of $\beta$) and $\beta$ is the \textbf{right mate} (of $\alpha$).

		\item If $H=1_{\mathcal{C}}$ and $K=1_{\mathcal{D}}$ are the identity functor,
		we say that $\alpha\colon F^\prime\to F$ and $\beta\colon G\to G^\prime$ are \textbf{conjugate $2$-morphisms} if they are mates under \eqref{diag:morph_adj} (see \cite[p. 100]{MacLane:1978}, \cite[(3.3.5)]{Lipman:2009}, or \cite[comments before Remark 4.3.9]{Riehl:26}).
		We also say $\alpha$ is the \textbf{left conjugate} (of $\beta$) and $\beta$ is the \textbf{right conjugate} (of $\alpha$).
	\end{enumerate}
\end{definition}
Diagramatically, we will picture a pair of mates under \eqref{diag:morph_adj} as
\begin{equation}
	\label{eq:mate_correspondence_2-cell}
	% https://q.uiver.app/#q=WzAsOCxbMCwwLCJcXG1hdGhjYWx7Q30iXSxbMSwwLCJcXG1hdGhjYWx7Q31eXFxwcmltZSJdLFswLDEsIlxcbWF0aGNhbHtEfSJdLFsxLDEsIlxcbWF0aGNhbHtEfV5cXHByaW1lIl0sWzMsMCwiXFxtYXRoY2Fse0N9Il0sWzQsMCwiXFxtYXRoY2Fse0N9XlxccHJpbWUiXSxbMywxLCJcXG1hdGhjYWx7RH0iXSxbNCwxLCJcXG1hdGhjYWx7RH1eXFxwcmltZSJdLFswLDEsIkgiXSxbMiwzLCJLIiwyXSxbMCwyLCJGIiwyXSxbMSwzLCJGXlxccHJpbWUiXSxbMSwyLCJcXGFscGhhIiwyLHsic2hvcnRlbiI6eyJzb3VyY2UiOjIwLCJ0YXJnZXQiOjIwfSwibGV2ZWwiOjJ9XSxbNCw1LCJIIl0sWzYsNywiSyIsMl0sWzYsNCwiRyJdLFs3LDUsIkdeXFxwcmltZSIsMl0sWzQsNywiXFxiZXRhIiwwLHsic2hvcnRlbiI6eyJzb3VyY2UiOjIwLCJ0YXJnZXQiOjIwfSwibGV2ZWwiOjJ9XSxbMTEsMTUsIlxcbGVmdHJpZ2h0c3F1aWdhcnJvdyIsMSx7InNob3J0ZW4iOnsic291cmNlIjoyMCwidGFyZ2V0IjoyMH0sInN0eWxlIjp7ImJvZHkiOnsibmFtZSI6Im5vbmUifSwiaGVhZCI6eyJuYW1lIjoibm9uZSJ9fX1dXQ==
	\begin{tikzcd}
		{\mathcal{C}} & {\mathcal{C}^\prime} && {\mathcal{C}} & {\mathcal{C}^\prime} \\
		{\mathcal{D}} & {\mathcal{D}^\prime} && {\mathcal{D}} & {\mathcal{D}^\prime}
		\arrow["H", from=1-1, to=1-2]
		\arrow["F"', from=1-1, to=2-1]
		\arrow["\alpha"',  shorten <=2pt, shorten >=2pt, Rightarrow,  from=1-2, to=2-1]
		\arrow[""{name=0, anchor=center, inner sep=0}, "{F^\prime}", from=1-2, to=2-2]
		\arrow["H", from=1-4, to=1-5]
		\arrow["\beta",  shorten <=2pt, shorten >=2pt, Rightarrow,  from=1-4, to=2-5]
		\arrow["K"', from=2-1, to=2-2]
		\arrow[""{name=1, anchor=center, inner sep=0}, "G", from=2-4, to=1-4]
		\arrow["K"', from=2-4, to=2-5]
		\arrow["{G^\prime}"', from=2-5, to=1-5]
		\arrow["\leftrightsquigarrow"{description}, draw=none, from=0, to=1]
	\end{tikzcd}
\end{equation}
Accordingly, the left mate $^\alpha$\hspace{-0.2em}\rotatebox[origin=c]{225}{$\Rightarrow$} points to the left
and that the right mate \rotatebox[origin=c]{-45}{$\Rightarrow$}$^\beta$ points to the right.
\begin{lemma}
	\label[lemma]{lem:morph_adj}
	Let $\alpha\colon F^\prime H\to KF$ and $\beta\colon HG\to G^\prime K$ be $2$-morphisms.
	Then the following conditions are equivalent:
	\begin{enumerate}[
		label=(\arabic*),
		ref=\thelemma(\arabic*)
		]
		\item
		\label[lemma]{lem:morph_adj1}
		$\alpha$ and $\beta$ are mates.

		\item
		\label[lemma]{lem:morph_adj2}
		The following diagram commutes:
		\begin{displaymath}
			% https://q.uiver.app/#q=WzAsNCxbMCwwLCJGJ0hHIl0sWzAsMSwiRidHJ0siXSxbMSwxLCJLLiJdLFsxLDAsIktGRyJdLFswLDMsIlxcYWxwaGEgSCJdLFswLDEsIkYnXFxiZXRhIiwyXSxbMSwyLCJcXHZhcmVwc2lsb24nSyIsMl0sWzMsMiwiS1xcdmFyZXBzaWxvbiJdXQ==
			\begin{tikzcd}
				{F^\prime HG} & KFG \\
				{F^\prime G^\prime K} & {K.}
				\arrow["{\alpha G}", from=1-1, to=1-2]
				\arrow["{F^\prime \beta}"', from=1-1, to=2-1]
				\arrow["{K\epsilon}", from=1-2, to=2-2]
				\arrow["{\epsilon^\prime K}"', from=2-1, to=2-2]
			\end{tikzcd}
		\end{displaymath}

		\item
		\label[lemma]{lem:morph_adj3}
		The following diagram commutes:
		\begin{displaymath}
			% https://q.uiver.app/#q=WzAsNCxbMSwxLCJHXlxccHJpbWUgS0YuIl0sWzAsMCwiSCJdLFsxLDAsIkdeXFxwcmltZSBGXlxccHJpbWUgSCJdLFswLDEsIkhHRiJdLFsyLDAsIkdeXFxwcmltZSBcXGFscGhhIl0sWzMsMCwiXFxiZXRhIEYiLDJdLFsxLDMsIkhcXGV0YSIsMl0sWzEsMiwiXFxldGFeXFxwcmltZSBIIl1d
			\begin{tikzcd}
				H & {G^\prime F^\prime H} \\
				HGF & {G^\prime KF.}
				\arrow["{\eta^\prime H}", from=1-1, to=1-2]
				\arrow["{H\eta}"', from=1-1, to=2-1]
				\arrow["{G^\prime \alpha}", from=1-2, to=2-2]
				\arrow["{\beta F}"', from=2-1, to=2-2]
			\end{tikzcd}
		\end{displaymath}
	\end{enumerate}
	Moreover, if $\mathcal{T}$ is the $2$-category of categories, then the previous conditions are also equivalent to
	\begin{enumerate}[
		label=(\arabic*),
		ref=\thelemma(\arabic*)
		]
		\setcounter{enumi}{3}
		\item
		\label{lem:morph_adj4}
		For all objects $c\in\mathcal{C}$ and $d\in\mathcal{D}$, the following diagram commutes:
		\begin{displaymath}
			% https://q.uiver.app/#q=WzAsNixbMCwwLCJcXG1hdGhjYWx7RH0oRmMsZCkiXSxbMCwxLCJcXG1hdGhjYWx7RH1eXFxwcmltZShLRmMsS2QpIl0sWzAsMiwiXFxtYXRoY2Fse0R9XlxccHJpbWUgKEZeXFxwcmltZSBIYyxLZCkiXSxbMSwwLCJcXG1hdGhjYWx7Q30oYyxHZCkiXSxbMSwxLCJcXG1hdGhjYWx7Q31eXFxwcmltZSAoSGMsSEdkKSJdLFsxLDIsIlxcbWF0aGNhbHtDfV5cXHByaW1lIChIYyxHXlxccHJpbWUgS2QpLiJdLFswLDMsIlxcY29uZyJdLFsyLDUsIlxcY29uZyJdLFswLDEsIksiLDJdLFsxLDIsIlxcYWxwaGFfY15cXGFzdCIsMl0sWzMsNCwiSCJdLFs0LDUsIlxcYmV0YV97ZCxcXGFzdH0iXV0=
			\begin{tikzcd}
				{\mathcal{D}(Fc,d)} & {\mathcal{C}(c,Gd)} \\
				{\mathcal{D}^\prime(KFc,Kd)} & {\mathcal{C}^\prime (Hc,HGd)} \\
				{\mathcal{D}^\prime (F^\prime Hc,Kd)} & {\mathcal{C}^\prime (Hc,G^\prime Kd).}
				\arrow["\cong", from=1-1, to=1-2]
				\arrow["K"', from=1-1, to=2-1]
				\arrow["H", from=1-2, to=2-2]
				\arrow["{\alpha_c^\ast}"', from=2-1, to=3-1]
				\arrow["{\beta_{d,\ast}}", from=2-2, to=3-2]
				\arrow["\cong", from=3-1, to=3-2]
			\end{tikzcd}
		\end{displaymath}
		where the top and bottom isomorphisms are respectively induced by $\eta$ and $\eta^\prime$ (equivalently, the inverse isomorphisms are respectively induced by $\epsilon$ and $\epsilon^\prime$).
	\end{enumerate}
\end{lemma}

Equivalence of \ref{lem:morph_adj4} with all other conditions is \cite[IV.7, Exercise 4]{MacLane:1978}.

\begin{proof}
	\ref{lem:morph_adj1} $\iff$ \ref{lem:morph_adj2}.
	The right mate $\overline{\alpha}$ of $\alpha$ is by definition the $2$-morphism
	\begin{displaymath}
		\overline{\alpha}\colon
		HG\xrightarrow{\eta^\prime HG}
		G^\prime F^\prime HG\xrightarrow{G^\prime \alpha G}
		G^\prime KFG\xrightarrow{G^\prime K\epsilon}G^\prime K.
	\end{displaymath}
	Therefore $\overline{\alpha}$ will equal $\beta$ if, and only if, the corresponding adjunct morphisms $F^\prime HG\to K$ (w.r.t.\ $F^\prime\dashv G^\prime$) are equal. This is exactly condition \ref{lem:morph_adj2}.

	\ref{lem:morph_adj1} $\iff$ \ref{lem:morph_adj3}.
	The left mate $\overline{\beta}$ of $\beta$ is by definition the $2$-morphism
	\begin{displaymath}
		\overline{\beta}\colon F^\prime H\xrightarrow{F^\prime H\eta}
		F^\prime HGF\xrightarrow{F^\prime \beta F}
		F^\prime G^\prime KF\xrightarrow{\epsilon^\prime KF}KF
	\end{displaymath}
	Therefore $\overline{\beta}$ will equal $\alpha$ if, and only if, the corresponding adjunct morphisms $H\to G^\prime KF$ (w.r.t.\ $F^\prime \dashv G^\prime$) are equal. This is exactly condition \ref{lem:morph_adj3}.

	\ref{lem:morph_adj3} $\iff$ \ref{lem:morph_adj4}.
	Suppose $\mathcal{T}$ is the $2$-category of categories.
	Let $c\in\mathcal{C}$, $d\in\mathcal{D}$, and $\phi\colon Fc\to d$.
	Going down the left column and then along the bottom row of the diagram in \ref{lem:morph_adj4}, $\phi$ is sent to $G^\prime K(\phi)\circ G^\prime (\alpha_c)\circ\eta^\prime_{Hc}$.
	Going along the top row and then down the right column, $\phi$ is sent to
	\begin{displaymath}
		\beta_d\circ HG(\phi)\circ H(\eta_c) = G^\prime K(\phi)\circ\beta_{Fc}\circ H(\eta_c),
	\end{displaymath}
	where the equality is naturality of $\beta$.
	Hence, \ref{lem:morph_adj3} implies \ref{lem:morph_adj4}.
	Conversely, taking $d=Fc$ and $\phi=1_{Fc}$ shows that \ref{lem:morph_adj4} implies \ref{lem:morph_adj3}.
\end{proof}

\begin{remark}\hfill
	\begin{enumerate}[label=(\arabic*),ref=\theremark(\arabic*)]
		\item \label{rem:unit_counit_is_mate} The unit of an adjunction is the right mate of the identity transformation of the left adjoint
		(dually, the analogous is true for the counit by swapping `left' and `right').
		More precisely, the two pairs $(\operatorname{id}_F,\eta)$ and $(\epsilon,\operatorname{id}_G)$ are mates respectively under the diagrams
		\begin{displaymath}
			% https://q.uiver.app/#q=WzAsOCxbMSwwLCJcXG1hdGhjYWx7Q30iXSxbMCwwLCJcXG1hdGhjYWx7Q30iXSxbMSwxLCJcXG1hdGhjYWx7RH0iXSxbMCwxLCJcXG1hdGhjYWx7Q30iXSxbNCwwLCJcXG1hdGhjYWx7RH0iXSxbNCwxLCJcXG1hdGhjYWx7RH0iXSxbNSwxLCJcXG1hdGhjYWx7RH0uIl0sWzUsMCwiXFxtYXRoY2Fse0N9Il0sWzEsMCwiMV97XFxtYXRoY2Fse0N9fSJdLFszLDIsIkYiLDJdLFsyLDAsIkciLDIseyJjdXJ2ZSI6MX1dLFswLDIsIkYiLDIseyJjdXJ2ZSI6MX1dLFszLDEsIjFfe1xcbWF0aGNhbHtDfX0iLDIseyJjdXJ2ZSI6MX1dLFsxLDMsIjFfe1xcbWF0aGNhbHtDfX0iLDIseyJjdXJ2ZSI6MX1dLFs0LDUsIjFfXFxtYXRoY2Fse0R9IiwyLHsiY3VydmUiOjF9XSxbNSw0LCIxX1xcbWF0aGNhbHtEfSIsMix7ImN1cnZlIjoxfV0sWzQsNywiRyJdLFs1LDYsIjFfXFxtYXRoY2Fse0R9IiwyXSxbNyw2LCJGIiwyLHsiY3VydmUiOjF9XSxbNiw3LCJHIiwyLHsiY3VydmUiOjF9XSxbMTAsMTQsIlxcdGV4dHtcXG5vcm1hbHNpemUgYW5kfSIsMSx7InNob3J0ZW4iOnsic291cmNlIjoyMCwidGFyZ2V0IjoyMH0sInN0eWxlIjp7ImJvZHkiOnsibmFtZSI6Im5vbmUifSwiaGVhZCI6eyJuYW1lIjoibm9uZSJ9fX1dLFsxMywxMiwiXFxkYXNodiIsMSx7InNob3J0ZW4iOnsic291cmNlIjoyMCwidGFyZ2V0IjoyMH0sInN0eWxlIjp7ImJvZHkiOnsibmFtZSI6Im5vbmUifSwiaGVhZCI6eyJuYW1lIjoibm9uZSJ9fX1dLFsxMSwxMCwiXFxkYXNodiIsMSx7InNob3J0ZW4iOnsic291cmNlIjoyMCwidGFyZ2V0IjoyMH0sInN0eWxlIjp7ImJvZHkiOnsibmFtZSI6Im5vbmUifSwiaGVhZCI6eyJuYW1lIjoibm9uZSJ9fX1dLFsxNCwxNSwiXFxkYXNodiIsMSx7InNob3J0ZW4iOnsic291cmNlIjoyMCwidGFyZ2V0IjoyMH0sInN0eWxlIjp7ImJvZHkiOnsibmFtZSI6Im5vbmUifSwiaGVhZCI6eyJuYW1lIjoibm9uZSJ9fX1dLFsxOCwxOSwiXFxkYXNodiIsMSx7InNob3J0ZW4iOnsic291cmNlIjoyMCwidGFyZ2V0IjoyMH0sInN0eWxlIjp7ImJvZHkiOnsibmFtZSI6Im5vbmUifSwiaGVhZCI6eyJuYW1lIjoibm9uZSJ9fX1dXQ==
			\begin{tikzcd}
				{\mathcal{C}} & {\mathcal{C}} &&& {\mathcal{D}} & {\mathcal{C}} \\
				{\mathcal{C}} & {\mathcal{D}} &&& {\mathcal{D}} & {\mathcal{D}.}
				\arrow["{1_{\mathcal{C}}}", from=1-1, to=1-2]
				\arrow[""{name=0, anchor=center, inner sep=0}, "{1_{\mathcal{C}}}"', curve={height=6pt}, from=1-1, to=2-1]
				\arrow[""{name=1, anchor=center, inner sep=0}, "F"', curve={height=6pt}, from=1-2, to=2-2]
				\arrow["G", from=1-5, to=1-6]
				\arrow[""{name=2, anchor=center, inner sep=0}, "{1_\mathcal{D}}"', curve={height=6pt}, from=1-5, to=2-5]
				\arrow[""{name=3, anchor=center, inner sep=0}, "F"', curve={height=6pt}, from=1-6, to=2-6]
				\arrow[""{name=4, anchor=center, inner sep=0}, "{1_{\mathcal{C}}}"', curve={height=6pt}, from=2-1, to=1-1]
				\arrow["F"', from=2-1, to=2-2]
				\arrow[""{name=5, anchor=center, inner sep=0}, "G"', curve={height=6pt}, from=2-2, to=1-2]
				\arrow[""{name=6, anchor=center, inner sep=0}, "{1_\mathcal{D}}"', curve={height=6pt}, from=2-5, to=1-5]
				\arrow["{1_\mathcal{D}}"', from=2-5, to=2-6]
				\arrow[""{name=7, anchor=center, inner sep=0}, "G"', curve={height=6pt}, from=2-6, to=1-6]
				\arrow["\dashv"{description}, draw=none, from=0, to=4]
				\arrow["\dashv"{description}, draw=none, from=1, to=5]
				\arrow["\dashv"{description}, draw=none, from=2, to=6]
				\arrow["\dashv"{description}, draw=none, from=3, to=7]
				\arrow["{\text{\normalsize and}}"{description}, draw=none, from=5, to=2]
			\end{tikzcd}
		\end{displaymath}

		\item \label{rem:mate_of_iso_needs_not_be_iso} It is not true in general that the mate of an isomorphism is an isomorphism.
		By \ref{rem:unit_counit_is_mate}, any adjunction whose unit or counit is not an isomorphism is an example of such behaviour (consider e.g.\ the free-forgetful adjunction of sets and groups).
		Nonetheless, see \Cref{lem:when_mate_correspondence_preserves_isos}.
		In the geometric setting, to show that a certain mate is an isomorphism, one requires extra geometric conditions, see e.g.\ \Cref{lem:base_change_relative_formula,lem:characterize_tor-independence}.
	\end{enumerate}
\end{remark}

Taking conjugates of a composed pair of adjoints amounts to taking mates twice.
Specifically:

\begin{lemma}[{{\cite[Lemma 2.12]{Shulman:2011}}}]
	\label{lem:taking_conjugates_is_taking_mates_twice}
	Suppose we have a $2$-category $\mathcal{T}$ and a diagram
	\begin{equation}
		\label{diag:square_conj_is_mates_twice}
		% https://q.uiver.app/#q=WzAsNCxbMCwwLCJYXlxccHJpbWUiXSxbMSwwLCJYIl0sWzEsMSwiWi4iXSxbMCwxLCJaXlxccHJpbWUiXSxbMSwyLCJwX1xcYXN0Il0sWzAsMSwicV5cXHByaW1lX1xcYXN0Il0sWzMsMiwicV9cXGFzdCIsMl0sWzAsMywicF5cXHByaW1lX1xcYXN0IiwyXV0=
		\begin{tikzcd}[cramped]
			{X^\prime} & X \\
			{Z^\prime} & {Z}
			\arrow["{q^\prime_\ast}", from=1-1, to=1-2]
			\arrow["{p^\prime_\ast}"', from=1-1, to=2-1]
			\arrow["{p_\ast}", from=1-2, to=2-2]
			\arrow["{q_\ast}"', from=2-1, to=2-2]
		\end{tikzcd}
	\end{equation}
	of $1$-morphisms in $\mathcal{T}$, not necessarily commutative.
	Assume $u_\ast$ as a left adjoint $u^\ast$, for $u\in\{p,p^\prime,q,q^\prime\}$.
	Then we have a commutative diagram of bijections:
	\begin{equation}
		\label{diag:triangle_conj_mates}
		% https://q.uiver.app/#q=WzAsMyxbMCwwLCJcXG9wZXJhdG9ybmFtZXtIb219KHBee1xccHJpbWVcXGFzdH1xXlxcYXN0LHFee1xccHJpbWVcXGFzdH1wXlxcYXN0KSJdLFsxLDEsIlxcb3BlcmF0b3JuYW1le0hvbX0ocV5cXGFzdCBwX1xcYXN0LHBfXFxhc3ReXFxwcmltZSBxXntcXHByaW1lXFxhc3R9KSJdLFsyLDAsIlxcb3BlcmF0b3JuYW1le0hvbX0ocF9cXGFzdCBxX1xcYXN0XlxccHJpbWUscV9cXGFzdCBwX1xcYXN0XlxccHJpbWUpIl0sWzAsMiwiXFx0ZXh0e2Nvbmp1Z2F0ZXN9IiwwLHsic3R5bGUiOnsidGFpbCI6eyJuYW1lIjoiYXJyb3doZWFkIn19fV0sWzAsMSwiXFx0ZXh0e21hdGVzIHcuci50Ln1cXFxccF57XFxwcmltZVxcYXN0fVxcZGFzaHYgcF5cXHByaW1lX1xcYXN0XFxcXHBeXFxhc3RcXGRhc2h2IHBfXFxhc3QiLDIseyJzdHlsZSI6eyJ0YWlsIjp7Im5hbWUiOiJhcnJvd2hlYWQifX19XSxbMSwyLCJcXHRleHR7bWF0ZXMgdy5yLnQufVxcXFxxXlxcYXN0XFxkYXNodiBxX1xcYXN0XFxcXHFee1xccHJpbWVcXGFzdH1cXGRhc2h2IHFfXFxhc3ReXFxwcmltZSIsMix7InN0eWxlIjp7InRhaWwiOnsibmFtZSI6ImFycm93aGVhZCJ9fX1dXQ==
		\begin{tikzcd}[cramped,row sep=2.25em]
			{\operatorname{Hom}(p^{\prime\ast}q^\ast,q^{\prime\ast}p^\ast)} && {\operatorname{Hom}(p_\ast q_\ast^\prime,q_\ast p_\ast^\prime)} \\
			& {\operatorname{Hom}(q^\ast p_\ast,p_\ast^\prime q^{\prime\ast})}
			\arrow["{\textrm{conjugates}}", leftrightarrow, from=1-1, to=1-3]
			\arrow["\def\arraystretch{0.8}\begin{array}{c} \textrm{mates w.r.t.}\\p^{\prime\ast}\dashv p^\prime_\ast\\p^\ast\dashv p_\ast \end{array}"', leftrightarrow, from=1-1, to=2-2]
			\arrow["\def\arraystretch{0.8}\begin{array}{c} \textrm{mates w.r.t.}\\q^\ast\dashv q_\ast\\q^{\prime\ast}\dashv q_\ast^\prime \end{array}"', leftrightarrow, from=2-2, to=1-3]
		\end{tikzcd}
	\end{equation}
\end{lemma}

\begin{proof}
	We will show that the top bijection equals the composite `top left$\to$bottom$\to$top right.'
	That is, taking the right conjugate of a $2$-morphism $p^{\prime\ast}q^\ast\to q^{\prime\ast}p^\ast$ amounts to first taking its right mate $q^\ast p_\ast\to p^\prime_\ast q^{\prime\ast}$ w.r.t.\ $p^{\prime\ast}\dashv p^\prime_\ast$ and $p^\ast\dashv p_\ast$ and then taking the subsequent right mate $p_\ast q^\prime_\ast\to q_\ast p^\prime_\ast$ w.r.t.\ $q^\ast\dashv q_\ast$ and $q^{\prime\ast}\dashv q^\prime_\ast$.
	Consider the following diagram:
	\begin{equation}
		\label{diag:conj_is_mates_twice}
		% https://q.uiver.app/#q=WzAsOSxbMCwwLCJcXG9wZXJhdG9ybmFtZXtIb219KHBee1xccHJpbWVcXGFzdH1xXlxcYXN0LHFee1xccHJpbWVcXGFzdH1wXlxcYXN0KSJdLFsxLDAsIlxcb3BlcmF0b3JuYW1le0hvbX0ocV5cXGFzdCxwX1xcYXN0XlxccHJpbWUgcV57XFxwcmltZVxcYXN0fXBeXFxhc3QpIl0sWzEsMSwiXFxvcGVyYXRvcm5hbWV7SG9tfShxXlxcYXN0IHBfXFxhc3QscF9cXGFzdF5cXHByaW1lIHFee1xccHJpbWVcXGFzdH0pIl0sWzAsMSwiXFxvcGVyYXRvcm5hbWV7SG9tfShwXntcXHByaW1lXFxhc3R9cV5cXGFzdCBwX1xcYXN0LHFee1xccHJpbWVcXGFzdH0pIl0sWzIsMCwiXFxvcGVyYXRvcm5hbWV7SG9tfShcXG9wZXJhdG9ybmFtZXtpZH0scV9cXGFzdCBwX1xcYXN0XlxccHJpbWUgcV57XFxwcmltZVxcYXN0fXBeXFxhc3QpIl0sWzIsMSwiXFxvcGVyYXRvcm5hbWV7SG9tfShwX1xcYXN0LHFfXFxhc3QgcF9cXGFzdF5cXHByaW1lIHFee1xccHJpbWVcXGFzdH0pIl0sWzIsMiwiXFxvcGVyYXRvcm5hbWV7SG9tfShwX1xcYXN0IHFfXFxhc3ReXFxwcmltZSxxX1xcYXN0IHBfXFxhc3ReXFxwcmltZSkiXSxbMCwyLCJcXG9wZXJhdG9ybmFtZXtIb219KHBee1xccHJpbWVcXGFzdH1xXlxcYXN0IHBfXFxhc3QgcV9cXGFzdF5cXHByaW1lLFxcb3BlcmF0b3JuYW1le2lkfSkiXSxbMSwyLCJcXG9wZXJhdG9ybmFtZXtIb219KHFeXFxhc3QgcF9cXGFzdCBxX1xcYXN0XlxccHJpbWUscF9cXGFzdF5cXHByaW1lKSJdLFswLDFdLFsxLDRdLFs0LDVdLFs1LDZdLFswLDNdLFszLDJdLFsyLDhdLFszLDddLFs3LDhdLFs4LDZdLFsyLDVdLFsxLDJdXQ==
		\begin{tikzcd}[cramped]
			{\operatorname{Hom}(p^{\prime\ast}q^\ast,q^{\prime\ast}p^\ast)} & {\operatorname{Hom}(q^\ast,p_\ast^\prime q^{\prime\ast}p^\ast)} & {\operatorname{Hom}(\operatorname{id},q_\ast p_\ast^\prime q^{\prime\ast}p^\ast)} \\
			{\operatorname{Hom}(p^{\prime\ast}q^\ast p_\ast,q^{\prime\ast})} & {\operatorname{Hom}(q^\ast p_\ast,p_\ast^\prime q^{\prime\ast})} & {\operatorname{Hom}(p_\ast,q_\ast p_\ast^\prime q^{\prime\ast})} \\
			{\operatorname{Hom}(p^{\prime\ast}q^\ast p_\ast q_\ast^\prime,\operatorname{id})} & {\operatorname{Hom}(q^\ast p_\ast q_\ast^\prime,p_\ast^\prime)} & {\operatorname{Hom}(p_\ast q_\ast^\prime,q_\ast p_\ast^\prime)}
			\arrow[from=1-1, to=1-2]
			\arrow[from=1-1, to=2-1]
			\arrow[from=1-2, to=1-3]
			\arrow[from=1-2, to=2-2]
			\arrow[from=1-3, to=2-3]
			\arrow[from=2-1, to=2-2]
			\arrow[from=2-1, to=3-1]
			\arrow[from=2-2, to=2-3]
			\arrow[from=2-2, to=3-2]
			\arrow[from=2-3, to=3-3]
			\arrow[from=3-1, to=3-2]
			\arrow[from=3-2, to=3-3]
		\end{tikzcd}
	\end{equation}
	Each inner square commutes by a suitable application of \Cref{lem:mates_bijection1} to each case.
	The claim follows by noting that in \eqref{diag:conj_is_mates_twice}:
	\begin{enumerate}
		\item Going from the top left corner to the right bottom one is `taking the right conjugate,' see \cite[Proposition 4.3.3]{Riehl:26}.
		\item Going from the top left corner to the central one is `taking the right mate with respect to $p^{\prime\ast}\dashv p^\prime_\ast$ and $p^\ast\dashv p_\ast$.'
		\item Going from the central object to the right bottom one is `taking the right mate with respect to $q^\ast\dashv q_\ast$ and $q^{\prime\ast}\dashv q^\prime_\ast$.'
		\qedhere
	\end{enumerate}
\end{proof}

We now explain how the mate correspondence is compatible with both horizontal and vertical composition.
More precisely, \Cref{lem:mates_bijection} upgrades to an isomorphism of double categories \Cref{lem:mates_double_cat_iso}.
We redirect the reader to \cite[\S 1.1]{Kelly/Street:1974} for the precise definition of a double category.
Roughly speaking, the data of a double category is made up of $0$-cells (its objects), two kinds of $1$-morphisms (horizontal and vertical ones) and $2$-cells
\begin{displaymath}
	% https://q.uiver.app/#q=WzAsNCxbMCwwLCJBIl0sWzEsMCwiQiJdLFswLDEsIkEnIl0sWzEsMSwiQiciXSxbMCwxLCJmIl0sWzIsMywiZiciLDJdLFswLDIsImEiLDIseyJjb2xvdXIiOlswLDYwLDYwXX0sWzAsNjAsNjAsMV1dLFsxLDMsImIiLDAseyJjb2xvdXIiOlswLDYwLDYwXX0sWzAsNjAsNjAsMV1dXQ==
	\begin{tikzcd}
		A & B \\
		{A^\prime} & {B^\prime}
		\arrow["f", from=1-1, to=1-2]
		\arrow["a"'{inner sep=.8ex}, "\shortmid"{marking}, from=1-1, to=2-1]
		\arrow["b"{inner sep=.8ex}, "\shortmid"{marking}, from=1-2, to=2-2]
		\arrow["{f^\prime}"', from=2-1, to=2-2]
	\end{tikzcd}
\end{displaymath}
involving two horizontal $1$-morphisms $f,f^\prime$ and two vertical $1$-morphisms $a,b$.\footnote{Typically, in a double category an horizontal arrow is depicted as $\to$,
whereas vertical arrows are depicted with some additional marker to distinguish them from the vertical ones.
Here we opted for $\slashedrightarrow$.}
The axioms of double categories allow to compose $2$-cells both vertically and horizontally; in fact, in a compatible way.

\begin{definition}[{{\cite[Definition 1.1]{Cheng/Gurski/Riehl:2014}}}]
	\label{def:mates_double_cat}
	Let $\mathcal{T}$ be a $2$-category.
	We have a pair of double categories $\mathbb{L}\mathbf{Adj}_{\mathcal{T}}$ and $\mathbb{R}\mathbf{Adj}_{\mathcal{T}}$ with same $0$-cells and $1$-cells but different $2$-cells.
	Namely, the objects and horizontal $1$-morphisms of $\mathbb{L}\mathbf{Adj}_{\mathcal{T}}$ and $\mathbb{R}\mathbf{Adj}_{\mathcal{T}}$ coincide with the objects and $1$-morphisms of $\mathcal{T}$, whereas a vertical morphism in $\mathbb{L}\mathbf{Adj}_{\mathcal{T}}$ and $\mathbb{R}\mathbf{Adj}_{\mathcal{T}}$ is an adjunction
	\begin{displaymath}
		% https://q.uiver.app/#q=WzAsMixbMCwwLCJcXG1hdGhjYWx7Q30iXSxbMCwxLCJcXG1hdGhjYWx7RH0iXSxbMCwxLCJGIiwyLHsib2Zmc2V0IjoyfV0sWzEsMCwiRyIsMix7Im9mZnNldCI6Mn1dLFsyLDMsIiIsMix7ImxldmVsIjoxLCJzdHlsZSI6eyJuYW1lIjoiYWRqdW5jdGlvbiJ9fV1d
		\begin{tikzcd}
			{\mathcal{C}} \\ {\mathcal{D}}
			\arrow[""{name=0, anchor=center, inner sep=0}, "F"', shift left=-2, from=1-1, to=2-1]
			\arrow[""{name=1, anchor=center, inner sep=0}, "G"', shift left=-2, from=2-1, to=1-1]
			\arrow["\dashv"{anchor=center}, draw=none, from=0, to=1]
		\end{tikzcd}
	\end{displaymath}
	in $\mathcal{T}$.
	Finally, a $2$-cell in $\mathbb{L}\mathbf{Adj}_{\mathcal{T}}$ and in $\mathbb{R}\mathbf{Adj}_{\mathcal{T}}$ is respectively a $2$-cell
	\begin{displaymath}
		%https://q.uiver.app/#q=WzAsOCxbMCwwLCJcXG1hdGhjYWx7Q30iXSxbMCwxLCJcXG1hdGhjYWx7RH0iXSxbMSwwLCJcXG1hdGhjYWx7Q31eXFxwcmltZSJdLFsxLDEsIlxcbWF0aGNhbHtEfV5cXHByaW1lIl0sWzQsMCwiXFxtYXRoY2Fse0N9Il0sWzUsMCwiXFxtYXRoY2Fse0N9XlxccHJpbWUiXSxbNCwxLCJcXG1hdGhjYWx7RH0iXSxbNSwxLCJcXG1hdGhjYWx7RH1eXFxwcmltZSJdLFswLDEsIkYiLDJdLFswLDIsIkgiXSxbMSwzLCJLIiwyXSxbMiwzLCJGJyJdLFswLDMsIlxcYWxwaGEiLDAseyJzaG9ydGVuIjp7InNvdXJjZSI6MjAsInRhcmdldCI6MjB9LCJsZXZlbCI6Mn1dLFs2LDcsIksiLDJdLFs0LDUsIkgiXSxbNiw0LCJHIl0sWzcsNSwiRyciLDJdLFs3LDQsIlxcYmV0YSIsMix7InNob3J0ZW4iOnsic291cmNlIjoyMCwidGFyZ2V0IjoyMH0sImxldmVsIjoyfV0sWzExLDE1LCJcXHRleHR7XFxub3JtYWxzaXplIGFuZH0iLDEseyJzaG9ydGVuIjp7InNvdXJjZSI6MjAsInRhcmdldCI6MjB9LCJzdHlsZSI6eyJib2R5Ijp7Im5hbWUiOiJub25lIn0sImhlYWQiOnsibmFtZSI6Im5vbmUifX19XV0=
		\begin{tikzcd}
			{\mathcal{C}} & {\mathcal{C}^\prime} && {\mathcal{C}} & {\mathcal{C}^\prime} \\
			{\mathcal{D}} & {\mathcal{D}^\prime} && {\mathcal{D}} & {\mathcal{D}^\prime}
			\arrow["H", from=1-1, to=1-2]
			\arrow["F"', from=1-1, to=2-1]
			\arrow["\alpha"', shorten <=2pt, shorten >=2pt, Rightarrow, from=1-2, to=2-1]
			\arrow[""{name=0, anchor=center, inner sep=0}, "{F^\prime}", from=1-2, to=2-2]
			\arrow["H", from=1-4, to=1-5]
			\arrow["\beta", shorten <=2pt, shorten >=2pt, Rightarrow, from=1-4, to=2-5]
			\arrow["K"', from=2-1, to=2-2]
			\arrow[""{name=1, anchor=center, inner sep=0}, "G", from=2-4, to=1-4]
			\arrow["K"', from=2-4, to=2-5]
			\arrow["{G^\prime}"', from=2-5, to=1-5]
			\arrow["{\textrm{\normalsize and}}"{description}, draw=none, from=0, to=1]
		\end{tikzcd}
	\end{displaymath}
	in $\mathcal{T}$.
	The horizontal and vertical composition laws of $2$-cells in $\mathbb{L}\mathbf{Adj}_{\mathcal{T}}$ and $\mathbb{R}\mathbf{Adj}_{\mathcal{T}}$
	are given by horizontal and vertical pasting of $2$-cells in $\mathcal{T}$.
\end{definition}

\begin{lemma}
	[{{\cite[Proposition 2.2]{Kelly/Street:1974}}}]
	\label{lem:mates_double_cat_iso}
	The mate bijection of \Cref{lem:mates_bijection} induces an isomorphism $\mathbb{L}\mathbf{Adj}_{\mathcal{T}}\cong\mathbb{R}\mathbf{Adj}_{\mathcal{T}}$ of double categories.\footnote{A \textbf{double functor} between double categories consists of assignments on objects, on horizontal morphisms, on vertical morphisms, and on
	squares, which are compatible with domains and codomains and preserve
	all compositions and identities strictly
	\cite[Definition 2.5]{Lyne:2026}.
	Double categories and double functors assemble into a $1$-category.
	Hence, the notion of isomorphism of double categories is defined.}
\end{lemma}

\begin{proof}
	The content of this lemma is that `the mate correspondence \eqref{eq:mate_correspondence_2-cell} commutes with horizontal and vertical composition of $2$-cells.'
	Namely, we define a double functor
	$\mathbb{L}\mathbf{Adj}_{\mathcal{T}}\to\mathbb{R}\mathbf{Adj}_{\mathcal{T}}$
	that acts as the identity on objects and vertical and horizontal $1$-morphisms,
	and that acts on $2$-cells as the mate correspondence of \Cref{lem:mates_bijection}.
	Since every assignment is bijective,
	it will be an isomorphism if it preserves horizontal and vertical composition.
	The explicit compatibility formulas are spelled out in \cite[pp. 344-345]{Cheng/Gurski/Riehl:2014}.
	Their proof is a $2$-categorical exercise whose details are in \cite[Proposition 2.2]{Kelly/Street:1974}.
\end{proof}

\begin{lemma}[{{functoriality of conjugacy \cite[Exercise 3.3.7 (c)]{Lipman:2009}}}]
	\label[lemma]{lem:conjugates_are_functorial}
	With the notation of \Cref{lem:mates_bijection}:
	For conjugates, the formulas $\overline{\operatorname{id}}=\operatorname{id}$ and
	$\overline{\alpha\circ\alpha^\prime}=\overline{\alpha^\prime}\circ\overline{\alpha}$ hold.
	Explicitly, the second one means that if
	\begin{displaymath}
		\begin{gathered}
			(\alpha\colon F^\prime\to F,
			\beta\colon G\to G^\prime),\\
			(\alpha^\prime\colon F^{\prime\prime}\to F^\prime,
			\beta^\prime\colon G^\prime\to G^{\prime\prime}),
		\end{gathered}
	\end{displaymath}
	are pairs of conjugates,
	then so is
	$(\alpha\circ\alpha^\prime\colon F^{\prime\prime}\to F,\beta^\prime\circ\beta\colon G\to G^{\prime\prime})$.
\end{lemma}

\begin{proof}
	The fact that the identity is auto-conjugate is trivial,
	whereas the second formula follows by considering horizontal composition and \eqref{lem:mates_double_cat_iso} for the following diagrams:
	\begin{displaymath}
		% https://q.uiver.app/#q=WzAsMTIsWzAsMCwiXFxtYXRoY2Fse0N9Il0sWzEsMCwiXFxtYXRoY2Fse0N9Il0sWzAsMSwiXFxtYXRoY2Fse0R9Il0sWzEsMSwiXFxtYXRoY2Fse0R9Il0sWzIsMCwiXFxtYXRoY2Fse0N9Il0sWzIsMSwiXFxtYXRoY2Fse0R9Il0sWzQsMCwiXFxtYXRoY2Fse0N9Il0sWzUsMCwiXFxtYXRoY2Fse0N9Il0sWzQsMSwiXFxtYXRoY2Fse0R9Il0sWzUsMSwiXFxtYXRoY2Fse0R9Il0sWzYsMCwiXFxtYXRoY2Fse0N9Il0sWzYsMSwiXFxtYXRoY2Fse0R9Il0sWzAsMSwiIiwwLHsibGV2ZWwiOjIsInN0eWxlIjp7ImhlYWQiOnsibmFtZSI6Im5vbmUifX19XSxbMiwzLCIiLDAseyJsZXZlbCI6Miwic3R5bGUiOnsiaGVhZCI6eyJuYW1lIjoibm9uZSJ9fX1dLFswLDIsIkYiLDJdLFsxLDMsIkZeXFxwcmltZSIsMV0sWzEsNCwiIiwxLHsibGV2ZWwiOjIsInN0eWxlIjp7ImhlYWQiOnsibmFtZSI6Im5vbmUifX19XSxbMyw1LCIiLDEseyJsZXZlbCI6Miwic3R5bGUiOnsiaGVhZCI6eyJuYW1lIjoibm9uZSJ9fX1dLFs0LDUsIkZee1xccHJpbWVcXHByaW1lfSJdLFsxLDIsIlxcYWxwaGEiLDIseyJsZXZlbCI6Mn1dLFs0LDMsIlxcYWxwaGFeXFxwcmltZSIsMix7ImxldmVsIjoyfV0sWzYsNywiIiwwLHsibGV2ZWwiOjIsInN0eWxlIjp7ImhlYWQiOnsibmFtZSI6Im5vbmUifX19XSxbOCw5LCIiLDAseyJsZXZlbCI6Miwic3R5bGUiOnsiaGVhZCI6eyJuYW1lIjoibm9uZSJ9fX1dLFs4LDYsIkciXSxbOSw3LCJHXlxccHJpbWUiLDFdLFs3LDEwLCIiLDEseyJsZXZlbCI6Miwic3R5bGUiOnsiaGVhZCI6eyJuYW1lIjoibm9uZSJ9fX1dLFs5LDExLCIiLDEseyJsZXZlbCI6Miwic3R5bGUiOnsiaGVhZCI6eyJuYW1lIjoibm9uZSJ9fX1dLFsxMSwxMCwiR157XFxwcmltZVxccHJpbWV9IiwyXSxbNiw5LCJcXGJldGEiLDAseyJsZXZlbCI6Mn1dLFs3LDExLCJcXGJldGFeXFxwcmltZSIsMCx7ImxldmVsIjoyfV0sWzE4LDIzLCJcXGxlZnRyaWdodHNxdWlnYXJyb3ciLDEseyJzaG9ydGVuIjp7InNvdXJjZSI6MjAsInRhcmdldCI6MjB9LCJzdHlsZSI6eyJib2R5Ijp7Im5hbWUiOiJub25lIn0sImhlYWQiOnsibmFtZSI6Im5vbmUifX19XV0=
		\begin{tikzcd}[/tikz/baseline=(tikz@f@1-2-1.base)]
			{\mathcal{C}} & {\mathcal{C}} & {\mathcal{C}} && {\mathcal{C}} & {\mathcal{C}} & {\mathcal{C}} \\
			{\mathcal{D}} & {\mathcal{D}} & {\mathcal{D}} && {\mathcal{D}} & {\mathcal{D}} & {\mathcal{D}}
			\arrow[equals, from=1-1, to=1-2]
			\arrow["F"', from=1-1, to=2-1]
			\arrow[equals, from=1-2, to=1-3]
			\arrow["\alpha"', shorten <=2pt, shorten >=2pt, Rightarrow, from=1-2, to=2-1]
			\arrow["{F^\prime}"{description}, from=1-2, to=2-2]
			\arrow["{\alpha^\prime}"', shorten <=2pt, shorten >=2pt, Rightarrow, from=1-3, to=2-2]
			\arrow[""{name=0, anchor=center, inner sep=0}, "{F^{\prime\prime}}", from=1-3, to=2-3]
			\arrow[equals, from=1-5, to=1-6]
			\arrow["\beta", shorten <=2pt, shorten >=2pt, Rightarrow, from=1-5, to=2-6]
			\arrow[equals, from=1-6, to=1-7]
			\arrow["{\beta^\prime}", shorten <=2pt, shorten >=2pt, Rightarrow, from=1-6, to=2-7]
			\arrow[equals, from=2-1, to=2-2]
			\arrow[equals, from=2-2, to=2-3]
			\arrow[""{name=1, anchor=center, inner sep=0}, "G", from=2-5, to=1-5]
			\arrow[equals, from=2-5, to=2-6]
			\arrow["{G^\prime}"{description}, from=2-6, to=1-6]
			\arrow[equals, from=2-6, to=2-7]
			\arrow["{G^{\prime\prime}}"', from=2-7, to=1-7]
			\arrow["\leftrightsquigarrow"{description}, draw=none, from=0, to=1]
		\end{tikzcd}
		\qedhere
	\end{displaymath}
\end{proof}

\begin{lemma-definition}
	\label{def:adj_equiv}
	An \textbf{adjoint equivalence} in a $2$-category $\mathcal{T}$ is an adjunction $(F,G,\eta,\epsilon)$ such that any of the equivalent conditions hold:
	\begin{enumerate}
		\item \label{cond:left_adj_is_equiv}
		$F$ is an equivalence.
		\item \label{cond:right_adj_is_equiv}
		$G$ is an equivalence.
		\item \label{cond:left_and_right_adj_are_rep_ff}
		$F$ and $G$ are representably fully faithful.\footnote{A $1$-morphism $F\colon \mathcal{C}\to\mathcal{C}^\prime$ in $\mathcal{T}$ is said to be \textbf{representably fully faithful}
		if $F_\ast\colon \mathcal{T}(\mathcal{D},\mathcal{C})\to\mathcal{T}(\mathcal{D},\mathcal{C}^\prime)$
		is a fully faithful functor for every object $\mathcal{D}\in\mathcal{T}$.
		In the case $\mathcal{T}=\mathbf{Cat}$ is the $2$-category of $1$-categories, then representable full faithfulness amounts to ordinary full faithfulness:
		On the one hand, the representable notion implies the ordinary one by considering $\mathcal{D}$ to be the category with a single object and a single morphism.
		Conversely, given functors $G_1,G_2\colon \mathcal{D}\rightrightarrows\mathcal{C}$,
		any natural transformation $FG_1\Rightarrow FG_2$ arises from a unique natural transformation $G_1\Rightarrow G_2$ if $F$ is fully faithful.}
		\item \label{cond:unit_counit_are_isos}
		$\eta$ and $\epsilon$ are isomorphisms.
	\end{enumerate}
\end{lemma-definition}

\begin{proof}
	It is clear that \eqref{cond:unit_counit_are_isos} implies \eqref{cond:left_adj_is_equiv} and \eqref{cond:right_adj_is_equiv}.
	The converses follow from \cite[Propositions 4.3.1, Corollary 4.3.6 and Remark 4.3.9]{Riehl:26}.
	The equivalence \eqref{cond:left_and_right_adj_are_rep_ff} $\iff$ \eqref{cond:unit_counit_are_isos} comes from the fact that $F$ (resp.\ $G$) is representably fully faithful if, and only if,
	$\eta$ (resp.\ $\epsilon$) is an isomorphism:
	Let us give the argument for the left adjoint.
	Given an object $\mathcal{D}\in\mathcal{T}$, since $2$-functors preserve adjunctions,
	we have an adjunction $F^\mathcal{D}\dashv G^\mathcal{D}$ in $\mathbf{Cat}$,
	where $H^\mathcal{D}$ means the image of a $1$-morphism $H$ in $\mathcal{T}$ under the $2$-functor $(-)^\mathcal{D}=\mathcal{T}(\mathcal{D},-)\colon \mathcal{T}\to\mathbf{Cat}$.
	Thus, if $\eta$ is an isomorphism, $\eta^\mathcal{D}$ is also an isomorphism and hence $F^\mathcal{D}$ is fully faithful \cite[\href{https://stacks.math.columbia.edu/tag/07RB}{Tag 07RB}]{StacksProject}.
	This means $F$ is representably fully faithful.
	The converse follows by reversing the argument using the $2$-Yoneda embedding,
	or by the argument in \cite[Proposition 2.19 (2)]{Weber:2007}.
\end{proof}
% In the case $\mathcal{T}=\mathbf{Cat}$ is the $2$-category of $1$-categories,
% then we might also replace `equivalence' by `fully faithful' in \eqref{cond:left_adj_is_equiv} and \eqref{cond:right_adj_is_equiv},
% see \cite[\href{https://stacks.math.columbia.edu/tag/07RB}{Tag 07RB}]{StacksProject}
% or \cite[Proposition 1.5.6]{Kashiwara/Schapira:2006}.

Any equivalence can be upgraded to an adjoint equivalence \cite[Corollary 4.3.6, Remark 4.3.9]{Riehl:26}.

\begin{lemma}
	\label{lem:when_mate_correspondence_preserves_isos}
	Assume \Cref{setup:mates}.
	Let $\alpha\colon F^\prime H\to KF$ and $\beta\colon HG\to G^\prime K$ be a pair of mates under \eqref{diag:morph_adj}.
	It holds that $\alpha$ is an isomorphism if, and only if, $\beta$ is an isomorphism in any of the following situations:
	\begin{enumerate}[label=(\arabic*),ref=\thelemma(\arabic*)]
		\item \label[lemma]{lem:when_mate_correspondence_preserves_isos1}
		\cite[Exercise 4.3.v]{Riehl:26}
		Conjugacy, i.e.\ $H=1_{\mathcal{C}}$ and $K=1_{\mathcal{D}}$.
		Moreover, in such case, $\overline{\alpha^{-1}}=\overline{\alpha}\,^{-1}$ (with the notation of \Cref{lem:mates_bijection}).

		\item \label[lemma]{lem:when_mate_correspondence_preserves_isos2}
		$(F,G,\eta,\epsilon)$ and $(F^\prime,G^\prime,\eta^\prime,\epsilon^\prime)$ are adjoint equivalences
		(see \Cref{def:adj_equiv}).
		\item \label[lemma]{lem:when_mate_correspondence_preserves_isos3}
		$H$ and $K$ are equivalences.
	\end{enumerate}
\end{lemma}

\begin{proof}
	\ref{lem:when_mate_correspondence_preserves_isos1} follows from the fact that conjugacy is functorial, \Cref{lem:conjugates_are_functorial}.
	On the other hand, \ref{lem:when_mate_correspondence_preserves_isos2} follows from the fact that the mate formulas in \Cref{lem:mates_bijection2} preserve $2$-isomorphisms: $\eta,\eta^\prime,\epsilon,\epsilon^\prime$ are $2$-isomorphisms and whiskering preserves $2$-isomorphisms.

	Let us show \ref{lem:when_mate_correspondence_preserves_isos3}.
	The $1$-morphisms $H$ and $K$ have a left adjoint quasi-inverse \cite[Corollary 4.3.6, Remark 4.3.9]{Riehl:26}.
	Hence, we may apply \Cref{lem:taking_conjugates_is_taking_mates_twice}
	by instantiating the square \eqref{diag:square_conj_is_mates_twice} as the square:
	\begin{displaymath}
		% https://q.uiver.app/#q=WzAsNCxbMCwwLCJcXG1hdGhjYWx7RH0iXSxbMCwxLCJcXG1hdGhjYWx7RH1eXFxwcmltZSJdLFsxLDEsIlxcbWF0aGNhbHtDfV5cXHByaW1lIl0sWzEsMCwiXFxtYXRoY2Fse0N9Il0sWzAsMywiRyJdLFsxLDIsIkdeXFxwcmltZSIsMl0sWzAsMSwiSyIsMl0sWzMsMiwiSCJdXQ==
		\begin{tikzcd}[cramped]
			{\mathcal{D}} & {\mathcal{C}} \\
			{\mathcal{D}^\prime} & {\mathcal{C}^\prime\mathrlap{.}}
			\arrow["G", from=1-1, to=1-2]
			\arrow["K"', from=1-1, to=2-1]
			\arrow["H", from=1-2, to=2-2]
			\arrow["{G^\prime}"', from=2-1, to=2-2]
		\end{tikzcd}
	\end{displaymath}
	The claim is that the bottom right bijection in \eqref{diag:triangle_conj_mates} preserves the property of being an isomorphism.
	It thus suffice to show that the other two bijections in \eqref{diag:triangle_conj_mates} preserve such property.
	This is \ref{lem:when_mate_correspondence_preserves_isos1} for the top bijection, and \ref{lem:when_mate_correspondence_preserves_isos2} for the bottom left bijection.
\end{proof}

\subsection{Base change morphisms}
\label{sec:base_change_morphisms}

Mates arise naturally within the operations of the six functor formalism.
For instance \cite[Example 4.5]{Bergh/Schnurer:2020},
if $f\colon X\to Y$ is a morphism of algebraic spaces and given $E\in D_{\operatorname{qc}}(Y)$, the natural transformations
\begin{displaymath}
	\begin{gathered}
		\mathbf{L}f^\ast E\otimes^{\mathbf{L}}f^\ast(-)\to
		f^\ast(E\otimes^{\mathbf{L}}-),\\
		f^\ast\operatorname{\mathbf{R}\mathcal{H}\! \mathit{om}}(E,-)\to
		\operatorname{\mathbf{R}\mathcal{H}\! \mathit{om}}(\mathbf{L}f^\ast E,\mathbf{L}f^\ast(-)),
	\end{gathered}
\end{displaymath}
are mates with respect to the adjunctions $E\otimes^{\mathbf{L}}(-)\dashv\operatorname{\mathbf{R}\mathcal{H}\! \mathit{om}}(E,-)$ and $\mathbf{L}f^\ast E\otimes^{\mathbf{L}}(-)\dashv\operatorname{\mathbf{R}\mathcal{H}\! \mathit{om}}(\mathbf{L}f^\ast E,-)$.

There are three different canonical maps regarded as `base change morphisms' in the theory of algebraic spaces, named as $^\ast$-$_\ast$, $_\ast$-$^\times$, and $^\ast$-$^\times$ base change morphism.\footnote{We read
$^\ast$-$_\ast$, $_\ast$-$^\times$, and $^\ast$-$^\times$ as `upper star-lower star,' `lower star-upper cross,' and `upper star-upper cross.'}
All three arise as mates.
In this section, we will explain how they do arise and prove functoriality with respect to pasting of commutative squares, see \Cref{lem:base_change_mates_pasting,lem:tor-independence_and_beta_pasting}.
Throughout this section, we fix a commutative square of quasi-compact and quasi-separated algebraic spaces:
\begin{equation}
	\label{diag:abstract_base_change}
	% https://q.uiver.app/#q=WzAsNCxbMCwwLCJYXlxccHJpbWUiXSxbMSwwLCJYIl0sWzEsMSwiWi4iXSxbMCwxLCJaXlxccHJpbWUiXSxbMSwyLCJwIl0sWzAsMSwicV5cXHByaW1lIl0sWzMsMiwicSIsMl0sWzAsMywicF5cXHByaW1lIiwyXV0=
	\begin{tikzcd}
		X^\prime
		\arrow[r,"q^\prime"]
		\arrow[d,"p^\prime"']
		& X
		\arrow[d,"p"]
		\\
		Z^\prime
		\arrow[r,"q"']
		& Z .
	\end{tikzcd}
\end{equation}

In the following lemmas and definitions, the decoration $(-)_\square$ means we are referring to this commutative square.

\begin{lemma}
	\label{lem:pseudofunctoriality_isos_are_conjugate}
	The pseudofunctoriality isomorphisms
	\begin{displaymath}
		\begin{gathered}
			\lambda_\square\colon \mathbf{L}(p^\prime)^\ast\mathbf{L}q^\ast
			\xrightarrow{\cong}\mathbf{L}(q^\prime)^\ast\mathbf{L}p^\ast,\\
			\rho_\square\colon \mathbf{R}p_\ast\mathbf{R}q^\prime_\ast
			\xrightarrow{\cong}\mathbf{R}q_\ast\mathbf{R}p^\prime_\ast,
		\end{gathered}
	\end{displaymath}
	are conjugate to each other (see \cite[Proposition 4.3.3]{Riehl:26}).
\end{lemma}

\begin{proof}
	We have two pairs of conjugate isomorphisms, respectively for the clockwise and counterclockwise composite in \eqref{diag:abstract_base_change}:
	\begin{displaymath}
		\begin{gathered}
			(\mathbf{L}(pq^\prime)^\ast\to\mathbf{L}(q^\prime)^\ast\mathbf{L}p^\ast,
			\mathbf{R}p_\ast\mathbf{R}q^\prime_\ast\to\mathbf{R}(pq^\prime)_\ast
			),
			\\
			(\mathbf{L}(p^\prime)^\ast\mathbf{L}q^\ast\to\mathbf{L}(qp^\prime)^\ast,
			\mathbf{R}(qp^\prime)_\ast\to\mathbf{R}q_\ast\mathbf{R}p^\prime_\ast).
		\end{gathered}
	\end{displaymath}
	Indeed, the fact that the second pair are conjugate transformations is exactly \cite[(3.6.2)]{Lipman:2009} coupled with \Cref{lem:morph_adj3} (or \cite[(3.6.2)\textsuperscript{\textsf{op}}]{Lipman:2009} coupled with \Cref{lem:morph_adj2}),
	thanks to \Cref{lem:spacey_Lipman_Hashimoto}.
	The fact that the first pair are conjugates is by same reasons combined with \Cref{lem:when_mate_correspondence_preserves_isos1}.
	The claim follows by composing these two pairs of conjugates (since $pq^\prime=qp^\prime$), which is possible by \Cref{lem:conjugates_are_functorial}.
\end{proof}

\begin{lemma-definition}[{{\cite[Definition and Remark 21.129]{Gortz/Wedhorn:2023}}}]
	\label[definition]{def:upper_star_lower_star_bc}
	The \textbf{$^\ast$-$_\ast$ base change morphism}
	\begin{equation}
		\label{eq:upper_star_lower_star}
		\vartheta_\square\colon
		\mathbf{L}q^\ast\mathbf{R}p_\ast
		\to \mathbf{R}p^\prime_\ast \mathbf{L}(q^\prime)^\ast
	\end{equation}
	is defined by any of the following descriptions, giving rise to the same morphism:
	\begin{enumerate}[label=(\arabic*), ref=\thedefinition(\arabic*)]
		\item
		\label[definition]{def:left_upper_star_lower_star_bc}
		The left mate of the pseudofunctoriality isomorphism
		\begin{displaymath}
			\rho_\square\colon \mathbf{R}p_\ast\mathbf{R}q^\prime_\ast
			\xrightarrow{\cong}\mathbf{R}q_\ast\mathbf{R}p^\prime_\ast
		\end{displaymath}
		with respect to the adjunctions $\mathbf{L}(q^\prime)^\ast\dashv\mathbf{R}q^\prime_\ast$ and $\mathbf{L}q^\ast\dashv\mathbf{R}q_\ast$.
		\item
		\label[definition]{def:right_upper_star_lower_star_bc}
		The right mate of the pseudofunctoriality isomorphism
		\begin{displaymath}
			\lambda_\square\colon \mathbf{L}(p^\prime)^\ast\mathbf{L}q^\ast
			\xrightarrow{\cong}\mathbf{L}(q^\prime)^\ast\mathbf{L}p^\ast
		\end{displaymath}
		with respect to the adjunctions $\mathbf{L}(p^\prime)^\ast\dashv\mathbf{R}p^\prime_\ast$ and $\mathbf{L}p^\ast\dashv\mathbf{R}p_\ast$.
	\end{enumerate}
\end{lemma-definition}

\begin{proof}
	Using \Cref{lem:pseudofunctoriality_isos_are_conjugate},
	this is a straight application of \Cref{lem:taking_conjugates_is_taking_mates_twice}.
	See also \cite[(3.7.2)]{Lipman:2009} or \cite[Proposition 2.1.3]{Deligne:1971}.\footnote{\label{fn:lipman_is_formal}These formal results apply in the algebraic space situation thanks to \Cref{lem:spacey_Lipman_Hashimoto}.}
\end{proof}

The functors in the domain and codomain of $\vartheta_\square$ have right adjoints, namely $\mathbf{L}q^\ast\mathbf{R}p_\ast \dashv p^\times\mathbf{R}q_\ast$ and $\mathbf{R}p^\prime_\ast \mathbf{L}(q^\prime)^\ast \dashv \mathbf{R}q^\prime_\ast (p^\prime)^\times$ \cite[Proposition 4.3.3]{Riehl:26}.
Indeed, for $M\in D_{\operatorname{qc}}(X)$ and $N\in D_{\operatorname{qc}}(Z^\prime)$, there are natural isomorphisms
\begin{displaymath}
	\begin{aligned}
		\operatorname{Hom}
		(\mathbf{L}q^\ast\mathbf{R}p_\ast M,N)
		&\cong \operatorname{Hom} (\mathbf{R}p_\ast M,\mathbf{R}q_\ast N)
		\\ &\cong \operatorname{Hom} (M,p^\times\mathbf{R}q_\ast N),
	\end{aligned}
\end{displaymath}
whereas
\begin{displaymath}
	\begin{aligned}
		\operatorname{Hom}
		(
		\mathbf{R}p^\prime_\ast
		\mathbf{L}(q^\prime)^\ast M,
		N
		)
		&\cong
		\operatorname{Hom}
		(
		\mathbf{L}(q^\prime)^\ast M,
		(p^\prime)^\times N
		)
		\\
		&\cong
		\operatorname{Hom}
		(
		M,
		\mathbf{R}q^\prime_\ast
		(p^\prime)^\times N
		).
	\end{aligned}
\end{displaymath}

\begin{definition}
	\label{def:lower_star_upper_cross_bc}
	The \textbf{$_\ast$-$^\times$ base change morphism}
	\begin{equation}
		\label{eq:lower_star_upper_cross}
		\xi_\square\colon
		\mathbf{R}q^\prime_\ast(p^\prime)^\times
		\to p^\times\mathbf{R}q_\ast
	\end{equation}
	is the right conjugate of $\vartheta_\square$.
\end{definition}

Therefore, by \Cref{lem:when_mate_correspondence_preserves_isos1}, $\vartheta_\square$ is an isomorphism if, and only if, $\xi_\square$ is an isomorphism.
This happens for instance if \eqref{diag:abstract_base_change} is a tor-independent fibered square \cite[\href{https://stacks.math.columbia.edu/tag/08IR}{Tag 08IR}]{StacksProject}.
Later we will show that in fact this is the only situation in which this happens, see \Cref{lem:characterize_tor-independence}.

\begin{lemma-definition}
	\label[definition]{def:upper_star_upper_cross_bc}
	Assume $\vartheta_\square$ is an isomorphism (equivalently,  $\xi_\square$ is an isomorphism).
	The \textbf{$^\ast$-$^\times$ base change morphism}
	\begin{equation}
		\label{eq:upper_star_upper_cross}
		\beta_\square\colon
		\mathbf{L}(q^\prime)^\ast p^\times
		\to (p^\prime)^\times\mathbf{L}q^\ast
	\end{equation}
	is defined by any of the following equivalent descriptions, giving rise to the same morphism:
	\begin{enumerate}[label=(\arabic*), ref=\thelemma(\arabic*)]
		\item
		\label[lemma-definition]{def:upper_star_upper_cross_bc1}
		The right mate of
		\begin{displaymath}
			\vartheta_\square^{-1}\colon
			\mathbf{R}p^\prime_\ast \mathbf{L}(q^\prime)^\ast
			\to \mathbf{L}q^\ast\mathbf{R}p_\ast
		\end{displaymath}
		with respect to the adjunctions $\mathbf{R}p_\ast\dashv p^\times$ and $\mathbf{R}p_\ast^\prime\dashv(p^\prime)^\times$.
		\item
		\label[lemma-definition]{def:upper_star_upper_cross_bc2}
		The left mate of
		\begin{displaymath}
			\xi_\square^{-1}\colon
			p^\times\mathbf{R}q_\ast
			\to \mathbf{R}q^\prime_\ast(p^\prime)^\times
		\end{displaymath}
		with respect to the adjunctions $\mathbf{L}q^\ast\dashv\mathbf{R}q_\ast$ and $\mathbf{L}(q^\prime)^\ast \dashv \mathbf{R}q^\prime_\ast$.
	\end{enumerate}
\end{lemma-definition}

\begin{proof}
	Using \Cref{lem:when_mate_correspondence_preserves_isos1},
	this is a straight application of \Cref{lem:taking_conjugates_is_taking_mates_twice}
	applied to the square:
	\begin{displaymath}
		% https://q.uiver.app/#q=WzAsNCxbMCwwLCJEX3tcXG9wZXJhdG9ybmFtZXtxY319KFpeXFxwcmltZSkiXSxbMCwxLCJEX3tcXG9wZXJhdG9ybmFtZXtxY319KFheXFxwcmltZSkiXSxbMSwwLCJEX3tcXG9wZXJhdG9ybmFtZXtxY319KFopIl0sWzEsMSwiRF97XFxvcGVyYXRvcm5hbWV7cWN9fShYKSJdLFswLDIsIlxcbWF0aGJme1J9cV9cXGFzdCJdLFsyLDMsInBeXFx0aW1lcyJdLFswLDEsIihwXlxccHJpbWUpXlxcdGltZXMiLDJdLFsxLDMsIlxcbWF0aGJme1J9cV5cXHByaW1lX1xcYXN0IiwyXV0=
		\begin{tikzcd}[/tikz/baseline=(tikz@f@1-2-1.base)]
			{D_{\operatorname{qc}}(Z^\prime)} & {D_{\operatorname{qc}}(Z)} \\
			{D_{\operatorname{qc}}(X^\prime)} & {D_{\operatorname{qc}}(X)}
			\arrow["{\mathbf{R}q_\ast}", from=1-1, to=1-2]
			\arrow["{(p^\prime)^\times}"', from=1-1, to=2-1]
			\arrow["{p^\times}", from=1-2, to=2-2]
			\arrow["{\mathbf{R}q^\prime_\ast}"', from=2-1, to=2-2]
		\end{tikzcd}
		\qedhere
	\end{displaymath}
\end{proof}

\begin{remark}
	In \cite[Example 4.4]{Bergh/Schnurer:2020} and \cite[\href{https://stacks.math.columbia.edu/tag/0E5D}{Tag 0E5D}]{StacksProject} we find two different definitions for \eqref{eq:upper_star_upper_cross}.
	\Cref{def:upper_star_upper_cross_bc} shows these two are the same.
	Indeed, \ref{def:upper_star_upper_cross_bc1} is exactly the definition of the $^\ast$-$^\times$ base change morphism given in \cite[Example 4.4]{Bergh/Schnurer:2020}.
	On the other hand,
	\ref{def:upper_star_upper_cross_bc2} made out explicit is the composite
	\begin{equation}
		\label{eq:base_change_mate_formula}
		\begin{aligned}
			\mathbf{L}(q^\prime)^\ast p^\times
			&\xrightarrow{ \mathbf{L}(q^\prime)^\ast p^\times(\eta_q)}
			\mathbf{L}(q^\prime)^\ast p^\times\mathbf{R}q_\ast\mathbf{L}q^\ast
			\\&\xrightarrow{\mathbf{L}(q^\prime)^\ast(\xi_\square^{-1}) \mathbf{L}q^\ast} \mathbf{L}(q^\prime)^\ast \mathbf{R}q^\prime_\ast (p^\prime)^\times \mathbf{L}q^\ast
			\\&\xrightarrow{ \epsilon_{q^\prime} (p^\prime)^\times \mathbf{L}q^\ast}
			(p^\prime)^\times\mathbf{L}q^\ast,
		\end{aligned}
	\end{equation}
	where $\eta_q\colon 1\to\mathbf{R}q_\ast\mathbf{L}q^\ast$ is the unit of
	$\mathbf{L}q^\ast\dashv\mathbf{R}q_\ast$, and $\epsilon_{q^\prime}\colon \mathbf{L}(q^\prime)^\ast \mathbf{R}q^\prime_\ast \to 1$ is the counit of
	$\mathbf{L}(q^\prime)^\ast \dashv\mathbf{R}q^\prime_\ast$.
	The composite \eqref{eq:base_change_mate_formula} is exactly \cite[\href{https://stacks.math.columbia.edu/tag/0E5D}{Tag 0E5D}]{StacksProject};
	this coincides with the description given in \cite[Exercise (3.10.4) (b)$^\prime$]{Lipman:2009}.
\end{remark}
We emphasize again that it does not formally follow that $\beta_{\square}$ is an isomorphism from the fact that $\xi_{\square}^{-1}$ is an isomorphism.
For this, one needs geometric conditions on $p$ and $q$, see \Cref{lem:base_change_relative_formula}.
% Thus, the unit and counit in \eqref{eq:base_change_mate_formula} are ingredients in the standard formula for the mate of $\xi_\square^{-1}$; neither is itself the mate.
% We emphasize that the fact that $\xi_\square$ is an isomorphism does not formally imply that $\beta_\square$ is an isomorphism.
% The mate correspondence does not in general preserve invertibility.
% Thus, \eqref{eq:base_change_mate_formula} defines the base change  morphism, whereas proving that this morphism is an isomorphism requires additional geometric hypotheses.

We now record compatibility of the base change morphisms with
horizontal pasting of fibered squares,
\Cref{lem:base_change_mates_pasting,lem:tor-independence_and_beta_pasting}.
Vertical pasting is obtained by a dual argument, see \Cref{rem:hor_pasting_bc2,rem:hor_pasting_beta2}.

\begin{lemma}
	\label{lem:base_change_mates_pasting}
	The formation of the base change morphisms of
	\Cref{def:upper_star_lower_star_bc,def:lower_star_upper_cross_bc}
	is compatible with horizontal pasting of commutative squares.
	More precisely,
	consider a commutative diagram of quasi-compact quasi-separated
	algebraic spaces
	\begin{equation}
		\label{eq:horizontal_cart_rect}
		\begin{tikzcd}
			X^{\prime\prime}
			\arrow[r,"{g^\prime}"]
			\arrow[d,"{f^{\prime\prime}}"']
			&
			X^\prime
			\arrow[r,"g"]
			\arrow[d,"{f^\prime}"']
			&
			X
			\arrow[d,"f"]
			\\
			Y^{\prime\prime}
			\arrow[r,"{h^\prime}"']
			&
			Y^\prime
			\arrow[r,"h"']
			&
			Y
		\end{tikzcd}
	\end{equation}
	The following holds:
	\begin{enumerate}[label=(\arabic*), ref=\thelemma(\arabic*)]
		\item
		\label[lemma]{lem:upper_star_lower_star_pastes_horizontally}
		Let
		\begin{displaymath}
			\vartheta_{l}\colon
			\mathbf{L}(h^\prime)^\ast \mathbf{R}f^\prime_\ast
			\to \mathbf{R}f^{\prime\prime}_\ast \mathbf{L}(g^\prime)^\ast,
			\qquad
			\vartheta_{r}\colon
			\mathbf{L}h^\ast\mathbf{R}f_\ast
			\to
			\mathbf{R}f^\prime_\ast \mathbf{L}g^\ast,
		\end{displaymath}
		be the base change morphisms \eqref{eq:upper_star_lower_star} for the left and right squares of \eqref{eq:horizontal_cart_rect}.
		Then, under the canonical pseudofunctorial identifications,
		the $^\ast$-$_\ast$ base change morphisms for the outer rectangle of \eqref{eq:horizontal_cart_rect} is the composite
		\begin{equation}
			\label{eq:upper_star_lower_star_hor_paste}
			\mathbf{L}(h^\prime)^\ast\mathbf{L}h^\ast\mathbf{R}f_\ast
			\xrightarrow{\mathbf{L}(h^\prime)^\ast\vartheta_r}
			\mathbf{L}(h^\prime)^\ast\mathbf{R}f_\ast^\prime\mathbf{L}g^\ast
			\xrightarrow{\vartheta_l\mathbf{L}g^\ast}
			\mathbf{R}f_\ast^{\prime\prime}\mathbf{L}(g^\prime)^\ast
			\mathbf{L}g^\ast.
		\end{equation}
		\item
		\label[lemma]{lem:lower_star_upper_cross_pastes_horizontally}
		Let
		\begin{displaymath}
			\xi_{l}\colon
			\mathbf{R}g^\prime_\ast  (f^{\prime\prime})^\times
			\to (f^\prime)^\times \mathbf{R}h^\prime_\ast,
			\qquad
			\xi_{r}\colon \mathbf{R}g_\ast(f^\prime)^\times
			\to f^\times\mathbf{R}h_\ast,
		\end{displaymath}
		be the base change morphisms \eqref{eq:lower_star_upper_cross} for the left and right squares of \eqref{eq:horizontal_cart_rect}.
		Then, under the canonical pseudofunctorial identifications,
		the $_\ast$-$^\times$ base change morphism for the outer rectangle of \eqref{eq:horizontal_cart_rect} is the composite
		\begin{equation}
			\label{eq:lower_star_upper_cross_hor_paste}
			\mathbf{R}g_\ast\mathbf{R}g_\ast^\prime(f^{\prime\prime})^\times
			\xrightarrow{\mathbf{R}g_\ast\xi_l}
			\mathbf{R}g_\ast(f^\prime)^\times\mathbf{R}h_\ast^\prime
			\xrightarrow{\xi_r\mathbf{R}h^\prime_\ast}
			f^\times\mathbf{R}h_\ast\mathbf{R}h^\prime_\ast.
		\end{equation}
	\end{enumerate}
\end{lemma}

\begin{proof}
	The proof of \ref{lem:lower_star_upper_cross_pastes_horizontally}
	will leverage \ref{lem:upper_star_lower_star_pastes_horizontally}.
	Throughout the proof, we will write several mate diagrams.
	When doing so, in such diagrams we will follow the following conventions (see e.g.\ \eqref{eq:pseudofunct_of_lower_star_and_upper_star_upper_cross_vert_paste}):
	\begin{itemize}
		\item Every functor is assumed to be derived.
		\item We will simply write $Z$ instead of $D_{\operatorname{qc}}(Z)$, for $Z\in\{X,X^\prime,Y^\prime,Y^{\prime\prime}\}$.
	\end{itemize}

	Let us show \ref{lem:upper_star_lower_star_pastes_horizontally}.
	By \Cref{def:left_upper_star_lower_star_bc} of this base change morphism, we have the following mate diagrams:
	\begin{equation}
		\label{eq:pseudofunct_of_lower_star_and_upper_star_upper_cross_vert_paste}
		% https://q.uiver.app/#q=WzAsMTIsWzAsMiwiWF57XFxwcmltZVxccHJpbWV9Il0sWzEsMiwiWV57XFxwcmltZVxccHJpbWV9Il0sWzEsMSwiWV5cXHByaW1lIl0sWzEsMCwiWSJdLFswLDAsIlgiXSxbMCwxLCJYXlxccHJpbWUiXSxbMywyLCJYXntcXHByaW1lXFxwcmltZX0iXSxbNCwyLCJZXntcXHByaW1lXFxwcmltZX0iXSxbNCwxLCJZXlxccHJpbWUiXSxbNCwwLCJZIl0sWzMsMCwiWCJdLFszLDEsIlheXFxwcmltZSJdLFs1LDAsImdee1xccHJpbWVcXGFzdH0iLDJdLFs0LDUsImdeXFxhc3QiLDJdLFs0LDMsImZfXFxhc3QiXSxbNSwyLCJmXlxccHJpbWVfXFxhc3QiLDFdLFszLDIsImheXFxhc3QiXSxbMCwxLCJmXntcXHByaW1lXFxwcmltZX1fXFxhc3QiLDJdLFsyLDEsImhee1xccHJpbWVcXGFzdH0iXSxbMyw1LCJcXHZhcnRoZXRhX3IiLDIseyJzaG9ydGVuIjp7InNvdXJjZSI6MjAsInRhcmdldCI6MjB9LCJsZXZlbCI6Mn1dLFsyLDAsIlxcdmFydGhldGFfbCIsMix7InNob3J0ZW4iOnsic291cmNlIjoyMCwidGFyZ2V0IjoyMH0sImxldmVsIjoyfV0sWzYsMTEsImdeXFxwcmltZV9cXGFzdCJdLFsxMSwxMCwiZ19cXGFzdCJdLFsxMCw5LCJmX1xcYXN0Il0sWzExLDgsImZeXFxwcmltZV9cXGFzdCIsMV0sWzgsOSwiaF9cXGFzdCIsMl0sWzYsNywiZl57XFxwcmltZVxccHJpbWV9X1xcYXN0IiwyXSxbNyw4LCJoXlxccHJpbWVfXFxhc3QiLDJdLFsxMSw3LCJcXGNvbmciLDAseyJzaG9ydGVuIjp7InNvdXJjZSI6MjAsInRhcmdldCI6MjB9LCJsZXZlbCI6Mn1dLFsxMCw4LCJcXGNvbmciLDAseyJzaG9ydGVuIjp7InNvdXJjZSI6MjAsInRhcmdldCI6MjB9LCJsZXZlbCI6Mn1dLFsyLDExLCJcXGxlZnRyaWdodHNxdWlnYXJyb3ciLDEseyJzdHlsZSI6eyJib2R5Ijp7Im5hbWUiOiJub25lIn0sImhlYWQiOnsibmFtZSI6Im5vbmUifX19XV0=
		\begin{tikzcd}[row sep=large]
			X & Y && X & Y \\
			{X^\prime} & {Y^\prime} && {X^\prime} & {Y^\prime} \\
			{X^{\prime\prime}} & {Y^{\prime\prime}} && {X^{\prime\prime}} & {Y^{\prime\prime}}
			\arrow["{f_\ast}", from=1-1, to=1-2]
			\arrow["{g^\ast}"', from=1-1, to=2-1]
			\arrow["{\vartheta_r}"',  shorten <=2pt, shorten >=2pt, Rightarrow,  from=1-2, to=2-1]
			\arrow["{h^\ast}", from=1-2, to=2-2]
			\arrow["{f_\ast}", from=1-4, to=1-5]
			\arrow["\cong",  shorten <=2pt, shorten >=2pt, Rightarrow,  from=1-4, to=2-5]
			\arrow["{f^\prime_\ast}"{description}, from=2-1, to=2-2]
			\arrow["{g^{\prime\ast}}"', from=2-1, to=3-1]
			\arrow["\leftrightsquigarrow"{description}, draw=none, from=2-2, to=2-4]
			\arrow["{\vartheta_l}"',  shorten <=2pt, shorten >=2pt, Rightarrow,  from=2-2, to=3-1]
			\arrow["{h^{\prime\ast}}", from=2-2, to=3-2]
			\arrow["{g_\ast}", from=2-4, to=1-4]
			\arrow["{f^\prime_\ast}"{description}, from=2-4, to=2-5]
			\arrow["\cong",  shorten <=2pt, shorten >=2pt, Rightarrow,  from=2-4, to=3-5]
			\arrow["{h_\ast}"', from=2-5, to=1-5]
			\arrow["{f^{\prime\prime}_\ast}"', from=3-1, to=3-2]
			\arrow["{g^\prime_\ast}", from=3-4, to=2-4]
			\arrow["{f^{\prime\prime}_\ast}"', from=3-4, to=3-5]
			\arrow["{h^\prime_\ast}"', from=3-5, to=2-5]
		\end{tikzcd}
	\end{equation}
	Where in the right mate the upper and lower $2$-cells are the $\mathbf{R}(-)_\ast$-pseudofunctoriality isomorphisms.
	The resulting pasted left mate in
	\eqref{eq:pseudofunct_of_lower_star_and_upper_star_upper_cross_vert_paste} is \eqref{eq:upper_star_lower_star_hor_paste}.
	Hence, by \Cref{lem:mates_double_cat_iso}, it suffices to show that the pasted right mate in \eqref{eq:pseudofunct_of_lower_star_and_upper_star_upper_cross_vert_paste}
	identifies with the pseudofunctoriality isomorphism $\mathbf{R}f_\ast\mathbf{R}(g\circ g)_\ast\cong\mathbf{R}(h\circ h^\prime)_\ast\mathbf{R}f^{\prime\prime}_\ast$ of the outer rectangle of \eqref{eq:horizontal_cart_rect}.
	This follows from the fact that $\mathbf{R}(-)_\ast$ is a pseudofunctor (this is one of the numerous claims contained in \Cref{lem:spacey_Lipman_Hashimoto}), so in particular,
	the pseudofunctoriality isomorphism is associative, i.e.\ \cite[(3.6.5.1)]{Lipman:2009} holds for $\mathbf{R}(-)_\ast$.
	This finishes the proof of the claim for $^\ast$-$_\ast$ base change morphism.
	Note a similar proof might be obtained by invoking \Cref{def:right_upper_star_lower_star_bc} instead.

	We next show \ref{lem:lower_star_upper_cross_pastes_horizontally}.
	This is the claim made at the end in the proof of \cite[\href{https://stacks.math.columbia.edu/tag/0E5F}{Tag 0E5F}]{StacksProject},
	where it is said that the argument is `dual' to \ref{lem:upper_star_lower_star_pastes_horizontally}.
	We develop this idea.
	Consider the following diagram of natural transformations:
	\begin{equation}
		\label{eq:upper_star_lower_star_left_mates}
		\begin{tikzcd}[sep=large]
			X & X & X & X \\
			{X^\prime} & {X^\prime} & {Y^\prime} & {Y^\prime} \\
			{Y^{\prime\prime}} & {Y^{\prime\prime}} & {Y^{\prime\prime}} & {Y^{\prime\prime}.}
			\arrow[equals,  from=1-1, to=1-2]
			\arrow["{g^\ast}"', from=1-1, to=2-1]
			\arrow[equals,  from=1-2, to=1-3]
			\arrow[ equals,  from=1-2, to=2-1]
			\arrow["{g^\ast}"{description}, from=1-2, to=2-2]
			\arrow[equals,  from=1-3, to=1-4]
			\arrow[ equals,  from=1-3, to=2-2]
			\arrow["{f_\ast^\prime g^\ast}"{description}, from=1-3, to=2-3]
			\arrow["{\vartheta_r}"',  shorten <=2pt, shorten >=2pt, Rightarrow,  from=1-4, to=2-3]
			\arrow["{h^\ast f_\ast}", from=1-4, to=2-4]
			\arrow[equals,  from=2-1, to=2-2]
			\arrow["{f_\ast^{\prime\prime}g^\prime_\ast}"', from=2-1, to=3-1]
			\arrow["{f_\ast^\prime}"{description}, from=2-2, to=2-3]
			\arrow["{\vartheta_l}"',  shorten <=2pt, shorten >=2pt, Rightarrow,  from=2-2, to=3-1]
			\arrow["{h^{\prime \ast}f_\ast^\prime}"{description}, from=2-2, to=3-2]
			\arrow[equals,  from=2-3, to=2-4]
			\arrow[ equals,  from=2-3, to=3-2]
			\arrow["{h^{\prime\ast}}"{description}, from=2-3, to=3-3]
			\arrow[ equals,  from=2-4, to=3-3]
			\arrow["{h^{\prime\ast}}", from=2-4, to=3-4]
			\arrow[equals,  from=3-1, to=3-2]
			\arrow[equals,  from=3-2, to=3-3]
			\arrow[equals,  from=3-3, to=3-4]
		\end{tikzcd}
	\end{equation}
	The resulting pasted $2$-cell equals \eqref{eq:upper_star_lower_star_hor_paste}, which we know identifies with
	\begin{displaymath}
		\vartheta_{\textrm{rect}}\colon
		\mathbf{L}(h\circ h^\prime)^\ast
		\mathbf{R}f_\ast
		\xrightarrow{\cong}
		\mathbf{R}f^{\prime\prime}_\ast
		\mathbf{L}(g\circ g^\prime)^\ast,
	\end{displaymath}
	the $^\ast$-$_\ast$ base change morphism for the outer rectangle in \eqref{eq:horizontal_cart_rect}.
	By definition, the right conjugate of $\vartheta_{\textrm{rect}}$ is
	\begin{displaymath}
		\xi_{\textrm{rect}}\colon
		\mathbf{R}(g\circ g^\prime)_\ast (f^{\prime\prime})^\times
		\to f^\times \mathbf{R}(h\circ h^\prime)_\ast,
	\end{displaymath}
	the $_\ast$-$^\times$ base change morphism for the outer rectangle.
	Taking mates in \eqref{eq:upper_star_lower_star_left_mates}, it follows from \Cref{lem:mates_double_cat_iso} that $\xi_{\textrm{rect}}$ identifies with the resulting pasted $2$-cell of:
	\begin{equation}
		\label{eq:lower_star_upper_cross_right_mates}
		\begin{tikzcd}[sep=large]
			X & X & X & X \\
			{X^\prime} & {X^\prime} & {Y^\prime} & {Y^\prime} \\
			{Y^{\prime\prime}} & {Y^{\prime\prime}} & {Y^{\prime\prime}} & {Y^{\prime\prime},   }
			\arrow[equals,  from=1-1, to=1-2]
			\arrow[ equals,  from=1-1, to=2-2]
			\arrow[equals,  from=1-2, to=1-3]
			\arrow["{g_\ast \eta}",  shorten <=2pt, shorten >=2pt, Rightarrow,  from=1-2, to=2-3]
			\arrow[equals,  from=1-3, to=1-4]
			\arrow["{\xi_r}",  shorten <=2pt, shorten >=2pt, Rightarrow,  from=1-3, to=2-4]
			\arrow["{g_\ast}", from=2-1, to=1-1]
			\arrow[equals,  from=2-1, to=2-2]
			\arrow["{\xi_l}",  shorten <=2pt, shorten >=2pt, Rightarrow,  from=2-1, to=3-2]
			\arrow["{g_\ast}"{description}, from=2-2, to=1-2]
			\arrow["{f_\ast^\prime}"{description}, from=2-2, to=2-3]
			\arrow["{\epsilon h^\prime_\ast}",  shorten <=2pt, shorten >=2pt, Rightarrow,  from=2-2, to=3-3]
			\arrow["{g_\ast f^{\prime\times}}"{description}, from=2-3, to=1-3]
			\arrow[equals,  from=2-3, to=2-4]
			\arrow[ equals,  from=2-3, to=3-4]
			\arrow["{f^\times h_\ast}"', from=2-4, to=1-4]
			\arrow["{g_\ast^\prime f^{\prime\prime\times}}", from=3-1, to=2-1]
			\arrow[equals,  from=3-1, to=3-2]
			\arrow["{f^{\prime\times}h_\ast^\prime}"{description}, from=3-2, to=2-2]
			\arrow[equals,  from=3-2, to=3-3]
			\arrow["{h^\prime_\ast}"{description}, from=3-3, to=2-3]
			\arrow[equals,  from=3-3, to=3-4]
			\arrow["{h_\ast^\prime}"', from=3-4, to=2-4]
		\end{tikzcd}
	\end{equation}
	where $\eta$, $\epsilon$ are the unit and counit of the adjunction $\mathbf{R}f_\ast^\prime\dashv (f^\prime)^\times$.
	It suffices to show that the vertically composed central $2$-cell
	\begin{equation}
		\label{eq:vert_2_cell_will_be_id}
		% https://q.uiver.app/#q=WzAsNixbMCwwLCJYIl0sWzEsMCwiWCJdLFswLDIsIllee1xccHJpbWVcXHByaW1lfSJdLFsxLDIsIllee1xccHJpbWVcXHByaW1lfSJdLFswLDEsIlheXFxwcmltZSJdLFsxLDEsIlleXFxwcmltZSJdLFs0LDUsImZfXFxhc3ReXFxwcmltZSIsMV0sWzAsMSwiIiwyLHsibGV2ZWwiOjIsInN0eWxlIjp7ImhlYWQiOnsibmFtZSI6Im5vbmUifX19XSxbMiwzLCIiLDIseyJsZXZlbCI6Miwic3R5bGUiOnsiaGVhZCI6eyJuYW1lIjoibm9uZSJ9fX1dLFs0LDAsImdfXFxhc3QiLDFdLFs1LDEsImdfXFxhc3QgZl57XFxwcmltZVxcdGltZXN9IiwxXSxbMyw1LCJoXlxccHJpbWVfXFxhc3QiLDFdLFsyLDQsImZee1xccHJpbWVcXHRpbWVzfWhfKl5cXHByaW1lIiwxXSxbMCw1LCJnXypcXGV0YSIsMCx7InNob3J0ZW4iOnsic291cmNlIjoyMCwidGFyZ2V0IjoyMH0sImxldmVsIjoyfV0sWzQsMywiXFxlcHNpbG9uIGheXFxwcmltZV9cXGFzdCIsMCx7InNob3J0ZW4iOnsic291cmNlIjoyMCwidGFyZ2V0IjoyMH0sImxldmVsIjoyfV1d
		\begin{tikzcd}[sep=large]
			X & X \\
			{X^\prime} & {Y^\prime} \\
			{Y^{\prime\prime}} & {Y^{\prime\prime}}
			\arrow[equals,  from=1-1, to=1-2]
			\arrow["{g_\ast \eta}",  shorten <=2pt, shorten >=2pt, Rightarrow,  from=1-1, to=2-2]
			\arrow["{g_\ast}"{description}, from=2-1, to=1-1]
			\arrow["{f_\ast^\prime}"{description}, from=2-1, to=2-2]
			\arrow["{\epsilon h^\prime_\ast}",  shorten <=2pt, shorten >=2pt, Rightarrow,  from=2-1, to=3-2]
			\arrow["{g_\ast f^{\prime\times}}"{description}, from=2-2, to=1-2]
			\arrow["{f^{\prime\times}h_\ast^\prime}"{description}, from=3-1, to=2-1]
			\arrow[equals,  from=3-1, to=3-2]
			\arrow["{h^\prime_\ast}"{description}, from=3-2, to=2-2]
		\end{tikzcd}
	\end{equation}
	is the identity $2$-cell of $\mathbf{R}g_\ast (f^\prime)^\times\mathbf{R}h^\prime_\ast$ (indeed, the pasted $2$-cell \eqref{eq:lower_star_upper_cross_right_mates} will thus equal \eqref{eq:lower_star_upper_cross_hor_paste}, as desired).
	Writing every functor implicitly as derived and using `$\star$' for whiskering \cite[Remark 1.7.6]{Riehl:26}, the pasted $2$-cell \eqref{eq:vert_2_cell_will_be_id} equals:
	\begin{displaymath}
		\begin{aligned}
			&[(g_\ast f^{\prime\times})\star(\epsilon\star h_\ast^\prime)]
			\circ[(g_\ast\star\eta)\star(f^{\prime\times}h_\ast^\prime)]\\
			&=[g_\ast \star(f^{\prime\times}\star\epsilon\star h_\ast^\prime)]
			\circ[(g_\ast\star\eta\star f^{\prime\times})\star h_\ast^\prime]
			&\textrm{(\cite[Appendix, \S 1, (I),(II),(III)]{Godement:1958})}\\
			&=[g_\ast \star(f^{\prime\times}\star\epsilon)\star h_\ast^\prime]
			\circ[g_\ast\star(\eta\star f^{\prime\times})\star h_\ast^\prime]
			&\textrm{(\cite[Appendix, \S 1, (III)]{Godement:1958})}\\
			&=g_\ast\star
			[(f^{\prime\times}\star\epsilon)\circ(\eta\star f^{\prime\times})]
			\star h^\prime_\ast
			&\textrm{(\cite[Appendix, \S 1, (IV)]{Godement:1958})}\\
			&=g_\ast \star\operatorname{id}_{f^{\prime\times}}\star h^\prime_\ast
			&\textrm{(triangle identity)}\\
			&=\operatorname{id}_{g_\ast f^{\prime\times}h_\ast^\prime},
		\end{aligned}
	\end{displaymath}
	since whiskering an identity $2$-cell gives an identity $2$-cell.
	This finishes the proof of \ref{lem:lower_star_upper_cross_pastes_horizontally}.

\end{proof}

\begin{remark}
	\hfill
	\begin{enumerate}[label=(\arabic*), ref=\theremark(\arabic*)]
		\item \label[remark]{rem:hor_pasting_bc1}
		\ref{lem:upper_star_lower_star_pastes_horizontally} is
		\cite[(3.7.2) (iii)]{Lipman:2009}
		and \cite[\href{https://stacks.math.columbia.edu/tag/0E47}{Tag 0E47}]{StacksProject},
		whereas \ref{lem:lower_star_upper_cross_pastes_horizontally} is claimed without proof at the end of the proof of \cite[\href{https://stacks.math.columbia.edu/tag/0E5F}{Tag 0E5F}]{StacksProject}.
		Our proofs for each claim make an explicit use of theory of mates, materialized in the form of \Cref{lem:mates_double_cat_iso}.
		\item \label[remark]{rem:hor_pasting_bc2}
		The analogous statement to \Cref{lem:upper_star_lower_star_pastes_horizontally} for vertical pasting of fibered squares holds, see
		\cite[\href{https://stacks.math.columbia.edu/tag/0E46}{Tag 0E46}]{StacksProject}
		or \cite[(3.7.2) (ii)]{Lipman:2009}.\footref{fn:lipman_is_formal}
		% and for \ref{lem:lower_star_upper_cross_pastes_horizontally} this is \cite[\href{https://stacks.math.columbia.edu/tag/0E5E}{0E5E}]{StacksProject}.
		% The vertical version of \Cref{lem:base_change_mates_pasting} for the flipped square is obtained from the horizontal one for the unflipped one in the following way.
		% Namely, if we denote $\square$ to the square \eqref{diag:abstract_base_change}, write $\square^{\textrm{f}}$ for the flipped square
		% \begin{displaymath}
		%     % https://q.uiver.app/#q=WzAsNCxbMCwwLCJYXlxccHJpbWUiXSxbMSwwLCJaXlxccHJpbWUiXSxbMCwxLCJYIl0sWzEsMSwiWi4iXSxbMSwzLCJxIl0sWzIsMywicCIsMl0sWzAsMSwicF5cXHByaW1lIl0sWzAsMiwicV5cXHByaW1lIiwyXV0=
		%     \begin{tikzcd}[cramped]
		%         {X^\prime} & {Z^\prime} \\
		%         X & {Z.}
		%         \arrow["{p^\prime}", from=1-1, to=1-2]
		%         \arrow["{q^\prime}"', from=1-1, to=2-1]
		%         \arrow["q", from=1-2, to=2-2]
		%         \arrow["p"', from=2-1, to=2-2]
		%     \end{tikzcd}
		% \end{displaymath}
		% Then horizontal pasting in the left and right squares $\square_l$ and $\square_r$ in \eqref{eq:horizontal_cart_rect} amounts to vertical pasting of the flipped squares $\square_l^\textrm{f}$ and $\square_r^\textrm{f}$.
	\end{enumerate}
\end{remark}

\begin{lemma}[\'{e}tale base change]
	\label{lem:etale_presen_bc}
	The base change morphism \eqref{eq:upper_star_lower_star} is an isomorphism if $q$ is \'{e}tale.
\end{lemma}

\begin{proof}
	See \cite[Tags \href{https://stacks.math.columbia.edu/tag/0D6G}{0D6G},
	\href{https://stacks.math.columbia.edu/tag/04LX}{04LX},
	\href{https://stacks.math.columbia.edu/tag/03LV}{03LV},
	\href{https://stacks.math.columbia.edu/tag/04JA}{04JA},
	\href{https://stacks.math.columbia.edu/tag/03DH}{03DH},
	\href{https://stacks.math.columbia.edu/tag/0735}{0735},
	\href{https://stacks.math.columbia.edu/tag/015F}{015F},
	\href{https://stacks.math.columbia.edu/tag/06YY}{06YY}, and
	\href{https://stacks.math.columbia.edu/tag/08FH}{08FH}]{StacksProject}.
\end{proof}

\begin{lemma}
	\label{lem:characterize_tor-independence}
	The following are equivalent:
	\begin{enumerate}
		\item \label[lemma]{lem:characterize_tor-independence1}
		The morphism $\vartheta_\square$ of \eqref{eq:upper_star_lower_star} is an isomorphism.
		\item \label[lemma]{lem:characterize_tor-independence2}
		The morphism $\xi_\square$ of \eqref{eq:lower_star_upper_cross} is an isomorphism.
		\item \label[lemma]{lem:characterize_tor-independence3}
		The square \eqref{diag:abstract_base_change} is tor-independent.
	\end{enumerate}
\end{lemma}

\begin{proof}
	\eqref{lem:characterize_tor-independence1} $\iff$ \eqref{lem:characterize_tor-independence2}.
	This is \Cref{lem:when_mate_correspondence_preserves_isos1}.

	\eqref{lem:characterize_tor-independence1} $\iff$ \eqref{lem:characterize_tor-independence3}.
	See \cite[\href{https://stacks.math.columbia.edu/tag/08IR}{Tag 08IR}]{StacksProject} for the backwards direction.
	We give a different proof by showing we can reduce to the situation when \eqref{diag:abstract_base_change} is a fibered square of affine schemes,
	in which case the result is \cite[Lemma 3.10.3.1]{Lipman:2009}, or \cite[proof of Theorem 22.99]{Gortz/Wedhorn:2023}.
	First, we reduce to the case when $Z$ is an affine scheme.
	Choose an \'{e}tale presentation $t\colon V\to Z$ with $V$ an affine scheme.
	We form the fibered cube (i.e.\ every face is a fibered square):
	\begin{equation}
		\label{eq:cube_affine_reduction}
		% https://q.uiver.app/#q=WzAsOCxbMywxLCJYIl0sWzEsMSwiWF5cXHByaW1lIl0sWzEsMywiWl5cXHByaW1lIl0sWzMsMywiWiJdLFswLDIsIlZeXFxwcmltZSJdLFsyLDIsIlYiXSxbMiwwLCJVIl0sWzAsMCwiVV5cXHByaW1lIl0sWzAsMywicCJdLFsxLDAsInFeXFxwcmltZSIsMSx7ImxhYmVsX3Bvc2l0aW9uIjozMH1dLFsxLDIsInBeXFxwcmltZSIsMCx7ImxhYmVsX3Bvc2l0aW9uIjoyMH1dLFsyLDMsInEiLDJdLFs3LDYsImdeXFxwcmltZSJdLFs2LDAsInMiXSxbNywxLCJzXlxccHJpbWUiXSxbNyw0LCJmXlxccHJpbWUiLDJdLFs0LDIsInReXFxwcmltZSIsMl0sWzUsMywidCJdLFs2LDUsImYiLDIseyJsYWJlbF9wb3NpdGlvbiI6MjB9XSxbNCw1LCJnIiwwLHsibGFiZWxfcG9zaXRpb24iOjIwfV1d
		\begin{tikzcd}[cramped]
			{U^\prime} && U & \\
			& {X^\prime} && X \\
			{V^\prime} && V \\
			& {Z^\prime} && Z.
			\arrow["{g^\prime}", from=1-1, to=1-3]
			\arrow["{s^\prime}", from=1-1, to=2-2]
			\arrow["{f^\prime}"', from=1-1, to=3-1]
			\arrow["s", from=1-3, to=2-4]
			\arrow["f"'{pos=0.2}, from=1-3, to=3-3]
			\arrow["{q^\prime}"{description, pos=0.3}, from=2-2, to=2-4]
			\arrow["{p^\prime}"{pos=0.2}, from=2-2, to=4-2]
			\arrow["p", from=2-4, to=4-4]
			\arrow["g"{pos=0.2}, from=3-1, to=3-3]
			\arrow["{t^\prime}"', from=3-1, to=4-2]
			\arrow["t", from=3-3, to=4-4]
			\arrow["q"', from=4-2, to=4-4]
		\end{tikzcd}
	\end{equation}
	On the one hand, by \cite[\href{https://stacks.math.columbia.edu/tag/08IQ}{Tag 08IQ}]{StacksProject}, the front face in \eqref{eq:cube_affine_reduction} is tor-independent if, and only, if, so is the back face.
	Write $\vartheta_{\textrm{rect}}$ for the $^\ast$-$_\ast$ base change morphism for
	\begin{displaymath}
		% https://q.uiver.app/#q=WzAsNCxbMCwwLCJVXlxccHJpbWUiXSxbMCwxLCJWXlxccHJpbWUiXSxbMSwwLCJYIl0sWzEsMSwiWiJdLFswLDIsInNfMCJdLFsyLDMsInAiXSxbMCwxLCJmXlxccHJpbWUiLDJdLFsxLDMsInRfMCJdXQ==
		\begin{tikzcd}[cramped]
			{U^\prime} & X \\
			{V^\prime} & Z.
			\arrow["{s_0}", from=1-1, to=1-2]
			\arrow["{f^\prime}"', from=1-1, to=2-1]
			\arrow["p", from=1-2, to=2-2]
			\arrow["{t_0}", from=2-1, to=2-2]
		\end{tikzcd}
	\end{displaymath}
	We can decompose this square into two rectangle diagrams: one involving the back and right faces of \eqref{eq:cube_affine_reduction} and another one involving the left and front faces of \eqref{eq:cube_affine_reduction}.
	By \Cref{lem:upper_star_lower_star_pastes_horizontally}, these decompositions induce a commutative diagram
	\begin{equation}
		\label{diag:two_rect_factorizations}
		% https://q.uiver.app/#q=WzAsOCxbMCwwLCJnXlxcYXN0IHReXFxhc3QgcF9cXGFzdCJdLFsxLDAsImdeXFxhc3QgZl9cXGFzdCBzXlxcYXN0Il0sWzIsMCwiZl9cXGFzdF5cXHByaW1lIGdee1xccHJpbWVcXGFzdH1zXlxcYXN0Il0sWzAsMSwidF5cXGFzdF8wIHBfXFxhc3QiXSxbMiwxLCJmX1xcYXN0XlxccHJpbWUgc18wXlxcYXN0Il0sWzAsMiwidF57XFxwcmltZVxcYXN0fXFeXFxhc3QgcF9cXGFzdCJdLFsxLDIsInRee1xccHJpbWVcXGFzdH1wX1xcYXN0XlxccHJpbWUgcV57XFxwcmltZVxcYXN0fSJdLFsyLDIsImZfXFxhc3ReXFxwcmltZSBzXntcXHByaW1lXFxhc3R9cV57XFxwcmltZVxcYXN0fSJdLFswLDEsImdeXFxhc3RcXHZhcnRoZXRhX3IiXSxbMSwyLCJcXHZhcnRoZXRhX2Igc15cXGFzdCJdLFswLDMsIlxcY29uZyIsMl0sWzIsNCwiXFxjb25nIl0sWzUsNiwidF57XFxwcmltZVxcYXN0fVxcdmFydGhldGFfXFxzcXVhcmUiLDJdLFs2LDcsIlxcdmFydGhldGFfbHFee1xccHJpbWVcXGFzdH0iLDJdLFs1LDMsIlxcY29uZyJdLFs3LDQsIlxcY29uZyIsMl0sWzMsNCwiXFx2YXJ0aGV0YV97XFx0ZXh0e3JlY3R9fSIsMV1d
		\begin{tikzcd}[cramped]
			{g^\ast t^\ast p_\ast} & {g^\ast f_\ast s^\ast} & {f_\ast^\prime g^{\prime\ast}s^\ast} \\
			{t^\ast_0 p_\ast} && {f_\ast^\prime s_0^\ast} \\
			{t^{\prime\ast}q^\ast p_\ast} & {t^{\prime\ast}p_\ast^\prime q^{\prime\ast}} & {f_\ast^\prime s^{\prime\ast}q^{\prime\ast}}
			\arrow["{g^\ast\vartheta_r}", from=1-1, to=1-2]
			\arrow["\cong"', from=1-1, to=2-1]
			\arrow["{\vartheta_b s^\ast}", from=1-2, to=1-3]
			\arrow["\cong", from=1-3, to=2-3]
			\arrow["{\vartheta_{\textrm{rect}}}"{description}, from=2-1, to=2-3]
			\arrow["\cong", from=3-1, to=2-1]
			\arrow["{t^{\prime\ast}\vartheta_\square}"', from=3-1, to=3-2]
			\arrow["{\vartheta_lq^{\prime\ast}}"', from=3-2, to=3-3]
			\arrow["\cong"', from=3-3, to=2-3]
		\end{tikzcd}
	\end{equation}
	where:
	\begin{itemize}
		\item every functor appearing is assumed to be derived,
		\item the vertical morphisms are pseudofunctoriality of derived pullback, and
		\item the morphisms $\vartheta_r$, $\vartheta_b$, $\vartheta_\square$, and $\vartheta_l$ are the $^\ast$-$_\ast$ base change morphism respectively for the right, back, front and left faces of \eqref{eq:cube_affine_reduction}.
	\end{itemize}
	By \Cref{lem:etale_presen_bc}, $\vartheta_r$ and $\vartheta_l$ are isomorphisms.
	Therefore, we have:
	\begin{align*}
		\vartheta_\square\textrm{ is an iso}
		&\iff\mathbf{L}(t^\prime)^\ast\vartheta_\square\textrm{ is an iso, by \Cref{lem:faithfully_flat_is_conservative}}\\
		&\iff\vartheta_b\mathbf{L}s^\ast\textrm{ is an iso,
		by commutativity of \eqref{diag:two_rect_factorizations}}\\
		&\iff\vartheta_b\textrm{ is an iso, by \Cref{lem:prewhiskering_by_ess_dense_reflects_isos}
		and \cite[Proposition 1.4]{Lank:2026c}.}\footnote{}
	\end{align*}
	\footnotetext{A map from a separated algebraic space (like an affine scheme) is always separated.
	This follows from \cite[\href{https://stacks.math.columbia.edu/tag/03KR}{Tag 03KR}]{StacksProject}.}%
	This proves that it suffices to show the claim when $Z$ is an affine scheme.

	Next, we reduce to the case $Z^\prime$ is an affine scheme.
	Choose an \'{e}tale presentation $h^\prime\colon Z^{\prime\prime}\to Z^\prime$ with $Z^{\prime\prime}$ affine,
	and form the fiber product
	$X^{\prime\prime}:=Z^{\prime\prime}\times_{Z^\prime}X^\prime$.
	Consider the following rectangle:
	\begin{displaymath}
		% https://q.uiver.app/#q=WzAsNixbMCwwLCJYXntcXHByaW1lXFxwcmltZX0iXSxbMCwxLCJaXntcXHByaW1lXFxwcmltZX0iXSxbMSwwLCJYXlxccHJpbWUiXSxbMSwxLCJaXlxccHJpbWUiXSxbMiwwLCJYIl0sWzIsMSwiWiJdLFsyLDRdLFswLDJdLFsxLDMsImheXFxwcmltZSIsMl0sWzAsMV0sWzIsM10sWzMsNV0sWzQsNV0sWzEwLDEyLCJcXHNxdWFyZSIsMSx7InNob3J0ZW4iOnsic291cmNlIjoyMCwidGFyZ2V0IjoyMH0sInN0eWxlIjp7ImJvZHkiOnsibmFtZSI6Im5vbmUifSwiaGVhZCI6eyJuYW1lIjoibm9uZSJ9fX1dLFs5LDEwLCJsIiwxLHsic2hvcnRlbiI6eyJzb3VyY2UiOjIwLCJ0YXJnZXQiOjIwfSwic3R5bGUiOnsiYm9keSI6eyJuYW1lIjoibm9uZSJ9LCJoZWFkIjp7Im5hbWUiOiJub25lIn19fV1d
		\begin{tikzcd}[cramped]
			{X^{\prime\prime}} & {X^\prime} & X \\
			{Z^{\prime\prime}} & {Z^\prime} & Z
			\arrow[from=1-1, to=1-2]
			\arrow[""{name=0, anchor=center, inner sep=0}, from=1-1, to=2-1]
			\arrow[from=1-2, to=1-3]
			\arrow[""{name=1, anchor=center, inner sep=0}, from=1-2, to=2-2]
			\arrow[""{name=2, anchor=center, inner sep=0}, from=1-3, to=2-3]
			\arrow["{h^\prime}"', from=2-1, to=2-2]
			\arrow[from=2-2, to=2-3]
			\arrow["l"{description}, draw=none, from=0, to=1]
			\arrow["\square"{description}, draw=none, from=1, to=2]
		\end{tikzcd}
	\end{displaymath}
	since $\vartheta_l$ is an isomorphism by \Cref{lem:etale_presen_bc}, from \Cref{lem:upper_star_lower_star_pastes_horizontally} we deduce that the $^\ast$-$_\ast$ base change morphism for the outer rectangle is an isomorphism if, and only if, so is $\mathbf{L}(h^\prime)^\ast\vartheta_\square$, which in turn happens if, and only if, $\vartheta_\square$ is an isomorphism (\Cref{lem:faithfully_flat_is_conservative}).
	On the other hand, by \cite[\href{https://stacks.math.columbia.edu/tag/08IQ}{Tag 08IQ}]{StacksProject} the outer rectangle is tor-independent if, and only if, the right square is tor-independent.
	This shows we may assume both $Z^\prime$ and $Z$ are affine schemes.
	The reduction to $X$ an affine scheme is the flipped argument.
	Namely, one takes an \'{e}tale presentation $U\to X$ with $U$ an affine scheme.
	Then one performs the flipped argument as before by invoking the vertical analogue to
	\Cref{lem:upper_star_lower_star_pastes_horizontally}
	(see \Cref{rem:hor_pasting_bc2}).
	This shows we may assume that \eqref{diag:abstract_base_change} is a diagram of affine schemes and thus the claim is \cite[Lemma 3.10.3.1]{Lipman:2009}, or \cite[proof of Theorem 22.99]{Gortz/Wedhorn:2023}.
\end{proof}

The following result is used extensively in the proof of \Cref{lem:base_change_relative_formula}.

\begin{lemma}
	\label[lemma]{lem:tor-independence_and_beta_pasting}
	The formation of the base change morphism of
	\Cref{def:upper_star_upper_cross_bc}
	is compatible with horizontal pasting of fibered squares.
	More precisely,
	consider a commutative diagram of quasi-compact quasi-separated
	algebraic spaces
	\begin{equation}
		\label{eq:horizontal_cart_rect_beta}
		\begin{tikzcd}
			X^{\prime\prime}
			\arrow[r,"{g^\prime}"]
			\arrow[d,"{f^{\prime\prime}}"']
			&
			X^\prime
			\arrow[r,"g"]
			\arrow[d,"{f^\prime}"']
			&
			X
			\arrow[d,"f"]
			\\
			Y^{\prime\prime}
			\arrow[r,"{h^\prime}"']
			&
			Y^\prime
			\arrow[r,"h"']
			&
			Y
		\end{tikzcd}
	\end{equation}
	in which the two squares are fibered and tor-independent.
	Then:
	\begin{enumerate}[label=(\arabic*), ref=\thelemma(\arabic*)]
		\item \label[lemma]{lem:tor-independence_pastes}
		The outer rectangle is fibered and tor-independent.
		\item \label[lemma]{lem:upper_star_upper_cross_pastes_horizontally}
		Let
		\begin{displaymath}
			\beta_{l}\colon
			\mathbf{L}(g^\prime)^\ast
			(f^\prime)^\times
			\to
			(f^{\prime\prime})^\times
			\mathbf{L}(h^\prime)^\ast,
			\qquad
			\beta_{r}\colon
			\mathbf{L}g^\ast f^\times
			\to
			(f^\prime)^\times\mathbf{L}h^\ast,
		\end{displaymath}
		be the base change morphisms \eqref{eq:base_change_mate_formula} (same as \cite[\href{https://stacks.math.columbia.edu/tag/0E5D}{Tag 0E5D}]{StacksProject}) for the left and right squares of \eqref{eq:horizontal_cart_rect_beta}.
		Then, under the canonical pseudofunctorial identifications,
		the $^\ast$-$^\times$ base change morphism for the outer rectangle of \eqref{eq:horizontal_cart_rect_beta} is the composite
		\begin{equation}
			\label{eq:upper_star_upper_cross_hor_paste}
			\mathbf{L}(g^\prime)^\ast
			\mathbf{L}g^\ast f^\times
			\xrightarrow{
			\mathbf{L}(g^\prime)^\ast
			\beta_{r}}
			\mathbf{L}(g^\prime)^\ast
			(f^\prime)^\times
			\mathbf{L}h^\ast
			\xrightarrow{
			\beta_{l}
			\mathbf{L}h^\ast}
			(f^{\prime\prime})^\times
			\mathbf{L}(h^\prime)^\ast
			\mathbf{L}h^\ast.
		\end{equation}
	\end{enumerate}
\end{lemma}

\begin{proof}
	Since the two squares are tor-independent, $\zeta_{r}$, $\zeta_{l}$, and $\zeta_{\textrm{rect}}$ are isomorphisms, for $\zeta\in\{\vartheta,\xi\}$.
	Hence the rectangle is tor-independent, by \Cref{lem:characterize_tor-independence}.

	We will give two proofs for
	\ref{lem:upper_star_upper_cross_pastes_horizontally},
	one which uses \ref{lem:upper_star_lower_star_pastes_horizontally}
	and another using \ref{lem:lower_star_upper_cross_pastes_horizontally}.
	We use the same diagramatic conventions explained at the beginning of the proof of \Cref{lem:base_change_mates_pasting}.

	\underline{First proof}:
	By \Cref{def:upper_star_upper_cross_bc1}, we have the following mate diagrams:
	\begin{displaymath}
		% https://q.uiver.app/#q=WzAsMTIsWzIsMCwiWF57XFxwcmltZVxccHJpbWV9Il0sWzIsMSwiWV57XFxwcmltZVxccHJpbWV9Il0sWzEsMSwiWV5cXHByaW1lIl0sWzAsMSwiWSJdLFswLDAsIlgiXSxbMSwwLCJYXlxccHJpbWUiXSxbNiwwLCJYXntcXHByaW1lXFxwcmltZX0iXSxbNiwxLCJZXntcXHByaW1lXFxwcmltZX0iXSxbNSwxLCJZXlxccHJpbWUiXSxbNCwxLCJZIl0sWzQsMCwiWCJdLFs1LDAsIlheXFxwcmltZSJdLFs1LDAsImdee1xccHJpbWVcXGFzdH0iXSxbNCw1LCJnXlxcYXN0Il0sWzQsMywiZl9cXGFzdCIsMl0sWzUsMiwiZl5cXHByaW1lX1xcYXN0IiwxXSxbMywyLCJoXlxcYXN0IiwyXSxbMCwxLCJmXntcXHByaW1lXFxwcmltZX1fXFxhc3QiXSxbMiwxLCJoXntcXHByaW1lXFxhc3R9IiwyXSxbNSwzLCJcXHZhcnRoZXRhX3Jeey0xfSIsMix7InNob3J0ZW4iOnsic291cmNlIjoyMCwidGFyZ2V0IjoyMH0sImxldmVsIjoyfV0sWzAsMiwiXFx2YXJ0aGV0YV9sXnstMX0iLDIseyJzaG9ydGVuIjp7InNvdXJjZSI6MjAsInRhcmdldCI6MjB9LCJsZXZlbCI6Mn1dLFsxMSw2LCJnXntcXHByaW1lXFxhc3R9Il0sWzEwLDExLCJnXlxcYXN0Il0sWzksMTAsImZeXFx0aW1lcyJdLFs4LDExLCJmXntcXHByaW1lXFx0aW1lc30iLDFdLFs5LDgsImheXFxhc3QiLDJdLFs3LDYsImZee1xccHJpbWVcXHByaW1lXFx0aW1lc30iLDJdLFs4LDcsImhee1xccHJpbWVcXGFzdH0iLDJdLFsxMCw4LCJcXGJldGFfciIsMCx7InNob3J0ZW4iOnsic291cmNlIjoyMCwidGFyZ2V0IjoyMH0sImxldmVsIjoyfV0sWzExLDcsIlxcYmV0YV9sIiwwLHsic2hvcnRlbiI6eyJzb3VyY2UiOjIwLCJ0YXJnZXQiOjIwfSwibGV2ZWwiOjJ9XSxbMTcsMjMsIlxcbGVmdHJpZ2h0c3F1aWdhcnJvdyIsMSx7InNob3J0ZW4iOnsic291cmNlIjoyMCwidGFyZ2V0IjoyMH0sInN0eWxlIjp7ImJvZHkiOnsibmFtZSI6Im5vbmUifSwiaGVhZCI6eyJuYW1lIjoibm9uZSJ9fX1dXQ==
		\begin{tikzcd}[sep=large]
			X & {X^\prime} & {X^{\prime\prime}} && X & {X^\prime} & {X^{\prime\prime}} \\
			Y & {Y^\prime} & {Y^{\prime\prime}} && Y & {Y^\prime} & {Y^{\prime\prime}}
			\arrow["{g^\ast}", from=1-1, to=1-2]
			\arrow["{f_\ast}"', from=1-1, to=2-1]
			\arrow["{g^{\prime\ast}}", from=1-2, to=1-3]
			\arrow["{\vartheta_r^{-1}}"', shorten <=2pt, shorten >=2pt, Rightarrow, from=1-2, to=2-1]
			\arrow["{f^\prime_\ast}"{description}, from=1-2, to=2-2]
			\arrow["{\vartheta_l^{-1}}"', shorten <=2pt, shorten >=2pt, Rightarrow, from=1-3, to=2-2]
			\arrow[""{name=0, anchor=center, inner sep=0}, "{f^{\prime\prime}_\ast}", from=1-3, to=2-3]
			\arrow["{g^\ast}", from=1-5, to=1-6]
			\arrow["{\beta_r}", shorten <=2pt, shorten >=2pt, Rightarrow, from=1-5, to=2-6]
			\arrow["{g^{\prime\ast}}", from=1-6, to=1-7]
			\arrow["{\beta_l}", shorten <=2pt, shorten >=2pt, Rightarrow, from=1-6, to=2-7]
			\arrow["{h^\ast}"', from=2-1, to=2-2]
			\arrow["{h^{\prime\ast}}"', from=2-2, to=2-3]
			\arrow[""{name=1, anchor=center, inner sep=0}, "{f^\times}", from=2-5, to=1-5]
			\arrow["{h^\ast}"', from=2-5, to=2-6]
			\arrow["{f^{\prime\times}}"{description}, from=2-6, to=1-6]
			\arrow["{h^{\prime\ast}}"', from=2-6, to=2-7]
			\arrow["{f^{\prime\prime\times}}"', from=2-7, to=1-7]
			\arrow["\leftrightsquigarrow"{description}, draw=none, from=0, to=1]
		\end{tikzcd}
	\end{displaymath}
	The pasted left mate is
	\begin{displaymath}
		(\mathbf{L}(h^\prime)^\ast\vartheta_r^{-1})
		\circ(\vartheta_l^{-1}\mathbf{L}g^\ast).
	\end{displaymath}
	This is the inverse of \eqref{eq:upper_star_lower_star_hor_paste}.
	Thus, by \ref{lem:upper_star_lower_star_pastes_horizontally}, it identifies with $\vartheta_{\textrm{rect}}^{-1}$.
	On the other hand, the pasted right mate equals \eqref{eq:upper_star_upper_cross_hor_paste}.
	The claim follows by \Cref{lem:mates_double_cat_iso}.

	\underline{Second proof}:
	Taking inverses in \eqref{eq:lower_star_upper_cross_hor_paste}, from \ref{lem:lower_star_upper_cross_pastes_horizontally} we get that $\xi_{\textrm{rect}}^{-1}$ identifies with
	\begin{equation}
		\label{eq:inverse_xi_base_change_pasting}
		(\mathbf{R}g_\ast\xi_{l}^{-1}) \circ(\xi_{r}^{-1} \mathbf{R}h^\prime_\ast).
	\end{equation}
	By \Cref{def:upper_star_upper_cross_bc2}, we have the following mate diagrams:
	\begin{equation}
		\label{eq:upper_cross_lower_star_and_upper_star_upper_cross_vert_paste}
		% https://q.uiver.app/#q=WzAsMTIsWzMsMCwiWSJdLFs0LDAsIlgiXSxbMywxLCJZXntcXHByaW1lfSJdLFs0LDEsIlhee1xccHJpbWV9Il0sWzQsMiwiWF57XFxwcmltZVxccHJpbWV9Il0sWzMsMiwiWV57XFxwcmltZVxccHJpbWV9Il0sWzAsMCwiWSJdLFsxLDAsIlgiXSxbMCwxLCJZXntcXHByaW1lfSJdLFsxLDEsIlhee1xccHJpbWV9Il0sWzEsMiwiWF57XFxwcmltZVxccHJpbWV9Il0sWzAsMiwiWV57XFxwcmltZVxccHJpbWV9Il0sWzAsMSwiZl5cXHRpbWVzIl0sWzIsMCwiaF9cXGFzdCJdLFszLDEsImdfXFxhc3QiLDJdLFsyLDMsImZee1xccHJpbWVcXHRpbWVzfSIsMV0sWzUsMiwiaF9cXGFzdF5cXHByaW1lIiwyXSxbNCwzLCJnX1xcYXN0XlxccHJpbWUiLDJdLFs1LDQsImZee1xccHJpbWVcXHByaW1lXFx0aW1lc30iLDJdLFswLDMsIlxceGlfcl57LTF9IiwwLHsic2hvcnRlbiI6eyJzb3VyY2UiOjIwLCJ0YXJnZXQiOjIwfSwibGV2ZWwiOjJ9XSxbMiw0LCJcXHhpX2xeey0xfSIsMCx7InNob3J0ZW4iOnsic291cmNlIjoyMCwidGFyZ2V0IjoyMH0sImxldmVsIjoyfV0sWzYsNywiZl5cXHRpbWVzIl0sWzYsOCwiaF5cXGFzdCIsMl0sWzcsOSwiZ15cXGFzdCJdLFs4LDksImZee1xccHJpbWVcXHRpbWVzfSIsMV0sWzgsMTEsImhee1xccHJpbWVcXGFzdH0iLDJdLFs5LDEwLCJnXntcXHByaW1lXFxhc3R9Il0sWzExLDEwLCJmXntcXHByaW1lXFxwcmltZVxcdGltZXN9IiwyXSxbNyw4LCJcXGJldGFfciIsMix7InNob3J0ZW4iOnsic291cmNlIjoyMCwidGFyZ2V0IjoyMH0sImxldmVsIjoyfV0sWzksMTEsIlxcYmV0YV9sIiwyLHsic2hvcnRlbiI6eyJzb3VyY2UiOjIwLCJ0YXJnZXQiOjIwfSwibGV2ZWwiOjJ9XSxbOSwyLCJcXGxlZnRyaWdodHNxdWlnYXJyb3ciLDEseyJzdHlsZSI6eyJib2R5Ijp7Im5hbWUiOiJub25lIn0sImhlYWQiOnsibmFtZSI6Im5vbmUifX19XV0=
		\begin{tikzcd}[sep=large]
			Y & X && Y & X \\
			{Y^\prime} & {X^\prime} && {Y^\prime} & {X^\prime} \\
			{Y^{\prime\prime}} & {X^{\prime\prime}} && {Y^{\prime\prime}} & {X^{\prime\prime}}
			\arrow["{f^\times}", from=1-1, to=1-2]
			\arrow["{h^\ast}"', from=1-1, to=2-1]
			\arrow["{\beta_r}"',  shorten <=2pt, shorten >=2pt, Rightarrow,  from=1-2, to=2-1]
			\arrow["{g^\ast}", from=1-2, to=2-2]
			\arrow["{f^\times}", from=1-4, to=1-5]
			\arrow["{\xi_r^{-1}}",  shorten <=2pt, shorten >=2pt, Rightarrow,  from=1-4, to=2-5]
			\arrow["{f^{\prime\times}}"{description}, from=2-1, to=2-2]
			\arrow["{h^{\prime\ast}}"', from=2-1, to=3-1]
			\arrow["\leftrightsquigarrow"{description}, draw=none, from=2-2, to=2-4]
			\arrow["{\beta_l}"',  shorten <=2pt, shorten >=2pt, Rightarrow,  from=2-2, to=3-1]
			\arrow["{g^{\prime\ast}}", from=2-2, to=3-2]
			\arrow["{h_\ast}", from=2-4, to=1-4]
			\arrow["{f^{\prime\times}}"{description}, from=2-4, to=2-5]
			\arrow["{\xi_l^{-1}}",  shorten <=2pt, shorten >=2pt, Rightarrow,  from=2-4, to=3-5]
			\arrow["{g_\ast}"', from=2-5, to=1-5]
			\arrow["{f^{\prime\prime\times}}"', from=3-1, to=3-2]
			\arrow["{h_\ast^\prime}"', from=3-4, to=2-4]
			\arrow["{f^{\prime\prime\times}}"', from=3-4, to=3-5]
			\arrow["{g_\ast^\prime}"', from=3-5, to=2-5]
		\end{tikzcd}
	\end{equation}
	The composite right mate is \eqref{eq:inverse_xi_base_change_pasting}, whereas the composite left mate is \eqref{eq:upper_star_upper_cross_hor_paste}.
	Using the notation of \Cref{lem:mates_bijection},
	\begin{displaymath}
		\begin{aligned}
			\beta_{\textrm{rect}}
			&=\overline{\xi_{\textrm{rect}}^{-1}}\\
			&=\overline{\xi_l^{-1}\diamond\xi_r^{-1}}\\
			&=\overline{\xi_l^{-1}}\diamond\overline{\xi_r^{-1}}
			&\qquad\textrm{(\Cref{lem:mates_double_cat_iso})}\\
			&=\beta_l\diamond\beta_r,
		\end{aligned}
	\end{displaymath}
	where `$\diamond$' denotes vertical paste of $2$-cells in \eqref{eq:upper_cross_lower_star_and_upper_star_upper_cross_vert_paste}.
\end{proof}

\begin{remark}
	\hfill
	\begin{enumerate}[label=(\arabic*), ref=\theremark(\arabic*)]
		\item \label[remark]{rem:hor_pasting_beta1}
		\Cref{lem:upper_star_upper_cross_pastes_horizontally} is the compatibility proven in \cite[\href{https://stacks.math.columbia.edu/tag/0E5F}{Tag 0E5F}]{StacksProject}, whose proof strategy is our second proof.
		The proof above makes explicit the mate operations
		entering the construction,
		building each subsequent pair of mates upon the previous one.
		\item \label[remark]{rem:hor_pasting_beta2}
		The analogous statement to \Cref{lem:upper_star_upper_cross_pastes_horizontally} for vertical pasting of fibered squares holds
		\cite[\href{https://stacks.math.columbia.edu/tag/0E5E}{Tag 0E5E}]{StacksProject}.
		% The vertical version of \Cref{lem:base_change_mates_pasting} for the flipped square is obtained from the horizontal one for the unflipped one in the following way.
		% Namely, if we denote $\square$ to the square \eqref{diag:abstract_base_change}, write $\square^{\textrm{f}}$ for the flipped square
		% \begin{displaymath}
		%     % https://q.uiver.app/#q=WzAsNCxbMCwwLCJYXlxccHJpbWUiXSxbMSwwLCJaXlxccHJpbWUiXSxbMCwxLCJYIl0sWzEsMSwiWi4iXSxbMSwzLCJxIl0sWzIsMywicCIsMl0sWzAsMSwicF5cXHByaW1lIl0sWzAsMiwicV5cXHByaW1lIiwyXV0=
		%     \begin{tikzcd}[cramped]
		%         {X^\prime} & {Z^\prime} \\
		%         X & {Z.}
		%         \arrow["{p^\prime}", from=1-1, to=1-2]
		%         \arrow["{q^\prime}"', from=1-1, to=2-1]
		%         \arrow["q", from=1-2, to=2-2]
		%         \arrow["p"', from=2-1, to=2-2]
		%     \end{tikzcd}
		% \end{displaymath}
		% Then horizontal pasting in the left and right squares $\square_l$ and $\square_r$ in \eqref{eq:horizontal_cart_rect} amounts to vertical pasting of the flipped squares $\square_l^\textrm{f}$ and $\square_r^\textrm{f}$.
	\end{enumerate}
\end{remark}

\begin{remark}
	\label{rem:abstract_nonsense_for_base_change}
	\Cref{def:upper_star_lower_star_bc,def:lower_star_upper_cross_bc,lem:pseudofunctoriality_isos_are_conjugate,lem:base_change_mates_pasting,lem:tor-independence_and_beta_pasting,def:upper_star_upper_cross_bc} are formal,
	and therefore replicable in higher generality.
	Let $\mathbf{S}$ be a category.
	\begin{enumerate}
		\item \label{rem:abstract_nonsense_for_base_change1}
		A pair $(^{\boldsymbol{\ast}},\, _{\boldsymbol{\ast}})$ with $^{\boldsymbol{\ast}}$ and $_{\boldsymbol{\ast}}$ respectively a contravariant and covariant pseudofunctor on $\mathbf{S}$
		is said to be adjoint if \cite[(3.6.7) (d)]{Lipman:2009} holds.
		Note that \cite[(3.6.7) (d) (iii)]{Lipman:2009} just means that the isomorphisms $f^\ast g^\ast \to (gf)^\ast$ and $(gf)_\ast\to g_\ast f_\ast$ are conjugate to each other,
		for \cite[(3.6.2) and (3.6.2)\textsuperscript{\textsf{op}}]{Lipman:2009}) corresponds to \ref{lem:morph_adj3} and \ref{lem:morph_adj2}, respectively.
		\item \label{rem:abstract_nonsense_for_base_change2}
		A pair $(_{\boldsymbol{\ast}},\, ^{\boldsymbol{\times}})$ with $_{\boldsymbol{\ast}}$ and $^{\boldsymbol{\times}}$ respectively a covariant and contravariant pseudofunctor on $\mathbf{S}$ is said to be adjoint if:
		\begin{enumerate}[label=(\roman*)]
			\item $\mathbf{X}_{\boldsymbol{\ast}}=\mathbf{X}^{\boldsymbol{\times}}$ for all objects $X$ in $\mathbf{S}$.
			\item For every $f\colon X\to Y$ in $\mathbf{S}$ there are bifunctorial isomorphisms
			\begin{displaymath}
				\operatorname{Hom}_{\mathbf{X}_{\boldsymbol{\ast}}}(f_\ast C,D)
				\cong
				\operatorname{Hom}_{\mathbf{Y}^{\boldsymbol{\times}}}(C, f^\times D)
				\qquad
				(C\in \mathbf{X}_{\boldsymbol{\ast}},
				D\in \mathbf{Y}^{\boldsymbol{\times}}),
			\end{displaymath}
			i.e.\ the functor
			$f_\ast\colon \mathbf{X}_{\boldsymbol{\ast}}\to \mathbf{Y}_{\boldsymbol{\ast}}$
			is left adjoint to
			$f^\times\colon \mathbf{Y}^{\boldsymbol{\times}}\to \mathbf{X}^{\boldsymbol{\times}}$.
			\item The isomorphisms $(gf)_\ast\to g_\ast f_\ast$ and $f^\times g^\times\to (gf)^\times$ are conjugate to each other.
		\end{enumerate}
		\item \label{rem:abstract_nonsense_for_base_change3}
		Define an \textbf{adjoint triple of pseudofunctors} on a category $\mathbf{S}$ is a triple $(^{\boldsymbol{\ast}},\, _{\boldsymbol{\ast}},\, ^{\boldsymbol{\times}})$ such that:
		\begin{enumerate}
			\item $^{\boldsymbol{\ast}}$ and $^{\boldsymbol{\times}}$ are contravariant pseudofunctors on $\mathbf{S}$.
			\item $_{\boldsymbol{\ast}}$ is a covariant pseudofunctor on $\mathbf{S}$.
			\item $(^{\boldsymbol{\ast}},\, _{\boldsymbol{\ast}})$ and $(_{\boldsymbol{\ast}},\, ^{\boldsymbol{\times}})$ are adjoint pairs of pseudofunctors on $\mathbf{S}$ in the senses \eqref{rem:abstract_nonsense_for_base_change1} and \eqref{rem:abstract_nonsense_for_base_change2}, respectively.
		\end{enumerate}
		There is the following existence and uniqueness result:
		Suppose we have a contravariant $^{\boldsymbol{\ast}}$ or a covariant $_{\boldsymbol{\ast}}$ pseudofunctor on $\mathbf{S}$.
		Suppose furthermore that for each morphism $f\colon X\to Y$ on $\mathbf{S}$, there is an adjoint triple $f^\ast\dashv f_\ast\dashv f^\times$, that we fix for each $f$.
		Then the data of these adjoint triples extend uniquely to an adjoint triple of pseudofunctors $(^{\boldsymbol{\ast}},\, _{\boldsymbol{\ast}},\, ^{\boldsymbol{\times}})$.
		Indeed, this is a straight application of the argument in \cite[(3.6.8.1) and its dual]{Lipman:2009}.
		Two instances where this result applies to are $D_{\operatorname{qc}}$ for quasi-compact separated schemes \cite[(3.10), first paragraph]{Lipman:2009} or for quasi-compact quasi-separated algebraic spaces, \Cref{lem:spacey_Lipman_Hashimoto}.

		\item \label{rem:abstract_nonsense_for_base_change4}
		(\textsc{Duality principle for adjoint triples of pseudofunctors}):\footnote{Although similar in spirit to \cite[(3.6.9)]{Lipman:2009}, this is not the same as `applying twice' the duality principle for adjoint pairs in loc.\ cit.,
		one for $(^{\boldsymbol{\ast}},\, _{\boldsymbol{\ast}})$ and another time for $(_{\boldsymbol{\ast}},\, ^{\boldsymbol{\times}})$.
		The latter procedure spits out a triple of signature covariant-contravariant-covariant,
		but we defined \eqref{rem:abstract_nonsense_for_base_change3} to be of signature contravariant-covariant-contravariant.
		Also, \cite[(3.6.9)]{Lipman:2009} for a pair spits out pseudofunctors on $\mathbf{S}^\mathrm{op}$, whereas our duality principle for triples spits out pseudofunctors on $\mathbf{S}$.}
		Given an adjoint triple of pseudofunctors $(^{\boldsymbol{\ast}},\, _{\boldsymbol{\ast}},\, ^{\boldsymbol{\times}})$ on $\mathbf{S}$,
		we define another adjoint triple of pseudofunctors
		$(^{\boldsymbol{\star}},\, _{\boldsymbol{\star}},\, ^{\boldsymbol{+}})$ on $\mathbf{S}$
		in the following way:
		\begin{displaymath}
			\mathbf{X}^{\boldsymbol{\star}}:=(\mathbf{X}^{\boldsymbol{\times}})^\mathrm{op},
			\qquad
			\mathbf{X}_{\boldsymbol{\star}}:=(\mathbf{X}_{\boldsymbol{\ast}})^\mathrm{op},
			\qquad
			\mathbf{X}^{\boldsymbol{+}}:=(\mathbf{X}^{\boldsymbol{\ast}})^\mathrm{op},
		\end{displaymath}
		for objects $X\in\mathbf{S}$, and
		\begin{displaymath}
			\begin{gathered}
				f^\star:=(f^\times)^\mathrm{op}\colon
				(\mathbf{Y}^{\boldsymbol{\times}})^\mathrm{op}\to
				(\mathbf{X}^{\boldsymbol{\times}})^\mathrm{op},\\
				f_\star:=(f_\ast)^\mathrm{op}\colon
				(\mathbf{X}_{\boldsymbol{\ast}})^\mathrm{op}\to
				(\mathbf{Y}_{\boldsymbol{\ast}})^\mathrm{op},\\
				f^+:=(f^\ast)^\mathrm{op}\colon
				(\mathbf{Y}^{\boldsymbol{\ast}})^\mathrm{op}\to
				(\mathbf{X}^{\boldsymbol{\ast}})^\mathrm{op},
			\end{gathered}
		\end{displaymath}
		for each map $f\colon X\to Y$ in $\mathbf{S}$.
		It follows that:

		\emph{Each diagram built up from the basic data defining adjoint
		triples can be interpreted in the dual context, giving rise to a `dual' diagram,
		obtained by swapping $^{\boldsymbol{\ast}}$ with $^{\boldsymbol{\times}}$ (while keeping $_{\boldsymbol{\ast}}$ fixed) and reversing arrows accordingly, etc., etc.}

		\item The notion \eqref{rem:abstract_nonsense_for_base_change1} has monoidal and triangulated refinements \cite[(3.6.7) (d), (iv, v)]{Lipman:2009}.
		Similarly, so do \eqref{rem:abstract_nonsense_for_base_change2}, \eqref{rem:abstract_nonsense_for_base_change3} and \eqref{rem:abstract_nonsense_for_base_change4}.

		\item \label{rem:abstract_nonsense_for_base_change6}
		Let $(^{\boldsymbol{\ast}},\, _{\boldsymbol{\ast}},\, ^{\boldsymbol{\times}})$ be an adjoint triple of pseudofunctors.
		Then \Cref{def:upper_star_lower_star_bc,def:lower_star_upper_cross_bc,def:upper_star_upper_cross_bc} can be defined in this case.
		Under duality, \eqref{rem:abstract_nonsense_for_base_change4}, it happens that $\vartheta_\square$ and $\xi_\square$ get swapped, while $\beta_\square$ is self-dual.
		The proof of \Cref{lem:upper_star_lower_star_pastes_horizontally} works verbatim in this generality.
		Therefore, by duality, the proof of \Cref{lem:lower_star_upper_cross_pastes_horizontally} follows on the nose in this more general setting.\footnote{This is probably what it is meant towards the end of proof \cite[\href{https://stacks.math.columbia.edu/tag/0E5F}{Tag 0E5F}]{StacksProject} when it is said that the pasting result for $\xi_\square$ is dual to that of $\vartheta_\square$.}
		Similarly, the proof of \Cref{lem:upper_star_upper_cross_pastes_horizontally} also works verbatim in this general case.
	\end{enumerate}
\end{remark}

%%%%%%%%%%%%%%%%%%%%%%%%%%%%%%%%%%%
\section{Derived reflexivity}
\label{sec:derived_reflexive}
%%%%%%%%%%%%%%%%%%%%%%%%%%%%%%%%%%

\begin{definition}
	\label{def:homologically_reflexive}
	Let $X$ be an algebraic space and $E,A\in D(X)$.
	We say that $E$ is \textbf{derived $A$-reflexive} if the following hold:
	\begin{enumerate}
		\item the canonical morphism
		\begin{displaymath}
			E\to \operatorname{\mathbb{R}\mathcal{H}\! \mathit{om}} ( \operatorname{\mathbb{R}\mathcal{H}\! \mathit{om}}(E, A) , A)
		\end{displaymath}
		is an isomorphism
		\item $E$ and $\operatorname{\mathbb{R}\mathcal{H}\! \mathit{om}}(E, A)$ are bounded pseudocoherent complexes.
	\end{enumerate}
\end{definition}

\begin{remark}
	\Cref{def:homologically_reflexive} is an analog of \cite[Definition 1.3.1]{Avramov/Iyengar/Lipman:2011}.
	For convenience, we recall that the canonical morphism is obtained via the following sequence of adjunctions:
	\begin{displaymath}
		\begin{aligned}
			\operatorname{Hom} ( \operatorname{\mathbb{R}\mathcal{H}\! \mathit{om}} (E , A) , \operatorname{\mathbb{R}\mathcal{H}\! \mathit{om}} (E , A) )
			&\cong \operatorname{Hom} ( \operatorname{\mathbb{R}\mathcal{H}\! \mathit{om}} (E , A) \otimes^{\mathbf{L}} E , A )
			\\&\cong \operatorname{Hom} ( E \otimes^{\mathbf{L}} \operatorname{\mathbb{R}\mathcal{H}\! \mathit{om}} (E , A) , A)
			\\&\cong \operatorname{Hom}(E, \operatorname{\mathbb{R}\mathcal{H}\! \mathit{om}} ( \operatorname{\mathbb{R}\mathcal{H}\! \mathit{om}}(E, A ) , A )).
		\end{aligned}
	\end{displaymath}
	Take the morphism obtained from the identity of $\operatorname{\mathbb{R}\mathcal{H}\! \mathit{om}} (E , A)$.
\end{remark}

\begin{lemma}
	\label{lem:derived_reflexive_by_covers}
	Let $X$ be a Noetherian algebraic space.
	Consider $E,A\in D^b_{\operatorname{coh}}(X)$ such that $E$ and $\operatorname{\mathbb{R}\mathcal{H}\! \mathit{om}}(E, A)$ are bounded pseudocoherent complexes.
	Then the following are equivalent:
	\begin{enumerate}
		\item \label{lem:derived_reflexive_by_covers1} $E$ is derived $A$-reflexive
		\item \label{lem:derived_reflexive_by_covers2} there exists an isomorphism
		\begin{displaymath}
			E\to \operatorname{\mathbb{R}\mathcal{H}\! \mathit{om}} ( \operatorname{\mathbb{R}\mathcal{H}\! \mathit{om}}(E, A) , A)
		\end{displaymath}
		\item \label{lem:derived_reflexive_by_covers3} $\mathbf{L}s^\ast E$ is derived $\mathbf{L}s^\ast A$-reflexive for every \'{e}tale morphism $s\colon Y \to X$ from a scheme
		\item \label{lem:derived_reflexive_by_covers4} there exists an \'{e}tale presentation $s\colon Y \to X$ such that $\mathbf{L}s^\ast E$ is derived $\mathbf{L}s^\ast A$-reflexive.
	\end{enumerate}
	Moreover, if $A$ and $B$ are isomorphic objects of $D(X)$, then $E$ is derived $A$-reflexive if, and only if, it is derived $B$-reflexive.
\end{lemma}

\begin{proof}
	It is immediate that $\eqref{lem:derived_reflexive_by_covers1} \implies \eqref{lem:derived_reflexive_by_covers2}$.
	We prove $\eqref{lem:derived_reflexive_by_covers2}\implies \eqref{lem:derived_reflexive_by_covers1}$.
	Suppose there exists an isomorphism
	\begin{displaymath}
		\alpha\colon E\to \operatorname{\mathbb{R}\mathcal{H}\! \mathit{om}} ( \operatorname{\mathbb{R}\mathcal{H}\! \mathit{om}}(E, A) , A).
	\end{displaymath}
	Choose an \'{e}tale presentation $s\colon U \to X$.
	Then it follows that $\mathbf{L}s^\ast E$ is isomorphic to $\operatorname{\mathbb{R}\mathcal{H}\! \mathit{om}} ( \operatorname{\mathbb{R}\mathcal{H}\! \mathit{om}}(\mathbf{L}s^\ast E, \mathbf{L}s^\ast A) , \mathbf{L}s^\ast A)$ \cite[\href{https://stacks.math.columbia.edu/tag/04LX}{Tag 04LX} and \href{https://stacks.math.columbia.edu/tag/08JB}{08JB}]{StacksProject}.
	By \cite[Proposition 1.3.3]{Avramov/Iyengar/Lipman:2011}, the canonical morphism
	\begin{displaymath}
		\mathbf{L}s^\ast E \to \operatorname{\mathbb{R}\mathcal{H}\! \mathit{om}} ( \operatorname{\mathbb{R}\mathcal{H}\! \mathit{om}}(\mathbf{L}s^\ast E, \mathbf{L}s^\ast A) , \mathbf{L}s^\ast A)
	\end{displaymath}
	is an isomorphism, and so,
	\begin{displaymath}
		\mathbf{L} s^\ast \operatorname{cone}(E\xrightarrow{ntrl.} \operatorname{\mathbb{R}\mathcal{H}\! \mathit{om}} ( \operatorname{\mathbb{R}\mathcal{H}\! \mathit{om}}(E, A) , A)) \cong 0.
	\end{displaymath}
	Hence, \Cref{lem:faithfully_flat_is_conservative} implies
	\begin{displaymath}
		\operatorname{cone}(E\xrightarrow{ntrl.} \operatorname{\mathbb{R}\mathcal{H}\! \mathit{om}} ( \operatorname{\mathbb{R}\mathcal{H}\! \mathit{om}}(E, A) , A)) \cong 0,
	\end{displaymath}
	and so $E$ is derived $A$-reflexive.

	To see that $\eqref{lem:derived_reflexive_by_covers2}\implies \eqref{lem:derived_reflexive_by_covers3}$, pullback the isomorphism, use \cite[\href{https://stacks.math.columbia.edu/tag/04LX}{Tag 04LX} and \href{https://stacks.math.columbia.edu/tag/08JB}{08JB}]{StacksProject}, and apply $\eqref{lem:derived_reflexive_by_covers2}$ on $Y$.
	Clearly, $\eqref{lem:derived_reflexive_by_covers3}\implies \eqref{lem:derived_reflexive_by_covers4}$, whereas $\eqref{lem:derived_reflexive_by_covers4}\implies \eqref{lem:derived_reflexive_by_covers1}$ is argued like above.

	Finally, we prove the last claim.
	Let $\alpha\colon A \to B$ be an isomorphism.
	If $E$ is derived $A$-reflexive, then the canonical morphism
	\begin{displaymath}
		E \to \operatorname{\mathbb{R}\mathcal{H}\! \mathit{om}} ( \operatorname{\mathbb{R}\mathcal{H}\! \mathit{om}}(E, A) , A)
	\end{displaymath}
	is an isomorphism.
	Moreover, $\alpha$ induces an isomorphism
	\begin{displaymath}
		\operatorname{\mathbb{R}\mathcal{H}\! \mathit{om}}(E, A) \xrightarrow{ \operatorname{\mathbb{R}\mathcal{H}\! \mathit{om}}(E, \alpha)}\operatorname{\mathbb{R}\mathcal{H}\! \mathit{om}}(E, B),
	\end{displaymath}
	and hence, an isomorphism
	\begin{displaymath}
		\operatorname{\mathbb{R}\mathcal{H}\! \mathit{om}} ( \operatorname{\mathbb{R}\mathcal{H}\! \mathit{om}}(E, A) , A) \to \operatorname{\mathbb{R}\mathcal{H}\! \mathit{om}} ( \operatorname{\mathbb{R}\mathcal{H}\! \mathit{om}}(E, B) , B).
	\end{displaymath}
	Consequently, we obtain an isomorphism
	\begin{displaymath}
		E \to \operatorname{\mathbb{R}\mathcal{H}\! \mathit{om}} ( \operatorname{\mathbb{R}\mathcal{H}\! \mathit{om}}(E, B) , B).
	\end{displaymath}
	By \eqref{lem:derived_reflexive_by_covers1}, $E$ is derived $B$-reflexive. The same argument gives the converse.
\end{proof}

\begin{example}
	\label{ex:relative_dualizing_is_bounded}
	Let $f\colon Y \to X$ be a proper flat morphism of Noetherian algebraic spaces.
	Then $f^! \mathcal{O}_X \in D^b_{\operatorname{coh}}(Y)$.
	Indeed, flatness of $f$ implies $\mathcal{O}_Y$ is $f$-perfect.
	Hence, by \Cref{prop:stacky_all_coincide}, $\mathbf{R}f_\ast \operatorname{Perf}(Y)\subseteq \operatorname{Perf}(X)$.
	By \Cref{lem:upper_shriek_coherent_cohomology}, $f^! \mathcal{O}_X \in D^+_{\operatorname{coh}}(Y)$.
	It suffices to check that $f^! \mathcal{O}_X \in D^-_{\operatorname{coh}}(Y)$.
	Choose a compact generator $G\in\operatorname{Perf}(Y)$.
	By adjunction, for any $n\in \mathbb{Z}$,
	\begin{displaymath}
		\operatorname{Hom}(G[n],f^! \mathcal{O}_X)
		\cong \operatorname{Hom}(\mathbf{R}f_\ast G[n] , \mathcal{O}_X).
	\end{displaymath}
	Applying \cite[\href{https://stacks.math.columbia.edu/tag/0GFH}{Tag 0GFH}]{StacksProject}, it follows that $\operatorname{Hom}(G[n],f^! \mathcal{O}_X) \cong 0$ for $n \ll 0$, and hence, $f^! \mathcal{O}_X \in D^-_{\operatorname{coh}}(Y)$.
\end{example}

\begin{lemma}
	\label{lem:involution_for_f_perfect}
	Let $f\colon Y\to X$ be a proper flat morphism of Noetherian algebraic spaces.
	If $E$ is an $f$-perfect pseudocoherent complex on $Y$, then there exists a canonical isomorphism
	\begin{displaymath}
		E\to \operatorname{\mathbf{R}\mathcal{H}\! \mathit{om}} ( \operatorname{\mathbf{R}\mathcal{H}\! \mathit{om}}(E, f^! \mathcal{O}_{X}) , f^! \mathcal{O}_{X}).
	\end{displaymath}
	In fact, $E$ is derived $f^! \mathcal{O}_{X}$-reflexive.
\end{lemma}

\begin{proof}
	By \Cref{prop:stacky_all_coincide}, $E$ is $f$-quasi-perfect.
	Applying \Cref{lem:quasi-perfect_bounded}, $\operatorname{\mathbb{R}\mathcal{H}\! \mathit{om}}(E, f^\times \mathcal{O}_X) \in D^b_{\operatorname{coh}}(Y)$.
	Then \Cref{prop:ballard_quasi_perfect_involution} yields an isomorphism
	\begin{displaymath}
		\nu\colon E \to \operatorname{\mathbf{R}\mathcal{H}\! \mathit{om}} (\operatorname{\mathbf{R}\mathcal{H}\! \mathit{om}} (E,f^\times \mathcal{O}_{X}) , f^\times \mathcal{O}_{X}).
	\end{displaymath}
	By \Cref{lem:neeman187}, there exists a natural isomorphism $f^\times \to f^!$.
	Hence, the first claim follows.

	By \Cref{lem:upper_shriek_coherent_cohomology}, $f^! D^+_{\operatorname{coh}}(X)\subseteq D^+_{\operatorname{coh}}(Y)$.
	In fact, \Cref{ex:relative_dualizing_is_bounded} implies $f^! \mathcal{O}_X\in D^b_{\operatorname{coh}}(Y)$.
	Then \Cref{lem:internal_hom} and \cite[\href{https://stacks.math.columbia.edu/tag/0A8A}{Tag 0A8A}]{StacksProject} yields an isomorphism
	\begin{displaymath}
		E \to \operatorname{\mathbb{R}\mathcal{H}\! \mathit{om}} (\operatorname{\mathbb{R}\mathcal{H}\! \mathit{om}} (E,f^\times \mathcal{O}_{X}) , f^\times \mathcal{O}_{X})
	\end{displaymath}
	Thus, by \Cref{lem:derived_reflexive_by_covers}, the claim for derived $f^! \mathcal{O}_X$-reflexivity follows.
\end{proof}

%%%%%%%%%%%%%%%%%%%%%%%%%%%%%%%%%%%
\section{Base change for adjoints}
\label{sec:base_change_for_adjoints}
%%%%%%%%%%%%%%%%%%%%%%%%%%%%%%%%%%

\begin{lemma}
	[Upper shriek base change]
	\label{lem:base_change_relative_formula}
	Consider a fibered square of Noetherian algebraic spaces
	\begin{displaymath}
		% https://q.uiver.app/#q=WzAsNCxbMSwwLCJZIl0sWzEsMSwiUyJdLFswLDAsIlleXFxwcmltZSJdLFswLDEsIlNeXFxwcmltZSJdLFswLDEsImYiXSxbMiwzLCJmXlxccHJpbWUiLDJdLFszLDEsImciLDJdLFsyLDAsImdeXFxwcmltZSJdXQ==
		\begin{tikzcd}
			{Y^\prime} & Y \\
			{S^\prime} & S
			\arrow["{g^\prime}", from=1-1, to=1-2]
			\arrow["{f^\prime}"', from=1-1, to=2-1]
			\arrow["f", from=1-2, to=2-2]
			\arrow["g"', from=2-1, to=2-2]
		\end{tikzcd}
	\end{displaymath}
	where $f$ is proper, finitely presented, and flat.
	Then the base change morphism \cite[\href{https://stacks.math.columbia.edu/tag/0E5D}{Tag 0E5D}]{StacksProject}
	\begin{displaymath}
		\beta_{f,g}(E)\colon
		\mathbf{L}(g^\prime)^\ast f^!E
		\to (f^\prime)^!\mathbf{L}g^\ast E
	\end{displaymath}
	is an isomorphism for every $E\in D_{\operatorname{qc}}(S)$.
\end{lemma}

\begin{proof}
	Let $s\colon U \to S$ be an \'{e}tale presentation from an affine scheme.
	Consider the fibered cube
	\begin{equation}
		\label{eq:base_change_relative_formula_cube}
		% https://q.uiver.app/#q=WzAsOCxbMywxLCJVIl0sWzMsMywiUy4iXSxbMSwzLCJZIl0sWzIsMiwiU15cXHByaW1lIl0sWzEsMSwiVV5cXHByaW1lIl0sWzIsMCwiViJdLFswLDAsIlZeXFxwcmltZSJdLFswLDIsIlleXFxwcmltZSJdLFswLDEsInMiXSxbMiwxLCJmIiwyXSxbMywxLCJnIl0sWzQsMCwiYiIsMCx7ImxhYmVsX3Bvc2l0aW9uIjo3MH1dLFs1LDAsImEiXSxbNCwyLCJzXzEiLDAseyJsYWJlbF9wb3NpdGlvbiI6NzB9XSxbNSwzLCJzXzIiLDIseyJsYWJlbF9wb3NpdGlvbiI6ODAsInN0eWxlIjp7ImJvZHkiOnsibmFtZSI6ImRhc2hlZCJ9fX1dLFs2LDQsImFeXFxwcmltZSIsMl0sWzYsNSwiYl5cXHByaW1lIiwyXSxbNywzLCJmXlxccHJpbWUiLDAseyJsYWJlbF9wb3NpdGlvbiI6MzAsInN0eWxlIjp7ImJvZHkiOnsibmFtZSI6ImRhc2hlZCJ9fX1dLFs3LDIsImdeXFxwcmltZSIsMl0sWzYsNywic15cXHByaW1lIiwyXV0=
		\begin{tikzcd}
			{V^\prime} && V & \\
			& {U^\prime} && U \\
			{Y^\prime} && {S^\prime} \\
			& Y && {S.}
			\arrow["{b^\prime}"', from=1-1, to=1-3]
			\arrow["{a^\prime}"', from=1-1, to=2-2]
			\arrow["{s^\prime}"', from=1-1, to=3-1]
			\arrow["a", from=1-3, to=2-4]
			\arrow["{s_2}"'{pos=0.8}, dashed, from=1-3, to=3-3]
			\arrow["b"{pos=0.7}, from=2-2, to=2-4]
			\arrow["{s_1}"{pos=0.7}, from=2-2, to=4-2]
			\arrow["s", from=2-4, to=4-4]
			\arrow["{f^\prime}"{pos=0.3}, dashed, from=3-1, to=3-3]
			\arrow["{g^\prime}"', from=3-1, to=4-2]
			\arrow["g", from=3-3, to=4-4]
			\arrow["f"', from=4-2, to=4-4]
		\end{tikzcd}
	\end{equation}
	As $s$ is finitely presented, every algebraic space in the fibered cube is Noetherian.
	By \Cref{lem:quasi_affine_diagonal} and \cite[\href{https://stacks.math.columbia.edu/tag/03HA}{Tag 03HA} and \href{https://stacks.math.columbia.edu/tag/01SL}{01SL}]{StacksProject}, $s$ is separated.
	Hence, its base changes $s_1,s_2,s^\prime$ are separated, \'{e}tale, and surjective.
	As $f,f^\prime,b,b^\prime$ are proper, \Cref{lem:neeman187} gives natural
	isomorphisms $f^\times\to f^!$, $(f^\prime)^\times \to (f^\prime)^!$, $b^\times \to b^!$, and $(b^\prime)^\times \to (b^\prime)^!$.
	Throughout the proof, the base change morphisms involving upper
	shriek functors are the base change morphisms of \cite[\href{https://stacks.math.columbia.edu/tag/0E5D}{Tag 0E5D}]{StacksProject} induced from these natural isomorphisms.

	We now compare the two decompositions of the outer fibered
	rectangle
	\begin{equation}
		\label{diag:outer_rectangle_relative_dualizing}
		% https://q.uiver.app/#q=WzAsNCxbMCwwLCJWXlxccHJpbWUiXSxbMiwwLCJZIl0sWzAsMSwiViJdLFsyLDEsIlMuIl0sWzAsMiwiYl5cXHByaW1lIiwyXSxbMSwzLCJmIl0sWzAsMSwiZ15cXHByaW1lIFxcY2lyYyBzXlxccHJpbWUgPSBzXzEgXFxjaXJjIGFeXFxwcmltZSJdLFsyLDMsImcgXFxjaXJjIHNfMj0gcyBcXGNpcmMgYSIsMl1d
		\begin{tikzcd}
			{V^\prime} && Y \\
			V && {S.}
			\arrow["{g^\prime \circ s^\prime = s_1 \circ a^\prime}", from=1-1, to=1-3]
			\arrow["{b^\prime}"', from=1-1, to=2-1]
			\arrow["f", from=1-3, to=2-3]
			\arrow["{g \circ s_2= s \circ a}"', from=2-1, to=2-3]
		\end{tikzcd}
	\end{equation}
	All the fibered squares which occur in \eqref{eq:base_change_relative_formula_cube} are tor-independent.
	Indeed, $f$ is flat, and hence by base change, $f^\prime,b,b^\prime$ are flat.
	The first decomposition of
	\eqref{diag:outer_rectangle_relative_dualizing}
	is given by the squares $(f,g)$ and $(f^\prime,s_2)$ in \eqref{eq:base_change_relative_formula_cube}.
	By \Cref{lem:upper_star_upper_cross_pastes_horizontally}, the base change morphism
	for the outer rectangle is the composite
	\begin{equation}
		\label{eq:first_decomposition_relative_dualizing}
		\begin{aligned}
			\mathbf{L}(s^\prime)^\ast \mathbf{L}(g^\prime)^\ast f^\times
			&\xrightarrow{\mathbf{L}(s^\prime)^\ast\beta_{f,g}} \mathbf{L}(s^\prime)^\ast (f^\prime)^\times\mathbf{L}g^\ast
			\\&\xrightarrow{\beta_{f^\prime,s_2}\mathbf{L}g^\ast} (b^\prime)^\times \mathbf{L}s_2^\ast\mathbf{L}g^\ast
			\\&\xrightarrow{\cong} (b^\prime)^\times \mathbf{L}a^\ast\mathbf{L}s^\ast .
		\end{aligned}
	\end{equation}
	The last arrow is the canonical pseudofunctoriality isomorphism induced by $g\circ s_2=s\circ a$.
	Moreover, the second decomposition of \eqref{diag:outer_rectangle_relative_dualizing} is given by the squares $(f,s)$ and $(b,a)$ in \eqref{eq:base_change_relative_formula_cube}.
	Again, by \Cref{lem:upper_star_upper_cross_pastes_horizontally}, the same outer base change morphism is the composite
	\begin{equation}
		\label{eq:second_decomposition_relative_dualizing}
		\begin{aligned}
			\mathbf{L}(s^\prime)^\ast \mathbf{L}(g^\prime)^\ast f^\times
			&\xrightarrow{\cong} \mathbf{L}(a^\prime)^\ast
			\mathbf{L}s_1^\ast f^\times
			\\&\xrightarrow{\mathbf{L}(a^\prime)^\ast\beta_{f,s}}
			\mathbf{L}(a^\prime)^\ast b^\times\mathbf{L}s^\ast
			\\&\xrightarrow{\beta_{b,a}\mathbf{L}s^\ast} (b^\prime)^\times \mathbf{L}a^\ast\mathbf{L}s^\ast .
		\end{aligned}
	\end{equation}
	Here the first arrow is the canonical pseudofunctoriality
	isomorphism induced by $g^\prime\circ s^\prime=s_1\circ a^\prime$.
	Since \eqref{eq:first_decomposition_relative_dualizing} and \eqref{eq:second_decomposition_relative_dualizing} are both the base change morphism for the same outer rectangle, the two composites are equal.
	In particular, the following diagram
	\begin{equation}
		\label{eq:decomposition_relative_dualizing}
		% https://q.uiver.app/#q=WzAsNixbMCwxLCJcXG1hdGhiZntMfShzXlxccHJpbWUpXlxcYXN0IFxcbWF0aGJme0x9KGdeXFxwcmltZSleXFxhc3QgZl5cXHRpbWVzIl0sWzAsMCwiXFxtYXRoYmZ7TH0oc15cXHByaW1lKV5cXGFzdCAoZl5cXHByaW1lKV5cXHRpbWVzXFxtYXRoYmZ7TH1nXlxcYXN0Il0sWzIsMCwiKGJeXFxwcmltZSleXFx0aW1lcyBcXG1hdGhiZntMfXNfMl5cXGFzdFxcbWF0aGJme0x9Z15cXGFzdCJdLFsyLDEsIihiXlxccHJpbWUpXlxcdGltZXMgXFxtYXRoYmZ7TH1hXlxcYXN0XFxtYXRoYmZ7TH1zXlxcYXN0ICJdLFswLDIsIlxcbWF0aGJme0x9KGFeXFxwcmltZSleXFxhc3QgICAgICAgICAgICAgXFxtYXRoYmZ7TH1zXzFeXFxhc3QgZl5cXHRpbWVzIl0sWzIsMiwiXFxtYXRoYmZ7TH0oYV5cXHByaW1lKV5cXGFzdCBiXlxcdGltZXNcXG1hdGhiZntMfXNeXFxhc3QiXSxbMCwxLCJcXG1hdGhiZntMfShzXlxccHJpbWUpXlxcYXN0XFxiZXRhX3tmLGd9Il0sWzEsMiwiXFxiZXRhX3tmXlxccHJpbWUsc18yfVxcbWF0aGJme0x9Z15cXGFzdCJdLFsyLDMsIlxcY29uZyJdLFswLDQsIlxcY29uZyIsMl0sWzQsNSwiXFxtYXRoYmZ7TH0oYV5cXHByaW1lKV5cXGFzdFxcYmV0YV97ZixzfSIsMl0sWzUsMywiXFxiZXRhX3tiLGF9XFxtYXRoYmZ7TH1zXlxcYXN0IiwyXV0=
		\begin{tikzcd}
			{\mathbf{L}(s^\prime)^\ast (f^\prime)^\times\mathbf{L}g^\ast} && {(b^\prime)^\times \mathbf{L}s_2^\ast\mathbf{L}g^\ast} \\
			{\mathbf{L}(s^\prime)^\ast \mathbf{L}(g^\prime)^\ast f^\times} && {(b^\prime)^\times \mathbf{L}a^\ast\mathbf{L}s^\ast } \\
			{\mathbf{L}(a^\prime)^\ast             \mathbf{L}s_1^\ast f^\times} && {\mathbf{L}(a^\prime)^\ast b^\times\mathbf{L}s^\ast}
			\arrow["{\beta_{f^\prime,s_2}\mathbf{L}g^\ast}", from=1-1, to=1-3]
			\arrow["\cong", from=1-3, to=2-3]
			\arrow["{\mathbf{L}(s^\prime)^\ast\beta_{f,g}}", from=2-1, to=1-1]
			\arrow["\cong"', from=2-1, to=3-1]
			\arrow["{\mathbf{L}(a^\prime)^\ast\beta_{f,s}}"', from=3-1, to=3-3]
			\arrow["{\beta_{b,a}\mathbf{L}s^\ast}"', from=3-3, to=2-3]
		\end{tikzcd}
	\end{equation}
	is commutative.

	Fix $E\in D_{\operatorname{qc}}(S)$.
	Since $f,f^\prime,b,b^\prime$ are proper, \Cref{lem:neeman187} gives natural
	isomorphisms $f^\times\to f^!$, $(f^\prime)^\times \to (f^\prime)^!$, $b^\times \to b^!$, and $(b^\prime)^\times \to (b^\prime)^!$.
	Then coupling these natural isomorphisms with \eqref{eq:decomposition_relative_dualizing} gives a commutative diagram
	\begin{equation}
		\label{diag:relative_dualizing_base_change}
		% https://q.uiver.app/#q=WzAsNCxbMCwwLCJcXG1hdGhiZntMfShzXlxccHJpbWUpXlxcYXN0IFxcbWF0aGJme0x9KGdeXFxwcmltZSleXFxhc3QgZl4hIEUiXSxbMiwxLCIoYl5cXHByaW1lKV4hIFxcbWF0aGJme0x9YV5cXGFzdFxcbWF0aGJme0x9c15cXGFzdCBFLiJdLFswLDEsIlxcbWF0aGJme0x9KGFeXFxwcmltZSleXFxhc3QgYl4hIFxcbWF0aGJme0x9c15cXGFzdCBFIl0sWzIsMCwiXFxtYXRoYmZ7TH0oc15cXHByaW1lKV5cXGFzdCAoZl5cXHByaW1lKV4hIFxcbWF0aGJme0x9Z15cXGFzdCBFIl0sWzIsMSwiXFxiZXRhX3tiLGF9XFxtYXRoYmZ7TH1zXlxcYXN0IEUiLDJdLFswLDMsIlxcbWF0aGJme0x9KHNeXFxwcmltZSleXFxhc3RcXGJldGFfe2YsZ30oRSkiXSxbMCwyLCJcXGNvbmciLDJdLFszLDEsIlxcY29uZyJdXQ==
		\begin{tikzcd}
			{\mathbf{L}(s^\prime)^\ast \mathbf{L}(g^\prime)^\ast f^! E} && {\mathbf{L}(s^\prime)^\ast (f^\prime)^! \mathbf{L}g^\ast E} \\
			{\mathbf{L}(a^\prime)^\ast b^! \mathbf{L}s^\ast E} && {(b^\prime)^! \mathbf{L}a^\ast\mathbf{L}s^\ast E.}
			\arrow["{\mathbf{L}(s^\prime)^\ast\beta_{f,g}(E)}", from=1-1, to=1-3]
			\arrow["\cong"', from=1-1, to=2-1]
			\arrow["\cong", from=1-3, to=2-3]
			\arrow["{\beta_{b,a}\mathbf{L}s^\ast E}"', from=2-1, to=2-3]
		\end{tikzcd}
	\end{equation}

	Let us spell out the two vertical arrows.
	The left vertical morphism is the composite
	\begin{displaymath}
		\begin{aligned}
			\mathbf{L}(s^\prime)^\ast \mathbf{L}(g^\prime)^\ast f^!E
			&\xrightarrow{\cong} \mathbf{L}(a^\prime)^\ast \mathbf{L}s_1^\ast f^!E
			\\&\xrightarrow{\mathbf{L}(a^\prime)^\ast \beta_{f,s}(E)} \mathbf{L}(a^\prime)^\ast b^!\mathbf{L}s^\ast E.
		\end{aligned}
	\end{displaymath}
	The right vertical morphism is
	\begin{displaymath}
		\begin{aligned}
			\mathbf{L}(s^\prime)^\ast (f^\prime)^!\mathbf{L}g^\ast E
			&\xrightarrow{\beta_{f^\prime,s_2}(\mathbf{L}g^\ast E)}
			(b^\prime)^!\mathbf{L}s_2^\ast\mathbf{L}g^\ast E
			\\&\xrightarrow{\cong} (b^\prime)^! \mathbf{L}a^\ast\mathbf{L}s^\ast E.
		\end{aligned}
	\end{displaymath}
	Both are isomorphisms.
	Indeed, there are fibered squares
	\begin{displaymath}
		% https://q.uiver.app/#q=WzAsNCxbMCwxLCJVIl0sWzEsMSwiUyJdLFsxLDAsIlkiXSxbMCwwLCJVXlxccHJpbWUiXSxbMCwxLCJzIiwyXSxbMiwxLCJmIl0sWzMsMCwiYiIsMix7ImxhYmVsX3Bvc2l0aW9uIjo3MH1dLFszLDIsInNfMSIsMCx7ImxhYmVsX3Bvc2l0aW9uIjo3MH1dXQ==
		\begin{tikzcd}
			{U^\prime} & Y \\
			U & S
			\arrow["{s_1}", from=1-1, to=1-2]
			\arrow["b"', from=1-1, to=2-1]
			\arrow["f", from=1-2, to=2-2]
			\arrow["s"', from=2-1, to=2-2]
		\end{tikzcd}
	\end{displaymath}
	and
	\begin{displaymath}
		% https://q.uiver.app/#q=WzAsNCxbMSwxLCJTXlxccHJpbWUiXSxbMCwxLCJWIl0sWzAsMCwiVl5cXHByaW1lIl0sWzEsMCwiWV5cXHByaW1lIl0sWzEsMCwic18yIiwyXSxbMiwxLCJiXlxccHJpbWUiLDJdLFszLDAsImZeXFxwcmltZSIsMCx7ImxhYmVsX3Bvc2l0aW9uIjozMH1dLFsyLDMsInNeXFxwcmltZSJdXQ==
		\begin{tikzcd}
			{V^\prime} & {Y^\prime} \\
			V & {S^\prime .}
			\arrow["{s^\prime}", from=1-1, to=1-2]
			\arrow["{b^\prime}"', from=1-1, to=2-1]
			\arrow["{f^\prime}"{pos=0.3}, from=1-2, to=2-2]
			\arrow["{s_2}"', from=2-1, to=2-2]
		\end{tikzcd}
	\end{displaymath}
	Applying \Cref{lem:neeman188} to these fibered squares yields respectively natural isomorphisms
	\begin{displaymath}
		\mathbf{L}s_1^\ast f^\times E
		\to b^\times\mathbf{L}s^\ast E
	\end{displaymath}
	and
	\begin{displaymath}
		\mathbf{L}(s^\prime)^\ast (f^\prime)^\times\mathbf{L}g^\ast E
		\to (b^\prime)^\times \mathbf{L}s_2^\ast\mathbf{L}g^\ast E.
	\end{displaymath}
	Here $s$ and $s_2$ are flat, $b$ and $b^\prime$ are proper and flat.
	By \Cref{lem:neeman187}, these isomorphisms induce the vertical isomorphisms in \eqref{diag:relative_dualizing_base_change}.
	Indeed, we can appeal to the case of proper morphisms of finite tor-dimension, which allows for all objects of $D_{\operatorname{qc}}$.

	It remains to consider the lower horizontal morphism.
	Choose an \'{e}tale presentation $r\colon W\to V$ from an affine scheme.
	Set $W^\prime:=V^\prime\times_V W$.
	Denote by $r^\prime\colon W^\prime\to V^\prime$ the natural morphism.
	Write $c:=a\circ r\colon W\to U$ and $b_W\colon W^\prime\to W$ for the natural projection.
	Thus, there is a diagram of fibered squares
	\begin{displaymath}
		\begin{tikzcd}
			{W^\prime} & {V^\prime} & {U^\prime} \\
			W & V & U.
			\arrow["{r^\prime}", from=1-1, to=1-2]
			\arrow["{b_W}"', from=1-1, to=2-1]
			\arrow["{a^\prime}", from=1-2, to=1-3]
			\arrow["{b^\prime}"', from=1-2, to=2-2]
			\arrow["b", from=1-3, to=2-3]
			\arrow["r"', from=2-1, to=2-2]
			\arrow["a"', from=2-2, to=2-3]
		\end{tikzcd}
	\end{displaymath}
	Since $W$ and $U$ are affine schemes and $b\colon U^\prime\to U$ is proper, flat, and finitely presented, \cite[\href{https://stacks.math.columbia.edu/tag/0E5I}{Tag 0E5I}]{StacksProject} implies that $\beta_{b,c}(\mathbf{L}s^\ast E)$ is an isomorphism.
	Moreover, $r$ is flat and $b^\prime$ is proper and flat.
	Hence, \Cref{lem:neeman188} shows that $\beta_{b^\prime,r} (\mathbf{L}a^\ast\mathbf{L}s^\ast E)$ is an isomorphism.
	By \Cref{lem:upper_star_upper_cross_pastes_horizontally}, under the canonical pseudofunctorial identifications, one has
	\begin{displaymath}
		\beta_{b,c}(\mathbf{L}s^\ast E)
		=
		\beta_{b^\prime,r}
		(\mathbf{L}a^\ast\mathbf{L}s^\ast E)
		\circ
		\mathbf{L}(r^\prime)^\ast
		(\beta_{b,a}(\mathbf{L}s^\ast E)).
	\end{displaymath}
	It follows that $\mathbf{L}(r^\prime)^\ast (\beta_{b,a}(\mathbf{L}s^\ast E))$ is an isomorphism.
	Since $r^\prime$ is \'{e}tale and surjective, \Cref{lem:zero_object_via_covers} implies that $\beta_{b,a}(\mathbf{L}s^\ast E)$ is an isomorphism.
	The commutativity of \eqref{diag:relative_dualizing_base_change} yields that $\mathbf{L}(s^\prime)^\ast\beta_{f,g}(E)$ is an isomorphism.

	Now we finish the proof.
	By base change, $s^\prime$ is \'{e}tale and surjective, and hence, faithfully flat.
	Therefore, by \Cref{lem:zero_object_via_covers}, $\beta_{f,g}(E)$ is an isomorphism.
	Indeed, $\mathbf{L}(s^\prime)^\ast\beta_{f,g}(E)$ being an isomorphism implies
	\begin{displaymath}
		0\cong \mathbf{L}(s^\prime)^\ast \operatorname{cone}(\beta_{f,g}(E)) \cong \operatorname{cone}(\mathbf{L}(s^\prime)^\ast\beta_{f,g}(E)).
	\end{displaymath}
	This completes the proof.
\end{proof}

\begin{lemma}
	[Base change for relative duals]
	\label{lem:base_change_relative_dualizing_complex}
	Consider a fibered square of Noetherian algebraic spaces
	\begin{displaymath}
		\begin{tikzcd}
			{X^\prime} & {S^\prime} \\
			X & S
			\arrow["{f^\prime}", from=1-1, to=1-2]
			\arrow["{g^\prime}"', from=1-1, to=2-1]
			\arrow["g", from=1-2, to=2-2]
			\arrow["f"', from=2-1, to=2-2]
		\end{tikzcd}
	\end{displaymath}
	where $f$ is proper and flat.
	If $g$ is quasi-affine and $E\in D^b_{\operatorname{coh}}(X)$ is
	$f$-perfect, then the canonical morphism
	\begin{equation}
		\label{eq:base_change_relative_dual}
		\begin{aligned}
			\theta_E
			&\colon
			\mathbf{L}(g^\prime)^\ast
			\operatorname{\mathbf{R}\mathcal{H}\! \mathit{om}}
			(E,f^!\mathcal{O}_S)
			\to
			\operatorname{\mathbf{R}\mathcal{H}\! \mathit{om}}
			(\mathbf{L}(g^\prime)^\ast E,
			\mathbf{L}(g^\prime)^\ast f^!\mathcal{O}_S)
			\\
			&\xrightarrow{
			\operatorname{\mathbf{R}\mathcal{H}\! \mathit{om}}
			(1,\beta_{f,g}(\mathcal{O}_S))}
			\operatorname{\mathbf{R}\mathcal{H}\! \mathit{om}}
			(\mathbf{L}(g^\prime)^\ast E,
			(f^\prime)^!\mathcal{O}_{S^\prime})
		\end{aligned}
	\end{equation}
	is an isomorphism, where $\beta_{f,g}(\mathcal{O}_S)$ is the
	morphism of \Cref{lem:base_change_relative_formula}.
\end{lemma}

\begin{proof}
	Since $f$ and $f^\prime$ are proper, \Cref{lem:neeman187} gives natural isomorphisms $f^\times\cong f^!$ and $(f^\prime)^\times\cong(f^\prime)^!$.
	We use these identifications throughout the proof.

	We first explain the morphisms occurring in \eqref{eq:base_change_relative_dual}.
	The second morphism is defined, and so we discuss the first morphism.
	By \Cref{lem:upper_shriek_preserves_Dplusqc}, $f^!\mathcal{O}_S \in D^+_{\operatorname{qc}}(X)$ and $(f^\prime)^!\mathcal{O}_{S^\prime} \in D^+_{\operatorname{qc}}(X^\prime)$.
	Moreover, \Cref{lem:base_change_relative_formula} gives an isomorphism
	\begin{displaymath}
		\mathbf{L}(g^\prime)^\ast f^!\mathcal{O}_S
		\xrightarrow{\beta_{f,g}(\mathcal{O}_S)} (f^\prime)^!\mathbf{L}g^\ast\mathcal{O}_S
		\cong (f^\prime)^!\mathcal{O}_{S^\prime}.
	\end{displaymath}
	Hence, $\mathbf{L}(g^\prime)^\ast f^!\mathcal{O}_S$ also has bounded below and quasi-coherent cohomology.
	Derived pullback preserves pseudocoherence \cite[\href{https://stacks.math.columbia.edu/tag/08H4}{Tag 08H4}]{StacksProject}.
	Hence, \cite[\href{https://stacks.math.columbia.edu/tag/0A8A}{Tag 0A8A}]{StacksProject} and \Cref{lem:internal_hom} identify the internal Homs occurring in \eqref{eq:base_change_relative_dual} with the corresponding
	derived sheaf Homs.
	Under these identifications, the first morphism in \eqref{eq:base_change_relative_dual} is the canonical pullback-Hom morphism of
	\cite[\href{https://stacks.math.columbia.edu/tag/08JF}{Tag 08JF}]{StacksProject}.

	Choose a compact generator $P$ of $D_{\operatorname{qc}}(X)$.
	Since $g^\prime$ is quasi-affine, $\mathbf{L}(g^\prime)^\ast P$ is a compact generator of $D_{\operatorname{qc}}(X^\prime)$ by \Cref{lem:pullback_compact_generator}.
	By \Cref{lem:zero_iff_pushforward_zero}, it is enough to show:
	\begin{equation}
		\label{eq:pushforward_generator_test}
		\mathbf{R}f^\prime_\ast
		(
		\operatorname{\mathbf{R}\mathcal{H}\! \mathit{om}}
		(\mathbf{L}(g^\prime)^\ast P,\mathcal{O}_{X^\prime})
		\otimes^{\mathbf{L}}
		\operatorname{cone}(\theta_E)
		)
		\cong 0.
	\end{equation}
	%%NOTE:
	% for every $n\in\mathbb{Z}$,
	% \begin{displaymath}
	%     \begin{aligned}
	%         &\operatorname{Hom}(\mathbf{L}(g^\prime)^\ast P, \operatorname{cone}(\theta_E)[n])
	%         \\&\cong \operatorname{Hom}(\mathcal{O}_{X^\prime}, \operatorname{\mathbf{R}\mathcal{H}\! \mathit{om}}(\mathbf{L}(g^\prime)^\ast P,\mathcal{O}_{X^\prime})\otimes^{\mathbf{L}}  \operatorname{cone}(\theta_E)[n])
	%         \\&\cong \operatorname{Hom}_{S^\prime} (\mathcal{O}_{S^\prime}, \mathbf{R}f^\prime_\ast (  \operatorname{\mathbf{R}\mathcal{H}\! \mathit{om}} (\mathbf{L}(g^\prime)^\ast P,\mathcal{O}_{X^\prime}) \otimes^{\mathbf{L}} \operatorname{cone}(\theta_E))[n])
	%         \\&\cong 0.
	%     \end{aligned}
	% \end{displaymath}
	% Since $\mathbf{L}(g^\prime)^\ast P$ compactly generates $D_{\operatorname{qc}}(X^\prime)$, this implies $\operatorname{cone}(\theta_E)\cong0$.

	By \cite[\href{https://stacks.math.columbia.edu/tag/08J9}{Tag 08J9} and \href{https://stacks.math.columbia.edu/tag/08JJ}{08JJ}]{StacksProject}
	and \Cref{lem:gortz_wedhorn_internal_hom}, the morphism
	\begin{displaymath}
		\operatorname{\mathbf{R}\mathcal{H}\! \mathit{om}}
		(\mathbf{L}(g^\prime)^\ast P,\mathcal{O}_{X^\prime})
		\otimes^{\mathbf{L}}\theta_E
	\end{displaymath}
	identifies with $\theta_{P\otimes^{\mathbf{L}}E}$.
	Consequently, it suffices to show that $\mathbf{R}f^\prime_\ast (\theta_{P\otimes^{\mathbf{L}}E})$ is an isomorphism.
	Since $E$ is $f$-perfect and $P$ is perfect, \Cref{lem:f_perf_implies_perf_preserved} gives $\mathbf{R}f_\ast(P\otimes^{\mathbf{L}}E) \in\operatorname{Perf}(S)$.
	For any $A\in D^-_{\operatorname{coh}}(X)$, denote by
	\begin{displaymath}
		\xi_f(A)\colon
		\mathbf{R}f_\ast
		\operatorname{\mathbf{R}\mathcal{H}\! \mathit{om}} (A,f^!\mathcal{O}_S)
		\to \operatorname{\mathbf{R}\mathcal{H}\! \mathit{om}}(\mathbf{R}f_\ast A,\mathcal{O}_S)
	\end{displaymath}
	the duality morphism from \cite[\href{https://stacks.math.columbia.edu/tag/0E58}{Tag 0E58} and \href{https://stacks.math.columbia.edu/tag/0GG4}{0GG4}]{StacksProject}.
	Define $\xi_{f^\prime}$ similarly.
	These are isomorphisms by \cite[\href{https://stacks.math.columbia.edu/tag/0GG5}{Tag 0GG5}]{StacksProject}.
	By the construction of $\xi_f(P\otimes^{\mathbf{L}}E)$ \cite[\href{https://stacks.math.columbia.edu/tag/0E57}{Tag 0E57}, \href{https://stacks.math.columbia.edu/tag/0B6C}{0B6C}, and \href{https://stacks.math.columbia.edu/tag/0E5K}{0E5K}]{StacksProject}, its adjunct is the composite
	\begin{equation}
		\label{eq:duality_via_trace}
		\begin{aligned}
			&\mathbf{R}f_\ast \operatorname{\mathbf{R}\mathcal{H}\! \mathit{om}}
			(P\otimes^{\mathbf{L}}E,f^!\mathcal{O}_S) \otimes^{\mathbf{L}}
			\mathbf{R}f_\ast(P\otimes^{\mathbf{L}}E)
			\\&\to\mathbf{R}f_\ast(  \operatorname{\mathbf{R}\mathcal{H}\! \mathit{om}} (P\otimes^{\mathbf{L}}E,f^!\mathcal{O}_S) \otimes^{\mathbf{L}} (P\otimes^{\mathbf{L}}E))
			\\&\to \mathbf{R}f_\ast f^!\mathcal{O}_S \xrightarrow{\operatorname{Tr}_f} \mathcal{O}_S.
		\end{aligned}
	\end{equation}
	The analogous description holds for $\xi_{f^\prime} (\mathbf{L}(g^\prime)^\ast(P\otimes^{\mathbf{L}}E))$.

	For $A\in D_{\operatorname{qc}}(X)$, denote by
	\begin{displaymath}
		\alpha_A\colon
		\mathbf{L}g^\ast\mathbf{R}f_\ast A
		\to \mathbf{R}f^\prime_\ast \mathbf{L}(g^\prime)^\ast A
	\end{displaymath}
	the canonical pullback-pushforward base change morphism \cite[\href{https://stacks.math.columbia.edu/tag/07A7}{Tag 07A7}]{StacksProject}.
	Since $f$ is flat, $\alpha_A$ is an isomorphism by \cite[\href{https://stacks.math.columbia.edu/tag/08IR}{Tag 08IR}]{StacksProject}.

	We claim that the following diagram commutes:
	\begin{equation}
		\label{diag:relative_dual_base_change_compatibility}
		\begin{tikzcd}
			{\mathbf{L}g^\ast\mathbf{R}f_\ast \operatorname{\mathbf{R}\mathcal{H}\! \mathit{om}} (P\otimes^{\mathbf{L}}E,f^!\mathcal{O}_S)} & \\
			{\mathbf{R}f^\prime_\ast\mathbf{L}(g^\prime)^\ast \operatorname{\mathbf{R}\mathcal{H}\! \mathit{om}} (P\otimes^{\mathbf{L}}E,f^!\mathcal{O}_S)} \\
			& {\mathbf{L}g^\ast \operatorname{\mathbf{R}\mathcal{H}\! \mathit{om}} (     \mathbf{R}f_\ast(P\otimes^{\mathbf{L}}E),     \mathcal{O}_S )} \\
			{\mathbf{R}f^\prime_\ast \operatorname{\mathbf{R}\mathcal{H}\! \mathit{om}} (     \mathbf{L}(g^\prime)^\ast     (P\otimes^{\mathbf{L}}E),     (f^\prime)^!\mathcal{O}_{S^\prime} )} & {\operatorname{\mathbf{R}\mathcal{H}\! \mathit{om}} (     \mathbf{L}g^\ast     \mathbf{R}f_\ast(P\otimes^{\mathbf{L}}E),     \mathcal{O}_{S^\prime} )} \\
			\\
			{\operatorname{\mathbf{R}\mathcal{H}\! \mathit{om}} (     \mathbf{R}f^\prime_\ast     \mathbf{L}(g^\prime)^\ast     (P\otimes^{\mathbf{L}}E),     \mathcal{O}_{S^\prime} ).}
			\arrow["{\alpha_{\operatorname{\mathbf{R}\mathcal{H}\! \mathit{om}}     (P\otimes^{\mathbf{L}}E,f^!\mathcal{O}_S)}}"', from=1-1, to=2-1]
			\arrow["{\mathbf{L}g^\ast     \xi_f(P\otimes^{\mathbf{L}}E)}", from=1-1, to=3-2]
			\arrow["{\mathbf{R}f^\prime_\ast     (\theta_{P\otimes^{\mathbf{L}}E})}"', from=2-1, to=4-1]
			\arrow[from=3-2, to=4-2]
			\arrow["{\xi_{f^\prime}     (\mathbf{L}(g^\prime)^\ast     (P\otimes^{\mathbf{L}}E))}"', from=4-1, to=6-1]
			\arrow["{\operatorname{\mathbf{R}\mathcal{H}\! \mathit{om}}     (\alpha^{-1}_{P\otimes^{\mathbf{L}}E},1)}", from=4-2, to=6-1]
		\end{tikzcd}
	\end{equation}
	The unlabeled right vertical morphism is the canonical pullback-Hom morphism \cite[\href{https://stacks.math.columbia.edu/tag/08JF}{Tag 08JF}]{StacksProject}.
	It is an isomorphism by  \Cref{lem:gortz_wedhorn_internal_hom}, since
	$\mathbf{R}f_\ast(P\otimes^{\mathbf{L}}E)$ is perfect.

	We verify the commutativity of \eqref{diag:relative_dual_base_change_compatibility}.
	By tensor-Hom adjunction \cite[\href{https://stacks.math.columbia.edu/tag/08J7}{Tag 08J7}]{StacksProject}, it suffices to show that the two boundary morphisms give the same morphism
	\begin{displaymath}
		\mathbf{L}g^\ast\mathbf{R}f_\ast \operatorname{\mathbf{R}\mathcal{H}\! \mathit{om}} (P\otimes^{\mathbf{L}}E,f^!\mathcal{O}_S) \otimes^{\mathbf{L}} \mathbf{R}f^\prime_\ast \mathbf{L}(g^\prime)^\ast
		(P\otimes^{\mathbf{L}}E)
		\to\mathcal{O}_{S^\prime}.
	\end{displaymath}
	The morphism obtained from right boundary is
	\begin{displaymath}
		\begin{aligned}
			& \mathbf{L}g^\ast\mathbf{R}f_\ast \operatorname{\mathbf{R}\mathcal{H}\! \mathit{om}} (P\otimes^{\mathbf{L}}E,f^!\mathcal{O}_S)
			\otimes^{\mathbf{L}} \mathbf{R}f^\prime_\ast
			\mathbf{L}(g^\prime)^\ast
			(P\otimes^{\mathbf{L}}E)
			\\&\xrightarrow{1\otimes  \alpha_{P\otimes^{\mathbf{L}}E}^{-1}}\mathbf{L}g^\ast\mathbf{R}f_\ast \operatorname{\mathbf{R}\mathcal{H}\! \mathit{om}} (P\otimes^{\mathbf{L}}E,f^!\mathcal{O}_S) \otimes^{\mathbf{L}} \mathbf{L}g^\ast\mathbf{R}f_\ast (P\otimes^{\mathbf{L}}E)
			\\&\to \mathbf{L}g^\ast\mathbf{R}f_\ast f^!\mathcal{O}_S \xrightarrow{\mathbf{L}g^\ast\operatorname{Tr}_f} \mathcal{O}_{S^\prime},
		\end{aligned}
	\end{displaymath}
	where the unlabeled morphism is the pullback of the relative cup
	product followed by evaluation.
	The morphism obtained from the left boundary is
	\begin{displaymath}
		\begin{aligned}
			&\mathbf{L}g^\ast\mathbf{R}f_\ast \operatorname{\mathbf{R}\mathcal{H}\! \mathit{om}} (P\otimes^{\mathbf{L}}E,f^!\mathcal{O}_S) \otimes^{\mathbf{L}} \mathbf{R}f^\prime_\ast \mathbf{L}(g^\prime)^\ast (P\otimes^{\mathbf{L}}E)
			\\&\xrightarrow{\alpha_{\operatorname{\mathbf{R}\mathcal{H}\! \mathit{om}}(P\otimes^{\mathbf{L}}E,f^!\mathcal{O}_S)} \otimes1}
			\mathbf{R}f^\prime_\ast \mathbf{L}(g^\prime)^\ast
			\operatorname{\mathbf{R}\mathcal{H}\! \mathit{om}}(P\otimes^{\mathbf{L}}E,f^!\mathcal{O}_S) \otimes^{\mathbf{L}} \mathbf{R}f^\prime_\ast \mathbf{L}(g^\prime)^\ast (P\otimes^{\mathbf{L}}E)
			\\&\to\mathbf{R}f^\prime_\ast ( \mathbf{L}(g^\prime)^\ast \operatorname{\mathbf{R}\mathcal{H}\! \mathit{om}} (P\otimes^{\mathbf{L}}E,f^!\mathcal{O}_S) \otimes^{\mathbf{L}} \mathbf{L}(g^\prime)^\ast (P\otimes^{\mathbf{L}}E))
			\\&\to\mathbf{R}f^\prime_\ast(\operatorname{\mathbf{R}\mathcal{H}\! \mathit{om}}(\mathbf{L}(g^\prime)^\ast (P\otimes^{\mathbf{L}}E), (f^\prime)^!\mathcal{O}_{S^\prime}) \otimes^{\mathbf{L}} \mathbf{L}(g^\prime)^\ast (P\otimes^{\mathbf{L}}E))
			\\&\to\mathbf{R}f^\prime_\ast (f^\prime)^!\mathcal{O}_{S^\prime} \xrightarrow{\operatorname{Tr}_{f^\prime}} \mathcal{O}_{S^\prime}.
		\end{aligned}
	\end{displaymath}
	Here the second arrow is the relative cup product, the third is
	induced by $\theta_{P\otimes^{\mathbf{L}}E}$, and the fourth is
	evaluation.

	These two morphisms agree.
	Indeed, the relative cup product is compatible with base change by \cite[\href{https://stacks.math.columbia.edu/tag/0H9A}{Tag 0H9A}]{StacksProject}.
	The first arrow in $\theta_{P\otimes^{\mathbf{L}}E}$ is constructed by pulling back the evaluation morphism, and the required compatibility is the
	construction of \cite[\href{https://stacks.math.columbia.edu/tag/08JG}{Tag 08JG}]{StacksProject}.
	Finally, \cite[\href{https://stacks.math.columbia.edu/tag/0E5L}{Tag 0E5L}]{StacksProject} gives the commutative diagram
	\begin{displaymath}
		\begin{tikzcd}
			{\mathbf{L}g^\ast\mathbf{R}f_\ast
			f^!\mathcal{O}_S} && {\mathcal{O}_{S^\prime}}
			\\{\mathbf{R}f^\prime_\ast \mathbf{L}(g^\prime)^\ast f^!\mathcal{O}_S} && {\mathbf{R}f^\prime_\ast (f^\prime)^!\mathcal{O}_{S^\prime}.}
			\arrow["{\mathbf{L}g^\ast\operatorname{Tr}_f}", from=1-1,to=1-3]
			\arrow["{\alpha_{f^!\mathcal{O}_S}}"', from=1-1,to=2-1]
			\arrow["{\mathbf{R}f^\prime_\ast(\beta_{f,g}(\mathcal{O}_S))}"',from=2-1,to=2-3]
			\arrow["{\operatorname{Tr}_{f^\prime}}"', from=2-3,to=1-3]
		\end{tikzcd}
	\end{displaymath}
	expressing compatibility with the trace morphisms.
	Hence, \eqref{diag:relative_dual_base_change_compatibility} commutes.
	Every morphism in \eqref{diag:relative_dual_base_change_compatibility}, except possibly $\mathbf{R}f^\prime_\ast(\theta_{P\otimes^{\mathbf{L}}E})$, is an isomorphism.
	Thus, $\mathbf{R}f^\prime_\ast (\theta_{P\otimes^{\mathbf{L}}E})$ is an isomorphism (e.g.\ use commutativity of \eqref{diag:relative_dual_base_change_compatibility}).
	This completes the proof.
\end{proof}

\begin{proposition}
	\label{prop:right_adjoint_pullback}
	Consider \Cref{setup:fm_cube_pullback} with $t$ quasi-affine.
	Let $K\in D_{\operatorname{qc}}(Y_1\times_{S} Y_2)$ be pseudocoherent and relatively perfect over $Y_2$.
	%%NOTE: So $K$ is bounded
	Then $\Phi_{\mathbf{L}(t^\prime)^\ast K^\prime}$ is right adjoint to $\Phi_{\mathbf{L} (t^\prime)^\ast K}$ on $D_{\operatorname{qc}}$ where
	\begin{displaymath}
		K^\prime := \operatorname{\mathbf{R}\mathcal{H}\! \mathit{om}}(K, (f^\prime_1)^! \mathcal{O}_{Y_2}).
	\end{displaymath}
	In particular, the kernel of the right adjoint is a bounded pseudocoherent complex.
\end{proposition}

\begin{proof}
	We check the first claim.
	By \Cref{lem:neeman187}, we have $b^\times \cong b^!$ on $D_{\operatorname{qc}}$ for any proper morphisms $b$ of Noetherian algebraic spaces.
	For any $E\in D_{\operatorname{qc}}(Y_1\times_{S} T)$ and $A\in D_{\operatorname{qc}}(Y_2\times_{S} T)$, there is a string of natural isomorphisms obtained by adjunctions:
	\begin{displaymath}
		\begin{aligned}
			\operatorname{Hom}
			&(\Phi_{\mathbf{L}(t^\prime)^\ast K}(E), A)
			\\&= \operatorname{Hom}(\mathbf{R}(g^\prime_1)_\ast(\mathbf{L}(g^\prime_2)^\ast E \otimes^{\mathbf{L}} \mathbf{L}(t^\prime)^\ast K), A)
			\\&\cong \operatorname{Hom}(\mathbf{L}(g^\prime_2)^\ast E \otimes^{\mathbf{L}} \mathbf{L}(t^\prime)^\ast K, (g^\prime_1)^! A) && \textrm{(right adjoint of push)}
			\\&\cong \operatorname{Hom}(\mathbf{L}(g^\prime_2)^\ast E,
			\operatorname{\mathbf{R}\mathcal{H}\! \mathit{om}}(\mathbf{L}(t^\prime)^\ast K, (g^\prime_1)^! A)) && \textrm{(tensor/hom)}
			\\&\cong \operatorname{Hom}(E,
			\mathbf{R}(g^\prime_2)_\ast \operatorname{\mathbf{R}\mathcal{H}\! \mathit{om}}(\mathbf{L}(t^\prime)^\ast K, (g^\prime_1)^! A)) && \textrm{(push/pull)}.
		\end{aligned}
	\end{displaymath}
	The desired claim follows if we can find an isomorphism
	\begin{displaymath}
		\begin{aligned}
			\mathbf{R}(g^\prime_2)_\ast
			& \operatorname{\mathbf{R}\mathcal{H}\! \mathit{om}}(\mathbf{L}(t^\prime)^\ast K, (g^\prime_1)^! A)
			\\&\cong \mathbf{R}(g^\prime_2)_\ast(
			\mathbf{L}(t^\prime)^\ast \operatorname{\mathbf{R}\mathcal{H}\! \mathit{om}}(K,(f^\prime_1)^! \mathcal{O}_{Y_2})
			\otimes^{\mathbf{L}} \mathbf{L}(g^\prime_1)^\ast A
			).
		\end{aligned}
	\end{displaymath}
	Moreover, by \Cref{lem:upper_shriek_coherent_cohomology}, $(f^\prime_1)^! \mathcal{O}_{Y_2} \in D^+_{\operatorname{coh}}(Y_1\times_{S} Y_2)$.
	Then \Cref{lem:base_change_relative_dualizing_complex} asserts that
	\begin{equation}
		\label{eq:right_adjoint_pullback}
		\begin{aligned}
			\mathbf{L}(t^\prime)^\ast
			\operatorname{\mathbf{R}\mathcal{H}\! \mathit{om}}(K,(f^\prime_1)^! \mathcal{O}_{Y_2})
			&\cong \operatorname{\mathbf{R}\mathcal{H}\! \mathit{om}}(\mathbf{L}(t^\prime)^\ast K, (g^\prime_1)^! \mathcal{O}_{Y_2\times_{S} T}).
			%%NOTE: There is an isomorphism $(g^\prime_1)^! \mathcal{O}_{Y_2\times_{S} T} \cong \mathbf{L}(t^\prime)^\ast (f^\prime_1)^! \mathcal{O}_{Y_2}$, and so, the latter object is in $D^+_{\operatorname{qc}}$.
		\end{aligned}
	\end{equation}
	From \Cref{lem:induced_dbcoh_perf_preservation_upon_pullback_with_t_affine}, $\Phi_{\mathbf{L}(t^\prime)^\ast K}$ restricts to an exact functor
	\begin{displaymath}
		\operatorname{Perf}(Y_1\times_{S} T)
		\to
		\operatorname{Perf}(Y_2\times_{S} T).
	\end{displaymath}
	Then \Cref{cor:preservation} implies
	\begin{displaymath}
		\mathbf{R}(g_1^\prime)_\ast
		(\mathbf{L}(t^\prime)^\ast K \otimes^{\mathbf{L}}
		\operatorname{Perf}(Y_1\times_{S} Y_2 \times_{S} T))
		\subseteq \operatorname{Perf}(Y_2\times_{S} T).
	\end{displaymath}
	There is another string of isomorphisms for all $A \in D_{\operatorname{qc}}$:
	\begin{displaymath}
		\begin{aligned}
			&\mathbf{R}(g^\prime_2)_\ast\operatorname{\mathbf{R}\mathcal{H}\! \mathit{om}}(\mathbf{L}(t^\prime)^\ast K, (g^\prime_1)^! A)
			\\&\cong \mathbf{R}(g^\prime_2)_\ast
			(
			\operatorname{\mathbf{R}\mathcal{H}\! \mathit{om}}(\mathbf{L}(t^\prime)^\ast K,(g^\prime_1)^! \mathcal{O}_{Y_2 \times_{S} T})
			\otimes^{\mathbf{L}} \mathbf{L}(g^\prime_1)^\ast A
			) && (\textrm{\Cref{prop:stacky_f_perf_gives_iso_between_uppershriek_and_pullback_up_to_tensor}})
			\\&\cong \mathbf{R}(g^\prime_2)_\ast
			(
			\mathbf{L}(t^\prime)^\ast \operatorname{\mathbf{R}\mathcal{H}\! \mathit{om}}(K,(f^\prime_1)^! \mathcal{O}_{Y_2})
			\otimes^{\mathbf{L}} \mathbf{L}(g^\prime_1)^\ast A
			) && \textrm{(use \eqref{eq:right_adjoint_pullback})}.
		\end{aligned}
	\end{displaymath}
	In particular, the second uses the fact that $\mathbf{L}(t^\prime)^\ast K$ is $g_1^\prime$-quasi-perfect via \Cref{thm:bounded_pseudocoherence_perfectness_faithfully_flat_affine,thm:stacky_all_coincide}.
	Lastly, the second claim follows from
	\Cref{lem:upper_shriek_coherent_cohomology,lem:quasi-perfect_bounded}.
\end{proof}

\begin{corollary}
	\label{cor:relatively_perfect_y2_implies_right_adjoint_restrict_to_dbcoh}
	Consider \Cref{setup:fm_cube_pullback} with $t=1_S$.
	Suppose $K\in D_{\operatorname{qc}}(Y_1\times_{S} Y_2)$ is pseudocoherent and relatively perfect over each $Y_i$.
	Denote by $K^\prime$ the kernel of the integral transform obtained in \Cref{prop:right_adjoint_pullback} which is right adjoint to $\Phi_K$ on $D_{\operatorname{qc}}$.
	Then $\Phi_{K^\prime}$ restricts to $D^b_{\operatorname{coh}}$.
	In particular, $\Phi_K$ and $\Phi_{K^\prime}$ form an adjoint pair on $D^b_{\operatorname{coh}}$.
\end{corollary}

\begin{proof}
	Since $K$ is $f^\prime_1$-quasi-perfect, \Cref{lem:Ballard_relative_perf_implies_dual_is_such} implies that $K^\prime$ is $f^\prime_1$-quasi-perfect.
	Hence, by \Cref{cor:preservation}, $\Phi_{K^\prime}(D^b_{\operatorname{coh}}(Y_2))\subseteq D^b_{\operatorname{coh}}(Y_1)$.
\end{proof}

\begin{lemma}
	\label{lem:left_adjoint}
	Consider \Cref{setup:fm_cube_pullback}. Suppose $K\in D_{\operatorname{qc}}(Y_1\times_{S} Y_2)$ is pseudocoherent and relatively perfect over $Y_1$.
	Then $\Phi_K$ admits a left adjoint on $D_{\operatorname{qc}}$.
	In particular, the left adjoint is of the form $\Phi_{K^\prime}$ where $K^\prime := \operatorname{\mathbf{R}\mathcal{H}\! \mathit{om}} (K, (f^\prime_2)^! \mathcal{O}_{Y_1})$.
\end{lemma}

\begin{proof}
	By \Cref{lem:neeman187}, $(f^\prime_i)^\times \cong (f^\prime_i)^!$ for each $i$.
	Also, from \Cref{cor:preservation}, $K$ is $f^\prime_2$-quasi-perfect.
	Thus, the desired claim follows from the string of natural isomorphisms for all $E\in D_{\operatorname{qc}}(Y_2)$ and $G\in D_{\operatorname{qc}}(Y_1)$:
	\begin{displaymath}
		\begin{aligned}
			\operatorname{Hom}& (E,\Phi_K (G))
			\cong \operatorname{Hom}( E ,\mathbf{R}(f_1^\prime)_\ast (K \otimes^{\mathbf{L}} \mathbf{L}(f^\prime_2)^\ast G)) && \textrm{(definition)}
			\\&\cong \operatorname{Hom}( \mathbf{L}(f_1^\prime)^\ast E , K \otimes^{\mathbf{L}} \mathbf{L}(f_2^\prime)^\ast G) && \textrm{(adjunction)}
			\\&\cong \operatorname{Hom}( \mathbf{R}(f_2^\prime)_\ast ( \operatorname{\mathbf{R}\mathcal{H}\! \mathit{om}} (K, (f^\prime_2)^! \mathcal{O}_{Y_1}) \otimes^{\mathbf{L}} \mathbf{L}(f_1^\prime)^\ast E ) , G) && (\textrm{\Cref{thm:adjoint_dqc_for_f-quasi-perfect}}).
		\end{aligned}
	\end{displaymath}
\end{proof}

\begin{proposition}
	\label{prop:left_adjoint_pullback}
	Consider \Cref{setup:fm_cube_pullback} with $t$ quasi-affine.
	Let $K\in D_{\operatorname{qc}}(Y_1\times_{S} Y_2)$ be pseudocoherent and relatively perfect over $Y_1$.
	Then $\Phi_{\mathbf{L}(t^\prime)^\ast K^\prime}$ is left adjoint to $\Phi_{\mathbf{L} (t^\prime)^\ast K}$ on $D_{\operatorname{qc}}$ where $K^\prime:= \operatorname{\mathbf{R}\mathcal{H}\! \mathit{om}} (  K,  (f^\prime_2)^! \mathcal{O}_{Y_1})$.
	In such a case, $K^\prime\in D^b_{\operatorname{coh}}(Y_1\times_{S} Y_2)$.
\end{proposition}

\begin{proof}
	The second claim follows from \Cref{lem:upper_shriek_coherent_cohomology,lem:quasi-perfect_bounded}.
	We check the first claim.
	By \Cref{lem:neeman187}, we have $b^\times \cong b^!$ on $D_{\operatorname{qc}}$ for any proper morphism $b$ of Noetherian algebraic spaces.
	By \Cref{thm:bounded_pseudocoherence_perfectness_faithfully_flat_affine}, $\mathbf{L} (t^\prime)^\ast K$ is relatively perfect over $Y_1 \times_{S} T$.
	Moreover, from \Cref{lem:left_adjoint}, the kernel of the integral transform which is left adjoint to $\Phi_{\mathbf{L} (t^\prime)^\ast K}$ is given by the object $\operatorname{\mathbf{R}\mathcal{H}\! \mathit{om}} ( \mathbf{L} (t^\prime)^\ast K, (g^\prime_2)^! \mathcal{O}_{Y_1\times_{S} T})$.
	It suffices to show there is an isomorphism
	\begin{displaymath}
		\begin{aligned}
			\mathbf{R}(g_2^\prime)_\ast & ( \operatorname{\mathbf{R}\mathcal{H}\! \mathit{om}} ( \mathbf{L} (t^\prime)^\ast K, (g^\prime_2)^! \mathcal{O}_{Y_1\times_{S} T}) \otimes^{\mathbf{L}} \mathbf{L}(g_1^\prime)^\ast E )
			\\&\to  \mathbf{R}(g_2^\prime)_\ast ( \mathbf{L} (t^\prime)^\ast  \operatorname{\mathbf{R}\mathcal{H}\! \mathit{om}} ( K, (f^\prime_2)^! \mathcal{O}_{Y_1}) \otimes^{\mathbf{L}} \mathbf{L}(g_1^\prime)^\ast E ).
		\end{aligned}
	\end{displaymath}
	This follows if we can find an isomorphism
	\begin{displaymath}
		\operatorname{\mathbf{R}\mathcal{H}\! \mathit{om}} ( \mathbf{L} (t^\prime)^\ast K, (g^\prime_2)^! \mathcal{O}_{Y_1\times_{S} T})
		\to  \mathbf{L} (t^\prime)^\ast ( \operatorname{\mathbf{R}\mathcal{H}\! \mathit{om}} ( K, (f^\prime_2)^! \mathcal{O}_{Y_1})).
	\end{displaymath}
	By \Cref{lem:upper_shriek_coherent_cohomology}, $(g^\prime_2)^! \mathcal{O}_{Y_1\times_{S} T} \in D^+_{\operatorname{coh}}(Y_1\times_{S} Y_2\times_{S} T)$.
	Then \Cref{lem:base_change_relative_dualizing_complex} asserts that
	\begin{displaymath}
		\begin{aligned}
			\mathbf{L}(t^\prime)^\ast \operatorname{\mathbf{R}\mathcal{H}\! \mathit{om}}( K, (f^\prime_2)^! \mathcal{O}_{Y_1})
			&\cong \operatorname{\mathbf{R}\mathcal{H}\! \mathit{om}}(\mathbf{L}(t^\prime)^\ast K, (g^\prime_2)^! \mathcal{O}_{Y_1\times_{S} T}).
			%%NOTE: There is an isomorphism $(g^\prime_2)^! \mathcal{O}_{Y_1\times_{S} T} \cong \mathbf{L}(t^\prime)^\ast (f^\prime_2)^! \mathcal{O}_{Y_1}$, and so, the latter object is in $D^+_{\operatorname{qc}}$.
		\end{aligned}
	\end{displaymath}
	This completes the proof.
\end{proof}

%%%%%%%%%%%%%%%%%%%%%%%%%%%%%%%%%%%%%
\section{Full faithfulness \& equivalences}
\label{sec:fully_faithful_equivalence}
%%%%%%%%%%%%%%%%%%%%%%%%%%%%%%%%%%%%%

\begin{lemma}
	\label{lem:fm_cube_pullback_natural_isomorphism}
	Consider \Cref{setup:fm_cube_pullback}. Let $K\in D_{\operatorname{qc}}(Y_1\times_{S} Y_2)$.
	Then on $D_{\operatorname{qc}}$ there is a natural isomorphism
	\begin{displaymath}
		\beta^K\colon \Phi_K \circ \mathbf{R}(t_1)_\ast
		\to \mathbf{R}(t_2)_\ast \circ \Phi_{\mathbf{L}(t^\prime)^\ast K}.
	\end{displaymath}
\end{lemma}

\begin{proof}
	This follows from the string of natural isomorphisms for each $E\in D_{\operatorname{qc}}(Y_1\times_{S} T)$:
	\begin{displaymath}
		\begin{aligned}
			\Phi_K \circ \mathbf{R} (t_1)_\ast (E)
			&= \mathbf{R} (f_1^\prime)_\ast
			( \mathbf{L} (f^\prime_2)^\ast \mathbf{R} (t_1)_\ast E \otimes^{\mathbf{L}} K )
			\\
			&\cong
			\mathbf{R} (f_1^\prime)_\ast
			(
			\mathbf{R} t^\prime_\ast \mathbf{L}(g^\prime_2)^\ast E
			\otimes^{\mathbf{L}} K
			)
			&
			\textrm{(flat base change)}
			\\
			&\cong
			\mathbf{R} (f_1^\prime)_\ast
			\mathbf{R} t^\prime_\ast
			(
			\mathbf{L}(g^\prime_2)^\ast E
			\otimes^{\mathbf{L}} \mathbf{L}(t^\prime)^\ast K
			)
			&
			\textrm{(projection formula)}
			\\
			&\cong
			\mathbf{R} (t_2)_\ast \mathbf{R} (g^\prime_1)_\ast
			(
			\mathbf{L}(g^\prime_2)^\ast E
			\otimes^{\mathbf{L}} \mathbf{L}(t^\prime)^\ast K
			)
			&
			\textrm{(pseudofunctoriality)}
			\\&=
			\mathbf{R} (t_2)_\ast
			\circ
			\Phi_{\mathbf{L}(t^\prime)^\ast K} (E).
		\end{aligned}
	\end{displaymath}
	We conclude the proof.
\end{proof}

The next two results, \Cref{lem:integral_transform_is_morphism_of_adjoints} and \Cref{prop:pullbacks_is_a_morphism_of_adjoint_integral_transforms}, deal with mates.
We redirect the reader to \Cref{sec:mates} for definitions.

\begin{lemma}
	\label{lem:integral_transform_is_morphism_of_adjoints}
	Consider \Cref{setup:fm_cube_pullback}.
	Let $K\in D_{\operatorname{qc}}(Y_1\times_SY_2)$.
	The isomorphisms $\alpha^K$ and $\beta^K$ of \Cref{lem:fm_cube_pullback_natural_isomorphism,lem:fm_cube_flat_base_change_natural_isomorphism} are respectively left and right mates under the diagram:
	\begin{displaymath}
		% https://q.uiver.app/#q=WzAsNCxbMCwwLCJEX1xccWMoWV8xKSJdLFswLDEsIkRfXFxxYyhZXzFcXHRpbWVzX1NUKSJdLFsxLDAsIkRfXFxxYyhZXzIpIl0sWzEsMSwiRF9cXHFjKFlfMlxcdGltZXNfU1QpIl0sWzAsMiwiXFxQaGlfSyJdLFsxLDMsIlxcUGhpX3tcXG1hdGhiZntMfSh0JyleKkt9IiwyXSxbMCwxLCJcXG1hdGhiZntMfXRfMV4qIiwyLHsib2Zmc2V0IjoyfV0sWzEsMCwiXFxtYXRoYmZ7Un0odF8xKV8qIiwyLHsib2Zmc2V0IjoyfV0sWzIsMywiXFxtYXRoYmZ7TH10XzJeKiIsMix7Im9mZnNldCI6Mn1dLFszLDIsIlxcbWF0aGJme1J9KHRfMilfKiIsMix7Im9mZnNldCI6Mn1dLFs2LDcsIiIsMix7ImxldmVsIjoxLCJzdHlsZSI6eyJuYW1lIjoiYWRqdW5jdGlvbiJ9fV0sWzgsOSwiIiwyLHsibGV2ZWwiOjEsInN0eWxlIjp7Im5hbWUiOiJhZGp1bmN0aW9uIn19XV0=
		\begin{tikzcd}
			{D_{\operatorname{qc}}(Y_1)} & {D_{\operatorname{qc}}(Y_2)} \\
			{D_{\operatorname{qc}}(Y_1\times_ST)} & {D_{\operatorname{qc}}(Y_2\times_ST).}
			\arrow["{\Phi_K}", from=1-1, to=1-2]
			\arrow[""{name=0, anchor=center, inner sep=0}, "{\mathbf{L}t_1^\ast }"', shift right=2, from=1-1, to=2-1]
			\arrow[""{name=1, anchor=center, inner sep=0}, "{\mathbf{L}t_2^\ast }"', shift right=2, from=1-2, to=2-2]
			\arrow[""{name=2, anchor=center, inner sep=0}, "{\mathbf{R}(t_1)_\ast }"', shift right=2, from=2-1, to=1-1]
			\arrow["{\Phi_{\mathbf{L}(t^\prime)^\ast K}}"', from=2-1, to=2-2]
			\arrow[""{name=3, anchor=center, inner sep=0}, "{\mathbf{R}(t_2)_\ast }"', shift right=2, from=2-2, to=1-2]
			\arrow["\dashv"{anchor=center}, draw=none, from=0, to=2]
			\arrow["\dashv"{anchor=center}, draw=none, from=1, to=3]
		\end{tikzcd}
	\end{displaymath}
\end{lemma}

\begin{remark}
	The compatibility of \Cref{lem:integral_transform_is_morphism_of_adjoints} also appears in \cite[Lemma 5.6(3)]{Hall/Priver:2024}.
	We provide a detailed proof here.
\end{remark}

\begin{proof}
	[Proof of \Cref{lem:integral_transform_is_morphism_of_adjoints}]
	Denote by $\xi^i\colon \mathbf{L}t_i^\ast \mathbf{R}(t_i)_\ast \to1$ the counit of the adjunction $\mathbf{L}t_i^\ast \dashv\mathbf{R}(t_i)_\ast $ with $i=1,2$. To prove the desired claim, it suffices to verify that Condition \ref{lem:morph_adj2} of \Cref{lem:morph_adj} holds; that is, that the following diagram commutes:
	\begin{displaymath}
		% https://q.uiver.app/#q=WzAsNCxbMCwwLCJcXG1hdGhiZntMfXRfMl5cXGFzdFxcUGhpX0tcXG1hdGhiZntSfSh0XzEpX1xcYXN0Il0sWzEsMCwiXFxQaGlfe1xcbWF0aGJme0x9KHReXFxwcmltZSleXFxhc3QgS31cXG1hdGhiZntMfXRfMV5cXGFzdCBcXG1hdGhiZntSfSh0XzEpX1xcYXN0Il0sWzAsMSwiXFxtYXRoYmZ7TH10XzJeXFxhc3RcXG1hdGhiZntSfSh0XzIpX1xcYXN0XFxQaGlfe1xcbWF0aGJme0x9KHReXFxwcmltZSleXFxhc3QgS30iXSxbMSwxLCJcXFBoaV97XFxtYXRoYmZ7TH0odF5cXHByaW1lKV5cXGFzdCBLfS4iXSxbMCwyLCJcXG1hdGhiZntMfXRfMl5cXGFzdChcXGJldGFeSykiLDJdLFswLDEsIlxcYWxwaGFeS197XFxtYXRoYmZ7Un0odF8xKV9cXGFzdH0iXSxbMSwzLCJcXFBoaV97XFxtYXRoYmZ7TH0odF5cXHByaW1lKV5cXGFzdCBLfShcXHhpXjEpIl0sWzIsMywiXFx4aV4yX3tcXFBoaV97XFxtYXRoYmZ7TH0odF5cXHByaW1lKV5cXGFzdCBLfX0iLDJdXQ==
		\begin{tikzcd}
			{\mathbf{L}t_2^\ast\Phi_K\mathbf{R}(t_1)_\ast} & {\Phi_{\mathbf{L}(t^\prime)^\ast K}\mathbf{L}t_1^\ast \mathbf{R}(t_1)_\ast} \\
			{\mathbf{L}t_2^\ast\mathbf{R}(t_2)_\ast\Phi_{\mathbf{L}(t^\prime)^\ast K}} & {\Phi_{\mathbf{L}(t^\prime)^\ast K}.}
			\arrow["{\alpha^K_{\mathbf{R}(t_1)_\ast}}", from=1-1, to=1-2]
			\arrow["{\mathbf{L}t_2^\ast(\beta^K)}"', from=1-1, to=2-1]
			\arrow["{\Phi_{\mathbf{L}(t^\prime)^\ast K}(\xi^1)}", from=1-2, to=2-2]
			\arrow["{\xi^2_{\Phi_{\mathbf{L}(t^\prime)^\ast K}}}"', from=2-1, to=2-2]
		\end{tikzcd}
	\end{displaymath}
	For brevity, in the remainder of the proof, we drop the $\mathbf{L}$'s and $\mathbf{R}$'s in the notation for the derived functors (i.e.\
	all functors now are understood to be derived).	Hence, the diagram above is now
	\begin{equation}
		\label{eq:square_without_L_and_R}
		% https://q.uiver.app/#q=WzAsNCxbMCwxLCJ0XzJeXFxhc3QgXFxQaGlfSyAodF8xKV9cXGFzdCJdLFsxLDIsIlxcUGhpX3sodF5cXHByaW1lKV5cXGFzdCBLfXRfMV5cXGFzdCAodF8xKV5cXGFzdCJdLFsxLDAsInRfMl5cXGFzdCAodF8yKV9cXGFzdCBcXFBoaV97KHReXFxwcmltZSleXFxhc3QgS30iXSxbMiwxLCJcXFBoaV97KHReXFxwcmltZSleXFxhc3QgS30uIl0sWzAsMiwidF8yXlxcYXN0KFxcYmV0YV5LKSJdLFsyLDMsIlxceGleMl97XFxQaGlfeyh0XlxccHJpbWUpXlxcYXN0IEt9fSJdLFsxLDMsIlxcUGhpX3sodF5cXHByaW1lKV5cXGFzdCBLfShcXHhpXjEpIiwyXSxbMCwxLCJcXGFscGhhXktfeyh0XzEpX1xcYXN0fSIsMl1d
		\begin{tikzcd}
			& {t_2^\ast (t_2)_\ast \Phi_{(t^\prime)^\ast K}} & \\
			{t_2^\ast \Phi_K (t_1)_\ast} && {\Phi_{(t^\prime)^\ast K}.} \\
			& {\Phi_{(t^\prime)^\ast K}t_1^\ast (t_1)^\ast}
			\arrow["{\xi^2_{\Phi_{(t^\prime)^\ast K}}}", from=1-2, to=2-3]
			\arrow["{t_2^\ast(\beta^K)}", from=2-1, to=1-2]
			\arrow["{\alpha^K_{(t_1)_\ast}}"', from=2-1, to=3-2]
			\arrow["{\Phi_{(t^\prime)^\ast K}(\xi^1)}"', from=3-2, to=2-3]
		\end{tikzcd}
	\end{equation}
	Choose $E\in D_{\operatorname{qc}}(Y_1\times_ST)$.
	If we expand the definition of $\beta^K$ and $\alpha^K$ in \eqref{eq:square_without_L_and_R}
	(see the proofs of
	\Cref{lem:fm_cube_pullback_natural_isomorphism,lem:fm_cube_flat_base_change_natural_isomorphism}),
	we obtain a large diagram consisting of various faces.
	In what follows, we make this a bit more explicit.
	In particular, we spell out each face needed to understand \eqref{eq:square_without_L_and_R}.
	To describe said diagram more explicitly, we use the following abbreviations for the (natural) isomorphisms:
	\begin{itemize}
		\item $\operatorname{BC}_i$, $i\in\{1,2\}$ for the tor-independent base change isomorphisms
		\item $\operatorname{PF}$ for the projection formula
		\item $F_{(-)_\ast }$ for functoriality of $(-)_\ast $
		\item $F_{(-)^\ast }$ for functoriality of $(-)^\ast $
		\item $M_{(-)^\ast }$ for monoidality of $(-)^\ast $
		\item $\xi^\prime $ for the counit of $(t^\prime)^\ast \dashv (t^\prime)_\ast $.
	\end{itemize}
	Now, consider the following diagrams:
	\begin{equation}
		\label{eq:big_boyA}
		% https://q.uiver.app/#q=WzAsNCxbMCwwLCJ0XzJeXFxhc3QgKGZfMSlfXFxhc3ReXFxwcmltZSAoKGZfMilfXFxhc3ReXFxwcmltZSAgKHRfMSlfXFxhc3QgRVxcb3RpbWVzIEspIl0sWzAsMiwiKGdfMSlfXFxhc3ReXFxwcmltZSAodF5cXHByaW1lKV5cXGFzdCAoKGZfMl5cXHByaW1lKV5cXGFzdCAgICh0XzEpX1xcYXN0IEVcXG90aW1lcyBLKSJdLFsyLDAsInRfMl5cXGFzdCAoZl8xKV9cXGFzdF5cXHByaW1lICgodF5cXHByaW1lKV9cXGFzdCAoZ18yXlxccHJpbWUpXlxcYXN0IEVcXG90aW1lcyBLKSJdLFsyLDIsIihnXzEpX1xcYXN0XlxccHJpbWUgKHReXFxwcmltZSleXFxhc3QgKHRfXFxhc3QgXlxccHJpbWUgKGdfMl5cXHByaW1lKV5cXGFzdCBFXFxvdGltZXMgSykiXSxbMCwxLCJcXG9wZXJhdG9ybmFtZXtCQ31fMiIsMV0sWzAsMiwiXFxvcGVyYXRvcm5hbWV7QkN9XzEiLDFdLFsyLDMsIlxcb3BlcmF0b3JuYW1le0JDfV8yIiwxXSxbMSwzLCJcXG9wZXJhdG9ybmFtZXtCQ31fMSIsMV1d
		\begin{tikzcd}
			{t_2^\ast (f_1^\prime)_\ast ((f_2^\prime)_\ast  (t_1)_\ast E\otimes K)} && {t_2^\ast (f_1^\prime)_\ast ((t^\prime)_\ast (g_2^\prime)^\ast E\otimes K)} \\
			\\
			{(g_1^\prime)_\ast (t^\prime)^\ast ((f_2^\prime)^\ast   (t_1)_\ast E\otimes K)} && {(g_1^\prime)_\ast (t^\prime)^\ast (t_\ast ^\prime (g_2^\prime)^\ast E\otimes K)}
			\arrow["{\operatorname{BC}_1}"{description}, from=1-1, to=1-3]
			\arrow["{\operatorname{BC}_2}"{description}, from=1-1, to=3-1]
			\arrow["{\operatorname{BC}_2}"{description}, from=1-3, to=3-3]
			\arrow["{\operatorname{BC}_1}"{description}, from=3-1, to=3-3]
		\end{tikzcd}
	\end{equation}

	\begin{equation}
		\label{eq:big_boyB}
		% https://q.uiver.app/#q=WzAsNCxbMCwwLCJ0XzJeXFxhc3QgKGZfMSlfXFxhc3ReXFxwcmltZSAoKHReXFxwcmltZSlfXFxhc3QgKGdfMl5cXHByaW1lKV5cXGFzdCBFXFxvdGltZXMgSykiXSxbMCwyLCIoZ18xKV9cXGFzdF5cXHByaW1lICh0XlxccHJpbWUpXlxcYXN0ICh0X1xcYXN0IF5cXHByaW1lIChnXzJeXFxwcmltZSleXFxhc3QgRVxcb3RpbWVzIEspIl0sWzIsMCwidF8yXlxcYXN0IChmXzEpX1xcYXN0XlxccHJpbWUgKHReXFxwcmltZSlfXFxhc3QgKChnXzJeXFxwcmltZSleXFxhc3QgRVxcb3RpbWVzICh0XlxccHJpbWUpXlxcYXN0IEspIl0sWzIsMiwiKGdfMSlfXFxhc3ReXFxwcmltZSAodF5cXHByaW1lKV5cXGFzdCB0X1xcYXN0IF5cXHByaW1lICgoZ18yXlxccHJpbWUpXlxcYXN0IEVcXG90aW1lcyAodF5cXHByaW1lKV5cXGFzdCBLKSJdLFswLDEsIlxcb3BlcmF0b3JuYW1le0JDfV8yIiwxXSxbMCwyLCJcXG9wZXJhdG9ybmFtZXtQRn0iLDFdLFsxLDMsIlxcb3BlcmF0b3JuYW1le1BGfSIsMV0sWzIsMywiXFxvcGVyYXRvcm5hbWV7QkN9XzIiLDFdXQ==
		\begin{tikzcd}
			{t_2^\ast (f_1^\prime)_\ast ((t^\prime)_\ast (g_2^\prime)^\ast E\otimes K)} && {t_2^\ast (f_1^\prime)_\ast (t^\prime)_\ast ((g_2^\prime)^\ast E\otimes (t^\prime)^\ast K)} \\
			\\
			{(g_1^\prime)_\ast (t^\prime)^\ast (t_\ast ^\prime (g_2^\prime)^\ast E\otimes K)} && {(g_1^\prime)_\ast (t^\prime)^\ast t_\ast ^\prime ((g_2^\prime)^\ast E\otimes (t^\prime)^\ast K)}
			\arrow["{\operatorname{PF}}"{description}, from=1-1, to=1-3]
			\arrow["{\operatorname{BC}_2}"{description}, from=1-1, to=3-1]
			\arrow["{\operatorname{BC}_2}"{description}, from=1-3, to=3-3]
			\arrow["{\operatorname{PF}}"{description}, from=3-1, to=3-3]
		\end{tikzcd}
	\end{equation}

	\begin{equation}
		\label{eq:big_boyC}
		% https://q.uiver.app/#q=WzAsNCxbMCwwLCIoZ18xKV9cXGFzdF5cXHByaW1lICh0XlxccHJpbWUpXlxcYXN0ICgoZl8yXlxccHJpbWUpXlxcYXN0ICAgKHRfMSlfXFxhc3QgRVxcb3RpbWVzIEspIl0sWzIsMCwiKGdfMSlfXFxhc3ReXFxwcmltZSAodF5cXHByaW1lKV5cXGFzdCAodF9cXGFzdCBeXFxwcmltZSAoZ18yXlxccHJpbWUpXlxcYXN0IEVcXG90aW1lcyBLKSJdLFswLDIsIihnXzEpX1xcYXN0XlxccHJpbWUgKCh0XlxccHJpbWUpXlxcYXN0IChmXzJeXFxwcmltZSleXFxhc3QgICAodF8xKV9cXGFzdCBFXFxvdGltZXMgKHReXFxwcmltZSleXFxhc3QgSykiXSxbMiwyLCIoZ18xKV9cXGFzdF5cXHByaW1lICgodF5cXHByaW1lKV5cXGFzdCB0X1xcYXN0IF5cXHByaW1lIChnXzJeXFxwcmltZSleXFxhc3QgRVxcb3RpbWVzICh0XlxccHJpbWUpXlxcYXN0IEspIl0sWzAsMSwiXFxvcGVyYXRvcm5hbWV7QkN9XzEiLDFdLFswLDIsIk1feygtKV5cXGFzdH0iLDFdLFsyLDMsIlxcb3BlcmF0b3JuYW1le0JDfV8xIiwxXSxbMSwzLCJNX3soLSl9XlxcYXN0IiwxXSxbMSwyLCIiLDEseyJzdHlsZSI6eyJib2R5Ijp7Im5hbWUiOiJub25lIn0sImhlYWQiOnsibmFtZSI6Im5vbmUifX19XV0=
		\begin{tikzcd}
			{(g_1^\prime)_\ast (t^\prime)^\ast ((f_2^\prime)^\ast   (t_1)_\ast E\otimes K)} && {(g_1^\prime)_\ast (t^\prime)^\ast (t_\ast ^\prime (g_2^\prime)^\ast E\otimes K)} \\
			\\
			{(g_1^\prime)_\ast ((t^\prime)^\ast (f_2^\prime)^\ast   (t_1)_\ast E\otimes (t^\prime)^\ast K)} && {(g_1^\prime)_\ast ((t^\prime)^\ast t_\ast ^\prime (g_2^\prime)^\ast E\otimes (t^\prime)^\ast K)}
			\arrow["{\operatorname{BC}_1}"{description}, from=1-1, to=1-3]
			\arrow["{M_{(-)^\ast}}"{description}, from=1-1, to=3-1]
			\arrow[draw=none, from=1-3, to=3-1]
			\arrow["{M_{(-)}^\ast}"{description}, from=1-3, to=3-3]
			\arrow["{\operatorname{BC}_1}"{description}, from=3-1, to=3-3]
		\end{tikzcd}
	\end{equation}

	\begin{equation}
		\label{eq:big_boyD}
		% https://q.uiver.app/#q=WzAsNCxbMiwwLCIoZ18xKV9cXGFzdF5cXHByaW1lICh0XlxccHJpbWUpXlxcYXN0ICh0X1xcYXN0IF5cXHByaW1lIChnXzJeXFxwcmltZSleXFxhc3QgRVxcb3RpbWVzIEspIl0sWzAsMCwiKGdfMSlfXFxhc3ReXFxwcmltZSAoKHReXFxwcmltZSleXFxhc3QgdF9cXGFzdCBeXFxwcmltZSAoZ18yXlxccHJpbWUpXlxcYXN0IEVcXG90aW1lcyAodF5cXHByaW1lKV5cXGFzdCBLKSJdLFsyLDIsIihnXzEpX1xcYXN0XlxccHJpbWUgKHReXFxwcmltZSleXFxhc3QgdF9cXGFzdCBeXFxwcmltZSAoKGdfMl5cXHByaW1lKV5cXGFzdCBFXFxvdGltZXMgKHReXFxwcmltZSleXFxhc3QgSykiXSxbMCwyLCIoZ18xKV9cXGFzdF5cXHByaW1lICgoZ18yXlxccHJpbWUpXlxcYXN0IEVcXG90aW1lcyAodF5cXHByaW1lKV5cXGFzdCBLKSJdLFswLDEsIk1feygtKX1eXFxhc3QiLDFdLFswLDIsIlxcb3BlcmF0b3JuYW1le1BGfSIsMV0sWzIsMywiXFx4aV5cXHByaW1lIiwxXSxbMSwzLCJcXHhpXlxccHJpbWUiLDFdXQ==
		\begin{tikzcd}
			{(g_1^\prime)_\ast ((t^\prime)^\ast t_\ast ^\prime (g_2^\prime)^\ast E\otimes (t^\prime)^\ast K)} && {(g_1^\prime)_\ast (t^\prime)^\ast (t_\ast ^\prime (g_2^\prime)^\ast E\otimes K)} \\
			\\
			{(g_1^\prime)_\ast ((g_2^\prime)^\ast E\otimes (t^\prime)^\ast K)} && {(g_1^\prime)_\ast (t^\prime)^\ast t_\ast ^\prime ((g_2^\prime)^\ast E\otimes (t^\prime)^\ast K)}
			\arrow["{\xi^\prime}"{description}, from=1-1, to=3-1]
			\arrow["{M_{(-)}^\ast}"{description}, from=1-3, to=1-1]
			\arrow["{\operatorname{PF}}"{description}, from=1-3, to=3-3]
			\arrow["{\xi^\prime}"{description}, from=3-3, to=3-1]
		\end{tikzcd}
	\end{equation}

	\begin{equation}
		\label{eq:big_boyE}
		% https://q.uiver.app/#q=WzAsNCxbMCwwLCJ0XzJeXFxhc3QgKGZfMSlfXFxhc3ReXFxwcmltZSAodF5cXHByaW1lKV9cXGFzdCAoKGdfMl5cXHByaW1lKV5cXGFzdCBFXFxvdGltZXMgKHReXFxwcmltZSleXFxhc3QgSykiXSxbMCwyLCIoZ18xKV9cXGFzdF5cXHByaW1lICh0XlxccHJpbWUpXlxcYXN0IHRfXFxhc3QgXlxccHJpbWUgKChnXzJeXFxwcmltZSleXFxhc3QgRVxcb3RpbWVzICh0XlxccHJpbWUpXlxcYXN0IEspIl0sWzIsMiwiKGdfMSlfXFxhc3ReXFxwcmltZSAoKGdfMl5cXHByaW1lKV5cXGFzdCBFXFxvdGltZXMgKHReXFxwcmltZSleXFxhc3QgSykiXSxbMiwwLCJ0XzJeXFxhc3QgKHRfMilfXFxhc3QgKGdfMSlfXFxhc3ReXFxwcmltZSAoKGdfMl5cXHByaW1lKV5cXGFzdCBFXFxvdGltZXMgKHReXFxwcmltZSleXFxhc3QgSykiXSxbMCwxLCJcXG9wZXJhdG9ybmFtZXtCQ31fMiIsMV0sWzEsMiwiXFx4aV5cXHByaW1lIiwxXSxbMCwzLCJGX3soLSleXFxhc3R9IiwxXSxbMywyLCJcXHhpXjIiLDFdXQ==
		\begin{tikzcd}
			{t_2^\ast (f_1^\prime)_\ast (t^\prime)_\ast ((g_2^\prime)^\ast E\otimes (t^\prime)^\ast K)} && {t_2^\ast (t_2)_\ast (g_1^\prime)_\ast ((g_2^\prime)^\ast E\otimes (t^\prime)^\ast K)} \\
			\\
			{(g_1^\prime)_\ast (t^\prime)^\ast t_\ast ^\prime ((g_2^\prime)^\ast E\otimes (t^\prime)^\ast K)} && {(g_1^\prime)_\ast ((g_2^\prime)^\ast E\otimes (t^\prime)^\ast K)}
			\arrow["{F_{(-)^\ast}}"{description}, from=1-1, to=1-3]
			\arrow["{\operatorname{BC}_2}"{description}, from=1-1, to=3-1]
			\arrow["{\xi^2}"{description}, from=1-3, to=3-3]
			\arrow["{\xi^\prime}"{description}, from=3-1, to=3-3]
		\end{tikzcd}
	\end{equation}

	\begin{equation}
		\label{eq:big_boyF}
		% https://q.uiver.app/#q=WzAsNCxbMCwwLCIoZ18xKV9cXGFzdF5cXHByaW1lICgodF5cXHByaW1lKV5cXGFzdCAoZl8yXlxccHJpbWUpXlxcYXN0ICAgKHRfMSlfXFxhc3QgRVxcb3RpbWVzICh0XlxccHJpbWUpXlxcYXN0IEspIl0sWzIsMCwiKGdfMSlfXFxhc3ReXFxwcmltZSAoKHReXFxwcmltZSleXFxhc3QgdF9cXGFzdCBeXFxwcmltZSAoZ18yXlxccHJpbWUpXlxcYXN0IEVcXG90aW1lcyAodF5cXHByaW1lKV5cXGFzdCBLKSJdLFsyLDIsIihnXzEpX1xcYXN0XlxccHJpbWUgKChnXzJeXFxwcmltZSleXFxhc3QgRVxcb3RpbWVzICh0XlxccHJpbWUpXlxcYXN0IEspIl0sWzAsMiwiKGdeXFxwcmltZV8xKV9cXGFzdCAoKGdfMl5cXHByaW1lKV5cXGFzdCB0XzFeXFxhc3QgICh0XzEpX1xcYXN0IEVcXG90aW1lcyAodF5cXHByaW1lKV5cXGFzdCBLKSJdLFswLDEsIlxcb3BlcmF0b3JuYW1le0JDfV8xIiwxXSxbMSwyLCJcXHhpXlxccHJpbWUiLDFdLFswLDMsIkZfeygtKV5cXGFzdH0iLDFdLFszLDIsIlxceGleMSIsMV1d
		\begin{tikzcd}
			{(g_1^\prime)_\ast ((t^\prime)^\ast (f_2^\prime)^\ast   (t_1)_\ast E\otimes (t^\prime)^\ast K)} && {(g_1^\prime)_\ast ((t^\prime)^\ast t_\ast ^\prime (g_2^\prime)^\ast E\otimes (t^\prime)^\ast K)} \\
			\\
			{(g^\prime_1)_\ast ((g_2^\prime)^\ast t_1^\ast  (t_1)_\ast E\otimes (t^\prime)^\ast K)} && {(g_1^\prime)_\ast ((g_2^\prime)^\ast E\otimes (t^\prime)^\ast K)}
			\arrow["{\operatorname{BC}_1}"{description}, from=1-1, to=1-3]
			\arrow["{F_{(-)^\ast}}"{description}, from=1-1, to=3-1]
			\arrow["{\xi^\prime}"{description}, from=1-3, to=3-3]
			\arrow["{\xi^1}"{description}, from=3-1, to=3-3]
		\end{tikzcd}
	\end{equation}
	Assume that we have shown \eqref{eq:big_boyA}, \eqref{eq:big_boyB}, \eqref{eq:big_boyC}, \eqref{eq:big_boyD}, \eqref{eq:big_boyE}, \eqref{eq:big_boyF} are commutative. Then \eqref{eq:square_without_L_and_R} (in the case of $E$) is the pasting of these diagrams, i.e.\ gives the desired morphism
	\begin{displaymath}
		% https://q.uiver.app/#q=WzAsMixbMCwwLCJ0XzJeXFxhc3QgKGZfMSlfXFxhc3ReXFxwcmltZSAoKGZfMilfXFxhc3ReXFxwcmltZSAgKHRfMSlfXFxhc3QgRVxcb3RpbWVzIEspIl0sWzIsMCwiKGdfMSlfXFxhc3ReXFxwcmltZSAoKGdfMl5cXHByaW1lKV5cXGFzdCBFXFxvdGltZXMgKHReXFxwcmltZSleXFxhc3QgSykuIl0sWzAsMV1d
		\begin{tikzcd}
			{t_2^\ast (f_1^\prime)_\ast ((f_2^\prime)_\ast  (t_1)_\ast E\otimes K)} && {(g_1^\prime)_\ast ((g_2^\prime)^\ast E\otimes (t^\prime)^\ast K).}
			\arrow[from=1-1, to=1-3]
		\end{tikzcd}
	\end{displaymath}

	We explain why \eqref{eq:big_boyA}, \eqref{eq:big_boyB}, \eqref{eq:big_boyC}, \eqref{eq:big_boyD}, \eqref{eq:big_boyE}, \eqref{eq:big_boyF} are commutative. As for \eqref{eq:big_boyA}, one may use naturality of $\operatorname{BC}_2$, whereas \eqref{eq:big_boyB} is due to the naturality of $\operatorname{BC}_2$. Also, for \eqref{eq:big_boyC}, one can use the naturality of $M_{(-)^\ast }$.

	Next, we explain \eqref{eq:big_boyD}. Observe that it is the functor $(g_1^\prime)_\ast $ applied to the following diagram evaluated at $((g_2^\prime)^\ast E,K)$:
	\begin{displaymath}
		% https://q.uiver.app/#q=WzAsNCxbMCwwLCIodF5cXHByaW1lKV5cXGFzdCAodF9cXGFzdF5cXHByaW1lICgtKVxcb3RpbWVzIC0pIl0sWzEsMCwiKHReXFxwcmltZSleXFxhc3QgdF9cXGFzdF5cXHByaW1lICgtXFxvdGltZXMgKHReXFxwcmltZSleXFxhc3QgKC0pKSJdLFswLDEsIih0XlxccHJpbWUpXlxcYXN0IHRfXFxhc3ReXFxwcmltZSAoLSlcXG90aW1lcyAodF5cXHByaW1lKV5cXGFzdCgtKSJdLFsxLDEsIigtKVxcb3RpbWVzICh0XlxccHJpbWUpXlxcYXN0KC0pIl0sWzAsMSwiXFx0ZXh0e1BGfSJdLFswLDIsIlxcbWF0aHJte019X3soLSleXFxhc3R9IiwyXSxbMiwzLCJcXHhpXlxccHJpbWUiLDJdLFsxLDMsIlxceGleXFxwcmltZSJdXQ==
		\begin{tikzcd}
			{(t^\prime)^\ast (t_\ast^\prime (-)\otimes -)} & {(t^\prime)^\ast t_\ast^\prime (-\otimes (t^\prime)^\ast (-))} \\
			{(t^\prime)^\ast t_\ast^\prime (-)\otimes (t^\prime)^\ast(-)} & {(-)\otimes (t^\prime)^\ast(-)}
			\arrow["{\operatorname{PF}}", from=1-1, to=1-2]
			\arrow["{M_{(-)^\ast}}"', from=1-1, to=2-1]
			\arrow["{\xi^\prime}", from=1-2, to=2-2]
			\arrow["{\xi^\prime}"', from=2-1, to=2-2]
		\end{tikzcd}
	\end{displaymath}
	which commutes by definition of the projection formula, see e.g.\ proof of \cite[Proposition 22.81]{Gortz/Wedhorn:2023}.

	Now, for \eqref{eq:big_boyE}. Note that it is the following diagram evaluated at $(g_2^\prime)^\ast E\otimes (t^\prime)^\ast K$:
	\begin{displaymath}
		% https://q.uiver.app/#q=WzAsNCxbMCwwLCJ0XzJeXFxhc3QgKGZeXFxwcmltZV8xKV9cXGFzdCB0XlxccHJpbWVfXFxhc3QiXSxbMCwxLCIoZ15cXHByaW1lXzEpX1xcYXN0ICh0XlxccHJpbWUpXlxcYXN0IHReXFxwcmltZV9cXGFzdCJdLFsxLDEsIihnXzFeXFxwcmltZSlfXFxhc3QiXSxbMSwwLCJ0XzJeXFxhc3QgKHRfMilfXFxhc3QgKGdeXFxwcmltZV8xKV9cXGFzdCJdLFswLDMsIlxcbWF0aHJte0Z9X3soLSlfXFxhc3R9Il0sWzAsMSwiXFx0ZXh0e0JDfV8yIiwyXSxbMSwyLCJcXHhpXlxccHJpbWUiLDJdLFszLDIsIlxceGleMiJdXQ==
		\begin{tikzcd}
			{t_2^\ast (f^\prime_1)_\ast t^\prime_\ast} & {t_2^\ast (t_2)_\ast (g^\prime_1)_\ast} \\
			{(g^\prime_1)_\ast (t^\prime)^\ast t^\prime_\ast} & {(g_1^\prime)_\ast}
			\arrow["{F_{(-)_\ast}}", from=1-1, to=1-2]
			\arrow["{\operatorname{BC}_2}"', from=1-1, to=2-1]
			\arrow["{\xi^2}", from=1-2, to=2-2]
			\arrow["{\xi^\prime}"', from=2-1, to=2-2]
		\end{tikzcd}
	\end{displaymath}
	which commutes by definition of the base change morphism.

	Lastly, we check \eqref{eq:big_boyF}. However, it is the functor $(g_1^\prime)_\ast (-\otimes (t^\prime)^\ast K)$ applied to the following diagram evaluated at $E$:
	\begin{displaymath}
		% https://q.uiver.app/#q=WzAsNCxbMCwwLCIodF5cXHByaW1lKV5cXGFzdCB0XlxccHJpbWVfXFxhc3QgKGdfMl5cXHByaW1lKV5cXGFzdCJdLFswLDEsIih0XlxccHJpbWUpXlxcYXN0IChmXzJeXFxwcmltZSleXFxhc3QgKHRfMSlfXFxhc3QiXSxbMSwxLCIoZ18yXlxccHJpbWUpXlxcYXN0IHRfMV5cXGFzdCAodF8xKV9cXGFzdCJdLFsxLDAsIihnXzJeXFxwcmltZSleXFxhc3QiXSxbMCwzLCJcXHhpXlxccHJpbWUiXSxbMSwwLCJcXHRleHR7QkN9XzEiXSxbMSwyLCJcXG1hdGhybXtGfV97KC0pXlxcYXN0fSIsMl0sWzIsMywiXFx4aV4xIiwyXV0=
		\begin{tikzcd}
			{(t^\prime)^\ast t^\prime_\ast (g_2^\prime)^\ast} & {(g_2^\prime)^\ast} \\
			{(t^\prime)^\ast (f_2^\prime)^\ast (t_1)_\ast} & {(g_2^\prime)^\ast t_1^\ast (t_1)_\ast}
			\arrow["{\xi^\prime}", from=1-1, to=1-2]
			\arrow["{\operatorname{BC}_1}", from=2-1, to=1-1]
			\arrow["{F_{(-)^\ast}}"', from=2-1, to=2-2]
			\arrow["{\xi^1}"', from=2-2, to=1-2]
		\end{tikzcd}
	\end{displaymath}
	which commutes by definition of the base change morphism.
\end{proof}

\begin{remark}
	Assume the hypothesis of \Cref{lem:integral_transform_is_morphism_of_adjoints}.
	Denote by $\zeta^i$ the unit of $t_i^\ast \dashv (t_i)_\ast$ if $i\in\{1,2\}$.
	From \Cref{lem:morph_adj,lem:integral_transform_is_morphism_of_adjoints}, we have that the following diagram is commutative:
	\begin{equation}
		\label{eq:integral_transform_is_compatible_with_units_of_upper_and_lower_stars}
		% https://q.uiver.app/#q=WzAsNCxbMCwwLCJcXFBoaV9LIl0sWzAsMSwiXFxQaGlfSyAodF8xKV9cXGFzdCB0XzFeXFxhc3QiXSxbMSwwLCIodF8yKV9cXGFzdCB0XzJeXFxhc3QgXFxQaGlfSyJdLFsxLDEsIih0XzIpX1xcYXN0IFxcUGhpX3sodF5cXHByaW1lKV5cXGFzdCBLfXRfMV5cXGFzdCAuIl0sWzAsMSwiXFxQaGlfSyhcXHpldGFeMSkiLDJdLFswLDIsIlxcemV0YV4yX3tcXFBoaV9LfSJdLFsyLDMsIih0XzIpX1xcYXN0IChcXGJldGFeSykiXSxbMSwzLCJcXGFscGhhXktfe3RfMV5cXGFzdH0iXV0=
		\begin{tikzcd}
			{\Phi_K} & {(t_2)_\ast t_2^\ast \Phi_K} \\
			{\Phi_K (t_1)_\ast t_1^\ast} & {(t_2)_\ast \Phi_{(t^\prime)^\ast K}t_1^\ast .}
			\arrow["{\zeta^2_{\Phi_K}}", from=1-1, to=1-2]
			\arrow["{\Phi_K(\zeta^1)}"', from=1-1, to=2-1]
			\arrow["{(t_2)_\ast (\alpha^K)}", from=1-2, to=2-2]
			\arrow["{\beta^K_{t_1^\ast}}", from=2-1, to=2-2]
		\end{tikzcd}
	\end{equation}
\end{remark}

\begin{proposition}
	\label{prop:pullbacks_is_a_morphism_of_adjoint_integral_transforms}
	Consider \Cref{setup:fm_cube_pullback} where $t$ is quasi-affine.
	Let $K\in D^b_{\operatorname{coh}}(Y_1\times_S Y_2)$ be relatively perfect over $Y_2$.
	Denote by $K^\prime$ the kernel obtained in \Cref{prop:right_adjoint_pullback}.
	Consider the diagram:
	\begin{displaymath}
		% https://q.uiver.app/#q=WzAsNCxbMCwwLCJEX3tcXG9wZXJhdG9ybmFtZXtxY319KFlfMSkiXSxbMCwxLCJEX3tcXG9wZXJhdG9ybmFtZXtxY319KFlfMikiXSxbMSwwLCJEX3tcXG9wZXJhdG9ybmFtZXtxY319KFlfMVxcdGltZXNfU1QpIl0sWzEsMSwiRF97XFxvcGVyYXRvcm5hbWV7cWN9fShZXzJcXHRpbWVzX1MgVCkiXSxbMCwyLCJcXG1hdGhiZntMfXRfMV4qIl0sWzEsMywiXFxtYXRoYmZ7TH10XzJeXFxhc3QiLDJdLFswLDEsIlxcUGhpX0siLDIseyJvZmZzZXQiOjJ9XSxbMSwwLCJcXFBoaV97S15cXHByaW1lfSIsMix7Im9mZnNldCI6Mn1dLFsyLDMsIlxcUGhpX3tcXG1hdGhiZntMfSh0XlxccHJpbWUpXlxcYXN0IEt9IiwyLHsib2Zmc2V0IjoyfV0sWzMsMiwiXFxQaGlfe1xcbWF0aGJme0x9KHReXFxwcmltZSleXFxhc3QgS15cXHByaW1lfSIsMix7Im9mZnNldCI6Mn1dLFs2LDcsIiIsMix7ImxldmVsIjoxLCJzdHlsZSI6eyJuYW1lIjoiYWRqdW5jdGlvbiJ9fV0sWzgsOSwiIiwyLHsibGV2ZWwiOjEsInN0eWxlIjp7Im5hbWUiOiJhZGp1bmN0aW9uIn19XV0=
		\begin{tikzcd}
			{D_{\operatorname{qc}}(Y_1)} & {D_{\operatorname{qc}}(Y_1\times_ST)} \\
			{D_{\operatorname{qc}}(Y_2)} & {D_{\operatorname{qc}}(Y_2\times_S T).}
			\arrow["{\mathbf{L}t_1^\ast }", from=1-1, to=1-2]
			\arrow[""{name=0, anchor=center, inner sep=0}, "{\Phi_K}"', shift right=2, from=1-1, to=2-1]
			\arrow[""{name=1, anchor=center, inner sep=0}, "{\Phi_{\mathbf{L}(t^\prime)^\ast K}}"', shift right=2, from=1-2, to=2-2]
			\arrow[""{name=2, anchor=center, inner sep=0}, "{\Phi_{K^\prime}}"', shift right=2, from=2-1, to=1-1]
			\arrow["{\mathbf{L}t_2^\ast}"', from=2-1, to=2-2]
			\arrow[""{name=3, anchor=center, inner sep=0}, "{\Phi_{\mathbf{L}(t^\prime)^\ast K^\prime}}"', shift right=2, from=2-2, to=1-2]
			\arrow["\dashv"{anchor=center}, draw=none, from=0, to=2]
			\arrow["\dashv"{anchor=center}, draw=none, from=1, to=3]
		\end{tikzcd}
	\end{displaymath}
	Then, under this diagram, the isomorphism
	coming from \Cref{lem:fm_cube_flat_base_change_natural_isomorphism}
	\begin{displaymath}
		\alpha^{K^\prime}\colon \mathbf{L}t_1^\ast\Phi_{K^\prime}
		\to \Phi_{\mathbf{L}(t^\prime)^\ast K^\prime} \mathbf{L}t_2^\ast
	\end{displaymath}
	has a left mate
	\begin{displaymath}
		\lambda^K\colon \Phi_{\mathbf{L}(t^\prime)^\ast K} \mathbf{L}t_1^\ast \to \mathbf{L}t_2^\ast\Phi_K
	\end{displaymath}
	that is an isomorphism.
\end{proposition}

\begin{remark}
	\Cref{prop:pullbacks_is_a_morphism_of_adjoint_integral_transforms} does not follow formally merely because $\alpha^{K^\prime}$ is an isomorphism.
	Indeed, a mate correspondence need not preserve invertibility, see \Cref{rem:mate_of_iso_needs_not_be_iso}.
	A similar result is \cite[Proposition 4.9]{Bergh/Schnurer:2020},
	stated for $t$ flat and $K$ perfect.
	The argument in loc.\ cit.\ breaks down if any of these two hypotheses is removed.
\end{remark}

\begin{proof}
	[Proof of \Cref{prop:pullbacks_is_a_morphism_of_adjoint_integral_transforms}]
	Denote by $\eta$ and $\epsilon$ the unit and counit of $\Phi_K\dashv\Phi_{K^\prime}$.
	Denote by $\eta^\prime$ and $\epsilon^\prime$ the unit and
	counit of $\Phi_{\mathbf{L}(t^\prime)^\ast K} \dashv \Phi_{\mathbf{L}(t^\prime)^\ast K^\prime}$.
	By \Cref{prop:right_adjoint_pullback}, the latter adjunction exists.
	By \Cref{lem:mates_bijection}, the component of $\lambda^K$ at $E\in D_{\operatorname{qc}}(Y_1)$ is the composite
	\begin{displaymath}
		\begin{aligned}
			\Phi_{\mathbf{L}(t^\prime)^\ast K}
			\mathbf{L}t_1^\ast E
			&\xrightarrow{\Phi_{\mathbf{L}(t^\prime)^\ast K} \mathbf{L}t_1^\ast(\eta_E)}
			\Phi_{\mathbf{L}(t^\prime)^\ast K}
			\mathbf{L}t_1^\ast\Phi_{K^\prime}\Phi_K(E)
			\\&\xrightarrow{\Phi_{\mathbf{L}(t^\prime)^\ast K}(\alpha^{K^\prime}_{\Phi_K(E)})}
			\Phi_{\mathbf{L}(t^\prime)^\ast K}
			\Phi_{\mathbf{L}(t^\prime)^\ast K^\prime}
			\mathbf{L}t_2^\ast\Phi_K(E)
			\\&\xrightarrow{\epsilon^\prime_{\mathbf{L}t_2^\ast\Phi_K(E)}}
			\mathbf{L}t_2^\ast\Phi_K(E).
		\end{aligned}
	\end{displaymath}
	Choose $E\in D_{\operatorname{qc}}(Y_1)$ and $A\in D_{\operatorname{qc}}(Y_2\times_S T)$.
	Precomposition with $\lambda^K_E$ gives a morphism
	\begin{displaymath}
		\operatorname{Hom}(\mathbf{L}t_2^\ast\Phi_K(E),A)
		\to \operatorname{Hom}
		(\Phi_{\mathbf{L}(t^\prime)^\ast K} \mathbf{L}t_1^\ast E, A ).
	\end{displaymath}
	We claim that this morphism is the composite
	\begin{equation}
		\label{eq:adjunctions}
		\begin{aligned}
			\operatorname{Hom}
			(\mathbf{L}t_2^\ast\Phi_K(E),A)
			&\cong \operatorname{Hom} (\Phi_K(E),\mathbf{R}(t_2)_\ast A)
			\\&\cong\operatorname{Hom} (E,\Phi_{K^\prime}\mathbf{R}(t_2)_\ast A)
			\\&\xrightarrow{\beta^{K^\prime}_A\circ-}\operatorname{Hom} ( E,\mathbf{R}(t_1)_\ast \Phi_{\mathbf{L}(t^\prime)^\ast K^\prime}(A))
			\\&\cong\operatorname{Hom} ( \mathbf{L}t_1^\ast E, \Phi_{\mathbf{L}(t^\prime)^\ast K^\prime}(A))
			\\&\cong \operatorname{Hom} ( \Phi_{\mathbf{L}(t^\prime)^\ast K} \mathbf{L}t_1^\ast E, A).
		\end{aligned}
	\end{equation}
	Here $\beta^{K^\prime}$ is the natural isomorphism of \Cref{lem:fm_cube_pullback_natural_isomorphism} applied to $K^\prime$.

	We verify the claim.
	Choose a morphism $h\colon \mathbf{L}t_2^\ast\Phi_K(E)\to A$.
	Let $h^\flat\colon \Phi_K(E)\to\mathbf{R}(t_2)_\ast A$ be its adjunct with respect to  $\mathbf{L}t_2^\ast\dashv\mathbf{R}(t_2)_\ast$.
	Under the adjunction $\Phi_{\mathbf{L}(t^\prime)^\ast K} \dashv \Phi_{\mathbf{L}(t^\prime)^\ast K^\prime}$, the morphism $h\circ\lambda^K_E$ corresponds to
	\begin{equation}
		\label{eq:adjunctions2}
		\begin{aligned}
			\mathbf{L}t_1^\ast E
			&\xrightarrow{\mathbf{L}t_1^\ast(\eta_E)}
			\mathbf{L}t_1^\ast \Phi_{K^\prime}\Phi_K(E)
			\xrightarrow{\alpha^{K^\prime}_{\Phi_K(E)}} \Phi_{\mathbf{L}(t^\prime)^\ast K^\prime} \mathbf{L}t_2^\ast\Phi_K(E)
			\\&\xrightarrow{ \Phi_{\mathbf{L}(t^\prime)^\ast K^\prime}(h)}
			\Phi_{\mathbf{L}(t^\prime)^\ast K^\prime}(A).
		\end{aligned}
	\end{equation}
	Moreover, \Cref{lem:integral_transform_is_morphism_of_adjoints} applied to
	$K^\prime$ says that
	\begin{displaymath}
		\alpha^{K^\prime}\colon
		\mathbf{L}t_1^\ast\Phi_{K^\prime}
		\to \Phi_{\mathbf{L}(t^\prime)^\ast K^\prime} \mathbf{L}t_2^\ast
	\end{displaymath}
	and
	\begin{displaymath}
		\beta^{K^\prime}\colon
		\Phi_{K^\prime}\mathbf{R}(t_2)_\ast
		\to\mathbf{R}(t_1)_\ast
		\Phi_{\mathbf{L}(t^\prime)^\ast K^\prime}
	\end{displaymath}
	are mates.
	Applying the Hom-set formulation of
	\Cref{lem:morph_adj} to $h\colon \mathbf{L}t_2^\ast\Phi_K(E)\to A$ shows that the adjunct of $\Phi_{\mathbf{L}(t^\prime)^\ast K^\prime}(h) \circ \alpha^{K^\prime}_{\Phi_K(E)}$ is
	\begin{displaymath}
		\beta^{K^\prime}_A
		\circ \Phi_{K^\prime}(h^\flat)
		\colon \Phi_{K^\prime}\Phi_K(E)
		\to \mathbf{R}(t_1)_\ast \Phi_{\mathbf{L}(t^\prime)^\ast K^\prime}(A).
	\end{displaymath}
	After precomposing with the unit $\eta_E\colon E\to\Phi_{K^\prime}\Phi_K(E)$, the adjunct of \eqref{eq:adjunctions2} is
	\begin{displaymath}
		E
		\xrightarrow{\eta_E} \Phi_{K^\prime}\Phi_K(E)
		\xrightarrow{\Phi_{K^\prime}(h^\flat)} \Phi_{K^\prime}\mathbf{R}(t_2)_\ast A
		\xrightarrow{\beta^{K^\prime}_A} \mathbf{R}(t_1)_\ast \Phi_{\mathbf{L}(t^\prime)^\ast K^\prime}(A).
	\end{displaymath}
	This is the image of $h$ under the composite \eqref{eq:adjunctions}.
	Hence, precomposition with $\lambda^K_E$ agrees with \eqref{eq:adjunctions}.

	Every arrow in \eqref{eq:adjunctions} is a bijection.
	Indeed, the first, second, fourth, and fifth arrows are adjunction isomorphisms, while the third is induced by the natural isomorphism $\beta^{K^\prime}$.
	Consequently,
	\begin{displaymath}
		\operatorname{Hom}
		(\mathbf{L}t_2^\ast\Phi_K(E),A)
		\to \operatorname{Hom}( \Phi_{\mathbf{L}(t^\prime)^\ast K} \mathbf{L}t_1^\ast E, A)
	\end{displaymath}
	induced by precomposition with $\lambda^K_E$ is a bijection for every $A\in D_{\operatorname{qc}}(Y_2\times_S T)$.
	By \cite[\href{https://stacks.math.columbia.edu/tag/001P}{Tag 001P}]{StacksProject}, $\lambda^K_E$ is an isomorphism.
	Since this holds for every $E\in D_{\operatorname{qc}}(Y_1)$, $\lambda^K$ is a natural isomorphism.
\end{proof}

\begin{remark}
	\label{rem:ff_descent}
	Consider the situation as in \Cref{prop:pullbacks_is_a_morphism_of_adjoint_integral_transforms}.
	Denote $\epsilon$ and $\epsilon^\prime $ respectively for the counits of the adjunctions $\Phi_K\dashv\Phi_{K^\prime }$ and $\Phi_{\mathbf{L}(t^\prime)^\ast K}\dashv\Phi_{\mathbf{L}(t^\prime)^\ast K^\prime }$.
	By \Cref{lem:morph_adj}, it follows that
	\begin{equation}
		\label{eq:ff_descent1}
		\Phi_{\mathbf{L}(t^\prime)^\ast K^\prime} (\lambda^K) \circ \eta^\prime \mathbf{L}t^\ast_1 = \alpha^{K^\prime} \Phi_K \circ \mathbf{L}t^\ast_1 (\eta)
	\end{equation}
	and
	\begin{equation}
		\label{eq:ff_descent2}
		\mathbf{L}t^\ast_2 (\epsilon) \circ \lambda^K\Phi_{K^\prime} = \epsilon^\prime \mathbf{L}t_2^\ast \circ \Phi_{\mathbf{L}(t^\prime)^\ast K} (\alpha^{K^\prime}).
	\end{equation}
	In other words, for $E\in D_{\operatorname{qc}}(Y_1)$ and $A\in D_{\operatorname{qc}}(Y_2)$, the diagrams
	\begin{displaymath}
		\begin{tikzcd}[column sep=large,row sep=large]
			\mathbf{L}t_1^\ast E
			\arrow[r,"\eta^\prime_{\mathbf{L}t_1^\ast E}"]
			\arrow[d,"\mathbf{L}t_1^\ast (\eta_E)"']
			&
			\Phi_{\mathbf{L}(t^\prime)^\ast K^\prime}
			\Phi_{\mathbf{L}(t^\prime)^\ast K}
			(\mathbf{L}t_1^\ast E)
			\arrow[d,"\Phi_{\mathbf{L}(t^\prime)^\ast K^\prime}(\lambda^K_E)"]
			\\
			\mathbf{L}t_1^\ast \Phi_{K^\prime}\Phi_K(E)
			\arrow[r,"\alpha^{K^\prime}_{\Phi_K(E)}"']
			&
			\Phi_{\mathbf{L}(t^\prime)^\ast K^\prime}
			\mathbf{L}t_2^\ast \Phi_K(E)
		\end{tikzcd}
	\end{displaymath}
	and
	\begin{displaymath}
		\begin{tikzcd}[column sep=large,row sep=large]
			\Phi_{\mathbf{L}(t^\prime)^\ast K}
			\mathbf{L}t_1^\ast \Phi_{K^\prime}(A)
			\arrow[r,"\lambda^K_{\Phi_{K^\prime}(A)}"]
			\arrow[d,"\Phi_{\mathbf{L}(t^\prime)^\ast K}(\alpha^{K^\prime}_A)"']
			&
			\mathbf{L}t_2^\ast \Phi_K\Phi_{K^\prime}(A)
			\arrow[d,"\mathbf{L}t_2^\ast (\epsilon_A)"]
			\\
			\Phi_{\mathbf{L}(t^\prime)^\ast K}
			\Phi_{\mathbf{L}(t^\prime)^\ast K^\prime}
			(\mathbf{L}t_2^\ast A)
			\arrow[r,"\epsilon^\prime_{\mathbf{L}t_2^\ast A}"']
			&
			\mathbf{L}t_2^\ast A
		\end{tikzcd}
	\end{displaymath}
	commute.
\end{remark}

\begin{remark}
	On a locally Noetherian algebraic space, any closed point can be represented by a morphism of finite type.
	See \cite[\href{https://stacks.math.columbia.edu/tag/0H1U}{Tag 0H1U} and \href{https://stacks.math.columbia.edu/tag/06QX}{06QX}]{StacksProject}.
\end{remark}

\begin{lemma}
	\label{lem:equivalences_induced}
	Let $f_1\colon Y_1\to S$ and $f_2\colon Y_2\to S$ be proper flat morphisms of Noetherian algebraic spaces.
	Suppose $K\in D^b_{\operatorname{coh}}(Y_1\times_S Y_2)$ is relatively perfect over each $Y_i$.
	Then the following are equivalent:
	\begin{enumerate}
		\item \label{lem:equivalences_induced1} $\Phi_K$ is fully faithful (resp.\ an equivalence) on $D_{\operatorname{qc}}$
		\item \label{lem:equivalences_induced2} $\Phi_K$ restricts to a fully faithful functor (resp.\ an equivalence) on $D^b_{\operatorname{coh}}$
		\item \label{lem:equivalences_induced3} $\Phi_K$ restricts to a fully faithful functor (resp.\ an equivalence) on $\operatorname{Perf}$.
	\end{enumerate}
	In particular, these conditions for full faithfulness or equivalences can be tested on a single compact generator (see proof for a precise statement).
\end{lemma}

\begin{remark}
	The equivalence case of \Cref{lem:equivalences_induced} follows from \Cref{rmk:equivalences_induced},
	even without relative perfection assumptions on $K$.
	Nonetheless, \Cref{rmk:equivalences_induced} requires weak approximability.
	The following proof only uses generation techniques (plus the relative perfection assumptions) and can be seen as an alternative to the weak approximability argument.
\end{remark}

\begin{proof}
	[Proof of \Cref{lem:equivalences_induced}]
	Set $K^\prime := \operatorname{\mathbf{R}\mathcal{H}\! \mathit{om}} (K,(f_1^\prime)^! \mathcal{O}_{Y_2})$ to be the kernel of the integral transform obtained in \Cref{prop:right_adjoint_pullback} which is right adjoint to $\Phi_K$ on $D_{\operatorname{qc}}$.
	Also, let $K^{\prime\prime}=\operatorname{\mathbf{R}\mathcal{H}\! \mathit{om}} (K,(f_2^\prime)^! \mathcal{O}_{Y_1})$ be the kernel of the integral transform obtained in \Cref{prop:left_adjoint_pullback} which is left adjoint to $\Phi_K$ on $D_{\operatorname{qc}}$.
	We start with full faithfulness.

	As $K$ is relatively perfect over both $Y_i$, \Cref{lem:induced_dbcoh_perf_preservation_upon_pullback_with_t_affine} implies $\Phi_K$ restricts to a functor on $D^b_{\operatorname{coh}}$ and $\operatorname{Perf}$.
	This shows $\eqref{lem:equivalences_induced1} \implies \eqref{lem:equivalences_induced2} \implies \eqref{lem:equivalences_induced3}$.

	Assume \eqref{lem:equivalences_induced3} holds.
	We are in the situation where $\Phi_K$ restricts to a fully faithful functor on $\operatorname{Perf}$.
	Set $e\colon \Phi_{K^{\prime\prime}} \circ \Phi_K \to 1$ to be the counit of the adjoint pair $\Phi_{K^{\prime\prime}}$ and $\Phi_K$.
	Then, by \Cref{lem:Ballard_relative_perf_implies_dual_is_such}, $K$ is relatively perfect over each $Y_i$ implies $K^{\prime\prime}$ is relatively perfect over $Y_1$.
	Hence, \Cref{lem:induced_dbcoh_perf_preservation_upon_pullback_with_t_affine} says $\Phi_{K^{\prime\prime}}$ restricts to a functor on $\operatorname{Perf}$.
	Moreover, from $\Phi_K$ restricting to a fully faithful functor on $\operatorname{Perf}$, \cite[\href{https://stacks.math.columbia.edu/tag/07RB}{Tag 07RB}]{StacksProject} ensures that $e_P$ is an isomorphism for all $P\in \operatorname{Perf}(Y_1)$.

	Consider the strictly full subcategory $\mathcal{T}$ of objects $B\in D_{\operatorname{qc}}(Y_1)$ such that $e_B$ is an isomorphism.
	It is verified this forms a localizing subcategory.
	Coupled with the above, $\operatorname{Perf}(Y_1)\subseteq \mathcal{T}$.
	Yet, $D_{\operatorname{qc}}(Y_1)$ is compactly generated, and so \Cref{lem:localizing_iff_cpt_gen} implies $D_{\operatorname{qc}}(Y_1) = \mathcal{T}$.
	Then \cite[\href{https://stacks.math.columbia.edu/tag/07RB}{Tag 07RB}]{StacksProject} says $\Phi_K$ is fully faithful on $D_{\operatorname{qc}}$.

	Next, we prove the case for equivalence.
	By \Cref{lem:equivalence_or_fully_faithful_via_compacts}, $\eqref{lem:equivalences_induced1} \iff \eqref{lem:equivalences_induced3}$.
	Note that $\Phi_{K^{\prime\prime}}$ and $\Phi_K$ restrict to adjoint pair on $\operatorname{Perf}$.
	Then, from \Cref{cor:relatively_perfect_y2_implies_right_adjoint_restrict_to_dbcoh}, $\Phi_K$ and $\Phi_{K^\prime}$ restrict to an adjoint pair on $D^b_{\operatorname{coh}}$.

	Assume $\eqref{lem:equivalences_induced1}$ holds.
	The unit and counit the adjoint pair $\Phi_K$ and $\Phi_{K^\prime}$ for objects in $D^b_{\operatorname{coh}}$ remain bounded and pseudocoherent.
	Hence, \cite[\href{https://stacks.math.columbia.edu/tag/07RB}{Tag 07RB}]{StacksProject} yields $\eqref{lem:equivalences_induced1} \implies \eqref{lem:equivalences_induced2}$.

	Next, suppose \eqref{lem:equivalences_induced2} holds.
	Since $\Phi_K$ restricts to a functor on $\operatorname{Perf}$, it is a fully faithful functor on $\operatorname{Perf}$.
	By the fully faithful case, we know that $\Phi_K$ is fully faithful on $D_{\operatorname{qc}}$.
	Set $\mathcal{S}$ to be the strictly full subcategory of $B\in D_{\operatorname{qc}}(Y_2)$ such the counit $\epsilon_B\colon (\Phi_K \circ \Phi_{K^\prime})(B) \to B$ is an isomorphism.
	Then $\mathcal{S}$ is localizing.
	Note that $\epsilon_B$ is a morphism of $D^b_{\operatorname{coh}}(Y_2)$ for all $B\in D^b_{\operatorname{coh}}(Y_2)$.
	It follows that $D^b_{\operatorname{coh}}(Y_2) \subseteq \mathcal{S}$.
	Thus, $\operatorname{Perf}(Y_2)\subseteq \mathcal{S}$, and so \Cref{lem:localizing_iff_cpt_gen} implies $\mathcal{S}=D_{\operatorname{qc}}(Y_2)$.
	Consequently, $\eqref{lem:equivalences_induced2} \implies \eqref{lem:equivalences_induced1}$.

	We prove the last claim.
	Set $\eta\colon 1\to\Phi_{K^\prime}\circ\Phi_K$ and $\epsilon\colon \Phi_K\circ \Phi_{K^\prime} \to 1$ respectively for the unit and counit.
	For $i=1,2$, choose a compact generator $G_i$ of $D_{\operatorname{qc}}(Y_i)$.
	We claim that $\Phi_K$ is fully faithful (resp.\ an equivalence) if, and only if, $\eta_{G_1}$ is an isomorphism (resp.\ both $\eta_{G_1}$ and $\epsilon_{G_2}$ are isomorphisms).
	To see, let $\mathcal{T}^\prime\subseteq\mathcal{T}$ and $\mathcal{S}^\prime\subseteq\mathcal{S}$ be the strictly full subcategories where $\eta$ and $\epsilon$ are isomorphisms.
	Again, one can check that $\mathcal{T}^\prime$ and $\mathcal{S}^\prime$ are closed under shifts, iterated cones, and small coproducts (hence, homotopy colimits).
	Hence, \cite[\href{https://stacks.math.columbia.edu/tag/09SN}{Tag 09SN}]{StacksProject} says the condition implies $\mathcal{T} \subseteq \mathcal{T}^\prime$
	(resp.\ $\mathcal{T} \subseteq \mathcal{T}^\prime$ and $\mathcal{S} \subseteq \mathcal{S}^\prime$).
\end{proof}

\begin{proposition}
	\label{prop:descent}
	Consider \Cref{setup:fm_cube_pullback}.
	Let $K\in D_{\operatorname{qc}}(Y_1 \times_{S}  Y_2 )$ be pseudocoherent and relatively perfect over each $Y_i$.
	If $\Phi_{\mathbf{L}(t^\prime)^\ast K}$ is fully faithful (resp.\ an equivalence) on $D^b_{\operatorname{coh}}$ where $t$ is one of the following:
	\begin{enumerate}
		\item \label{prop:descent1} $t$ is quasi-affine and faithfully flat
		%\item \label[proposition]{prop:descent2} $t\colon \operatorname{Spec}(k)\to S$ is any morphism from a field (so here, we want the condition to all for all such morphisms)
		\item \label{prop:descent2_closed} for every finite type representative $t\colon \operatorname{Spec}(k)\to S$ of a closed point $p\in S$
	\end{enumerate}
	then so is $\Phi_{K}$.
\end{proposition}

\begin{proof}
	We prove the claim for full faithfulness.
	We first recall some properties satisfied by any morphism $t$ appearing in the statement.
	In \eqref{prop:descent2_closed}, every morphism from a field to $S$ is affine \cite[\href{https://stacks.math.columbia.edu/tag/09TF}{Tag 09TF}]{StacksProject}, and hence quasi-affine.
	By \Cref{thm:bounded_pseudocoherence_perfectness_faithfully_flat_affine}, $\mathbf{L}(t^\prime)^\ast K$ is relatively perfect over each $Y_i \times_{S} T$. By \Cref{prop:right_adjoint_pullback}, $\Phi_K$ and $\Phi_{K^\prime}$, as well as $\Phi_{\mathbf{L}(t^\prime)^\ast K}$ and $\Phi_{\mathbf{L}(t^\prime)^\ast K^\prime}$, form adjoint pairs on $D_{\operatorname{qc}}$. Denote by $\eta$ (resp.\ $\eta^\prime$) the unit of the adjunction $\Phi_K$ and $\Phi_{K^\prime}$ (resp.\ $\Phi_{\mathbf{L}(t^\prime)^\ast K}$ and $\Phi_{\mathbf{L}(t^\prime)^\ast K^\prime}$) on $D_{\operatorname{qc}}$.
	By \Cref{cor:relatively_perfect_y2_implies_right_adjoint_restrict_to_dbcoh}, these adjoint pairs restrict to $D^b_{\operatorname{coh}}$. Moreover, \Cref{cor:preservation} implies that $\Phi_K$ and $\Phi_{\mathbf{L}(t^\prime)^\ast K}$ restrict to functors on $\operatorname{Perf}$.

	We start by proving \eqref{prop:descent1}.
	Let $E\in D^b_{\operatorname{coh}}(Y_1)$.
	Since $t$ is faithfully flat, $t_1$ is as well, and so $\mathbf{L}t^\ast_1 E\in D^b_{\operatorname{coh}}(Y_1 \times_{S} T)$.
	By \eqref{eq:ff_descent1}, it follows that
	\begin{displaymath}
		\operatorname{cone}(\eta^\prime_{\mathbf{L}t_1^\ast E}) \cong \mathbf{L}t_1^\ast \operatorname{cone}(\eta_E).
	\end{displaymath}
	%%NOTE: Apply octahedral axiom to the equation
	However, $\Phi_{\mathbf{L}(t^\prime)^\ast K}$ restricts to a fully faithful functor on $D^b_{\operatorname{coh}}$, and so $\operatorname{cone}(\eta^\prime_{\mathbf{L}t_1^\ast E})\cong 0$.
	Then \Cref{lem:faithfully_flat_is_conservative} implies that $\operatorname{cone}(\eta_E)\cong 0$.
	Hence, $\Phi_K$ restricts to a fully faithful functor on $D^b_{\operatorname{coh}}$ because $E$ was arbitrary.

	Lastly we check \eqref{prop:descent2_closed}. By \Cref{lem:equivalences_induced}, it suffices to check $\Phi_K$ restricts to a fully faithful functor on $D^b_{\operatorname{coh}}(Y_1)$.
	Choose any $P\in D^b_{\operatorname{coh}}(Y_1)$.
	Since $\Phi_K$ restricts to a functor on $D^b_{\operatorname{coh}}$, we know that $\operatorname{cone}(\eta_P) \in D^b_{\operatorname{coh}}(Y_1)$.
	Let $p\in |Y_1|$ be a closed point. Suppose $s\colon \operatorname{Spec}(k)\to Y_1$ represents $p$ and is of finite type.
	Since $f_1$ is proper, \Cref{lem:factor_residual_gerbes} says $t:= f_1 \circ s$ represents the closed point $f_1 (p)\in |S|$ and is of finite type.

	Consider the following commutative diagram
	\begin{displaymath}
		% https://q.uiver.app/#q=WzAsNSxbMiwxLCJcXG9wZXJhdG9ybmFtZXtTcGVjfShrKSJdLFsyLDIsIlxcbWF0aGNhbHtTfS4iXSxbMSwyLCJcXG1hdGhjYWx7WX1fMSJdLFsxLDEsIlxcbWF0aGNhbHtZfV8xIFxcdGltZXNfUyBcXG9wZXJhdG9ybmFtZXtTcGVjfShrKSJdLFswLDAsIlxcb3BlcmF0b3JuYW1le1NwZWN9KGspIl0sWzAsMSwidCJdLFszLDIsInRfMSJdLFszLDAsImZfMV5cXHByaW1lIl0sWzQsMywiaCIsMV0sWzQsMCwiMV97XFxvcGVyYXRvcm5hbWV7U3BlY30oayl9Il0sWzQsMiwicyIsMl0sWzIsMSwiZl8xIiwyXV0=
		\begin{tikzcd}
			{\operatorname{Spec}(k)} && \\
			& {Y_1 \times_S \operatorname{Spec}(k)} & {\operatorname{Spec}(k)} \\
			& {Y_1} & {S.}
			\arrow["h"{description}, from=1-1, to=2-2]
			\arrow["{1_{\operatorname{Spec}(k)}}", bend left=30, from=1-1, to=2-3]
			\arrow["s"', bend right=30, from=1-1, to=3-2]
			\arrow["{f_1^\prime}", from=2-2, to=2-3]
			\arrow["{t_1}", from=2-2, to=3-2]
			\arrow["t", from=2-3, to=3-3]
			\arrow["{f_1}"', from=3-2, to=3-3]
		\end{tikzcd}
	\end{displaymath}
	By \eqref{eq:ff_descent1}, we obtain that
	\begin{displaymath}
		\operatorname{cone}(\eta^\prime_{\mathbf{L}t_1^\ast P}) \cong \mathbf{L}t_1^\ast \operatorname{cone}(\eta_P).
	\end{displaymath}
	By hypothesis, $\Phi_{\mathbf{L}(t^\prime)^\ast K}$ is fully faithful, and so $\operatorname{cone}(\eta^\prime_{\mathbf{L}t_1^\ast P})\cong 0$ (see \Cref{lem:equivalences_induced}).
	Then
	\begin{displaymath}
		\mathbf{L}s^\ast \operatorname{cone}(\eta_P) \cong \mathbf{L}h^\ast \mathbf{L}t_1^\ast \operatorname{cone}(\eta_P) \cong 0.
	\end{displaymath}
	By \Cref{lem:support_is_cohomological_for_finite_type_cohomology}, we see that $p\not\in \operatorname{Supp}(\operatorname{cone}(\eta_P))$.
	Since $p$ was an arbitrary closed point of $Y_1$, it follows that $\operatorname{Supp}(\operatorname{cone}(\eta_P))=\emptyset$.
	In other words, $\operatorname{cone}(\eta_P)$ is the zero object, which completes the proof for full faithfulness.

	To prove the last claim for an equivalence, the same argument above can be adapted with the counit of the adjoint pair.
	This completes the proof.
\end{proof}

\begin{proposition}
	\label{prop:ascending}
	Consider \Cref{setup:fm_cube_pullback} where $t$ is quasi-affine.
	Let $K\in D_{\operatorname{qc}}(Y_1 \times_{S}  Y_2 )$ be pseudocoherent and relatively perfect over each  $Y_i$.
	If $\Phi_K$ is fully faithful (resp.\ an equivalence) on $D^b_{\operatorname{coh}}$, then so is $\Phi_{\mathbf{L}(t^\prime)^\ast K}$.
\end{proposition}

\begin{proof}
	Let $G_i$ be a compact generator of $D_{\operatorname{qc}}(Y_i)$.
	By \Cref{lem:equivalences_induced}, full faithfulness of $\Phi_K$ is equivalent to $\eta_{G_1}$ being an isomorphism.
	Then \eqref{eq:ff_descent1}, together with
	the fact that $\lambda^K$ and $\alpha^{K^\prime}$ are isomorphisms,
	implies that $\eta^\prime_{\mathbf{L}t_1^\ast G_1}$ is an isomorphism.
	Since $t_1$ is quasi-affine, \Cref{lem:pullback_compact_generator} implies $\mathbf{L}t_1^\ast G_1$ is a compact generator.
	Hence, \Cref{lem:equivalences_induced} implies that
	$\Phi_{\mathbf{L}(t^\prime)^\ast K}$ is fully faithful.

	If $\Phi_K$ is an equivalence, then $\epsilon_{G_2}$ is an isomorphism.
	Then \eqref{eq:ff_descent2} implies that $\epsilon^\prime_{\mathbf{L}t_2^\ast G_2}$ is an isomorphism.
	By \Cref{lem:pullback_compact_generator}, $\mathbf{L}t_2^\ast G_2$ is a compact generator, and so \Cref{lem:equivalences_induced} implies that $\Phi_{\mathbf{L}(t^\prime)^\ast K}$ is an equivalence. 
\end{proof}

\begin{theorem}
	\label{thm:descent_ascent}
	Consider \Cref{setup:fm_cube_pullback}.
	Let $K\in D_{\operatorname{qc}}(Y_1 \times_{S}  Y_2 )$ be pseudocoherent and  relatively perfect over each $Y_i$.
	Then the following are equivalent:
	\begin{enumerate}
		\item \label{thm:descent_ascent1} $\Phi_{K}$ is fully faithful (resp.\ an equivalence) on $D^b_{\operatorname{coh}}$
		\item \label{thm:descent_ascent2} $\Phi_{\mathbf{L}(t^\prime)^\ast K}$ is fully faithful (resp.\ an equivalence) on $D^b_{\operatorname{coh}}$ for any quasi-affine morphism $t$
		%\item \label{thm:descent_ascent3} $\Phi_{\mathbf{L}(t^\prime)^\ast K}$ is fully faithful (resp.\ an equivalence) on $D^b_{\operatorname{coh}}$ for every morphism $t\colon \operatorname{Spec}(k)\to S$ from a field
		\item \label{thm:descent_ascent3} $\Phi_{\mathbf{L}(t^\prime)^\ast K}$ is fully faithful (resp.\ an equivalence) on $D^b_{\operatorname{coh}}$ for every $t\colon \operatorname{Spec}(k)\to S$ of finite type that represents any closed point $p\in |S|$
		\item \label{thm:descent_ascent4} for any closed point $p\in |S|$ there exists a representative $t\colon \operatorname{Spec}(k)\to S$ of $p$ such that $\Phi_{\mathbf{L}(t^\prime)^\ast K}$ is fully faithful (resp.\ an equivalence) on $D^b_{\operatorname{coh}}$.
		%%NOTE: $t$ need not be of finite type for this condition
	\end{enumerate}
\end{theorem}

\begin{proof}
	By \Cref{prop:ascending}, $\eqref{thm:descent_ascent1}\implies \eqref{thm:descent_ascent2}$.
	By \cite[\href{https://stacks.math.columbia.edu/tag/09TF}{Tag 09TF}]{StacksProject}, $\eqref{thm:descent_ascent2}\implies \eqref{thm:descent_ascent3}$.
	It is obvious that $\eqref{thm:descent_ascent3}\implies \eqref{thm:descent_ascent4}$.
	By \Cref{prop:descent}, $\eqref{thm:descent_ascent3} \implies \eqref{thm:descent_ascent1}$.
	We prove that $\eqref{thm:descent_ascent4}\implies \eqref{thm:descent_ascent3}$.

	For each closed point $p\in S$, let $t\colon \operatorname{Spec}(k)\to S$ be a representative of $p$ such that $\Phi_{\mathbf{L}(t^\prime)^\ast K}$ is fully faithful (resp.\ an equivalence) on $D^b_{\operatorname{coh}}$.
	Fix an arbitrary representative $\operatorname{Spec}(k^\prime)\to S$ of $p$.
	There exist a field $\ell$ and a commutative diagram
	\begin{displaymath}
		% https://q.uiver.app/#q=WzAsNCxbMSwwLCJcXG9wZXJhdG9ybmFtZXtTcGVjfShrKSJdLFsxLDEsIlMuIl0sWzAsMSwiXFxvcGVyYXRvcm5hbWV7U3BlY30oa15cXHByaW1lKSJdLFswLDAsIlxcb3BlcmF0b3JuYW1le1NwZWN9KFxcZWxsKSJdLFswLDEsInQiXSxbMiwxLCJ0XlxccHJpbWUiLDJdLFszLDIsImEiLDJdLFszLDAsImIiXV0=
		\begin{tikzcd}
			{\operatorname{Spec}(\ell)} & {\operatorname{Spec}(k)} \\
			{\operatorname{Spec}(k^\prime)} & {S.}
			\arrow["b", from=1-1, to=1-2]
			\arrow["a"', from=1-1, to=2-1]
			\arrow["t", from=1-2, to=2-2]
			\arrow["{t^\prime}"', from=2-1, to=2-2]
		\end{tikzcd}
	\end{displaymath}
	By \Cref{prop:ascending}, we can ascend full faithfulness or equivalences from $\operatorname{Spec}(k)$ to $\operatorname{Spec}(\ell)$.
	The canonical morphism $a\colon \operatorname{Spec}(\ell) \to \operatorname{Spec}(k^\prime)$ is affine and faithfully flat.
	Hence, \Cref{prop:descent} allows us to descend full faithfulness or equivalences from $\operatorname{Spec}(\ell)$ to $\operatorname{Spec}(k^\prime)$.
	This completes the proof.
\end{proof}

%%%%%%%%%%%%%%%%%%%%%%%%%%%%%%%%%%%
\part{Ramifications}
\label{sec:part3}
%%%%%%%%%%%%%%%%%%%%%%%%%%%%%%%%%%%

%%%%%%%%%%%%%%%%%%%%%%%%%%%%%%%%%%%
\section{Elliptic fibrations}
\label{sec:autoequivalence}
%%%%%%%%%%%%%%%%%%%%%%%%%%%%%%%%%%%

\begin{lemma}
	\label{lem:ideal_sheaf_relatively_perfect}
	Let $f\colon Y \to X$ be proper flat morphism of Noetherian algebraic spaces.
	Denote by $C_f$ the cone of the unit $\mathcal{O}_{Y \times_{X} Y}\to \mathbf{R} (\Delta_f)_\ast \mathcal{O}_{Y}$ where $\Delta_f\colon Y \to Y \times_{X} Y$ is the diagonal.
	Set $p,q\colon Y \times_{X} Y \to Y$ the canonical morphisms.
	Then $C_f$ is relatively perfect for $p$ and $q$.
\end{lemma}

\begin{proof}
	We prove the case for $q$ because the other is similar (e.g.\ switch notation).
	There exists a fibered square
	\begin{displaymath}
		% https://q.uiver.app/#q=WzAsNCxbMCwxLCJcXG1hdGhjYWx7WX0iXSxbMSwxLCJcXG1hdGhjYWx7WH0uIl0sWzEsMCwiXFxtYXRoY2Fse1l9Il0sWzAsMCwiXFxtYXRoY2Fse1l9IFxcdGltZXNfe1xcbWF0aGNhbHtYfX0gXFxtYXRoY2Fse1l9Il0sWzAsMSwiZiIsMl0sWzIsMSwiZiJdLFszLDIsInEiXSxbMywwLCJwIiwyXV0=
		\begin{tikzcd}
			{Y \times_{X} Y} & {Y} \\
			{Y} & {X.}
			\arrow["q", from=1-1, to=1-2]
			\arrow["p"', from=1-1, to=2-1]
			\arrow["f", from=1-2, to=2-2]
			\arrow["f"', from=2-1, to=2-2]
		\end{tikzcd}
	\end{displaymath}
	Consider the distinguished triangle
	\begin{displaymath}
		C_f [-1] \to \mathcal{O}_{Y \times_{X} Y} \to \mathbf{R} (\Delta_f)_\ast \mathcal{O}_{Y} \to C_f.
	\end{displaymath}
	Since $f$ has affine and proper diagonal, the canonical morphism $(\Delta_f)_\ast \mathcal{O}_{Y} \to \mathbf{R} (\Delta_f)_\ast \mathcal{O}_{Y}$ is an isomorphism.
	Now, if we extend to the long exact sequence in cohomology with the distinguished triangle above, we see that $C_f [-1]$ is a coherent sheaf on $Y \times_{X} Y$ because $\mathcal{O}_{Y \times_{X} Y} \to (\Delta_f)_\ast \mathcal{O}_{Y}$ is surjective since $\Delta_f$ is a closed immersion.
	Let $E\in D^b_{\operatorname{qc}}(Y)$. Since $q$ is flat, we know that $\mathbf{L} q^\ast E \in D^b_{\operatorname{qc}}(Y \times_{X} Y)$.
	Consider the distinguished triangle obtained by tensoring with $\mathbf{L} q^\ast E$,
	\begin{displaymath}
		C_f [-1] \otimes^{\mathbf{L}} \mathbf{L} q^\ast E \to \mathcal{O}_{Y \times_{X} Y} \otimes^{\mathbf{L}} \mathbf{L} q^\ast E  \to \mathbf{R} (\Delta_f)_\ast \mathcal{O}_{Y} \otimes^{\mathbf{L}} \mathbf{L} q^\ast E  \to C_f \otimes^{\mathbf{L}} \mathbf{L} q^\ast E .
	\end{displaymath}
	If we can show that $\mathbf{R} (\Delta_f)_\ast \mathcal{O}_{Y} \otimes^{\mathbf{L}} \mathbf{L} q^\ast E \in D^b_{\operatorname{qc}}(Y \times_{X} Y)$, then we are done.
	%%NOTE: Indeed, since $p$ is proper, flat, and concentrated, we can apply \Cref{cor:preservation}.
	Now, by projection formula, we have
	\begin{displaymath}
		\begin{aligned}
			\mathbf{R} (\Delta_f)_\ast \mathcal{O}_{Y} \otimes^{\mathbf{L}} \mathbf{L} q^\ast E
			&\cong \mathbf{R} (\Delta_f)_\ast \mathbf{L} \Delta_f^\ast \mathbf{L} q^\ast E
			\\&\cong \mathbf{R} (\Delta_f)_\ast \mathbf{L} (q\circ \Delta_f)^\ast E
			\\&\cong \mathbf{R} (\Delta_f)_\ast E && (q\circ \Delta_f= 1_{Y}).
		\end{aligned}
	\end{displaymath}
	Since $\Delta_f$ is affine, it follows that $\mathbf{R}(\Delta_f)_\ast$ is $t$-exact and conservative.
	Hence, $\mathbf{R} (\Delta_f)_\ast E \in D^b_{\operatorname{qc}}(Y \times_{X} Y)$, which completes the proof.
\end{proof}

\begin{reminder}
	\label{rem:elliptic_fibration}
	Recall that an \textbf{elliptic $X$-fibration} is a proper Gorenstein morphism $Y \to X$ between Noetherian algebraic spaces such that for each $p\in |X|$ there is a representative $\operatorname{Spec}(k)\to X$ of $p$ such that $Y \times_X \operatorname{Spec}(k)$ is a geometrically integral $1$-dimensional $k$-scheme of genus one\footnote{We follow \cite[\href{https://stacks.math.columbia.edu/tag/0BY7}{Tag 0BY7}]{StacksProject} for the definition of `genus'.}.
\end{reminder}

\begin{example}
	\label{ex:elliptic_fibrations}
	Let $k$ be a field.
	Consider a geometrically integral $1$-dimensional $k$-scheme $Y$ of genus one with trivial dualizing sheaf.
	Then the canonical morphism $Y\times_k Y \to Y$ is an elliptic fibration.

	To see, consider a morphism $\operatorname{Spec}(L)\to Y$ from a field. Note that the base change of $Y\times_k Y$ along $\operatorname{Spec}(L)\to Y$ is the base change $Y$ along $L/k$.
	As $Y \times_k \operatorname{Spec}(L)$ is geometrically integral (and hence, connected), \cite[\href{https://stacks.math.columbia.edu/tag/0FD2}{Tag 0FD2}]{StacksProject} says that
	\begin{displaymath}
		H^0(Y \times_k \operatorname{Spec}(L), \mathcal{O}_{Y \times_k  \operatorname{Spec}(L) }) = L.
	\end{displaymath}
	From base change, we see that $Y \times_k \operatorname{Spec}(L)$ is proper and Gorenstein over $L$.
	Hence, by \cite[\href{https://stacks.math.columbia.edu/tag/0A26}{Tag 0A26}]{StacksProject}, $Y \times_k \operatorname{Spec}(L)$ is projective over $L$.
	Moreover, \cite[\href{https://stacks.math.columbia.edu/tag/0BY9}{Tag 0BY9}]{StacksProject} tells us that $Y \times_k \operatorname{Spec}(L)$ has genus one.
	Applying \cite[\href{https://stacks.math.columbia.edu/tag/0FW1}{Tag 0FW1}]{StacksProject}, it follows that $Y \times_k \operatorname{Spec}(L)$ has trivial dualizing sheaf.
	Consequently, the base change of $Y\times_k Y$ over $\operatorname{Spec}(L)\to Y$ is projective, geometrically integral, of Krull dimension one, genus one, and with trivial dualizing sheaf.

	More generally, if $Y\to S$ is an elliptic fibration, then $Y\times_S Y\to Y$ is an elliptic fibration.
	This can be argued in a similar fashion to the case $S=\operatorname{Spec}(k)$ for a field $k$ above.
	In fact, if $Y\to S$ and $X\to S$ are elliptic fibrations, then the natural projections of $X\times_S Y$ to $X$ and $Y$ are elliptic fibrations.
	Again, the argument follows analogously to that above.
\end{example}

\begin{proposition}
	\label{prop:derived_pullback_ideal_sheaf}
	Let $X$ be a Noetherian algebraic space.
	Let $T\to X$ be an affine morphism from a Noetherian algebraic space and $f\colon Y \to X$ a proper flat morphism from a Noetherian algebraic space.
	Denote by $\mathcal{I}_{\Delta_f}$ (resp.\ $\mathcal{I}_{\Delta_f^\prime}$) the ideal sheaf of the diagonal morphism $\Delta_f\colon Y \to Y \times_X Y$ (resp.\ $\Delta_{f^\prime}\colon Y \times_X T \to Y \times_X Y\times_X T$).
	Then $\mathbf{L}(t^\prime)^\ast \mathcal{I}_{\Delta_f}\cong \mathcal{I}_{\Delta_f^\prime}$ where $t^\prime\colon Y \times_X Y\times_X T \to Y \times_X Y$ is the canonical morphism.
\end{proposition}

\begin{proof}
	Since $f$ is proper, the diagonal $\Delta_f$ is a closed immersion. Consider the commutative diagram obtained by base changes,
	\begin{displaymath}
		% https://q.uiver.app/#q=WzAsOCxbNCwxLCJUIl0sWzQsMywiXFxtYXRoY2Fse1N9LiJdLFsxLDMsIlkiXSxbMiwyLCJZIl0sWzEsMSwiWSBcXHRpbWVzX1ggIFQiXSxbMiwwLCJZIFxcdGltZXNfWCBUIl0sWzAsMCwiWVxcdGltZXNfWCBZIFxcdGltZXNfWCBUIl0sWzAsMiwiWSBcXHRpbWVzX1ggWSJdLFswLDEsInQiXSxbMiwxLCJmIiwyXSxbMywxLCJmIl0sWzQsMCwiZl5cXHByaW1lIl0sWzUsMCwiZl5cXHByaW1lIl0sWzQsMiwicyIsMCx7ImxhYmVsX3Bvc2l0aW9uIjo3MH1dLFs1LDMsInMiLDIseyJsYWJlbF9wb3NpdGlvbiI6ODAsInN0eWxlIjp7ImJvZHkiOnsibmFtZSI6ImRhc2hlZCJ9fX1dLFs2LDQsInBeXFxwcmltZSJdLFs2LDUsInFeXFxwcmltZSJdLFs3LDMsInEiLDAseyJsYWJlbF9wb3NpdGlvbiI6MzAsInN0eWxlIjp7ImJvZHkiOnsibmFtZSI6ImRhc2hlZCJ9fX1dLFs3LDIsInAiLDJdLFs2LDcsInReXFxwcmltZSIsMl1d
		\begin{tikzcd}
			{Y\times_X Y \times_X T} && {Y \times_X T} && \\
			& {Y \times_X  T} &&& T \\
			{Y \times_X Y} && Y \\
			& Y &&& {X.}
			\arrow["{q^\prime}", from=1-1, to=1-3]
			\arrow["{p^\prime}", from=1-1, to=2-2]
			\arrow["{t^\prime}"', from=1-1, to=3-1]
			\arrow["{f^\prime}", from=1-3, to=2-5]
			\arrow["s"'{pos=0.8}, dashed, from=1-3, to=3-3]
			\arrow["{f^\prime}", from=2-2, to=2-5]
			\arrow["s"{pos=0.7}, from=2-2, to=4-2]
			\arrow["t", from=2-5, to=4-5]
			\arrow["q"{pos=0.3}, dashed, from=3-1, to=3-3]
			\arrow["p"', from=3-1, to=4-2]
			\arrow["f", from=3-3, to=4-5]
			\arrow["f"', from=4-2, to=4-5]
		\end{tikzcd}
	\end{displaymath}
	Note that base change shows that $p,q$ are flat.
	We obtain a diagram of fibered squares
	\begin{displaymath}
		% https://q.uiver.app/#q=WzAsOCxbMywwLCJUIl0sWzMsMiwiXFxtYXRoY2Fse1N9Il0sWzIsMiwiWSJdLFsyLDAsIlkgXFx0aW1lc19YIFQiXSxbMSwwLCJZXFx0aW1lc19YIFkgXFx0aW1lc19YIFQiXSxbMSwyLCJZIFxcdGltZXNfWCBZIl0sWzAsMiwiWSJdLFswLDAsIllcXHRpbWVzX1ggVCJdLFswLDEsInQiXSxbMiwxLCJmIiwyXSxbMywwLCJmXlxccHJpbWUiXSxbMywyLCJzIl0sWzQsMywicF5cXHByaW1lIl0sWzUsMiwicCIsMl0sWzQsNSwidF5cXHByaW1lIiwyXSxbNiw1LCJcXERlbHRhX2YiLDJdLFs3LDQsIlxcRGVsdGFfe2ZeXFxwcmltZX0iXSxbNyw2LCJzIiwyXV0=
		\begin{tikzcd}
			{Y\times_X T} & {Y\times_X Y \times_X T} & {Y \times_X T} & T \\
			\\
			Y & {Y \times_X Y} & Y & {X.}
			\arrow["{\Delta_{f^\prime}}", from=1-1, to=1-2]
			\arrow["s"', from=1-1, to=3-1]
			\arrow["{p^\prime}", from=1-2, to=1-3]
			\arrow["{t^\prime}"', from=1-2, to=3-2]
			\arrow["{f^\prime}", from=1-3, to=1-4]
			\arrow["s", from=1-3, to=3-3]
			\arrow["t", from=1-4, to=3-4]
			\arrow["{\Delta_f}"', from=3-1, to=3-2]
			\arrow["p"', from=3-2, to=3-3]
			\arrow["f"', from=3-3, to=3-4]
		\end{tikzcd}
	\end{displaymath}
	The leftmost square occurs from the fact that $\Delta_{f^\prime}$ is the base change of $\Delta_f$ along $t^\prime$ (see e.g.\ proof of \cite[\href{https://stacks.math.columbia.edu/tag/03KL}{Tag 03KL}]{StacksProject}).
	%%NOTE: See Lem 2 of 'Magic Squares' in note folder.
	Since $t$ is affine, base change tells us both $s$ and $t^\prime$ are affine, and hence representable by schemes.
	By \cite[\href{https://stacks.math.columbia.edu/tag/07U8}{Tag 07U8}]{StacksProject}, there is an isomorphism
	\begin{displaymath}
		\begin{aligned}
			\mathbf{R} (\Delta_{f^\prime})_\ast \mathbf{L} s^\ast \mathcal{O}_Y
			&\cong \mathbf{R} (\Delta_{f^\prime})_\ast s^\ast \mathcal{O}_Y
			\\&\cong \mathbf{R} (\Delta_{f^\prime})_\ast \mathcal{O}_{Y\times_X T}
			\\&\cong (t^\prime)^\ast \mathbf{R} (\Delta_f)_\ast \mathcal{O}_Y \\&\cong (t^\prime)^\ast (\Delta_f)_\ast \mathcal{O}_Y
		\end{aligned}
	\end{displaymath}
	where $(t^\prime)^\ast$ is underived.

	We want to show that the canonical morphism $\mathbf{L} (t^\prime)^\ast \mathbf{R} (\Delta_f)_\ast \mathcal{O}_Y \to (t^\prime)^\ast \mathbf{R} (\Delta_f)_\ast \mathcal{O}_Y$ is an isomorphism.
	This amounts to checking that $\mathbf{L} (t^\prime)^\ast \mathbf{R} (\Delta_f)_\ast \mathcal{O}_Y$ is concentrated in degree zero.
	Since $\mathbf{R} t^\prime_\ast$ is $t$-exact and conservative, it suffices to check that $\mathbf{R} t^\prime_\ast \mathbf{L} (t^\prime)^\ast \mathbf{R} (\Delta_f)_\ast \mathcal{O}_Y$ is concentrated in degree zero.
	Indeed, $t$-exactness gives
	\begin{displaymath}
		\mathcal{H}^i (\mathbf{R} t^\prime_\ast \mathbf{L} (t^\prime)^\ast \mathbf{R} (\Delta_f)_\ast \mathcal{O}_Y) \cong t^\prime_\ast \mathcal{H}^i (\mathbf{L} (t^\prime)^\ast \mathbf{R} (\Delta_f)_\ast \mathcal{O}_Y),
	\end{displaymath}
	whereas conservativeness shows
	\begin{displaymath}
		t^\prime_\ast \mathcal{H}^i (\mathbf{L} (t^\prime)^\ast \mathbf{R} (\Delta_f)_\ast \mathcal{O}_Y)\cong 0
		\implies \mathcal{H}^i (\mathbf{L} (t^\prime)^\ast \mathbf{R} (\Delta_f)_\ast \mathcal{O}_Y)\cong 0.
	\end{displaymath}
	By projection formula, we have
	\begin{displaymath}
		\mathbf{R} t^\prime_\ast \mathbf{L} (t^\prime)^\ast \mathbf{R} (\Delta_f)_\ast \mathcal{O}_Y \cong \mathbf{R} t^\prime_\ast \mathcal{O}_{Y\times_X Y \times_X T} \otimes^{\mathbf{L}}\mathbf{R} (\Delta_f)_\ast \mathcal{O}_Y.
	\end{displaymath}
	Note that
	\begin{displaymath}
		\mathcal{O}_{Y\times_X Y \times_X T} \cong \mathbf{L}(f^\prime \circ p^\prime)^\ast \mathcal{O}_{T}.
	\end{displaymath}
	Using flat base change,
	\begin{displaymath}
		\begin{aligned}
			\mathbf{R} t^\prime_\ast \mathcal{O}_{Y\times_X Y \times_X T}
			&\cong \mathbf{R} t^\prime_\ast \mathbf{L}(f^\prime \circ p^\prime)^\ast \mathcal{O}_{T}
			\\&\cong \mathbf{L}(f\circ p)^\ast \mathbf{R}t_\ast \mathcal{O}_{T}.
		\end{aligned}
	\end{displaymath}
	Then we obtain:
	\begin{displaymath}
		\begin{aligned}
			\mathbf{R} t^\prime_\ast
			& \mathcal{O}_{Y\times_X Y \times_X T} \otimes^{\mathbf{L}}\mathbf{R} (\Delta_f)_\ast \mathcal{O}_Y
			\\&\cong \mathbf{L}(f\circ p)^\ast \mathbf{R}t_\ast \mathcal{O}_{T} \otimes^{\mathbf{L}}\mathbf{R} (\Delta_f)_\ast \mathcal{O}_Y && \textrm{(use above)}
			\\&\cong \mathbf{R} (\Delta_f)_\ast \mathbf{L} \Delta_f^\ast \mathbf{L}(f\circ p)^\ast \mathbf{R}t_\ast \mathcal{O}_{T} && \textrm{(projection formula)}
			\\&\cong \mathbf{R} (\Delta_f)_\ast \mathbf{L}(f\circ p \circ \Delta_f)^\ast \mathbf{R}t_\ast \mathcal{O}_{T} && \textrm{(pseudofunctoriality)}
			\\&\cong \mathbf{R} (\Delta_f)_\ast \mathbf{L}f^\ast \mathbf{R}t_\ast \mathcal{O}_{T} && (p \circ \Delta_f = 1_Y).
		\end{aligned}
	\end{displaymath}
	Since $t$ is affine, the canonical morphism $t_\ast \mathcal{O}_{T} \to \mathbf{R}t_\ast \mathcal{O}_{T}$ is an isomorphism.
	Moreover, $f$ is flat, and so the canonical morphism $f^\ast \mathbf{R}t_\ast \mathcal{O}_{T} \to \mathbf{L}f^\ast \mathbf{R}t_\ast \mathcal{O}_{T}$ is an isomorphism.
	Hence, $\mathbf{L}f^\ast \mathbf{R}t_\ast \mathcal{O}_{T}$ is concentrated in degree zero.
	As $\Delta_f$ is affine, the canonical morphism $(\Delta_f)_\ast \mathbf{L}f^\ast \mathbf{R}t_\ast \mathcal{O}_{T} \to \mathbf{R} (\Delta_f)_\ast \mathbf{L}f^\ast \mathbf{R}t_\ast \mathcal{O}_{T}$ is an isomorphism.
	Thus, $\mathbf{R}^j (\Delta_f)_\ast \mathbf{L}f^\ast \mathbf{R}t_\ast \mathcal{O}_{T}\cong 0$ for all $j\not=0$.

	% Consider the distinguished triangle
	% \begin{displaymath}
	%     \mathcal{I}_{\Delta_f} \to \mathcal{O}_{Y\times_X Y}\to \mathbf{R} (\Delta_f)_\ast \mathcal{O}_Y \to \mathcal{I}_{\Delta_f} [1].
	% \end{displaymath}
	% Applying $\mathbf{L}(t^\prime)^\ast$ yields a distinguished triangle
	% \begin{displaymath}
	%     \mathbf{L}(t^\prime)^\ast \mathcal{I}_{\Delta_f} \to \mathbf{L}(t^\prime)^\ast \mathcal{O}_{Y\times_X Y}\to \mathbf{L}(t^\prime)^\ast \mathbf{R} (\Delta_f)_\ast \mathcal{O}_Y \to \mathbf{L}(t^\prime)^\ast \mathcal{I}_{\Delta_f} [1].
	% \end{displaymath}
	% This gives a long exact sequence in cohomology
	% \begin{displaymath}
	%     \cdots \to \mathcal{H}^i ( \mathbf{L}(t^\prime)^\ast \mathcal{I}_{\Delta_f}) \to \mathcal{H}^i (\mathbf{L}(t^\prime)^\ast \mathcal{O}_{Y\times_X Y}) \to \mathcal{H}^i (\mathbf{L}(t^\prime)^\ast \mathbf{R} (\Delta_f)_\ast \mathcal{O}_Y) \to \cdots.
	% \end{displaymath}
	Recall the natural isomorphism $\mathbf{L}(t^\prime)^\ast \mathcal{O}_{Y\times_X Y} \cong \mathcal{O}_{Y\times_X Y \times_X T}$.
	There is a morphism of distinguished triangles
	\begin{displaymath}
		% https://q.uiver.app/#q=WzAsOCxbMSwwLCJcXG1hdGhiZntMfSh0XlxccHJpbWUpXlxcYXN0IFxcbWF0aGNhbHtPfV97WVxcdGltZXNfWCBZfSJdLFsxLDEsIlxcbWF0aGNhbHtPfV97WVxcdGltZXNfWCBZIFxcdGltZXNfWCBUfSJdLFsyLDAsIlxcbWF0aGJme0x9KHReXFxwcmltZSleXFxhc3QgXFxtYXRoYmZ7Un0gKFxcRGVsdGFfZilfXFxhc3QgXFxtYXRoY2Fse099X1kiXSxbMiwxLCJcXG1hdGhiZntSfSAoXFxEZWx0YV9mXlxccHJpbWUpX1xcYXN0IFxcbWF0aGNhbHtPfV97WVxcdGltZXNfWCBUfSJdLFswLDEsIlxcbWF0aGNhbHtJfV97XFxEZWx0YV9mXlxccHJpbWV9Il0sWzAsMCwiXFxtYXRoYmZ7TH0odF5cXHByaW1lKV5cXGFzdCBcXG1hdGhjYWx7SX1fe1xcRGVsdGFfZn0iXSxbMywwLCJcXG1hdGhiZntMfSh0XlxccHJpbWUpXlxcYXN0IFxcbWF0aGNhbHtJfV97XFxEZWx0YV9mfVsxXSJdLFszLDEsIlxcbWF0aGNhbHtJfV97XFxEZWx0YV9mXlxccHJpbWV9WzFdLiJdLFswLDFdLFswLDJdLFsxLDNdLFsyLDNdLFs0LDFdLFs1LDBdLFs1LDRdLFsyLDZdLFszLDddLFs2LDddXQ==
		\begin{tikzcd}
			{\mathbf{L}(t^\prime)^\ast \mathcal{I}_{\Delta_f}} & {\mathbf{L}(t^\prime)^\ast \mathcal{O}_{Y\times_X Y}} & {\mathbf{L}(t^\prime)^\ast \mathbf{R} (\Delta_f)_\ast \mathcal{O}_Y} & {\mathbf{L}(t^\prime)^\ast \mathcal{I}_{\Delta_f}[1]} \\
			{\mathcal{I}_{\Delta_f^\prime}} & {\mathcal{O}_{Y\times_X Y \times_X T}} & {\mathbf{R} (\Delta_f^\prime)_\ast \mathcal{O}_{Y\times_X T}} & {\mathcal{I}_{\Delta_f^\prime}[1].}
			\arrow[from=1-1, to=1-2]
			\arrow[from=1-1, to=2-1]
			\arrow[from=1-2, to=1-3]
			\arrow[from=1-2, to=2-2]
			\arrow[from=1-3, to=1-4]
			\arrow[from=1-3, to=2-3]
			\arrow[from=1-4, to=2-4]
			\arrow[from=2-1, to=2-2]
			\arrow[from=2-2, to=2-3]
			\arrow[from=2-3, to=2-4]
		\end{tikzcd}
	\end{displaymath}
	Since we have shown that the two inner morphisms are isomorphisms, it follows that the morphism $\mathbf{L}(t^\prime)^\ast \mathcal{I}_{\Delta_f} \to \mathcal{I}_{\Delta_f^\prime}$ is an isomorphism.
\end{proof}

\begin{proof}
	[Proof of \Cref{thm:elliptic_fibration_autoequivalence}]
	We start with an observation.
	In this paragraph, fix $p\in |X|$. Choose a representative $\operatorname{Spec}(k)\to X$ of $p$ such that $Y \times_X \operatorname{Spec}(k)$ is geometrically connected, of Krull dimension one, genus one, and with trivial dualizing sheaf.
	Let $\overline{k}$ be an algebraic closure of $k$. Since $Y \times_X \operatorname{Spec}(k)$ is geometrically integral, \cite[\href{https://stacks.math.columbia.edu/tag/0FD2}{Tag 0FD2}]{StacksProject} shows that
	\begin{displaymath}
		H^0(Y \times_X \operatorname{Spec}(\overline{k}), \mathcal{O}_{Y \times_X \operatorname{Spec}(\overline{k})}) = \overline{k}.
	\end{displaymath}
	By base change, $Y \times_X \operatorname{Spec}(\overline{k})$ is proper and Gorenstein over $\overline{k}$, and so \cite[\href{https://stacks.math.columbia.edu/tag/0A26}{Tag 0A26}]{StacksProject} implies $Y \times_X \operatorname{Spec}(\overline{k})$ is also projective over $\overline{k}$. Moreover, from \cite[\href{https://stacks.math.columbia.edu/tag/0BY9}{Tag 0BY9}]{StacksProject}, $Y \times_X \operatorname{Spec}(\overline{k})$ has genus one. Using \cite[\href{https://stacks.math.columbia.edu/tag/0FW1}{Tag 0FW1}]{StacksProject}, we see that $Y \times_X \operatorname{Spec}(\overline{k})$ has trivial dualizing sheaf.
	Thus, for every $p\in |X|$, we have a representative of $p$ from an algebraically closed field such that the fiber of $f$ is projective, geometrically connected, of Krull dimension one, genus one, and with trivial dualizing sheaf.
	Denote by $K^\prime$ the kernel of the integral transform from \Cref{prop:right_adjoint_pullback} which is right adjoint to $\Phi_{\mathcal{I}_{\Delta_f}}$ on $D_{\operatorname{qc}}$.

	Now, back to the proof.
	Since $f$ is proper, the diagonal $\Delta_f$ is a closed immersion.
	Hence, $\mathcal{I}_{\Delta_f}$ is the negative shift of the cone of the unit $\mathcal{O}_{Y \times_X Y}\to \mathbf{R} (\Delta_f)_\ast \mathcal{O}_Y$.
	Note that \Cref{lem:ideal_sheaf_relatively_perfect} shows that $\mathcal{I}_{\Delta_f}$ is relatively perfect over $Y$.
	By \Cref{prop:right_adjoint_pullback}, $K^\prime$ is bounded and pseudocoherent.
	Fix a closed point $p\in |X|$.
	Choose a representative $t\colon \operatorname{Spec}(k)\to X$ of $p$ with $k$ algebraically closed representatives constructed in the first paragraph.
	From \Cref{prop:right_adjoint_pullback}, we know that $\Phi_{\mathbf{L}(t^\prime)^\ast \mathcal{I}_{\Delta_f}}$ and $\Phi_{\mathbf{L}(t^\prime)^\ast K^\prime}$ remain an adjoint pair.

	By \Cref{prop:derived_pullback_ideal_sheaf}, we have $\mathbf{L}(t^\prime)^\ast \mathcal{I}_{\Delta_f} \cong \mathcal{I}_{\Delta_f^\prime}$.
	Consequently, we have that $\Phi_{\mathbf{L}(t^\prime)^\ast \mathcal{I}_{\Delta_f}}$ is naturally isomorphic to $\Phi_{\mathcal{I}_{\Delta_f^\prime}}$ as endofunctors on $D^b_{\operatorname{coh}}(Y\times_X \operatorname{Spec}(k))$.
	Since geometric integrality of the fibers implies irreducible components have the same Krull dimension (i.e.\ equidimensional fibers), \cite[Proposition 2.16 \& \S 3.4.1]{Ruiperez/Hernandez/Martin/SanchodeSalas:2009} implies that $\Phi_{\mathcal{I}_{\Delta_f^\prime}}$ is an autoequivalence of $D^b_{\operatorname{coh}}(Y\times_X \operatorname{Spec}(k))$.
	Now, using \Cref{thm:descent_ascent}, the desired claim follows since we have checked $\Phi_{\mathbf{L}(t^\prime)^\ast \mathcal{I}_{\Delta_f}}$ is an equivalence for all closed points.
\end{proof}

%%%%%%%%%%%%%%%%%%%%%%%%%%%%%%%%%%%

%%%%%%%%%%%%%%%%%%%%%%%%%%%%%%%%%%%

%%%%%%%%%%%%%%%%%%%%%%%%%%%%%%%%%%%

\section{Openness of loci}
\label{sec:loci}
%%%%%%%%%%%%%%%%%%%%%%%%%%%%%%%%%%%

\subsection{With relative perfect assumptions}
\label{sec:loci_rel_perf}

\begin{definition}
	\label{def:fm_locus}
	Consider \Cref{setup:fm_cube_pullback}.
	Let $K\in D_{\operatorname{qc}}(Y_1 \times_{S}  Y_2 )$.
	We say that $\Phi_K$ is \textbf{fully faithful} (resp.\ an \textbf{equivalence}) \textbf{at} $p\in |S|$
	if $\Phi_{\mathbf{L}(t^\prime)^\ast K}$ is such on $D_{\operatorname{qc}}$ for some representative $t$ of $p$.
	Denote by $\operatorname{fm}(K; Y_1,Y_2,S)$ (resp.\ $\operatorname{FM}(K; Y_1,Y_2,S)$)
	the collection of $p\in |S|$ where $\Phi_K$ is fully faithful (resp.\ an equivalence) at $p$.
	If $Y_1=Y_2=Y$,
	we write this as $\operatorname{fm}(K; Y,S)$
	(resp.\ $\operatorname{FM}(K; Y,S)$). In the case $S=\operatorname{Spec}(R)$ is an affine scheme,
	then we simply write $\operatorname{fm}(K; Y_1,Y_2,R)$ and $\operatorname{fm}(K; Y,R)$
	(resp.\ $\operatorname{FM}(K; Y_1,Y_2,R)$ and $\operatorname{FM}(K; Y,R)$).
	We will just write $\operatorname{fm}(K)$ and $\operatorname{FM}(K)$ whenever $Y_1,Y_2,S$ are understood.
\end{definition}

\begin{lemma}
	\label{lem:independence_for_rep_with_ff_or_eq}
	Consider \Cref{setup:fm_cube_pullback}.
	Let $K\in D_{\operatorname{qc}}(Y_1 \times_{S}  Y_2 )$ be pseudocoherent and relatively perfect over each $Y_i$.
	The condition that $\Phi_K$ is fully faithful (resp.\ an equivalence) at some $p\in |S|$ is independent of the representative of $p$.
\end{lemma}

\begin{proof}
	Choose some $p\in \operatorname{fm}(K; Y_1,Y_2,S)$, and let $t\colon \operatorname{Spec}(k)\to S$ be a representative where full faithfulness holds.
	Pick any representative $h\colon \operatorname{Spec}(k^\prime)\to S$ of $p$.
	There exist a field $\ell$ and a commutative diagram,
	\begin{displaymath}
		% https://q.uiver.app/#q=WzAsNCxbMSwwLCJcXG9wZXJhdG9ybmFtZXtTcGVjfShrXlxccHJpbWUpIl0sWzEsMSwiUy4iXSxbMCwxLCJcXG9wZXJhdG9ybmFtZXtTcGVjfShrKSJdLFswLDAsIlxcb3BlcmF0b3JuYW1le1NwZWN9KFxcZWxsKSJdLFswLDEsInQiXSxbMiwxLCJoIiwyXSxbMywyLCJiIiwyXSxbMywwLCJhIl1d
		\begin{tikzcd}
			{\operatorname{Spec}(\ell)} & {\operatorname{Spec}(k^\prime)} \\
			{\operatorname{Spec}(k)} & {S.}
			\arrow["a", from=1-1, to=1-2]
			\arrow["b"', from=1-1, to=2-1]
			\arrow["h", from=1-2, to=2-2]
			\arrow["t"', from=2-1, to=2-2]
		\end{tikzcd}
	\end{displaymath}
	Note that $a,b$ are affine and faithfully flat.
	The hypothesis is that full faithfulness occurs after base change along $t$.
	We can apply \Cref{prop:ascending} to realize full faithfulness occurs after base change along $t\circ b$.
	Since $a$ is affine and faithfully flat, \Cref{prop:descent} implies full faithfulness occurs after base change along $h$ (e.g.\ use that $h\circ a = t \circ b$).

	Lastly, let $p\in \operatorname{FM}(K; Y_1,Y_2,S)$.
	Choose a representative $t\colon \operatorname{Spec}(k)\to S$ for which the adjoint equivalence occurs.
	By \Cref{sec:base_change_for_adjoints}, there exists $K^\prime\in D^b_{\operatorname{coh}}(Y_1 \times_{S} Y_2)$ such that $\Phi_K$ is a left adjoint of $\Phi_{K^\prime}$ on $D_{\operatorname{qc}}$.
	As $\Phi_{\mathbf{L}(t^\prime)^\ast K}$ and $\Phi_{\mathbf{L}(t^\prime)^\ast K^\prime}$ form an adjoint equivalence, $\Phi_{\mathbf{L}(t^\prime)^\ast K^\prime}$ is an exact equivalence.
	Hence, $\Phi_{\mathbf{L}(t^\prime)^\ast K^\prime}$ restricts to an exact equivalence on the compacts (i.e.\ $\operatorname{Perf}$).
	This means $\mathbf{L}(t^\prime)^\ast K^\prime$ is relatively perfect over $Y_1 \times_S \operatorname{Spec}(k)$.
	Pick any representative $h\colon \operatorname{Spec}(k^\prime)\to S$ of $p$.
	There exist a field $\ell$ and a commutative diagram,
	\begin{displaymath}
		% https://q.uiver.app/#q=WzAsNCxbMSwwLCJcXG9wZXJhdG9ybmFtZXtTcGVjfShrXlxccHJpbWUpIl0sWzEsMSwiUy4iXSxbMCwxLCJcXG9wZXJhdG9ybmFtZXtTcGVjfShrKSJdLFswLDAsIlxcb3BlcmF0b3JuYW1le1NwZWN9KFxcZWxsKSJdLFswLDEsInQiXSxbMiwxLCJoIiwyXSxbMywyLCJiIiwyXSxbMywwLCJhIl1d
		\begin{tikzcd}
			{\operatorname{Spec}(\ell)} & {\operatorname{Spec}(k^\prime)} \\
			{\operatorname{Spec}(k)} & {S.}
			\arrow["a", from=1-1, to=1-2]
			\arrow["b"', from=1-1, to=2-1]
			\arrow["h", from=1-2, to=2-2]
			\arrow["t"', from=2-1, to=2-2]
		\end{tikzcd}
	\end{displaymath}
	Note that $a,b$ are affine and faithfully flat.
	By \Cref{lem:induced_dbcoh_perf_preservation_upon_pullback_with_t_affine}, the derived pullback of $\mathbf{L}(t^\prime)^\ast K^\prime$ along natural base change obtained from $b$ is relatively perfect over $Y_1 \times_S \operatorname{Spec}(\ell)$.
	Applying \Cref{prop:reflecting_bounded_pseudocoherence_perfectness1}, the derived pullback of $K^\prime$ along the canonical morphism obtained by base change along $h$ is relatively perfect over $Y_1 \times_S \operatorname{Spec}(k^\prime)$.
	This shows the derived pullback of $K^\prime$ is relatively $Y_1\times_S \operatorname{Spec}(k^\prime)$-perfect for all representatives $\operatorname{Spec}(k^\prime) \to S$ of $p$.
	Applying \Cref{lem:equivalences_induced}, we can argue like fully faithful case above to finish the proof.
\end{proof}

\begin{theorem}
	\label{thm:openness_for_equivalence_fullfaithful_dbcoh}
	Consider \Cref{setup:fm_cube_pullback}.
	Let $K\in D_{\operatorname{qc}}(Y_1 \times_{S}  Y_2 )$ be pseudocoherent and relatively perfect over each $Y_i$.
	Then the collection of $p\in |S|$, for which $\Phi_{\mathbf{L}(t^\prime)^\ast K}$ along some representative $t$ of $p$ is an equivalence (resp.\ fully faithful) on $D_{\operatorname{qc}}$, is an open subset of $|S|$.
\end{theorem}

\begin{proof}
	By \Cref{sec:base_change_for_adjoints}, there exists $K^\prime\in D^b_{\operatorname{coh}}(Y_1 \times_{S} Y_2)$ such that $\Phi_K$ is a left adjoint of $\Phi_{K^\prime}$ on $D_{\operatorname{qc}}$.
	In fact, the adjoint pair restricts to $D^b_{\operatorname{coh}}$.
	Denote by $\epsilon\colon \Phi_K \circ \Phi_{K^\prime}\to 1$ the counit and $\eta\colon 1 \to \Phi_{K^\prime} \circ \Phi_K$ the unit for this adjunction.
	Choose compact generators $G_i$ for $D_{\operatorname{qc}}(Y_i)$.
	Set $C_1 := \operatorname{cone}(\eta_{G_1})$ and $C_2 := \operatorname{cone}(\epsilon_{G_2})$.
	We claim that
	\begin{displaymath}
		f_1 (\operatorname{Supp}(C_1)) = |S|\setminus \operatorname{fm}(K; Y_1,Y_2,S)
	\end{displaymath}
	and
	\begin{displaymath}
		|S| \setminus \operatorname{FM}(K; Y_1,Y_2,S) = f_1 (\operatorname{Supp}(C_1))  \cup f_2 (\operatorname{Supp}(C_2)) .
	\end{displaymath}
	Since $\Phi_K$ and $\Phi_{K^\prime}$ restrict to an adjoint pair on $D^b_{\operatorname{coh}}$, each $C_1$ and $C_2$ are bounded pseudocoherent.
	Thus, as each $f_i$ is a closed morphism, the desired claim would follow.

	We prove the desired claim. For any $p\in |S|$ and representative $t\colon \operatorname{Spec}(k)\to S$ of $p$, \eqref{eq:ff_descent1} and \eqref{eq:ff_descent2} imply
	\begin{displaymath}
		\mathbf{L} t_1^\ast C_1 \cong \operatorname{cone}(\eta^\prime_{\mathbf{L} t_1^\ast G_1})
	\end{displaymath}
	and
	\begin{displaymath}
		\mathbf{L} t_2^\ast C_2 \cong \operatorname{cone}(\epsilon^\prime_{\mathbf{L} t_2^\ast G_2}) .
	\end{displaymath}
	where $\epsilon^\prime\colon \Phi_{\mathbf{L} (t^\prime)^\ast K} \circ \Phi_{\mathbf{L} (t^\prime)^\ast K^\prime}\to 1$ is the counit and $\eta^\prime\colon 1 \to \Phi_{\mathbf{L} (t^\prime)^\ast K^\prime} \circ \Phi_{\mathbf{L} (t^\prime)^\ast K}$ the unit for the adjoint pair $\Phi_{\mathbf{L} (t^\prime)^\ast K}$ and $\Phi_{\mathbf{L} (t^\prime)^\ast K^\prime}$.
	See \Cref{prop:right_adjoint_pullback}.
	Since $\mathbf{L}t_i^\ast G_i$ are compact generators of $D_{\operatorname{qc}}(Y_i \times_S \operatorname{Spec}(k))$, \Cref{lem:equivalences_induced} gives
	\begin{displaymath}
		\mathbf{L} t_1^\ast C_1 \cong 0 \iff \Phi_{\mathbf{L} (t^\prime)^\ast K} \textrm{ is fully faithful}
	\end{displaymath}
	and
	\begin{displaymath}
		\mathbf{L} t_1^\ast C_1 \cong 0 \textrm{ and }\mathbf{L} t_2^\ast C_2 \cong 0 \iff \Phi_{\mathbf{L} (t^\prime)^\ast K} \textrm{ is an equivalence}.
	\end{displaymath}
	Consequently, the claim follows from \Cref{lem:independence_for_rep_with_ff_or_eq} and \Cref{lem:support_is_cohomological_for_finite_type_cohomology}.
\end{proof}

\begin{corollary}
	\label{cor:deform}
	Consider \Cref{setup:fm_cube_pullback}.
	Let $K\in D_{\operatorname{qc}}(Y_1 \times_{S}  Y_2 )$ be pseudocoherent and relatively perfect over each $Y_i$.
	Then
	\begin{displaymath}
		t^{-1}(\operatorname{fm}(K)) = \operatorname{fm}(\mathbf{L}(t^\prime)^\ast K)
	\end{displaymath}
	and
	\begin{displaymath}
		t^{-1}(\operatorname{FM}(K)) = \operatorname{FM}(\mathbf{L}(t^\prime)^\ast K).
	\end{displaymath}
\end{corollary}

\begin{proof}
	By \Cref{cor:noetherian_base_change_for_relative_perfectness}, $\mathbf{L}(t^\prime)^\ast K$ is relatively perfect over each $Y_i \times_S T$.
	Moreover, for each $p\in |T|$ there exists a commutative diagram
	\begin{displaymath}
		% https://q.uiver.app/#q=WzAsNCxbMCwwLCJcXG9wZXJhdG9ybmFtZXtTcGVjfShcXGthcHBhKHApKSJdLFsxLDAsIlxcb3BlcmF0b3JuYW1le1NwZWN9KFxca2FwcGEodChwKSkpIl0sWzEsMSwiUyJdLFswLDEsIlQiXSxbMCwxXSxbMSwyXSxbMCwzXSxbMywyLCJ0IiwyXV0=
		\begin{tikzcd}
			{\operatorname{Spec}(\kappa(p))} & {\operatorname{Spec}(\kappa(t(p)))} \\
			T & S
			\arrow[from=1-1, to=1-2]
			\arrow[from=1-1, to=2-1]
			\arrow[from=1-2, to=2-2]
			\arrow["t"', from=2-1, to=2-2]
		\end{tikzcd}
	\end{displaymath}
	where $\kappa(-)$ denotes the residue field of a point in the algebraic space.
	See \cite[\href{https://stacks.math.columbia.edu/tag/0EMV}{Tag 0EMV}]{StacksProject}.
	Denote by $K_{\kappa(t(p))}$ and $(\mathbf{L}(t^\prime)^\ast K)_{\kappa(p)}$ respectively the base changes of $K$ and $\mathbf{L}(t^\prime)^\ast K$ along $\kappa(t(p)) \to S$ and $\kappa(p)\to T$.
	By \Cref{thm:descent_ascent}, if $\Phi_{K_{\kappa(t(p))}}$ is fully faithful on $D^b_{\operatorname{coh}}$, then $\Phi_{(\mathbf{L}(t^\prime)^\ast K)_{\kappa(p)}}$is fully faithful on $D^b_{\operatorname{coh}}$.
	Hence, $t(p) \in \operatorname{fm}(K)$ implies $p\in \operatorname{fm}(\mathbf{L}(t^\prime)^\ast K)$.
	On the other hand, if $p\in \operatorname{fm}(\mathbf{L}(t^\prime)^\ast K)$, then \Cref{prop:descent} implies $t(p)\in \operatorname{fm}(K)$ because $\operatorname{Spec}(\kappa(p)) \to \operatorname{Spec}(\kappa(t(p)))$ is affine and faithfully flat.
	Thus, the claim follows for full faithfulness.
	The proof for equivalences is similar.
\end{proof}

\begin{corollary}
	\label{cor:FM_partners_detects_surjectivity}
	Consider \Cref{setup:fm_cube_pullback}.
	Let $K\in D_{\operatorname{qc}}(Y_1 \times_{S}  Y_2 )$ be pseudocoherent and relatively perfect over each $Y_i$.
	\begin{enumerate}
		\item \label{cor:FM_partners_detects_surjectivity1} If $\Phi_K$ is fully faithful (resp.\ an equivalence) on $D^b_{\operatorname{coh}}$, then so is $\Phi_{\mathbf{L}(t^\prime)^\ast K}$.
		In particular, if $Y_1$ and $Y_2$ are Fourier--Mukai $S$-partners, then $Y_1\times_S T$ and $Y_2\times_S T$ are Fourier--Mukai $T$-partners.
		\item \label{cor:FM_partners_detects_surjectivity2} Suppose $t$ is surjective.
		If $\Phi_{\mathbf{L}(t^\prime)^\ast K}$ is fully faithful (resp.\ an equivalence) on $D^b_{\operatorname{coh}}$, then so is $\Phi_K$.
	\end{enumerate}
\end{corollary}

\begin{proof}
	We prove the claims for full faithfulness.
    The case of equivalences is similar.
	By \Cref{cor:noetherian_base_change_for_relative_perfectness}, $\mathbf{L}(t^\prime)^\ast K$ is pseudocoherent and relatively perfect over each $Y_i\times_S T$.
	We first show that $\Phi_K$ is fully faithful on $D^b_{\operatorname{coh}}$ if, and only if, $\operatorname{fm}(K)=|S|$.
	Suppose $\Phi_K$ is fully faithful on $D^b_{\operatorname{coh}}$.
	Pick $p\in |S|$.
    Choose $s\colon \operatorname{Spec}(k)\to S$ a representative of $p$.
	Since $s$ is affine \cite[\href{https://stacks.math.columbia.edu/tag/09TF}{Tag 09TF}]{StacksProject}, \Cref{prop:ascending} shows that the fiber of $\Phi_K$ at $s$ is fully faithful on $D^b_{\operatorname{coh}}$, and hence on $D_{\operatorname{qc}}$ by \Cref{lem:equivalences_induced}.
	Thus, $p\in \operatorname{fm}(K)$.
	Conversely, suppose $\operatorname{fm}(K)=|S|$.
	Then every closed point of $|S|$ admits a representative such that the corresponding fiber of $\Phi_K$ is fully faithful on $D_{\operatorname{qc}}$, and hence on $D^b_{\operatorname{coh}}$ by \Cref{lem:equivalences_induced}.
	By \Cref{thm:descent_ascent}, $\Phi_K$ is fully faithful on $D^b_{\operatorname{coh}}$.
	The same argument applies to $\mathbf{L}(t^\prime)^\ast K$ over $T$.

	By \Cref{cor:deform}, $t^{-1}(\operatorname{fm}(K)) = \operatorname{fm}(\mathbf{L}(t^\prime)^\ast K)$.
	If $\operatorname{fm}(K)=|S|$, then $\operatorname{fm}(\mathbf{L}(t^\prime)^\ast K) = |T|$.
    This proves the first claim of \eqref{cor:FM_partners_detects_surjectivity1}.
	If $t$ is surjective and $\operatorname{fm}(\mathbf{L}(t^\prime)^\ast K) = |T|$, then $\operatorname{fm}(K)=|S|$.
    This proves \eqref{cor:FM_partners_detects_surjectivity2}.
	Lastly, suppose $E\in D^b_{\operatorname{coh}}(Y_1\times_S Y_2)$ is such that $\Phi_E$ restricts to an equivalence on $D^b_{\operatorname{coh}}$.
	By \Cref{rem:rel_perfect_is_free}, $E$ is relatively perfect over each $Y_i$.
	Hence, the first claim of \eqref{cor:FM_partners_detects_surjectivity1} applies to $E$, and so $\Phi_{\mathbf{L}(t^\prime)^\ast E}$ restricts to an equivalence on $D^b_{\operatorname{coh}}$.
\end{proof}

\begin{remark}
	In forthcoming work, the authors will prove variations of \Cref{thm:openness_for_equivalence_fullfaithful_dbcoh} and \Cref{cor:deform} for singularity categories.
\end{remark}

\subsection{Without relative perfect assumptions}
\label{sec:loci__weak_approx}

Let $S$ be a Noetherian algebraic space.
Suppose $f_i\colon Y_i \to S$ are proper flat morphisms of algebraic spaces.
For $p\in |S|$, denote by $t_p\colon \operatorname{Spec}(\kappa(p)) \to S$ the canonical morphism from the residue field of $p$ \cite[\href{https://stacks.math.columbia.edu/tag/0EMV}{Tag 0EMV}]{StacksProject}.
The morphisms $t_p$ yield commutative cubes with fibered faces,
\begin{displaymath}
	% https://q.uiver.app/#q=WzAsOCxbNCwxLCJcXG9wZXJhdG9ybmFtZXtTcGVjfShcXGthcHBhKHApKSJdLFs0LDMsIlMiXSxbMSwzLCJZXzEiXSxbMywyLCJZXzIiXSxbMSwxLCJZXnBfMSJdLFszLDAsIllecF8yIl0sWzAsMCwiWV8xXnAgXFx0aW1lc197XFxrYXBwYShwKX0gWV5wXzIiXSxbMCwyLCJZXzFcXHRpbWVzX1MgWV8yIl0sWzAsMSwidF9wIl0sWzIsMSwiZl8xIiwyXSxbMywxLCJmXzIiXSxbNCwwLCJmXnBfMSJdLFs1LDAsImZecF8yIl0sWzQsMiwidF5wXzEiLDAseyJsYWJlbF9wb3NpdGlvbiI6NzB9XSxbNSwzLCJ0XzJecCIsMix7ImxhYmVsX3Bvc2l0aW9uIjo4MCwic3R5bGUiOnsiYm9keSI6eyJuYW1lIjoiZGFzaGVkIn19fV0sWzYsNCwiZ18xXnAiLDJdLFs2LDUsImdecF8yIiwyXSxbNywzLCJnXzIiLDAseyJsYWJlbF9wb3NpdGlvbiI6MTAsInN0eWxlIjp7ImJvZHkiOnsibmFtZSI6ImRhc2hlZCJ9fX1dLFs3LDIsImdfMSIsMl0sWzYsNywidF5cXHByaW1lX3AiLDJdXQ==
	\begin{tikzcd}
		{Y_1^p \times_{\kappa(p)} Y^p_2} &&& {Y^p_2} & \\
		& {Y^p_1} &&& {\operatorname{Spec}(\kappa(p))} \\
		{Y_1\times_S Y_2} &&& {Y_2} \\
		& {Y_1} &&& S
		\arrow["{g^p_2}"', from=1-1, to=1-4]
		\arrow["{g_1^p}"', from=1-1, to=2-2]
		\arrow["{t^\prime_p}"', from=1-1, to=3-1]
		\arrow["{f^p_2}", from=1-4, to=2-5]
		\arrow["{t_2^p}"'{pos=0.8}, dashed, from=1-4, to=3-4]
		\arrow["{f^p_1}", from=2-2, to=2-5]
		\arrow["{t^p_1}"{pos=0.7}, from=2-2, to=4-2]
		\arrow["{t_p}", from=2-5, to=4-5]
		\arrow["{g_2}"{pos=0.1}, dashed, from=3-1, to=3-4]
		\arrow["{g_1}"', from=3-1, to=4-2]
		\arrow["{f_2}", from=3-4, to=4-5]
		\arrow["{f_1}"', from=4-2, to=4-5]
	\end{tikzcd}
\end{displaymath}
where $Y^p_i:= Y_i \times_S \operatorname{Spec}(\kappa(p))$ for each $i$.
Given $E\in D_{\operatorname{qc}}(Y_1\times_S Y_2)$, let $E_p:= \mathbf{L}(t^\prime_p)^\ast E$.

\begin{lemma}
	Each face of the cube above is a tor-independent square.
\end{lemma}

\begin{proof}
	This follows from the fact that $f_i$ and $g_i$ are flat.
	This can be checked after taking \'{e}tale presentations, applying \cite[\href{https://stacks.math.columbia.edu/tag/08IQ}{Tag 08IQ}]{StacksProject} and reducing to the scheme case (see e.g.\ \cite[Lemma 22.93 \& Remark 22.94]{Gortz/Wedhorn:2020}).
	%%NOTE: In Tag 08IQ, apply (3) and use that base change along projections are flat at times.
\end{proof}

\begin{proposition}
	\label{prop:relatively_perfect_locus}
	Let $K\in D^b_{\operatorname{coh}}(Y_1 \times_S Y_2)$.
	Then the collection of $p\in |S|$ such that $\mathbf{L}(t^\prime_p)^\ast K$ is relatively perfect over $Y_1^p$ (resp.\ $Y_2^p$) forms an open subset.
\end{proposition}

\begin{proof}
	We prove the case for $Y_1$ because the other is symmetric.
	Denote by $S_1$ the collection of $p\in |S|$ such that $\mathbf{L}(t^\prime_p)^\ast K$ is relatively perfect over $Y_1^p$.
	Choose a compact generator $G_i$ of $D_{\operatorname{qc}}(Y_i)$ for each $i$.
	Since $\mathbf{L}(g_1)^\ast G_1 \otimes^{\mathbf{L}} \mathbf{L}(g_2)^\ast G_2 \in \operatorname{Perf}(Y_1\times_S Y_2)$ and $K$ is bounded pseudocoherent, it follows that
	\begin{displaymath}
		K\otimes^{\mathbf{L}} \mathbf{L} (g_1)^\ast G_1 \otimes^{\mathbf{L}} \mathbf{L}(g_2)^\ast G_2 \in D^b_{\operatorname{coh}}(Y_1 \times_S Y_2).
	\end{displaymath}
	Properness of $g_1$ yields that
	\begin{displaymath}
		\mathbf{R}(g_1)_\ast (K\otimes^{\mathbf{L}} \mathbf{L} (g_1)^\ast G_1 \otimes^{\mathbf{L}} \mathbf{L}(g_2)^\ast G_2) \in D^b_{\operatorname{coh}}(Y_1).
	\end{displaymath}
	By \Cref{lem:perfect_locus_open}, the collection $U$ of $y\in Y_1$ such that
	\begin{displaymath}
		\mathbf{L}h^\ast \mathbf{R}(g_1)_\ast (K\otimes^{\mathbf{L}} \mathbf{L} (g_1)^\ast G_1 \otimes^{\mathbf{L}} \mathbf{L}(g_2)^\ast G_2) \in D^b_{\operatorname{coh}}(k),
	\end{displaymath}
	where $h$ represents $y$, is an open subset of $|Y_1|$.
	Moreover, $f_1$ being proper implies $|f_1|(|Y_1| \setminus U) \subseteq |S|$ is closed.
	It suffices to show that $|f_1|(|Y_1|\setminus U) = |S| \setminus S_1$.

	Let $p\in |S| \setminus S_1$.
	Consider the fibered square
	\begin{displaymath}
		% https://q.uiver.app/#q=WzAsNCxbMywwLCJcXG9wZXJhdG9ybmFtZXtTcGVjfShcXGthcHBhKHApKSJdLFszLDIsIlMuIl0sWzAsMiwiWV8xIl0sWzAsMCwiWV5wXzEiXSxbMCwxLCJ0X3AiXSxbMiwxLCJmXzEiLDJdLFszLDAsImZecF8xIl0sWzMsMiwidF5wXzEiLDAseyJsYWJlbF9wb3NpdGlvbiI6NzB9XV0=
		\begin{tikzcd}
			{Y^p_1} &&& {\operatorname{Spec}(\kappa(p))} \\
			\\
			{Y_1} &&& {S.}
			\arrow["{f^p_1}", from=1-1, to=1-4]
			\arrow["{t^p_1}"{pos=0.7}, from=1-1, to=3-1]
			\arrow["{t_p}", from=1-4, to=3-4]
			\arrow["{f_1}"', from=3-1, to=3-4]
		\end{tikzcd}
	\end{displaymath}
	Since $p\not\in S_1$, $\mathbf{L}(t^\prime_p)^\ast K$ is not relatively perfect over $Y_1 \times_{S} \operatorname{Spec}(\kappa(p))$.
	Indeed, since $\mathbf{L}(t^p_1 \circ g^p_1)^\ast G_1 \otimes^{\mathbf{L}} \mathbf{L}(t^p_2 \circ g^p_2)^\ast G_2$ is a compact generator (see \Cref{lem:neeman2023cor5_10}), \Cref{thm:stacky_all_coincide} implies $\mathbf{L}(t^\prime_p)^\ast K$ is relatively perfect over $Y^p_1$ if, and only if,
	\begin{displaymath}
		\mathbf{R}(g^p_1)_\ast  \mathbf{L} (t^\prime_p)^\ast (K\otimes^{\mathbf{L}} \mathbf{L} (g_1)^\ast G_1 \otimes^{\mathbf{L}} \mathbf{L}(g_2)^\ast G_2)
	\end{displaymath}
	is perfect.
	By tor-independent base change \cite[\href{https://stacks.math.columbia.edu/tag/08IQ}{Tag 08IQ} and \href{https://stacks.math.columbia.edu/tag/08IR}{08IR}]{StacksProject},
	\begin{displaymath}
		\begin{aligned}
			\mathbf{R}(g^p_1)_\ast
			& \mathbf{L} (t^\prime_p)^\ast (K\otimes^{\mathbf{L}} \mathbf{L} (g_1)^\ast G_1 \otimes^{\mathbf{L}} \mathbf{L}(g_2)^\ast G_2)
			\\&\cong \mathbf{L}(t^p_1)^\ast \mathbf{R}(g_1)_\ast (K\otimes^{\mathbf{L}} \mathbf{L} (g_1)^\ast G_1 \otimes^{\mathbf{L}} \mathbf{L}(g_2)^\ast G_2)
			\\&\not\in \operatorname{Perf}(Y_1 \times_{S} \operatorname{Spec}(\kappa(p))).
		\end{aligned}
	\end{displaymath}
	By \Cref{prop:perfectness}, there exist  $q\in |Y_1 \times_S \operatorname{Spec}(\kappa(p))|$ and a representative $s\colon \operatorname{Spec}(\ell) \to Y_1 \times_S \operatorname{Spec}(\kappa(p))$ of $q$ such that
	\begin{displaymath}
		\mathbf{L}(t_1^p \circ s)^\ast \mathbf{R}(g_1)_\ast (K\otimes^{\mathbf{L}} \mathbf{L} (g_1)^\ast G_1 \otimes^{\mathbf{L}} \mathbf{L}(g_2)^\ast G_2) \not\in D^b_{\operatorname{coh}}(\ell).
	\end{displaymath}
	By \Cref{lem:factor_residual_gerbes}, $t^p_1 \circ s$ represents $|t^p_1| (q)$.
	Hence, $|t^p_1| (q) \in |Y_1|\setminus U$.
	Thus, $p\in |f_1| (|Y_1|\setminus U)$, and so $|S| \setminus S_1 \subseteq |f_1|(|Y_1|\setminus U)$.

	Conversely, let $p\in |f_1|(|Y_1|\setminus U)$.
	Choose $q\in |Y_1|\setminus U$ such that $|f_1| (q)=p$.
	Consider the fibered square
	\begin{displaymath}
		% https://q.uiver.app/#q=WzAsNCxbMywwLCJcXG9wZXJhdG9ybmFtZXtTcGVjfShcXGthcHBhKHApKSJdLFszLDIsIlMuIl0sWzAsMiwiWV8xIl0sWzAsMCwiWV5wXzEiXSxbMCwxLCJ0X3AiXSxbMiwxLCJmXzEiLDJdLFszLDAsImZecF8xIl0sWzMsMiwidF5wXzEiLDAseyJsYWJlbF9wb3NpdGlvbiI6NzB9XV0=
		\begin{tikzcd}
			{Y^p_1} &&& {\operatorname{Spec}(\kappa(p))} \\
			\\
			{Y_1} &&& {S.}
			\arrow["{f^p_1}", from=1-1, to=1-4]
			\arrow["{t^p_1}"{pos=0.7}, from=1-1, to=3-1]
			\arrow["{t_p}", from=1-4, to=3-4]
			\arrow["{f_1}"', from=3-1, to=3-4]
		\end{tikzcd}
	\end{displaymath}
	By \cite[\href{https://stacks.math.columbia.edu/tag/03H4}{Tag 03H4}]{StacksProject}, there exists $q^\prime \in |Y_1 \times_S \operatorname{Spec}(\kappa(p))|$ such that $|t^p_1| (q^\prime)=q$.
	Indeed, the morphism of the fiber product in category of sets says $(q,\{pt\})\in |Y_1| \times_S |\operatorname{Spec}(\kappa(p))|$ is a lift of $q$ along the natural morphism.
	Yet, \cite[\href{https://stacks.math.columbia.edu/tag/03H4}{Tag 03H4}]{StacksProject} gives us a lift of $(q,\{pt\})$ and clearly $|t^p_1|$ factors through the natural morphism on sets $|Y_1| \times_{|S|} |\operatorname{Spec}(\kappa(p))|\to |Y_1|$.
	Set $s\colon \operatorname{Spec}(\ell) \to Y_1 \times_S \operatorname{Spec}(\kappa(p))$ to be a representative of $q^\prime$.
	By \Cref{lem:factor_residual_gerbes}, $t^p_1 \circ s$ represents $q$.
	Since $q\not\in U$, tor-independent base change \cite[\href{https://stacks.math.columbia.edu/tag/08IQ}{Tag 08IQ} and \href{https://stacks.math.columbia.edu/tag/08IR}{08IR}]{StacksProject} implies
	\begin{displaymath}
		\begin{aligned}
			\mathbf{L}s^\ast \mathbf{R}(g^p_1)_\ast
			& \mathbf{L} (t^\prime_p)^\ast (K\otimes^{\mathbf{L}} \mathbf{L} (g_1)^\ast G_1 \otimes^{\mathbf{L}} \mathbf{L}(g_2)^\ast G_2)
			\\&\cong \mathbf{L}s^\ast \mathbf{L}(t^p_1)^\ast \mathbf{R}(g_1)_\ast (K\otimes^{\mathbf{L}} \mathbf{L} (g_1)^\ast G_1 \otimes^{\mathbf{L}} \mathbf{L}(g_2)^\ast G_2)
			\\&\not\in D^b_{\operatorname{coh}}(\ell).
		\end{aligned}
	\end{displaymath}
	Consequently, by \Cref{prop:perfectness},
	\begin{displaymath}
		\mathbf{R}(g^p_1)_\ast
		\mathbf{L} (t^\prime_p)^\ast (K\otimes^{\mathbf{L}} \mathbf{L} (g_1)^\ast G_1 \otimes^{\mathbf{L}} \mathbf{L}(g_2)^\ast G_2)
		\not\in \operatorname{Perf}(Y_1 \times_{S} \operatorname{Spec}(\kappa(p))).
	\end{displaymath}
	Thus, $p\not\in S_1$, and hence $|f_1|(|Y_1|\setminus U) \subseteq |S| \setminus S_1$.
\end{proof}

\begin{theorem}
	\label{thm:fm_locus_weak_approx}
	Consider proper flat morphisms $f_1\colon Y_1 \to S$ and $f_2\colon Y_2 \to S$ of Noetherian algebraic spaces.
	Let $E\in D_{\operatorname{coh}}^b(Y_1 \times_{S}  Y_2 )$.
	Then:
	\begin{enumerate}
		\item \label{cor:fm_locus1} $\operatorname{FM}(E)$ is an open subset of $|S|$
		\item \label{cor:fm_locus2} $\operatorname{FM}(E)$ commutes with arbitrary Noetherian base change.
	\end{enumerate}
\end{theorem}

\begin{proof}
	By \Cref{prop:relatively_perfect_locus}, the collection $U_1$ (resp.\ $U_2$) of $p\in |S|$ such that $\mathbf{L}(t^\prime_p)^\ast E$ is relatively perfect over $Y_1^p$ (resp.\ $Y_2^p$) forms an open subset.
	Set $U:=U_1\cap U_2$.
	Denote by $j\colon U\to S$ the associated open immersion.

	If $p\in \operatorname{FM}(E)$, then by \Cref{rem:rel_perfect_is_free} we have $E_p$ is relatively perfect over each $Y_i^p$, and so $p\in U$.
	Therefore, $\operatorname{FM}(E)\subseteq U$.

	There exists a fibered cube
	\begin{displaymath}
		% https://q.uiver.app/#q=WzAsOCxbNCwxLCJVIl0sWzQsMywiUy4iXSxbMSwzLCJZXzEiXSxbMywyLCJZXzIiXSxbMSwxLCJZXzEgXFx0aW1lc19TIFUiXSxbMywwLCJZXzIgXFx0aW1lc19TIFUiXSxbMCwwLCJZXzEgXFx0aW1lc19TIFlfMlxcdGltZXNfUyBVIl0sWzAsMiwiWV8xXFx0aW1lc19TIFlfMiJdLFswLDEsImoiXSxbMiwxLCJmXzEiLDJdLFszLDEsImZfMiJdLFs0LDAsImZeVV8xIl0sWzUsMCwiZl5VXzIiXSxbNCwyLCJ0XlVfMSIsMCx7ImxhYmVsX3Bvc2l0aW9uIjo3MH1dLFs1LDMsInRfM15VIiwyLHsibGFiZWxfcG9zaXRpb24iOjgwLCJzdHlsZSI6eyJib2R5Ijp7Im5hbWUiOiJkYXNoZWQifX19XSxbNiw0LCJnXzFeVSIsMl0sWzcsMywiZ18yIiwwLHsibGFiZWxfcG9zaXRpb24iOjEwLCJzdHlsZSI6eyJib2R5Ijp7Im5hbWUiOiJkYXNoZWQifX19XSxbNywyLCJnXzEiLDJdLFs2LDcsInReXFxwcmltZV9VIiwyXSxbNiw1LCJnXlVfMiIsMl1d
		\begin{tikzcd}
			{Y_1 \times_S Y_2\times_S U} &&& {Y_2 \times_S U} & \\
			& {Y_1 \times_S U} &&& U \\
			{Y_1\times_S Y_2} &&& {Y_2} \\
			& {Y_1} &&& {S.}
			\arrow["{g^U_2}"', from=1-1, to=1-4]
			\arrow["{g_1^U}"', from=1-1, to=2-2]
			\arrow["{t^\prime_U}"', from=1-1, to=3-1]
			\arrow["{f^U_2}", from=1-4, to=2-5]
			\arrow["{t_2^U}"'{pos=0.8}, dashed, from=1-4, to=3-4]
			\arrow["{f^U_1}", from=2-2, to=2-5]
			\arrow["{t^U_1}"{pos=0.7}, from=2-2, to=4-2]
			\arrow["j", from=2-5, to=4-5]
			\arrow["{g_2}"{pos=0.1}, dashed, from=3-1, to=3-4]
			\arrow["{g_1}"', from=3-1, to=4-2]
			\arrow["{f_2}", from=3-4, to=4-5]
			\arrow["{f_1}"', from=4-2, to=4-5]
		\end{tikzcd}
	\end{displaymath}
	Set $E_j := \mathbf{L}(t^\prime_U)^\ast E$.
	By \Cref{thm:bounded_pseudocoherence_perfectness_faithfully_flat_affine}, $E_j$ is relatively perfect over each $Y_i\times_S U$.
	Applying \Cref{thm:openness_for_equivalence_fullfaithful_dbcoh}, $\operatorname{FM}(E_j)$ is open in $U$.

	We claim that $\operatorname{FM}(E)=|j|(\operatorname{FM}(E_j))$.
	Let $p\in \operatorname{FM}(E_j)$.
	Choose a representative $a\colon \operatorname{Spec}(k)\to U$ of $p$ such that the corresponding fiber transform is an equivalence.
	Composing $a$ with $j$, pseudofunctoriality of derived pullback identifies the resulting fiber kernel with the corresponding fiber of $E$.
	Hence, $|j|(p)\in \operatorname{FM}(E)$.

	Conversely, let $p\in \operatorname{FM}(E)$.
	We have already shown that $p\in U$.
	Choose a representative $a\colon \operatorname{Spec}(k)\to S$ of $p$ such that the corresponding fiber transform is an equivalence.
	Since the image of $a$ is contained in $U$, the morphism $a$ factors through $j$.
	Again, pseudofunctoriality identifies the resulting fiber kernel of $E_j$ with the fiber kernel of $E$.
	Thus, $p\in |j|(\operatorname{FM}(E_j))$.
	This proves the claim.

	If $j^\prime\colon \operatorname{FM}(E_j)\to U$ is the associated open immersion, then $\operatorname{FM}(E)$ is the topological image of the open immersion $j\circ j^\prime$.
	Therefore, $\operatorname{FM}(E)$ is open in $|S|$, proving \eqref{cor:fm_locus1}.

	We prove \eqref{cor:fm_locus2}.
	Let $t\colon T\to S$ be a morphism of Noetherian algebraic spaces.
	Set $V:=T\times_S U$.
	Denote by $t_U\colon V\to U$ and $j_T\colon V\to T$ the natural morphisms.
	Denote by $E_t$ and $(E_j)_{t_U}$ the derived pullbacks of $E$ and $E_j$ to $Y_1\times_S Y_2\times_S T$ and $Y_1\times_S Y_2\times_S V$ respectively, and by $j_T^\prime\colon Y_1\times_S Y_2\times_S V \to Y_1\times_S Y_2\times_S T$ the base change of $j_T$.
	By pseudofunctoriality, there is an isomorphism $\mathbf{L}(j_T^\prime)^\ast E_t \cong (E_j)_{t_U}$.
	Since $E_j$ is relatively perfect over each $Y_i\times_S U$, \Cref{cor:deform} gives $t_U^{-1}(\operatorname{FM}(E_j)) = \operatorname{FM}((E_j)_{t_U})$.

	Choose $q\in |T|$.
	Set $p:=t(q)$.
	First, assume that $p\in \operatorname{FM}(E)$.
	Since $\operatorname{FM}(E)\subseteq U$, we have $q\in V$.
	By the fact that $\operatorname{FM}(E)=|j|(\operatorname{FM}(E_j))$, the point $p\in U$ belongs to $\operatorname{FM}(E_j)$.
	Hence, \Cref{cor:deform} implies $q\in \operatorname{FM}((E_j)_{t_U})$.
	By the isomorphism above, this gives $q\in \operatorname{FM}(E_t)$.
	Hence, $t^{-1}(\operatorname{FM}(E)) \subseteq \operatorname{FM}(E_t)$.

	Conversely, suppose $q\in \operatorname{FM}(E_t)$.
	Denote by $a\colon \operatorname{Spec}(\kappa(q))\to \operatorname{Spec}(\kappa(p))$ the canonical morphism induced by $t$.
	The fiber of $E_t$ at $q$ is the base change of $E_p$ along $a$.
	Since $q\in \operatorname{FM}(E_t)$, \Cref{rem:rel_perfect_is_free} implies that this fiber kernel is relatively perfect over both factors.
	As $a$ is faithfully flat, \Cref{prop:reflecting_bounded_pseudocoherence_perfectness1} implies that $E_p$ is relatively perfect over each $Y_i^p$.
	Hence, $p\in U$, and therefore $q\in V$.

	By the isomorphism $\mathbf{L}(j_T^\prime)^\ast E_t \cong (E_j)_{t_U}$,
	$q\in \operatorname{FM}((E_j)_{t_U})$.
	Applying \Cref{cor:deform}, the point $p$, considered as a point of $U$, belongs to $\operatorname{FM}(E_j)$.
	Hence, $p\in \operatorname{FM}(E)$.
	Consequently, $\operatorname{FM}(E_t) \subseteq t^{-1}(\operatorname{FM}(E))$.
	Therefore, $t^{-1}(\operatorname{FM}(E)) = \operatorname{FM}(E_t)$, which proves \eqref{cor:fm_locus2}.
\end{proof}

\begin{example}
	\label{ex:fm_locus_counterexample}
	Let $S=\operatorname{Spec}(\mathbb{Z})$.
	Set $Y_1=Y_2=S$ with $f_1$ and $f_2$ the identity.
	Choose an infinite sequence of prime numbers $p_1 < p_2 < \cdots$.
	For each $n\geq 1$, denote by $t_n\colon \operatorname{Spec}(\mathbb{F}_{p_n}) \to S$ the associated closed immersion.
	Consider
	\begin{displaymath}
		K := \mathcal{O}_S \oplus \bigoplus_{n\geq 1} (t_n)_\ast \mathcal{O}_{\operatorname{Spec}(\mathbb{F}_{p_n})} [n] \in D^-_{\operatorname{coh}}(S).
	\end{displaymath}
	Note that $\Phi_K (\mathcal{O}_S)\cong K$ is not bounded below, and so $\Phi_K$ does not preserve $D^b_{\operatorname{coh}}$.
	Let $t\colon \operatorname{Spec}(k)\to S$ be a representative of a point $p\in |S|$.
	If $p\not\in \{p_1,p_2,\dots\}$, then $\mathbf{L}t^\ast K\cong \mathcal{O}_{\operatorname{Spec}(k)}$, and so the fiber of $\Phi_K$ at $t$ is the identity.
	If $p=p_n$, then $\mathbf{L}t^\ast K\cong \mathcal{O}_{\operatorname{Spec}(k)} \oplus \mathcal{O}_{\operatorname{Spec}(k)}[n] \oplus \mathcal{O}_{\operatorname{Spec}(k)}[n+1]$, and so the fiber of $\Phi_K$ at $t$ is the endofunctor $E\mapsto E\oplus E[n] \oplus E[n+1]$.
	This endofunctor is not fully faithful (e.g.\ compare the dimensions of the endomorphism algebras of $\mathcal{O}_{\operatorname{Spec}(k)}$ and of its image).
	Thus, $\operatorname{FM}(K)=\operatorname{fm}(K)=|S|\setminus \{p_1,p_2,\dots\}$, which is not open.
\end{example}

%%%%%%%%%%%%%%%%%%%%%%%%%%%%%%%%%%%
\part{Appendices}
\label{sec:part4}
%%%%%%%%%%%%%%%%%%%%%%%%%%%%%%%%%%%

%%%%%%%%%%%%%%%%%%%%%%%%%%%%%%%%%%%
\appendix
%%%%%%%%%%%%%%%%%%%%%%%%%%%%%%%%%%%

%%%%%%%%%%%%%%%%%%%%%%%%%%%%%%%%%%%
\section{Preferred equivalence classes}
\label{app:preferred}
%%%%%%%%%%%%%%%%%%%%%%%%%%%%%%%%%%%

We prove versions of \cite[Proposition 6.6]{DeDeyn/Lank/ManaliRahul/Peng:2025}, \cite[Proposition 3.1]{Lank:2026a}, and \cite[Proposition 4.2]{Hall/Rydh:2023} for algebraic spaces.
The difference is that loc.\ cit.\ takes place on the lisse-\'{e}tale site, and so we prove it on the small \'{e}tale site for the sake of completeness.
This requires a placeholder:
\begin{equation}
	\begin{minipage}{13cm}
		An algebraic space $X$ satisfies \textbf{P.E.C.} if for any quasi-compact open immersion $U \to X$ with complement $Z$, the standard $t$-structure on $D_{\operatorname{qc},Z}(X)$ is in the preferred equivalence class.
	\end{minipage}
\end{equation}

\begin{lemma}
	\label{lem:BNP_weak_version_for_glueing}
	Consider a recollement of triangulated categories
	\begin{displaymath}
		\begin{tikzcd}[ampersand replacement=\&]
			{\mathcal{D}_Z} \&\& {\mathcal{D}} \&\& {\mathcal{D}_U}\rlap{ .}
			\arrow["{i_\ast}" description, from=1-1, to=1-3]
			\arrow["{i^!}", bend left=30, from=1-3, to=1-1]
			\arrow["{i^\ast}"', bend right=30, from=1-3, to=1-1]
			\arrow["{j^\ast}" description, from=1-3, to=1-5]
			\arrow["{j_!}"', bend right=30, from=1-5, to=1-3]
			\arrow["{j_\ast}", bend left=30, from=1-5, to=1-3]
		\end{tikzcd}
	\end{displaymath}
	Let $G_U$ be a compact generator $\mathcal{D}_U$. Assume that $\mathcal{D}$ admits a compact generator $G^\prime$ such that $\operatorname{Hom}(G^\prime [-n],G^\prime)=0$ for $n\gg 0$ and $j^\ast G^\prime\in \mathcal{D}_U^{\leq m}$ for some $m>0$ (here, $\mathcal{D}_U^{\leq 0}$ is the aisle of the $t$-structure on $\mathcal{D}_U$ compactly generated by $G_U$).
	There exists a compact generator $G\in \mathcal{D}$, which can be taken to be $G^\prime[m] \oplus j_! G_U$, such that the $t$-structure on $\mathcal{D}$ compactly generated by $G$ is the gluing of the $t$-structure on $\mathcal{D}_Z$ compactly generated by $i^\ast G$ and the $t$-structure on $\mathcal{D}_U$ compactly generated by $G_U$.
	In particular, the glued $t$-structure is in the preferred equivalence class.
\end{lemma}

\begin{proof}
	We follow the argument of \cite[Lemma 3.11]{Burke/Neeman/Pauwels:2023}.
	Set $G:=G^\prime[m]\oplus j_!G_U$.
	We first check that $G$ is a compact generator of $\mathcal{D}$.
	The object $G^\prime[m]$ is a compact generator.
	Moreover, $j_!G_U$ is compact.
	Indeed, the right adjoint of $j_!$ is $j^\ast$, and $j^\ast$ preserves small coproducts since it is also a left adjoint in the recollement.
	Hence, $G$ is compact, and it generates $\mathcal{D}$ because $G^\prime[m]$ is a direct summand of $G$.

	The object $i^\ast G$ is a compact generator of $\mathcal{D}_Z$.
	Indeed, $i_\ast$ preserves small coproducts since it admits both a left and a right adjoint, so $i^\ast$ preserves compact objects.
	If $A\in\mathcal{D}_Z$ satisfies $\operatorname{Hom}(i^\ast G,A[n])=0$ for all $n\in\mathbb{Z}$, then adjunction gives $\operatorname{Hom}(G,i_\ast A[n])=0$ for all $n\in\mathbb{Z}$.
	Since $G$ generates $\mathcal{D}$, this implies $i_\ast A\cong0$, and hence $A\cong0$.

	Let $(\mathcal{D}_Z^{\leq0},\mathcal{D}_Z^{\geq0})$ and $(\mathcal{D}_U^{\leq0},\mathcal{D}_U^{\geq0})$ respectively be the $t$-structures compactly generated by $i^\ast G$ and $G_U$.
	Define $(\mathcal{D}_{\operatorname{gl}}^{\leq0}, \mathcal{D}_{\operatorname{gl}}^{\geq0})$ to be their glued $t$-structure on $\mathcal{D}$.
	Then
	\begin{displaymath}
		\mathcal{D}_{\operatorname{gl}}^{\leq0}
		= \left\{ A\in\mathcal{D} \ : i^\ast A\in\mathcal{D}_Z^{\leq0}
		\textrm{ and } j^\ast A\in\mathcal{D}_U^{\leq0}
		\right\}.
	\end{displaymath}
	We show that this is the aisle compactly generated by $G$.

	By hypothesis, $j^\ast G^\prime\in\mathcal{D}_U^{\leq m}$.
	Hence, $j^\ast G^\prime[m]\in\mathcal{D}_U^{\leq0}$.
	Since $j^\ast j_!\cong1$, we obtain
	\begin{displaymath}
		j^\ast G
		\cong j^\ast G^\prime[m]\oplus G_U
		\in\mathcal{D}_U^{\leq0}.
	\end{displaymath}
	By definition, $i^\ast G\in\mathcal{D}_Z^{\leq0}$.
	Since $i^\ast$ and $j^\ast$ are exact and preserve small coproducts, they respectively carry the cocomplete aisle generated by $G$ into $\mathcal{D}_Z^{\leq0}$ and $\mathcal{D}_U^{\leq0}$.
	Therefore, $\overline{\langle G\rangle}^{(-\infty,0]} \subseteq \mathcal{D}_{\operatorname{gl}}^{\leq0}$.

	We prove the reverse inclusion.
	Since $j_!G_U$ is a direct summand of $G$, $j_!G_U \in \overline{\langle G\rangle}^{(-\infty,0]}$.
	The functor $j_!$ is exact and preserves small coproducts.
	Since $\mathcal{D}_U^{\leq0}$ is the cocomplete aisle generated by $G_U$, it follows that $j_!\mathcal{D}_U^{\leq0} \subseteq \overline{\langle G\rangle}^{(-\infty,0]}$.
	Consider the recollement triangle
	\begin{displaymath}
		j_!j^\ast G
		\to G
		\to i_\ast i^\ast G
		\to j_!j^\ast G[1].
	\end{displaymath}
	Since $j^\ast G\in\mathcal{D}_U^{\leq0}$, the preceding inclusion gives $j_!j^\ast G \in \overline{\langle G\rangle}^{(-\infty,0]}$.
	As an aisle is closed under extensions and positive shifts, the
	triangle implies $i_\ast i^\ast G \in \overline{\langle G\rangle}^{(-\infty,0]}$.
	The functor $i_\ast$ is exact and preserves small coproducts.
	Since $\mathcal{D}_Z^{\leq0}$ is the cocomplete aisle generated by $i^\ast G$, it follows that $i_\ast\mathcal{D}_Z^{\leq0} \subseteq \overline{\langle G\rangle}^{(-\infty,0]}$.

	Finally, let $A\in\mathcal{D}_{\operatorname{gl}}^{\leq0}$.
	Then $j^\ast A\in\mathcal{D}_U^{\leq0}$ and $i^\ast A\in\mathcal{D}_Z^{\leq0}$.
	Thus, both $j_!j^\ast A$ and $i_\ast i^\ast A$ belong to $\overline{\langle G\rangle}^{(-\infty,0]}$.
	The recollement triangle
	\begin{displaymath}
		j_!j^\ast A
		\to A
		\to i_\ast i^\ast A
		\to j_!j^\ast A[1]
	\end{displaymath}
	implies $A\in \overline{\langle G\rangle}^{(-\infty,0]}$.
	Consequently, $\mathcal{D}_{\operatorname{gl}}^{\leq0}  = \overline{\langle G\rangle}^{(-\infty,0]}$.
	Hence, the $t$-structure compactly generated by $G=G^\prime[m]\oplus j_!G_U$ is precisely the glued $t$-structure.
	Since $G$ is a compact generator of $\mathcal{D}$, the glued $t$-structure belongs to the preferred equivalence class.
\end{proof}

\begin{lemma}
	\label{lem:support_for_pullback}
	Let $f\colon Y\to X$ be a morphism of quasi-compact quasi-separated algebraic spaces.
	Consider a quasi-compact open immersion $j\colon U\to X$.
	Set $Z:= |X|\setminus |U|$.
	Then
	\begin{displaymath}
		\mathbf{L}f^\ast (D_{\operatorname{qc},Z}(X))\subseteq D_{\operatorname{qc},f^{-1}(Z)}(Y).
	\end{displaymath}
\end{lemma}

\begin{proof}
	Choose $E\in D_{\operatorname{qc},Z}(X)$.
	Set $f^\prime\colon f^{-1}(U) \to U$ to be the base change of $j$ along $f$.
	As $\mathbf{L}j^\ast E\cong 0$, one has $\mathbf{L}(j\circ f^\prime)^\ast E\cong0$.
	Define $j^\prime\colon f^{-1}(U) \to Y$ to be the open immersion obtained by base change.
	Hence, $\mathbf{L}(f\circ j^\prime)^\ast E\cong0$, implying that $\mathbf{L}f^\ast E \in D_{\operatorname{qc},f^{-1}(Z)}(Y)$.
\end{proof}

\begin{lemma}
	\label{lem:support_for_pushforward}
	Let $f\colon Y\to X$ be a morphism of quasi-compact quasi-separated algebraic spaces.
	Consider a quasi-compact open immersion $j\colon U\to X$.
	Set $Z:= |X|\setminus |U|$. Then
	\begin{displaymath}
		\mathbf{R}f_\ast (D_{\operatorname{qc},f^{-1}(Z)}(Y))\subseteq D_{\operatorname{qc},Z}(X).
	\end{displaymath}
\end{lemma}

\begin{proof}
	By \cite[\href{https://stacks.math.columbia.edu/tag/0AEC}{Tag 0AEC}]{StacksProject}, $D_{\operatorname{qc},f^{-1}(Z)} (Y)$ is compactly generated by a single object $G$.
	Recall that $\mathbf{R}f_\ast$ preserves small coproducts.
	Applying \cite[\href{https://stacks.math.columbia.edu/tag/09SN}{Tag 09SN}]{StacksProject}, we can reduce to checking that $\mathbf{R}f_\ast G \in D_{\operatorname{qc},Z}(X)$.
	As $f^{-1} (Z) = |Y|\setminus |f^{-1}(U)|$, it follows that $\mathbf{L}(j^\prime)^\ast G\cong 0$ where $j^\prime\colon Y\times_X U\to Y$ is the projection from the fiber product obtained base changing $f$ along $j$.
	Denote $f^\prime\colon Y\times_X U\to U$ for the other projection.
	By flat base change, it follows that
	\begin{displaymath}
		\mathbf{L}j^\ast \mathbf{R}f_\ast G \cong \mathbf{R}f^\prime_\ast \mathbf{L}(j^\prime)^\ast G \cong 0.
	\end{displaymath}
	This completes the proof.
\end{proof}

\begin{lemma}
	\label{lem:spacey_recollement}
	Let $X$ be a quasi-compact quasi-separated algebraic space.
	Let $j\colon U\to X$ be a quasi-compact open immersion.
	Let $Z:=|X|\setminus |U|$.
	There exists a recollement
	\begin{displaymath}
		% https://q.uiver.app/#q=WzAsMyxbNCwwLCJEX3tcXG9wZXJhdG9ybmFtZXtxY30sWn0oXFxtYXRoY2Fse1h9KSJdLFsyLDAsIkRfe1xcb3BlcmF0b3JuYW1le3FjfX0oXFxtYXRoY2Fse1h9KSJdLFswLDAsIkRfe1xcb3BlcmF0b3JuYW1le3FjfX0oXFxtYXRoY2Fse1V9KSJdLFsxLDAsImleISIsMV0sWzIsMSwiXFxtYXRoYmZ7Un1qX1xcYXN0IiwxXSxbMCwxLCJpX1xcYXN0IiwyLHsib2Zmc2V0IjozLCJjdXJ2ZSI6Mn1dLFswLDEsImleXFxhc3QiLDAseyJvZmZzZXQiOi0zLCJjdXJ2ZSI6LTJ9XSxbMSwyLCJcXG1hdGhiZntMfWpeXFxhc3QiLDIseyJvZmZzZXQiOjMsImN1cnZlIjoyfV0sWzEsMiwial5cXHRpbWVzIiwwLHsib2Zmc2V0IjotMywiY3VydmUiOi0yfV1d
		\begin{tikzcd}
			{D_{\operatorname{qc}}(U)} && {D_{\operatorname{qc}}(X)} && {D_{\operatorname{qc},Z}(X)}
			\arrow["{\mathbf{R}j_\ast}"{description}, from=1-1, to=1-3]
			\arrow["{\mathbf{L}j^\ast}"', shift right=3, bend right=30, from=1-3, to=1-1]
			\arrow["{j^\times}", shift left=3, bend left=30, from=1-3, to=1-1]
			\arrow["{i^!}"{description}, from=1-3, to=1-5]
			\arrow["{i_\ast}"', shift right=3, bend right=30, from=1-5, to=1-3]
			\arrow["{i^\ast}", shift left=3, bend left=30, from=1-5, to=1-3]
		\end{tikzcd}
	\end{displaymath}
	where $i_\ast$ is the natural inclusion and $j^\times$ the right adjoint of $\mathbf{R}j_\ast$.
	Additionally, given a closed $W\subseteq |X|$ subset with quasi-compact complement, there exists a recollement
	\begin{displaymath}
		% https://q.uiver.app/#q=WzAsMyxbNCwwLCJEX3tcXG9wZXJhdG9ybmFtZXtxY30sV1xcY2FwIFp9KFxcbWF0aGNhbHtYfSkiXSxbMiwwLCJEX3tcXG9wZXJhdG9ybmFtZXtxY30sV30oXFxtYXRoY2Fse1h9KSJdLFswLDAsIkRfe1xcb3BlcmF0b3JuYW1le3FjfSx8XFxtYXRoY2Fse1V9fFxcY2FwIFd9KFxcbWF0aGNhbHtVfSkiXSxbMSwwLCJpXiEiLDFdLFsyLDEsIlxcbWF0aGJme1J9al9cXGFzdCIsMV0sWzAsMSwiaV9cXGFzdCIsMix7Im9mZnNldCI6MywiY3VydmUiOjJ9XSxbMCwxLCJpXlxcYXN0IiwwLHsib2Zmc2V0IjotMywiY3VydmUiOi0yfV0sWzEsMiwiXFxtYXRoYmZ7TH1qXlxcYXN0IiwyLHsib2Zmc2V0IjozLCJjdXJ2ZSI6Mn1dLFsxLDIsImpeISIsMCx7Im9mZnNldCI6LTMsImN1cnZlIjotMn1dXQ==
		\begin{tikzcd}
			{D_{\operatorname{qc},|U|\cap W}(U)} && {D_{\operatorname{qc},W}(X)} && {D_{\operatorname{qc},W\cap Z}(X)}
			\arrow["{\mathbf{R}j_\ast}"{description}, from=1-1, to=1-3]
			\arrow["{\mathbf{L}j^\ast}"', shift right=3, bend right=30, from=1-3, to=1-1]
			\arrow["{j^!_W}", shift left=3, bend left=30, from=1-3, to=1-1]
			\arrow["{i^!}"{description}, from=1-3, to=1-5]
			\arrow["{i_\ast}"', shift right=3, bend right=30, from=1-5, to=1-3]
			\arrow["{i^\ast}", shift left=3, bend left=30, from=1-5, to=1-3]
		\end{tikzcd}
	\end{displaymath}
	where $i_\ast$ is the natural inclusion.
\end{lemma}

\begin{proof}
	By \cite[\href{https://stacks.math.columbia.edu/tag/0AEC}{Tag 0AEC}]{StacksProject}, $D_{\operatorname{qc},|U|\cap W}(U)$ is compactly generated.
	We show the second claim because the first claim is similar (in the absolute case $W=|X|$, note that one obtains $j^!$ is naturally isomorphic to $j^\times$ by uniqueness of adjoints).
	By \Cref{lem:support_for_pullback}, the restriction of $\mathbf{L}j^\ast D_{\operatorname{qc},W}(X) \subseteq D_{\operatorname{qc},|U|\cap W}(U)$.
	Moreover, \Cref{lem:support_for_pushforward} says $\mathbf{R}j_\ast D_{\operatorname{qc},|U|\cap W}(U) \subseteq D_{\operatorname{qc},W}(X)$.
	Hence, the adjoint pair $\mathbf{L}j^\ast$ and $\mathbf{R}j_\ast$ restricts to $D_{\operatorname{qc},W}(X)$ and $D_{\operatorname{qc},|U|\cap W}(U)$.
	There exists a Verdier localization sequence,
	\begin{displaymath}
		D_{\operatorname{qc},W\cap Z}(X) \xrightarrow{i_\ast} D_{\operatorname{qc},W}(X) \xrightarrow{\mathbf{L}j^\ast} D_{\operatorname{qc},|U|\cap W}(U).
	\end{displaymath}
	To see, use that the counit of $\mathbf{L}j^\ast$ and $\mathbf{R}j_\ast$ is a natural isomorphism on $D_{\operatorname{qc}}$, which implies that $\mathbf{L}j^\ast\colon D_{\operatorname{qc},W}(X) \to D_{\operatorname{qc},|U|\cap W}(U)$ is essentially surjective.
	To see that full faithfulness holds, use \cite[\S I.\ Proposition 1.3]{Gabriel/Zisman:1967} and \cite[Lemma 4.3.1]{Krause:2010}.
	As $\mathbf{R}j_\ast$ is right adjoint to $\mathbf{L}j^\ast$, \cite[Lemma 2.2(ii)]{Gao/Psaroudakis:2018} implies $i_\ast$ admits a right adjoint. We denote it by $i^!$ (see also \cite[Theorem 1.1]{Cline/Parshall/Scott:1988a} and \cite[Theorem 2.1]{Cline/Parshall/Scott:1988b}).
	Hence, we obtain another Verdier localization sequence,
	\begin{displaymath}
		D_{\operatorname{qc},Z\cap W}(X) \xleftarrow{i^!} D_{\operatorname{qc},W}(X) \xleftarrow{\mathbf{R}j_\ast} D_{\operatorname{qc},|U|\cap W}(U).
	\end{displaymath}
	Since $D_{\operatorname{qc},|U|\cap W}(U)$ is compactly generated and $\mathbf{R}j_\ast$ preserves small coproducts, \cite[Corollary 2.3]{Balmer/Dell'Ambrogio/Sanders:2016} implies the restriction of $\mathbf{R}j_\ast$ admits a right adjoint.
	Applying \cite[Lemma 2.2(ii)]{Gao/Psaroudakis:2018}, $i^!$ admits a right adjoint as well, which we denote by $i^\ast$.
	Thus, we have the required data for a recollement; i.e.\ a Verdier localization sequence which is a localization and colocalization sequence in the sense of \cite[\S 4]{Krause:2010}.
\end{proof}

\begin{remark}
	\label{rmk:stacky_recollement_compacts_preserved_by_ishriek}
	In \Cref{lem:spacey_recollement}, $i^!$ preserves small coproducts, and so, $i_\ast$ preserves compact objects (see e.g.\ the proof of $\implies$ in \cite[Theorem 5.1]{Neeman:1996}; which this fact does not require compact generation).
\end{remark}

\begin{lemma}
	\label{lem:etale_nbhd_preferred_class_equivalences}
	Consider an elementary distinguished square $j\colon U \to X$ and $f\colon V \to X$ where $j$ is an open immersion of quasi-compact quasi-separated algebraic spaces.
	%%NOTE: So $j$ is quasi-compact
	Set $Z:= |X|\setminus |U|$.
	Then the restrictions of $\mathbf{L}f^\ast$ and $\mathbf{R}f_\ast$ yield $t$-exact equivalences:
	\begin{enumerate}
		\item $D_{\operatorname{qc},f^{-1}(Z)}(V)$ and $D_{\operatorname{qc},Z}(X)$
		\item $D_{\operatorname{qc},f^{-1}(Z)}^b(V)$ and $D_{\operatorname{qc},Z}^b(X)$.
	\end{enumerate}
\end{lemma}

\begin{proof}
	Denote by $\epsilon$ and $\eta$ respectively the counit and unit of the adjoint pair $\mathbf{L} f^\ast$ and $\mathbf{R}f_\ast$.
	Applying \cite[\href{https://stacks.math.columbia.edu/tag/08GG}{Tag 08GG}]{StacksProject}, $\epsilon$ and $\eta$ are isomorphisms when restricted to objects with the prescribed support.
	It follows that $\mathbf{L} f^\ast$ and $\mathbf{R}f_\ast$ restrict to an equivalence $D_{\operatorname{qc},f^{-1}(Z)}(V)$ and $D_{\operatorname{qc},Z}(X)$.
	Since $f$ is \'{e}tale, $\mathbf{L} f^\ast$ is $t$-exact, whereas the restricted adjoint pair being an equivalence implies the restriction of $\mathbf{R} f_\ast$ is $t$-exact.
	By $t$-exactness, the restricted functors preserve bounded cohomology, and so the second claim follows.
\end{proof}

\begin{lemma}
	\label{lem:pullback_compact_generator}
	Let $f\colon Y\to X$ be a quasi-affine morphism of quasi-compact quasi-separated algebraic spaces.
	Fix a closed subset $W\subseteq |X|$ with quasi-compact complement.
	Then $\mathbf{L}f^\ast G$ compactly generates $D_{\operatorname{qc},f^{-1}(W)}(Y)$ whenever $G\in D_{\operatorname{qc},W}(X)$ does such for $D_{\operatorname{qc},W}(X)$.
\end{lemma}

\begin{proof}
	By \Cref{lem:quasi-affine_is_conservative,lem:support_for_pushforward}, the restriction of $\mathbf{R}f_\ast$ is conservative.
	Then apply \cite[Corollary 3.6]{Lank:2026a}.
\end{proof}

\begin{lemma}
	\label{lem:etale_nbhd_preferred_class}
	Consider an elementary distinguished square $j\colon U \to X$ and $f\colon Y \to X$ of quasi-compact quasi-separated algebraic spaces.
	If $U$ and $Y$ satisfy P.E.C., then so does $X$.
\end{lemma}

\begin{proof}
	Define $Z:=|X|\setminus|U|$. Let $v\colon V\to X$ be a
	quasi-compact open immersion and set $W:=|X|\setminus|V|$.
	It suffices to prove that the standard $t$-structure on
	$D_{\operatorname{qc},W}(X)$ belongs to the preferred
	equivalence class.
	By \Cref{lem:spacey_recollement}, there is a recollement
	\begin{displaymath}
		\begin{tikzcd}
			{D_{\operatorname{qc},|U|\cap W}(U)}
			&&
			{D_{\operatorname{qc},W}(X)}
			&&
			{D_{\operatorname{qc},Z\cap W}(X)}
			\arrow["{\mathbf{R}j_\ast}"{description}, from=1-1,to=1-3]
			\arrow["{\mathbf{L}j^\ast}"', bend right=30, from=1-3,to=1-1]
			\arrow["{j^!_W}", bend left=30, from=1-3,to=1-1]
			\arrow["{i^!}"{description}, from=1-3,to=1-5]
			\arrow["{i_\ast}"', bend right=30, from=1-5,to=1-3]
			\arrow["{i^\ast}", bend left=30, from=1-5,to=1-3]
		\end{tikzcd}
	\end{displaymath}
	where $i_\ast$ is the natural inclusion.

	We first record a bound for $i^!$ with respect to the standard
	$t$-structures. There exists $a\geq0$ such that
	\begin{equation}
		\label{eq:i_shriek_right_amplitude}
		i^! D_{\operatorname{qc},W}^{\leq0}(X)
		\subseteq D_{\operatorname{qc},Z\cap W}^{\leq a+1}(X).
	\end{equation}
	Indeed, if $E\in D_{\operatorname{qc},W}^{\leq0}(X)$, then flatness of $j$
	gives $\mathbf{L}j^\ast E \in D_{\operatorname{qc},|U|\cap W}^{\leq0}(U)$.
	Moreover, there exists $a\geq0$ which is independent of $E$ such that $\mathbf{R}j_\ast D_{\operatorname{qc},|U|\cap W}^{\leq0}(U) \subseteq D_{\operatorname{qc},W}^{\leq a}(X)$.
	This can be checked after base change along an \'{e}tale presentation, using that open immersions are representable by schemes and \cite[Proposition 3.9.2]{Lipman:2009}.
	The recollement triangle
	\begin{equation}
		\label{eq:rec}
		i_\ast i^!E
		\to E
		\to \mathbf{R}j_\ast\mathbf{L}j^\ast E
		\to i_\ast i^!E[1]
	\end{equation}
	then shows that $i_\ast i^!E[1]\in D_{\operatorname{qc},W}^{\leq a}(X)$.
	Since $i_\ast$ is the natural inclusion, \eqref{eq:i_shriek_right_amplitude} follows.

	Choose a compact generator $P\in D_{\operatorname{qc},W}(X)\cap\operatorname{Perf}(X)$.
	Since $j$ is quasi-affine, \Cref{lem:pullback_compact_generator} shows that
	$\mathbf{L}j^\ast P$ compactly generates $D_{\operatorname{qc},|U|\cap W}(U)$.
	Choose a compact generator $P_Z\in D_{\operatorname{qc},Z\cap W}(X)$.
	Since $i^!$ preserves small coproducts, $i_\ast$ preserves compact objects, so $i_\ast P_Z\in D_{\operatorname{qc},W}(X)\cap\operatorname{Perf}(X)$.

	We next compare the standard $t$-structure on $D_{\operatorname{qc},Z\cap W}(X)$ with the $t$-structure compactly generated by $P_Z$.
	Since $f\colon Y\to X$ occurs in an elementary distinguished square, its restriction over $Z$ is an isomorphism.
	Hence \Cref{lem:etale_nbhd_preferred_class_equivalences}, applied to $Z\cap W$, gives a $t$-exact equivalence $\mathbf{L}f^\ast\colon D_{\operatorname{qc},Z\cap W}(X) \to D_{\operatorname{qc},f^{-1}(Z\cap W)}(Y)$.
	Under this equivalence, $\mathbf{L}f^\ast P_Z$ is a compact generator.
	Since $Y$ satisfies P.E.C., the standard $t$-structure on $D_{\operatorname{qc},f^{-1}(Z\cap W)}(Y)$ belongs to the preferred equivalence class.
	Transporting along the $t$-exact equivalence shows that the standard $t$-structure on $D_{\operatorname{qc},Z\cap W}(X)$ is equivalent to the
	$t$-structure compactly generated by $P_Z$.

	In particular, there exists $c_Z\geq0$ such that
	\begin{displaymath}
		D_{\operatorname{qc},Z\cap W}^{\leq-c_Z}(X)
		\subseteq D_{P_Z}^{\leq0}
		\subseteq D_{\operatorname{qc},Z\cap W}^{\leq c_Z}(X),
	\end{displaymath}
	where $D_{P_Z}^{\leq0}$ denotes the aisle compactly generated by
	$P_Z$.
	Since $P$ is perfect, after shifting if necessary we may assume $P\in D_{\operatorname{qc},W}^{\leq0}(X)$.
	Then \eqref{eq:i_shriek_right_amplitude} gives $i^!P\in D_{\operatorname{qc},Z\cap W}^{\leq a+1}(X)$.
	By equivalence of the two $t$-structures, there exists $m>0$ such that
	$i^!P\in D_{P_Z}^{\leq m}$.
	Also, $\operatorname{Hom}(P[-n],P)=0$ for $n\gg0$.

	Apply \Cref{lem:BNP_weak_version_for_glueing} to the recollement
	above, viewed in the order
	\begin{displaymath}
		D_{\operatorname{qc},|U|\cap W}(U)
		\to D_{\operatorname{qc},W}(X)
		\to D_{\operatorname{qc},Z\cap W}(X).
	\end{displaymath}
	In the notation of that result, the functor denoted $j^\ast$ there is $i^!$ here, and the functor denoted $j_!$ there is $i_\ast$ here.
	Hence, $P[m]\oplus i_\ast P_Z$ is a compact generator of $D_{\operatorname{qc},W}(X)$, and its compactly generated $t$-structure is obtained by gluing the $t$-structure on $D_{\operatorname{qc},|U|\cap W}(U)$ compactly generated by $\mathbf{L}j^\ast P[m]$ with the $t$-structure on $D_{\operatorname{qc},Z\cap W}(X)$ compactly generated by $P_Z$.

	Since $U$ satisfies P.E.C. and $\mathbf{L}j^\ast P[m]$ is a compact generator, there exists $c_U\geq0$ such that
	\begin{displaymath}
		D_{\operatorname{qc},|U|\cap W}^{\leq-c_U}(U)
		\subseteq D_{\mathbf{L}j^\ast P[m]}^{\leq0}
		\subseteq D_{\operatorname{qc},|U|\cap W}^{\leq c_U}(U).
	\end{displaymath}
	Let $D_{P[m]\oplus i_\ast P_Z}^{\leq0}$
	denote the glued aisle on $D_{\operatorname{qc},W}(X)$.

	Choose $N\geq\max\{c_U,a+1+c_Z\}$.
	If $E\in D_{\operatorname{qc},W}^{\leq-N}(X)$, then $t$-exactness of $\mathbf{L}j^\ast$ gives
	\begin{displaymath}
		\mathbf{L}j^\ast E
		\in D_{\operatorname{qc},|U|\cap W}^{\leq-N}(U)
		\subseteq D_{\mathbf{L}j^\ast P[m]}^{\leq0}.
	\end{displaymath}
	Moreover, \eqref{eq:i_shriek_right_amplitude} and compatibility with shifts
	give
	\begin{displaymath}
		i^!E
		\in D_{\operatorname{qc},Z\cap W}^{\leq a+1-N}(X)
		\subseteq D_{\operatorname{qc},Z\cap W}^{\leq-c_Z}(X)
		\subseteq D_{P_Z}^{\leq0}.
	\end{displaymath}
	By the definition of the glued aisle,
	\begin{equation}
		\label{eq:standard_below_glued}
		D_{\operatorname{qc},W}^{\leq-N}(X)
		\subseteq D_{P[m]\oplus i_\ast P_Z}^{\leq0}.
	\end{equation}

	Conversely, let $E\in D_{P[m]\oplus i_\ast P_Z}^{\leq0}$.
	Then $\mathbf{L}j^\ast E \in D_{\operatorname{qc},|U|\cap W}^{\leq c_U}(U)$ and $i^!E\in D_{\operatorname{qc},Z\cap W}^{\leq c_Z}(X)$.
	By the choice of $a$, $\mathbf{R}j_\ast \mathbf{L}j^\ast E \in  D_{\operatorname{qc},W}^{\leq a+c_U}(X)$.
	Since $i_\ast$ is the natural inclusion, $i_\ast i^!E\in D_{\operatorname{qc},W}^{\leq c_Z}(X)$.
	Applying the recollement triangle \eqref{eq:rec} yields $E\in D_{\operatorname{qc},W}^{\leq \max\{a+c_U,c_Z\}}(X)$.
	Hence,
	\begin{equation}
		\label{eq:glued_below_standard}
		D_{P[m]\oplus i_\ast P_Z}^{\leq0}
		\subseteq D_{\operatorname{qc},W}^{\leq \max\{a+c_U,c_Z\}}(X).
	\end{equation}

	The inclusions \eqref{eq:standard_below_glued} and \eqref{eq:glued_below_standard} show that the standard $t$-structure on
	$D_{\operatorname{qc},W}(X)$ is equivalent to the $t$-structure compactly generated by $P[m]\oplus i_\ast P_Z$.
	Since this object is a compact generator, the latter lies in the
	preferred equivalence class.
	Thus, the standard $t$-structure on $D_{\operatorname{qc},W}(X)$ belongs to the preferred equivalence class.
	Since $V\to X$ was arbitrary, $X$ satisfies P.E.C.
\end{proof}

\begin{lemma}
	\label{lem:preferred_eq_class}
	Any quasi-compact quasi-separated algebraic space satisfies P.E.C.
\end{lemma}

\begin{proof}
	The claim is known for schemes \cite[Theorem 3.2(iii)]{Neeman:2024}. Hence, the claim holds for the schematic locus of a quasi-compact quasi-separated algebraic space. We apply the induction principle \cite[\href{https://stacks.math.columbia.edu/tag/09IT}{Tag 09IT}]{StacksProject}. By \Cref{lem:etale_nbhd_preferred_class}, we finish the proof.
\end{proof}

%%%%%%%%%%%%%%%%%%%%%%%%%%%%%%%%%%%

%%%%%%%%%%%%%%%%%%%%%%%%%%%%%%%%%%%
\section{The \texorpdfstring{$2$}{2}-category of triangulated categories}
\label{app:triang_is_2-cat}
%%%%%%%%%%%%%%%%%%%%%%%%%%%%%%%%%%%

By `$2$-category' we mean so in the strict sense \cite[Definition 1.7.8]{Riehl:26}.
$2$-categorically speaking, triangulated categories are contained in categories with translation, which in turn are contained in $1$-categories.
Symbolically,
\begin{displaymath}
	\mathbf{Triang}\subseteq\mathbf{Trans}\subseteq\mathbf{Cat},
\end{displaymath}
where the first containment is full on $2$-morphisms, and none of the containments is full on $1$-morphisms.
Let us explain this.

We extend \cite[Definition 10.1.1]{Kashiwara/Schapira:2006}:

\begin{definition}
	\label{def:cat_with_trans}
	\hfill
	\begin{enumerate}
		\item A \textbf{category with translation} is a category $\mathcal{C}$ equipped with an autoequivalence $T\colon \mathcal{C}\to\mathcal{C}$.
		\item A $1$-morphism of categories with translation $(\mathcal{C},T)\to(\mathcal{C}^\prime,T^\prime)$ is a pair $(F,\xi)$, where $F\colon \mathcal{C}\to\mathcal{C}^\prime$ is a functor and $\xi\colon FT\cong T^\prime F$ is a natural isomorphism.
		\item The unit $1$-morphism $(\mathcal{C},T)\to(\mathcal{C},T)$ is $1_{(\mathcal{C},T)}=(1_\mathcal{C},1_T)$.
		\item The composite of $1$-morphisms of categories with translation
		\begin{displaymath}
			(\mathcal{C},T)
			\xrightarrow{(F,\xi)}(\mathcal{C}^\prime,T^\prime)
			\xrightarrow{(G,\zeta)}(\mathcal{C}^{\prime\prime},T^{\prime\prime})
		\end{displaymath}
		is $(G,\zeta)(F,\xi)=(GF,\zeta_F\cdot G\xi)$.
		\item Given morphisms of categories with translation $(F,\xi),(G,\zeta)\colon (\mathcal{C},T)\rightrightarrows(\mathcal{C}^\prime,T^\prime)$, a $2$-morphism (also called a \textbf{transnatural transformation}) $(F,\xi)\Rightarrow(G,\zeta)$ is a natural transformation $\eta\colon F\Rightarrow G$ compatible with $\xi$ and $\zeta$ in the sense that
		\begin{equation}
			\label{diag:transnat}
			% https://q.uiver.app/#q=WzAsNCxbMCwwLCJGVCJdLFswLDEsIlQnRiJdLFsxLDAsIkdUIl0sWzEsMSwiVCdHIl0sWzAsMSwiXFx4aSIsMl0sWzIsMywiXFx6ZXRhIl0sWzAsMiwiXFxldGFfVCJdLFsxLDMsIlQnXFxldGEiLDJdXQ==
			\begin{tikzcd}
				FT & GT \\
				{T^\prime F} & {T^\prime G}
				\arrow["{\eta_T}", from=1-1, to=1-2]
				\arrow["\xi"', from=1-1, to=2-1]
				\arrow["\zeta", from=1-2, to=2-2]
				\arrow["{T^\prime\eta}"', from=2-1, to=2-2]
			\end{tikzcd}
		\end{equation}
		commutes.
	\end{enumerate}
\end{definition}

Write $\mathbf{Trans}$ for the collection of categories with translation with their $1$-morphisms and $2$-morphisms.
There is an obvious assignment $\mathbf{Trans}\to \mathbf{Cat}$ which forgets the translation of a category with translation, forgets the translation compatibility isomorphism of a $1$-morphism and acts as the identity on $2$-morphisms.
Categories with translation become a $2$-subcategory of $\mathbf{Cat}$ in a unique way:

\begin{lemma}
	\label{lem:trans_is_2-cat}
	Categories with translation along with the data from \Cref{def:cat_with_trans} assemble uniquely into a $2$-category $\mathbf{Trans}$ such that $\mathbf{Trans}\to \mathbf{Cat}$ is a $2$-functor.
\end{lemma}

\begin{proof}
	The assignment being faithful on the level of $2$-morphisms specifies the unit $2$-cell and the vertical and horizontal composition rules in $\mathbf{Trans}$.
	Namely, they are the ones from $\mathbf{Cat}$.
	By pasting the corresponding transnaturality squares, it is easy to see that vertical composition of $2$-morphisms satisfy the transnaturality condition \eqref{diag:transnat}.
	The proof will conclude if we verify that horizonal composition of $2$-cells in $\mathbf{Trans}$ satisfy the transnaturality condition \eqref{diag:transnat}.
	Note that this condition amounts to the equality of $2$-cells:
	\begin{displaymath}
		% https://q.uiver.app/#q=WzAsOCxbMCwwLCJcXG1hdGhjYWx7Q30iXSxbMCwyLCJcXG1hdGhjYWx7Q30iXSxbMiwwLCJcXG1hdGhjYWx7Q30nIl0sWzIsMiwiXFxtYXRoY2Fse0N9JyJdLFs0LDIsIlxcbWF0aGNhbHtDfSJdLFs2LDIsIlxcbWF0aGNhbHtDfSciXSxbNCwwLCJcXG1hdGhjYWx7Q30iXSxbNiwwLCJcXG1hdGhjYWx7Q30nIl0sWzAsMiwiRiJdLFsxLDAsIlQiXSxbMywyLCJUJyIsMl0sWzEsMywiRiIsMSx7ImN1cnZlIjotMn1dLFsxLDMsIkciLDIseyJjdXJ2ZSI6Mn1dLFswLDMsIlxceGkiLDAseyJsYWJlbF9wb3NpdGlvbiI6NDAsIm9mZnNldCI6LTEsInNob3J0ZW4iOnsic291cmNlIjoyMCwidGFyZ2V0IjozMH0sImxldmVsIjoyfV0sWzYsNywiRiIsMCx7ImN1cnZlIjotMn1dLFs0LDYsIlQiXSxbNSw3LCJUJyIsMl0sWzQsNSwiRiciLDJdLFs2LDcsIkciLDEseyJjdXJ2ZSI6Mn1dLFs2LDUsIlxcemV0YSIsMCx7Im9mZnNldCI6Mywic2hvcnRlbiI6eyJzb3VyY2UiOjMwLCJ0YXJnZXQiOjIwfSwibGV2ZWwiOjJ9XSxbMTEsMTIsIlxcZXRhIiwwLHsic2hvcnRlbiI6eyJzb3VyY2UiOjIwLCJ0YXJnZXQiOjIwfX1dLFsxNCwxOCwiXFxldGEiLDAseyJzaG9ydGVuIjp7InNvdXJjZSI6MjAsInRhcmdldCI6MjB9fV0sWzEwLDE1LCI9IiwxLHsic2hvcnRlbiI6eyJzb3VyY2UiOjIwLCJ0YXJnZXQiOjIwfSwic3R5bGUiOnsiYm9keSI6eyJuYW1lIjoibm9uZSJ9LCJoZWFkIjp7Im5hbWUiOiJub25lIn19fV1d
		\begin{tikzcd}[cramped]
			{\mathcal{C}} && {\mathcal{C}^\prime} && {\mathcal{C}} && {\mathcal{C}^\prime} \\
			\\
			{\mathcal{C}} && {\mathcal{C}^\prime} && {\mathcal{C}} && {\mathcal{C}^\prime}
			\arrow["F", from=1-1, to=1-3]
			\arrow["\xi"{pos=0.4}, shift left, shorten <=2pt, shorten >=2pt, Rightarrow, from=1-1, to=3-3]
			\arrow[""{name=0, anchor=center, inner sep=0}, "F", bend left=30, from=1-5, to=1-7]
			\arrow[""{name=1, anchor=center, inner sep=0}, "G"{description}, bend right=30, from=1-5, to=1-7]
			\arrow["\zeta", shift right=3, shorten <=2pt, shorten >=2pt, Rightarrow, from=1-5, to=3-7]
			\arrow["T", from=3-1, to=1-1]
			\arrow[""{name=2, anchor=center, inner sep=0}, "F"{description}, bend left=30, from=3-1, to=3-3]
			\arrow[""{name=3, anchor=center, inner sep=0}, "G"', bend right=30, from=3-1, to=3-3]
			\arrow[""{name=4, anchor=center, inner sep=0}, "{T^\prime}"', from=3-3, to=1-3]
			\arrow[""{name=5, anchor=center, inner sep=0}, "T", from=3-5, to=1-5]
			\arrow["{F^\prime}"', from=3-5, to=3-7]
			\arrow["{T^\prime}"', from=3-7, to=1-7]
			\arrow["\eta", shorten <=2pt, shorten >=2pt, Rightarrow, from=0, to=1]
			\arrow["\eta", shorten <=2pt, shorten >=2pt, Rightarrow, from=2, to=3]
			\arrow["{=}"{description}, draw=none, from=4, to=5]
		\end{tikzcd}
	\end{displaymath}
	Consider the following diagram in $\mathbf{Trans}$:
	\begin{displaymath}
		% https://q.uiver.app/#q=WzAsMyxbMCwwLCJcXG1hdGhjYWx7QX0iXSxbMiwwLCJcXG1hdGhjYWx7Qn0iXSxbNCwwLCJcXG1hdGhjYWx7Q30iXSxbMCwxLCIoRixcXHhpKSIsMCx7ImN1cnZlIjotMn1dLFswLDEsIihGJyxcXHhpJykiLDIseyJjdXJ2ZSI6Mn1dLFsxLDIsIihHLFxcemV0YSkiLDAseyJjdXJ2ZSI6LTJ9XSxbMSwyLCIoRycsXFx6ZXRhJykiLDIseyJjdXJ2ZSI6Mn1dLFszLDQsIlxcYWxwaGEiLDAseyJzaG9ydGVuIjp7InNvdXJjZSI6MjAsInRhcmdldCI6MjB9fV0sWzUsNiwiXFxiZXRhIiwwLHsic2hvcnRlbiI6eyJzb3VyY2UiOjIwLCJ0YXJnZXQiOjIwfX1dXQ==
		\begin{tikzcd}[cramped]
			{\mathcal{A}} && {\mathcal{B}} && {\mathcal{C}}
			\arrow[""{name=0, anchor=center, inner sep=0}, "{(F,\xi)}", bend left=30, from=1-1, to=1-3]
			\arrow[""{name=1, anchor=center, inner sep=0}, "{(F^\prime,\xi^\prime)}"', bend right=30, from=1-1, to=1-3]
			\arrow[""{name=2, anchor=center, inner sep=0}, "{(G,\zeta)}", bend left=30, from=1-3, to=1-5]
			\arrow[""{name=3, anchor=center, inner sep=0}, "{(G^\prime,\zeta^\prime)}"', bend right=30, from=1-3, to=1-5]
			\arrow["\alpha", shorten <=2pt, shorten >=2pt, Rightarrow, from=0, to=1]
			\arrow["\beta", shorten <=2pt, shorten >=2pt, Rightarrow, from=2, to=3]
		\end{tikzcd}
	\end{displaymath}
	For simplicity, we will abuse notation by denoting with same letter $T$ to the translation of $\mathcal{A}$, $\mathcal{B}$ and $\mathcal{C}$. We have to verify that the horizontal composition $\beta \ast\alpha$ (performed in $\mathbf{Cat}$) is transnatural (i.e.\ it defines a $2$-morphism $\beta \ast\alpha\colon (GF,\zeta_F\circ G\xi)\Rightarrow(G^\prime F^\prime,\zeta^\prime_{F^\prime}\circ G^\prime\xi^\prime)$ in $\mathbf{Trans}$).

	By the middle four interchange law and transnaturality of $\alpha$ and $\beta$, we get the sequence of equalities:
	% https://q.uiver.app/#q=WzAsOSxbNSwwLCJcXG1hdGhjYWx7QX0iXSxbNSwyLCJcXG1hdGhjYWx7QX0iXSxbNywwLCJcXG1hdGhjYWx7Qn0iXSxbNywyLCJcXG1hdGhjYWx7Qn0iXSxbOSwwLCJcXG1hdGhjYWx7Q30iXSxbOSwyLCJcXG1hdGhjYWx7Q30iXSxbNCwxLCI9Il0sWzAsMSwiXFxtYXRoY2Fse0F9Il0sWzMsMSwiXFxtYXRoY2Fse0N9Il0sWzAsMiwiRiJdLFsxLDAsIlQiXSxbMywyLCJUIiwyXSxbMSwzLCJGIiwxLHsiY3VydmUiOi0yfV0sWzEsMywiRiciLDIseyJjdXJ2ZSI6Mn1dLFswLDMsIlxceGkiLDAseyJsYWJlbF9wb3NpdGlvbiI6NDAsIm9mZnNldCI6LTEsInNob3J0ZW4iOnsic291cmNlIjoyMCwidGFyZ2V0IjozMH0sImxldmVsIjoyfV0sWzIsNCwiRyJdLFszLDUsIkciLDEseyJjdXJ2ZSI6LTJ9XSxbMyw1LCJHJyIsMix7ImN1cnZlIjoyfV0sWzUsNCwiVCIsMl0sWzIsNSwiXFx6ZXRhIiwwLHsibGFiZWxfcG9zaXRpb24iOjQwLCJvZmZzZXQiOi0xLCJzaG9ydGVuIjp7InNvdXJjZSI6MjAsInRhcmdldCI6MzB9LCJsZXZlbCI6Mn1dLFs3LDgsIkdGVCIsMCx7ImN1cnZlIjotNX1dLFs3LDgsIlRHJ0YnIiwyLHsiY3VydmUiOjV9XSxbNyw4LCJUR0YiLDFdLFsxMiwxMywiXFxhbHBoYSIsMCx7InNob3J0ZW4iOnsic291cmNlIjoyMCwidGFyZ2V0IjoyMH19XSxbMTYsMTcsIlxcYmV0YSIsMCx7InNob3J0ZW4iOnsic291cmNlIjoyMCwidGFyZ2V0IjoyMH19XSxbMjAsMjIsIlxcemV0YV9GXFxjZG90IEdcXHhpIiwwLHsib2Zmc2V0IjoyLCJzaG9ydGVuIjp7InNvdXJjZSI6MjAsInRhcmdldCI6MjB9fV0sWzIyLDIxLCJUKFxcYmV0YSpcXGFscGhhKSIsMCx7Im9mZnNldCI6Miwic2hvcnRlbiI6eyJzb3VyY2UiOjIwLCJ0YXJnZXQiOjIwfX1dXQ==
	\begin{displaymath}
		\begin{tikzcd}[cramped]
		&&&&& {\mathcal{A}} && {\mathcal{B}} && {\mathcal{C}} \\
		{\mathcal{A}} &&& {\mathcal{C}} & {=} \\
		&&&&& {\mathcal{A}} && {\mathcal{B}} && {\mathcal{C}}
		\arrow["F", from=1-6, to=1-8]
		\arrow["\xi"{pos=0.4}, shift left, between={0.2}{0.7}, Rightarrow, from=1-6, to=3-8]
		\arrow["G", from=1-8, to=1-10]
		\arrow["\zeta"{pos=0.4}, shift left, between={0.2}{0.7}, Rightarrow, from=1-8, to=3-10]
		\arrow[""{name=0, anchor=center, inner sep=0}, "GFT", curve={height=-30pt}, from=2-1, to=2-4]
		\arrow[""{name=1, anchor=center, inner sep=0}, "{TG'F'}"', curve={height=30pt}, from=2-1, to=2-4]
		\arrow[""{name=2, anchor=center, inner sep=0}, "TGF"{description}, from=2-1, to=2-4]
		\arrow["T", from=3-6, to=1-6]
		\arrow[""{name=3, anchor=center, inner sep=0}, "F"{description}, curve={height=-12pt}, from=3-6, to=3-8]
		\arrow[""{name=4, anchor=center, inner sep=0}, "{F'}"', curve={height=12pt}, from=3-6, to=3-8]
		\arrow["T"', from=3-8, to=1-8]
		\arrow[""{name=5, anchor=center, inner sep=0}, "G"{description}, curve={height=-12pt}, from=3-8, to=3-10]
		\arrow[""{name=6, anchor=center, inner sep=0}, "{G'}"', curve={height=12pt}, from=3-8, to=3-10]
		\arrow["T"', from=3-10, to=1-10]
		\arrow["{\zeta_F\cdot G\xi}", shift right=2, between={0.2}{0.8}, Rightarrow, from=0, to=2]
		\arrow["{T(\beta*\alpha)}", shift right=2, between={0.2}{0.8}, Rightarrow, from=2, to=1]
		\arrow["\alpha", between={0.2}{0.8}, Rightarrow, from=3, to=4]
		\arrow["\beta", between={0.2}{0.8}, Rightarrow, from=5, to=6]
		\end{tikzcd}
	\end{displaymath}

	% https://q.uiver.app/#q=WzAsOCxbMSwwLCJcXG1hdGhjYWx7QX0iXSxbMywwLCJcXG1hdGhjYWx7Qn0iXSxbNSwwLCJcXG1hdGhjYWx7Q30iXSxbNiwwLCI9Il0sWzcsMCwiXFxtYXRoY2Fse0F9Il0sWzksMCwiXFxtYXRoY2Fse0J9Il0sWzExLDAsIlxcbWF0aGNhbHtDfSJdLFswLDAsIj0iXSxbMCwxLCJGVCIsMCx7ImN1cnZlIjotNH1dLFswLDEsIlRGJyIsMix7ImN1cnZlIjo0fV0sWzAsMSwiVEYiLDFdLFsxLDIsIkdUIiwwLHsiY3VydmUiOi00fV0sWzEsMiwiVEcnIiwyLHsiY3VydmUiOjR9XSxbMSwyLCJURyIsMV0sWzQsNSwiRlQiLDAseyJjdXJ2ZSI6LTR9XSxbNCw1LCJURiciLDIseyJjdXJ2ZSI6NH1dLFs1LDYsIkdUIiwwLHsiY3VydmUiOi00fV0sWzUsNiwiVEcnIiwyLHsiY3VydmUiOjR9XSxbNCw1LCJGJ1QiLDFdLFs1LDYsIkcnVCIsMV0sWzgsMTAsIlxceGkiLDAseyJzaG9ydGVuIjp7InNvdXJjZSI6MjAsInRhcmdldCI6MjB9fV0sWzEwLDksIlRcXGFscGhhIiwwLHsic2hvcnRlbiI6eyJzb3VyY2UiOjIwLCJ0YXJnZXQiOjIwfX1dLFsxMSwxMywiXFx6ZXRhIiwwLHsic2hvcnRlbiI6eyJzb3VyY2UiOjIwLCJ0YXJnZXQiOjIwfX1dLFsxMywxMiwiVFxcYmV0YSIsMCx7InNob3J0ZW4iOnsic291cmNlIjoyMCwidGFyZ2V0IjoyMH19XSxbMTQsMTgsIlxcYWxwaGFfVCIsMCx7InNob3J0ZW4iOnsic291cmNlIjoyMCwidGFyZ2V0IjoyMH19XSxbMTgsMTUsIlxceGknIiwwLHsic2hvcnRlbiI6eyJzb3VyY2UiOjIwLCJ0YXJnZXQiOjIwfX1dLFsxNiwxOSwiXFxiZXRhX1QiLDAseyJzaG9ydGVuIjp7InNvdXJjZSI6MjAsInRhcmdldCI6MjB9fV0sWzE5LDE3LCJcXHpldGEnIiwwLHsic2hvcnRlbiI6eyJzb3VyY2UiOjIwLCJ0YXJnZXQiOjIwfX1dXQ==
	\begin{displaymath}
		\begin{tikzcd}[cramped]
		{=} & {\mathcal{A}} && {\mathcal{B}} && {\mathcal{C}} & {=} & {\mathcal{A}} && {\mathcal{B}} && {\mathcal{C}}
		\arrow[""{name=0, anchor=center, inner sep=0}, "FT", curve={height=-24pt}, from=1-2, to=1-4]
		\arrow[""{name=1, anchor=center, inner sep=0}, "{TF'}"', curve={height=24pt}, from=1-2, to=1-4]
		\arrow[""{name=2, anchor=center, inner sep=0}, "TF"{description}, from=1-2, to=1-4]
		\arrow[""{name=3, anchor=center, inner sep=0}, "GT", curve={height=-24pt}, from=1-4, to=1-6]
		\arrow[""{name=4, anchor=center, inner sep=0}, "{TG'}"', curve={height=24pt}, from=1-4, to=1-6]
		\arrow[""{name=5, anchor=center, inner sep=0}, "TG"{description}, from=1-4, to=1-6]
		\arrow[""{name=6, anchor=center, inner sep=0}, "FT", curve={height=-24pt}, from=1-8, to=1-10]
		\arrow[""{name=7, anchor=center, inner sep=0}, "{TF'}"', curve={height=24pt}, from=1-8, to=1-10]
		\arrow[""{name=8, anchor=center, inner sep=0}, "{F'T}"{description}, from=1-8, to=1-10]
		\arrow[""{name=9, anchor=center, inner sep=0}, "GT", curve={height=-24pt}, from=1-10, to=1-12]
		\arrow[""{name=10, anchor=center, inner sep=0}, "{TG'}"', curve={height=24pt}, from=1-10, to=1-12]
		\arrow[""{name=11, anchor=center, inner sep=0}, "{G'T}"{description}, from=1-10, to=1-12]
		\arrow["\xi", between={0.2}{0.8}, Rightarrow, from=0, to=2]
		\arrow["{T\alpha}", between={0.2}{0.8}, Rightarrow, from=2, to=1]
		\arrow["\zeta", between={0.2}{0.8}, Rightarrow, from=3, to=5]
		\arrow["{T\beta}", between={0.2}{0.8}, Rightarrow, from=5, to=4]
		\arrow["{\alpha_T}", between={0.2}{0.8}, Rightarrow, from=6, to=8]
		\arrow["{\xi'}", between={0.2}{0.8}, Rightarrow, from=8, to=7]
		\arrow["{\beta_T}", between={0.2}{0.8}, Rightarrow, from=9, to=11]
		\arrow["{\zeta'}", between={0.2}{0.8}, Rightarrow, from=11, to=10]
		\end{tikzcd}
	\end{displaymath}

	% https://q.uiver.app/#q=WzAsMTAsWzAsMSwiPSJdLFsxLDAsIlxcbWF0aGNhbHtBfSJdLFszLDAsIlxcbWF0aGNhbHtCfSJdLFs1LDAsIlxcbWF0aGNhbHtDfSJdLFsxLDIsIlxcbWF0aGNhbHtBfSJdLFszLDIsIlxcbWF0aGNhbHtCfSJdLFs1LDIsIlxcbWF0aGNhbHtDfSJdLFs2LDEsIj0iXSxbNywxLCJcXG1hdGhjYWx7QX0iXSxbMTAsMSwiXFxtYXRoY2Fse0N9Il0sWzQsMSwiVCJdLFs1LDIsIlQiLDFdLFs2LDMsIlQiLDJdLFs0LDUsIkYnIiwyXSxbNSw2LCJHJyIsMl0sWzEsMiwiRiIsMCx7ImN1cnZlIjotMn1dLFsxLDIsIkYnIiwxLHsiY3VydmUiOjJ9XSxbMiwzLCJHIiwwLHsiY3VydmUiOi0yfV0sWzIsMywiRyciLDEseyJjdXJ2ZSI6Mn1dLFsxLDUsIlxceGknIiwyLHsic2hvcnRlbiI6eyJzb3VyY2UiOjMwLCJ0YXJnZXQiOjIwfSwibGV2ZWwiOjJ9XSxbMiw2LCJcXHpldGEnIiwyLHsic2hvcnRlbiI6eyJzb3VyY2UiOjMwLCJ0YXJnZXQiOjIwfSwibGV2ZWwiOjJ9XSxbOCw5LCJHRlQiLDAseyJjdXJ2ZSI6LTV9XSxbOCw5LCJURydGJyIsMix7ImN1cnZlIjo1fV0sWzgsOSwiRydGJ1QiLDFdLFsxNSwxNiwiXFxhbHBoYSIsMCx7InNob3J0ZW4iOnsic291cmNlIjoyMCwidGFyZ2V0IjoyMH19XSxbMTcsMTgsIlxcYmV0YSIsMCx7InNob3J0ZW4iOnsic291cmNlIjoyMCwidGFyZ2V0IjoyMH19XSxbMjEsMjMsIihcXGJldGEqXFxhbHBoYSlfVCIsMCx7Im9mZnNldCI6Miwic2hvcnRlbiI6eyJzb3VyY2UiOjIwLCJ0YXJnZXQiOjIwfX1dLFsyMywyMiwiXFx6ZXRhJ197Rid9XFxjaXJjIEcnXFx4aSciLDAseyJvZmZzZXQiOjIsInNob3J0ZW4iOnsic291cmNlIjoyMCwidGFyZ2V0IjoyMH19XV0=
	\begin{displaymath}
		\begin{tikzcd}[cramped]
		& {\mathcal{A}} && {\mathcal{B}} && {\mathcal{C}} &&&&& \\
		{=} &&&&&& {=} & {\mathcal{A}} &&& {\mathcal{C}} \\
		& {\mathcal{A}} && {\mathcal{B}} && {\mathcal{C},}
		\arrow[""{name=0, anchor=center, inner sep=0}, "F", curve={height=-12pt}, from=1-2, to=1-4] 
		\arrow[""{name=1, anchor=center, inner sep=0}, "{F'}"{description}, curve={height=12pt}, from=1-2, to=1-4]
		\arrow["{\xi'}"', between={0.3}{0.8}, Rightarrow, from=1-2, to=3-4]
		\arrow[""{name=2, anchor=center, inner sep=0}, "G", curve={height=-12pt}, from=1-4, to=1-6]
		\arrow[""{name=3, anchor=center, inner sep=0}, "{G'}"{description}, curve={height=12pt}, from=1-4, to=1-6]
		\arrow["{\zeta'}"', between={0.3}{0.8}, Rightarrow, from=1-4, to=3-6]
		\arrow[""{name=4, anchor=center, inner sep=0}, "GFT", curve={height=-30pt}, from=2-8, to=2-11]
		\arrow[""{name=5, anchor=center, inner sep=0}, "{TG'F'}"', curve={height=30pt}, from=2-8, to=2-11]
		\arrow[""{name=6, anchor=center, inner sep=0}, "{G'F'T}"{description}, from=2-8, to=2-11]
		\arrow["T", from=3-2, to=1-2]
		\arrow["{F'}"', from=3-2, to=3-4]
		\arrow["T"{description}, from=3-4, to=1-4]
		\arrow["{G'}"', from=3-4, to=3-6]
		\arrow["T"', from=3-6, to=1-6]
		\arrow["\alpha", between={0.2}{0.8}, Rightarrow, from=0, to=1]
		\arrow["\beta", between={0.2}{0.8}, Rightarrow, from=2, to=3]
		\arrow["{(\beta*\alpha)_T}", shift right=2, between={0.2}{0.8}, Rightarrow, from=4, to=6]
		\arrow["{\zeta'_{F'}\circ G'\xi'}", shift right=2, between={0.2}{0.8}, Rightarrow, from=6, to=5]
		\end{tikzcd}
	\end{displaymath}
	as desired.
\end{proof}

\begin{definition}
	\hfill
	\begin{enumerate}
		\item A \textbf{transadjunction} is an adjunction in the $2$-category $\mathbf{Trans}$.
		\item A \textbf{transequivalence} is an equivalence in the $2$-category $\mathbf{Trans}$.
	\end{enumerate}
\end{definition}

\begin{lemma}
	\label{lem:equiv_transadj}
	Let $(F,\xi)\colon (\mathcal{C},T)\to(\mathcal{C}^\prime,T^\prime)$ and
	$(G,\zeta)\colon (\mathcal{C}^\prime,T^\prime)\to(\mathcal{C},T)$ be functors of categories with translation.
	Suppose $F\dashv G$ in $\mathbf{Cat}$ with unit and counit $\eta$ and $\epsilon$, respectively.
	The following are equivalent:
	\begin{enumerate}[label=(\arabic*), ref=\thelemma(\arabic*)]
		\item \label[lemma]{lem:equiv_transadj1}
		$(\xi,\zeta^{-1})$ are mates under the diagram
		\begin{equation}
			\label{diag:mates_transadj}
			% https://q.uiver.app/#q=WzAsNCxbMCwwLCJcXG1hdGhjYWx7Q30iXSxbMSwwLCJcXG1hdGhjYWx7Q30iXSxbMCwxLCJcXG1hdGhjYWx7Q31eXFxwcmltZSJdLFsxLDEsIlxcbWF0aGNhbHtDfV5cXHByaW1lIl0sWzAsMSwiVCJdLFsyLDMsIlReXFxwcmltZSIsMl0sWzIsMCwiRyIsMix7ImN1cnZlIjoxfV0sWzAsMiwiRiIsMix7ImN1cnZlIjoxfV0sWzMsMSwiRyIsMix7ImN1cnZlIjoxfV0sWzEsMywiRiIsMix7ImN1cnZlIjoxfV0sWzksOCwiXFxkYXNodiIsMyx7InNob3J0ZW4iOnsic291cmNlIjozMCwidGFyZ2V0IjozMH0sImxldmVsIjoxLCJzdHlsZSI6eyJib2R5Ijp7Im5hbWUiOiJub25lIn0sImhlYWQiOnsibmFtZSI6Im5vbmUifX19XSxbNyw2LCJcXGRhc2h2IiwzLHsic2hvcnRlbiI6eyJzb3VyY2UiOjMwLCJ0YXJnZXQiOjMwfSwibGV2ZWwiOjEsInN0eWxlIjp7ImJvZHkiOnsibmFtZSI6Im5vbmUifSwiaGVhZCI6eyJuYW1lIjoibm9uZSJ9fX1dXQ==
			\begin{tikzcd}[cramped]
				{\mathcal{C}} & {\mathcal{C}} \\
				{\mathcal{C}^\prime} & {\mathcal{C}^\prime\mathrlap{.}}
				\arrow["T", from=1-1, to=1-2]
				\arrow[""{name=0, anchor=center, inner sep=0}, "F"', bend right=30, from=1-1, to=2-1]
				\arrow[""{name=1, anchor=center, inner sep=0}, "F"', bend right=30, from=1-2, to=2-2]
				\arrow[""{name=2, anchor=center, inner sep=0}, "G"', bend right=30, from=2-1, to=1-1]
				\arrow["{T^\prime}"', from=2-1, to=2-2]
				\arrow[""{name=3, anchor=center, inner sep=0}, "G"', bend right=30, from=2-2, to=1-2]
				\arrow["\dashv"{marking, allow upside down}, draw=none, from=0, to=2]
				\arrow["\dashv"{marking, allow upside down}, draw=none, from=1, to=3]
			\end{tikzcd}
		\end{equation}
		\item \label[lemma]{lem:equiv_transadj2}
		$\epsilon$ is transnatural.
		\item \label[lemma]{lem:equiv_transadj3}
		$\eta$ is transnatural.
		\item \label[lemma]{lem:equiv_transadj4}
		The following diagram commutes:
		\begin{displaymath}
			% https://q.uiver.app/#q=WzAsNixbMCwwLCJcXG1hdGhjYWx7Q31eXFxwcmltZShGYyxkKSJdLFswLDEsIlxcbWF0aGNhbHtDfV5cXHByaW1lKFReXFxwcmltZSBGYyxUXlxccHJpbWUgZCkiXSxbMCwyLCJcXG1hdGhjYWx7Q31eXFxwcmltZSAoRl5cXHByaW1lIFRjLFReXFxwcmltZSBkKSJdLFsxLDAsIlxcbWF0aGNhbHtDfShjLEdkKSJdLFsxLDEsIlxcbWF0aGNhbHtDfShUYyxUR2QpIl0sWzEsMiwiXFxtYXRoY2Fse0N9KFRjLEdeXFxwcmltZSBUXlxccHJpbWUgZCkuIl0sWzAsMywiXFxjb25nIl0sWzIsNSwiXFxjb25nIl0sWzAsMSwiVF5cXHByaW1lIiwyXSxbMSwyLCJcXHhpX2NeXFxhc3QiLDJdLFszLDQsIlQiXSxbNCw1LCJcXHpldGFfe2QsXFxhc3R9XnstMX0iXV0=
			\begin{tikzcd}[cramped]
				{\mathcal{C}^\prime(Fc,d)} & {\mathcal{C}(c,Gd)} \\
				{\mathcal{C}^\prime(T^\prime Fc,T^\prime d)} & {\mathcal{C}(Tc,TGd)} \\
				{\mathcal{C}^\prime (F^\prime Tc,T^\prime d)} & {\mathcal{C}(Tc,G^\prime T^\prime d).}
				\arrow["\cong", from=1-1, to=1-2]
				\arrow["{T^\prime}"', from=1-1, to=2-1]
				\arrow["T", from=1-2, to=2-2]
				\arrow["{\xi_c^\ast}"', from=2-1, to=3-1]
				\arrow["{\zeta_{d,\ast}^{-1}}", from=2-2, to=3-2]
				\arrow["\cong", from=3-1, to=3-2]
			\end{tikzcd}
		\end{displaymath}
		\item \label[lemma]{lem:equiv_transadj5} $\eta$ and $\epsilon$ are transnatural and $(F,G,\eta,\epsilon)$ is an adjunction in $\mathbf{Trans}$.
	\end{enumerate}
\end{lemma}

\begin{proof}
	The equivalence of the first four points is \Cref{lem:morph_adj}.
	The equivalence \ref{lem:equiv_transadj5}$\iff$(\ref{lem:equiv_transadj2}$+$\ref{lem:equiv_transadj3}) is clear.
\end{proof}

Recall that given a $2$-category $\mathcal{T}$, an equivalence in $\mathcal{T}$ is a $1$-morphism $F\colon \mathcal{C}\to\mathcal{C}^\prime$ for which there exists a $1$-morphism $G\colon \mathcal{C}^\prime\to\mathcal{C}$ along with $2$-isomorphisms $\operatorname{id}_{\mathcal{C}}\cong GF$ and $FG\cong\operatorname{id}_{\mathcal{C}^\prime}$.
The $1$-morphism $G$ is called a quasi-inverse for $F$.
Quasi-inverses are unique up to isomorphism.
For a $1$-morphism in $\mathbf{Trans}$, testing adjointness and equivalence in $\mathbf{Triang}$ amounts to doing so in $\mathbf{Cat}$:

\begin{lemma}
	\label{lem:adj_and_equiv_in_trans}
	Let $(F,\xi)\colon (\mathcal{C},T)\to(\mathcal{C}^\prime,T^\prime)$ be a functor of categories with translation.
	\begin{enumerate}[label=(\arabic*), ref=\thelemma(\arabic*)]
		\item \label[lemma]{lem:adj_in_trans}
		$F$ has a right adjoint in $\mathbf{Cat}$ if, and only if, it has a right adjoint in $\mathbf{Trans}$.
		Actually, a right adjoint to $F$ in $\mathbf{Cat}$ promotes to a right adjoint in $\mathbf{Trans}$ in a unique way. Namely,
		suppose $F$ has a right adjoint $G$ in $\mathbf{Cat}$, with unit and counit $\eta$ and $\epsilon$.
		Then there is a unique natural transformation $\zeta\colon GT^\prime\cong TG$ such that $\eta$ and $\epsilon$ are transnatural.
		Thus $(F,\xi)\dashv(G,\zeta)$ in $\mathbf{Trans}$ with unit and counit $\eta$ and $\epsilon$.
		There is a similar statement for left adjoints.
		\item \label[lemma]{lem:equiv_in_trans}
		$F$ is an equivalence in $\mathbf{Cat}$ if, and only if, it is an equivalence in $\mathbf{Trans}$.
	\end{enumerate}
\end{lemma}

\begin{proof}
	In \ref{lem:adj_in_trans} and \ref{lem:equiv_in_trans} the converses should be clear.

	We prove the claim in \ref{lem:adj_in_trans}.
	\Cref{lem:equiv_transadj} gives uniqueness and tells us what $\zeta$ should be.
	Indeed, the right mate $TG\to GT^\prime$ of $\xi$ under the diagram \eqref{diag:mates_transadj} is an isomorphism by \Cref{lem:when_mate_correspondence_preserves_isos3}.
	One defines $\zeta$ to be the inverse isomorphism $GT^\prime\xrightarrow{\cong}TG$.

	The proof of \ref{lem:equiv_in_trans} follows from the fact that any equivalence can be promoted to an adjoint equivalence \cite[Corollary 4.3.6]{Riehl:26},
	so one reduces to \ref{lem:adj_in_trans}.
\end{proof}

\begin{definition}
	\label{def:trinat_trans}
	Let $(F,\xi),(G,\zeta)\colon (\mathcal{D},T)\to(\mathcal{D}^\prime,T^\prime)$ be triangulated functors between triangulated categories. A \textbf{trinatural transformation} $(F,\xi)\Rightarrow (G,\zeta)$ is, for each object $A\in\mathcal{D}$, a morphism $\eta_A\colon FA\to GA$, subject to the following property: For every distinguished triangle $A\to B\to C\xrightarrow{h} TA$ in $\mathcal{D}$, the following diagram commutes:
	\begin{equation}
		\label{diag:trinat}
		% https://q.uiver.app/#q=WzAsOCxbMCwwLCJGQSJdLFsxLDAsIkZCIl0sWzIsMCwiRkMiXSxbMywwLCJUJ0ZBIl0sWzAsMSwiR0EiXSxbMSwxLCJHQiJdLFsyLDEsIkdDIl0sWzMsMSwiVCdHQSJdLFswLDFdLFsxLDJdLFsyLDMsIlxceGlfQVxcY2lyYyBGKGgpIl0sWzAsNCwiXFxldGFfQSJdLFs0LDVdLFs1LDZdLFs2LDcsIlxcemV0YV9BXFxjaXJjIEcoaCkiLDJdLFszLDcsIlReXFxwcmltZVxcZXRhX0EiXSxbMiw2LCJcXGV0YV9DIl0sWzEsNSwiXFxldGFfQiJdXQ==
		\begin{tikzcd}
			FA & FB & FC &[5pt] {T^\prime FA} \\
			GA & GB & GC & {T^\prime GA}.
			\arrow[from=1-1, to=1-2]
			\arrow["{\eta_A}", from=1-1, to=2-1]
			\arrow[from=1-2, to=1-3]
			\arrow["{\eta_B}", from=1-2, to=2-2]
			\arrow["{\xi_A\circ F(h)}"{yshift=4pt}, from=1-3, to=1-4]
			\arrow["{\eta_C}", from=1-3, to=2-3]
			\arrow["{T^\prime\eta_A}", from=1-4, to=2-4]
			\arrow[from=2-1, to=2-2]
			\arrow[from=2-2, to=2-3]
			\arrow["{\zeta_A\circ G(h)}"'{yshift=-4pt}, from=2-3, to=2-4]
		\end{tikzcd}
	\end{equation}
\end{definition}

\begin{lemma}
	\label{lem:transnat_is_trinat}
	Let $(F,\xi),(G,\zeta)\colon (\mathcal{D},T)\to(\mathcal{D}^\prime,T^\prime)$ be triangulated functors between triangulated categories.
	A collection of morphisms $\eta_A\colon FA\to GA$, $A\in\mathcal{D}$, defines a trinatural transformation $(F,\xi)\Rightarrow (G,\zeta)$ if, and only if,
	it defines a transnatural transformation $(F,\xi)\Rightarrow (G,\zeta)$.
\end{lemma}

Therefore a trinatural transformation of triangulated functors coincides with the notion of `morphism of $\Delta$-functors' in the sense of \cite[\S 1.5]{Lipman:2009}.

\begin{proof}
	Write $\eta=\{\eta_A\}_{A\in\mathcal{D}}$.

	Suppose first we have a trinatural transformation. Let $A\to B$ be any morphism in $\mathcal{D}$. Extend it to a distinguished triangle $A\to B\to C\to TA$. Commutativity of the left square of \eqref{diag:trinat} implies naturality of $\eta$. On the other hand, the triangle $A\to 0\to TA\xrightarrow{\operatorname{id}}TA$ is distinguished in $\mathcal{D}$. Hence commutativity of the right square of \eqref{diag:trinat} realized in this triangle becomes the transnaturality condition \eqref{diag:transnat}.

	Conversely, assume $\eta$ is transnatural and let $A\to B\to C\to TA$ be a distinguished triangle in $\mathcal{D}$. Then the first and middle square of \eqref{diag:trinat} commutes, by naturality of $\eta$. The right square of \eqref{diag:trinat} decomposes as
	\begin{displaymath}
		% https://q.uiver.app/#q=WzAsNixbMCwwLCJGQyJdLFsxLDAsIkZUQSJdLFsyLDAsIlQnRkEiXSxbMCwxLCJHQyJdLFsxLDEsIkdUQSJdLFsyLDEsIlQnR0EiXSxbMCwxLCJGaCJdLFszLDQsIkdoIl0sWzAsM10sWzEsNF0sWzIsNV0sWzQsNSwiXFx6ZXRhX0EiXSxbMSwyLCJcXHhpX0EiXV0=
		\begin{tikzcd}[cramped]
			FC & FTA & {T^\prime FA} \\
			GC & GTA & {T^\prime GA}
			\arrow["Fh", from=1-1, to=1-2]
			\arrow[from=1-1, to=2-1]
			\arrow["{\xi_A}", from=1-2, to=1-3]
			\arrow[from=1-2, to=2-2]
			\arrow[from=1-3, to=2-3]
			\arrow["Gh", from=2-1, to=2-2]
			\arrow["{\zeta_A}", from=2-2, to=2-3]
		\end{tikzcd}
	\end{displaymath}
	In this diagram, the left (resp.\ right) square commutes by naturality (resp.\ transnaturality) of $\eta$.
\end{proof}

\begin{definition}
	\label{def:triang}
	The category of triangulated categories $\mathbf{Triang}$ has as objects the triangulated categories, as $1$-morphisms the triangulated functors and as $2$-morphisms the trinatural transformations.
	By \Cref{lem:trans_is_2-cat,lem:transnat_is_trinat}, it follows that $\mathbf{Triang}$ is a $2$-category, and actually a $2$-subcategory of $\mathbf{Trans}$.
\end{definition}

\begin{definition}
	\hfill
	\begin{enumerate}
		\item A \textbf{triadjunction} is an adjunction $F\dashv G$ in the $2$-category $\mathbf{Triang}$.
		We say $F$ is \textbf{left triadjoint} to $G$ and that $G$ is \textbf{right triadjoint} to $F$.
		\item A \textbf{triequivalence} (or \textbf{triangulated equivalence}) is an equivalence in the $2$-category $\mathbf{Triang}$.
	\end{enumerate}
\end{definition}

For triadjunctions, the equivalence of \eqref{lem:equiv_transadj2}, \eqref{lem:equiv_transadj3} and \eqref{lem:equiv_transadj4} is
\cite[Lemma-Definition 3.3.1]{Lipman:2009}.
For a $1$-morphism in $\mathbf{Triang}$, testing adjointness and equivalence in $\mathbf{Triang}$ amounts to doing so in $\mathbf{Cat}$:

\begin{lemma}
	Let $(F,\xi)$ be a triangulated functor of triangulated categories.
	\begin{enumerate}[label=(\arabic*), ref=\thelemma(\arabic*)]
		\item \label[lemma]{lem:triadjoint_iff_adjoint}
		$F$ has a right adjoint in $\mathbf{Cat}$ if, and only if, it has a right adjoint in $\mathbf{Triang}$.
		Actually, a right adjoint to $F$ in $\mathbf{Cat}$ promotes to a right adjoint in $\mathbf{Triang}$ in a unique way.
		Namely, let $G$ be a right adjoint to $F$ in $\mathbf{Cat}$ with unit $\eta\colon \operatorname{id}\to GF$ and counit $\epsilon\colon FG\to\operatorname{id}$.
		Then there is a unique natural isomorphism $\zeta\colon GT\cong TG$ such that $(G,\zeta)$ is a triangulated functor and $\eta$ and $\epsilon$ are trinatural.
		Therefore, $(G,\zeta)\dashv(F,\xi)$ in $\mathbf{Triang}$ with unit and counit $\eta$ and $\epsilon$.
		There is a similar statement for left adjoints.
		\item \label[lemma]{lem:triequiv_iff_equiv}
		$F$ is an equivalence in $\mathbf{Cat}$ if, and only if, it is an equivalence in $\mathbf{Triang}$.
	\end{enumerate}
\end{lemma}

\begin{proof}
	In \ref{lem:triadjoint_iff_adjoint} and \ref{lem:triequiv_iff_equiv} the converses should be clear.

	We prove the claim in \ref{lem:triadjoint_iff_adjoint}.
	Define $\zeta$ to be the isomorphism arising in \Cref{lem:adj_in_trans}.
	If one shows $(G,\zeta)$ is triangulated, then we are done by \Cref{lem:transnat_is_trinat}.
	The fact that $(G,\zeta)$ is triangulated is shown in the proof of \cite[(3.3.8)]{Lipman:2009}.

	The proof of \ref{lem:triequiv_iff_equiv} follows from the fact that any equivalence can be promoted to an adjoint equivalence \cite[Corollary 4.3.6]{Riehl:26},
	so one reduces to \ref{lem:triadjoint_iff_adjoint}.
\end{proof}

%%%%%%%%%%%%%%%%%%%%%%%%%%%%%%%%%%%
\section{Comparison with enhanced approaches}
\label{app:higher_categories}
%%%%%%%%%%%%%%%%%%%%%%%%%%%%%%%%%%%

We compare our results with kernel theorems for enhanced categories \cite{BenZvi/Francis/Nadler:2010,Ben-Zvi/Nadler/Preygel:2017}.
Moreover, we give an adjunction between integral transforms whose unit is not induced by a morphism of kernels.
The following is straightforward.

\begin{lemma}
	\label{cmp:lem:units}
	Let $\mathcal{T}$ be a $2$-category.
	Let $F\dashv G$ be an adjunction in $\mathcal{T}$ with unit $\eta$ and counit $\epsilon$.
	\begin{enumerate}[label=(\arabic*),ref=\arabic*]
		\item \label{cmp:lem:units1}
		If $\tau\colon G\Rightarrow G$ is an invertible $2$-morphism, then $F\dashv G$ is an adjunction in $\mathcal{T}$ with unit $(\tau F)\circ\eta$ and counit $\epsilon\circ(F\tau^{-1})$.
		\item \label{cmp:lem:units2}
		If $F\dashv G$ also has unit $\eta^\prime$ and counit $\epsilon^\prime$, then $\tau=(G\epsilon)\circ(\eta^\prime G)$ is the unique invertible $2$-morphism $G\Rightarrow G$ satisfying $\eta^\prime=(\tau F)\circ\eta$.
		Moreover, $\epsilon^\prime=\epsilon\circ(F\tau^{-1})$.
	\end{enumerate}
	In particular this holds in $\mathcal{T}=\mathbf{Triang}$, see \Cref{def:triang}.
\end{lemma}

\begin{proof}
	We use the interchange law \cite[Appendix, \S 1]{Godement:1958}.
	For \eqref{cmp:lem:units1}, the triangle identities give
	\begin{displaymath}
		\begin{aligned}
			(\epsilon F)\circ(F\tau^{-1}F)\circ(F\tau F)\circ(F\eta)
			&=(\epsilon F)\circ(F\eta)=1_F
		\end{aligned}
	\end{displaymath}
	and
	\begin{displaymath}
		\begin{aligned}
			(G\epsilon)\circ(GF\tau^{-1})\circ(\tau FG)\circ(\eta G)
			&=\tau\circ(G\epsilon)\circ(\eta G)\circ\tau^{-1}=1_G.
		\end{aligned}
	\end{displaymath}
	For \eqref{cmp:lem:units2}, set $\sigma:=(G\epsilon^\prime)\circ(\eta G)$.
	Then
	\begin{displaymath}
		\begin{aligned}
			\tau\circ\sigma
			&=(G\epsilon)\circ(\eta^\prime G)\circ(G\epsilon^\prime)\circ(\eta G)\\
			&=(G\epsilon^\prime)\circ(G\epsilon FG)\circ(GF\eta G)\circ(\eta^\prime G)\\
			&=(G\epsilon^\prime)\circ(\eta^\prime G)=1_G.
		\end{aligned}
	\end{displaymath}
	By symmetry, $\sigma\circ\tau=1_G$.
	Moreover,
	\begin{displaymath}
		\begin{aligned}
			(\tau F)\circ\eta
			&=(G\epsilon F)\circ(\eta^\prime GF)\circ\eta
			=(G\epsilon F)\circ(GF\eta)\circ\eta^\prime=\eta^\prime
		\end{aligned}
	\end{displaymath}
	and
	\begin{displaymath}
		\begin{aligned}
			\epsilon\circ(F\sigma)
			&=\epsilon\circ(FG\epsilon^\prime)\circ(F\eta G)
			=\epsilon^\prime\circ(\epsilon FG)\circ(F\eta G)=\epsilon^\prime.
		\end{aligned}
	\end{displaymath}
	If $\tau_0\colon G\Rightarrow G$ is invertible with $\eta^\prime=(\tau_0 F)\circ\eta$, then
	\begin{displaymath}
		\begin{aligned}
			(G\epsilon)\circ(\eta^\prime G)
			&=(G\epsilon)\circ(\tau_0 FG)\circ(\eta G)
			=\tau_0\circ(G\epsilon)\circ(\eta G)=\tau_0.
		\end{aligned}
	\end{displaymath}
	This completes the proof.
\end{proof}

\begin{definition}
	\label{cmp:def:induced}
	Let $f_i\colon Y_i\to S$ be proper flat morphisms of Noetherian algebraic spaces for $i=1,2$.
	Denote by $f_j^\prime \colon Y_1\times_S Y_2 \to Y_i$ be the canonical projection for $j\not= i \in \{1,2\}$.
	A morphism $w\colon M\to N$ in $D_{\operatorname{qc}}(Y_1\times_S Y_2)$ induces a $2$-morphism $\Phi_w\colon\Phi_M\Rightarrow\Phi_N$ in $\mathbf{Triang}$ with components
	\begin{displaymath}
		(\Phi_w)_A=\mathbf{R}(f_1^\prime)_\ast(\mathbf{L}(f_2^\prime)^\ast A\otimes^{\mathbf{L}}w),\qquad A\in D_{\operatorname{qc}}(Y_1).
	\end{displaymath}
	A $2$-morphism $\Phi_M\Rightarrow\Phi_N$ is \textbf{induced by a morphism of kernels} if it is of the form $\Phi_w$.
\end{definition}

\begin{example-claim}
	\label{cmp:ex:unit}
	Let $k$ be an algebraically closed field.
	Suppose $X$ is a smooth projective curve over $k$.
	Set $S=\operatorname{Spec}(k)$ and $Y_1=Y_2=X$.
	Denote by $\Delta\colon X\to X\times_k X$ the diagonal.
	Consider the kernels
	\begin{displaymath}
		K=\mathbf{R}\Delta_\ast\mathcal{O}_X\oplus\mathbf{R}\Delta_\ast\mathcal{O}_X[1],
		\qquad
		K^\prime=\mathbf{R}\Delta_\ast\mathcal{O}_X\oplus\mathbf{R}\Delta_\ast\mathcal{O}_X[-1]
	\end{displaymath}
	in $D^b_{\operatorname{coh}}(X\times_k X)$.
	Since $X\times_k X$ is regular, $K$ is perfect.
	Hence, by \Cref{rmk:relative_perfects}\eqref{rmk:relative_perfects3}, $K$ is $f^\prime_i$-perfect for each canonical projection $f_j^\prime \colon X\times_k X \to X$.
	By the projection formula, $\Phi_{\mathbf{R}\Delta_\ast\mathcal{O}_X[n]}$ is naturally isomorphic to the shift $[n]$ for every integer $n$.
	Hence, it follows that
	\begin{equation}
		\label{cmp:eq:composite}
		\begin{aligned}
			\Phi_{K^\prime}\circ\Phi_K
			&\cong \operatorname{id}\oplus[1]\oplus[-1]\oplus \operatorname{id}\\
			&\cong\Phi_{\mathbf{R}\Delta_\ast\mathcal{O}_X\oplus\mathbf{R}\Delta_\ast\mathcal{O}_X[1]\oplus\mathbf{R}\Delta_\ast\mathcal{O}_X[-1]\oplus\mathbf{R}\Delta_\ast\mathcal{O}_X}.
		\end{aligned}
	\end{equation}
	The kernel $K^\prime$ is isomorphic to the right-adjoint kernel of \Cref{prop:right_adjoint_pullback}.
	Indeed, it follows from the computation:
	\begin{displaymath}
		\begin{aligned}
			&\operatorname{\mathbf{R}\mathcal{H}\! \mathit{om}}(K ,(f_1^\prime)^!\mathcal{O}_X)
			\\&\cong\operatorname{\mathbf{R}\mathcal{H}\! \mathit{om}}(\mathbf{R}\Delta_\ast\mathcal{O}_X\oplus\mathbf{R}\Delta_\ast\mathcal{O}_X[1] ,(f_1^\prime)^!\mathcal{O}_X)
			\\&\cong\operatorname{\mathbf{R}\mathcal{H}\! \mathit{om}}(\mathbf{R}\Delta_\ast (\mathcal{O}_X\oplus \mathcal{O}_X[1] ),(f_1^\prime)^!\mathcal{O}_X)
			\\&\cong \mathbf{R}\Delta_\ast \operatorname{\mathbf{R}\mathcal{H}\! \mathit{om}}(\mathcal{O}_X \oplus \mathcal{O}_X[1] ,\Delta^! (f_1^\prime)^!\mathcal{O}_X) && \textrm{(\Cref{lem:neeman23lem5_3,lem:neeman187})}
			\\&\cong \mathbf{R}\Delta_\ast \operatorname{\mathbf{R}\mathcal{H}\! \mathit{om}}(\mathcal{O}_X \oplus \mathcal{O}_X[1] ,(f_1^\prime \circ \Delta)^!\mathcal{O}_X) && \textrm{(\Cref{lem:neeman186})}
			\\&\cong \mathbf{R}\Delta_\ast \operatorname{\mathbf{R}\mathcal{H}\! \mathit{om}}(\mathcal{O}_X \oplus \mathcal{O}_X[1] , \mathcal{O}_X) && (f_1^\prime \circ\Delta=1_X)
			\\&\cong \mathbf{R}\Delta_\ast (\operatorname{\mathbf{R}\mathcal{H}\! \mathit{om}}(\mathcal{O}_X, \mathcal{O}_X) \oplus \operatorname{\mathbf{R}\mathcal{H}\! \mathit{om}}(\mathcal{O}_X[1] ,\mathcal{O}_X))
			\\&\cong \mathbf{R}\Delta_\ast (\mathcal{O}_X \oplus \mathcal{O}_X[-1])
			\\&\cong K^\prime.
		\end{aligned}
	\end{displaymath}
	We claim that there exist an adjunction $\Phi_K\dashv\Phi_{K^\prime}$ on $D^b_{\operatorname{coh}}(X)$ in $\mathbf{Triang}$ whose unit is not induced by a morphism of kernels under \eqref{cmp:eq:composite}.
\end{example-claim}

\begin{proof}[Proof of \Cref{cmp:ex:unit}]
	In the proof of \cite[Proposition 2.3]{Canonaco/Stellari:2012}, Canonaco--Stellari construct natural transformations $\zeta_p\colon \operatorname{id} \Rightarrow[1]$ on $D^b_{\operatorname{coh}}(X)$. These are indexed by closed points $p\in |X|$.
	Loc.\ cit.\ used these to show that the functor 
	\begin{displaymath}
		\Phi_\# \colon D^b_{\operatorname{coh}}(X\times_k X) \to \mathbf{Triang}(D^b_{\operatorname{coh}}(X),D^b_{\operatorname{coh}}(X)) , \quad E\mapsto \Phi_E
	\end{displaymath}
	is neither full or faithful.
	Here each $\zeta_p$ satisfies $(\zeta_p)_{\mathcal{O}_p}\neq0$ and $(\zeta_p)_F=0$ for every indecomposable coherent sheaf $F$ not isomorphic to $\mathcal{O}_p$.
	The construction uses that every object of $D^b_{\operatorname{coh}}(X)$ is a finite direct sum of shifts of indecomposable coherent sheaves.
	We bootstrap the dimension argument of \cite[Proposition 2.3]{Canonaco/Stellari:2012}.

	There exist natural isomorphisms of $\Phi_{\mathbf{R}\Delta_\ast\mathcal{O}_X}$ and $\Phi_{\mathbf{R}\Delta_\ast\mathcal{O}_X[1]}$ respectively with $\operatorname{id}$ and $[1]$.
	By \Cref{cmp:def:induced}, there exists a $k$-linear map
	\begin{displaymath}
		\begin{aligned}
			\Theta \colon \operatorname{Hom}(
				\mathbf{R}\Delta_\ast\mathcal{O}_X,
				\mathbf{R}\Delta_\ast\mathcal{O}_X[1])
			&\longrightarrow
			\mathbf{Triang}(
				D^b_{\operatorname{coh}}(X),
				D^b_{\operatorname{coh}}(X))
				(\operatorname{id},[1]),\\
			w&\longmapsto\Phi_w.
		\end{aligned}
	\end{displaymath}
	Its image are the $2$-morphisms $\operatorname{id}\Rightarrow[1]$ induced by morphisms between these diagonal kernels.
	As $\Theta$ is $k$-linear, its image is finite-dimensional.
	Indeed,	apply \cite[Lemma 2.7]{Lank/Venkatesh:2026} with the rank-nullity theorem.

	We show that the $\zeta_p$ are linearly independent over $k$.
	It is indexed by the closed points $p\in |X|$.
	Pick $p_1,\ldots,p_r$ distinct closed points of $|X|$ such that $\sum_{i=1}^r a_i\zeta_{p_i}=0$ for some $a_1,\ldots,a_r\in k$.
	For each $j\in\{1,\ldots,r\}$, the construction of the $\zeta_q$ yields
	%%NOTE: evaluating at $\mathcal{O}_{p_j}$
	\begin{displaymath}
		\sum_{i=1}^r a_i(\zeta_{p_i})_{\mathcal{O}_{p_j}}
		=a_j(\zeta_{p_j})_{\mathcal{O}_{p_j}}=0.
	\end{displaymath}
	Since $(\zeta_{p_j})_{\mathcal{O}_{p_j}}\neq0$, it follows that $a_j=0$.
	As $j$ is arbitrary, the $\zeta_p$ are linearly independent.

	Since $X$ is a curve over an algebraically closed field, $|X|$ has infinitely many closed points.
	Consequently, the finite-dimensional image of $\Theta$ can contain only finitely many of $\zeta_p$.
	Hence, there exists a closed point $p\in |X|$ such that $\zeta_p$ is not induced by any $w\colon\mathbf{R}\Delta_\ast\mathcal{O}_X \to \mathbf{R}\Delta_\ast\mathcal{O}_X[1]$.
	%%NOTE: Although the identity and shift functors are represented by diagonal kernels, the transformation $\zeta_p$ between them is not induced by a morphism between those kernels.

	Recall that the direct sum of the adjunctions $\operatorname{id}\dashv\operatorname{id}$ and $[1]\dashv[-1]$ yields an adjunction $\Phi_K\dashv\Phi_{K^\prime}$.
	Under \eqref{cmp:eq:composite}, its unit $\eta$ has components
	\begin{displaymath}
		\eta_E=(1_E,0,0,1_E)\colon E\to E\oplus E[1]\oplus E[-1]\oplus E,\qquad E\in D^b_{\operatorname{coh}}(X).
	\end{displaymath}
	Thus, $\eta$ is induced by a morphism of kernels.
	Denote by $\tau$ the automorphism of $\Phi_{K^\prime}\cong \operatorname{id} \oplus[-1]$ with components given by
	\begin{displaymath}
		\tau_A=\begin{pmatrix}
			1_A&(\zeta_p)_{A[-1]}\\
			0&1_{A[-1]}
		\end{pmatrix}\colon A\oplus A[-1]\to A\oplus A[-1],\qquad A\in D^b_{\operatorname{coh}}(X).
	\end{displaymath}
	Since $\zeta_p$ is a $2$-morphism of $\mathbf{Triang}$, so is $\tau$.
	Its inverse is obtained by replacing $\zeta_p$ with $-\zeta_p$.
	By \Cref{cmp:lem:units}\eqref{cmp:lem:units1}, $\Phi_K\dashv\Phi_{K^\prime}$ in $\mathbf{Triang}$ with unit and counit
	\begin{displaymath}
		\eta^\prime=(\tau\Phi_K)\circ\eta\quad\textrm{and}\quad\epsilon^\prime=\epsilon\circ(\Phi_K\tau^{-1}).
	\end{displaymath}

	Now we show that the unit $\eta^\prime$ is not induced by kernels.
	For $E\in D^b_{\operatorname{coh}}(X)$, $\Phi_K(E)=E\oplus E[1]$.
	By the construction of $\tau$,
	\begin{displaymath}
		\eta^\prime_E=(1_E,(\zeta_p)_E,0,1_E)\colon E\to E\oplus E[1]\oplus E[-1]\oplus E.
	\end{displaymath}
	Assume that $\eta^\prime=\Phi_u$ via \eqref{cmp:eq:composite} for some morphism
	\begin{displaymath}
		u\colon\mathbf{R}\Delta_\ast\mathcal{O}_X
		\to \mathbf{R}\Delta_\ast\mathcal{O}_X\oplus\mathbf{R}\Delta_\ast\mathcal{O}_X[1] \oplus\mathbf{R}\Delta_\ast\mathcal{O}_X[-1]\oplus\mathbf{R}\Delta_\ast\mathcal{O}_X.
	\end{displaymath}
	By composing $u$ with the projection onto $\mathbf{R}\Delta_\ast\mathcal{O}_X[1]$, we obtain some $w$ with $\Phi_w=\zeta_p$.
	This contradicts the choice of $p$, and hence proves the claim.
\end{proof}

\begin{remark}
	\label{cmp:rem:enhanced}
	Under the given identifications on homotopy categories, the adjunction of \Cref{cmp:ex:unit} does not lift to an adjunction between the corresponding enhanced integral transforms.
	Indeed, $X$ is smooth and proper, and so \cite[Theorem 1.2(2)]{BenZvi/Francis/Nadler:2010} identifies the stable $\infty$-category of perfect kernels on $X\times_k X$ with $k$-linear exact endofunctors of perfect complexes on $X$.
	Since the identification respects composition, the unit of an enhanced adjunction is induced (up to homotopy) by a morphism of kernels.
	However, this is contrary to \Cref{cmp:ex:unit}.

	Additionally, the counit $\epsilon^\prime$ cannot be induced by a morphism of kernels.
	Under $\Phi_K\Phi_{K^\prime}\cong\operatorname{id}\oplus[-1]\oplus[1]\oplus\operatorname{id}$, its components are $\epsilon^\prime_A=(1_A,-(\zeta_p)_{A[-1]},0,1_A)$. 
	A morphism of kernels inducing $\epsilon^\prime$ would also induce $\zeta_p$, and hence this contradicts the choice of $p$.
	However, this does not contradict \cite[Theorems 3.1 and 3.2]{Anno/Logvinenko:2012}.
	Loc.\ cit.\ construct morphisms of kernels inducing the counits.
	In $\mathbf{Triang}$, the adjunctions defined by $(\eta,\epsilon)$ and $(\eta^\prime,\epsilon^\prime)$ are isomorphic.
	Yet, the automorphism $\tau$ relating them is not induced by a morphism of kernels.
	Thus, the obstruction concerns the choice of unit and counit, not the underlying adjoint functors.

	By \Cref{cmp:lem:units}\eqref{cmp:lem:units2}, changing the adjunction data preserves the cones of its unit and counit up to isomorphism.
	Hence, the criteria for full faithfulness and equivalence, as well as \Cref{thm:openness_for_equivalence_fullfaithful_dbcoh}, are independent of this choice.
	Moreover, kernels alone do not determine units, counits, or mates in $\mathbf{Triang}$.
	Hence, the compatibilities in \Cref{prop:pullbacks_is_a_morphism_of_adjoint_integral_transforms,rem:ff_descent} cannot be determined strictly from enhancements.
	See \Cref{sec:mates_in_geometry,sec:fully_faithful_equivalence}.

	We do not know if $\zeta_p$ extends to a $2$-morphism on $D_{\operatorname{qc}}(X)$.
	\Cref{cmp:ex:unit} is formulated on $D^b_{\operatorname{coh}}(X)$, and the adjunctions of \Cref{cor:relatively_perfect_y2_implies_right_adjoint_restrict_to_dbcoh} exist.
\end{remark}

\bibliographystyle{alpha}
\bibliography{mainbib}

@article{Hall/Rydh:2017,
    author = {Hall, Jack and Rydh, David},
    title = {Perfect complexes on algebraic stacks},
    fjournal = {Compositio Mathematica},
    journal = {Compos. Math.},
    issn = {0010-437X},
    volume = {153},
    number = {11},
    pages = {2318--2367},
    year = {2017},
    language = {English},
    doi = {10.1112/S0010437X17007394},
    zbMATH = {6810455},
    Zbl = {1390.14057}
}

@article{Lank/Venkatesh:2026,
    author = {Lank, Pat and Venkatesh, Sridhar},
    title = {Triangulated characterizations of singularities},
    fjournal = {Nagoya Mathematical Journal},
    journal = {Nagoya Math. J.},
    issn = {0027-7630},
    volume = {261},
    pages = {15},
    note = {Id/No e3},
    year = {2026},
    language = {English},
    doi = {10.1017/nmj.2025.11},
    zbMATH = {8135530}
}

@article{AlonsoTarrio/JeremiasLopez/SanchodeSalas:2023,
    author = {{Alonso Tarr{\'{\i}}o}, Leovigildo and {Jerem{\'{\i}}as L{\'o}pez}, Ana and {Sancho de Salas}, Fernando},
    title = {Relative perfect complexes},
    fjournal = {Mathematische Zeitschrift},
    journal = {Math. Z.},
    issn = {0025-5874},
    volume = {304},
    number = {3},
    pages = {25},
    note = {Id/No 42},
    year = {2023},
    language = {English},
    doi = {10.1007/s00209-023-03294-7},
    zbMATH = {7700874},
    Zbl = {1520.18013}
}

@book{MacLane:1978,
    title = {Categories for the Working Mathematician},
    DOI = {10.1007/978-1-4757-4721-8},
    journal = {Graduate Texts in Mathematics},
    publisher = {Springer New York},
    author = {Mac\space{}Lane,  Saunders},
    year = {1978}
}

@article{Lyne:2026,
     author = {Moser, Lyne},
     title = {A double $(\infty ,1)$-categorical nerve for double categories},
     journal = {Annales de l'Institut Fourier},
     year = {2026},
     publisher = {Association des Annales de l{\textquoteright}institut Fourier},
     doi = {10.5802/aif.3750},
     language = {en},
     note = {Online first},
}

@article{Neeman:2026,
    author = {Neeman, Amnon},
    title = {Triangulated categories with a single compact generator, and two {Brown} representability theorems},
    fjournal = {Inventiones Mathematicae},
    journal = {Invent. Math.},
    issn = {0020-9910},
    volume = {244},
    number = {2},
    pages = {531--616},
    year = {2026},
    language = {English},
    doi = {10.1007/s00222-025-01401-5},
    zbMATH = {8188002}
}

@article{Canonaco/Stellari:2012,
    author = {Canonaco, Alberto and Stellari, Paolo},
    title = {Non-uniqueness of {Fourier}--{Mukai} kernels},
    fjournal = {Mathematische Zeitschrift},
    journal = {Math. Z.},
    issn = {0025-5874},
    volume = {272},
    number = {1-2},
    pages = {577--588},
    year = {2012},
    language = {English},
    doi = {10.1007/s00209-011-0950-3},
    zbMATH = {6097118},
    Zbl = {1282.14033}
}

@Article{Mukai:1981,
    Author = {Mukai, Shigeru},
    Title = {Duality between {{\(D(X)\)}} and {{\(D(\hat X)\)}} with its application to {Picard} sheaves},
    FJournal = {Nagoya Mathematical Journal},
    Journal = {Nagoya Math. J.},
    ISSN = {0027-7630},
    Volume = {81},
    Pages = {153--175},
    Year = {1981},
    Language = {English},
    DOI = {10.1017/S002776300001922X},
    zbMATH = {3648900},
    Zbl = {0417.14036}
}

@article{Lieblich:2006,
    author = {Lieblich, Max},
    title = {Moduli of complexes on a proper morphism},
    fjournal = {Journal of Algebraic Geometry},
    journal = {J. Algebr. Geom.},
    issn = {1056-3911},
    volume = {15},
    number = {1},
    pages = {175--206},
    year = {2006},
    language = {English},
    doi = {10.1090/S1056-3911-05-00418-2},
    zbMATH = {5001156},
    Zbl = {1085.14015}
}

@incollection{Kelly/Street:1974,
    author = {Kelly, G. M. and Street, Ross},
    fseries = {Lecture Notes in Mathematics},
    series = {Lect. Notes Math.},
    issn = {0075-8434},
    volume = {420},
    year = {1974},
    publisher = {Springer, Cham},
    language = {English},
    doi = {10.1007/BFb0063096},
    url = {link.springer.com/978-3-540-37270-7},
    zbMATH = {3445176},
    Zbl = {0284.00002},
    title = {Review of the elements of 2-categories},
    booktitle = {Category seminar. Proceedings Sydney category theory seminar 1972/1973}
}

@article{Ben-Zvi/Nadler/Preygel:2017,
    author = {{Ben-Zvi}, David and Nadler, David and Preygel, Anatoly},
    title = {Integral transforms for coherent sheaves},
    fjournal = {Journal of the European Mathematical Society (JEMS)},
    journal = {J. Eur. Math. Soc. (JEMS)},
    issn = {1435-9855},
    volume = {19},
    number = {12},
    pages = {3763--3812},
    year = {2017},
    language = {English},
    doi = {10.4171/JEMS/753},
    zbMATH = {6821388},
    Zbl = {1402.14020}
}

@misc{GuisadoVillalgordo/Lank:2026a,
    title={Full faithful loci for integral transforms on algebraic stacks}, 
    author={El\'{i}as {Guisado Villalgordo} and Pat Lank},
    year={2026},
    url={https://lankp.github.io/assets/PDFs/OlanderCheng.pdf},
    howpublished    = {\href{https://lankp.github.io/assets/PDFs/OlanderCheng.pdf}{Note on author webpage}},
    note = {see also Appendix E of \href{https://arxiv.org/abs/2609.01298v4}{arXiv:2609.01298 v4}}
}

@misc{GuisadoVillalgordo/Lank:2026b,
    title={Remarks on {Fourier--Mukai} invariances}, 
    author={El\'{i}as {Guisado Villalgordo} and Pat Lank},
    year={2026},
    note = {in preparation}
}

@misc{Lank:2026d,
    title={Gluing weak approximability along recollements}, 
    author={Pat Lank},
    year={2026},
    url={https://arxiv.org/abs/2610.22161},
    eprint={2610.22161},
    archivePrefix={arXiv},
    primaryClass={math.AG},
    howpublished    = {\href{https://arxiv.org/abs/2610.22161}{arXiv:2610.22161}},
    publisher     = {arXiv}
}

@article{BenZvi/Francis/Nadler:2010,
    author = {{Ben-Zvi}, David and Francis, John and Nadler, David},
    title = {Integral transforms and {Drinfeld} centers in derived algebraic geometry},
    fjournal = {Journal of the American Mathematical Society},
    journal = {J. Am. Math. Soc.},
    issn = {0894-0347},
    volume = {23},
    number = {4},
    pages = {909--966},
    year = {2010},
    language = {English},
    doi = {10.1090/S0894-0347-10-00669-7},
    zbMATH = {5797887},
    Zbl = {1202.14015}
}

@article{Iyengar/Lipman/Neeman:2015,
    author = {Iyengar, Srikanth B. and Lipman, Joseph and Neeman, Amnon},
    title = {Relation between two twisted inverse image pseudofunctors in duality theory},
    fjournal = {Compositio Mathematica},
    journal = {Compos. Math.},
    issn = {0010-437X},
    volume = {151},
    number = {4},
    pages = {735--764},
    year = {2015},
    language = {English},
    doi = {10.1112/S0010437X14007672},
    zbMATH = {6437571},
    Zbl = {1348.13022}
}

@Article{Burke/Neeman/Pauwels:2023,
    Author = {Burke, Jesse and Neeman, Amnon and Pauwels, Bregje},
    Title = {Gluing approximable triangulated categories},
    FJournal = {Forum of Mathematics, Sigma},
    Journal = {Forum Math. Sigma},
    ISSN = {2050-5094},
    Volume = {11},
    Pages = {18},
    Note = {Id/No e110},
    Year = {2023},
    Language = {English},
    DOI = {10.1017/fms.2023.97},
    zbMATH = {7781650},
    Zbl = {1536.18009}
}

@incollection{Krause:2010,
    author = {Krause, Henning},
    title = {Localization theory for triangulated categories},
    booktitle = {Triangulated categories. Based on a workshop, Leeds, UK, August 2006},
    isbn = {978-0-521-74431-7},
    pages = {161--235},
    year = {2010},
    publisher = {Cambridge: Cambridge University Press},
    language = {English},
    zbMATH = {5831327},
    Zbl = {1232.18012}
}

@article{Cline/Parshall/Scott:1988b,
    author = {Cline, E. and Parshall, B. and Scott, L.},
    title = {Algebraic stratification in representation categories},
    fjournal = {Journal of Algebra},
    journal = {J. Algebra},
    issn = {0021-8693},
    volume = {117},
    number = {2},
    pages = {504--521},
    year = {1988},
    language = {English},
    doi = {10.1016/0021-8693(88)90123-8},
    zbMATH = {4077518},
    Zbl = {0659.18011}
}

@book{Gabriel/Zisman:1967,
    author = {Gabriel, P. and Zisman, M.},
    title = {Calculus of fractions and homotopy theory.},
    fseries = {Ergebnisse der Mathematik und ihrer Grenzgebiete},
    series = {Ergeb. Math. Grenzgeb.},
    volume = {35},
    year = {1967},
    publisher = {Springer-Verlag, Berlin},
    language = {English},
    zbMATH = {3297895},
    Zbl = {0186.56802}
}

@article{Neeman:2024,
    author = {Neeman, Amnon},
    title = {Bounded $t$-structures on the category of perfect complexes},
    fjournal = {Acta Mathematica},
    journal = {Acta Math.},
    issn = {0001-5962},
    volume = {233},
    number = {2},
    pages = {239--284},
    year = {2024},
    language = {English},
    doi = {10.4310/ACTA.2024.v233.n2.a2},
    zbMATH = {7983480}
}

@article {Hall/Rydh:2023,
    AUTHOR = {Hall, Jack and Rydh, David},
    TITLE = {Mayer--{V}ietoris squares in algebraic geometry},
    JOURNAL = {J. Lond. Math. Soc. (2)},
    FJOURNAL = {Journal of the London Mathematical Society. Second Series},
    VOLUME = {107},
    YEAR = {2023},
    NUMBER = {5},
    PAGES = {1583--1612},
    ISSN = {0024-6107,1469-7750},
    MRCLASS = {14A20 (13F40 14F08 14F20)},
    MRNUMBER = {4585296},
    MRREVIEWER = {Ariyan\ Javanpeykar},
    DOI = {10.1112/jlms.12575},
    URL = {https://doi.org/10.1112/jlms.12575},
}

@article{Cline/Parshall/Scott:1988a,
    author = {Cline, E. and Parshall, B. and Scott, L.},
    title = {Finite dimensional algebras and highest weight categories},
    fjournal = {Journal f{\"u}r die Reine und Angewandte Mathematik},
    journal = {J. Reine Angew. Math.},
    issn = {0075-4102},
    volume = {391},
    pages = {85--99},
    year = {1988},
    language = {English},
    url = {https://eudml.org/doc/153071},
    zbMATH = {4073249},
    Zbl = {0657.18005}
}

@article{Gao/Psaroudakis:2018,
    author = {Gao, Nan and Psaroudakis, Chrysostomos},
    title = {Ladders of compactly generated triangulated categories and preprojective algebras},
    fjournal = {Applied Categorical Structures},
    journal = {Appl. Categ. Struct.},
    issn = {0927-2852},
    volume = {26},
    number = {4},
    pages = {657--679},
    year = {2018},
    language = {English},
    doi = {10.1007/s10485-017-9508-9},
    zbMATH = {6916413},
    Zbl = {1428.18021}
}

@article{Schwede/Shipley:2003,
    author = {Schwede, Stefan and Shipley, Brooke},
    title = {Stable model categories are categories of modules},
    fjournal = {Topology},
    journal = {Topology},
    issn = {0040-9383},
    volume = {42},
    number = {1},
    pages = {103--153},
    year = {2003},
    language = {English},
    doi = {10.1016/S0040-9383(02)00006-X},
    zbMATH = {1895949},
    Zbl = {1013.55005}
}

@Article{Letz:2021,
    Author = {Letz, Janina C.},
    Title = {Local to global principles for generation time over commutative {Noetherian} rings},
    FJournal = {Homology, Homotopy and Applications},
    Journal = {Homology Homotopy Appl.},
    ISSN = {1532-0073},
    Volume = {23},
    Number = {2},
    Pages = {165--182},
    Year = {2021},
    Language = {English},
    DOI = {10.4310/HHA.2021.v23.n2.a10},
    zbMATH = {7420249},
    Zbl = {1485.18021}
}

@Article{Neeman:1996,
    Author = {Neeman, Amnon},
    Title = {The {Grothendieck} duality theorem via {Bousfield}'s techniques and {Brown} representability},
    FJournal = {Journal of the American Mathematical Society},
    Journal = {J. Am. Math. Soc.},
    ISSN = {0894-0347},
    Volume = {9},
    Number = {1},
    Pages = {205--236},
    Year = {1996},
    Language = {English},
    DOI = {10.1090/S0894-0347-96-00174-9},
    zbMATH = {868352},
    Zbl = {0864.14008}
}

@article {Rouquier:2008,
    AUTHOR = {Rouquier, Rapha\"el},
    TITLE = {Dimensions of triangulated categories},
    JOURNAL = {J. K-Theory},
    FJOURNAL = {Journal of K-Theory. K-Theory and its Applications in
            Algebra, Geometry, Analysis &amp; Topology},
    VOLUME = {1},
    YEAR = {2008},
    PAGES = {193--256},
    ISSN = {1865-2433},
    MRCLASS = {18E30 (14F05)},
    DOI = {10.1017/is007011012jkt010},
    URL = {https://doi.org/10.1017/is007011012jkt010},
}

@article{Hall/Priver:2024,
    author = {Hall, Jack and Priver, Kyle},
    title = {A generalized {Bondal}--{Orlov} full faithfulness criterion for {Deligne}--{Mumford} stacks},
    fjournal = {Selecta Mathematica. New Series},
    journal = {Sel. Math., New Ser.},
    issn = {1022-1824},
    volume = {32},
    number = {3},
    pages = {47},
    note = {Id/No 51},
    year = {2026},
    language = {English},
    doi = {10.1007/s00029-026-01155-9},
    zbMATH = {8211357}
}

@book{Kashiwara/Schapira:2006,
    author    = {Kashiwara, Masaki and Schapira, Pierre},
    title     = {Categories and Sheaves},
    series    = {Grundlehren der mathematischen Wissenschaften},
    volume    = {332},
    publisher = {Springer},
    address   = {Berlin},
    year      = {2006},
    isbn      = {978-3-540-27949-5},
    doi       = {10.1007/3-540-27950-4}
}

@article{Weber:2007,
    author  = {Mark Weber},
    title   = {Yoneda structures from 2-toposes},
    journal = {Applied Categorical Structures},
    year    = {2007},
    volume  = {15},
    number  = {3},
    pages   = {259--323},
    doi     = {10.1007/s10485-007-9079-2}
}

@article{Shulman:2011,
    author = {Shulman, Michael},
    journal = {The New York Journal of Mathematics [electronic only]},
    language = {eng},
    pages = {75-125},
    publisher = {University at Albany, Deptartment of Mathematics and Statistics},
    title = {Comparing composites of left and right derived functors},
    volume = {17},
    year = {2011}
}

@article{Bokstedt/Neeman:1993,
     author = {B\"okstedt, Marcel and Neeman, Amnon},
     title = {Homotopy limits in triangulated categories},
     journal = {Compositio Mathematica},
     pages = {209--234},
     year = {1993},
     publisher = {Kluwer Academic Publishers},
     volume = {86},
     number = {2},
     mrnumber = {1214458},
     zbl = {0802.18008},
     language = {en},
     url = {https://www.numdam.org/item/CM_1993__86_2_209_0/}
}

@article{Neeman:1992,
    author = {Neeman, Amnon},
    title = {The chromatic tower for {{\(D(R)\)}}. {With} an appendix by {Marcel} {B{\"o}kstedt}},
    fjournal = {Topology},
    journal = {Topology},
    issn = {0040-9383},
    volume = {31},
    number = {3},
    pages = {519--532},
    year = {1992},
    language = {English},
    doi = {10.1016/0040-9383(92)90047-L},
    zbMATH = {166187},
    Zbl = {0793.18008}
}

@article{Raoult:1971,
    author = {Raoult, Jean-Claude},
    title = {Compactification des espaces alg{\'e}briques normaux. ({Compactification} of algebraic normal spaces)},
    fjournal = {Comptes Rendus Hebdomadaires des S{\'e}ances de l'Acad{\'e}mie des Sciences, S{\'e}rie A},
    journal = {C. R. Acad. Sci., Paris, S{\'e}r. A},
    issn = {0366-6034},
    volume = {273},
    pages = {766--767},
    year = {1971},
    language = {French},
    zbMATH = {3358629},
    Zbl = {0226.14002}
}

@book{Gortz/Wedhorn:2020,
    author = {G{\"o}rtz, Ulrich and Wedhorn, Torsten},
    title = {Algebraic geometry {I}. {Schemes}. {With} examples and exercises},
    edition = {2nd edition},
    fseries = {Springer Studium Mathematik -- Master},
    series = {Springer Stud. Math. -- Master},
    issn = {2509-9310},
    isbn = {978-3-658-30732-5; 978-3-658-30733-2},
    year = {2020},
    publisher = {Wiesbaden: Springer Spektrum},
    language = {English},
    doi = {10.1007/978-3-658-30733-2},
    zbMATH = {7261253},
    Zbl = {1444.14001}
}

@book{Gortz/Wedhorn:2023,
    author = {G{\"o}rtz, Ulrich and Wedhorn, Torsten},
    title = {Algebraic geometry {II}: cohomology of schemes. {With} examples and exercises},
    fseries = {Springer Studium Mathematik -- Master},
    series = {Springer Stud. Math. -- Master},
    issn = {2509-9310},
    isbn = {978-3-658-43030-6; 978-3-658-43031-3},
    year = {2023},
    publisher = {Wiesbaden: Springer Spektrum},
    language = {English},
    doi = {10.1007/978-3-658-43031-3},
    zbMATH = {7802900}
}

@article{Ruiperez/Martin/deSalas:2007,
    author = {{Hern{\'a}ndez  Ruip{\'e}rez}, Daniel and {L{\'o}pez Mart{\'{\i}}n}, Ana Cristina and {Sancho de Salas}, Fernando},
    title = {Fourier--{Mukai} transforms for {Gorenstein} schemes},
    fjournal = {Advances in Mathematics},
    journal = {Adv. Math.},
    issn = {0001-8708},
    volume = {211},
    number = {2},
    pages = {594--620},
    year = {2007},
    language = {English},
    doi = {10.1016/j.aim.2006.09.006},
    zbMATH = {5148885},
    Zbl = {1118.14022}
}

@article{Bondal/VandenBergh:2003,
    author = {Bondal, A. and {Van den Bergh}, M.},
    title = {Generators and representability of functors in commutative and noncommutative geometry},
    fjournal = {Moscow Mathematical Journal},
    journal = {Mosc. Math. J.},
    issn = {1609-3321},
    volume = {3},
    number = {1},
    pages = {1--36},
    year = {2003},
    language = {English},
    zbMATH = {2069670},
    Zbl = {1135.18302}
}

@incollection{Deligne:1971,
    author    = {Deligne, P.},
    title     = {Cohomologie \`{a} supports propres},
    booktitle = {Th\'{e}orie de Topos et Cohomologie {\'{E}}tale des Sch\'{e}mas ({SGA}4)},
    volume    = {3},
    publisher = {Springer},
    year      = {1971},
    series    = {Lecture Notes in Mathematics},
    pages     = {250--461},
    note = {Exp. XVII}
}

@incollection{Lipman:2009,
    author = {Lipman, Joseph},
    title = {Notes on {D}erived {F}unctors and {G}rothendieck {D}uality},
    booktitle = {Foundations of {Grothendieck} duality for diagrams of schemes},
    fseries = {Lecture Notes in Mathematics},
    series = {Lect. Notes Math.},
    issn = {0075-8434},
    volume = {1960},
    isbn = {978-3-540-85419-7; 978-3-540-85420-3},
    year = {2009},
    publisher = {Berlin: Springer},
    language = {English},
    doi = {10.1007/978-3-540-85420-3},
    zbMATH = {5490400},
    Zbl = {1163.14001}
}

@book{Godement:1958,
    author    = {Godement, Roger},
    title     = {Topologie alg{\'e}brique et th{\'e}orie des faisceaux},
    publisher = {Hermann},
    year      = {1958},
    series    = {Actualit{\'e}s scientifiques et industrielles},
    number    = {1252},
    address   = {Paris},
    isbn      = {9782705612528}
}

@article{Neeman:2023,
    author = {Neeman, Amnon},
    title = {An improvement on the base-change theorem and the functor {{\(f^!\)}}},
    fjournal = {Bulletin of the Iranian Mathematical Society},
    journal = {Bull. Iran. Math. Soc.},
    issn = {1017-060X},
    volume = {49},
    number = {3},
    pages = {163},
    note = {Id/No 25},
    year = {2023},
    language = {English},
    doi = {10.1007/s41980-023-00768-6},
    zbMATH = {7700381},
    Zbl = {1527.14038}
}

@misc{DeDeyn/Lank/ManaliRahul/Peng:2025,
    title={Categorical characterizations of regularity for algebraic stacks}, 
    author={Timothy {De Deyn} and Pat Lank and Kabeer {Manali Rahul} and Fei Peng},
    year={2025},
    url={https://arxiv.org/abs/2504.02813},
    eprint={2504.02813},
    archivePrefix={arXiv},
    primaryClass={math.AG},
    howpublished    = {\href{https://arxiv.org/abs/2504.02813}{arXiv:2504.02813}},
    publisher     = {arXiv},
}

@misc{GuisadoVillalgordo/Lank/ManaliRahul/Pavic:2025,
    title={Fiberwise criteria for {Fourier}--{Mukai} equivalences}, 
    author={El\'{i}as {Guisado Villalgordo} and Pat Lank and Kabeer {Manali Rahul} and Nebojsa Pavic},
    year={2025},
    url={https://arxiv.org/abs/2512.16503},
    eprint={2512.16503},
    archivePrefix={arXiv},
    primaryClass={math.AG},
    howpublished    = {\href{https://arxiv.org/abs/2512.16503}{arXiv:2512.16503}},
    publisher     = {arXiv},
}

@misc{DeDeyn/Lank/ManaliRahul:2025b,
    title={A note on quasi-perfect morphisms}, 
    author={Timothy {De Deyn} and Pat Lank and Kabeer {Manali Rahul}},
    year={2025},
    url={https://arxiv.org/abs/2508.10845},
    eprint={2508.10845},
    archivePrefix={arXiv},
    primaryClass={math.AG},
    howpublished    = {\href{https://arxiv.org/abs/2508.10845}{arXiv:2508.10845}},
    publisher     = {arXiv},
}

@article{Lipman/Neeman:2007,
    author = {Lipman, Joseph and Neeman, Amnon},
    title = {Quasi-perfect scheme-maps and boundedness of the twisted inverse image functor},
    fjournal = {Illinois Journal of Mathematics},
    journal = {Ill. J. Math.},
    issn = {0019-2082},
    volume = {51},
    number = {1},
    pages = {209--236},
    year = {2007},
    language = {English},
    zbMATH = {5197677},
    Zbl = {1124.14003}
}

@misc{StacksProject,
    shorthand    = {Stacks},
    author       = {The {Stacks Project Authors}},
    title        = {\textit{Stacks Project}},
    howpublished = {\url{https://stacks.math.columbia.edu}},
    year         = {2026},
}

@misc{Lank:2026a,
    title={Perfect generation for regular algebraic stacks}, 
    author={Pat Lank},
    year={2026},
    url={https://arxiv.org/abs/2601.04053},
    eprint={2601.04053},
    archivePrefix={arXiv},
    primaryClass={math.AG},
    howpublished    = {\href{https://arxiv.org/abs/2601.04053}{arXiv:2601.04053}},
    publisher     = {arXiv},
}

@misc{Dutta/Lank/ManaliRahul:2025,
    title={Integral transforms on singularity categories for Noetherian schemes}, 
    author={Uttaran Dutta and Pat Lank and Kabeer {Manali Rahul}},
    year={2025},
    eprint={2501.13834},
    archivePrefix={arXiv},
    primaryClass={math.AG},
    url={https://arxiv.org/abs/2501.13834}, 
    howpublished	= {\href{https://arxiv.org/abs/2501.13834}{arXiv:2501.13834}},
    publisher     = {arXiv},
    note = {to appear in Michigan Math. J.}
}

@misc{Lank:2026c,
    title={Dualizing complexes and {{\(t\)}}-structures for algebraic spaces}, 
    author={Pat Lank},
    year={2026},
    url={https://arxiv.org/abs/2602.20742},
    eprint={2602.20742},
    archivePrefix={arXiv},
    primaryClass={math.AG},
    howpublished    = {\href{https://arxiv.org/abs/2602.20742}{arXiv:2602.20742}},
    publisher     = {arXiv},
}

@article{Orlov:2002,
    author = {Orlov, D. O.},
    title = {Derived categories of coherent sheaves on abelian varieties and equivalences between them},
    fjournal = {Izvestiya: Mathematics},
    journal = {Izv. Math.},
    issn = {1064-5632},
    volume = {66},
    number = {3},
    pages = {569--594},
    year = {2002},
    language = {English},
    doi = {10.1070/IM2002v066n03ABEH000389},
    zbMATH = {1930697},
    Zbl = {1031.18007}
}

@article{Cheng/Gurski/Riehl:2014,
    pages = {337--396},
    volume = {13},
    doi = {10.1017/is013012007jkt250},
    number = {02},
    year = {2014},
    publisher = {Cambridge University Press},
    title = {Cyclic multicategories, multivariable adjunctions and mates},
    journal = {Journal of K-Theory},
    author = {Eugenia Cheng and Nick Gurski and Emily Riehl},
}

@article{Raoult:1974,
    author = {Raoult, Jean-Claude},
    title = {Compactification des espaces alg{\'e}briques},
    fjournal = {Comptes Rendus Hebdomadaires des S{\'e}ances de l'Acad{\'e}mie des Sciences, S{\'e}rie A},
    journal = {C. R. Acad. Sci., Paris, S{\'e}r. A},
    issn = {0366-6034},
    volume = {278},
    pages = {867--869},
    year = {1974},
    language = {French},
    zbMATH = {3437306},
    Zbl = {0278.14002}
}

@article{Conrad/Lieblich/Olsson:2012,
    author = {Conrad, Brian and Lieblich, Max and Olsson, Martin},
    title = {Nagata compactification for algebraic spaces},
    fjournal = {Journal of the Institute of Mathematics of Jussieu},
    journal = {J. Inst. Math. Jussieu},
    issn = {1474-7480},
    volume = {11},
    number = {4},
    pages = {747--814},
    year = {2012},
    language = {English},
    doi = {10.1017/S1474748011000223},
    zbMATH = {6101597},
    Zbl = {1255.14003}
}

@article{Avramov/Iyengar/Lipman:2011,
    author = {Avramov, Luchezar and Iyengar, Srikanth B. and Lipman, Joseph},
    title = {Reflexivity and rigidity for complexes. {II}: {Schemes}},
    fjournal = {Algebra \& Number Theory},
    journal = {Algebra Number Theory},
    issn = {1937-0652},
    volume = {5},
    number = {3},
    pages = {379--429},
    year = {2011},
    language = {English},
    doi = {10.2140/ant.2011.5.379},
    url = {msp.berkeley.edu/ant/2011/5-3/p04.xhtml},
    zbMATH = {6000310},
    Zbl = {1246.14005}
}

@misc{Ballard:2009,
    title={Equivalences of derived categories of sheaves on quasi-projective schemes}, 
    author={Matthew Robert Ballard},
    year={2009},
    url={https://arXiv.org/abs/0905.3148}, 
    eprint={0905.3148},
    archivePrefix={arXiv},
    howpublished	= {\href{https://arXiv.org/abs/0905.3148}{arXiv:0905.3148}},
    publisher     = {arXiv},
}

@Article{Rizzardo:2017,
    Author = {Rizzardo, Alice},
    Title = {Adjoints to a {Fourier}--{Mukai} functor},
    FJournal = {Advances in Mathematics},
    Journal = {Adv. Math.},
    ISSN = {0001-8708},
    Volume = {322},
    Pages = {83--96},
    Year = {2017},
    Language = {English},
    DOI = {10.1016/j.aim.2017.10.015},
    zbMATH = {6806880},
    Zbl = {1386.13045}
}

@article{Ruiperez/Hernandez/Martin/SanchodeSalas:2009,
    author = {Hern{\'a}ndez\space{}Ruip{\'e}rez, Daniel and {L{\'o}pez Mart{\'{\i}}n}, Ana Cristina and Sancho\space{}de\space{}Salas, Fernando},
    title = {Relative integral functors for singular fibrations and singular partners},
    fjournal = {Journal of the European Mathematical Society (JEMS)},
    journal = {J. Eur. Math. Soc. (JEMS)},
    issn = {1435-9855},
    volume = {11},
    number = {3},
    pages = {597--625},
    year = {2009},
    language = {English},
    doi = {10.4171/JEMS/162},
    url = {www.ems-ph.org/journals/show_pdf.php?issn=1435-9855&vol=11&iss=3&rank=7},
    zbMATH = {5565393},
    Zbl = {1221.18010}
}

@book{Neeman:2001,
    author = {Neeman, Amnon},
    title = {Triangulated categories},
    fseries = {Annals of Mathematics Studies},
    series = {Ann. Math. Stud.},
    volume = {148},
    isbn = {0-691-08685-0; 0-691-08686-9},
    year = {2001},
    publisher = {Princeton, NJ: Princeton University Press},
    language = {English},
    doi = {10.1515/9781400837212},
    zbMATH = {1573275},
    Zbl = {0974.18008}
}

@misc{Riehl:26,
	title = {Category {T}heory in {C}ontext},
	author = {Riehl, Emily},
	year = {2026},
    note = {Upcoming 2nd ed. Available at author's webpage: \url{https://emilyriehl.github.io/files/context.pdf}}
}

@misc{Canonaco/Haesemeyer/Neeman/Stellari:2024,
    title={The passage among the subcategories of weakly approximable triangulated categories},
    author={Alberto Canonaco and Christian Haesemeyer and Amnon Neeman and Paolo Stellari},
    year={2024},
    url={https://arxiv.org/abs/2402.04605},
    eprint={2402.04605},
    archivePrefix={arXiv},
    howpublished    = {\href{https://arxiv.org/abs/2402.04605}{arXiv:2402.04605}},
    publisher     = {arXiv},
}

@article{Bergh/Schnurer:2020,
    author = {Bergh, Daniel and Schn{\"u}rer, Olaf M.},
    title = {Conservative descent for semi-orthogonal decompositions},
    fjournal = {Advances in Mathematics},
    journal = {Adv. Math.},
    issn = {0001-8708},
    volume = {360},
    pages = {39},
    note = {Id/No 106882},
    year = {2020},
    language = {English},
    doi = {10.1016/j.aim.2019.106882},
    zbMATH = {7146109},
    Zbl = {1453.14048}
}

@article{Anno/Logvinenko:2012,
    author = {Anno, Rina and Logvinenko, Timothy},
    title = {On adjunctions for {Fourier}--{Mukai} transforms},
    fjournal = {Advances in Mathematics},
    journal = {Adv. Math.},
    issn = {0001-8708},
    volume = {231},
    number = {3-4},
    pages = {2069--2115},
    year = {2012},
    language = {English},
    doi = {10.1016/j.aim.2012.06.007},
    zbMATH = {6094108},
    Zbl = {1316.14033}
}

@Book{Liu:2002,
    Author = {Liu, Qing},
    Title = {Algebraic geometry and arithmetic curves},
    FSeries = {Oxford Graduate Texts in Mathematics},
    Series = {Oxf. Grad. Texts Math.},
    Volume = {6},
    ISBN = {0-19-850284-2},
    Year = {2002},
    Publisher = {Oxford: Oxford University Press},
    Language = {English},
    zbMATH = {1748084},
    Zbl = {0996.14005}
}

\end{document}